\documentclass{amsbook}

\usepackage{hyperref}
\usepackage{cleveref}

\crefname{section}{Section}{Sections}
\crefname{chapter}{Chapter}{Chapters}
\crefname{part}{Part}{Parts}
\crefname{corollary}{Corollary}{Corollaries}
\crefname{definition}{Definition}{Definitions}
\crefname{example}{Example}{Examples}
\crefname{explanation}{Explanation}{Explanations}
\crefname{figure}{Figure}{Figures}
\crefname{lemma}{Lemma}{Lemmas}
\crefname{proposition}{Proposition}{Propositions}
\crefname{question}{Question}{Questions}
\crefname{remark}{Remark}{Remarks}
\crefname{theorem}{Theorem}{Theorems}

\crefname{equation}{}{} 
\Crefname{equation}{}{}
\crefname{enumi}{}{} 
\Crefname{enumi}{}{} 
\crefname{enumii}{}{} 
\Crefname{enumii}{}{}

\usepackage{sty/preamble}

\title[Grothendieck Construction and Inverse $K$-Theory]{The Grothendieck Construction of Bipermutative-Indexed Categories and Pseudo Symmetric Inverse $K$-Theory}

\authorinfoDY

\date{\today}
\subjclass[2020]{}

\begin{document}
\frontmatter

\begin{abstract}
  The Grothendieck construction is a fundamental link between indexed categories and opfibrations.  This work is a detailed study of the Grothendieck construction over a small tight bipermutative category in the context of $\Cat$-enriched multicategories, with applications to inverse $K$-theory and pseudo symmetric $E_\infty$-algebras.
  
The ordinary Grothendieck construction over a small category $\C$ is a 2-equivalence that sends a $\C$-indexed category to an opfibration over $\C$.  We show that the Grothendieck construction over a small tight bipermutative category $\cD$ is a pseudo symmetric $\Cat$-multifunctor that is generally not a $\Cat$-multifunctor in the symmetric sense.  When the projection to $\cD$ is taken into account, we prove that the Grothendieck construction over $\cD$ lifts to a non-symmetric $\Cat$-multiequivalence whose codomain is a non-symmetric $\Cat$-multicategory with small permutative opfibrations over $\cD$ as objects.

As applications we show that inverse $K$-theory, from $\Ga$-categories to small permutative categories, is a pseudo symmetric $\Cat$-multifunctor but not a $\Cat$-multifunctor in the symmetric sense.  It follows that inverse $K$-theory preserves algebraic structures parametrized by non-symmetric and pseudo symmetric $\Cat$-multifunctors but not $\Cat$-multifunctors in general.  As a special case, we observe that inverse $K$-theory sends pseudo symmetric $E_\infty$-algebras in $\Ga$-categories to pseudo symmetric $E_\infty$-algebras in small permutative categories.

\end{abstract}

\maketitle

\cleardoublepage
\thispagestyle{empty}
\vspace*{13.5pc}
\begin{center}
To Jacqueline
\end{center}
\cleardoublepage

\tableofcontents

\newcommand{\sect}[1]{\section*{#1}}

\newcommand{\prefacepartNumName}[1]{\Cref{#1}: \nameref{#1}}
\newcommand{\prefacechapNumName}[1]{\medskip\begin{center}\Cref{#1}: \nameref{#1} \end{center}}

\chapter*{Preface}

\chapquote{``That's the worst joke I ever heard."}{Alexander Grothendieck}{Lectures on topos theory in Buffalo, July 12, 1973}

The Grothendieck construction is a fundamental link between indexed categories and opfibrations.  With $\Cat$ denoting the category of small categories, for each small category $\C$, the Grothendieck construction associates to each pseudo functor $F \cn \C \to \Cat$ a category $\groc F$ together with an opfibration over $\C$:
\[\begin{tikzcd}[column sep=large]
\groc F \ar{r}{U_F} & \C.
\end{tikzcd}\]
Extended to strong transformations and modifications, this construction is a 2-equivalence between 2-categories \cite[10.6.16]{johnson-yau}
\begin{equation}\label{groiiequivalence}
\begin{tikzpicture}[xscale=1,yscale=1,baseline={(x1.base)}]
\draw[0cell=1]
(0,0) node (x1) {[\C,\Cat]}
(x1)++(3,0) node (x2) {\Fibop(\C).};
\draw[1cell=.9]  
(x1) edge node {\groc} node[swap] {\sim} (x2);
\end{tikzpicture}
\end{equation}
There are many variations and extensions of the Grothendieck construction in the literature.  See, for example, \cite{beardsley_wong} for an enriched version, \cite{moeller_vas} for a monoidal version, \cite{shulman-framed} for a version involving double categories, and \cite{hermida-fibrations} where $F$ is a multifunctor.

\sect{Purpose and Main Results}

This work is a detailed study of the Grothendieck construction over a bipermutative category $(\cD,\oplus,\otimes)$ in the context of $\Cat$-enriched multicategories, with applications to inverse $K$-theory and pseudo symmetric $E_\infty$-algebras.  Next we describe the main results of this work.

\subsection*{Characterization of $\Cat$-Multiequivalences}

A key ingredient of the proof of the 2-equivalence $\groc$ in \Cref{groiiequivalence} is a characterization of a 2-equivalence as a 2-functor that is
\begin{itemize}
\item essentially surjective on objects and
\item an isomorphism on each hom category.  
\end{itemize} 
See \cite[7.5.8]{johnson-yau} for a detailed proof.  Just as a multicategory is a multiple-input extension of a category, a $\Cat$-multicategory is a multiple-input generalization of a 2-category.  From this view point, our first main result extends the characterization of 2-equivalences to $\Cat$-\emph{multiequivalences}, which are equivalences in the 2-category $\catmulticat$ of small $\Cat$-multicategories, $\Cat$-multifunctors, and $\Cat$-multinatural transformations.  The Multiequivalence \Cref{thm:multiwhitehead} says that a $\Cat$-multifunctor is a $\Cat$-multiequivalence if and only if it is
\begin{itemize}
\item essentially surjective on objects and
\item an isomorphism on each multimorphism category.
\end{itemize}
The non-symmetric variant, where there is no symmetric group action, is also true.  Moreover, this result recovers the characterization of 2-equivalences by restricting to 1-ary multimorphisms; see \Cref{thm:iiequivalence}.

\subsection*{Pseudo Symmetric $\Cat$-Multifunctors}

We introduce a weaker variant of a $\Cat$-multifunctor, which is called a \emph{pseudo symmetric $\Cat$-multifunctor}, that preserves the symmetric group action up to coherent natural isomorphisms 
\begin{equation}\label{Fsigma-preface}
\begin{tikzpicture}[xscale=3.5,yscale=1.3,vcenter]
\draw[0cell=.9]
(0,0) node (x11) {\M\scmap{\angc;c'}}
(x11)++(1,0) node (x12) {\N\scmap{F\angc;Fc'}}
(x11)++(0,-1) node (x21) {\M\scmap{\angc\sigma;c'}}
(x12)++(0,-1) node (x22) {\N\scmap{F\angc\sigma;Fc'}}
;
\draw[1cell=.9]
(x11) edge node (f1) {F} (x12)
(x21) edge node[swap] (f2) {F} (x22)
(x11) edge node[swap] {\sigma} (x21)
(x12) edge node {\sigma} (x22)
;
\draw[2cell=.9]
node[between=f1 and f2 at .55, 2label={above,F_{\sigma}}, 2label={below,\iso}] {\Longrightarrow}
;
\end{tikzpicture}
\end{equation}
called \emph{pseudo symmetry isomorphisms} \Cref{Fsigmaangccp}.  As we will discuss shortly, pseudo symmetric $\Cat$-multifunctors are important in the context of the Grothendieck construction over a bipermutative category and inverse $K$-theory.  With an appropriate notion of \emph{pseudo symmetric} $\Cat$-multinatural transformations, we prove in \Cref{thm:catmulticatps} that there is a 2-category $\catmulticatps$ with
\begin{itemize}
\item small $\Cat$-multicategories as objects,
\item pseudo symmetric $\Cat$-multifunctors as 1-cells, and
\item pseudo symmetric $\Cat$-multinatural transformations as 2-cells.
\end{itemize}   
In this context, we extend \Cref{thm:multiwhitehead} to characterize pseudo symmetric $\Cat$-multiequivalences; see \Cref{thm:psmultiwhitehead}.

\subsection*{$\Cat$-Multicategory of Bipermutative-Indexed Categories}

A \emph{bipermutative category} $(\cD,\oplus,\otimes)$ is equipped with two permutative---that is, strict symmetric monoidal---structures, $\oplus$ and $\otimes$, that interact via factorization natural transformations
\begin{equation}\label{falfar-preface}
\begin{tikzcd}[row sep=tiny,column sep=huge]
(A \otimes C) \oplus (B \otimes C) \ar{r}{\fal_{A,B,C}} & (A \oplus B) \otimes C\\
(A \otimes B) \oplus (A \otimes C) \ar{r}{\far_{A,B,C}} & A \otimes (B \oplus C)
\end{tikzcd}
\end{equation}
and satisfy appropriate axioms (\Cref{def:embipermutativecat}).  Bipermutative categories arise naturally in multiplicative infinite loop space theory and algebraic $K$-theory spectra; see \cite{elmendorf-mandell} and \cite[Ch.\! 11]{cerberusIII}.  In this work bipermutative categories are mostly small and \emph{tight}, which means that $\fal$ and $\far$ in \Cref{falfar-preface} are natural isomorphisms.  

In the original Grothendieck construction $\groc$ in \Cref{groiiequivalence} with $\C$ a small category, the domain $[\C,\Cat]$ is a 2-category.  For a small tight bipermutative category $(\cD,\oplus,\otimes)$, we prove in \Cref{thm:dcatcatmulticat} that there is a $\Cat$-multicategory 
\[\DCat\]
with symmetric monoidal functors 
\[(\cD,\oplus) \to (\Cat,\times)\]
as objects.  The multimorphism categories of $\DCat$ (\Cref{def:additivenattr,def:additivemodification}) make full use of the bipermutative structure of $\cD$.  We emphasize that the $\Cat$-multicategory structure on $\DCat$ is \emph{not} induced by a monoidal structure (\Cref{expl:dcatbipermzero}).

\subsection*{Pseudo Symmetric $\Cat$-Multifunctorial Grothendieck Construction}

We prove in \Cref{thm:grocatmultifunctor} that the Grothendieck construction over a small tight bipermutative category $(\cD,\oplus,\otimes)$ is a \emph{pseudo symmetric} $\Cat$-multifunctor
\begin{equation}\label{groddcatpermcatsg-preface}
\begin{tikzpicture}[xscale=2.75,yscale=1,baseline={(dc.base)}]
\draw[0cell=1]
(0,0) node (dc) {\DCat}
(dc)++(1,0) node (pc) {\permcatsg.};
\draw[1cell=.9]  
(dc) edge node[pos=.5] {\grod} (pc);
\end{tikzpicture}
\end{equation}
The codomain $\permcatsg$ is a $\Cat$-multicategory with small permutative categories as objects (\Cref{thm:permcatenrmulticat}).  Moreover, if the multiplicative braiding $\betate$ in $\cD$ is the identity natural transformation, then $\grod$ is a $\Cat$-multifunctor, which means that it strictly preserves the symmetric group action.  A weaker converse is true.  If $\grod$ is a $\Cat$-multifunctor, then 
\[x \otimes y = y \otimes x\]
for all objects $x,y \in \cD$.  An important implication of this theorem is that inverse $K$-theory is a pseudo symmetric $\Cat$-multifunctor but \emph{not} a $\Cat$-multifunctor in the symmetric sense.

\subsection*{Permutative Opfibrations from the Grothendieck Construction}

For a symmetric monoidal functor 
\[\begin{tikzcd}[column sep=large]
(\cD,\oplus) \ar{r}{X} & (\Cat,\times),
\end{tikzcd}\] 
the object assignment of $\grod$ in \Cref{groddcatpermcatsg-preface} produces a small permutative category $\grod X$.  On the other hand, the Grothendieck construction $\groc$ in \Cref{groiiequivalence} produces not only the category $\groc F$ but also an opfibration $\groc F \to \C$ over $\C$.  To take into account the projection 
\[\grod X \to \cD\]
to $\cD$, in \Cref{thm:pfibdmulticat} we prove that there is a non-symmetric $\Cat$-multicategory $\pfibd$ with \emph{small permutative opfibrations over $\cD$} as objects.  Similar to $\DCat$, the objects in $\pfibd$ only involve the additive structure $(\cD,\oplus)$, but the multimorphism categories in $\pfibd$ make full use of the bipermutative structure of $\cD$.  On the other hand, unlike $\DCat$, the non-symmetric $\Cat$-multicategory $\pfibd$ is rarely a $\Cat$-multicategory in the symmetric sense (\Cref{pfibdsymwelldefobj,pfibdsymwelldefmor}).  

To extend the 2-equivalence $\groc$ in \Cref{groiiequivalence} to the bipermutative context, first we construct a non-symmetric $\Cat$-multifunctor \Cref{pfibdtopermcatsg}
\[\begin{tikzcd}[column sep=large]
\pfibd \ar{r}{\U} & \permcatsg
\end{tikzcd}\]
that forgets the projection down to $\cD$.  In \Cref{thm:dcatpfibd} we prove that $\grod$ in \Cref{groddcatpermcatsg-preface} lifts as a non-symmetric $\Cat$-multifunctor to $\pfibd$ as in the following commutative diagram.
\begin{equation}\label{grodgrodu}
\begin{tikzpicture}[xscale=2.75,yscale=1.2,vcenter]
\draw[0cell=.9]
(0,0) node (dc) {\DCat}
(dc)++(1,0) node (pc) {\permcatsg}
(pc)++(0,1) node (pf) {\pfibd}
;
\draw[1cell=.9]  
(dc) edge node[pos=.6] {\grod} (pc)
(pf) edge node {\U} (pc)
(dc) edge[bend left=20] node[pos=.7] {\grod} (pf)
;
\end{tikzpicture}
\end{equation}
Utilizing the characterization \Cref{thm:multiwhitehead}, in \Cref{thm:dcatpfibdeq} we prove that the lifted Grothendieck construction
\begin{equation}\label{groddcatpfibd-preface}
\begin{tikzpicture}[xscale=2.75,yscale=1.2,baseline={(dc.base)}]
\draw[0cell=1]
(0,0) node (dc) {\DCat}
(dc)++(1,0) node (pf) {\pfibd};
\draw[1cell=.9]  
(dc) edge node {\grod} node[swap] {\sim} (pf);
\end{tikzpicture}
\end{equation}
is a non-symmetric $\Cat$-\emph{multiequivalence}.  Here is a summary table.
\smallskip
\begin{center}
\resizebox{0.85\width}{!}{
{\renewcommand{\arraystretch}{1.4}%
{\setlength{\tabcolsep}{1ex}
\begin{tabular}{|c|c|c|}\hline
& structure & ref \\ \hline
$\cD$ & small tight bipermutative category & \ref{def:embipermutativecat}\\ \hline
$\DCat$ & $\Cat$-multicategory & \ref{thm:dcatcatmulticat} \\ \hline
$\permcatsg$ & $\Cat$-multicategory & \ref{thm:permcatenrmulticat} \\ \hline
$\pfibd$ & non-symmetric $\Cat$-multicategory & \ref{thm:pfibdmulticat} \\ \hline\hline
$\grod \cn \DCat \to \permcatsg$ & pseudo symmetric $\Cat$-multifunctor & \ref{thm:grocatmultifunctor} \\ \hline
$\U \cn \pfibd \to \permcatsg$ & non-symmetric $\Cat$-multifunctor & \ref{cor:pfibdtopermcatsg} \\ \hline
$\grod \cn \DCat \fto{\sim} \pfibd$ & non-symmetric $\Cat$-multiequivalence & \ref{thm:dcatpfibdeq} \\ \hline
\end{tabular}}}}
\end{center}
\medskip

\subsection*{Pseudo Symmetric $\Cat$-Multifunctorial Inverse $K$-Theory}

The theorems described above have major implications in inverse $K$-theory.  In \cite{mandell_inverseK} Mandell constructed a functor
\begin{equation}\label{invKfunctor}
\begin{tikzcd}[column sep=large]
\Gacat \ar{r}{\cP} & \permcatsg,
\end{tikzcd}
\end{equation}
from the category of $\Ga$-categories to the category of small permutative categories, and proved that $\cP$ is a homotopy inverse of Segal $K$-theory.  It is proved in \cite{johnson-yau-invK} that the inverse $K$-theory functor $\cP$ is a non-symmetric $\Cat$-multifunctor.  We emphasize that $\cP$ is \emph{not} a $\Cat$-multifunctor in the symmetric sense because it does not strictly preserve the symmetric group action.  

By construction $\cP$ is a composite $\groa \circ A$ of two functors.  We prove in \Cref{thm:Acatmultifunctor} that the first functor $A$ extends to a $\Cat$-multifunctor
\begin{equation}\label{catmultifunctorA-preface}
\begin{tikzcd}[column sep=large]
\Gacat \ar{r}{A} & \ACat.
\end{tikzcd}
\end{equation}
The codomain $\ACat$ is the $\Cat$-multicategory $\DCat$ discussed above for the small tight bipermutative category $\cA$ in \Cref{ex:mandellcategory}.  Moreover, by \Cref{groddcatpermcatsg-preface} the Grothendieck construction $\groa$ is a pseudo symmetric $\Cat$-multifunctor.  As a result, inverse $K$-theory
$\cP$
\begin{equation}\label{factorP}
\begin{tikzpicture}[xscale=2.75,yscale=1,vcenter]
\draw[0cell=.9]
(0,0) node (x11) {\Gacat}
(x11)++(1,0) node (x12) {\permcatsg}
(x11)++(.42,-1) node (x21) {\ACat}
;
\draw[1cell=.9] 
(x11) edge node {\cP} (x12)
(x11) edge node[swap,pos=.3] {A} (x21)
(x21) edge node[swap,pos=.7] {\groa} (x12)
;
\end{tikzpicture}
\end{equation}
is a pseudo symmetric $\Cat$-multifunctor but \emph{not} a $\Cat$-multifunctor; see \Cref{thm:invKpseudosym}.  Combining \Cref{grodgrodu,groddcatpfibd-preface,factorP} yields the further factorization of $\cP$
\begin{equation}\label{furtherfactorP}
\begin{tikzpicture}[xscale=3.3,yscale=1,vcenter]
\def\h{.1} \def\v{-1.2}
\draw[0cell=.9]
(0,0) node (x11) {\Gacat}
(x11)++(1,0) node (x12) {\permcatsg}
(x11)++(\h,\v) node (x21) {\ACat}
(x12)++(-\h,\v) node (x22) {\pfiba}
;
\draw[1cell=.85] 
(x11) edge[bend left=15] node {\cP} (x12)
(x11) edge node[swap,pos=.3] {A} (x21)
(x21) edge[transform canvas={xshift=-1ex}] node[pos=.4] {\groa} (x12)
(x21) edge[bend right=20] node[pos=.6] {\groa} node[swap,pos=.6] {\sim} (x22)
(x22) edge node[swap,pos=.7] {\U} (x12)
;
\end{tikzpicture}
\end{equation}
into
\begin{enumerate}[label=(\roman*)]
\item a $\Cat$-multifunctor $A$,
\item a non-symmetric $\Cat$-multiequivalence $\groa$, and
\item a non-symmetric $\Cat$-multifunctor $\U$.
\end{enumerate}
Here is a summary table.
\smallskip
\begin{center}
\resizebox{0.85\width}{!}{
{\renewcommand{\arraystretch}{1.4}%
{\setlength{\tabcolsep}{1ex}
\begin{tabular}{|c|c|c|}\hline
& structure & ref \\ \hline
$\cA$ & small tight bipermutative category & \ref{ex:mandellcategory}\\ \hline
$\ACat$ & $\Cat$-multicategory & \ref{ex:dcatexamples} \\ \hline
$\pfiba$ & non-symmetric $\Cat$-multicategory & \ref{ex:pfibdcatmulticat} \\ \hline
$\Gacat$ & $\Cat$-multicategory & \ref{thm:Gacatsmc}\\ \hline\hline
$\groa \cn \ACat \to \permcatsg$ & pseudo symmetric $\Cat$-multifunctor & \ref{ex:grodpseudosymmetric}\\ \hline
$\U \cn \pfiba \to \permcatsg$ & non-symmetric $\Cat$-multifunctor & \ref{ex:pfibdcatmulticat} \\ \hline
$\groa \cn \ACat \fto{\sim} \pfiba$ & non-symmetric $\Cat$-multiequivalence & \ref{ex:dcatgrodpfibdeq} \\ \hline
$A \cn \Gacat \to \ACat$ & $\Cat$-multifunctor & \ref{thm:Acatmultifunctor}\\ \hline
$\cP \cn \Gacat \to \permcatsg$ & pseudo symmetric $\Cat$-multifunctor & \ref{thm:invKpseudosym}\\ \hline
\end{tabular}}}}
\end{center}
\medskip

\subsection*{Application: Preservation of Pseudo Symmetric Algebras}

Since inverse $K$-theory $\cP$ is not a $\Cat$-multifunctor, it does not preserve algebraic structures parametrized by $\Cat$-multifunctors in general.  However, the pseudo symmetric $\Cat$-multifunctoriality of $\cP$ \Cref{factorP} implies that it preserves algebraic structures parametrized by
\begin{itemize}
\item non-symmetric $\Cat$-multifunctors (\Cref{cor:Ppreservesnonsym}) and
\item pseudo symmetric $\Cat$-multifunctors (\Cref{cor:Ppreservespseudosym}).
\end{itemize}
As a special case of the latter, we observe that $\cP$ preserves \emph{pseudo symmetric $E_\infty$-algebras} (\Cref{cor:Pefinitypsalg}), which are pseudo symmetric $\Cat$-multifunctors from the categorical Barratt-Eccles operad $\BE$.  Pseudo symmetric $E_\infty$-algebras in $\Gacat$ and $\permcatsg$ are described explicitly in \Cref{sec:apppseudo}.

\sect{Chapter Summaries}

This work consists of the following three parts.
\begin{itemize}
\item \prefacepartNumName{part:multibipermcat}
\item \prefacepartNumName{part:multiequivalence}
\item \prefacepartNumName{part:invK}
\end{itemize}
Below is a brief summary of each chapter.  In the main text, each chapter starts with an introduction that describes the content of that chapter more thoroughly and its connection with other chapters.

\prefacechapNumName{ch:prelim}

To keep the prerequisite to an absolute minimum, in this preliminary chapter we review basic concepts related to monoidal categories, enriched categories, and 2-categories.  We recall the symmetric monoidal closed structure on a diagram category $\DV$ when the domain category $\cD$ and the codomain category $\V$ have the relevant structure.  Due to its importance in later chapters, the special case $\V = \Cat$ is explained further.  More detailed discussion of the material in this chapter is in \cite{johnson-yau,cerberusIII,cerberusI}.

\prefacechapNumName{ch:bipermutative}

This chapter has two main purposes.  First we review the closely related notions of \emph{symmetric bimonoidal categories} and \emph{bipermutative categories}.  We discuss several important examples of small tight bipermutative categories, including the indexing categories $\Fskel$ and $\cA$ for, respectively, $\Ga$-categories and inverse $K$-theory.  The second purpose of this chapter is to review Laplaza's Coherence \Cref{thm:laplaza-coherence-1} for symmetric bimonoidal categories.  Tight bipermutative categories, which have invertible factorization morphisms \Cref{falfar-preface}, are special cases of tight symmetric bimonoidal categories.  Thus \Cref{thm:laplaza-coherence-1} also applies to tight bipermutative categories.  More detailed discussion of the material in this chapter is in \cite{cerberusI,cerberusII}.

\prefacechapNumName{ch:multicat}

This chapter is about enriched multicategories and contains the first main theorem of this work.  First we recall the 2-category $\VMulticat$ of small $\V$-multicategories, $\V$-multifunctors, and $\V$-multinatural transformations for a symmetric monoidal category $\V$.  The most relevant case for this work is $\V = \Cat$.  More detailed discussion of enriched multicategories is in \cite{johnson-yau,cerberusIII,yau-operad}.  The second half of this chapter is devoted to the proof of \Cref{thm:multiwhitehead}, which characterizes (adjoint) $\Cat$-multiequivalences locally.  From this theorem we recover the analogous local characterization of 2-equivalences (\Cref{thm:iiequivalence}).

\prefacechapNumName{ch:pseudosymmetry}

This chapter introduces the notion of \emph{pseudo symmetric $\Cat$-multifunctors}.  They preserve the symmetric group action up to coherent natural isomorphisms \Cref{Fsigma-preface} that satisfy four coherence axioms.  They are required to strictly preserve the colored units and the composition.  In later chapters we show that 
\begin{itemize}
\item the Grothendieck construction $\grod$ for each small tight bipermutative category $\cD$ and
\item inverse $K$-theory $\cP$
\end{itemize}
are pseudo symmetric $\Cat$-multifunctors.  Further lax variants of $\Cat$-multifunctors are discussed in \Cref{expl:laxsmf,qu:multibicat}.  \Cref{thm:catmulticatps} proves the existence of a 2-category $\catmulticatps$ with small $\Cat$-multicategories as objects, pseudo symmetric $\Cat$-multifunctors as 1-cells, and pseudo symmetric $\Cat$-multinatural transformations as 2-cells.  \Cref{thm:psmultiwhitehead} provides a local characterization of (adjoint) equivalences in the 2-category $\catmulticatps$.  This concludes \Cref{part:multibipermcat}.

\prefacechapNumName{ch:diagram}

\Cref{part:multiequivalence} is about the Grothendieck construction $\grod$ over a small tight bipermutative category $\cD$.  In \Cref{sec:diagrampermutative} we first discuss in detail the $\Cat$-multicategory $\Dcat$, which is also denoted $\Dtecat$, of $\cD$-indexed categories for a small permutative category $(\cD,\otimes)$.  The rest of this chapter discusses a different $\Cat$-multicategory $\DCat$, which serves as the domain of $\grod$, for a small tight bipermutative category $(\cD,\oplus,\otimes)$.  The objects of $\DCat$ are symmetric monoidal functors $(\cD,\oplus) \to \Cat$, instead of just $\cD$-indexed categories.  Moreover, the multimorphism categories of $\DCat$ involve the full extent of the bipermutative structure of $\cD$ and Laplaza's Coherence \Cref{thm:laplaza-coherence-1}.  Furthermore, while the $\Cat$-multicategory structure on $\Dtecat$ is induced by its symmetric monoidal closed structure, the $\Cat$-multicategory structure on $\DCat$ is \emph{not} induced by a monoidal structure in general (\Cref{expl:dcatbipermzero}).

\prefacechapNumName{ch:multigro}

\Cref{sec:permcatmulticat} discusses the $\Cat$-multicategories $\permcat$ and its strong variant $\permcatsg$, which serves as the codomain of $\grod$ when the projection down to $\cD$ is not taken into account.  The $\Cat$-multicategory $\permcatsg$ is also the codomain of inverse $K$-theory.  In each of $\permcat$ and $\permcatsg$, the objects are small permutative categories, and the $n$-ary 2-cells are $n$-linear natural transformations.  The $n$-ary 1-cells in $\permcat$ and $\permcatsg$ are, respectively, $n$-linear functors and their strong variants.  The rest of this chapter is devoted to the proof of \Cref{thm:grocatmultifunctor}, which proves that the Grothendieck construction $\grod$ in \Cref{groddcatpermcatsg-preface} is a pseudo symmetric $\Cat$-multifunctor in the sense of \Cref{ch:pseudosymmetry}.  It is, furthermore, a $\Cat$-multifunctor if the multiplicative braiding $\betate$ in $\cD$ is the identity.  Conversely, if $\grod$ is a $\Cat$-multifunctor, then $x \otimes y = y \otimes x$ for all objects $x,y \in \cD$.

\prefacechapNumName{ch:permfib}

This chapter studies the Grothendieck construction $\grod$ on a small tight bipermutative category $\cD$ with the projection down to $\cD$ taken into account.  \Cref{thm:pfibdmulticat} proves the existences of the non-symmetric $\Cat$-multicategory $\pfibd$ with small permutative opfibrations over $\cD$ as objects.  Its multimorphism categories are more subtle and involve the bipermutative structure on $\cD$ and Laplaza's Coherence \Cref{thm:laplaza-coherence-1}.  The non-symmetry of $\pfibd$ in general is once again related to the multiplicative braiding $\betate$ in $\cD$.  It is explained in \Cref{pfibdsymwelldefobj,pfibdsymwelldefmor}.  The second main result of this chapter is \Cref{thm:dcatpfibd}, which lifts the pseudo symmetric $\Cat$-multifunctor $\grod$ to a non-symmetric $\Cat$-multifunctor with codomain $\pfibd$ as in \Cref{grodgrodu}.

\prefacechapNumName{ch:gromultequiv}

The main result of this chapter is \Cref{thm:dcatpfibdeq}, which proves that the lifted Grothendieck construction $\grod$ in \Cref{groddcatpfibd-preface} is a non-symmetric $\Cat$-\emph{multiequivalence}.  This proof uses the local characterization for (adjoint) equivalences in \Cref{thm:multiwhitehead}.  In other words, instead of explicitly constructing an adjoint inverse of $\grod$, we show that $\grod$ is (i) essentially surjective on objects and (ii) an isomorphism on each multimorphism category.  This concludes \Cref{part:multiequivalence}.

\prefacechapNumName{ch:multifunctorA}

\Cref{part:invK} is about
\begin{itemize}
\item the pseudo symmetric $\Cat$-multifunctoriality of inverse $K$-theory $\cP$ and
\item its ramification in terms of preservation of pseudo symmetric algebraic structures.
\end{itemize}
Inverse $K$-theory $\cP = \groa \circ A$ involves two constructions:
\begin{itemize}
\item $A$ in \Cref{catmultifunctorA-preface} and
\item the Grothendieck construction $\groa$ for the small tight bipermutative category $\cA$ in \Cref{ex:mandellcategory}.
\end{itemize}   
This chapter discusses the first construction $A$.  For the domain of $\cP$ and $A$, \Cref{sec:Gammacategories,sec:Gacatmulticat} discuss the symmetric monoidal closed category $\Gacat$ of $\Ga$-categories and its $\Cat$-multicategory structure.  Then we discuss the underlying functor of $A$ and its extension to a $\Cat$-multifunctor.

\prefacechapNumName{ch:invK}

This chapter contains our main \Cref{thm:invKpseudosym} on inverse $K$-theory.  It constructs
\begin{itemize}
\item the pseudo symmetric $\Cat$-multifunctorial inverse $K$-theory $\cP = \groa \circ A$ as in \Cref{factorP} and
\item a further factorization of $\cP$ as in \Cref{furtherfactorP} that involves the non-symmetric $\Cat$-multiequivalence $\groa \cn \ACat \fto{\sim} \pfiba$.
\end{itemize} 
In \Cref{sec:unravelingP} we completely unravel the pseudo symmetric $\Cat$-multifunctorial structure of $\cP$.  In the second half of this section, we discuss related work in the literature \cite{elmendorf-multi-inverse,johnson-yau-invK} about the multifunctoriality of $\cP$.  To illustrate the utility of \Cref{thm:invKpseudosym}, the rest of this chapter discusses the fact that $\cP$ preserves pseudo symmetric $E_\infty$-algebras.  As preparation in \Cref{sec:barratteccles} we discuss the categorical Barratt-Eccles operad $\BE$, its Coherence \Cref{mapfrombarratteccles}, and $E_\infty$-algebras in $\Gacat$.  In \Cref{sec:apppseudo} we define \emph{pseudo symmetric $E_\infty$-algebras} as pseudo symmetric $\Cat$-multifunctors from $\BE$.  Inverse $K$-theory sends pseudo symmetric $E_\infty$-algebras in $\Gacat$ to those in $\permcatsg$ (\Cref{cor:Pefinitypsalg}).  Then we completely unravel these structures in $\Gacat$ and $\permcatsg$.  This concludes \Cref{part:invK}.

\prefacechapNumName{ch:questions}

In this appendix we discuss a number of open questions related to the topics of this work.  They provide further motivation for the main text.

\sect{Audience}

This work is aimed at graduate students and researchers with an interest in category theory, algebraic $K$-theory, and homotopy theory.  Our highly detailed exposition is designed to make this work accessible to a wide audience.

\subsection*{Related  Literature}

The following references provide background and further discussion.
\begin{description}
\item[\textnormal{2-dimensional categories}] 
\cite{johnson-yau}
\item[\textnormal{Monoidal categories and enriched multicategories}]
\cite{cerberusIII,yau-operad}
\item[\textnormal{Bimonoidal and bipermutative categories}]
\cite{cerberusI,cerberusII}
\item[\textnormal{The Grothendieck construction}]
\cite{beardsley_wong,hermida-fibrations,johnson-yau,moeller_vas,shulman-framed}
\item[\textnormal{Algebraic $K$-theory spectra}]
\cite{bohmann_osorno,cerberusIII,johnson-yau-permmult,johnson-yau-Fmulti,johnson-yau-multiK,bousfield_friedlander,may-permutative,quillenKI,segal,thomason,waldhausen}
\item[\textnormal{Multifunctorial $K$-theory}]
\cite{elmendorf-mandell,elmendorf-mandell-perm,cerberusIII}
\item[\textnormal{Inverse $K$-theory}]
\cite{gjo1,mandell_inverseK,johnson-yau-invK}
\end{description}

\mainmatter

\part{Bipermutative Categories, Enriched Multicategories, and Pseudo Symmetry}
\label{part:multibipermcat}

\chapter{Preliminaries on Enriched Categories and 2-Categories}
\label{ch:prelim}
To prepare for later chapters, in this chapter we review some relevant definitions and theorems about
\begin{itemize}
\item monoidal categories (\cref{sec:monoidalcat}),
\item enriched categories (\cref{sec:enrichedcat}),
\item 2-categories (\cref{sec:twocategories}), and
\item the Day convolution for functors between symmetric monoidal categories (\cref{sec:dayconvolution}).
\end{itemize}   
Here is a summary table of some of the definitions and results in this chapter.
\smallskip
\begin{center}
\resizebox{0.8\width}{!}{
{\renewcommand{\arraystretch}{1.4}%
{\setlength{\tabcolsep}{1ex}
\begin{tabular}{|c|c|}\hline
monoidal categories (braided, symmetric, closed) & \cref{def:monoidalcategory} (\ref{def:braidedmoncat}, \ref{def:symmoncat}, \ref{def:closedcat}) \\ \hline
diagram categories & \cref{def:diagramcat} \\ \hline
monoidal functors and natural transformations & \cref{def:monoidalfunctor,def:monoidalnattr} \\ \hline\hline
enriched categories, functors, and natural transformations & \cref{def:enriched-category,def:enriched-functor,def:enriched-natural-transformation} \\ \hline
tensor products of enriched categories & \cref{definition:vtensor-0} \\ \hline\hline
2-categories, 2-functors, and 2-natural transformations & \cref{def:twocategory,def:twofunctor,def:twonaturaltr} \\ \hline
modifications & \cref{def:modification} \\ \hline
adjoint equivalences in a 2-category & \cref{def:equivalences} \\ \hline
2-equivalences & \cref{def:iiequivalence} \\ \hline\hline
Day convolution  & \cref{def:dayconvolution} \\ \hline
symmetric monoidal closed diagram categories & \cref{thm:Day,thm:DV} \\ \hline
permutative-indexed categories & \cref{expl:dcatsmclosed} \\ \hline
\end{tabular}}}}
\end{center}
\smallskip
A detailed discussion of enriched categories and the Day convolution is in \cite{cerberusIII}.  All the necessary concepts of 2-dimensional category theory, including pasting diagrams, for this work are discussed in detail in \cite{johnson-yau}.

\subsection*{Connection with Subsequent Chapters}

The notion of a symmetric monoidal category in \cref{sec:monoidalcat} is central to this work.  One of our main objectives is to study the Grothendieck construction from $\DCat$, whose objects are additive symmetric monoidal functors $\Dplus \to \Cat$, with $(\cD,\oplus,\otimes)$ a small tight bipermutative category (\cref{def:embipermutativecat}).  In a bipermutative category, each of $\oplus$ and $\otimes$ is a permutative structure, which means a strict symmetric monoidal structure.  The codomains of the Grothendieck constructions in \cref{thm:grocatmultifunctor,thm:dcatpfibd} also involve small permutative categories.  The codomain of inverse $K$-theory (\cref{thm:invKpseudosym}) has small permutative categories as objects.

From \cref{ch:multicat} onward we discuss enriched multicategories, especially $\Cat$-multicategories, which are multicategories enriched in the symmetric monoidal category $(\Cat, \times, \boldone)$ of small categories with the Cartesian product.  \cref{sec:enrichedcat,sec:twocategories} on enriched categories and 2-categories serve as motivation for enriched multicategories and the 2-categories that they form.  The language of 2-categories also provides a convenient framework to define (adjoint) $\Cat$-multiequivalences for $\Cat$-multicategories (\cref{def:catmultiequivalence}).  One of our main observations in this work is that the Grothendieck construction on bipermutative-indexed categories is a non-symmetric $\Cat$-multiequivalence (\cref{thm:dcatpfibdeq}).

The Day convolution in \cref{sec:dayconvolution} provides a symmetric monoidal structure, hence also a $\Cat$-multicategory structure, on $\cD$-indexed categories (\cref{thm:permindexedcat}).  This, in turn, allows us to define the $\Cat$-multicategory $\DCat$ (\cref{def:dcatcatmulticat}).  However, $\DCat$ is \emph{not} induced by a symmetric monoidal structure due to the non-existence of a potential monoidal unit in general  (\cref{expl:dcatbipermzero}).

\section{Monoidal Categories}
\label{sec:monoidalcat}

In this section we review the definitions of monoidal categories, monoidal functors, and monoidal natural transformations.  For further discussion of monoidal categories, the reader is referred to  \cite{joyal-street,maclane,cerberusIII,cerberusI,cerberusII}.

\begin{definition}\label{def:universe}
A \emph{Grothendieck universe}\index{universe} is a set\label{notation:universe} $\calu$ with the following properties.
\begin{enumerate}
\item\label{univ1} If $x \in \calu$ and $y \in x$, then $y\in\calu$.
\item\label{univ2} If $x\in\calu$, then $\P(x)\in\calu$, where $\P(x)$ is the set of subsets of $x$.
\item\label{univ3} If $J\in\calu$ and $x_j\in\calu$ for each $j\in J$, then the union $\bigcup_{j\in J} x_j \in \calu$.
\item\label{univ4} $\bbN \in \calu$, where $\bbN$ is the set of finite ordinals.\defmark
\end{enumerate}
\end{definition}

\begin{convention}[Grothendieck Universe]\label{conv:universe}\index{Grothendieck Universe}\index{universe}\index{Axiom of Universes}
We assume Grothendieck's \emph{Axiom of Universes}: every set belongs to some universe.  We fix a universe $\cU$\label{not:universe} and call an element in $\cU$ a \emph{set}.  A subset of $\cU$ is called a \emph{class}.  Any categorical structure is called \emph{small} if its objects form a set.  Whenever necessary we tacitly replace $\cU$ by a larger universe $\cU'$ in which $\cU$ is a set.  For more discussion of universes, see \cite[Section 1.1]{johnson-yau}, \cite[I.6]{maclane}, and \cite{agv,maclane-foundation}.
\end{convention}

\begin{definition}\label{def:monoidalcategory}\index{monoidal category}\index{category!monoidal}
A \emph{monoidal category} is a sextuple 
\[(\C,\otimes,\tu,\alpha,\lambda,\rho)\]
consisting of the following data and axioms.
\begin{itemize}
\item $\C$ is a category.
\item \label{notation:monoidal-product}$\otimes \cn \C \times \C \to \C$ is a functor, which is called the \index{monoidal category!monoidal product}\emph{monoidal product}.
\item \label{not:monoidalunit}$\tu \in \C$ is an object, which is called the \index{monoidal category!monoidal unit}\emph{monoidal unit}.
\item \label{not:associativityiso}\label{not:unitisos}$\alpha$, $\lambda$, and $\rho$ are natural isomorphisms
\[\begin{tikzcd}[column sep=huge,row sep=tiny]
(X \otimes Y) \otimes Z \ar{r}{\alpha_{X,Y,Z}}[swap]{\iso} & X \otimes (Y \otimes Z)
\end{tikzcd}\]\vspace{-1em}
\[\begin{tikzcd}[column sep=large]
\tu \otimes X \ar{r}{\lambda_X}[swap]{\iso} & X & X \otimes \tu \ar{l}{\iso}[swap]{\rho_X}
\end{tikzcd}\]
for objects $X,Y,Z \in \C$.  They are called, respectively, the \index{associativity isomorphism}\index{monoidal category!associativity isomorphism}\emph{associativity isomorphism}, the \index{left unit isomorphism}\index{monoidal category!left unit isomorphism}\emph{left unit isomorphism}, and the \index{right unit isomorphism}\index{monoidal category!right unit isomorphism}\emph{right unit isomorphism}.
\end{itemize}
The following middle unity\index{monoidal category!middle unity axiom} and pentagon\index{monoidal category!pentagon axiom}\index{pentagon axiom} diagrams are required to commute for objects $W,X,Y,Z \in \C$, where all $\otimes$ symbols are omitted to save space.
\begin{equation}\label{monoidalcataxioms}
\begin{tikzpicture}[xscale=1,yscale=1,vcenter]
\tikzset{0cell/.append style={nodes={scale=.85}}}
\tikzset{1cell/.append style={nodes={scale=.85}}}
\def\h{1.6} \def\g{1.5} \def\v{.9} \def\u{2}
\draw[0cell]
(0,0) node (a) {(X \tu) Y}
(a)++(0,-2) node (b) {X (\tu Y)}
(a)++(1.5,-1) node (c) {XY}
;
\draw[1cell]  
(a) edge node[pos=.4] {\rho_X 1_Y} (c)
(a) edge node[swap] {\alpha_{X,\tu,Y}} (b)
(b) edge node[swap,pos=.4] {1_X \lambda_Y} (c)
;
\begin{scope}[shift={(6,0)}]
\draw[0cell]
(0,0) node (x0) {(WX)(YZ)}
(x0)++(-\h,-\v) node (x11) {((WX)Y)Z}
(x0)++(\h,-\v) node (x12) {W(X(YZ))}
(x0)++(-\g,-\u) node (x21) {(W(XY))Z}
(x0)++(\g,-\u) node (x22) {W((XY)Z)}
;
\draw[1cell]
(x11) edge node[pos=.2] {\al_{WX,Y,Z}} (x0)
(x0) edge node[pos=.8] {\al_{W,X,YZ}} (x12)
(x11) edge node[swap,pos=.2] {\al_{W,X,Y} 1_Z} (x21)
(x21) edge node {\al_{W,XY,Z}} (x22)
(x22) edge node[swap,pos=.8] {1_W \al_{X,Y,Z}} (x12)
;
\end{scope}
\end{tikzpicture}
\end{equation}
This finishes the definition of a monoidal category.  Moreover, we define the following.
\begin{itemize}
\item A monoidal category is \emph{strict}\index{monoidal category!strict}\index{strict monoidal category} if $\alpha$, $\lambda$, and $\rho$ are identity natural transformations.
\item A (monoidal) category is also called a \index{monoidal 1-category}\index{1-category}\emph{(monoidal) 1-category}.\defmark
\end{itemize}
\end{definition}

In each monoidal category, the morphism equality holds:\index{monoidal category!unity properties} 
\begin{equation}\label{lambda=rho}
\lambda_{\tensorunit} = \rho_{\tensorunit} \cn \tensorunit \otimes \tensorunit \to \tensorunit.
\end{equation}
Moreover, the following left and right unity diagrams commute \cite{kelly-coherence}.  
\begin{equation}\label{moncat-other-unit-axioms}
\begin{tikzcd}[column sep=normal]
(\tensorunit \otimes X) \otimes Y \dar[swap]{\lambda_X \otimes 1_Y} \rar{\alpha_{\tensorunit,X,Y}}
& \tensorunit \otimes (X\otimes Y) \dar{\lambda_{X\otimes Y}}\\ X \otimes Y \rar[equal]& X \otimes Y\end{tikzcd}\qquad
\begin{tikzcd}[column sep=normal]
(X \otimes Y) \otimes \tensorunit \dar[swap]{\rho_{X \otimes Y}} \rar{\alpha_{X,Y,\tensorunit}}
& X \otimes (Y\otimes \tensorunit ) \dar{1_X \otimes \rho_Y}\\ X \otimes Y \rar[equal]& X \otimes Y\end{tikzcd}
\end{equation} 
Proofs of these unity properties in the more general context of bicategories are in \cite[Section 2.2]{johnson-yau}.

\begin{definition}\label{def:monoid}
A \emph{monoid}\index{monoid} in a monoidal category $\C$ is a triple\label{notation:monoid} $(X,\mu,\eta)$ consisting of the following data and axioms.
\begin{itemize}
\item $X \in \C$ is an object.
\item $\mu \cn X \otimes X \to X$ is a morphism, which is called the \index{multiplication}\emph{multiplication}.
\item $\eta \cn \tensorunit \to X$ is a morphism, which is called the \index{unit}\emph{unit}.
\end{itemize}
The following associativity and unity diagrams are required to commute.
\[\begin{tikzcd}[column sep=large]
(X\otimes X) \otimes X \arrow{dd}[swap]{\mu\otimes 1_X} \rar{\alpha} & X \otimes (X \otimes X) \dar{1_X\otimes \mu}\\ & X \otimes X \dar{\mu}\\  
X \otimes X \arrow{r}{\mu} & X
\end{tikzcd}\qquad
\begin{tikzcd}[column sep=large]
\tensorunit \otimes X \ar{d}[swap]{\eta \otimes 1_X} \ar{r}{\lambda_X} & X \ar[equal]{d}\\ 
X \otimes X \ar{r}{\mu} & X \ar[equal]{d}\\
X \otimes \tensorunit \ar{u}{1_X \otimes \eta} \ar{r}{\rho_X} & X
\end{tikzcd}\]
A \emph{morphism} of monoids
\[\begin{tikzcd}[column sep=large]
(X,\mu^X,\eta^X) \ar{r}{f} & (Y,\mu^Y,\eta^Y)
\end{tikzcd}\] 
is a morphism $f \cn X \to Y$ in $\C$ such that the diagrams
\[\begin{tikzcd}[column sep=large]
X\otimes X \dar[swap]{\mu^X} \rar{f\otimes f} & Y \otimes Y \dar{\mu^Y}\\ X \rar{f} & Y\end{tikzcd} \qquad
\begin{tikzcd}[column sep=large]
\tensorunit \dar[equal] \rar{\eta^X} & X \dar{f}\\ \tensorunit \rar{\eta^Y} & Y\end{tikzcd}\]
commute.
\end{definition}

\begin{definition}\label{def:braidedmoncat}\index{braided monoidal category}\index{monoidal category!braided}\index{category!braided monoidal}
A \emph{braided monoidal category} is a pair $(\C,\xi)$ consisting of the following data and axioms.
\begin{itemize}
\item $\C = (\C,\otimes,\tu,\al,\lambda,\rho)$ is a monoidal category.
\item $\xi$ is a natural isomorphism\label{notation:symmetry-iso} 
\[X \otimes Y \fto[\iso]{\xi_{X,Y}} Y \otimes X\]
for $X,Y \in \C$, which is called the \index{braiding}\emph{braiding}.
\end{itemize}
The hexagon diagrams\index{hexagon diagram}
\begin{equation}\label{hexagon-braided}
\def\sb{\scalebox{.7}}
\begin{tikzpicture}[scale=.75,commutative diagrams/every diagram]
\node (P0) at (0:2cm) {\sb{$Y \otimes (Z\otimes X)$}};
\node (P1) at (60:2cm) {\makebox[3ex][l]{\sb{$Y \otimes (X \otimes Z)$}}};
\node (P2) at (120:2cm) {\makebox[3ex][r]{\sb{$(Y \otimes X) \otimes Z$}}};
\node (P3) at (180:2cm) {\sb{$(X \otimes Y) \otimes Z$}};
\node (P4) at (240:2cm) {\makebox[3ex][r]{\sb{$X \otimes (Y \otimes Z)$}}};
\node (P5) at (300:2cm) {\makebox[3ex][l]{\sb{$(Y \otimes Z) \otimes X$}}};
\path[commutative diagrams/.cd, every arrow, every label]
(P3) edge node[pos=.25] {\sb{$\xitimes_{X,Y}\otimes 1_Z$}} (P2)
(P2) edge node {\sb{$\alpha$}} (P1)
(P1) edge node[pos=.75] {\sb{$1_Y \otimes \xitimes_{X,Z}$}} (P0)
(P3) edge node[swap,pos=.4] {\sb{$\alpha$}} (P4)
(P4) edge node {\sb{$\xitimes_{X,Y \otimes Z}$}} (P5)
(P5) edge node[swap,pos=.6] {\sb{$\alpha$}} (P0);
\end{tikzpicture}
\qquad
\begin{tikzpicture}[scale=.75, commutative diagrams/every diagram]
\node (P0) at (0:2cm) {\sb{$(Z\otimes X) \otimes Y$}};
\node (P1) at (60:2cm) {\makebox[3ex][l]{\sb{$(X \otimes Z) \otimes Y$}}};
\node (P2) at (120:2cm) {\makebox[3ex][r]{\sb{$X \otimes (Z \otimes Y)$}}};
\node (P3) at (180:2cm) {\sb{$X \otimes (Y \otimes Z)$}};
\node (P4) at (240:2cm) {\makebox[3ex][r]{\sb{$(X \otimes Y) \otimes Z$}}};
\node (P5) at (300:2cm) {\makebox[3ex][l]{\sb{$Z \otimes (X \otimes Y)$}}};
\path[commutative diagrams/.cd, every arrow, every label]
(P3) edge node[pos=.25] {\sb{$1_X \otimes \xitimes_{Y,Z}$}} (P2)
(P2) edge node {\sb{$\alpha^\inv$}} (P1)
(P1) edge node[pos=.75] {\sb{$\xitimes_{X,Z}\otimes 1_Y$}} (P0)
(P3) edge node[swap,pos=.4] {\sb{$\alpha^\inv$}} (P4)
(P4) edge node {\sb{$\xi_{X \otimes Y, Z}$}} (P5)
(P5) edge node[swap,pos=.6] {\sb{$\alpha^\inv$}} (P0);
\end{tikzpicture}
\end{equation}
are required to commute for objects $X,Y,Z \in \C$. 
\end{definition}

The following two unity diagrams\index{braided monoidal category!unity properties} commute for each object $X$ in a braided monoidal category; see \cite[1.3.21]{cerberusII} for the proof. 
\begin{equation}\label{braidedunity}
\begin{tikzcd}[column sep=large]
X \otimes \tensorunit \dar[swap]{\rho_X} \rar{\xi_{X,\tensorunit}} & \tensorunit \otimes X \dar{\lambda_X} & \tu \otimes X \dar[swap]{\lambda_X} \rar{\xi_{\tu,X}} & X \otimes \tu \dar{\rho_X}\\ 
X \rar[equal] & X & X \rar[equal] & X
\end{tikzcd}
\end{equation}

\begin{definition}\label{def:symmoncat}\index{symmetric monoidal category}\index{monoidal category!symmetric}\index{category!symmetric monoidal}
A \emph{symmetric monoidal category} is a pair $(\C,\xi)$ consisting of the following data and axioms.
\begin{itemize}
\item $\C = (\C,\otimes,\tu,\al,\lambda,\rho)$ is a monoidal category.
\item $\xi$ is a natural isomorphism 
\[X \otimes Y \fto[\iso]{\xi_{X,Y}} Y \otimes X\]
for $X,Y \in \C$, which is called the \emph{symmetry isomorphism} or the \emph{braiding}.
\end{itemize}
The following symmetry\index{symmetry axiom} and \index{hexagon diagram}hexagon diagrams are required to commute for objects $X,Y,Z \in \C$. 
\begin{equation}\label{symmoncatsymhexagon}
\begin{tikzpicture}[xscale=3,yscale=1.2,vcenter]
\def\h{.1}
\draw[0cell=.8] 
(0,0) node (a) {X \otimes Y}
(a)++(.7,0) node (c) {X \otimes Y}
(a)++(.35,-1) node (b) {Y \otimes X}
;
\draw[1cell=.8] 
(a) edge node {1_{X \otimes Y}} (c)
(a) edge node [swap,pos=.3] {\xi_{X,Y}} (b)
(b) edge node [swap,pos=.7] {\xi_{Y,X}} (c)
;
\begin{scope}[shift={(1.5,.5)}]
\draw[0cell=.8] 
(0,0) node (x11) {(Y \otimes X) \otimes Z}
(x11)++(.9,0) node (x12) {Y \otimes (X \otimes Z)}
(x11)++(-\h,-1) node (x21) {(X \otimes Y) \otimes Z}
(x12)++(\h,-1) node (x22) {Y \otimes (Z \otimes X)}
(x11)++(0,-2) node (x31) {X \otimes (Y \otimes Z)}
(x12)++(0,-2) node (x32) {(Y \otimes Z) \otimes X}
;
\draw[1cell=.8]
(x21) edge node[pos=.25] {\xi_{X,Y} \otimes 1_Z} (x11)
(x11) edge node {\alpha} (x12)
(x12) edge node[pos=.75] {1_Y \otimes \xi_{X,Z}} (x22)
(x21) edge node[swap,pos=.25] {\alpha} (x31)
(x31) edge node {\xi_{X,Y \otimes Z}} (x32)
(x32) edge node[swap,pos=.75] {\alpha} (x22)
;
\end{scope}
\end{tikzpicture}
\end{equation}
A \emph{permutative category}\index{permutative category}\index{category!permutative} is a strict symmetric monoidal category.  For a permutative category, we sometimes denote its monoidal product by \label{not:opluse}$\oplus$ and its monoidal unit by $e$.
\end{definition}

\begin{remark}\label{rk:smcat}
The symmetry axiom, $\xi_{Y,X} \xi_{X,Y} = 1$, implies that the left and right hexagons \cref{hexagon-braided} are equivalent.  Therefore, a symmetric monoidal category is precisely a braided monoidal category that satisfies the symmetry axiom.  As a consequence, in each symmetric monoidal category the unity diagrams \cref{braidedunity} commute.
\end{remark}

\begin{definition}\label{def:closedcat}\index{closed category}\index{category!closed}\index{internal hom}
A symmetric monoidal category $(\C,\otimes,\tu,\xi)$ is \emph{closed} if, for each object $X \in \C$, the functor $- \otimes X \cn \C \to \C$ admits a right adjoint, which is called an \emph{internal hom}.  A right adjoint of $- \otimes X$ is denoted by $\Hom(X,-)$ or \label{notation:internal-hom}$[X,-]$.
\end{definition}

\begin{definition}[Diagrams]\label{def:diagramcat}\index{diagram category}\index{category!diagram}
For a small category $\cD$ and a category $\C$, the \emph{diagram category} \label{not:DC}$\DC$ has 
\begin{itemize}
\item functors $\cD \to \C$ as objects, 
\item natural transformations between such functors as morphisms,
\item identity natural transformations as identities, and
\item vertical composition of natural transformations as composition.  
\end{itemize}
Moreover, we define the following.
\begin{itemize}
\item A functor $\cD \to \C$ is also called a \emph{$\cD$-diagram in $\C$}.
\item For the category $\Cat$ of small categories and functors, a functor $\cD \to \Cat$ is called a \index{indexed category}\index{category!indexed}\emph{$\cD$-indexed category}.\defmark
\end{itemize}
\end{definition}

\begin{example}[Small Categories]\label{ex:cat}\index{category!of small categories}
\label{not:cat}$(\Cat, \times, \boldone, [,])$ is a symmetric monoidal closed category, with $\Cat$ the category of small categories and functors.
\begin{itemize}
\item The monoidal product is the Cartesian product $\times$.  
\item The monoidal unit $\boldone$ is the terminal category with only one object $*$ and its identity morphism $1_*$.  
\item The closed structure $[,]$ is given by diagram categories (\cref{def:diagramcat}). 
\end{itemize} 
The category $\Cat$ is complete and cocomplete \cite[Section 1.4]{yau-involutive}.
\end{example}

\subsection*{Monoidal Functors and Natural Transformations}

\begin{definition}\label{def:monoidalfunctor}\index{monoidal functor}\index{functor!monoidal}
Suppose $\C$ and $\D$ are monoidal categories.  A \emph{monoidal functor} 
\[(F, F^2, F^0) \cn \C \to \D\]
consists of the following data and axioms.
\begin{itemize}
\item $F \cn \C \to \D$ is a functor.
\item $F^0 \cn \tu \to F\tu$ is a morphism in $\D$, which is called the \index{unit constraint}\index{constraint!unit}\emph{unit constraint}.
\item $F^2$ is a natural transformation
\[\begin{tikzcd}[column sep=large]
FX \otimes FY \ar{r}{F^2_{X,Y}} & F(X \otimes Y)
\end{tikzcd}
\forspace X,Y \in \C,\]
which is called the \index{monoidal constraint}\index{constraint!monoidal}\emph{monoidal constraint}.
\end{itemize}
The following unity and associativity diagrams are required to commute for objects $X,Y,Z \in \C$.  
\begin{equation}\label{monoidalfunctorunity}
\begin{tikzcd}
\tensorunit \otimes FX \dar[swap]{F^0 \otimes 1_{FX}} \rar{\lambda_{FX}} & FX \\ 
F\tensorunit \otimes FX \rar{F^2} & F(\tensorunit \otimes X)
\uar[swap]{F\lambda_X}
\end{tikzcd}
\qquad
\begin{tikzcd}
FX \otimes \tensorunit \dar[swap]{1_{FX} \otimes F^0} \rar{\rho_{FX}} & FX \\ 
FX \otimes F\tensorunit \rar{F^2} & F(X \otimes \tensorunit)
\uar[swap]{F\rho_X}
\end{tikzcd}
\end{equation}
\begin{equation}\label{monoidalfunctorassoc}
\begin{tikzcd}[column sep=large]
\bigl(FX \otimes FY\bigr) \otimes FZ \rar{\alpha} \dar[swap]{F^2 \otimes 1_{FZ}} 
& FX \otimes \bigl(FY \otimes FZ\bigr) \dar{1_{FX} \otimes F^2}\\
F(X \otimes Y) \otimes FZ \dar[swap]{F^2} & FX \otimes F(Y \otimes Z) \dar[d]{F^2}\\
F\bigl((X \otimes Y) \otimes Z\bigr) \rar{F\alpha} &
F\bigl(X \otimes (Y \otimes Z)\bigr)
\end{tikzcd}
\end{equation}
A monoidal functor is said to be 
\begin{itemize}
\item \index{monoidal functor!strong}\index{strong monoidal functor}\emph{strong} if $F^0$ and $F^2$ are isomorphisms;
\item \index{monoidal functor!strictly unital}\index{strictly unital monoidal functor}\emph{strictly unital} if $F^0$ is the identity morphism; and
\item \index{monoidal functor!strict}\index{strict monoidal functor}\emph{strict} if $F^0$ and $F^2$ are identities.
\end{itemize}
For braided monoidal categories $\C$ and $\D$, a monoidal functor $(F,F^2,F^0)$ as above is \index{braided monoidal functor}\index{monoidal functor!braided}\emph{braided} if the diagram
\begin{equation}\label{monoidalfunctorbraiding}
\begin{tikzcd}[column sep=large]
FX \otimes FY \dar[swap]{F^2} \rar{\xi_{FX,FY}}[swap]{\cong} & FY \otimes FX \dar{F^2} \\ 
F(X \otimes Y) \rar{F\xi_{X,Y}}[swap]{\cong} & F(Y \otimes X)
\end{tikzcd}
\end{equation}
commutes for all objects $X,Y \in \C$.  If $\C$ and $\D$ are symmetric monoidal categories, then a braided monoidal functor as above is called a \index{symmetric monoidal functor}\index{monoidal functor!symmetric}\emph{symmetric monoidal functor}.
\end{definition}

\begin{definition}\label{def:monoidalnattr}
For monoidal functors 
\[(F, F^2, F^0) \andspace (G,G^2,G^0) \cn \C \to \D\]
between monoidal categories $\C$ and $\D$, a \index{monoidal natural transformation}\index{natural transformation!monoidal}\emph{monoidal natural transformation} $\theta \cn F \to G$ is a natural transformation of the underlying functors such that the unit constraint and monoidal constraint diagrams 
\begin{equation}\label{monnattr}
\begin{tikzpicture}[xscale=2,yscale=1.3,vcenter]
\draw[0cell=.9]
(0,0) node (x11) {\tu}
(x11)++(1,0) node (x12) {F\tu}
(x12)++(0,-1) node (x2) {G\tu}
(x12)++(1.2,0) node (y11) {FX \otimes FY}
(y11)++(1.5,0) node (y12) {GX \otimes GY}
(y11)++(0,-1) node (y21) {F(X \otimes Y)}
(y12)++(0,-1) node (y22) {G(X \otimes Y)}
;
\draw[1cell=.9]  
(x11) edge node {F^0} (x12)
(x11) edge node[swap] {G^0} (x2)
(x12) edge node {\theta_{\tu}} (x2)
(y11) edge node {\theta_X \otimes \theta_Y} (y12)
(y12) edge node {G^2} (y22)
(y11) edge node[swap] {F^2} (y21)
(y21) edge node {\theta_{X \otimes Y}} (y22)
;
\end{tikzpicture}
\end{equation}
commute for all objects $X,Y \in \C$.  
\end{definition}

Such a (monoidal) natural transformation is also denoted by the following \index{2-cell!notation}2-cell notation.
\begin{equation}\label{twocellnotation}
\begin{tikzpicture}[xscale=2,yscale=1.7,baseline={(x1.base)}]
\draw[0cell=.9]
(0,0) node (x1) {\C}
(x1)++(1,0) node (x2) {\D}
;
\draw[1cell=.9]  
(x1) edge[bend left] node {F} (x2)
(x1) edge[bend right] node[swap] {G} (x2)
;
\draw[2cell]
node[between=x1 and x2 at .47, rotate=-90, 2label={above,\theta}] {\Rightarrow}
;
\end{tikzpicture}
\end{equation}
Pasting diagrams involving 2-cells are discussed in \cite[Ch.\! 3]{johnson-yau}.

\begin{convention}[Left Normalized Bracketing]\label{expl:leftbracketing}
Unless otherwise specified, an iterated monoidal product is \index{left normalized}\emph{left normalized} with the left half of each pair of parentheses at the far left.
\end{convention}  

For example, we denote
\[a \otimes b \otimes c \otimes d = \big((a \otimes b) \otimes c\big) \otimes d.\]
We omit parentheses for iterated monoidal products and tacitly insert the necessary associativity and unit isomorphisms.  This is valid because each monoidal category $\C$ admits a strong monoidal adjoint equivalence $\C \to \C_\st$ with $\C_\st$ strict monoidal \cite[XI.3.1]{maclane}.  The braided and symmetric versions are in \cite[21.3.1 and 21.6.1]{yau-inf-operad}.  In each case, the strict diagrams commute if and only if their preimages in $\C$ commute.

\section{Enriched Categories}
\label{sec:enrichedcat}

In this section we briefly review the definitions of enriched categories, enriched functors, enriched natural transformations, and their tensor products.  More detailed discussion is in \cite[Section 1.3]{johnson-yau}, \cite[Ch.\! 1]{cerberusIII}, and \cite{kelly-enriched}. 

\begin{definition}\label{def:enriched-category}
For a monoidal category $(\V,\otimes,\tu,\alpha,\lambda,\rho)$ (\cref{def:monoidalcategory}), a \emph{$\V$-category} $\C$, which is also called a \index{category!enriched}\index{enriched category}\emph{category enriched in $\V$}, consists of the following data and axioms.
\begin{itemize}
\item It is equipped with a class \index{object!enriched category}$\Ob(\C)$ of \emph{objects}.
\item For each pair of objects $X,Y$ in $\C$, it is equipped with a \index{hom object}\index{object!hom}\emph{hom object} $\C(X,Y) \in \V$ with domain $X$ and codomain $Y$.
\item For objects $X,Y,Z$ in $\C$, it is equipped with a morphism\label{not:enrcomposition}
\[\begin{tikzcd}[column sep=huge] 
\C(Y,Z) \otimes \C(X,Y) \rar{m_{X,Y,Z}} & \C(X,Z)
\end{tikzcd} \inspace \V,\]
called the \index{composition!enriched category}\emph{composition}.
\item For each object $X$ in $\C$, it is equipped with a morphism\label{not:enridentity}
\[\begin{tikzcd}[column sep=large] 
\tensorunit \rar{i_X} & \C(X,X)
\end{tikzcd} \inspace \V,\]
called the \emph{identity} of $X$.
\end{itemize}
The following \index{associativity!enriched category}\emph{associativity diagram} and \index{unity!enriched category}\emph{unity diagram} are required to commute for objects $W,X,Y,Z$ in $\C$.
\begin{equation}\label{enriched-cat-associativity}
\begin{tikzcd}[cells={nodes={scale=.9}}]
  \big(\C(Y,Z) \otimes \C(X,Y)\big) \otimes \C(W,X) \arrow{dd}[swap]{m \otimes 1}
  \rar{\alpha}
  & \C(Y,Z) \otimes \big(\C(X,Y) \otimes \C(W,X)\big)  \dar{1 \otimes m} \\
& \C(Y,Z) \otimes \C(W,Y) \dar{m}\\
\C(X,Z) \otimes \C(W,X) \rar{m} & \C(W,Z)  
\end{tikzcd}
\end{equation}
\begin{equation}\label{enriched-cat-unity}
\begin{tikzcd}[cells={nodes={scale=.9}}]
\tensorunit \otimes \C(X,Y) \dar[swap]{i_Y \otimes 1} \rar{\lambda} & \C(X,Y) \dar[equal] & \C(X,Y) \otimes \tensorunit \lar[swap]{\rho} \dar{1 \otimes i_X}\\
\C(Y,Y) \otimes \C(X,Y) \rar{m} & \C(X,Y) & \C(X,Y) \otimes \C(X,X) \lar[swap]{m}
\end{tikzcd}
\end{equation}
This finishes the definition of a $\V$-category.  Moreover, a $\V$-category $\C$ is \emph{small}\index{small!enriched category}\index{enriched category!small} if $\Ob(\C)$ is a set.
\end{definition}

\begin{definition}\label{def:enriched-functor}
For $\V$-categories $\C$ and $\D$, a \index{functor!enriched}\index{enriched!functor}\emph{$\V$-functor} $F \cn \C \to \D$ consists of the following data and axioms.
\begin{itemize}
\item It is equipped with an object assignment $F \cn \Ob(\C) \to \Ob(\D)$.
\item For each pair of objects $X,Y$ in $\C$, it is equipped with a component morphism
\[\begin{tikzcd}[column sep=large]
\C(X,Y) \rar{F_{X,Y}} & \D\bigl(FX,FY\bigr)\end{tikzcd} \inspace \V.\] 
\end{itemize}
The following two diagrams are required to commute for objects $X,Y,Z$ in $\C$.\index{composition!enriched functor}\index{identities!enriched functor} 
\begin{equation}\label{eq:enriched-composition}
\begin{tikzcd}[cells={nodes={scale=.9}}]
\C(Y,Z) \otimes \C(X,Y) \rar{m} \dar[swap]{F \otimes F} & \C(X,Z) \dar{F}\\
\D(FY,FZ) \otimes \D(FX,FY) \rar{m} & \D(FX,FZ)\end{tikzcd}
\qquad
\begin{tikzcd}[cells={nodes={scale=.9}}]
\tensorunit \dar[equal] \rar{i_X} & \C(X,X) \dar{F}\\
\tensorunit \rar{i_{FX}} & \D(FX,FX)
\end{tikzcd}
\end{equation}
Moreover, identity $\V$-functors and composition of $\V$-functors are defined separately on objects and component morphisms.
\end{definition}

\begin{definition}\label{def:enriched-natural-transformation}
For $\V$-functors $F,G \cn \C\to\D$ between $\V$-categories $\C$ and $\D$, a \index{enriched!natural transformation}\index{natural transformation!enriched}\emph{$\V$-natural transformation} $\theta \cn F\to G$ consists of, for each object $X$ in $\C$, a morphism
\[\begin{tikzcd}[column sep=large]
\tensorunit \ar{r}{\theta_X} & \D(FX,GX)
\end{tikzcd} \inspace \V,\]
called the \index{component}\emph{$X$-component} of $\theta$, such that the following naturality diagram commutes for objects $X,Y$ in $\C$.
\[\begin{tikzpicture}[x=40mm,y=12mm,vcenter]
\def\s{.85}
  \draw[0cell=\s] 
  (0,0) node (a) {
    \C(X,Y)
  }
  (.25,1) node (b1) {
    \tensorunit \otimes \C(X,Y)    
  }
  (b1)++(1,0) node (c1) {
    \D(FY,GY) \otimes \D(FX,FY)
  }
  (.25,-1) node (b2) {
    \C(X,Y) \otimes \tensorunit
  }
  (b2)++(1,0) node (c2) {
    \D(GX,GY) \otimes \D(FX,GX)
  }
  (c1)++(.25,-1) node (d) {
    \D(FX,GY)
  }
  ;
  \draw[1cell=\s] 
  (a) edge node[pos=.3] {\la^\inv} node['] {\iso}(b1)
  (a) edge node[swap,pos=.2] {\rho^\inv} node {\iso} (b2)
  (b1) edge node {\theta_Y \otimes F} (c1)
  (b2) edge node {G \otimes \theta_X} (c2)
  (c1) edge node[pos=.6] {m} (d)
  (c2) edge node[swap,pos=.7] {m} (d)
  ;
\end{tikzpicture}\]
This finishes the definition of a $\V$-natural transformation.  Moreover, the following notions are defined componentwise:
\begin{itemize}
\item identity $\V$-natural transformations,
\item vertical composition of $\V$-natural transformations, and
\item horizontal composition of $\V$-natural transformations.\defmark
\end{itemize}
\end{definition}

The following definition is from \cite[Section 1.3]{cerberusIII}.

\begin{definition}\label{definition:vtensor-0}
For $\V$-categories $\C$ and $\D$ with $(\V,\otimes,\tu,\xi)$ a braided monoidal category (\cref{def:braidedmoncat}), the \emph{tensor product}\index{tensor product!enriched category}\index{enriched category!tensor product} $\C \otimes \D$ is the $\V$-category defined by the following data.
\begin{description}
\item[Objects] $\Ob(\C \otimes \D) = \Ob(\C) \times \Ob(\D)$, with objects denoted $X \otimes Y$ for $X \in \C$ and $Y \in \D$.
\item[Hom Objects] For objects $X \otimes Y$ and $X' \otimes Y'$, the hom object is given by
\[(\C \otimes \D)(X \otimes Y, X' \otimes Y') = \C(X,X') \otimes \D(Y,Y').\]
\item[Composition] For objects $X \otimes Y$, $X'\otimes Y'$, and $X''\otimes Y''$, the composition is given by the following composite in $\V$, with $\ximid$ interchanging the middle two factors using the associativity isomorphism and braiding in $\V$.\label{not:ximid}
\[\begin{tikzpicture}[x=30mm,y=15mm]
\draw[0cell=.8]
(0,0) node (a) {\big(\C(X',X'')\otimes\D(Y',Y'')\big) \otimes \big(\C(X,X')\otimes\D(Y,Y')\big)}
(1,-1) node (b) {\big(\C(X',X'')\otimes\C(X,X')\big) \otimes \big(\D(Y',Y'')\otimes\D(Y,Y')\big)}
(2,0) node (c) {\C(X,X'') \otimes \D(Y,Y'')}
;
\draw[1cell=.9]
(a) edge[transform canvas={xshift=-1em}] node[swap,pos=.2] {\ximid} (b)
(b) edge[transform canvas={xshift=1em}] node[swap,pos=.7] {m \otimes m} (c);
\end{tikzpicture}\]
\item[Identities] The identity of an object $X \otimes Y$ is the following composite.
    \[\tensorunit \fto[\iso]{\la^\inv} \tensorunit \otimes \tensorunit \fto{i_X \otimes i_Y} \C(X,X)
    \otimes \D(Y,Y)\]
\end{description}
This finishes the definition of the $\V$-category $\C \otimes \D$.  Moreover, the tensor product $\otimes$ is extended to $\V$-functors and $\V$-natural transformations componentwise.
\end{definition}

\section{2-Categories}
\label{sec:twocategories}

In this section we review the definitions of 2-categories, 2-functors, 2-natural transformations, modifications, (adjoint) equivalences in 2-categories, and (adjoint) 2-equivalences between 2-categories.  
The reader is referred to \cite{johnson-yau} for an introduction to 2-categories.  

\begin{definition}\label{def:twocategory}\index{2-category}\index{category!2-}
A \emph{2-category} $\B$ consists of the following data and axioms.
\begin{description}
\item[Objects] It is equipped with a class $\B_0$ of \emph{objects}.
\item[1-Cells] For each pair of objects $X,Y\in\B_0$, it is equipped with a class $\B_1(X,Y)$ of \index{1-cell}\emph{1-cells} from $X$ to $Y$.  Such a 1-cell is denoted $X \to Y$.
\item[2-Cells] For 1-cells $f,f' \in \B_1(X,Y)$, it is equipped with a set $\B_2(f,f')$ of \index{2-cell}\emph{2-cells} from $f$ to $f'$.  Such a 2-cell is denoted $f \to f'$.
\item[Identities]  It is equipped with
\begin{itemize}
\item an \emph{identity 1-cell} $1_X \in \B_1(X,X)$ for each object $X$ and
\item an \emph{identity 2-cell} $1_f \in \B_2(f,f)$ for each 1-cell $f \in \B_1(X,Y)$.
\end{itemize}
\item[Compositions]
$X$, $Y$, and $Z$ denote objects in $\B$ in the following compositions.
\begin{itemize}
\item For 1-cells $f,f',f'' \in \B_1(X,Y)$, it is equipped with an assignment\label{not:vcompiicell}
\[\begin{tikzcd}
\B_2(f',f'') \times \B_2(f,f') \rar{v} & \B_2(f,f'')
\end{tikzcd},\qquad v(\alpha',\alpha) = \alpha'\alpha\] 
called the \index{vertical composition!2-category}\emph{vertical composition} of 2-cells.
\item It is equipped with an assignment\label{not:hcompicell}
\[\begin{tikzcd}
\B_1(Y,Z) \times \B_1(X,Y) \rar{c_1} & \B_1(X,Z)
\end{tikzcd},\qquad c_1(g,f) = gf\]
called the \index{horizontal composition!2-category}\emph{horizontal composition} of 1-cells.
\item For 1-cells $f,f' \in \B_1(X,Y)$ and $g,g' \in \B_1(Y,Z)$, it is equipped with an assignment\label{not:hcompiicell}
\[\begin{tikzcd}
\B_2(g,g') \times \B_2(f,f') \rar{c_2} & \B_2(gf,g'f')
\end{tikzcd},\quad c_2(\beta,\alpha) = \beta * \alpha\]
called the \emph{horizontal composition} of 2-cells.
\end{itemize}
\end{description}
The data above are required to satisfy the following four axioms:
\begin{enumerate}[label=(\roman*)]
\item\label{twocat-i} Vertical composition is associative and unital for identity 2-cells.
\item\label{twocat-ii} Horizontal composition preserves identity 2-cells and vertical composition.
\item\label{twocat-iii} Horizontal composition of 1-cells is associative and unital for identity 1-cells.
\item\label{twocat-iv} Horizontal composition of 2-cells is associative and unital for identity 2-cells of identity 1-cells.
\end{enumerate}
This finishes the definition of a 2-category.

Moreover, we define the following.
\begin{itemize}
\item For objects $X$ and $Y$ in a 2-category $\B$, the \index{hom category!2-category}\emph{hom category} \label{not:homcat}$\B(X,Y)$ is the category defined by the following data.
\begin{itemize}
\item Its objects are 1-cells from $X$ to $Y$.
\item Morphisms from a 1-cell $f \cn X \to Y$ to a 1-cell $f' \cn X \to Y$ are 2-cells $f \to f'$.
\item Its composition is vertical composition of 2-cells.
\item Its identities are identity 2-cells.
\end{itemize} 
\item A 2-category is \index{locally small}\emph{locally small} if each hom category is a small category.
\item A 2-category is \emph{small} if it is locally small and has a set of objects.
\item The \index{underlying 1-category}\index{1-category!underlying}\emph{underlying 1-category} of a 2-category is the category defined by the following data.
\begin{itemize}
\item It has the same class of objects.
\item Its morphisms are 1-cells.
\item Its composition is horizontal composition of 1-cells.
\item Its identities are identity 1-cells.
\end{itemize} 
\end{itemize} 
The 2-cell notation \cref{twocellnotation} also applies to 2-cells in a 2-category.
\end{definition}

\begin{example}[Small Categories]\label{ex:catastwocategory}\index{2-category!of small categories}
$\Cat$ in \cref{ex:cat} is the underlying 1-category of a 2-category, also denoted by $\Cat$, where the 2-cells are natural transformations.  Horizontal, respectively vertical, composition of 2-cells is defined as that of natural transformations. 
\end{example}

\begin{example}[Small Enriched Categories]\label{ex:vcatastwocategory}\index{2-category!of small enriched categories}
Each monoidal category $\V$ has an associated 2-category $\Catv$ defined by the following data.
\begin{itemize}
\item Its objects are small $\V$-categories.
\item Its 1-cells are $\V$-functors.
\item Its 2-cells are $\V$-natural transformations.
\end{itemize}
The 2-category $\Cat$ in \cref{ex:catastwocategory} is $\Catv$ for $\V = (\Set, \times, *)$.
\end{example}

A $\Cat$-category---that is, a category enriched in $(\Cat, \times, \boldone)$---is a 2-category.  The converse is also true up to a set-theoretic condition.

\begin{proposition}\label{locallysmalltwocat}
A locally small 2-category is precisely a $\Cat$-category.
\end{proposition}

\subsection*{2-Functors}

Up to the set-theoretic condition in \cref{locallysmalltwocat}, a 2-functor is a $\Cat$-functor between $\Cat$-categories.  We organize the precise definition as follows.

\begin{definition}\label{def:twofunctor}
For 2-categories $\B$ and $\C$, a \index{2-functor}\index{functor!2-}\emph{2-functor} $F \cn \B \to \C$ consists of the following data and axioms.
\begin{itemize}
\item It is equipped with an \emph{object assignment} 
\[\begin{tikzcd}[column sep=large]
\B_0 \ar{r}{F_0} & \C_0.
\end{tikzcd}\]
\item It is equipped with a \emph{1-cell assignment} 
\[\begin{tikzcd}[column sep=large]
\B_1(X,Y) \ar{r}{F_1} & \C_1(F_0 X,F_0 Y)
\end{tikzcd}\] 
for each pair of objects $X,Y$ in $\B$.
\item It is equipped with a \emph{2-cell assignment}
\[\begin{tikzcd}[column sep=large]
\B_2(f,f') \ar{r}{F_2} & \C_2(F_1 f,F_1 f')
\end{tikzcd}\] 
for each pair of objects $X,Y$ and 1-cells $f,f' \in \B_1(X,Y)$. 
\end{itemize}
The data above are required to satisfy the following three axioms, with $F_0$, $F_1$, and $F_2$ all abbreviated to $F$:
\begin{enumerate}[label=(\roman*)]
\item\label{twofunctor-i} The object and 1-cell assignments of $F$ form a functor between the underlying 1-categories of $\B$ and $\C$.
\item\label{twofunctor-ii} For each pair of objects $X,Y$ in $\B$, the 1-cell and 2-cell assignments of $F$ form a functor 
\[\begin{tikzcd}[column sep=large]
\B(X,Y) \ar{r}{F} & \C(FX,FY)
\end{tikzcd}\] 
between hom categories.
\item\label{twofunctor-iii} $F$ preserves horizontal composition of 2-cells.
\end{enumerate}
This finishes the definition of a 2-functor.

Moreover, given a 2-functor $G \cn \C \to \D$, the \emph{composite 2-functor} 
\[\begin{tikzcd}[column sep=large]
\B \ar{r}{GF} & \D
\end{tikzcd}\]
is defined by separately composing the object, 1-cell, and 2-cell assignments.
\end{definition}

\subsection*{2-Natural Transformations}

Up to the set-theoretic condition in \cref{locallysmalltwocat}, a 2-natural transformation is a $\Cat$-natural transformation between $\Cat$-functors.  We organize the precise definition as follows.

\begin{definition}\label{def:twonaturaltr}\index{2-natural!transformation}\index{natural transformation!2-}
For 2-functors $F,G \cn \B\to\C$ between 2-categories $\B$ and $\C$, a \emph{2-natural transformation} $\phi \cn F \to G$ consists of, for each object $X$ in $\B$, a \emph{component 1-cell} 
\[\begin{tikzcd}[column sep=large]
FX \ar{r}{\phi_X} & GX \inspace \C
\end{tikzcd}\]
such that the following two axioms hold:
\begin{description}
\item[1-Cell Naturality] For each 1-cell $f \cn X \to Y$ in $\B$, the following two composite 1-cells in $\C(FX,GY)$ are equal.
\begin{equation}\label{onecellnaturality}
\begin{tikzcd}[column sep=large]
FX \ar{d}[swap]{Ff} \ar{r}{\phi_X} & GX \ar{d}{Gf}\\
FY \ar{r}{\phi_Y} & GY
\end{tikzcd}
\end{equation}
\item[2-Cell Naturality] 
For each 2-cell $\theta \cn f \to g$ in $\B(X,Y)$, the following two whiskered 2-cells are equal.
\begin{equation}\label{twocellnaturality}
\begin{tikzpicture}[xscale=3,yscale=1.4,vcenter]
\def\a{30} \def\s{.85}
\draw[0cell=\s]
(0,0) node (x11) {FX}
(x11)++(1,0) node (x12) {GX}
(x11)++(0,-1) node (x21) {FY}
(x12)++(0,-1) node (x22) {GY}
;
\draw[1cell=\s]  
(x11) edge node {\phi_X} (x12)
(x21) edge node {\phi_Y} (x22)
(x11) edge[bend right] node[swap] {Ff} (x21)
(x11) edge[bend left] node {Fg} (x21)
(x12) edge[bend right] node[swap] {Gf} (x22)
(x12) edge[bend left] node {Gg} (x22)
;
\draw[2cell=.9]
node[between=x11 and x21 at .6, 2label={above,F\theta}] {\Rightarrow}
node[between=x12 and x22 at .6, 2label={above,G\theta}] {\Rightarrow}
;
\end{tikzpicture}
\end{equation}
This means that the following 2-cell equality holds in $\C(FX,GY)$.
\[G\theta * 1_{\phi_X} = 1_{\phi_Y} * F\theta\]
\end{description}
This finishes the definition of a 2-natural transformation.  The 2-cell notation \cref{twocellnotation} applies to 2-natural transformations.  Moreover, a \index{2-natural!isomorphism}\index{natural isomorphism!2-}\emph{2-natural isomorphism} is a 2-natural transformation with each component 1-cell an isomorphism in the underlying 1-category.
\end{definition}

\begin{definition}\label{def:twonatcomposition}
Suppose $\phi \cn F \to G$ and $\varphi \cn G \to H$ are 2-natural transformations for 2-functors $F,G,H \cn \B \to \C$ between 2-categories $\B$ and $\C$.  The \emph{horizontal composite}\index{horizontal composition!2-natural transformation} 2-natural transformation
\[\begin{tikzcd}[column sep=large]
F \ar{r}{\varphi \phi} & H
\end{tikzcd}\]
is defined by, for each object $X$ in $\B$, the horizontal composite component 1-cell 
\[\begin{tikzcd}[column sep=large]
FX \ar{r}{\phi_X} & GX \ar{r}{\varphi_X} & HX
\end{tikzcd}\]
in $\C(FX,HX)$.
\end{definition}

\subsection*{Modifications}

\begin{definition}\label{def:modification}
Suppose $\phi, \varphi \cn F \to G$ are 2-natural transformations for 2-functors $F,G \cn \B \to \C$ between 2-categories $\B$ and $\C$.  A \index{modification}\emph{modification} $\Phi \cn \phi \to \varphi$ consists of, for each object $X$ in $\B$, a \emph{component 2-cell}\label{not:componentiicell}
\[\begin{tikzcd}[column sep=large]
\phi_X \ar{r}{\Phi_X} & \varphi_X \inspace \C(FX,GX)
\end{tikzcd}\]
such that the following two whiskered 2-cells in $\C(FX,GY)$ are equal for each 1-cell $f \cn X \to Y$ in $\B$.
\begin{equation}\label{modificationaxiom}
\begin{tikzpicture}[xscale=2.5,yscale=1.5,vcenter]
\def\a{35} \def\s{.85}
\draw[0cell=\s]
(0,0) node (x11) {FX}
(x11)++(1,0) node (x12) {GX}
(x11)++(0,-1) node (x21) {FY}
(x12)++(0,-1) node (x22) {GY}
;
\draw[1cell=\s]  
(x11) edge[bend left=\a] node[pos=.4] {\phi_X} (x12)
(x11) edge[bend right=\a] node[swap,pos=.6] {\varphi_X} (x12)
(x21) edge[bend left=\a] node[pos=.4] {\phi_Y} (x22)
(x21) edge[bend right=\a] node[swap,pos=.6] {\varphi_Y} (x22)
(x11) edge node[swap] {Ff} (x21)
(x12) edge node {Gf} (x22)
;
\draw[2cell=.9]
node[between=x11 and x12 at .44, rotate=-90, 2label={above,\Phi_X}] {\Rightarrow}
node[between=x21 and x22 at .44, rotate=-90, 2label={above,\Phi_Y}] {\Rightarrow}
;
\end{tikzpicture}
\end{equation}
This means the equality
\[1_{Gf} * \Phi_X = \Phi_Y * 1_{Ff}.\]
This is called the \index{modification!axiom}\emph{modification axiom}.
\end{definition}

\begin{definition}\label{def:modcomposition}
Suppose $F,G,H \cn \B \to \C$ are 2-functors between 2-categories $\B$ and $\C$.
\begin{enumerate}[label=(\roman*)]
\item Suppose $\phi, \varphi, \psi \cn F \to G$ are 2-natural transformations, and $\Phi \cn \phi \to \varphi$ and $\Psi \cn \varphi \to \psi$ are modifications.  The \emph{vertical composite modification}\index{vertical composition!modification}
\begin{equation}\label{modificationvcomp}
\begin{tikzcd}[column sep=large]
\phi \ar{r}{\Psi \Phi} & \psi
\end{tikzcd}
\end{equation}
is defined by, for each object $X$ in $\B$, the vertical composite 2-cell
\[\begin{tikzcd}[column sep=large]
\phi_X \ar{r}{\Phi_X} & \varphi_X \ar{r}{\Psi_X} & \psi_X
\end{tikzcd}\]
in $\C(FX,GX)$.
\item Suppose $\Phi' \cn \phi' \to \varphi'$ is a modification for 2-natural transformations $\phi', \varphi' \cn G \to H$.  The \emph{horizontal composite modification}\index{horizontal composition!modification}
\begin{equation}\label{modificationhcomp}
\begin{tikzcd}[column sep=huge]
\phi' \phi \ar{r}{\Phi' * \Phi} & \varphi' \varphi
\end{tikzcd}
\end{equation}
is defined by, for each object $X$ in $\B$, the horizontal composite 2-cell
\[\begin{tikzcd}[column sep=huge]
\phi'_X \phi_X \ar{r}{\Phi'_X * \Phi_X} & \varphi'_X \varphi_X
\end{tikzcd}\]
in $\C(FX,HX)$.\defmark
\end{enumerate}
\end{definition}

\subsection*{Adjoint Equivalences}

The following definitions are \cite[5.1.18, 6.1.1, and 6.2.1]{johnson-yau} restricted to 2-categories.

\begin{definition}\label{def:equivalences}
Suppose $\B$ is a 2-category.
\begin{enumerate}[label=(\roman*)]
\item\label{def:equivalences-i} A 1-cell $f \cn X \to Y$ in $\B$ is called an \index{2-category!equivalence}\index{equivalence!2-category}\emph{equivalence} if there exist
\begin{itemize}
\item a 1-cell $h \cn Y \to X$, which is called an \emph{inverse} of $f$, and
\item 2-cell isomorphisms $\eta \cn 1_X \fto{\iso} hf$ and $\epz \cn fh \fto{\iso} 1_Y$.
\end{itemize}
\item\label{def:equivalences-ii} An \index{2-category!adjunction}\index{adjunction!2-category}\emph{adjunction} in $\B$ is a quadruple 
\[(f,h,\eta,\epz)\]
consisting of the following data and axioms.
\begin{itemize}
\item $f \cn X \to Y$ and $h \cn Y \to X$ are 1-cells, which are called, respectively, the \index{left adjoint}\emph{left adjoint} and the \index{right adjoint}\emph{right adjoint}.
\item $\eta \cn 1_X \to hf$ and $\epz \cn fh \to 1_Y$ are 2-cells, which are called, respectively, the \index{unit!adjunction}\emph{unit} and the \index{counit!adjunction}\emph{counit}.
\end{itemize}
The following two diagrams of 2-cells, which are called, respectively, the \index{triangle identities}\emph{left triangle identity} and the \emph{right triangle identity}, are required to commute.
\begin{equation}\label{triangleidentities}
\begin{tikzpicture}[xscale=1,yscale=1.3,vcenter]
\def\s{.9} \def\h{1.8}
\draw[0cell=\s]
(0,0) node (x11) {f}
(x11)++(\h,0) node (x12) {fhf}
(x12)++(0,-1) node (x22) {f}
(x12)++(\h,0) node (y11) {h}
(y11)++(\h,0) node (y12) {hfh}
(y12)++(0,-1) node (y22) {h}
;
\draw[1cell=\s]  
(x11) edge node {1_f * \eta} (x12)
(x12) edge node {\epz * 1_f} (x22)
(x11) edge node[swap] {1_f} (x22)
(y11) edge node {\eta * 1_h} (y12)
(y12) edge node {1_h * \epz} (y22)
(y11) edge node[swap] {1_h} (y22)
;
\end{tikzpicture}
\end{equation}
\item\label{def:equivalences-iii} An adjunction is called an \index{2-category!adjoint equivalence}\index{adjoint equivalence!2-category}\emph{adjoint equivalence} if the unit and the counit are isomorphisms.
\end{enumerate}
This finishes the definition.
\end{definition}

\begin{example}\label{ex:adjointequivalences}
Consider \cref{def:equivalences}.
\begin{itemize}
\item If $f$ is an equivalence with $(h,\eta,\epz)$ as in \cref{def:equivalences-i}, then $h$ is also an equivalence with inverse 1-cell $f$ and 2-cell isomorphisms $(\epz^\inv, \eta^\inv)$.
\item A quadruple $(f,h,\eta,\epz)$ is an adjoint equivalence as in \cref{def:equivalences-iii} if and only if the quadruple $(h,f,\epz^\inv,\eta^\inv)$ is an adjoint equivalence.
\item If $(f,h,\eta,\epz)$ is an adjoint equivalence, then $f$ and $h$ are equivalences.  Conversely, each equivalence is the left adjoint of an adjoint equivalence by \cite[6.2.4]{johnson-yau}.  It is also the right adjoint of an adjoint equivalence by the previous item.\defmark
\end{itemize}
\end{example}


\begin{definition}\label{def:iiequivalence}
A 2-functor $F \cn \B \to \C$ is called a \index{2-equivalence}\emph{2-equivalence} if there exist
\begin{itemize}
\item a 2-functor $G \cn \C \to \B$, which is called an \emph{inverse} of $F$, and
\item 2-natural isomorphisms $\eta \cn 1_{\B} \fto{\iso} GF$ and $\epz \cn FG \fto{\iso} 1_{\C}$.
\end{itemize}
Moreover, if the quadruple $(F,G,\eta,\epz)$ satisfies the triangle identities \cref{triangleidentities}, then it is called an \index{adjoint 2-equivalence}\emph{adjoint 2-equivalence}. 
\end{definition}

A functor between categories is an equivalence if and only if it is 
\begin{itemize}
\item essentially surjective on objects and
\item fully faithful on morphism sets.
\end{itemize}
A similar characterization of 2-equivalences is \cite[7.5.8]{johnson-yau}; see \cref{thm:iiequivalence}.   \cref{thm:multiwhitehead} is the extension to $\Cat$-multiequivalences.

\section{The Day Convolution}
\label{sec:dayconvolution}

In this section we review the Day convolution for functors between symmetric monoidal categories and associated constructions.  The special case where the codomain is the symmetric monoidal closed category $\Cat$ is discussed in detail in  \cref{expl:dcatsmclosed}.  The following discussion is adapted from \cite[Ch.\! 3]{cerberusIII}, where much more detail can be found.  The \emph{pointed} Day convolution for $\Ga$-categories, which we need for the inverse $K$-theory discussion in \cref{ch:invK}, is in \cref{def:Gacatsymmetricmonoidal}.

\begin{definition}\label{def:dayconvolution}
Suppose $(\V,\otimes,\tu,\xi,[,])$ is a symmetric monoidal closed category that is complete and cocomplete, and $(\cD,\Dtimes,e)$ is a small symmetric monoidal category.  Suppose $X, Y \cn \cD \to \V$ are $\cD$-diagrams in $\V$.  
\begin{itemize}
\item The \index{Day convolution}\emph{Day convolution} of $X$ and $Y$ is the $\cD$-diagram in $\V$ defined as the coend
\begin{equation}\label{eq:dayconvolution}
X \otimes Y = \int^{(a,b) \in \cD \times \cD} \coprod_{\cD(a \Dtimes b, -)} (Xa \otimes Yb).
\end{equation}
For natural transformations $\phi \cn X \to X'$ and $\psi \cn Y \to Y'$ between $\cD$-diagrams in $\V$, the \emph{Day convolution} $\phi \otimes \psi$ is the natural transformation
\begin{equation}\label{daynattr}
\begin{tikzcd}[column sep=large]
X \otimes Y \ar{r}{\phi \otimes \psi} & X' \otimes Y'
\end{tikzcd}
\end{equation}
induced by the morphisms in $\V$
\[\begin{tikzcd}[column sep=large]
Xa \otimes Yb \ar{r}{\phi_a \otimes \psi_b} & X'a \otimes Y'b
\end{tikzcd}
\forspace a,b \in \cD\]
via the universal properties of coproducts and coends. 
\item The \emph{braiding}\index{braiding!Day convolution}
\begin{equation}\label{daybraiding}
\begin{tikzcd}[column sep=large]
X \otimes Y \ar{r}{\xi_{X,Y}}[swap]{\iso} & Y \otimes X
\end{tikzcd}
\end{equation}
is the natural transformation induced by the braiding in $\V$
\[\begin{tikzcd}[column sep=large]
Xa \otimes Yb \ar{r}{\xi_{Xa,Yb}}[swap]{\iso} & Yb \otimes Xa
\end{tikzcd}
\forspace a,b \in \cD,\]
the braiding in $\cD$, and the universal properties of coproducts and coends. 
\item The \index{hom diagram}\index{diagram!hom}\index{Day convolution!hom diagram}\emph{hom diagram} from $X$ to $Y$ is defined as the end
\begin{equation}\label{eq:dayhom}
\Dhom(X,Y) = \int_{(b,c) \in \cD \times \cD} \bigg[\coprod_{\cD(- \Dtimes b, c)} Xb \scs Yc\bigg].
\end{equation}
This is extended to natural transformations between $\cD$-diagrams in $\V$ via the universal properties of coproducts and ends.
\item The \index{unit diagram}\index{diagram!unit}\index{Day convolution!unit diagram}\emph{unit diagram} $J$ is the $\cD$-diagram
\begin{equation}\label{eq:dayunit}
J = \coprod_{\cD(e,-)} \tu \cn \cD \to \V.
\end{equation}
\end{itemize}
This finishes the definition.
\end{definition}

\begin{explanation}[Universal Property]\label{expl:dayconvolution}
The Day convolution $X \otimes Y$ in \cref{eq:dayconvolution} is characterized by the following universal property.  For each pair of objects $(a,b) \in \cD \times \cD$, it is equipped with a structure morphism
\begin{equation}\label{daystructure}
\begin{tikzcd}[column sep=large]
Xa \otimes Yb \ar{r}{\zeta_{a,b}} & (X \otimes Y)(a \Dtimes b) \inspace \V
\end{tikzcd}
\end{equation}
such that the following two conditions hold.
\begin{description}
\item[Naturality] For each pair of morphisms $f \cn a \to a'$ and $g \cn b \to b'$ in $\cD$, the following diagram in $\V$ commutes.
\begin{equation}\label{dayconvnaturality}
\begin{tikzpicture}[xscale=5,yscale=1.3,vcenter]
\draw[0cell=.9]
(0,0) node (x11) {Xa \otimes Yb}
(x11)++(1,0) node (x12) {Xa' \otimes Yb'}
(x11)++(0,-1) node (x21) {(X \otimes Y)(a \Dtimes b)}
(x12)++(0,-1) node (x22) {(X \otimes Y)(a' \Dtimes b')}
;
\draw[1cell=.9]  
(x11) edge node {Xf \otimes Yg} (x12)
(x12) edge node {\zeta_{a',b'}} (x22)
(x11) edge node[swap] {\zeta_{a,b}} (x21)
(x21) edge node {(X \otimes Y)(f \Dtimes g)} (x22)
;
\end{tikzpicture}
\end{equation}
\item[Universality] The pair $(X \otimes Y,\zeta)$ is initial with respect to the previous condition.
\end{description}
We refer to \cref{dayconvnaturality} as the \emph{naturality} of $\zeta$.
\end{explanation}

\begin{explanation}[Iterated Day Convolution]\label{expl:nday}
For $n \geq 2$ suppose $X_1, \ldots, X_n$ are $\cD$-diagrams in $\V$.  By general coend calculus their iterated Day convolution is the $\cD$-diagram given by the coend
\begin{equation}\label{iteratedday}
\bigotimes_{j=1}^n X_j \iso \int^{\{a_j\}_{j=1}^n \,\in\, \cD^n} \coprod_{\cD(a_1 \Dtimes \cdots \Dtimes\, a_n \scs -)} \bigg( \bigotimes_{j=1}^n X_j a_j \bigg).
\end{equation}
It is characterized by the universal property involving the $n$-variable variants of the structure morphisms \cref{daystructure} and the naturality diagram \cref{dayconvnaturality}.  More explicitly, for each $n$-tuple of objects $\{a_j\}_{j=1}^n$ in $\cD$, it is equipped with a structure morphism
\begin{equation}\label{daystructuren}
\begin{tikzcd}[column sep=huge]
\bigotimes_{j=1}^n X_j a_j \ar{r}{\zeta_{a_1,\ldots,a_n}} & \big(\bigotimes_{j=1}^n X_j\big) (a_1 \Dtimes \cdots \Dtimes a_n) \inspace \V
\end{tikzcd}
\end{equation} 
such that the following two conditions hold.
\begin{description}
\item[Naturality] For each $n$-tuple of morphisms $\{f_j \cn a_j \to a_j'\}_{j=1}^n$ in $\cD$, denote by
\[a = a_1 \Dtimes \cdots \Dtimes a_n, \quad a' = a_1' \Dtimes \cdots \Dtimes a_n', \andspace f = f_1 \Dtimes \cdots \Dtimes f_n.\]
Then the following naturality diagram in $\V$ commutes.
\begin{equation}\label{dayconvnaturalityn}
\begin{tikzpicture}[xscale=5,yscale=1.5,vcenter]
\draw[0cell=.9]
(0,0) node (x11) {\textstyle\bigotimes_{j=1}^n X_j a_j}
(x11)++(1,0) node (x12) {\textstyle\bigotimes_{j=1}^n X_j a_j'}
(x11)++(0,-1) node (x21) {\textstyle\big(\bigotimes_{j=1}^n X_j\big) (a)}
(x12)++(0,-1) node (x22) {\textstyle\big(\bigotimes_{j=1}^n X_j\big) (a')}
;
\draw[1cell=.9]  
(x11) edge node {\otimes_{j=1}^n X_j f_j} (x12)
(x12) edge node {\zeta_{a_1',\ldots,a_n'}} (x22)
(x11) edge node[swap, shorten <=-1ex, shorten >=-1ex] {\zeta_{a_1,\ldots,a_n}} (x21)
(x21) edge node {(\otimes_{j=1}^n X_j)(f)} (x22)
;
\end{tikzpicture}
\end{equation}
\item[Universality] The pair $\big(\bigotimes_{j=1}^n X_j \scs \zeta\big)$ is initial with respect to the previous condition.
\end{description}
In \cref{iteratedday,daystructuren,dayconvnaturalityn} each iterated monoidal product is left normalized as in \cref{expl:leftbracketing}.
\end{explanation}

\subsection*{Symmetric Monoidal Closed Structure}

Recall the diagram category $\DC$ (\cref{def:diagramcat}).  The following theorem is due to Day \cite{day-convolution}.  A detailed proof is in \cite[3.7.22]{cerberusIII}, which covers the more general situation when $\cD$ is a small symmetric monoidal $\V$-category.  The simpler situation as stated below is sufficient for this work.

\begin{theorem}\label{thm:Day}\index{Day convolution!symmetric monoidal closed}\index{symmetric monoidal!Day convolution}
Suppose $\V$ is a symmetric monoidal closed category that is complete and cocomplete.  Suppose $\cD$ is a small symmetric monoidal category.  Then the diagram category $\DV$ is a symmetric monoidal closed category with the following data.
\begin{itemize}
\item The monoidal product is given by the Day convolution \cref{eq:dayconvolution}.
\item The braiding is given by \cref{daybraiding}.
\item The closed structure is given by the hom diagram \cref{eq:dayhom}.
\item The monoidal unit is given by the unit diagram $J$ \cref{eq:dayunit}.  
\end{itemize}
\end{theorem}

\begin{definition}[Tensored and Cotensored]\label{def:tensored}
Suppose $\C$ is a $\V$-category with $(\V,[,])$ a symmetric monoidal closed category.
\begin{description}
\item[Tensored] We say that $\C$ is \index{tensored}\emph{tensored over $\V$} if $\C$ is equipped with, for each $X \in \C$ and $A \in \V$, an object \label{not:tensored}$X \otimes A \in \C$ together with isomorphisms in $\V$
\[\C(X \otimes A, Y) \iso [A,\C(X,Y)]\]
for all $Y \in \C$.
\item[Cotensored] We say that $\C$ is \index{cotensored}\index{tensored!co-}\emph{cotensored over $\V$} if $\C$ is equipped with, for each $X \in \C$ and $A \in \V$, an object \label{not:cotensored}$X^A \in \C$ together with isomorphisms in $\V$
\[\C(Y,X^A) \iso [A,\C(Y,X)]\]
for all $Y \in \C$.\dqed
\end{description}
\end{definition}

The following result is a nontrivial consequence of Day's \cref{thm:Day}, as presented in \cite[3.9.15]{cerberusIII}.

\begin{theorem}\label{thm:DV}\index{diagram category!enriched, tensored, and cotensored}
Suppose $\V$ is a symmetric monoidal closed category that is complete and cocomplete.  Suppose $\cD$ is a small symmetric monoidal category.  Then the symmetric monoidal closed category $\DV$ is enriched, tensored, and cotensored over $\V$.
\end{theorem}

\begin{proof}
The reader is referred to \cite[Section 3.9]{cerberusIII} for the detailed proof.  Here we provide a brief outline.  Since $\DV$ is a symmetric monoidal closed category, it is enriched, tensored, and cotensored over itself.  Changing the enrichment along the symmetric monoidal functor\label{not:evaluation}
\[\begin{tikzcd}[column sep=large]
\DV \ar{r}{\ev_e} & \V,
\end{tikzcd}\]
given by evaluation at the monoidal unit $e$ in $\cD$, we conclude that $\DV$ is enriched, tensored, and cotensored over $\V$.
\end{proof}

\subsection*{Permutative-Indexed Categories}
Recall that a $\Cat$-category is a locally small 2-category (\cref{locallysmalltwocat}).

\begin{explanation}\label{expl:dcatsmclosed}
Let us explain the structure of $\Dcat$ in \cref{thm:Day,thm:DV} with
\begin{itemize}
\item $\V = (\Cat, \times, \boldone, [,])$ in \cref{ex:cat,ex:catastwocategory} and
\item $(\cD,\otimes,\tu,\beta)$ a small permutative category.
\end{itemize}
The 1-category $\cD$ is also regarded as a 2-category with only identity 2-cells.  We first describe the 2-category structure on $\Dcat$.
\begin{description}
\item[Objects] An object in $\Dcat$ is a 2-functor $\cD \to \Cat$ (\cref{def:twofunctor}).  This is the same thing as a functor $\cD \to \Cat$, since $\cD$ has only identity 2-cells.  In other words, an object in $\Dcat$ is a $\cD$-indexed category (\cref{def:diagramcat}).
\item[1-Cells]
For 2-functors $F,G \cn \cD \to \Cat$, the morphism object 
\[(\Dcat)(F,G)\]
is a small category.  An \emph{object} in $(\Dcat)(F,G)$ is a 2-natural transformation $\phi \cn F \to G$ (\cref{def:twonaturaltr}).  This means that $\phi$ is equipped with a component functor
\[\begin{tikzcd}[column sep=large]
Fa \ar{r}{\phi_a} & Ga
\end{tikzcd} \foreachspace a \in \cD.\]
For each morphism $f \cn a \to b$ in $\cD$, the following diagram of functors, which is the 1-cell naturality condition \cref{onecellnaturality}, is required to commute.
\begin{equation}\label{dcatonecellnaturality}
\begin{tikzcd}[column sep=large]
Fa \ar{r}{\phi_a} \ar{d}[swap]{Ff} & Ga \ar{d}{Gf}\\
Fb \ar{r}{\phi_b} & Gb
\end{tikzcd}
\end{equation}
The 2-cell naturality condition \cref{twocellnaturality} is automatically satisfied because $\cD$ has only identity 2-cells.  In other words, $\phi \cn F \to G$ is a natural transformation.

The \emph{identity 1-cell} $1_F \cn F \to F$ is the identity natural transformation.  It has component identity functors: 
\[(1_F)_a = 1_{Fa}.\]
\item[Horizontal Composition for 1-Cells]
For another 2-functor $H \cn \cD \to \Cat$ and a 2-natural transformation $\phi' \cn G \to H$, the \emph{horizontal composite} 
\[\begin{tikzcd}[column sep=large]
F \ar{r}{\phi' \phi} & H
\end{tikzcd}\]
is the 2-natural transformation as in \cref{def:twonatcomposition}.  It has component composite functor
\[\begin{tikzcd}[column sep=large]
Fa \ar{r}{\phi_a} & Ga \ar{r}{\phi'_a} & Ha
\end{tikzcd} \foreachspace a \in \cD.\]
\item[2-Cells] 
For 2-natural transformations $\phi, \varphi \cn F \to G$, a \emph{morphism} 
\[\begin{tikzcd}[column sep=large]
\phi \ar{r}{\Phi} & \varphi
\end{tikzcd} \inspace (\Dcat)(F,G)\]
is a modification (\cref{def:modification}).  This means that $\Phi$ is equipped with a component natural transformation $\Phi_a$ for each object $a \in \cD$, as displayed below.
\[\begin{tikzpicture}[xscale=2.3,yscale=1.7,baseline={(x1.base)}]
\draw[0cell=.9]
(0,0) node (x1) {Fa}
(x1)++(1,0) node (x2) {Ga}
;
\draw[1cell=.9]  
(x1) edge[bend left] node {\phi_a} (x2)
(x1) edge[bend right] node[swap] {\varphi_a} (x2)
;
\draw[2cell]
node[between=x1 and x2 at .44, rotate=-90, 2label={above,\Phi_a}] {\Rightarrow}
;
\end{tikzpicture}\]
The modification axiom \cref{modificationaxiom} in this case requires that, for each morphism $f \cn a \to b$ in $\cD$, the following two whiskered natural transformations are equal.
\begin{equation}\label{dcatmodification}
\begin{tikzpicture}[xscale=2.5,yscale=1.5,vcenter]
\def\a{35} \def\s{.9}
\draw[0cell=\s]
(0,0) node (x11) {Fa}
(x11)++(1,0) node (x12) {Ga}
(x11)++(0,-1) node (x21) {Fb}
(x12)++(0,-1) node (x22) {Gb}
;
\draw[1cell=\s]  
(x11) edge[bend left=\a] node[pos=.4] {\phi_a} (x12)
(x11) edge[bend right=\a] node[swap,pos=.6] {\varphi_a} (x12)
(x21) edge[bend left=\a] node[pos=.4] {\phi_b} (x22)
(x21) edge[bend right=\a] node[swap,pos=.6] {\varphi_b} (x22)
(x11) edge node[swap] {Ff} (x21)
(x12) edge node {Gf} (x22)
;
\draw[2cell=\s]
node[between=x11 and x12 at .44, rotate=-90, 2label={above,\Phi_a}] {\Rightarrow}
node[between=x21 and x22 at .44, rotate=-90, 2label={above,\Phi_b}] {\Rightarrow}
;
\end{tikzpicture}
\end{equation}
In other words, the modification axiom in this case is the equality
\[1_{Gf} * \Phi_a = \Phi_b * 1_{Ff}\]
of horizontal compositions of natural transformations for each morphism $f \cn a \to b$ in $\cD$.

The \emph{identity 2-cell} 
\[\begin{tikzcd}[column sep=large]
\phi \ar{r}{1_{\phi}} & \phi
\end{tikzcd} \inspace (\Dcat)(F,G)\]
is the identity modification.  It has component identity natural transformations: 
\[(1_{\phi})_a = 1_{\phi_a}.\] 
\item[Vertical Composition for 2-Cells]
For another modification 
\[\begin{tikzcd}[column sep=large]
\varphi \ar{r}{\Psi} & \psi 
\end{tikzcd}\inspace (\Dcat)(F,G)\]
the \emph{vertical composite modification} 
\[\begin{tikzcd}[column sep=large]
\phi \ar{r}{\Psi\Phi} & \psi 
\end{tikzcd}\inspace (\Dcat)(F,G)\]
is as in \cref{modificationvcomp}.  It has component vertical composite natural transformation
\[\begin{tikzpicture}[xscale=2.5,yscale=2.5,vcenter]
\def\a{50} \def\s{.9}
\draw[0cell=\s]
(0,0) node (x11) {Fa}
(x11)++(1,0) node (x12) {Ga}
;
\draw[1cell=\s]  
(x11) edge[bend left=\a] node[pos=.4] {\phi_a} (x12)
(x11) edge node [pos=.25] {\varphi_a} (x12)
(x11) edge[bend right=\a] node[swap,pos=.4] {\psi_a} (x12)
;
\draw[2cell=\s]
node[between=x11 and x12 at .46, shift={(0,.4)}, rotate=-90, 2label={above,\Phi_a}] {\Rightarrow}
node[between=x11 and x12 at .46, shift={(0,-.4)}, rotate=-90, 2label={above,\Psi_a}] {\Rightarrow}
;
\end{tikzpicture}\]
for each object $a \in \cD$.
\item[Horizontal Composition for 2-Cells]
For a modification 
\[\begin{tikzcd}[column sep=large]
\phi' \ar{r}{\Phi'} & \varphi' 
\end{tikzcd} \inspace (\Dcat)(G,H)\]
the \emph{horizontal composite modification}
\[\begin{tikzcd}[column sep=huge]
\phi' \phi \ar{r}{\Phi' * \Phi} & \varphi' \varphi
\end{tikzcd} \inspace (\Dcat)(F,H)\]
is as in \cref{modificationhcomp}.  It has component horizontal composite natural transformation
\[\begin{tikzpicture}[xscale=2.3,yscale=1.7,baseline={(x1.base)}]
\draw[0cell=.9]
(0,0) node (x1) {Fa}
(x1)++(1,0) node (x2) {Ga}
(x2)++(1,0) node (x3) {Ha}
;
\draw[1cell=.9]  
(x1) edge[bend left] node {\phi_a} (x2)
(x1) edge[bend right] node[swap] {\varphi_a} (x2)
(x2) edge[bend left] node {\phi'_a} (x3)
(x2) edge[bend right] node[swap] {\varphi'_a} (x3)
;
\draw[2cell]
node[between=x1 and x2 at .44, rotate=-90, 2label={above,\Phi_a}] {\Rightarrow}
node[between=x2 and x3 at .44, rotate=-90, 2label={above,\Phi'_a}] {\Rightarrow}
;
\end{tikzpicture}\]
for each object $a \in \cD$.  
\end{description}
This finishes the description of the 2-category $\Dcat$.

The (co)tensor structure (\cref{def:tensored}) on $\Dcat$ over $\Cat$ is induced by the symmetric monoidal closed structure on $\Cat$, that is, the Cartesian product $\times$ and diagram categories $[-,-]$.  More explicitly, suppose $F \cn \cD \to \Cat$ is a $\cD$-indexed category, and $\B$ is a small category.
\begin{description}
\item[Tensored]
The $\cD$-indexed category $F \otimes \B \cn \cD \to \Cat$ is the composite functor
\[\begin{tikzcd}[column sep=large]
\cD \ar{r}{F} & \Cat \ar{r}{- \times \B} & \Cat.
\end{tikzcd}\]
\item[Cotensored]
The $\cD$-indexed category $F^\B \cn \cD \to \Cat$ is the composite functor
\[\begin{tikzcd}[column sep=large]
\cD \ar{r}{F} & \Cat \ar{r}{[\B,-]} & \Cat.
\end{tikzcd}\]
\end{description}

Next we describe the monoidal structure on $\Dcat$.
\begin{description}
\item[Monoidal Unit]
It is the unit diagram
\begin{equation}\label{dcatunit}
J = \cD(\tu,-) \cn \cD \to \Cat.
\end{equation}
For each $a \in \cD$ the morphism set $\cD(\tu,a)$ is regarded as a discrete category with only identity morphisms.
\item[Monoidal Product]
On $\cD$-indexed categories and natural transformations, the monoidal product is the Day convolution in \cref{eq:dayconvolution,daynattr}, where $\otimes$ in $\Cat$ means the Cartesian product $\times$.  Suppose given modifications $\Phi \cn \phi \to \varphi$ and $\Phi' \cn \phi' \to \varphi'$, as displayed below, for $\cD$-indexed categories $F,G,F',G' \cn \cD \to \Cat$.
\[\begin{tikzpicture}[xscale=2,yscale=1.7,baseline={(x1.base)}]
\draw[0cell=.9]
(0,0) node (x1) {F}
(x1)++(1,0) node (x2) {G}
(x2)++(.7,0) node (y1) {F'}
(y1)++(1,0) node (y2) {G'}
;
\draw[1cell=.9]  
(x1) edge[bend left] node {\phi} (x2)
(x1) edge[bend right] node[swap] {\varphi} (x2)
(y1) edge[bend left] node {\phi'} (y2)
(y1) edge[bend right] node[swap] {\varphi'} (y2)
;
\draw[2cell]
node[between=x1 and x2 at .44, rotate=-90, 2label={above,\Phi}] {\Rightarrow}
node[between=y1 and y2 at .42, rotate=-90, 2label={above,\Phi'}] {\Rightarrow}
;
\end{tikzpicture}\]
The monoidal product $\Phi \otimes \Phi'$
\begin{equation}\label{daymodifications}
\begin{tikzpicture}[xscale=3,yscale=1.7,baseline={(x1.base)}]
\draw[0cell=.9]
(0,0) node (x1) {F \otimes F'}
(x1)++(1,0) node (x2) {G \otimes G'}
;
\draw[1cell=.9]  
(x1) edge[bend left] node {\phi \otimes \phi'} (x2)
(x1) edge[bend right] node[swap] {\varphi \otimes \varphi'} (x2)
;
\draw[2cell]
node[between=x1 and x2 at .34, rotate=-90, 2label={above,\Phi \otimes \Phi'}] {\Rightarrow}
;
\end{tikzpicture}
\end{equation}
is the modification whose component natural transformations are induced by the product natural transformations
\[\begin{tikzpicture}[xscale=3.5,yscale=1.7,baseline={(x1.base)}]
\draw[0cell=.9]
(0,0) node (x1) {Fa \times F'b}
(x1)++(1,0) node (x2) {Ga \times G'b}
;
\draw[1cell=.9]  
(x1) edge[bend left] node {\phi_a \times \phi'_b} (x2)
(x1) edge[bend right] node[swap] {\varphi_a \times \varphi'_b} (x2)
;
\draw[2cell]
node[between=x1 and x2 at .34, rotate=-90, 2label={above,\Phi_a \times \Phi'_b}] {\Rightarrow}
;
\end{tikzpicture}\]
for $a,b \in \cD$.
\end{description}
This completes the description of $\Dcat$.
\end{explanation}

\chapter{Symmetric Bimonoidal and Bipermutative Categories}
\label{ch:bipermutative}
In this chapter we discuss two closely related categorical structures, symmetric bimonoidal categories and bipermutative categories, several detailed examples, and Laplaza's Coherence Theorem for symmetric bimonoidal categories.  Here is a summary table of some of the definitions and results in this chapter.
\smallskip
\begin{center}
\resizebox{0.8\width}{!}{
{\renewcommand{\arraystretch}{1.4}%
{\setlength{\tabcolsep}{1ex}
\begin{tabular}{|c|c|}\hline
symmetric bimonoidal categories & \cref{def:sbc}\\ \hline
ring categories & \cref{def:ringcat} \\ \hline
bipermutative categories & \cref{def:embipermutativecat} \\ \hline
tight bipermutative categories & \cref{thm:smcbipermuative} \\ \hline
Laplaza's Coherence & \cref{thm:laplaza-coherence-1} \\ \hline
\end{tabular}}}}
\end{center}
\smallskip
In later chapters Laplaza's Coherence Theorem is a key part of several important definitions and constructions, including
\begin{itemize}
\item additive natural transformations \cref{laplazaapp},
\item additive modifications \cref{additivemodadditivity},
\item the Grothendieck construction of an additive natural transformation  (\cref{grodphiconstraintdom,grodphilinearity}), and
\item opcartesian $n$-linear functors \cref{laplazapx}.
\end{itemize}
As references, a detailed discussion of symmetric bimonoidal categories and Laplaza's Coherence Theorem is in \cite{cerberusI}.  Bipermutative categories and their relations to symmetric bimonoidal categories are discussed in detail in  \cite{cerberusII}.

\subsection*{Organization}

In \cref{sec:symbimonoidalcat} we review symmetric bimonoidal categories in the sense of Laplaza \cite{laplaza}.  A symmetric bimonoidal category is a categorical analog of a \emph{commutative rig}, which is a commutative ring without additive inverses.  

In \cref{sec:bipermutativecat} we review bipermutative categories in the sense of Elmendorf-Mandell \cite{elmendorf-mandell}.  The main difference between symmetric bimonoidal categories and bipermutative categories is that the former have \emph{distributivity morphisms} \cref{sbc-distributivity}, while the latter have \emph{factorization morphisms} \cref{ringcatfactorization} that go in the opposite direction as the distributivity morphisms.  Since none of these structure morphisms is required to be invertible in general, bipermutative categories do not form a subclass of symmetric bimonoidal categories, or vice versa.  On the other hand, a \emph{tight} bipermutative category, which has invertible factorization morphisms, is a \emph{tight} symmetric bimonoidal category, which has invertible distributivity morphisms (\cref{thm:smcbipermuative}).

In \cref{sec:examples} we discuss the following tight bipermutative categories.
\smallskip
\begin{center}
\resizebox{0.9\textwidth}{!}{
{\renewcommand{\arraystretch}{1.4}%
{\setlength{\tabcolsep}{1ex}
\begin{tabular}{|c|c|c|c|}\hline
category & objects & morphisms & ref \\ \hline
$\Finsk$ & finite sets & functions & \ref{ex:finsk}\\ \hline
$\Fset$ & finite sets & permutations & \ref{ex:Fset} \\ \hline
$\Fskel$ & pointed finite sets & pointed functions & \ref{ex:Fskel} \\ \hline
$\cA$ & finite sequences of positive integers & certain maps of unpointed finite sets & \ref{ex:mandellcategory} \\ \hline
\end{tabular}}}}
\end{center}
\smallskip
In particular, \cref{ex:mandellcategory} exhibits a bipermutative category structure on Mandell's indexing category $\cA$ for inverse $K$-theory \cite{mandell_inverseK}.  This bipermutative structure on $\cA$ is a crucial reason why inverse $K$-theory is a \emph{pseudo symmetric} $\Cat$-multifunctor, as we will discuss in \cref{ch:multigro,ch:invK}.

In \cref{sec:laplazacoherence} we review Laplaza's Coherence \cref{thm:laplaza-coherence-1} for symmetric bimonoidal categories.  Analogous to Mac Lane's Coherence Theorem for monoidal categories \cite[VII.2]{maclane}, Laplaza's Coherence Theorem states that, in a symmetric bimonoidal category that satisfies a certain monomorphism condition, each formal diagram built from the symmetric bimonoidal structure commutes, provided that the domain is \emph{regular}.  The regularity condition excludes generally non-commutative diagrams such as
\[\begin{tikzcd}
x \otimes x \ar[shift left]{r}{1} \ar[shift right]{r}[swap]{\xitimes} & x \otimes x\end{tikzcd}\]
with $\xitimes$ the multiplicative braiding.  The monomorphism condition in Laplaza's Coherence Theorem is satisfied if the distributivity morphisms are invertible.  In particular, by \cref{thm:smcbipermuative} Laplaza's Coherence Theorem applies to tight bipermutative categories, including those in \cref{sec:examples}.  In  \cref{ch:diagram} we use Laplaza's Coherence \cref{thm:laplaza-coherence-1} to construct the $\Cat$-enriched multicategory of certain indexed categories when the indexing category is a tight bipermutative category.

\section{Symmetric Bimonoidal Categories}
\label{sec:symbimonoidalcat}

In this section we review symmetric bimonoidal categories in the sense of Laplaza \cite{laplaza}.  The presentation below is adapted from \cite[Ch.\! 2]{cerberusI}, where much more detailed discussion can be found.

\begin{definition}\label{def:sbc}
A \emph{symmetric bimonoidal category}\index{symmetric bimonoidal category}\index{category!symmetric bimonoidal}\index{bimonoidal category!symmetric} is a tuple
\[\big(\C, (\oplus,\zero,\alphaplus,\lambdaplus,\rhoplus,\xiplus),(\otimes,\tensorunit,\alphatimes,\lambdatimes,\rhotimes,\xitimes), (\lambdadot,\rhodot), (\deltal,\deltar)\big)\]
consisting of the following data.
\begin{description}
\item[Additive Structure] $\Cplus = (\C,\oplus,\zero,\alphaplus,\lambdaplus,\rhoplus,\xiplus)$\label{not:sbcadditive} is a symmetric monoidal category (\cref{def:symmoncat}).\index{additive structure!symmetric bimonoidal category}
\item[Multiplicative Structure] $\Cte = (\C,\otimes,\tensorunit,\alphatimes,\lambdatimes,\rhotimes,\xitimes)$\label{not:sbcmultiplicative} is a symmetric monoidal category.\index{multiplicative structure!symmetric bimonoidal category}
\item[Multiplicative Zeros] $\lambdadot$ and $\rhodot$ are natural isomorphisms
\begin{equation}\label{sbc-multiplicative-zero}
\begin{tikzcd}[column sep=large]
\zero \otimes A \ar{r}{\lambdadot_A}[swap]{\cong} & \zero & A \otimes \zero \ar{l}{\cong}[swap]{\rhodot_A}
\end{tikzcd} \forspace A \in \C,
\end{equation}
which are called the \index{left multiplicative zero!symmetric bimonoidal category}\emph{left multiplicative zero} and the \index{right multiplicative zero!symmetric bimonoidal category}\emph{right multiplicative zero}, respectively.
\item[Distributivity] $\deltal$ and $\deltar$ are natural monomorphisms
\begin{equation}\label{sbc-distributivity}
\begin{tikzcd}[row sep=tiny,column sep=huge]
A \otimes (B \oplus C) \ar{r}{\deltal_{A,B,C}} & (A \otimes B) \oplus (A \otimes C)\\
(A \oplus B) \otimes C \ar{r}{\deltar_{A,B,C}} & (A \otimes C) \oplus (B \otimes C)
\end{tikzcd}
\end{equation}
for objects $A,B,C\in\C$, which are called the \index{left distributivity morphism!symmetric bimonoidal category}\emph{left distributivity morphism} and the \index{right distributivity morphism!symmetric bimonoidal category}\emph{right distributivity morphism}, respectively.
\end{description}
We abbreviate $\otimes$ using concatenation and let $\tensor$ take precedence over $\oplus$.  For example, the left distributivity morphism is abbreviated to $A(B \oplus C) \to AB \oplus AC$.  

The following 24 diagrams in $\C$, which are called \index{Laplaza's Axioms}\emph{Laplaza's Axioms}, are required to commute for all objects $A,B,C$, and $D$ in $\C$.
\begin{description}
\item[Distributivity and Multiplicative Braiding]
\begin{equation}\label{laplaza-II}
\begin{tikzcd}[column sep=large]
(A\oplus B)C \ar{d}[swap]{\xitimes_{A\oplus B,C}} \ar{r}{\deltar_{A,B,C}} & AC \oplus BC \ar{d}{\xitimes_{A,C}\oplus \xitimes_{B,C}}\\
C(A\oplus B) \ar{r}{\deltal_{C,A,B}} & CA \oplus CB
\end{tikzcd}
\end{equation}
\item[Distributivity and Additive Braiding]
\begin{equation}\label{laplaza-I}
\begin{tikzcd}[column sep=large]
A(B\oplus C) \ar{d}[swap]{1_A\xiplus_{B,C}} \ar{r}{\deltal_{A,B,C}} & AB \oplus AC \ar{d}{\xiplus_{AB,AC}}\\
A(C\oplus B) \ar{r}{\deltal_{A,C,B}} & AC \oplus AB
\end{tikzcd}
\end{equation}
\begin{equation}\label{laplaza-III}
\begin{tikzcd}[column sep=large]
(A\oplus B)C \ar{d}[swap]{\xiplus_{A,B} 1_C} \ar{r}{\deltar_{A,B,C}} & AC \oplus BC \ar{d}{\xiplus_{AC,BC}}\\
(B\oplus A)C \ar{r}{\deltar_{B,A,C}} & BC \oplus AC
\end{tikzcd}
\end{equation}
\item[Distributivity and Additive Associativity]
\begin{equation}\label{laplaza-IV}
\begin{tikzcd}[column sep=large, cells={nodes={scale=.75}}]
[(A\oplus B)\oplus C]D \ar{d}[swap]{\alphaplus_{A,B,C} 1_D} \ar{r}{\deltar_{A\oplus B,C,D}} & (A\oplus B)D \oplus CD \ar{r}[scale=.75]{\deltar_{A,B,D} \oplus 1_{CD}} & (AD \oplus BD) \oplus CD \ar{d}{\alphaplus_{AD,BD,CD}}\\
\left[A \oplus (B\oplus C)\right]D \ar{r}{\deltar_{A,B\oplus C,D}} & AD \oplus (B\oplus C)D \ar{r}[scale=.75]{1_{AD}\oplus \deltar_{B,C,D}} & AD \oplus (BD \oplus CD)
\end{tikzcd}
\end{equation}
\begin{equation}\label{laplaza-V}
\begin{tikzcd}[column sep=large, cells={nodes={scale=.75}}]
A[(B\oplus C)\oplus D] \ar{d}[swap]{1_A\alphaplus_{B,C,D}} \ar{r}{\deltal_{A,B\oplus C,D}} & A(B\oplus C) \oplus AD \ar{r}[scale=.75]{\deltal_{A,B,C} \oplus 1_{AD}} & (AB \oplus AC) \oplus AD \ar{d}{\alphaplus_{AB,AC,AD}}\\
A[B\oplus(C\oplus D)] \ar{r}{\deltal_{A,B,C\oplus D}} & AB \oplus A(C\oplus D) \ar{r}[scale=.75]{1_{AB}\oplus \deltal_{A,C,D}} & AB \oplus (AC \oplus AD)
\end{tikzcd}
\end{equation}
\item[Distributivity and Multiplicative Associativity]
\begin{equation}\label{laplaza-VI}
\begin{tikzcd}[column sep=large, cells={nodes={scale=.75}}]
(AB)(C\oplus D) \ar{d}[swap]{\alphatimes_{A,B,C\oplus D}} \ar{rr}{\deltal_{AB,C,D}} && (AB)C \oplus (AB)D \ar{d}{\alphatimes_{A,B,C}\oplus\alphatimes_{A,B,D}}\\
A[B(C\oplus D)] \ar{r}{1_A\deltal_{B,C,D}} & A(BC \oplus BD) \ar{r}{\deltal_{A,BC,BD}} & A(BC) \oplus A(BD)
\end{tikzcd}
\end{equation}
\begin{equation}\label{laplaza-VII}
\begin{tikzcd}[column sep=large, cells={nodes={scale=.75}}]
[(A\oplus B)C]D \ar{d}[swap]{\alphatimes_{A\oplus B,C,D}} \ar{r}{\deltar_{A,B,C}1_D} & (AC\oplus BC)D \ar{r}{\deltar_{AC,BC,D}} & (AC)D \oplus (BC)D \ar{d}{\alphatimes_{A,C,D}\oplus\alphatimes_{B,C,D}}\\
(A\oplus B)(CD) \ar{rr}{\deltar_{A,B,CD}} && A(CD) \oplus B(CD)
\end{tikzcd}
\end{equation}
\begin{equation}\label{laplaza-VIII}
\begin{tikzcd}[column sep=large, cells={nodes={scale=.75}}]
[A(B\oplus C)]D \ar{d}[swap]{\alphatimes_{A,B\oplus C,D}} \ar{r}{\deltal_{A,B,C}1_D} & (AB \oplus AC)D \ar{r}{\deltar_{AB,AC,D}} & (AB)D \oplus (AC)D \ar{d}{\alphatimes_{A,B,D}\oplus\alphatimes_{A,C,D}}\\
A[(B\oplus C)D] \ar{r}{1_A\deltar_{B,C,D}} & A(BD \oplus CD) \ar{r}{\deltal_{A,BD,CD}} & A(BD) \oplus A(CD)
\end{tikzcd}
\end{equation}
\item[2-By-2 Distributivity]
\begin{equation}\label{laplaza-IX}
\begin{tikzcd}
(A\oplus B)(C\oplus D) \ar{d}[swap]{\deltal_{A\oplus B,C,D}} \ar{r}{\deltar_{A,B,C\oplus D}} & A(C\oplus D) \oplus B(C \oplus D) \ar{d}{\deltal_{A,C,D}\oplus\deltal_{B,C,D}}\\
(A\oplus B)C \oplus (A\oplus B)D \ar{d}[swap]{\deltar_{A,B,C}\oplus\deltar_{A,B,D}} & (AC \oplus AD) \oplus (BC \oplus BD) \ar{d}{\alphaplus_{AC,AD,BC\oplus BD}}\\
(AC \oplus BC) \oplus (AD \oplus BD) \ar{d}[swap]{\alphaplus_{AC,BC,AD\oplus BD}} & AC \oplus \left[AD \oplus (BC \oplus BD)\right] \ar{d}{1_{AC}\oplus (\alphaplus)^{-1}}\\
AC \oplus \left[BC \oplus (AD \oplus BD)\right] \ar{d}[swap]{1_{AC}\oplus (\alphaplus)^{-1}} & AC \oplus \left[(AD \oplus BC) \oplus BD\right] \ar{d}{1_{AC}\oplus (\xiplus_{AD,BC}\oplus 1_{BD})}\\ 
AC \oplus \left[(BC \oplus AD) \oplus BD\right] \ar[equal]{r} & AC \oplus \left[(BC \oplus AD) \oplus BD\right] 
\end{tikzcd}
\end{equation}
\item[Multiplicative Zero of 0]
\begin{equation}\label{laplaza-X}
\begin{tikzcd}[column sep=large]
\zero \otimes \zero \ar[shift left]{r}{\lambdadot_\zero} \ar[shift right]{r}[swap]{\rhodot_\zero} & \zero
\end{tikzcd}
\end{equation}
\item[Multiplicative Zero of a Sum]
\begin{equation}\label{laplaza-XI}
\begin{tikzcd}[column sep=large]
\zero(A\oplus B) \ar{d}[swap]{\lambdadot_{A\oplus B}} \ar{r}{\deltal_{\zero,A,B}} & \zero A \oplus \zero B \ar{d}{\lambdadot_A \oplus \lambdadot_B}\\
\zero & \zero \oplus \zero \ar{l}[swap]{\lambdaplus_{\zero}}
\end{tikzcd}
\end{equation}
\begin{equation}\label{laplaza-XII}
\begin{tikzcd}[column sep=large]
(A\oplus B)\zero \ar{d}[swap]{\rhodot_{A\oplus B}} \ar{r}{\deltar_{A,B,\zero}} & A\zero \oplus B\zero \ar{d}{\rhodot_A \oplus \rhodot_B}\\
\zero & \zero \oplus \zero \ar{l}[swap]{\lambdaplus_{\zero}}
\end{tikzcd}
\end{equation}
\item[Multiplicative Zero and Multiplicative Unit]
\begin{equation}\label{laplaza-XIII}
\begin{tikzcd}[column sep=large]
\zero \otimes \tu \ar[shift left]{r}{\lambdadot_{\tu}} \ar[shift right]{r}[swap]{\rhotimes_{\zero}} & \zero
\end{tikzcd}
\end{equation}
\begin{equation}\label{laplaza-XIV}
\begin{tikzcd}[column sep=large]
\tu \otimes \zero \ar[shift left]{r}{\rhodot_{\tu}} \ar[shift right]{r}[swap]{\lambdatimes_{\zero}} & \zero
\end{tikzcd}
\end{equation}
\item[Braiding of Multiplicative Zero]
\begin{equation}\label{laplaza-XV}
\begin{tikzcd}[column sep=small]
A \otimes \zero \ar{dr}[swap]{\rhodot_A} \ar{rr}{\xitimes_{A,\zero}} && \zero \otimes A \ar{dl}{\lambdadot_A}\\ & \zero &
\end{tikzcd}
\end{equation}
\item[Multiplicative Zero and Multiplicative Associativity]
\begin{equation}\label{laplaza-XVIII}
\begin{tikzcd}[column sep=large]
(AB)\zero \ar{d}[swap]{\rhodot_{AB}} \ar{r}{\alphatimes_{A,B,\zero}} & A(B\zero) \ar{d}{1_A\rhodot_B}\\
\zero & A\zero \ar{l}[swap]{\rhodot_A}
\end{tikzcd}
\end{equation}
\begin{equation}\label{laplaza-XVII}
\begin{tikzcd}
(A\zero)B \ar{d}[swap]{\rhodot_A 1_B} \ar{rr}{\alphatimes_{A,\zero,B}} && A(\zero B) \ar{d}{1_A\lambdadot_B}\\
\zero B \ar{dr}[swap]{\lambdadot_B} && A\zero \ar{dl}{\rhodot_A}\\
& \zero 
\end{tikzcd}
\end{equation}
\begin{equation}\label{laplaza-XVI}
\begin{tikzcd}[column sep=large]
(\zero A)B \ar{d}[swap]{\alphatimes_{\zero,A,B}} \ar{r}{\lambdadot_{A}1_B} & \zero B \ar{d}{\lambdadot_B}\\
\zero(AB) \ar{r}{\lambdadot_{AB}} & \zero
\end{tikzcd}
\end{equation}
\item[Additive and Multiplicative Zero]
\begin{equation}\label{laplaza-XIX}
\begin{tikzcd}[column sep=large]
A(\zero \oplus B) \ar{d}[swap]{1_A\lambdaplus_B} \ar{r}{\deltal_{A,\zero,B}} & A\zero \oplus AB \ar{d}{\rhodot_A \oplus 1_{AB}}\\
AB & \zero \oplus AB \ar{l}[swap]{\lambdaplus_{AB}}
\end{tikzcd}
\end{equation}
\begin{equation}\label{laplaza-XX}
\begin{tikzcd}[column sep=large]
(\zero \oplus B)A \ar{d}[swap]{\lambdaplus_B 1_A} \ar{r}{\deltar_{\zero,B,A}} & \zero A \oplus BA \ar{d}{\lambdadot_A \oplus 1_{BA}}\\
BA & \zero \oplus BA \ar{l}[swap]{\lambdaplus_{BA}} 
\end{tikzcd}
\end{equation}
\begin{equation}\label{laplaza-XXI}
\begin{tikzcd}[column sep=large]
A(B \oplus \zero) \ar{d}[swap]{1_A\rhoplus_B} \ar{r}{\deltal_{A,B,\zero}} & AB \oplus A\zero \ar{d}{1_{AB}\oplus \rhodot_A}\\
AB & AB \oplus\zero \ar{l}[swap]{\rhoplus_{AB}}
\end{tikzcd}
\end{equation}
\begin{equation}\label{laplaza-XXII}
\begin{tikzcd}[column sep=large]
(B \oplus \zero)A \ar{d}[swap]{\rhoplus_B 1_A} \ar{r}{\deltar_{B,\zero,A}} & BA \oplus \zero A \ar{d}{1_{BA}\oplus \lambdadot_A}\\
BA & BA \oplus \zero \ar{l}[swap]{\rhoplus_{BA}}
\end{tikzcd}
\end{equation}
\item[Distributivity and Multiplicative Unit]
\begin{equation}\label{laplaza-XXIII}
\begin{tikzcd}[column sep=small]
\tu(A\oplus B) \ar{dr}[swap]{\lambdatimes_{A\oplus B}} \ar{rr}{\deltal_{\tu,A,B}} && \tu A \oplus \tu B \ar{dl}{\lambdatimes_A \oplus \lambdatimes_B}\\
& A\oplus B &
\end{tikzcd}
\end{equation}
\begin{equation}\label{laplaza-XXIV}
\begin{tikzcd}[column sep=small]
(A\oplus B)\tu \ar{dr}[swap]{\rhotimes_{A\oplus B}} \ar{rr}{\deltar_{A,B,\tu}} && A\tu \oplus B\tu  \ar{dl}{\rhotimes_A \oplus \rhotimes_B}\\
& A\oplus B &
\end{tikzcd}
\end{equation}
\end{description}
This finishes the definition of a symmetric bimonoidal category.

Moreover, we define the following.
\begin{itemize}
\item A symmetric bimonoidal category is \index{symmetric bimonoidal category!small}\index{small!symmetric bimonoidal category}\emph{small} if its class of objects is a set.
\item A symmetric bimonoidal category is \index{symmetric bimonoidal category!tight}\index{tight!symmetric bimonoidal category}\emph{tight} if both $\deltal$ and $\deltar$ are natural isomorphisms.
\item The objects $\zero$ and $\tu$ are called the \emph{additive zero}\index{additive!zero} and the \index{multiplicative!unit}\emph{multiplicative unit}, respectively.
\item $\oplus$ and $\otimes$ are called the \index{sum}\emph{sum} and the \index{product}\emph{product}, respectively.
\item $\alphaplus$, $\lambdaplus$, $\rhoplus$, and $\xiplus$ are called the \index{additive!associativity isomorphism}\index{associativity isomorphism!additive}\emph{additive associativity isomorphism}, the \index{left additive zero}\index{additive!left - zero}\emph{left additive zero}, the \index{right additive zero}\index{additive!right - zero}\emph{right additive zero}, and the \index{additive!braiding}\index{braiding!additive}\emph{additive braiding}, respectively.
\item $\alphatimes$, $\lambdatimes$, $\rhotimes$, and $\xitimes$ are called the \index{multiplicative!associativity isomorphism}\index{associativity isomorphism!multiplicative}\emph{multiplicative associativity isomorphism}, the \index{left multiplicative unit}\index{multiplicative!left - unit}\emph{left multiplicative unit}, the \index{right multiplicative unit}\index{multiplicative!right - unit}\emph{right multiplicative unit}, and the \index{multiplicative!braiding}\index{braiding!multiplicative}\emph{multiplicative braiding}, respectively.\defmark
\end{itemize}
\end{definition}

Due to the presence of the axioms 
\begin{itemize}
\item \eqref{laplaza-II} relating $\deltal$ and $\deltar$ and
\item \eqref{laplaza-XV} relating $\lambdadot$ and $\rhodot$,
\end{itemize}
12 of the 24 axioms in \cref{def:sbc} follow from the other 12 axioms.  The following observation is \cite[2.2.13]{cerberusI}.

\begin{theorem}\label{thm:sbc-redundant-axioms}
In \Cref{def:sbc} of a symmetric bimonoidal category, it is sufficient to assume the first axiom in each of the twelve groups of axioms, namely, \eqref{laplaza-II}, \eqref{laplaza-I}, \eqref{laplaza-IV}, \eqref{laplaza-VI}, \eqref{laplaza-IX}, \eqref{laplaza-X}, \eqref{laplaza-XI}, \eqref{laplaza-XIII}, \eqref{laplaza-XV}, \eqref{laplaza-XVIII}, \eqref{laplaza-XIX}, and \eqref{laplaza-XXIII}.
\end{theorem}

\section{Bipermutative Categories}
\label{sec:bipermutativecat}

In this section we review bipermutative categories in the sense of Elmendorf-Mandell \cite{elmendorf-mandell}.  The presentation below is adapted from \cite[Ch.\! 9]{cerberusII}.  A bipermutative category is a \emph{ring category}, to be defined below, with further structure and axioms.  The following concept of a ring category is due to \cite[Def.\! 3.3]{elmendorf-mandell}, as presented in \cite[Section 9.1]{cerberusII}.  The main difference between these concepts and symmetric bimonoidal categories (\cref{def:sbc}) is that ring and bipermutative categories have \emph{factorization morphisms} that go in the opposite direction as the distributivity morphisms \cref{sbc-distributivity}.

\begin{definition}\label{def:ringcat}\index{ring category}\index{category!ring}
A \emph{ring category} is a tuple
\[\big(\C,(\oplus,\zero,\xiplus),(\otimes,\tu),(\fal,\far)\big)\]
consisting of the following data.
\begin{description}
\item[Additive Structure]\index{additive structure!ring category}\index{ring category!additive structure} 
$\Cplus = (\C,\oplus,\zero,\xiplus)$\label{not:Cplusring} is a permutative category (\cref{def:symmoncat}), with $\oplus$, $\zero$, and $\xiplus$ called, respectively, the \index{ring category!sum}\index{sum}\emph{sum}, the \index{ring category!additive zero}\index{additive zero}\emph{additive zero}, and the \index{additive braiding}\emph{additive braiding}.
\item[Multiplicative Structure]\index{multiplicative structure!ring category}\index{ring category!multiplicative structure}
$\Cte = (\C,\otimes,\tu)$\label{not:Ctering} is a strict monoidal category (\cref{def:monoidalcategory}), with $\otimes$ and $\tu$ called, respectively, the \index{ring category!product}\index{product}\emph{product} and the \index{ring category!multiplicative unit}\index{multiplicative unit}\emph{multiplicative unit}.
\item[Factorization Morphisms] 
$\fal$ and $\far$ are natural transformations
\begin{equation}\label{ringcatfactorization}
\begin{tikzcd}[row sep=tiny,column sep=huge]
(A \otimes C) \oplus (B \otimes C) \ar{r}{\fal_{A,B,C}} & (A \oplus B) \otimes C\\
(A \otimes B) \oplus (A \otimes C) \ar{r}{\far_{A,B,C}} & A \otimes (B \oplus C)
\end{tikzcd}
\end{equation}
for objects $A,B,C\in\C$, which are called the \emph{left factorization morphism}\index{left factorization morphism!ring category}\index{ring category!left factorization morphism} and the \emph{right factorization morphism}\index{right factorization morphism!ring category}\index{ring category!right factorization morphism}, respectively.
\end{description}
We abbreviate $\otimes$ to concatenation, with $\tensor$ taking precedence over $\oplus$ in the absence of clarifying parentheses.  For example, the left factorization morphism is abbreviated to 
\[AC \oplus BC \to (A \oplus B)C.\]  
The subscripts in $\xiplus$, $\fal$, and $\far$ are sometimes omitted.

The data above are required to satisfy the following seven axioms for all objects $A$, $A'$, $A''$, $B$, $B'$, $B''$, $C$, and $C'$ in $\C$, with $\boldone$ denoting the terminal category with one object $*$ and its identity morphism.
\begin{description}
\item[The Multiplicative Zero Axiom] 
The following diagram of functors commutes.\index{ring category!multiplicative zero axiom}\index{multiplicative zero axiom}
\begin{equation}\label{ringcataxiommultzero}
\begin{tikzcd}[column sep=large]
\boldone \times \C \ar{d}[swap]{\zero \times 1_{\C}} \ar{r}{\iso} & \C \ar{d}{\zero} & \C \times \boldone \ar{l}[swap]{\iso} \ar{d}{1_{\C} \times \zero}\\
\C \times \C \ar{r}{\otimes} & \C & \C \times \C \ar{l}[swap]{\otimes}
\end{tikzcd}
\end{equation}
In \cref{ringcataxiommultzero} the top horizontal isomorphisms drop the $\boldone$ argument.  Each $\zero$ denotes the constant functor at $\zero \in \C$ and $1_\zero$.
\item[The Zero Factorization Axiom]\index{ring category!zero factorization axiom}\index{zero factorization axiom} 
\begin{equation}\label{ringcataxiomzerofact}
\begin{aligned}
\fal_{\zero,B,C} &= 1_{B \otimes C} &\qquad \far_{\zero,B,C} &= 1_{\zero}\\
\fal_{A,\zero,C} &= 1_{A \otimes C} & \far_{A,\zero,C} &= 1_{A \otimes C}\\
\fal_{A,B,\zero} &= 1_{\zero} & \far_{A,B,\zero} &= 1_{A \otimes B}
\end{aligned}
\end{equation}
The three equalities for $\fal$ are called the \emph{left zero factorization axioms}.  The three equalities for $\far$ are called the \emph{right zero factorization axioms}.
\item[The Unit Factorization Axiom]\index{ring category!unit factorization axiom}\index{unit factorization axiom} 
\begin{equation}\label{ringcataxiomunitfact}
\begin{split}
\fal_{A,B,\tu} &= 1_{A \oplus B}\\
\far_{\tu,B,C} &= 1_{B \oplus C}
\end{split}
\end{equation}
These are called, respectively, the \emph{left} and the \emph{right} unit factorization axioms.
\item[The Braiding Factorization Axiom]\index{ring category!symmetry factorization axiom}\index{braiding factorization axiom} 
The following two diagrams in $\C$ commute.
\begin{equation}\label{ringcataxiomfactsyml}
\begin{tikzcd}
AC \oplus BC \ar{d}[swap]{\xiplus} \ar{r}{\fal} & (A \oplus B)C \ar{d}{\xiplus 1_C}\\
BC \oplus AC \ar{r}{\fal} & (B \oplus A)C
\end{tikzcd}\qquad
\begin{tikzcd}
AB \oplus AC \ar{d}[swap]{\xiplus} \ar{r}{\far} & A(B \oplus C) \ar{d}{1_A \xiplus}\\
AC \oplus AB \ar{r}{\far} & A(C \oplus B)
\end{tikzcd}
\end{equation}
These are called, respectively, the \emph{left} and the \emph{right} symmetry factorization axioms.
\item[The Internal Factorization Axiom]\index{ring category!internal factorization axiom}\index{internal factorization axiom}  
The following two diagrams in $\C$ commute.
\begin{equation}\label{ringcataxiomifa}
\begin{tikzcd}[cells={nodes={scale=.8}},
every label/.append style={scale=.8}]
AB \oplus A'B \oplus A''B \ar{d}[swap]{1 \oplus \fal} \ar{r}{\fal \oplus 1} & (A \oplus A')B \oplus A''B \ar{d}{\fal}\\
AB \oplus (A' \oplus A'')B \ar{r}{\fal} & (A \oplus A' \oplus A'')B
\end{tikzcd}\quad
\begin{tikzcd}[cells={nodes={scale=.8}},
every label/.append style={scale=.8}]
AB \oplus AB' \oplus AB'' \ar{d}[swap]{1 \oplus \far} \ar{r}{\far \oplus 1} & A(B \oplus B') \oplus AB'' \ar{d}{\far}\\
AB \oplus A(B' \oplus B'') \ar{r}{\far} & A(B \oplus B' \oplus B")
\end{tikzcd}
\end{equation}
These are called, respectively, the \emph{left} and the \emph{right} internal factorization axioms.
\item[The External Factorization Axiom]\index{ring category!external factorization axiom}\index{external factorization axiom} 
The three diagrams in $\C$ below commute.
\begin{equation}\label{ringcataxiomefai} 
\begin{tikzcd}[column sep=huge]
ABC \oplus A'BC \ar{d}[swap]{\fal_{AB,A'B,C}} \ar{r}{\fal_{A,A',BC}} & (A \oplus A')BC \ar[equal]{d}\\
(AB \oplus A'B)C \ar{r}{\fal_{A,A',B} 1_C} & (A \oplus A')BC
\end{tikzcd}
\end{equation}
\begin{equation}\label{ringcataxiomefaii}
\begin{tikzcd}[column sep=huge]
ABC \oplus AB'C \ar{d}[swap]{\far_{A,BC,B'C}} \ar{r}{\fal_{AB,AB',C}} & (AB \oplus AB')C \ar{d}{\far 1_C}\\
A(BC \oplus B'C) \ar{r}{1_A \fal_{B,B',C}} & A(B \oplus B')C
\end{tikzcd}
\end{equation}
\begin{equation}\label{ringcataxiomefaiii}
\begin{tikzcd}[column sep=huge]
ABC \oplus ABC' \ar{d}[swap]{\far_{A,BC,BC'}} \ar{r}{\far_{AB,C,C'}} & AB(C \oplus C') \ar[equal]{d}\\
A(BC \oplus BC') \ar{r}{1_A \far_{B,C,C'}} & AB(C \oplus C')
\end{tikzcd}
\end{equation}
These are called, respectively, the \emph{left}, the \emph{middle}, and the \emph{right} external factorization axioms.
\item[The 2-By-2 Factorization Axiom]
The following diagram in $\C$ commutes.\index{ring category!2-by-2 factorization axiom}\index{2-by-2 factorization axiom}
\begin{equation}\label{ringcataxiomtwotwo}
\begin{tikzpicture}[xscale=3,yscale=.8,baseline={(x2.base)}]
\def\h{1} \def\v{-1}
\tikzset{0cell/.append style={nodes={scale=.8}}}
\tikzset{1cell/.append style={nodes={scale=.9}}}
\draw[0cell] 
(0,0) node (x0) {A(B \oplus B') \oplus A'(B \oplus B')}
(x0)++(-\h,\v) node (x1) {AB \oplus AB' \oplus A'B \oplus A'B'}
(x1)++(1.3*\h,\v) node (x2) {(A \oplus A')(B \oplus B')}
(x1)++(0,2*\v) node (x3) {AB \oplus A'B \oplus AB' \oplus A'B'}
(x3)++(\h,\v) node (x4) {(A \oplus A')B \oplus (A \oplus A')B'}
;
\draw[1cell] 
(x1) edge[transform canvas={xshift={-2ex}}] node[pos=.3] {\far \oplus \far} (x0)
(x0) edge node {\fal} (x2)
(x1) edge node[swap] {1 \oplus \xiplus \oplus 1} (x3)
(x3) edge[transform canvas={xshift={-2ex}}] node[swap,pos=.3] {\fal \oplus \fal} (x4)
(x4) edge node[swap] {\far} (x2)
; 
\end{tikzpicture}
\end{equation}
\end{description}
This finishes the definition of a ring category.

Moreover, a ring category as above is said to be
\begin{itemize}
\item \emph{small}\index{small!ring category}\index{ring category!small} if it has a set of objects and
\item \emph{tight}\index{tight!ring category}\index{ring category!tight} if $\fal$ and $\far$ in \eqref{ringcatfactorization} are natural isomorphisms.\defmark
\end{itemize}
\end{definition}

The following concept of a bipermutative category is due to \cite[Def.\! 3.6]{elmendorf-mandell}, as presented in \cite[Section 9.3]{cerberusII}.

\begin{definition}\label{def:embipermutativecat}\index{bipermutative category}\index{category!bipermutative}\index{ring category!bipermutative category}
A \emph{bipermutative category} is a tuple
\[\big(\C,(\oplus,\zero,\xiplus),(\otimes,\tu,\xitimes),(\fal,\far)\big)\]
consisting of the following data.
\begin{itemize}
\item The tuple
\[\big(\C,(\oplus,\zero,\xiplus),(\otimes,\tu),(\fal,\far)\big)\]
is a ring category (\Cref{def:ringcat}).
\item The tuple \index{multiplicative structure!bipermutative category}\index{bipermutative category!multiplicative structure}\label{not:Ctebiperm}$\Cte = (\C,\otimes,\tu,\xitimes)$ is a permutative category (\cref{def:symmoncat}), with $\xitimes$ called the \index{braiding!multiplicative}\index{multiplicative braiding}\emph{multiplicative braiding}.
\end{itemize}
The data above are required to satisfy the following two axioms for objects $A,B,C \in \C$.
\begin{description}
\item[The Zero Braiding Axiom] 
There is an equality of morphisms\index{zero braiding axiom}
\begin{equation}\label{embipermutativeaxiomi}
\xitimes_{A,\zero} = 1_{\zero} \cn A \otimes \zero = \zero \to \zero = \zero \otimes A.
\end{equation}
\item[The Multiplicative Braiding Factorization Axiom] 
The following diagram commutes.\index{multiplicative braiding!factorization axiom}
\begin{equation}\label{embipermutativeaxiomii}
\begin{tikzcd}[column sep=huge]
(A \otimes C) \oplus (B \otimes C) \ar{d}[swap]{\xitimes_{A,C} \oplus \xitimes_{B,C}} \ar{r}{\fal_{A,B,C}} & (A \oplus B) \otimes C \ar{d}{\xitimes_{A \oplus B,C}}\\
(C \otimes A) \oplus (C \otimes B) \ar{r}{\far_{C,A,B}} & C \otimes (A \oplus B)
\end{tikzcd}
\end{equation}
\end{description}
This finishes the definition of a bipermutative category.  A bipermutative category is \index{small!bipermutative category}\index{bipermutative category!small}\emph{small}, respectively \index{tight!bipermutative category}\index{bipermutative category!tight}\emph{tight}, if the underlying ring category is so.
\end{definition}

\begin{explanation}\label{expl:embipermutativecat}
Consider \Cref{def:embipermutativecat} of a bipermutative category.
\begin{enumerate}
\item A bipermutative category $\C$ has 
\begin{itemize}
\item a permutative category $(\C,\oplus,\zero,\xiplus)$ as its additive structure and
\item a permutative category $(\C,\otimes,\tu,\xitimes)$ as its multiplicative structure.
\end{itemize}
This justifies the name \emph{bipermutative} category.
\item By the symmetry axiom \cref{symmoncatsymhexagon} in the permutative category $(\C,\xitimes)$, the zero braiding axiom \cref{embipermutativeaxiomi} is equivalent to the equality
\[\xitimes_{\zero,A} = 1_{\zero} \cn \zero \otimes A = \zero \to \zero = A \otimes \zero.\]
of morphisms.\defmark
\end{enumerate}
\end{explanation}

Bipermutative categories do \emph{not} form a subclass of symmetric bimonoidal categories, or vice versa, because
\begin{itemize}
\item the distributivity morphisms $\deltal$ and $\deltar$ in \cref{sbc-distributivity} and
\item the factorization morphisms $\fal$ and $\far$ in \cref{ringcatfactorization}
\end{itemize}
are not invertible in general.  On the other hand, \emph{tight} bipermutative categories---that is, those with $\fal$ and $\far$ invertible---form a subclass of tight symmetric bimonoidal categories.  The following result is \cite[9.3.7]{cerberusII}.

\begin{theorem}\label{thm:smcbipermuative}\index{bipermutative category!tight!tight symmetric bimonoidal category}\index{tight!symmetric bimonoidal category!tight bipermutative categories}\index{tight!bipermutative category!tight symmetric bimonoidal category}
There is a canonical bijective correspondence between
\begin{enumerate}
\item\label{tightbipermutativesbcI} the class of tight bipermutative categories in \Cref{def:embipermutativecat} and
\item\label{tightbipermutativesbcII} the class of tight symmetric bimonoidal categories in \Cref{def:sbc} with
\begin{itemize}
\item a permutative category as the additive structure,
\item a permutative category as the multiplicative structure, and
\item $\lambdadot = 1$ and $\rhodot = 1$.
\end{itemize}
\end{enumerate}
The correspondence between
\begin{itemize}
\item the factorization morphisms $\fal$ and $\far$ in \cref{ringcatfactorization} and
\item the distributivity morphisms $\deltal$ and $\deltar$ in \cref{sbc-distributivity}
\end{itemize}  
is given by the equalities
\begin{equation}\label{distributivityfactorizationsy}
\deltal = (\far)^{\inv} \andspace \deltar = (\fal)^{\inv}.
\end{equation}
\end{theorem}

In the bijective correspondence in \cref{thm:smcbipermuative}, the additive structure $(\C,\oplus,\zero,\xiplus)$ and the multiplicative structure $(\C,\otimes, \tu,\xitimes)$ remain unchanged.  Using \cref{thm:smcbipermuative} we identify tight bipermutative categories with the corresponding subclass of tight symmetric bimonoidal categories with distributivity morphisms given by \cref{distributivityfactorizationsy}.

\section{Examples of Bipermutative Categories}
\label{sec:examples}

In this section we discuss several small tight symmetric bimonoidal categories that are also small tight bipermutative categories by \cref{thm:smcbipermuative}.
\begin{enumerate}
\item In \cref{ex:finsk} the underlying category is $\Finsk$, which consists of the finite sets $\{1,\ldots,n\}$ and functions.  It is the most basic example in this section.  Each subsequent example reuses its symmetric bimonoidal structure.
\item In \cref{ex:Fset} the underlying category is the subcategory $\Fset$ of $\Finsk$ consisting of permutations.
\item $\Fskel$ in \cref{ex:Fskel} is the pointed variant of $\Finsk$ with pointed finite sets and pointed functions.  This is the indexing category for $\Ga$-categories (\cref{def:gammacategory}) and is a fundamental part of Segal $K$-theory; see \cite[Ch.\! 8]{cerberusIII}.
\item \cref{ex:mandellcategory} describes Mandell's category $\cA$ from \cite{mandell_inverseK} as a symmetric bimonoidal category.  This structure on $\cA$ is the ultimate reason why inverse $K$-theory is a pseudo symmetric $\Cat$-multifunctor, as we discuss in \cref{ch:invK}.  
\end{enumerate}
Moreover, in each of these small tight symmetric bimonoidal categories the following statements hold:
\begin{itemize}
\item Both the additive structure and the multiplicative structure are permutative categories.
\item The left multiplicative zero, the right multiplicative zero, and the left distributivity morphism are identities.
\item The right distributivity morphism is a natural isomorphism that is not the identity.
\end{itemize}
By \cref{thm:smcbipermuative} each of these small tight symmetric bimonoidal categories is also a small tight bipermutative category with factorization morphisms given by
\begin{equation}\label{factdist}
\fal = (\deltar)^\inv \andspace \far = (\deltal)^\inv.
\end{equation}

The symmetric group\index{symmetric group} on $n$ letters is denoted $\Sigma_n$.

\begin{example}[Finite Sets]\label{ex:finsk}\index{symmetric bimonoidal category!finite sets}\index{bipermutative category!finite sets}
There is a small tight symmetric bimonoidal category $\Finsk$ given by the following data.  The superscript $\mathsf{sk}$ indicates that the category $\Finsk$ is a skeleton of the category of all finite sets and functions.
\begin{description}
\item[Category] Its set of objects is $\bbN$, the set of natural numbers.  For each $n \in \bbN$, denote by
\begin{equation}\label{unpointedfs}
\ufs{n} = \begin{cases}
\{1,\ldots,n\} & \text{if $n \geq 1$},\\
\emptyset & \text{if $n=0$}
\end{cases} 
\end{equation}
the unpointed finite set with $n$ elements.  To simplify the presentation, we identify the object $n \in \Finsk$ with $\ufs{n}$.  The morphisms in $\Finsk$ are functions of unpointed finite sets:
\[\Finsk(\ufs{m}, \ufs{n}) = \Set(\ufs{m}, \ufs{n}).\]
Composition of morphisms is given by composition of functions.
\item[Additive Structure] It is given on objects by taking coproducts:
\[\ufs{m} \oplus \ufs{n} = \ufs{m+n} = \ufs{m} \bincoprod \ufs{n}.\]
\begin{itemize}
\item It is extended to morphisms using the coproduct. 
\item The sum $\oplus$ is strictly associative with strict additive zero given by the empty set $\ufs{0} = \emptyset$.
\end{itemize}
The additive braiding
\begin{equation}\label{fset-xiplus}
\begin{tikzcd}[column sep=large]
\ufs{m} \oplus \ufs{n} \ar{r}{\xiplus_{m,n}} & \ufs{n} \oplus \ufs{m} 
\end{tikzcd}
\end{equation}
is the \index{block swapping}\emph{block swapping} in $\Sigma_{m+n}$ given by
\[\xiplus_{m,n}(k) = \begin{cases}
n+k & \text{if $1 \leq k \leq m$},\\
k-m & \text{if $m+1 \leq k \leq m+n$}.
\end{cases}\]
It swaps the first $m$ elements with the last $n$ elements.
\item[Multiplicative Structure] It is given on objects by taking products with the lexicographic ordering that first compares the second entry: 
\begin{equation}\label{fset-mtensorn}
\ufs{m} \otimes \ufs{n} = \ufs{mn} = \ang{\ang{i + (j-1)m}_{i=1}^m}_{j=1}^n.
\end{equation}
\begin{itemize}
\item It is extended to morphisms using the correspondence 
\[\begin{tikzcd}
\ufs{m} \times \ufs{n} \ni (i,j) \ar[<->]{r}{\iso} & i + (j-1)m \in \ufs{mn}.
\end{tikzcd}\]
\item The product $\otimes$ is strictly associative with strict multiplicative unit $\ufs{1}$.
\end{itemize}
The multiplicative braiding
\begin{equation}\label{fset-xitimes}\begin{tikzcd}[column sep=large]
\ufs{m} \otimes \ufs{n} \ar{r}{\xitimes_{m,n}} & \ufs{n} \otimes \ufs{m} 
\end{tikzcd}
\end{equation}
is the \index{transpose permutation}\index{permutation!transpose}\emph{transpose permutation} in $\Sigma_{mn}$ given by
\[\xitimes_{m,n}\big(i + (j-1)m\big) = j + (i-1)n\]
for $1 \leq i \leq m$ and $1 \leq j \leq n$.  Regarding $\ufs{m} \otimes \ufs{n}$ as an $n \times m$ matrix with $(j,i)$-entry $i + (j-1)m$, the multiplicative braiding $\xitimes_{m,n}$ in \cref{fset-xitimes} corresponds to taking the transpose. 
\item[Multiplicative Zeros] Both $\lambdadot$ and $\rhodot$ are identity natural transformations.
\item[Distributivity] The left distributivity morphism is the identity function
\begin{equation}\label{fset-deltal}
\begin{tikzcd}[column sep=large]
\ufs{m} \otimes (\ufs{n} \oplus \ufs{p}) \ar{r}{\deltal_{m,n,p}}[swap]{=} & (\ufs{m} \otimes \ufs{n}) \oplus (\ufs{m} \otimes \ufs{p}).
\end{tikzcd}
\end{equation}
The right distributivity morphism
\begin{equation}\label{fset-deltar}
\begin{tikzcd}[column sep=large]
(\ufs{m} \oplus \ufs{n}) \otimes \ufs{p} \ar{r}{\deltar_{m,n,p}}[swap]{\iso} & (\ufs{m} \otimes \ufs{p}) \oplus (\ufs{n} \otimes \ufs{p})
\end{tikzcd}
\end{equation}
is the \emph{tetris permutation} in $\Sigma_{(m+n)p}$ given by
\begin{equation}\label{fset-deltarformula}
\left\{\begin{split}
\deltar_{m,n,p}\big(i+(k-1)(m+n)\big) &= i+(k-1)m \\
\deltar_{m,n,p}\big(j+m+(k-1)(m+n)\big) &= j+(k-1)n+pm
\end{split}\right.
\end{equation}
for $1\leq i \leq m$, $1\leq j \leq n$, and $1\leq k \leq p$.  Regarding $(\ufs{m} \oplus \ufs{n}) \otimes \ufs{p}$ as a $p \times (m+n)$ matrix of the form $[A \,|\, B]$ with
\begin{itemize}
\item $A$ the $p \times m$ matrix consisting of the first $m$ columns and
\item $B$ the $p \times n$ matrix consisting of the last $n$ columns,
\end{itemize} 
$\deltar_{m,n,p}$ rearranges the matrix $[A \,|\, B]$ to the array $\left[\frac{A}{B}\right]$, leaving the order of the entries in each of $A$ and $B$ unchanged.  The two rows in \cref{fset-deltarformula} correspond to the action of $\deltar_{m,n,p}$ on, respectively, the $(k,i)$-entry in $A$ and the $(k,j)$-entry in $B$.
\end{description} 
This finishes the construction of $\Finsk$.

All the symmetric bimonoidal category axioms in \cref{def:sbc} for $\Finsk$ follow from 
\begin{itemize}
\item the fact that $\lambdadot$, $\rhodot$, and $\deltal$ are identity natural transformations,
\item the fact that $\xiplus$ is the block swapping, and
\item the matrix description of each of $\xitimes$ and $\deltar$.
\end{itemize}
For example, the axiom \cref{laplaza-II} for $\Finsk$ states the commutativity of the following diagram for $p,q,r$ in $\bbN$.
\begin{equation}\label{laplazaIIFinsk}
\begin{tikzpicture}[xscale=4,yscale=1.2,vcenter]
\draw[0cell=.9]
(0,0) node (x11) {(\ufs{p} \oplus \ufs{r}) \otimes \ufs{q}}
(x11)++(1,0) node (x12) {(\ufs{p} \otimes \ufs{q}) \oplus (\ufs{r} \otimes \ufs{q})}
(x11)++(0,-1) node (x21) {\ufs{q} \otimes (\ufs{p}\oplus \ufs{r})}
(x12)++(0,-1) node (x22) {(\ufs{q} \otimes \ufs{p}) \oplus (\ufs{q} \otimes \ufs{r})} 
;
\draw[1cell=.9]  
(x11) edge node {\deltar_{p,r,q}} (x12)
(x12) edge node {\xitimes_{p,q} \oplus \xitimes_{r,q}} (x22)
(x11) edge node[swap] {\xitimes_{p+r,q}} (x21)
(x21) edge node {\deltal_{q,p,r}} node[swap] {=} (x22)
;
\end{tikzpicture}
\end{equation}
As above we regard $(\ufs{p} \oplus \ufs{r}) \otimes \ufs{q}$ as a $q \times (p+r)$ matrix of the form $[A \,|\, B]$ with
\begin{itemize}
\item $A$ the $q \times p$ matrix consisting of the first $p$ columns and
\item $B$ the $q \times r$ matrix consisting of the last $r$ columns.
\end{itemize} 
Then $\xitimes_{p+r,q}$ in \cref{laplazaIIFinsk} transposes the matrix $[A \,|\, B]$ to the matrix $\left[\frac{A^T}{B^T}\right]$, with $A^T$ and $B^T$ denoting the transposes of, respectively, $A$ and $B$.  On the other hand, the top-right composite in \cref{laplazaIIFinsk} rearranges the matrix $[A \,|\, B]$  
\begin{itemize}
\item first to the array $\left[\frac{A}{B}\right]$ and
\item then to the matrix $\left[\frac{A^T}{B^T}\right]$ by taking the transpose of each of $A$ and $B$.
\end{itemize}
This proves that \cref{laplazaIIFinsk} is commutative.  The axioms \cref{laplaza-I,laplaza-III,laplaza-IV,laplaza-VII,laplaza-VIII,laplaza-IX} are checked in a similar manner.  In every other axiom in \cref{def:sbc} for $\Finsk$, each morphism is the identity morphism.

Both the additive structure and the multiplicative structure of $\Finsk$ are permutative categories, with both distributivity morphisms invertible and strict multiplicative zeros.  By \cref{thm:smcbipermuative} $\Finsk$ is a small tight bipermutative category with factorization morphisms as in \cref{factdist}.  
\end{example}

Restricting the morphisms in $\Finsk$ to just the permutations, we obtain the following example.

\begin{example}[Symmetric Groups]\label{ex:Fset}\index{symmetric bimonoidal category!symmetric groups}\index{bipermutative category!symmetric groups}
There is a small tight symmetric bimonoidal category $\Fset$ given by the following data.
\begin{description}
\item[Category] Its set of objects is $\bbN$, the set of natural numbers.  For $m,n \in \bbN$, the morphism set is
\[\Fset(m,n) = \begin{cases}
\Sigma_n & \text{if $m=n$},\\
\emptyset & \text{if $m \neq n$}.
\end{cases}\]
\item[Symmetric Bimonoidal Structure] It is the restriction of the symmetric bimonoidal structure of $\Finsk$ in \cref{ex:finsk} to permutations.  This is well defined because
\begin{itemize}
\item $\xiplus$ \cref{fset-xiplus}, $\xitimes$ \cref{fset-xitimes}, and $\deltar$ \cref{fset-deltar} are permutations, and
\item $\lambdadot$, $\rhodot$, and $\deltal$ \cref{fset-deltal} are identity functions.
\end{itemize} 
\end{description}
This example is discussed in detail in \cite[Section 2.4]{cerberusI}.  By \cref{thm:smcbipermuative} $\Fset$ is a small tight bipermutative category with factorization morphisms as in \cref{factdist}.
\end{example}

The next example is a pointed variant of $\Finsk$.  It is used to define $\Ga$-categories (\cref{def:gammacategory}) in the discussion of inverse $K$-theory in \cref{ch:multifunctorA,ch:invK}.

\begin{example}[Pointed Finite Sets]\label{ex:Fskel}\index{symmetric bimonoidal category!pointed finite sets}\index{bipermutative category!pointed finite sets}
There is a small tight symmetric bimonoidal category $\Fskel$ given by the following data.
\begin{description}
\item[Category] Its objects are the \index{pointed!finite set}\emph{pointed} finite sets
\begin{equation}\label{ordn}
\ord{n} = \{0,\ldots,n\} \forspace n \geq 0
\end{equation}
with basepoint 0.  Its morphisms are \index{pointed!function}\emph{pointed} functions, that is, functions preserving the basepoint 0.
\item[Additive Structure]
It is given on objects by the wedge
\[\ord{m} \vee \ord{n} = \ord{m+n}.\]
\begin{itemize}
\item It is the quotient of the disjoint union with the basepoints in $\ord{m}$ and $\ord{n}$ identified. 
\item It is extended to morphisms using the wedge.
\item The wedge is strictly associative with strict additive zero given by $\ord{0}$.
\end{itemize}
The additive braiding
\[\begin{tikzcd}[column sep=large]
\ord{m} \vee \ord{n} = \ord{m+n} \ar{r}{\xiwedge_{\ord{m},\ord{n}}} & \ord{n+m} = \ord{n} \vee \ord{m} 
\end{tikzcd}\]
sends 0 to 0 and is given by the block swapping $\xiplus_{m,n}$ in $\Finsk$ \cref{fset-xiplus} for nonzero elements in $\ord{m+n}$.
\item[Multiplicative Structure]
It is given on objects by the smash product
\[\ord{m} \sma \ord{n} = \ord{mn}.\]
\begin{itemize}
\item It is the quotient of the product $\ord{m} \times \ord{n}$ with $(x,0)$ and $(0,y)$ identified with $0 \in \ord{mn}$ for $x \in \ord{m}$ and $y \in \ord{n}$. 
\item For $1 \leq i \leq m$ and $1 \leq j \leq n$, we use the identification
\[\ord{m} \sma \ord{n} \ni (i,j)  = i + (j-1)m \in \ord{mn}\]
as in \cref{fset-mtensorn}.
\item It is extended to pointed functions using the basepoint-preserving condition and the expression $i + (j-1)m$ for non-basepoint elements.
\item The smash product is strictly associative with strict multiplicative unit given by $\ord{1}$.
\end{itemize}
The multiplicative braiding
\begin{equation}\label{xisma}
\begin{tikzcd}[column sep=large]
\ord{m} \sma \ord{n} = \ord{mn} \ar{r}{\xisma_{\ord{m},\ord{n}}} & \ord{nm} = \ord{n} \sma \ord{m} 
\end{tikzcd}
\end{equation}
sends 0 to 0 and is given by the transpose permutation $\xitimes_{m,n}$ in $\Finsk$ \cref{fset-xitimes} for nonzero elements in $\ord{mn}$.
\item[Multiplicative Zeros]
Both $\lambdadot$ and $\rhodot$ are identity natural transformations.
\item[Distributivity]
The left distributivity morphism is the identity function
\begin{equation}\label{Fskeldeltal}
\begin{tikzcd}[column sep=large]
\ord{m} \sma (\ord{n} \vee \ord{p}) \ar{r}{\deltal_{\ord{m},\ord{n},\ord{p}}}[swap]{=} 
& (\ord{m} \sma \ord{n}) \vee (\ord{m} \sma \ord{p}).
\end{tikzcd}
\end{equation}
The right distributivity morphism
\begin{equation}\label{Fskeldeltar}
\begin{tikzcd}[column sep=large]
(\ord{m} \vee \ord{n}) \sma \ord{p} \ar{r}{\deltar_{\ord{m},\ord{n},\ord{p}}}[swap]{\iso} 
& (\ord{m} \sma \ord{p}) \vee (\ord{n} \sma \ord{p})
\end{tikzcd}
\end{equation}
sends 0 to 0 and is given by $\deltar_{m,n,p}$ in $\Finsk$ \cref{fset-deltar} for nonzero elements in $\ord{(m+n)p}$.
\end{description}
The symmetric bimonoidal category axioms of $\Fskel$ follow from those of $\Finsk$ in \cref{ex:finsk}.  By \cref{thm:smcbipermuative} $\Fskel$ is a small tight bipermutative category with factorization morphisms as in \cref{factdist}.  In \cite[Ch.\! 8]{cerberusIII} the multiplicative structure $(\Fskel,\sma,\ord{1})$ is used to define Segal $K$-theory.  However, in \cite[8.1.5]{cerberusIII} the smash product $\ord{m} \sma \ord{n}$ uses the lexicographic order that compares the first entry first.  Here $\ord{m} \sma \ord{n}$ is defined with the lexicographic order that compares the second entry first.  
\end{example}

The next example is crucial for our discussion of inverse $K$-theory in \cref{ch:multifunctorA,ch:invK}.

\begin{example}[Mandell's Category for Inverse $K$-Theory]\label{ex:mandellcategory}\index{symmetric bimonoidal category!inverse $K$-theory}\index{bipermutative category!inverse $K$-theory}\index{inverse K-theory@inverse $K$-theory!bipermutative category}
There is a small tight symmetric bimonoidal category $\cA$ given by the following data.
\begin{description}
\item[Category]
An object in $\cA$ is a finite sequence of positive integers.  The empty sequence is denoted by either \label{not:emptyseq}$\ang{}$ or $\emptyset$.  With the notation in \cref{unpointedfs} for unpointed finite sets, a morphism 
\begin{equation}\label{Amorphismphi}
\begin{tikzcd}[column sep=large]
m = (m_1, \ldots, m_p) \ar{r}{\phi} & (n_1, \ldots, n_q) = n
\end{tikzcd} \inspace \cA
\end{equation}
is a function of unpointed finite sets
\begin{equation}\label{phiufs}
\begin{tikzcd}[column sep=large]
\coprod\limits_{i=1}^p \ufs{m_i} \ar{r}{\phi} & \coprod\limits_{j=1}^q \ufs{n_j}
\end{tikzcd}
\end{equation}
such that the preimage of each $\ufs{n_j}$ is either empty or contained in a single $\ufs{m_i}$, with the index $i$ determined by the index $j$.  Note that, however, each $\ufs{m_i}$ may hit several different $\ufs{n_j}$.  The category $\cA$ is Mandell's indexing category for inverse $K$-theory \cite{mandell_inverseK}; see also \cite{johnson-yau-invK} and \cref{ch:multifunctorA,ch:invK}.
\item[Additive Structure]
With the objects $m$ and $n$ as in \cref{Amorphismphi}, the sum is given by concatenation of finite sequences:
\[m \oplus n = (m_1,\ldots,m_p, n_1,\ldots,n_q).\]
\begin{itemize}
\item It is extended to morphisms by taking coproducts of functions. 
\item The sum $\oplus$ is strictly associative with strict additive zero given by the empty sequence $\ang{}$.
\end{itemize}
The additive braiding
\begin{equation}\label{Abetaplusmn}
\begin{tikzcd}[column sep=large]
m \oplus n \ar{r}{\betaplus_{m,n}} & n \oplus m
\end{tikzcd}
\end{equation}
is the block swapping in $\Sigma_{\mathtt{m}+\mathtt{n}}$ in \cref{fset-xiplus}
\[\begin{tikzcd}[column sep=large]
\ufs{\mathtt{m}+\mathtt{n}} = \Big(\coprod\limits_{i=1}^p \ufs{m_i}\Big) \coprod \Big(\coprod\limits_{j=1}^q \ufs{n_j}\Big) \ar{r}{\xiplus_{\mathtt{m},\mathtt{n}}} 
& \Big(\coprod\limits_{j=1}^q \ufs{n_j}\Big) \coprod \Big(\coprod\limits_{i=1}^p \ufs{m_i}\Big) = \ufs{\mathtt{n}+\mathtt{m}}
\end{tikzcd}\]
where 
\begin{equation}\label{summn}
\mathtt{m} = m_1 + \cdots + m_p \andspace \mathtt{n} = n_1 + \cdots + n_q.
\end{equation}
Equivalently, 
\[\betaplus_{m,n} = \ufs{\xiplus_{p,q}}\]
is the block permutation induced by the block swapping $\xiplus_{p,q} \in \Sigma_{p+q}$ in \cref{fset-xiplus} that permutes the $p+q$ blocks of lengths $\{m_1,\ldots,m_p,n_1,\ldots,n_q\}$.
\item[Multiplicative Structure]
With the objects $m$ and $n$ as in \cref{Amorphismphi}, the product is given by
\begin{equation}\label{Amtensorn}
m \otimes n = \ang{\ang{m_i n_j}_{i=1}^p}_{j=1}^q.
\end{equation}
We visualize $m \otimes n$ as the $q \times p$ matrix with $(j,i)$-entry given by $m_i n_j$, which is also regarded as the unpointed finite set $\ufs{m_i} \times \ufs{n_j} = \ufs{m_i n_j}$.  

To define $\otimes$ on morphisms, suppose given morphisms in $\cA$ as follows.
\begin{equation}\label{Aphipsi}
\begin{tikzcd}[column sep=large, row sep=0ex,
/tikz/column 1/.append style={anchor=base east},
/tikz/column 2/.append style={anchor=base west}]
m = (m_1, \ldots, m_p) \ar{r}{\phi} & (m_1', \ldots, m_r') = m'\\
n = (n_1, \ldots, n_q) \ar{r}{\psi} & (n_1', \ldots, n_s') = n'
\end{tikzcd}
\end{equation}
The morphism
\[\begin{tikzcd}[column sep=large]
\ang{\ang{m_i n_j}_i}_j = m \otimes n \ar{r}{\phi \otimes \psi} 
& m' \otimes n' = \ang{\ang{m_k' n_\ell'}_k}_\ell
\end{tikzcd} \inspace \cA\]
is the function
\[\begin{tikzcd}[column sep=large]
\coprod\limits_{j=1}^q \coprod\limits_{i=1}^p \ufs{m_i n_j} \ar{r}{(\phi,\psi)} &
\coprod\limits_{\ell=1}^s \coprod\limits_{k=1}^r \ufs{m_k' n_\ell'} 
\end{tikzcd}\]
given by
\[(\phi,\psi)(x,y) = (\phi x, \psi y) \forspace (x,y) \in \ufs{m_i n_j}.\]
In this definition, each of $\ufs{m_i n_j}$ and $\ufs{m_k' n_\ell'}$ is given the lexicographic ordering in \cref{fset-mtensorn} that compares the second entry first.

The product $\otimes$ is strictly associative with strict multiplicative unit given by the length-1 sequence $(1)$.  The multiplicative braiding
\begin{equation}\label{Abetate}
\begin{tikzcd}[column sep=large]
\ang{\ang{m_i n_j}_i}_j = m \otimes n \ar{r}{\betatimes_{m,n}} & n \otimes m = \ang{\ang{n_j m_i}_j}_i
\end{tikzcd}
\end{equation}
is the composite bijection
\[\begin{tikzpicture}[xscale=4,yscale=1.2,vcenter]
\draw[0cell=.9]
(0,0) node (x11) {\textstyle\coprod_{j=1}^q \coprod_{i=1}^p \ufs{m_i n_j}}
(x11)++(1,0) node (x12) {\textstyle\coprod_{i=1}^p \coprod_{j=1}^q \ufs{n_j m_i}}
(x11)++(.5,-1) node (x2) {\textstyle\coprod_{i=1}^p \coprod_{j=1}^q \ufs{m_i n_j}}
;
\draw[1cell=.9]  
(x11) edge node {\betatimes_{m,n}} (x12)
(x11) edge[transform canvas={xshift={-2ex}}] node[swap,pos=.2] {\ufs{\xitimes_{p,q}}} (x2)
(x2) edge[transform canvas={xshift={2ex}}] node[swap,pos=.7] {\textstyle \coprod_{i,j} \xitimes_{m_i,n_j}} (x12)
;
\end{tikzpicture}\]
given by
\[\betatimes_{m,n}(x,y) = (y,x) \in \ufs{n_j m_i} \forspace (x,y) \in \ufs{m_i n_j}.\]
\begin{itemize}
\item With $\mathtt{m} = m_1 + \cdots + m_p$ and $\mathtt{n} = n_1 + \cdots + n_q$ as in \cref{summn}, $\ufs{\xitimes_{p,q}} \in \Sigma_{\mathtt{m} \mathtt{n}}$ is the block permutation induced by the transpose permutation $\xitimes_{p,q} \in \Sigma_{pq}$ in \cref{fset-xitimes} that permutes the $pq$ blocks of lengths $\{\{m_i n_j\}_{i=1}^p\}_{j=1}^q$.
\item For each $1 \leq i \leq p$ and each $1 \leq j \leq q$, $\xitimes_{m_i,n_j} \in \Sigma_{m_i n_j}$ is the transpose permutation in \cref{fset-xitimes}.
\end{itemize}
\item[Multiplicative Zeros]
Both $\lambdadot$ and $\rhodot$ are identity natural transformations.  This is well defined because there are object equalities
\[m \otimes \ang{} = \ang{} = \ang{} \otimes m \forspace m \in \cA.\]
\item[Distributivity]
Suppose $m$, $m'$, $n$, and $n'$ are objects in $\cA$ as in \cref{Aphipsi}.  The left distributivity morphism is the identity morphism
\begin{equation}\label{Adeltal}
\begin{tikzcd}[column sep=large]
m \otimes (n \oplus n') \ar{r}{\deltal_{m,n,n'}}[swap]{=} & (m \otimes n) \oplus (m \otimes n').
\end{tikzcd}
\end{equation}
With $\mathtt{m}$ and $\mathtt{n}$ as in \cref{summn} and $\mathtt{m'} = m_1' + \cdots + m_r'$, the right distributivity morphism 
\begin{equation}\label{Adeltar}
\begin{tikzcd}[column sep=large]
(m \oplus m') \otimes n \ar{r}{\deltar_{m,m',n}}[swap]{\iso} & (m \otimes n) \oplus (m' \otimes n)
\end{tikzcd}
\end{equation}
is the block permutation in $\Sigma_{(\mathtt{m} + \mathtt{m'}) \mathtt{n}}$
\[\begin{tikzcd}[column sep=large,cells={nodes={scale=.9}}]
\coprod\limits_{j=1}^q \bigg[\Big( \coprod\limits_{i=1}^p \ufs{m_i n_j} \Big) \coprod \Big( \coprod\limits_{k=1}^r \ufs{m_k' n_j} \Big)\bigg] \ar{r}{\ufs{\deltar_{p,r,q}}} &
\bigg( \coprod\limits_{j=1}^q \coprod\limits_{i=1}^p \ufs{m_i n_j} \bigg) \coprod  \bigg( \coprod\limits_{j=1}^q \coprod\limits_{k=1}^r \ufs{m_k' n_j} \bigg)
\end{tikzcd}\]
induced by the tetris permutation \cref{fset-deltar} 
\[\begin{tikzcd}[column sep=large]
(\ufs{p} \oplus \ufs{r}) \otimes \ufs{q} \ar{r}{\deltar_{p,r,q}}[swap]{\iso} & (\ufs{p} \otimes \ufs{q}) \oplus (\ufs{r} \otimes \ufs{q})
\end{tikzcd}\]
that permutes the $(p+r)q$ blocks of lengths
\[\big\{ \{m_i n_j\}_{i=1}^p \scs \{m_k' n_j\}_{k=1}^r \big\}_{j=1}^q.\]
\end{description}
This finishes the construction of $\cA$.  

The symmetric bimonoidal category axioms in \cref{def:sbc} for $\cA$ follow from those for $\Finsk$ in \cref{ex:finsk}, applied to induced block permutations.  More explicitly, the axiom \cref{laplaza-II} for $\cA$ states that the following diagram is commutative for objects $m,m',n \in \cA$ as in \cref{Aphipsi}.
\[\begin{tikzpicture}[xscale=4,yscale=1.3,vcenter]
\draw[0cell=.9]
(0,0) node (x11) {(m \oplus m') \otimes n}
(x11)++(1,0) node (x12) {(m \otimes n) \oplus (m' \otimes n)}
(x11)++(0,-1) node (x21) {n \otimes (m \oplus m')}
(x12)++(0,-1) node (x22) {(n \otimes m) \oplus (n \otimes m')}
;
\draw[1cell=.9]  
(x11) edge node {\deltar_{m,m',n}} (x12)
(x12) edge node {\betatimes_{m,n} \oplus \betatimes_{m',n}} (x22)
(x11) edge node[swap] {\betatimes_{m \oplus m',n}} (x21)
(x21) edge node {\deltal_{n,m,m'}} node[swap] {=} (x22)
;
\end{tikzpicture}\]
The two composites in the previous diagram are the following two composite functions, where the coproduct indices run through $i \in \{1,\ldots,p\}$, $j \in \{1, \ldots,q\}$, and $k \in \{1,\ldots,r\}$.
\[\begin{tikzpicture}[xscale=5.5,yscale=1.3,vcenter]
\draw[0cell=.85]
(0,0) node (x11) {\textstyle\coprod_j \big[\big(\coprod_i \ufs{m_i n_j}\big) \coprod \big(\coprod_k \ufs{m_k' n_j}\big)\big]}
(x11)++(1,0) node (x12) {\textstyle \big(\coprod_j \coprod_i \ufs{m_i n_j}\big) \coprod \big(\coprod_j \coprod_k \ufs{m_k' n_j} \big)}
(x11)++(.5,-.8) node (x2) {\textstyle \big(\coprod_i \coprod_j \ufs{m_i n_j}\big) \coprod \big(\coprod_k \coprod_j \ufs{m_k' n_j} \big)}
(x2)++(0,-1) node (x3) {\textstyle \big(\coprod_i \coprod_j \ufs{n_j m_i}\big) \coprod \big(\coprod_k \coprod_j \ufs{n_j m_k'} \big)}
;
\draw[1cell=.85]  
(x11) edge node {\ufs{\deltar_{p,r,q}}} (x12)
(x12) edge[transform canvas={xshift={4ex}}] node[pos=.2] {\ufs{\xitimes_{p,q}} \sqcup \ufs{\xitimes_{r,q}}} (x2)
(x11) edge[transform canvas={xshift={-4ex}}] node[swap,pos=.2] {\ufs{\xitimes_{p+r,q}}} (x2)
(x2) edge node[swap] {\textstyle(\sqcup_{i,j}\, \xitimes_{m_i,n_j}) \,\sqcup} node {(\sqcup_{k,j}\, \xitimes_{m_k',n_j})} (x3)
;
\end{tikzpicture}\]
In the previous diagram, the top triangle commutes by \cref{laplazaIIFinsk}, which is the axiom \cref{laplaza-II} for $\Finsk$, applied to induced block permutations.  This proves the axiom \cref{laplaza-II} for $\cA$.  The only other axiom that involves the multiplicative braiding is \cref{laplaza-XV}, which holds in $\cA$ because $\lambdadot$ and $\rhodot$ are both the identity natural transformations.  Every other axiom in \cref{def:sbc} for $\cA$ holds as they do for $\Finsk$, applied to induced block permutations.  By \cref{thm:smcbipermuative} $\cA$ is a small tight bipermutative category with factorization morphisms as in \cref{factdist}.
\end{example}

\section{Laplaza's Coherence Theorem}
\label{sec:laplazacoherence}

In this section we recall the statement of Laplaza's Coherence \cref{thm:laplaza-coherence-1} for symmetric bimonoidal categories \cite{laplaza}.  A detailed proof, along with corrections of the original argument in \cite{laplaza}, is in \cite[Ch.\! 3]{cerberusI}.  The statement of this coherence theorem requires some preliminary definitions about graphs, which we first discuss.

\subsection*{Elementary Graph} 

\begin{definition}\label{def:free-algebra}
Suppose $S$ is a set.  The \index{algebra!free}\index{free!algebra}\emph{free $\plustimes$-algebra} of $S$ is the set $\frees$ defined inductively by the following two conditions.
\begin{itemize}
\item $S \subset \frees$.
\item If $a,b \in \frees$, then the symbols \[a \oplus b \andspace a \otimes b\] also belong to $\frees$.  They are called, respectively, the \index{sum}\emph{sum} and the \index{product}\emph{product} of $a$ and $b$.
\end{itemize}
To simplify the presentation, we sometimes abbreviate $a \otimes b$ to $ab$.  In the absence of clarifying parentheses, $\otimes$ takes precedence over $\oplus$. 
\end{definition}

\begin{definition}\label{def:graph}
A \index{graph}\emph{graph} $G = (V,E)$ is a pair consisting of the following data. 
\begin{itemize}
\item $V$ is a class.  An element in $V$ is called a \index{vertex}\emph{vertex} in $G$.
\item $E$ assigns to each ordered pair $(u,v)$ with $u,v\in V$ a set $E(u,v)$, an element of which is called an \index{edge}\emph{edge} with \index{domain!of an edge}\emph{domain} $u$ and \index{codomain!of an edge}\emph{codomain} $v$.  We also denote such an edge by\label{not:edge}
\[u \to v, \quad e \cn u \to v, \orspace \begin{tikzcd}u \ar{r}{e} & v\end{tikzcd}\] 
if $e$ is the name of the edge.
\end{itemize}
A \index{path}\emph{path} in such a graph is a nonempty finite sequence of edges \label{not:pathe}$(e_n,\ldots,e_1)$ as in
\[\begin{tikzcd}
v_0 \ar{r}{e_1} & v_1 \ar{r}{e_2} & \cdots \ar{r}{e_n} & v_n.\end{tikzcd}\]
Such a path is denoted by \label{not:pathv}$v_0 \apath v_n$ and is said to have \index{length!of a path}\emph{length} $n$, \index{domain!of a path}domain $v_0$, and \index{codomain!of a path}codomain $v_n$.
\end{definition}

\begin{example}\label{ex:category-as-graph}
Each category $\C$ has an associated graph $(V,E)$, with $V$ the class of objects in $\C$, and $E(u,v) = \C(u,v)$ for objects $u,v \in \C$.  A nonempty finite sequence of composable morphisms in $\C$ yields a path in the associated graph.
\end{example}

\begin{definition}\label{def:elementary-graph}
Suppose $X$ is a set with two distinguished elements $\zeroofx$ and $\oneofx$, called the \index{additive!zero}\emph{additive zero} and the \index{multiplicative!unit}\emph{multiplicative unit}, respectively.  The \index{elementary!graph}\index{graph!elementary}\emph{elementary graph of $X$}, denoted by \label{not:grelx}$\grelx$, is the graph defined as follows.
\begin{description}
\item[Vertices] The set of vertices in $\grelx$ is the free $\plustimes$-algebra $\freex$ of $X$.
\item[Edges]  Edges in $\grelx$ are of the following types for all $x,y,z \in \freex$.
\begin{description}
\item[The Additive Structure] 
\[\begin{tikzcd}[column sep=large]
(x \oplus y) \oplus z \ar[shift left]{r}{\alphaplus_{x,y,z}} & x \oplus (y \oplus z) \ar[shift left]{l}{\alphaplusinv_{x,y,z}} & x \oplus y \ar[shift left]{r}{\xiplus_{x,y}} & y \oplus x \ar[shift left]{l}{\xiplusinv_{x,y}}\end{tikzcd}\]
\[\begin{tikzcd}[column sep=large]
\zeroofx \oplus x \ar[shift left]{r}{\lambdaplus_{x}} & x \ar[shift left]{l}{\lambdaplusinv_{x}} \ar[shift right]{r}[swap]{\rhoplusinv_{x}} & x \oplus \zeroofx \ar[shift right]{l}[swap]{\rhoplus_{x}}
\end{tikzcd}\]
\item[The Multiplicative Structure] 
\[\begin{tikzcd}[column sep=large]
(x \otimes y) \otimes z \ar[shift left]{r}{\alphatimes_{x,y,z}} & x \otimes (y \otimes z) \ar[shift left]{l}{\alphatimesinv_{x,y,z}} & x \otimes y \ar[shift left]{r}{\xitimes_{x,y}} & y \otimes x \ar[shift left]{l}{\xitimesinv_{x,y}}\end{tikzcd}\]
\[\begin{tikzcd}[column sep=large]
\oneofx \otimes x \ar[shift left]{r}{\lambdatimes_{x}} & x \ar[shift left]{l}{\lambdatimesinv_{x}} \ar[shift right]{r}[swap]{\rhotimesinv_{x}} & x \otimes \oneofx \ar[shift right]{l}[swap]{\rhotimes_{x}}
\end{tikzcd}\]
\item[The Multiplicative Zeros]
\[\begin{tikzcd}[column sep=large]
\zeroofx \otimes x \ar[shift left]{r}{\lambdadot_{x}} & \zeroofx \ar[shift left]{l}{\lambdadotinv_{x}} \ar[shift right]{r}[swap]{\rhodotinv_{x}} & x \otimes \zeroofx \ar[shift right]{l}[swap]{\rhodot_{x}}
\end{tikzcd}\]
\item[Distributivity]
\[\begin{tikzcd}[column sep=large,row sep=tiny
    ,/tikz/column 1/.append style={anchor=base east}
    ,/tikz/column 2/.append style={anchor=base west}]
x \otimes (y \oplus z) \ar{r}{\deltal_{x,y,z}} & (x \otimes y) \oplus (x \otimes z)\\
(x \oplus y) \otimes z \ar{r}{\deltar_{x,y,z}} & (x \otimes z) \oplus (y \otimes z)
\end{tikzcd}\]
\item[Identity]
\[\begin{tikzcd}[column sep=large]
x \ar{r}{1_x} & x\end{tikzcd}\]
\end{description}
\end{description}
This finishes the definition of $\grelx$.  

Moreover, we define the following.
\begin{itemize}
\item The set of edges in $\grelx$ is denoted by \label{not:eelx}$\eelx$, the elements of which are called \index{elementary!edge}\index{edge!elementary}\emph{elementary edges}.  
\item $\alphaplus$ and $\alphaplusinv$ are \index{formal inverse}\emph{formal inverses} of each other, and similarly for the other 9 pairs of elementary edges in the first three groups above.
\item $1_x$ is called the \index{edge!identity}\index{identity}\emph{identity} of $x$.
\item The names in \Cref{def:sbc} are reused for elementary edges.  For example, $\lambdadot$ is called the left multiplicative zero, and $\deltar$ is called the right distributivity.\defmark
\end{itemize}
\end{definition}

\subsection*{Graph}

The next two definitions will be combined to yield the graph of $X$ in \Cref{def:graph-X}.

\begin{definition}\label{def:eelfr}
With $(X,\zeroofx,\oneofx)$ as in \Cref{def:elementary-graph}, consider the free $\plustimes$-algebra \label{not:eelfrx}$\eelfrx$ of the set $\eelx$ of elementary edges.  The \index{domain!of an edge}\emph{domain} and \index{codomain!of an edge}\emph{codomain} of an element $f \in \eelfrx$ are elements in $\freex$ defined inductively as follows.
\begin{itemize}
\item For an elementary edge $f \in \eelx$, its (co)domain are those of $f$ in the elementary graph $\grelx$.
\item Suppose $f_1,f_2 \in \eelfrx$ with
\begin{itemize}
\item $u_i \in \freex$ the domain of $f_i$ and 
\item $v_i \in \freex$ the codomain of $f_i$
\end{itemize}  
already defined for $i=1,2$.  Then
\begin{itemize}
\item $f_1 \oplus f_2$ has \emph{domain} $u_1 \oplus u_2$ and \emph{codomain} $v_1 \oplus v_2$, and
\item $f_1 \otimes f_2$ has \emph{domain} $u_1 \otimes u_2$ and \emph{codomain} $v_1 \otimes v_2$.\defmark
\end{itemize} 
\end{itemize}
\end{definition}

\begin{definition}\label{def:prime-edges}
Continuing \Cref{def:eelfr}, \index{prime edge!identity}\emph{identity prime edges} and \index{prime edge!nonidentity}\emph{nonidentity prime edges} are elements in $\eelfrx$ defined inductively by the following four conditions.
\begin{itemize}
\item Elementary edges of the type $1_x$ for $x\in\freex$ are identity prime edges.
\item Elementary edges not of the type $1_x$ for $x\in\freex$ are nonidentity prime edges.  
\item If $e_1,e_2 \in \eelfrx$ are identity prime edges, then so are $e_1\oplus e_2$ and $e_1 \otimes e_2$.
\item If $f$ is a nonidentity prime edge and if $e$ is an identity prime edge, then 
\[f\oplus e, \quad e\oplus f, \quad f\otimes e, \andspace e\otimes f\]
are nonidentity prime edges.
\end{itemize}
Moreover, we define the following.
\begin{itemize}
\item A \index{edge!prime}\index{prime edge}\emph{prime edge} means either an identity prime edge or a nonidentity prime edge.  The set of prime edges is denoted by \label{not:pedgex}$\pedgex$.
\item A \emph{$\delta$-prime edge} is a prime edge that involves either $\deltal$ or $\deltar$.\defmark
\end{itemize}
\end{definition}

\begin{definition}\label{def:graph-X}
With $(X,\zeroofx,\oneofx)$ as in \Cref{def:elementary-graph}, the \index{graph!of a set}\emph{graph of $X$}, which is denoted by \label{not:grx}$\grx$, is the graph defined as follows.
\begin{description}
\item[Vertices] The set of vertices in $\grx$ is the free $\plustimes$-algebra $\freex$ of $X$.
\item[Edges] The set of edges in $\grx$ is the set $\pedgex$ of prime edges as in \Cref{def:prime-edges}, with (co)domain as in \Cref{def:eelfr}.\defmark
\end{description}
\end{definition}

\subsection*{The Value of a Path}

The graph of $X$ is interpreted in a symmetric bimonoidal category via the following concept.

\begin{definition}\label{def:graph-morphism}
Suppose $G_1 = (V_1,E_1)$ and $G_2 = (V_2,E_2)$ are graphs (\cref{def:graph}).  A \index{morphism!graph}\index{graph!morphism}\emph{graph morphism} $f \cn G_1 \to G_2$ consists of functions
\begin{itemize}
\item $f_V \cn V_1 \to V_2$ and
\item $f_E \cn E_1(u,v) \to E_2\big(f_V u, f_V v\big)$ for $u,v\in V_1$.
\end{itemize}
To simplify the notation, both $f_V$ and $f_E$ are denoted by $f$ below.
\end{definition}

\begin{definition}\label{def:graphx-to-c}
Suppose given the data $(X,\C,\vphi)$ as follows.
\begin{itemize}
\item $X$ is a set with two distinguished elements $\zeroofx$ and $\oneofx$ as in \Cref{def:elementary-graph}.
\item $\C$ is a symmetric bimonoidal category as in \Cref{def:sbc}, equipped with the graph structure in \Cref{ex:category-as-graph}.
\item $\vphi \cn X \to \Ob(\C)$ is any function such that 
\begin{equation}\label{vphi-zero-one}
\vphi(\zeroofx) = \zero \andspace \vphi(\oneofx) = \tu.
\end{equation}
\end{itemize} 
Extend $\vphi$ to a graph morphism\label{not:vphigrxc}\index{graph!morphism!extension}
\[\begin{tikzcd}[column sep=large]
\grx \ar{r}{\vphi} & \C\end{tikzcd}\] 
as follows.
\begin{description}
\item[Vertices] For $x,y\in \freex$ such that $\vphi x,\vphi y \in \Ob(\C)$ are already defined, we define
\begin{equation}\label{vphi-objects}
\begin{split}
\vphi(x\oplus y) &= \vphi x \oplus \vphi y \andspace\\ 
\vphi(x\otimes y) &= \vphi x \otimes \vphi y.\end{split}
\end{equation}
\item[Elementary Edges] $\vphi$ sends each elementary edge to the structure morphism in $\C$ with the same name and with the subscripts replaced by their images under $\vphi$.
\item[Prime Edges] If $f_1,f_2\in\pedgex$ are prime edges with at most one of them nonidentity and with $\vphi(f_1)$ and $\vphi(f_2)$ already defined, then we define
\begin{equation}\label{vphi-morphisms}
\begin{split}
\vphi(f_1 \oplus f_2) &= \vphi(f_1) \oplus \vphi(f_2) \andspace\\ 
\vphi(f_1 \otimes f_2) &= \vphi(f_1) \otimes \vphi(f_2).\end{split}
\end{equation}
\end{description}
This finishes the definition of the graph morphism $\vphi$.

Moreover, for a path $P=(f_n,\ldots,f_1)$ in $\grx$ with domain $u$ and codomain $v$, its \index{value!path}\index{path!value}\emph{value} in $\C$ is the composite 
\begin{equation}\label{value}
\vphi P = \vphi(f_n) \circ \cdots \circ \vphi(f_1) 
\end{equation}
in $\C(\vphi u, \vphi v)$.
\end{definition}

\subsection*{Regularity}

Next we introduce a restriction, called regularity, on the elements in $\freex$ for the coherence theorem of symmetric bimonoidal categories.

\begin{definition}\label{def:strict-algebra}
Suppose $(X,\zeroofx,\oneofx)$ is as in \Cref{def:elementary-graph}.  Define its \index{algebra!strict}\index{strict algebra}\emph{strict $\plustimes$-algebra} \label{not:stx}$\stx$ as the quotient set of $\freex$ in \Cref{def:free-algebra} by the smallest relation that contains the following identifications for all elements $x,y,z \in \freex$.
\begin{description}
\item[Additive Structure]
\[\begin{split}
(x \oplus y) \oplus z &= x \oplus (y \oplus z)\\
\zeroofx \oplus x &= x = x \oplus \zeroofx\\
x \oplus y &= y \oplus x\\
\end{split}\]
\item[Multiplicative Structure]
\[\begin{split}
(x \otimes y) \otimes z &= x \otimes (y \otimes z)\\
\oneofx \otimes x &= x = x \otimes \oneofx\\
x \otimes y &= y \otimes x\\
\end{split}\]
\item[Multiplicative Zeros]
\[\zeroofx \otimes x = \zeroofx = x \otimes \zeroofx\]
\item[Distributivity]
\[\begin{split}
x \otimes (y \oplus z) &= (x \otimes y) \oplus (x \otimes z)\\
(x \oplus y) \otimes z &= (x \otimes z) \oplus (y \otimes z)
\end{split}\]
\end{description}
This finishes the definition of $\stx$.  Moreover, denote by
\begin{equation}\label{support}
\begin{tikzcd}[column sep=large]
\freex \ar{r}{\supp} & \stx\end{tikzcd}
\end{equation}
the quotient map, called the \index{support}\emph{support}.
\end{definition}

\begin{definition}\label{def:regular}
An element $x\in\freex$ is \index{regular}\emph{regular} if there exist elements $x^i_j \in X$ for $1\leq i \leq m$ and $1 \leq j \leq k_i$ for each $i$ with $m,k_1,\ldots,k_m>0$, such that the following three conditions hold.
\begin{enumerate}[label=(\roman*)]
\item\label{supporti} The equality
\begin{equation}\label{supp-sum-product}
\supp(x) = \supp\bigg(\bigoplus_{i=1}^m \big(x^i_1 \otimes \cdots \otimes x^i_{k_i}\big)\bigg)
\end{equation}
holds in $\stx$.  The iterated sum $\bigoplus_{i=1}^m$ and each of the $m$ iterated products $x^i_1 \otimes \cdots \otimes x^i_{k_i}$ have some bracketings.  By the definition of $\stx$, different bracketings yield the same support.
\item\label{supportii} For each $1\leq i \leq m$, the elements $x^i_1,\ldots,x^i_{k_i} \in X$ are distinct.
\item\label{supportiii} The $m$ elements
\[\supp\big(x^i_1 \otimes \cdots \otimes x^i_{k_i}\big) \in \stx\]
for $1\leq i \leq m$ are distinct.\defmark
\end{enumerate}
\end{definition}

We are now ready to state the following coherence theorem for symmetric bimonoidal categories.  See \cite[3.9.1]{cerberusI} for a detailed proof.

\begin{theorem}[Laplaza's Coherence]\label{thm:laplaza-coherence-1}\index{Laplaza's Coherence Theorem}\index{symmetric bimonoidal category!Laplaza's Coherence Theorem}
In the context of \cref{def:graphx-to-c}, suppose $\C$ is a symmetric bimonoidal category in which the value of each \index{delta-@$\delta$-!prime edge}\index{prime edge!delta-@$\delta$-}$\delta$-prime edge is a \index{monomorphism}monomorphism.  Suppose
\[\begin{tikzcd}[column sep=large]
a \ar[shift left,twoheadrightarrow]{r}{P_1} \ar[shift right,twoheadrightarrow]{r}[swap]{P_2} & b\end{tikzcd}\]
are two \index{path}paths in $\grx$ with $a\in\freex$ \index{regular}regular.  Then the \index{value}values of $P_1$ and $P_2$ in $\C$ are equal.
\end{theorem}

\begin{example}\label{ex:laplazacoherence}
If $\C$ is a \emph{tight} symmetric bimonoidal category---that is, both distributivity morphisms $\deltal$ and $\deltar$ in $\C$ are invertible---then the value of each $\delta$-prime edge is an isomorphism, which is, in particular, a monomorphism.  In this case \cref{thm:laplaza-coherence-1} applies to $\C$.  The value in $\C$ of each path in \cref{thm:laplaza-coherence-1} is called a \emph{Laplaza coherence isomorphism}.
\begin{itemize}
\item Via \cref{thm:smcbipermuative} each tight bipermutative category\index{bipermutative category!coherence} is also a tight symmetric bimonoidal category, so Laplaza's Coherence \cref{thm:laplaza-coherence-1} applies to tight bipermutative categories. 
\item In particular, \cref{thm:laplaza-coherence-1} applies to the tight bipermutative categories $\Finsk$, $\Fset$, $\Fskel$, and $\cA$ in, respectively, \cref{ex:finsk,ex:Fset,ex:Fskel,ex:mandellcategory}.\defmark
\end{itemize} 
\end{example}

There is another coherence theorem for symmetric bimonoidal categories that is also due to Laplaza \cite{laplaza2}.  Since we will not use that coherence theorem in this work, the interested reader is referred to \cite[Ch. 4]{cerberusI} for a detailed discussion.

\chapter{Enriched Multicategories and Multiequivalences}
\label{ch:multicat}
This chapter has two main purposes:
\begin{enumerate}
\item We review enriched multicategories and related concepts that are used in subsequent chapters.
\item We characterize categorically-enriched multiequivalences in terms of essential surjectivity on objects and invertibility on multimorphism categories (\cref{thm:multiwhitehead}).  This characterization holds in both the symmetric and the non-symmetric contexts.  \cref{thm:multiwhitehead} has a pseudo symmetric variant in \cref{thm:psmultiwhitehead}.  
\end{enumerate}
Here is a summary table of some of the definitions and results in this chapter.
\begin{center}
\resizebox{\textwidth}{!}{
{\renewcommand{\arraystretch}{1.4}%
{\setlength{\tabcolsep}{1ex}
\begin{tabular}{|c|c|}\hline
enriched multicategories and multifunctors & \cref{def:enr-multicategory,def:enr-multicategory-functor} \\ \hline
enriched multinatural transformations & \cref{def:enr-multicat-natural-transformation} \\ \hline
2-categories of enriched multicategories & \cref{v-multicat-2cat} \\ \hline
$\Cat$-multicategories and multinatural transformations & \cref{def:catmulticat,expl:catmultitransformation} \\ \hline\hline
$\Cat$-multiequivalences & \cref{def:catmultiequivalence} \\ \hline
characterizations of $\Cat$-multiequivalences and 2-equivalences & \cref{thm:multiwhitehead,thm:iiequivalence} \\ \hline 
\end{tabular}}}}
\end{center}
\smallskip
In \cref{thm:dcatpfibdeq} we use \cref{thm:multiwhitehead} to show that the Grothendieck construction on a small tight bipermutative category is a non-symmetric $\Cat$-multiequivalence.  An important consequence of this fact is that inverse $K$-theory factors through a non-symmetric $\Cat$-multiequivalence \cref{invKAgroaU}.

\subsection*{Organization}

In \cref{sec:enrmulticat} we review multicategories enriched in a symmetric monoidal category $\V$.  In \cref{sec:iicatofmulticat} we review $\V$-enriched multifunctors, $\V$-enriched multinatural transformations, and the 2-category $\VMulticat$ that they form with small $\V$-multicategories as objects.  Most of the material in \cref{sec:enrmulticat,sec:iicatofmulticat} is presented with much more detail in \cite[Ch.\! 5 and 6]{cerberusIII} and \cite{yau-operad}.  

In \cref{sec:multequiv} we show that $\Cat$-multiequivalences have the properties of being essentially surjective on objects and invertible on multimorphism categories.  \cref{thm:multiwhitehead} in \cref{sec:catmultieqtheorem} shows that these properties precisely characterize $\Cat$-multiequivalences.  Moreover, restricting this proof to the 1-ary multimorphism categories recovers an analogous characterization of 2-equivalences between 2-categories (\cref{thm:iiequivalence}).

Throughout this chapter we assume that 
\[\V = (\V,\otimes,\tu,\alpha,\lambda,\rho,\xi)\]
is a symmetric monoidal category (\cref{def:symmoncat}).  We remind the reader of our left normalized bracketing \cref{expl:leftbracketing}.

\section{Enriched Multicategories}
\label{sec:enrmulticat}

In this section we review the definition of an enriched multicategory and a few basic examples.  More elaborate examples appear in later chapters; see \cref{expl:permindexedcat,thm:dcatcatmulticat,thm:permcatenrmulticat,thm:pfibdmulticat,thm:Gacatsmc}.  We begin with some notation.

\begin{definition}\label{def:profile}
Suppose $D$\label{notation:s-class} is a class.  
\begin{itemize}
\item Denote by\index{profile}\label{notation:profs}
\[\Prof(D) = \coprodover{n \geq 0}\ D^{\times n}\] 
the class of finite tuples in $D$.  An element in $\Prof(D)$ is called a \emph{$D$-profile}.  
\item A $D$-profile of length\index{length!of a profile}
  $n=\len\angd$ is denoted by $\angd = (d_1, \ldots, d_n) \in
  D^{n}$\label{notation:us} or by $\ang{d_i}_{i=1}^n$.  The empty $D$-profile\index{empty profile}
  is denoted by $\ang{}$.
\item We let $\oplus$ denote the \label{not:concat}concatenation\index{concatenation} of profiles.  It is associative with unit given by the empty tuple $\ang{}$.
\item An element in $\Prof(D)\times D$ is denoted
  by\label{notation:duc} $\smscmap{\angd; d'}$ with $d'\in D$ and
  $\angd\in\Prof(D)$.\defmark
\end{itemize}
\end{definition}

The \index{symmetric group}symmetric group on $n$ letters is denoted $\Sigma_n$.

\begin{definition}\label{def:enr-multicategory}
A \emph{$\V$-multicategory}\index{enriched multicategory}\index{category!enriched multi-}\index{multicategory!enriched} $(\M, \gamma, \operadunit)$\label{notation:enr-multicategory} consists of the following data.
\begin{itemize}
\item $\M$ is equipped with a class $\ObM$
  of\index{object!enriched multicategory} \emph{objects}.  We abbreviate $\Prof(\Ob\M)$ to $\Prof(\M)$.
\item For $x'\in\ObM$ and $\angx= \ang{x_j}_{j=1}^n \in\ProfM$, $\M$ is
  equipped with an object\label{notation:enr-cduc}
  \[\M\scmap{\angx;x'} = \M\mmap{x'; x_1,\ldots,x_n} \in \V,\]
  called the \emph{$n$-ary operation object}\index{n-ary@$n$-ary!operation object} or \index{multimorphism}\emph{$n$-ary multimorphism object} with \emph{input profile}\index{input profile} $\angx$ and
  \emph{output}\index{output} $x'$.
\item
For $\smscmap{\angx;x'} \in \ProfMM$ as above and a permutation $\sigma \in
\Sigma_n$, $\M$ is equipped with an isomorphism in $\V$
\begin{equation}\label{rightsigmaaction}
\begin{tikzcd}[column sep=large]
\M\scmap{\angx;x'} \rar{\sigma}[swap]{\cong} & \M\scmap{\angx\sigma; x'},\end{tikzcd}
\end{equation}
called the \emph{right $\sigma$-action}\index{right action} or the \index{symmetric group!action}\emph{symmetric group action}, with\label{enr-notation:c-sigma}
\[\angx\sigma = \left(x_{\sigma(1)}, \ldots, x_{\sigma(n)}\right) \in \ProfM\]
the right permutation\index{right permutation} of $\angx$ by $\sigma$.
\item Each object $x$ in $\M$ is equipped with a morphism\label{notation:enr-unit-c}
\begin{equation}\label{ccoloredunit}
\begin{tikzcd}[column sep=large]
\tensorunit \ar{r}{\operadunit_x} & \M\scmap{x;x},
\end{tikzcd}
\end{equation}
called the \index{colored unit}\emph{$x$-colored unit}.
\item For
\begin{itemize}
\item $x'' \in \ObM$,
\item $\ang{x'} = \ang{x'_j}_{j=1}^n \in \ProfM$, 
\item $\ang{x_j} = \ang{x_{j,i}}_{i=1}^{k_j} \in \ProfM$ for each $j\in\{1,\ldots,n\}$ with $\angx = \oplus_{j=1}^n \ang{x_j}$,
\end{itemize}
$\M$ is equipped with a morphism in $\V$\label{notation:enr-multicategory-composition}
  \begin{equation}\label{eq:enr-defn-gamma}
    \begin{tikzcd}[column sep=large]
      \M\mmap{x'';\ang{x'}} \otimes
      \txotimes_{j=1}^n \M\mmap{x_j';\ang{x_j}}
      \rar{\gamma}
      &
      \M\mmap{x'';\ang{x}},
    \end{tikzcd}
  \end{equation}
called the \index{multicategory!composition}\index{composition!multicategory}\emph{composition} or \emph{multicategorical composition}. 
\end{itemize}
The data above are required to satisfy the following axioms.
\begin{description}
\item[Symmetric Group Action]
For $\mmap{x';\ang{x}}\in\ProfMM$ with $n=\len\ang{x}$ and
$\sigma,\tau\in\Sigma_n$, the following diagram in $\V$ commutes. 
\begin{equation}\label{enr-multicategory-symmetry}
\begin{tikzcd}
\M\scmap{\angx;x'} \arrow{rd}[swap]{\sigma\tau} \rar{\sigma} 
& \M\scmap{\angx\sigma;x'} \dar{\tau}\\
& \M\mmap{x';\angx\sigma\tau}
\end{tikzcd}
\end{equation}
The identity in $\Sigma_n$ acts as the identity morphism on $\M\scmap{\angx;x'}$.
\item[Associativity]
Suppose given
\begin{itemize}
\item $x''' \in \ObM$,
\item $\ang{x''} = \ang{x''_j}_{j=1}^n \in \ProfM$,
\item $\ang{x_j'} = \ang{x'_{j,i}}_{i=1}^{k_{j}} \in \ProfM$ for each
  $j \in \{1,\ldots,n\}$ with $\ang{x'} = \oplus_{j=1}^n \ang{x_j'}$ and $k_j > 0$ for at least one $j$, and
\item $\ang{x_{j,i}} = \ang{x_{j,i,p}}_{p=1}^{\ell_{j,i}} \in \ProfM$ for
  each $j\in\{1,\ldots,n\}$ and each $i \in \{1,\ldots,k_j\}$ with $\ang{x_j} = \oplus_{i=1}^{k_j}{\ang{x_{j,i}}}$ and $\ang{x} =
\oplus_{j=1}^n{\ang{x_{j}}}$.
\end{itemize}
Then the \index{associativity!enriched multicategory}\emph{associativity diagram} below commutes.
\begin{equation}\label{enr-multicategory-associativity}
\begin{tikzpicture}[x=40mm,y=15mm,vcenter]
  \draw[0cell=.85] 
  (0,0) node (a) {\textstyle
    \M\mmap{x''';\ang{x''}}
    \otimes
    \biggl[\bigotimes\limits_{j=1}^n \M\mmap{x''_j;\ang{x'_{j}}}\biggr]
    \otimes
    \bigotimes\limits_{j=1}^n \biggl[\bigotimes\limits_{i=1}^{k_j} \M\mmap{x'_{j,i};\ang{x_{j,i}}}\biggr] 
  }
  (1,.8) node (b) {\textstyle
    \M\mmap{x''';\ang{x'}}
    \otimes
    \bigotimes\limits_{j=1}^{n} \biggl[\bigotimes\limits_{i=1}^{k_j} \M\mmap{x'_{j,i};\ang{x_{j,i}}}\biggr]
  }
  (0,-1.2) node (a') {\textstyle
    \M\mmap{x''';\ang{x''}} \otimes
    \bigotimes\limits_{j=1}^n \biggl[\M\mmap{x_j'';\ang{x_j'}} \otimes \bigotimes\limits_{i=1}^{k_j} \M\mmap{x'_{j,i};\ang{x_{j,i}}}\biggr]
  }
  (1,-2) node (b') {\textstyle
    \M\mmap{x''';\ang{x''}} \otimes \bigotimes\limits_{j=1}^n \M\mmap{x_j'';\ang{x_{j}}}
  }
  (1.2,-.6) node (c) {\textstyle
    \M\mmap{x''';\ang{x}}
  }
  ;
  \draw[1cell=.85]
  (a) edge[shorten <=-1ex,shorten >=-1ex] node {\iso} node['] {\mathrm{permute}} (a')
  (a) edge[shorten >=-4ex,shorten <=-1ex, transform canvas={xshift=-2.5em}] node[pos=.6] {(\ga,1)} (b)
  (b) edge node {\ga} (c)
  (a') edge[shorten >=-3ex, shorten <=-1ex,transform canvas={xshift=-2.5em}] node[swap,pos=.6] {(1,\textstyle\bigotimes_j \ga)} (b')
  (b') edge['] node {\ga} (c)
  ;
\end{tikzpicture}
\end{equation}

\item[Unity]
Suppose $x' \in \ObM$.\index{unity!enriched multicategory}
\begin{enumerate}
\item If $\angx = \ang{x_j}_{j=1}^n \in \ProfM$ has length $n \geq 1$,
  then the following \emph{right unity diagram}\index{right unity!enriched multicategory} commutes, with $\tensorunit^n$ the $n$-fold monoidal product of $\tensorunit$ with itself.
  \begin{equation}\label{enr-multicategory-right-unity}
    \begin{tikzcd} \M\scmap{\angx;x'} \otimes \tensorunit^{n} \dar[swap]{1 \otimes (\otimes_j \operadunit_{x_j})} \rar{\rho} & \M\scmap{\angx;x'} \dar{1}\\
      \M\scmap{\angx;x'} \otimes \txotimes_{j=1}^n \M\scmap{x_j;x_j} \rar{\gamma} & \M\scmap{\angx;x'}
    \end{tikzcd}
  \end{equation}

\item
For any $\ang{x} \in \ProfM$, the \index{unity!enriched multicategory}\index{left unity!enriched multicategory}\emph{left unity diagram} below commutes.
\begin{equation}\label{enr-multicategory-left-unity}
\begin{tikzcd}
\tensorunit \otimes \M\scmap{\angx;x'} \dar[swap]{\operadunit_{x'} \otimes 1} \rar{\lambda} & 
\M\scmap{\angx;x'} \dar{1}\\
\M\mmap{x';x'} \otimes \M\scmap{\angx;x'} \rar{\gamma} & \M\scmap{\angx;x'}
\end{tikzcd}
\end{equation}
\end{enumerate}
\item[Equivariance]
Suppose that in the definition of $\gamma$ \eqref{eq:enr-defn-gamma}, $\len\ang{x_j} = k_j \geq 0$.\index{equivariance!enriched multicategory}
\begin{enumerate}
\item For each $\sigma \in \Sigma_n$, the following \index{top equivariance!enriched multicategory}\emph{top equivariance diagram} commutes.
\begin{equation}\label{enr-operadic-eq-1}
\begin{tikzcd}[column sep=large,cells={nodes={scale=.8}},
every label/.append style={scale=.8}]
\M\mmap{x'';\ang{x'}} \otimes \txotimes_{j=1}^n \M\mmap{x'_j;\ang{x_j}} 
\dar[swap]{\gamma} \rar{(\sigma, \sigma^{-1})}
& \M\mmap{x'';\ang{x'}\sigma} \otimes \txotimes_{j=1}^n \M\mmap{x'_{\sigma(j)};\ang{x_{\sigma(j)}}} \dar{\gamma}\\
\M\mmap{x'';\ang{x_1},\ldots,\ang{x_n}} \rar{\sigma\langle k_{\sigma(1)}, \ldots , k_{\sigma(n)}\rangle}
& \M\mmap{x'';\ang{x_{\sigma(1)}},\ldots,\ang{x_{\sigma(n)}}}
\end{tikzcd}
\end{equation}
Here\label{notation:enr-block-permutation}  
\begin{equation}\label{blockpermutation}
\sigma\langle k_{\sigma(1)}, \ldots , k_{\sigma(n)} \rangle \in \Sigma_{k_1+\cdots+k_n}
\end{equation}
is the block permutation\index{block!permutation} that permutes $n$ consecutive blocks of lengths $k_{\sigma(1)}$, $\ldots$, $k_{\sigma(n)}$ as $\sigma$ permutes $\{1,\ldots,n\}$, leaving the relative order within each block unchanged.
\item
Given permutations $\tau_j \in \Sigma_{k_j}$ for $1 \leq j \leq n$,
the following \index{bottom equivariance!enriched multicategory}\emph{bottom equivariance
  diagram} commutes.
\begin{equation}\label{enr-operadic-eq-2}
\begin{tikzcd}[cells={nodes={scale=.85}},
every label/.append style={scale=.85}]
\M\mmap{x'';\ang{x'}} \otimes \bigotimes_{j=1}^n \M\mmap{x'_j;\ang{x_j}}
\dar[swap]{\gamma} \rar{(1, \otimes_j \tau_j)} & 
\M\mmap{x'';\ang{x'}} \otimes \bigotimes_{j=1}^n \M\mmap{x'_j;\ang{x_j}\tau_j}\dar{\gamma} \\
\M\mmap{x'';\ang{x_1},\ldots,\ang{x_n}} \rar{\tau_1 \times \cdots \times \tau_n}
& \M\mmap{x'';\ang{x_1}\tau_1,\ldots,\ang{x_n}\tau_n}
\end{tikzcd}
\end{equation}
Here the block sum\index{block!sum}\label{notation:enr-block-sum} 
\begin{equation}\label{blocksum}
\tau_1 \times\cdots \times\tau_n \in \Sigma_{k_1+\cdots+k_n}
\end{equation} 
is the image of $(\tau_1, \ldots, \tau_n)$ under the canonical inclusion \[\Sigma_{k_1} \times \cdots \times \Sigma_{k_n} \to \Sigma_{k_1 + \cdots + k_n}.\]
\end{enumerate}
\end{description}
This finishes the definition of a $\V$-multicategory.  

Moreover, we define the following.
\begin{itemize}
\item A $\V$-multicategory is \emph{small}\index{multicategory!enriched!small}\index{small!enriched multicategory} if its class of objects is a set.
\item A \emph{non-symmetric $\V$-multicategory}\index{non-symmetric!multicategory}\index{multicategory!non-symmetric} is defined in the same way as a $\V$-multicategory but without the symmetric group action and the axioms \cref{enr-multicategory-symmetry,enr-operadic-eq-1,enr-operadic-eq-2}.
\item A \index{multicategory}\emph{(non-symmetric) multicategory} is a (non-symmetric) $\Set$-multicategory, where $(\Set,\times,*)$ is the symmetric monoidal category of sets and functions with the Cartesian product.
\item A \index{enriched operad}\index{enriched!operad}\index{operad!enriched}\emph{$\V$-operad} is a $\V$-multicategory with one object.  If $\M$ is a $\V$-operad, then its object of $n$-ary operations is denoted by \label{not:nthobject}$\M_n \in \V$.
\item An \emph{operad}\index{multicategory!one object} is a $\Set$-operad, that is, a multicategory with one object.\defmark
\end{itemize}
\end{definition}

\begin{remark}\label{rk:lambekmulticat}
A non-symmetric multicategory as in \cref{def:enr-multicategory} is the same as a multicategory in \cite{lambek}.
\end{remark}

\begin{example}\label{definition:terminal-operad-comm}  
The \index{terminal!multicategory}\index{multicategory!terminal}\emph{terminal multicategory} $\Mterm$ consists of a single object $*$ and a single $n$-ary operation $\iota_n$ for each $n \ge 0$.
\end{example}

\begin{example}[Endomorphism Operad]\label{example:enr-End}
For each $\V$-multicategory $\M$ and object $x$ of $\M$, the \emph{endomorphism $\V$-operad} \index{endomorphism!operad}\index{operad!endomorphism}$\End(x)$ consists of
  the single object $x$ and $n$-ary operation object
  \[\End(x)_n = \M\mmap{x;\ang{x}} \in \V,\]
  with $\ang{x}$ the constant $n$-tuple at $x$.  The
  symmetric group action, unit, and composition of $\End(x)$ are given
  by those of $\M$.
\end{example}

\begin{example}[Endomorphism Multicategory]\label{ex:endc}
Each permutative category $(\C,\oplus,e,\xi)$ has an associated \index{endomorphism!multicategory}\index{multicategory!endomorphism}\index{permutative category!endomorphism multicategory}\emph{endomorphism multicategory} $\End(\C)$ with object set $\Ob\C$ and $n$-ary multimorphism set
  \[\End(\C)\mmap{y;\ang{x}} 
  = \C(x_1 \oplus \cdots \oplus x_n , y)\]
  for $y \in \Ob\C$ and $\ang{x} = \ang{x_j}_{j=1}^n \in (\Ob\C)^{\times n}$.  An empty $\oplus$ is, by definition, the unit object $e$.
\end{example}

\begin{example}[Underlying $\V$-Category]\label{ex:unarycategory}\index{enriched multicategory!underlying enriched category}\index{enriched category!underlying - of an enriched multicategory}
Each non-symmetric $\V$-multicategory $(\M,\ga,\operadunit)$ has an underlying $\V$-category (\cref{def:enriched-category}) defined by the following data.
\begin{itemize}
\item It has the same class of objects as $\M$.
\item It has hom object $\M\smscmap{a;b}$ with domain $a$ and codomain $b$.
\item Its identities are given by the colored units in $\M$.
\item Its composition is given by 
\[\begin{tikzcd}[column sep=large]
\M\scmap{b;c} \otimes \M\scmap{a;b} \ar{r}{\gamma} & \M\scmap{a;c}
\end{tikzcd}\]
for objects $a, b, c \in \M$.
\end{itemize}   
The associativity diagram \cref{enriched-cat-associativity} and the unity diagram \cref{enriched-cat-unity} of a $\V$-category are the 1-ary special cases of, respectively, the associativity diagram \cref{enr-multicategory-associativity} and the unity diagrams \cref{enr-multicategory-right-unity,enr-multicategory-left-unity} of a $\V$-multicategory.
\end{example}

\section{The 2-Category of Enriched Multicategories}
\label{sec:iicatofmulticat}

As in \cref{sec:enrmulticat} we assume that $\V$ is a symmetric monoidal category.  In this section we  review the 2-category of small $\V$-enriched multicategories, $\V$-enriched multifunctors, and $\V$-enriched multinatural transformations.  In \cref{expl:catmultitransformation} we explicitly describe the data and axioms of a $\Cat$-enriched multinatural transformation.

\subsection*{Enriched Multifunctors and Multinatural Transformations}

\begin{definition}\label{def:enr-multicategory-functor}
A \emph{$\V$-multifunctor}\index{multifunctor!enriched}\index{enriched!multifunctor}\index{functor!multi-} $F \cn \M \to \N$ between $\V$-multicategories $\M$ and $\N$ consists of the following data and axioms.
\begin{itemize}
\item It is equipped with an object assignment $F \cn \ObM \to \ObN$.
\item For each $\smscmap{\angx;x'} \in \ProfMM$ with $\angx= \ang{x_j}_{j=1}^n$, it is equipped with a component morphism in $\V$
\[\begin{tikzcd}[column sep=large]
\M\mmap{x';\ang{x}} \ar{r}{F} & \N\mmap{Fx';F\ang{x}}
\end{tikzcd}\]
with $F\angx= \ang{Fx_j}_{j=1}^n$.
\end{itemize}
The data above are required to satisfy the following three axioms.
\begin{description}
\item[Symmetric Group Action] For each $\smscmap{\angx;x'}$ as above and permutation $\sigma \in \Sigma_n$, the following diagram\index{equivariance!enriched multifunctor} in $\V$ commutes.
\begin{equation}\label{enr-multifunctor-equivariance}
\begin{tikzcd}[column sep=large]
\M\mmap{x';\ang{x}} \ar{d}{\cong}[swap]{\sigma} \ar{r}{F} & \N\mmap{Fx';F\ang{x}} \ar{d}{\cong}[swap]{\sigma}\\
\M\mmap{x';\ang{x}\sigma} \ar{r}{F} & \N\mmap{x';F\ang{x}\sigma}\end{tikzcd}
\end{equation}
\item[Units] For each object $x$ in $\M$, the following diagram in $\V$ commutes.
\begin{equation}\label{enr-multifunctor-unit}
\begin{tikzpicture}[x=25mm,y=15mm,vcenter]
  \draw[0cell] 
  (0,0) node (a) {\tu}
  (1,.5) node (b) {\M\mmap{x;x}}
  (1,-.5) node (b') {\N\mmap{Fx;Fx}}
  ;
  \draw[1cell] 
  (a) edge node {\operadunit_x} (b)
  (a) edge node[swap,pos=.6] {\operadunit_{Fx}} (b')
  (b) edge node {F} (b')
  ;
\end{tikzpicture}
\end{equation} 
\item[Composition] For $x''$, $\ang{x'}$, and $\ang{x} =
  \oplus_{j=1}^n\ang{x_j}$ as in \eqref{eq:enr-defn-gamma}, the following diagram in $\V$ commutes.
\begin{equation}\label{v-multifunctor-composition}
\begin{tikzcd}[column sep=large,cells={nodes={scale=.8}},
every label/.append style={scale=.9}]
\M\mmap{x'';\ang{x'}} \otimes \bigotimes_{j=1}^n \M\mmap{x'_j;\ang{x_j}} \dar[swap]{\gamma} \ar{r}{(F,\otimes_j F)} & \N\mmap{Fx'';F\ang{x'}} \otimes \bigotimes_{j=1}^n \N\mmap{Fx'_j;F\ang{x_j}} \dar{\gamma}\\  
\M\mmap{x'';\ang{x}} \ar{r}{F} & \N\mmap{Fx'';F\ang{x}}
\end{tikzcd}
\end{equation}
\end{description}
This finishes the definition of a $\V$-multifunctor.  

Moreover, we define the following.
\begin{itemize}
\item For another $\V$-multifunctor $G \cn \N\to\P$ between $\V$-multicategories, where $\P$ has object class $\ObP$, the \index{composition!enriched multifunctor}\emph{composition} $GF \cn \M\to\P$ is the $\V$-multifunctor defined by composing the object assignments 
\[\begin{tikzcd} \ObM \ar{r}{F} & \ObN \ar{r}{G} & \ObP
\end{tikzcd}\]
and the morphisms on $n$-ary multimorphism objects
\[\begin{tikzcd}
\M\mmap{x';\ang{x}} \ar{r}{F} & \N\mmap{Fx';F\ang{x}} \ar{r}{G} & \P\mmap{GFx';GF\ang{x}}.
\end{tikzcd}\]
\item The \index{identity!enriched multifunctor}\emph{identity $\V$-multifunctor} $1_{\M} \cn \M\to\M$ is defined by the identity assignment on objects and the identity morphisms on $n$-ary multimorphism objects.
\item A \emph{$\V$-operad morphism}\index{operad!morphism}\index{enriched operad!morphism}\index{morphism!enriched operad} is a $\V$-multifunctor between two $\V$-multicategories with one object.
\item A \emph{non-symmetric $\V$-multifunctor}\index{non-symmetric!multifunctor}\index{multifunctor!non-symmetric} $F \cn \M \to \N$ between non-symmetric $\V$-multicategories $\M$ and $\N$ is defined in the same way as a $\V$-multifunctor but without the symmetric group action axiom \cref{enr-multifunctor-equivariance}.  Composition and identities are defined as above.
\item A \emph{(non-symmetric) multifunctor} is a (non-symmetric) $\Set$-multifunctor.\defmark
\end{itemize}
\end{definition}

\begin{definition}\label{def:enr-multicat-natural-transformation}
For $\V$-multifunctors $F,G \cn \M\to\N$, a \emph{$\V$-multinatural transformation}\index{enriched!multinatural transformation}\index{natural transformation!enriched multi-}\index{multinatural transformation!enriched} $\theta \cn F\to G$ consists of, for each object $x$ in $\M$, a component morphism in $\V$
\[\begin{tikzcd}[column sep=large]
\tu \ar{r}{\theta_x} & \N\mmap{Gx;Fx}
\end{tikzcd}\]
such that the following \emph{$\V$-naturality diagram} in $\V$
commutes for each $\smscmap{\ang{x};x'} \in \ProfMM$ with
$\angx= \ang{x_j}_{j=1}^n$.
\begin{equation}\label{enr-multinat}
\begin{tikzpicture}[x=25mm,y=12mm,vcenter]
  \draw[0cell=.8]
  (.5,0) node (a) {\M\mmap{x';\ang{x}}}
  (1,1) node (b) {\tu \otimes \M\mmap{x';\ang{x}}}
  (3,1) node (c) {\N\mmap{Gx';Fx'} \otimes \N\mmap{Fx';F\ang{x}}}
  (3.5,0) node (d) {\N\mmap{Gx';F\ang{x}}}
  (1,-1) node (b') {\M\mmap{x';\ang{x}}\otimes \txotimes_{j=1}^n \tu}
  (3,-1) node (c') {\N\mmap{Gx';G\ang{x}} \otimes \txotimes_{j=1}^n \N\mmap{Gx_j;Fx_j}}
  ;
  \draw[1cell=.8] 
  (a) edge node[pos=.3] {\la^\inv} (b)
  (a) edge node[swap,pos=.2] {\rho^\inv} (b')
  (b) edge node {\theta_{x'} \otimes F} (c)
  (c) edge node[pos=.7] {\ga} (d)
  (b') edge node {G \otimes \txotimes_{j=1}^n \theta_{x_j}} (c')
  (c') edge['] node[pos=.7] {\ga} (d)
  ;
\end{tikzpicture}
\end{equation}
This finishes the definition of a $\V$-multinatural transformation.  

Moreover, we define the following.
\begin{itemize}
\item The \emph{identity $\V$-multinatural transformation}\index{multinatural transformation!identity} $1_F \cn F\to F$ has components \[(1_F)_x = \operadunit_{Fx} \forspace x\in\ObM.\]
\item A \emph{multinatural transformation} is a $\Set$-multinatural transformation.
\item A \emph{non-symmetric $\V$-multinatural transformation}\index{multinatural transformation!non-symmetric}\index{non-symmetric!multinatural transformation} $\theta \cn F \to G$ between non-symmetric $\V$-multifunctors $F, G \cn \M \to \N$ is defined as above.\defmark
\end{itemize} 
\end{definition}

\begin{definition}\label{def:enr-multinatural-composition}
Suppose $\theta \cn F \to G$ is a $\V$-multinatural transformation between $\V$-multifunctors as in \cref{def:enr-multicat-natural-transformation}.
\begin{itemize}
\item Suppose $\psi \cn G \to H$ is a $\V$-multinatural transformation with $H \cn \M \to \N$ a $\V$-multifunctor.  The \emph{vertical composition}\index{vertical composition!enriched multinatural transformation}\label{notation:enr-operad-vcomp}
\begin{equation}\label{multinatvcomp}
\begin{tikzcd}[column sep=large]
F \ar{r}{\psi\theta} & H
\end{tikzcd}  
\end{equation} 
is the $\V$-multinatural transformation with component at each $x \in \ObM$ given by the following composite in $\V$.
  \[
  \begin{tikzpicture}[x=45mm,y=12mm]
    \draw[0cell=.9] 
    (0,0) node (a) {\tu}
    (0,-1) node (b) {\tu \otimes \tu}
    (1,-1) node (c) {\N\mmap{Hx;Gx} \otimes \N\mmap{Gx;Fx}}
    (1,0) node (d) {\N\mmap{Hx;Fx}}
    ;
    \draw[1cell=.9] 
    (a) edge node['] {\la^\inv} (b)
    (b) edge node {\psi_x \otimes \theta_x} (c)
    (c) edge['] node {\ga} (d)
    (a) edge node {(\psi\theta)_x} (d)
    ;
  \end{tikzpicture}
  \]
\item Suppose $F', G' \cn \N \to \P$ are $\V$-multifunctors, and $\theta' \cn F' \to G'$ is a $\V$-multinatural transformation.  The \emph{horizontal composition}\index{horizontal composition!enriched multinatural transformation}\label{notation:enr-operad-hcomp}
\begin{equation}\label{multinathcomp}
\begin{tikzcd}[column sep=large]
F'F \ar{r}{\theta' \ast \theta} & G'G
\end{tikzcd}
\end{equation} 
is the $\V$-multinatural transformation with component at each $x \in \ObM$ given by the following composite in $\V$.
\[
\begin{tikzpicture}[x=50mm,y=12mm]
  \draw[0cell=.9] 
  (0,0) node (a) {\tu}
  (0,-2) node (b) {\tu \otimes \tu}
  (1,-2) node (c) {\P\mmap{G'Gx;F'Gx} \otimes \N\mmap{Gx;Fx}}
  (1,-1) node (d) {\P\mmap{G'Gx;F'Gx} \otimes \P\mmap{F'Gx;F'Fx}}
  (1,0) node (e) {\P\mmap{G'Gx;F'Fx}}
  ;
  \draw[1cell=.9] 
  (a) edge node {(\theta' * \theta)_x} (e)
  (a) edge['] node {\la^\inv} (b)
  (b) edge node {\theta'_{Gx} \otimes \theta_x} (c)
  (c) edge['] node {1 \otimes F'} (d)
  (d) edge['] node {\ga} (e)
  ;
\end{tikzpicture}
\]
\end{itemize}
Vertical and horizontal compositions of non-symmetric $\V$-multinatural transformations are defined as above.
This finishes the definition.
\end{definition}

\subsection*{The 2-Category Structure}

\begin{theorem}\label{v-multicat-2cat}
For each symmetric monoidal category $\V$, there is a 2-category\index{2-category!of small enriched multicategories}\index{enriched multicategory!2-category}\index{multicategory!enriched!2-category} 
\[\VMulticat\]
consisting of the following data.
\begin{itemize}
\item Its objects are small $\V$-multicategories.
\item For small $\V$-multicategories $\M$ and $\N$, the hom category 
\[\VMulticat(\M,\N)\]
is defined as follows.
\begin{itemize}
\item Its objects are $\V$-multifunctors $\M \to \N$.
\item Its morphisms are $\V$-multinatural transformations.
\item Composition is vertical composition of $\V$-multinatural transformations.
\item Identity morphisms are identity $\V$-multinatural transformations.
\end{itemize}
\item The identity 1-cell $1_{\M}$ is the identity $\V$-multifunctor $1_{\M}$.
\item Horizontal composition of 1-cells is the composition of $\V$-multifunctors.
\item Horizontal composition of 2-cells is that of $\V$-multinatural transformations.
\end{itemize}
Moreover, there is a 2-category \label{not:vmulticatns}$\VMulticatns$ consisting of the following data.
\begin{itemize}
\item Its objects are small non-symmetric $\V$-multicategories.
\item Its 1-cells are non-symmetric $\V$-multifunctors.
\item Its 2-cells are non-symmetric $\V$-multinatural transformations.
\item The rest of the 2-category structure is as in $\VMulticat$.
\end{itemize} 
\end{theorem}

\begin{proof}
The special case $\V = \Set$ is \cite[2.4.26]{johnson-yau}.  That proof can be adapted to prove the $\V$-enriched version essentially without change.  The proof for the non-symmetric variant is obtained by omitting the parts involving the symmetric group action.
\end{proof}

Recall from \cref{def:twocategory} the \emph{underlying 1-category} of a 2-category.  Next is the $\V$-enriched analog of \cite[5.5.14]{cerberusIII}, whose proof naturally extends to the $\V$-enriched setting.  

\begin{theorem}\label{vmulticatbicomplete}
Suppose $\V$ is a complete and cocomplete symmetric monoidal closed category.  Then the underlying 1-categories of the 2-categories
\[\VMulticat \andspace \VMulticatns\]
are complete and cocomplete.
\end{theorem}

\begin{example}[Initial and Terminal Objects]\label{ex:vmulticatinitialterminal}\index{initial!enriched multicategory}\index{enriched multicategory!initial}\index{terminal!enriched multicategory}\index{enriched multicategory!terminal}
The initial $\V$-multicategory has an empty set of objects.  If $\V$ has a terminal object \label{not:bt}$\bt$, then a terminal $\V$-multicategory $\Mterm$ has a single object $*$ and $n$-ary multimorphism object
\[\Mterm_n = \Mterm(\overbracket[.5pt]{*, \ldots, *}^{\text{$n$ terms}} ; *) = \bt\]
for each $n \geq 0$.
\end{example}

\subsection*{Categorically-Enriched Multicategories}

Recall the symmetric monoidal closed category $(\Cat,\times,\boldone)$ of small categories and functors with the Cartesian product (\cref{ex:cat}).

\begin{definition}\label{def:catmulticat}
Suppose $\angx = \ang{x_j}_{j=1}^n, y, z$ are objects in a $\Cat$-multicategory $(\M,\ga,1)$.
\begin{itemize}
\item The category $\M\scmap{\angx;y}$ is called a \index{category!multimorphism}\index{multimorphism category}\emph{multimorphism category}.
\item An object in the multimorphism category $\M\scmap{\angx;y}$ is called an \index{1-cell!n-ary@$n$-ary}\index{n-ary@$n$-ary!1-cell}\emph{$n$-ary 1-cell} and is denoted \label{not:naryonecell}$\angx \to y$.
\item A 1-ary 1-cell $f \cn y \to z$ is an \emph{isomorphism} if there exists a 1-ary 1-cell $g \cn z \to y$ such that
\[\ga(g;f) = 1_y \andspace \ga(f;g) = 1_z.\]
In this case, we say that $f$ is \emph{invertible} with \emph{inverse} $g$ and denote it by $y \fto{\iso} z$.
\item A $\Cat$-multifunctor $F \cn \M \to \N$ is \emph{essentially surjective}\index{essentially surjective} on objects if, for each object $w \in \N$, there exist
\begin{itemize}
\item an object $v \in \M$ and
\item an isomorphism $h \cn Fv \fto{\iso} w$ in $\N$.
\end{itemize} 
\item A morphism $\theta \cn f \to g$ in the multimorphism category $\M\scmap{\angx;y}$ is called an \index{2-cell!n-ary@$n$-ary}\index{n-ary@$n$-ary!2-cell}\emph{$n$-ary 2-cell}.  It is also denoted by
\begin{equation}\label{twocellmulticat}
\begin{tikzpicture}[xscale=2,yscale=1.7,baseline={(x1.base)}]
\draw[0cell=.9]
(0,0) node (x1) {\angx}
(x1)++(1,0) node (x2) {y}
;
\draw[1cell=.9]  
(x1) edge[bend left] node {f} (x2)
(x1) edge[bend right] node[swap] {g} (x2)
;
\draw[2cell]
node[between=x1 and x2 at .47, rotate=-90, 2label={above,\theta}] {\Rightarrow}
;
\end{tikzpicture}
\end{equation}
which extends the 2-cell notation in \cref{twocellnotation}.
\end{itemize}
The same terminology and notation apply to non-symmetric $\Cat$-multicategories and $\Cat$-multifunctors.
\end{definition}

\begin{explanation}[$\Cat$-Multinatural Transformations]\label{expl:catmultitransformation}
Suppose $\M$ and $\N$ are $\Cat$-multicategories, and 
\[F,G \cn \M \to \N\]
are $\Cat$-multifunctors.  A \index{multinatural transformation!Cat-@$\Cat$-}\index{Cat-multinatural@$\Cat$-multinatural!transformation}\emph{$\Cat$-multinatural transformation} $\theta \cn F \to G$ consists of a component 1-ary 1-cell
\begin{equation}\label{thetaccomponent}
\begin{tikzcd}[column sep=large]
Fc \ar{r}{\theta_c} & Gc
\end{tikzcd}
\qquad \text{in $\N$ for each $c \in \Ob\M$}
\end{equation}
such that the following two \index{Cat-naturality conditions@$\Cat$-naturality conditions}\emph{$\Cat$-naturality conditions} hold:
\begin{description}
\item[Objects] For each $k$-ary 1-cell $p \cn \angc \to c'$ in $\M$ with
$\angc= \ang{c_j}_{j=1}^k$, denote by
\begin{equation}\label{Fangcthetaangc}
\left\{\begin{split}
F\ang{c} &= \ang{F c_j}_{j=1}^k \in (\Ob\N)^k\\
\theta_{\ang{c}} &= \ang{\theta_{c_j}}_{j=1}^k \in \txprod_{j=1}^k \N\mmap{G c_j; F c_j}.
\end{split}\right.
\end{equation}
Then the following $k$-ary 1-cell equality holds, with composition taken in the $\Cat$-multicategory $\N$:
\begin{equation}\label{catmultinaturality}
\gamma\scmap{Gp; \theta_{\angc}} =
\gamma\scmap{\theta_{c'}; Fp} \inspace \N\mmap{Gc';F\angc}.
\end{equation}
\item[Morphisms] For each $k$-ary 2-cell $f \cn p \to q$ in $\M\mmap{c';\ang{c}}$, the following $k$-ary 2-cell equality holds, with $1_{\theta_{\angc}} = \ang{1_{\theta_{c_j}}}_{j=1}^k$:
\begin{equation}\label{catmultinaturalityiicell}
\gamma\scmap{Gf; 1_{\theta_{\angc}}} =
\gamma\scmap{1_{\theta_{c'}}; Ff} \inspace \N\mmap{Gc';F\angc}.
\end{equation}
\end{description}
The conditions \cref{catmultinaturality,catmultinaturalityiicell} together constitute the $\V$-naturality condition \cref{enr-multinat} when $\V$ is $\Cat$.

The two sides of the object $\Cat$-naturality condition \cref{catmultinaturality} use the following two compositions in $\N$, applied to objects.
\begin{equation}\label{catmultinaturalitycomp}
\begin{tikzpicture}[xscale=2.5,yscale=1.2,vcenter]
\draw[0cell=.85]
(0,0) node (a) {\N\mmap{Gc';G\angc} \times \txprod_{j=1}^k \N\mmap{G c_j;F c_j}}
(a)++(2,0) node (b) {\N\mmap{Gc';Fc'} \times \N\mmap{Fc'; F\angc}}
(a)++(1,-1) node (c) {\N\mmap{Gc';F\angc}}
;
\draw[1cell=.9]  
(a) edge node[swap,pos=.5] {\gamma} (c)
(b) edge node[pos=.5] {\gamma} (c)
;
\end{tikzpicture}
\end{equation}
The object $\Cat$-naturality condition \cref{catmultinaturality} is the commutative diagram
\begin{equation}\label{catmultinaturalitydiagram}
\begin{tikzcd}[column sep=large]
F\angc \ar{d}[swap]{Fp} \ar{r}{\theta_{\angc}} & G \angc \ar{d}{Gp}\\
Fc' \ar{r}{\theta_{c'}} & Gc'
\end{tikzcd}
\end{equation}
with the understanding that it is a $k$-ary 1-cell equality involving the compositions \cref{catmultinaturalitycomp} in $\N$.

Similarly, the  two sides of the morphism $\Cat$-naturality condition \cref{catmultinaturalityiicell} use the compositions in \cref{catmultinaturalitycomp}, applied to morphisms.  The condition \cref{catmultinaturalityiicell} is the equality of \index{pasting diagram!multicategorical}multicategorical pasting diagrams
\begin{equation}\label{catmultinatiicellpasting}
\begin{tikzpicture}[xscale=2.5,yscale=1.5,vcenter]
\def\a{35} \def\s{.8} \def\h{.8}
\def\boundary{
\draw[0cell=\s]
(0,0) node (x11) {F\ang{c}}
(x11)++(\h,0) node (x12) {G\ang{c}}
(x11)++(0,-1) node (x21) {Fc'}
(x12)++(0,-1) node (x22) {Gc'}
;
\draw[1cell=\s]  
(x11) edge node {\theta_{\ang{c}}} (x12)
(x12) edge[bend left=\a] node[pos=.75] {Gq} (x22)
(x11) edge[bend right=\a] node[swap,pos=.25] {Fp} (x21)
(x21) edge node {\theta_{c'}} (x22)
;}
\draw (0,0) node[font=\LARGE] (equal) {=};
\begin{scope}[shift={(.6,.5)}]
\boundary
\draw[1cell=\s] 
(x11) edge[bend left=\a] node[pos=.4] {Fq} (x21)
;
\draw[2cell=\s]
node[between=x11 and x21 at .55, rotate=0, 2label={above,Ff}] {\Longrightarrow}
;
\end{scope}
\begin{scope}[shift={(-1.4,.5)}]
\boundary
\draw[1cell=\s] 
(x12) edge[bend right=\a] node[swap,pos=.4] {Gp} (x22)
;
\draw[2cell=\s]
node[between=x12 and x22 at .55, rotate=0, 2label={above,Gf}] {\Longrightarrow}
;
\end{scope}
\end{tikzpicture}
\end{equation}
with the understanding that it is a $k$-ary 2-cell equality involving the composition in $\N$.

By the object $\Cat$-naturality condition \cref{catmultinaturality} for $\theta' \cn F' \to G'$, there are two equal ways to write the $c$-component of the horizontal composite for each object $c$ in $\M$:
\begin{equation}\label{hcompexpressions}
\begin{split}
(\theta' * \theta)_c &= \ga\scmap{\theta'_{Gc}; F' \theta_c}\\
&= \ga\scmap{G' \theta_c; \theta'_{Fc}}.
\end{split}
\end{equation}

A $\Cat$-multinatural transformation $\theta \cn F \to G$ is called a \emph{$\Cat$-multinatural isomorphism}\index{multinatural isomorphism!Cat-@$\Cat$-}\index{Cat-multinatural@$\Cat$-multinatural!isomorphism} if $\theta$ is invertible with respect to the vertical composition.  This happens if and only if each component 1-ary 1-cell $\theta_c$ \cref{thetaccomponent} is invertible with respect to the composition in $\N$.  In other words, $\theta$ is a $\Cat$-multinatural isomorphism if and only if each component 1-ary 1-cell $\theta_c$ is an isomorphism in the underlying 1-ary 1-category of $\N$ (\cref{def:twocategory,locallysmalltwocat,ex:unarycategory}). 

A non-symmetric $\Cat$-multinatural transformation or isomorphism admits the same description as above.
\end{explanation}

\section{Categorically-Enriched Multiequivalences}
\label{sec:multequiv}

In this section we show that equivalences in the 2-categories $\catmulticat$ and $\catmulticatns$ are
\begin{enumerate}
\item essentially surjective on objects and
\item isomorphisms between multimorphism categories.
\end{enumerate}  
See \cref{catmultieqesff}.  In \cref{sec:catmultieqtheorem} we show that equivalences are actually characterized by these two properties.  

Recall from \cref{v-multicat-2cat} the 2-category \label{not:catmulticat}$\catmulticat$ with
\begin{itemize}
\item small $\Cat$-multicategories as objects,
\item $\Cat$-multifunctors as 1-cells, and
\item $\Cat$-multinatural transformations as 2-cells.
\end{itemize}
It has a non-symmetric variant \label{not:catmulticatns}$\catmulticatns$.  The following notions are restrictions of those in \cref{def:equivalences} to these 2-categories.

\begin{definition}\label{def:catmultiequivalence}
We define the following.
\begin{itemize}
\item A \emph{$\Cat$-multiequivalence}\index{Cat-multiequivalence@$\Cat$-multiequivalence}\index{multiequivalence!Cat-@$\Cat$-} is an equivalence in the 2-category $\catmulticat$.
\item An \emph{adjoint $\Cat$-multiequivalence}\index{adjoint Cat-multiequivalence@adjoint $\Cat$-multiequivalence} is an adjoint equivalence in the 2-category $\catmulticat$.
\end{itemize}
Moreover, their non-symmetric variants are defined in the same way in the 2-category $\catmulticatns$.
\end{definition}

\begin{explanation}\label{expl:catmultiequivalence}
Adjoint $\Cat$-multiequivalences can be defined directly without mentioning the 2-category $\catmulticat$.  In more detail, a $\Cat$-multifunctor $F \cn \M \to \N$ is a \emph{$\Cat$-multiequivalence} if there exist
\begin{itemize}
\item a $\Cat$-multifunctor $G \cn \N \to \M$ and
\item $\Cat$-multinatural isomorphisms (\cref{expl:catmultitransformation})
\begin{equation}\label{catmultiequnit}
\begin{tikzcd}[column sep=large]
1_{\M} \ar{r}{\eta}[swap]{\iso} & GF
\end{tikzcd} \andspace 
\begin{tikzcd}[column sep=large]
FG \ar{r}{\epz}[swap]{\iso} & 1_{\N}.
\end{tikzcd}
\end{equation}
\end{itemize}
The quadruple $(F,G,\eta,\epz)$ is, furthermore, an \emph{adjoint $\Cat$-multiequivalence} if it satisfies the two triangle identities \cref{triangleidentities}, which involve vertical and horizontal compositions of $\Cat$-multinatural transformations.  The same remarks apply to the non-symmetric variants.  In other words, in \cref{def:catmultiequivalence} it is \emph{not} necessary to assume that the $\Cat$-multicategories $\M$ and $\N$ are small.
\end{explanation}
 
Recall from \cref{def:catmulticat} that $n$-ary 1-cells and $n$-ary 2-cells are, respectively, the objects and morphisms in the $n$-ary multimorphism categories.

\begin{lemma}\label{catmultieqesff}\index{Cat-multiequivalence@$\Cat$-multiequivalence!properties}
Suppose
\[\begin{tikzcd}[column sep=large]
\N \ar{r}{G} & \M
\end{tikzcd}\]
is a $\Cat$-multiequivalence.  Then 
$G$ is
\begin{itemize}
\item essentially surjective on objects and
\item an isomorphism on each multimorphism category.
\end{itemize}
Moreover, the non-symmetric variant is also true.
\end{lemma}

\begin{proof}
Since $G$ is a $\Cat$-multiequivalence, by definition there exist
\begin{itemize}
\item a $\Cat$-multifunctor $F \cn \M \to \N$ and
\item $\Cat$-multinatural isomorphisms $\eta \cn 1_{\M} \fto{\iso} GF$ and $\epz \cn FG \fto{\iso} 1_{\N}$ as in \cref{catmultiequnit}.
\end{itemize} 
We consider essential surjectivity first and then the invertibility of the multimorphism functors of $G$.

\newcounter{catmultieqesff-step}
\newcommand{\catmultieqesffstep}[1]{\medskip\textbf{Step \stepcounter{catmultieqesff-step}\arabic{catmultieqesff-step}}: #1\medskip}

\catmultieqesffstep{Essential Surjectivity on Objects}

For each object $x \in \M$, there is an invertible component 1-ary 1-cell
\[\begin{tikzcd}[column sep=large]
x \ar{r}{\eta_x}[swap]{\iso} & G(Fx)
\end{tikzcd} \inspace \M\smscmap{x;GFx}.\]
This shows that $G$ is essentially surjective on the underlying 1-ary 1-categories.

\catmultieqesffstep{Injectivity on $n$-Ary 1-Cells}

Suppose given an $n$-tuple of objects $\angx = \ang{x_j}_{j=1}^n$ with $n \geq 0$ and an object $y$ in $\N$.  We want to show that the multimorphism functor
\begin{equation}\label{Gmultimorphismfunctor}
\begin{tikzcd}[column sep=large]
\N\scmap{\angx;y} \ar{r}{G} & \M\scmap{G\angx;Gy}
\end{tikzcd}
\end{equation}
is injective on objects.  For each $n$-ary 1-cell $f \cn \angx \to y$ in $\N$, the following computation shows that $f$ is uniquely determined by $Gf$, with $\ga$ denoting the composition \cref{eq:enr-defn-gamma} in $\N$.
\begin{equation}\label{fdeterminedbyGf}
\left\{\begin{aligned}
f &= \ga\scmap{f; \ang{1_{x_j}}} & \phantom{M} & \text{by right unity \cref{enr-multicategory-right-unity}}\\
&= \ga\scmap{f; \bang{\ga\smscmap{\epz_{x_j}; \epz_{x_j}^\inv}}} & \phantom{M} & \text{by invertibility of $\epz$}\\
&= \ga\scmap{\ga\smscmap{f; \epz_{\angx}}; \epz^\inv_{\angx}} & \phantom{M} & \text{by associativity \cref{enr-multicategory-associativity}}\\
&= \ga\scmap{\ga\smscmap{\epz_y;FGf}; \epz^\inv_{\angx}} & \phantom{M} & \text{by object $\Cat$-naturality \cref{catmultinaturality}}
\end{aligned}\right.
\end{equation}
This shows that the multimorphism functor $G$ in \cref{Gmultimorphismfunctor} is injective on objects.

The computation \cref{fdeterminedbyGf} with $\epz$ replaced by 
\[\eta^\inv \cn GF \fto{\iso} 1_{\M}\] 
also shows that the $\Cat$-multiequivalence $F \cn \M \to \N$ is injective on $n$-ary 1-cells.

\catmultieqesffstep{Injectivity on $n$-Ary 2-Cells}

Suppose $\theta \cn f \to f'$ is an $n$-ary 2-cell in $\N\scmap{\angx;y}$.  The computation \cref{fdeterminedbyGf}---with $f$, $\epz$, and \cref{catmultinaturality} replaced by, respectively, $\theta$, $1_{\epz}$, and morphism $\Cat$-naturality \cref{catmultinaturalityiicell}---shows that $\theta$ is uniquely determined by $G\theta$.  This shows that the multimorphism functor $G$ in \cref{Gmultimorphismfunctor} is injective on morphisms.  With $1_{\epz}$ replaced by $1_{\eta^\inv}$, the same computation also shows that the $\Cat$-multiequivalence $F \cn \M \to \N$ is injective on $n$-ary 2-cells.

\catmultieqesffstep{Surjectivity on $n$-Ary 1-Cells}

In the context of \cref{Gmultimorphismfunctor}, suppose given an $n$-ary 1-cell 
\[h \cn G\angx \to Gy \inspace \M.\]  
We must show that there exists an $n$-ary 1-cell $f \cn \angx \to y$ in $\N$ such that $Gf = h$.  We define the $n$-ary 1-cell
\[f = \ga\scmap{\ga\smscmap{\epz_y; Fh}; \epz^\inv_{\angx}} \inspace \N\scmap{\angx;y},\] 
which we visualize as the long composite below.
\[\begin{tikzpicture}[xscale=3, yscale=1.2,vcenter]
\draw[0cell=.9]
(0,0) node (x11) {FG\angx}
(x11)++(1,0) node (x12) {\angx}
(x11)++(0,-1) node (x21) {FGy}
(x12)++(0,-1) node (x22) {y}
;
\draw[1cell=.9]
(x12) edge node {\iso} node[swap] {\epz^\inv_{\angx}} (x11)
(x11) edge node[swap] {Fh} (x21)
(x21) edge node {\epz_y} node[swap] {\iso} (x22)
(x12) edge node {f} (x22)
;
\end{tikzpicture}\]
Similar to \cref{fdeterminedbyGf}, we compute as follows.
\begin{equation}\label{FhFGf}
\left\{\begin{aligned}
Fh 
&= \ga\scmap{\epz^\inv_y; \ga\smscmap{f;\epz_{\angx}}} & \phantom{=} & \text{by \cref{enr-multicategory-associativity,enr-multicategory-right-unity,enr-multicategory-left-unity}}\\
&= \ga\scmap{\epz^\inv_y; \ga\smscmap{\epz_y; FGf}} & \phantom{=} & \text{by \cref{catmultinaturality}}\\
&= FGf & \phantom{=} & \text{by \cref{enr-multicategory-associativity,enr-multicategory-left-unity}}
\end{aligned}\right.
\end{equation}
The injective of $F$ on $n$-ary 1-cells established above and \cref{FhFGf} imply the desired equality, $Gf = h$.

\catmultieqesffstep{Surjectivity on $n$-Ary 2-Cells}

In the context of \cref{Gmultimorphismfunctor}, suppose given an $n$-ary 2-cell 
\[\begin{tikzcd}[column sep=large]
h \ar{r}{\varphi} & h' 
\end{tikzcd}
\inspace \M\scmap{G\angx;Gy}.\]
We must show that there exists an $n$-ary 2-cell $\theta \cn f \to f'$ in $\N\smscmap{\angx;y}$ such that $G\theta = \varphi$.  We define the $n$-ary 2-cell
\[\theta = \ga\scmap{\ga\smscmap{1_{\epz_y}; F\varphi}; 1_{\epz^\inv_{\angx}}} \inspace \N\scmap{\angx;y},\] 
which we visualize as the long composite below.
\[\begin{tikzpicture}[xscale=3, yscale=1.3,vcenter]
\def\a{15} \def\b{20}
\draw[0cell=.9]
(0,0) node (x11) {FG\angx}
(x11)++(1,0) node (x12) {\angx}
(x11)++(0,-1) node (x21) {FGy}
(x12)++(0,-1) node (x22) {y}
;
\draw[1cell=.9]
(x12) edge node {\iso} node[swap] {\epz^\inv_{\angx}} (x11)
(x21) edge node {\epz_y} node[swap] {\iso} (x22)
(x11) edge[bend right=\a, transform canvas={xshift=-.5em}] node[swap] {Fh} (x21)
(x11) edge[bend left=\a, transform canvas={xshift=.6em}] node {Fh'} (x21)
(x12) edge[bend right=\b] (x22)
(x12) edge[bend left=\b] (x22)
;
\draw[2cell=.9]
node[between=x11 and x21 at .6, 2label={above,F\varphi}] {\Rightarrow}
node[between=x12 and x22 at .6, 2label={above,\theta}] {\Rightarrow}
;
\end{tikzpicture}\]
The computation \cref{FhFGf}---with $h$, $f$, $\epz$, and \cref{catmultinaturality} replaced by, respectively, $\varphi$, $\theta$, $1_{\epz}$, and morphism $\Cat$-naturality \cref{catmultinaturalityiicell}---shows that
\[F\varphi = FG\theta.\]
The injective of $F$ on $n$-ary 2-cells established above now implies the desired equality, $G\theta = \varphi$.  This finishes the proof that each multimorphism functor \cref{Gmultimorphismfunctor} of $G$ is an isomorphism.

For the non-symmetric variant, we reuse the proof above by noting that the symmetric group action is not used.
\end{proof}

\section{Characterization of $\Cat$-Multiequivalences}
\label{sec:catmultieqtheorem}

In \cref{catmultieqesff} we observe that (non-symmetric) $\Cat$-multiequivalences are essentially surjective on objects and isomorphisms on multimorphism categories.  In this section we prove that these properties precisely characterize (non-symmetric) $\Cat$-multiequivalences; see \cref{thm:multiwhitehead}.  In \cref{thm:dcatpfibdeq} we apply \cref{thm:multiwhitehead} to show that the Grothendieck construction on bipermutative-indexed categories is a non-symmetric $\Cat$-multiequivalence.

Here is an outline of this section.
\begin{enumerate}
\item For a $\Cat$-multifunctor $G$ that satisfies the conclusion of \cref{catmultieqesff}, we construct the left adjoint $F$ and the unit $\eta \cn 1 \to GF$ in \cref{def:Fandeta}.
\begin{itemize}
\item \cref{Fcomponentfunctorial,Fcatmultifunctor} show that $F$ is a $\Cat$-multifunctor, which is non-symmetric if $G$ is.
\item \cref{etacatmultinatural} shows that $\eta$ is a $\Cat$-multinatural isomorphism, which is non-symmetric if $G$ is.
\end{itemize}
\item The counit $\epz \cn FG \to 1$ is constructed in \cref{def:counitepz}.
\begin{itemize}
\item \cref{epzcatmultiiso} shows that $\epz$ is a $\Cat$-multinatural isomorphism, which is non-symmetric if $G$ is.
\item \cref{multieqtriangleids} verifies the triangle identities for the quadruple $(F,G,\eta,\epz)$, providing the last part of the proof of \cref{thm:multiwhitehead}.
\end{itemize} 
\item \cref{thm:iiequivalence} is the 1-ary special case that characterizes 2-equivalences between 2-categories.
\end{enumerate}

\subsection*{Left Adjoint and Unit}

\begin{definition}\label{def:Fandeta}
Suppose given a $\Cat$-multifunctor between $\Cat$-multicategories
\[\begin{tikzcd}[column sep=large]
\N \ar{r}{G} & \M
\end{tikzcd}\]
that is
\begin{itemize}
\item essentially surjective on objects and
\item an isomorphism on each multimorphism category.
\end{itemize}
Define the data of
\begin{itemize}
\item a $\Cat$-multifunctor $F \cn \M \to \N$ and
\item a $\Cat$-multinatural isomorphism $\eta \cn 1_{\M} \fto{\iso} GF$
\end{itemize}
as follows.
\begin{description}
\item[$F$ on objects and components of $\eta$]
Given an object $x \in \M$, by the essential surjectivity of $G$ and the Axiom of Choice, we choose an object $Fx \in \N$ and an isomorphism
\begin{equation}\label{etaxxgfx}
\begin{tikzcd}[column sep=large]
x \ar{r}{\eta_x}[swap]{\iso} & G(Fx)
\end{tikzcd} \inspace \M\scmap{x;GFx}.
\end{equation}
These choices define the object assignment, $x \mapsto Fx$, of $F$ and the component 1-ary 1-cells of $\eta$.
\item[$F$ on $n$-ary 1-cells]
Given an $n$-ary 1-cell $f \cn \angx \to y$ in $\M$ with $\angx = \ang{x_j}_{j=1}^n$, consider the $n$-ary 1-cell
\[\ga\lrscmap{\eta_y; \ga\lrscmap{f; \eta^\inv_{\angx}}} \inspace \M\scmap{GF\angx;GFy},\]
which we visualize as the long composite below.
\[\begin{tikzpicture}[xscale=3, yscale=1.2,vcenter]
\draw[0cell=.9]
(0,0) node (x11) {\angx}
(x11)++(1,0) node (x12) {GF\angx}
(x11)++(0,-1) node (x21) {y}
(x12)++(0,-1) node (x22) {GFy}
;
\draw[1cell=.9]
(x12) edge node {\iso} node[swap] {\eta^\inv_{\angx}} (x11)
(x11) edge node[swap] {f} (x21)
(x21) edge node {\eta_y} node[swap] {\iso} (x22)
(x12) edge[densely dashed] node {G(Ff)} (x22)
;
\end{tikzpicture}\]
Since the multimorphism functor
\begin{equation}\label{Gnfxfymgfxgfy}
\begin{tikzcd}[column sep=large]
\N\scmap{F\angx;Fy} \ar{r}{G}[swap]{\iso} & \M\scmap{GF\angx;GFy}
\end{tikzcd}
\end{equation}
is an isomorphism, there exists a unique $n$-ary 1-cell
\[\begin{tikzcd}[column sep=large]
F\angx \ar{r}{Ff} & Fy
\end{tikzcd} \inspace \N\scmap{F\angx;Fy}\]
such that
\begin{equation}\label{Ffdefinition}
G(Ff) = \ga\lrscmap{\eta_y; \ga\lrscmap{f; \eta^\inv_{\angx}}}. 
\end{equation}
This defines the $n$-ary 1-cell assignment, $f \mapsto Ff$, of $F$.
\item[$F$ on $n$-ary 2-cells]
Given an $n$-ary 2-cell $\theta \cn f \to f'$ in $\M\smscmap{\angx;y}$, consider the $n$-ary 2-cell
\[\ga\lrscmap{1_{\eta_y}; \ga\lrscmap{\theta; 1_{\eta^\inv_{\angx}}}} \cn GFf \to GFf' \inspace \M\scmap{GF\angx;GFy},\]
which we visualize as the long composite below.
\[\begin{tikzpicture}[xscale=3.5, yscale=1.3,vcenter]
\def\a{20}
\draw[0cell=.9]
(0,0) node (x11) {\angx}
(x11)++(1,0) node (x12) {GF\angx}
(x11)++(0,-1) node (x21) {y}
(x12)++(0,-1) node (x22) {GFy}
;
\draw[1cell=.9]
(x12) edge node {\iso} node[swap] {\eta^\inv_{\angx}} (x11)
(x11) edge[bend right=\a] node[swap,pos=.45] {f} (x21)
(x11) edge[bend left=\a] node[pos=.45] {f'} (x21)
(x21) edge node {\eta_y} node[swap] {\iso} (x22)
(x12) edge[bend right=\a, transform canvas={xshift=-1ex}] node[swap] {GFf} (x22)
(x12) edge[bend left=\a, transform canvas={xshift=1ex}] node {GFf'} (x22)
;
\draw[2cell=.9]
node[between=x11 and x21 at .6, 2label={above,\theta}] {\Rightarrow}
node[between=x12 and x22 at .65, 2label={above,G(F\theta)}] {\Rightarrow}
;
\end{tikzpicture}\]
By the isomorphism $G$ in \cref{Gnfxfymgfxgfy}, there exists a unique $n$-ary 2-cell
\[\begin{tikzcd}[column sep=large]
Ff \ar{r}{F\theta} & Ff'
\end{tikzcd} \inspace \N\scmap{F\angx;Fy}\]
such that
\begin{equation}\label{Fthetadefinition}
G(F\theta) = \ga\lrscmap{1_{\eta_y}; \ga\lrscmap{\theta; 1_{\eta^\inv_{\angx}}}}.
\end{equation}
This defines the $n$-ary 2-cell assignment, $\theta \mapsto F\theta$, of $F$.
\end{description}
This finishes the definitions of $F$ and $\eta$.  Moreover, we make the same definitions in the non-symmetric case with $G$, $F$, and $\eta$ non-symmetric.
\end{definition}

\begin{explanation}\label{expl:Fandeta}
Consider \cref{def:Fandeta}.
\begin{enumerate}[label=(\roman*)]
\item\label{expl:Fandeta-i} This definition still makes sense even if $G \cn \N \to \M$ is a non-symmetric $\Cat$-multifunctor between non-symmetric $\Cat$-multicategories because the symmetric group action on $\M$ and $\N$ and the symmetric group action axiom for $G$ are not used.
\item\label{expl:Fandeta-ii}
The definitions \cref{Ffdefinition,Fthetadefinition} are modeled after, respectively, the object and morphism $\Cat$-naturality conditions, \cref{catmultinaturality,catmultinaturalityiicell}, for $\eta \cn 1_{\M} \to GF$.  Once we prove that $F$ is a (non-symmetric) $\Cat$-multifunctor, these $\Cat$-naturality conditions for $\eta$ will then follow.\defmark
\end{enumerate}
\end{explanation}

To show that $F \cn \M \to \N$ is a (non-symmetric) $\Cat$-multifunctor, first we check its functoriality between multimorphism categories.  The following lemma holds in both the symmetric and the non-symmetric cases.

\begin{lemma}\label{Fcomponentfunctorial}
In the context of \cref{def:Fandeta}, for objects $\angx = \ang{x_j}_{j=1}^n$ and $y$ in $\M$, the assignments
\[\begin{tikzcd}[column sep=large]
\M\scmap{\angx;y} \ar{r}{F} & \N\scmap{F\angx;Fy}
\end{tikzcd}\]
form a functor.
\end{lemma}

\begin{proof}
Suppose $f \cn \angx \to y$ is an $n$-ary 1-cell in $\M$.  To show that $F$ preserves the identity $n$-ary 2-cell $1_f \cn f \to f$, we compute as follows.
\[\begin{aligned}
GF1_f &=\ga\scmap{1_{\eta_y}; \ga\smscmap{1_f; 1_{\eta^\inv_{\angx}}}} & \phantom{=} & \text{by \cref{Fthetadefinition}}\\
&= 1_{\ga\smscmap{\eta_y; \ga\smscmap{f; \eta^\inv_{\angx}}}} & \phantom{=} & \text{by functoriality of $\ga$ in $\M$}\\
&= 1_{GFf} & \phantom{=} & \text{by \cref{Ffdefinition}}\\
&= G1_{Ff} & \phantom{=} & \text{by functoriality of $G$ in \cref{Gnfxfymgfxgfy}}
\end{aligned}\]
By the injectivity of $G$ in \cref{Gnfxfymgfxgfy} on $n$-ary 2-cells, it follows that
\[F1_f = 1_{Ff} \inspace \N\scmap{F\angx;Fy}.\]
This shows that $F$ preserves identity $n$-ary 2-cells.

To show that $F$ preserves composition of $n$-ary 2-cells
\[\begin{tikzcd}[column sep=large]
f \ar{r}{\theta} & f' \ar{r}{\theta'} & f''
\end{tikzcd} \inspace \M\scmap{\angx;y},\]
we compute as follows.
\[\scalebox{.9}{$\begin{aligned}
&\phantom{==} GF(\theta' \circ \theta) &&\\
&= \ga\scmap{1_{\eta_y}; \ga\smscmap{\theta' \circ \theta; 1_{\eta^\inv_{\angx}}}} & \phantom{=} & \text{by \cref{Fthetadefinition}}\\
&= \left[\ga\scmap{1_{\eta_y}; \ga\smscmap{\theta'; 1_{\eta^\inv_{\angx}}}}\right] \circ \left[\ga\scmap{1_{\eta_y}; \ga\smscmap{\theta; 1_{\eta^\inv_{\angx}}}}\right] & \phantom{=} & \text{by functoriality of $\ga$}\\
&= (GF\theta') \circ (GF\theta) & \phantom{=} & \text{by \cref{Fthetadefinition}}\\
&= G\big( F\theta' \circ F\theta\big) & \phantom{=} & \text{by functoriality of $G$}
\end{aligned}$}\]
By the injectivity of $G$ in \cref{Gnfxfymgfxgfy} on $n$-ary 2-cells, it follows that
\[F(\theta' \circ \theta) = F\theta' \circ F\theta \inspace \N\scmap{F\angx;Fy}.\]
This finishes the proof that $F$ is a functor on each multimorphism category.
\end{proof}

\begin{lemma}\label{Fcatmultifunctor}
In the context of \cref{def:Fandeta}, 
\[\begin{tikzcd}[column sep=large]
\M \ar{r}{F} & \N
\end{tikzcd}\]
is a $\Cat$-multifunctor.  If $G$ is a non-symmetric $\Cat$-multifunctor, then $F$ is a non-symmetric $\Cat$-multifunctor.
\end{lemma}

\begin{proof}
We must show that $F$ preserves the colored units, composition, and, in the symmetric case only, symmetric group action as in \cref{def:enr-multicategory-functor}.

\newcounter{Fcatmultifunctor-step}
\newcommand{\Fcatmultifunctorstep}[1]{\medskip\textbf{Step \stepcounter{Fcatmultifunctor-step}\arabic{Fcatmultifunctor-step}}: #1\medskip}

\Fcatmultifunctorstep{Colored Units}

For each object $x \in \M$, its colored unit $1_x \cn x \to x$ is a 1-ary 1-cell in $\M$.  We compute as follows.
\[\begin{aligned}
GF1_x &= \ga\scmap{\eta_x; \ga\smscmap{1_x; \eta^\inv_x}} & \phantom{=} & \text{by \cref{Ffdefinition}}\\
&= \ga\scmap{\eta_x; \eta^\inv_x} & \phantom{=} & \text{by left unity \cref{enr-multicategory-left-unity}}\\ 
&= 1_{GFx} & \phantom{=} & \text{by \cref{etaxxgfx}}\\
&= G1_{Fx} & \phantom{=} & \text{by \cref{enr-multifunctor-unit} for $G$}
\end{aligned}\]
The injectivity of $G$ in \cref{Gnfxfymgfxgfy} on 1-ary 1-cells implies that
\[F1_x = 1_{Fx} \inspace \N\scmap{Fx;Fx}.\]
This shows that $F$ preserves colored units in the sense of \cref{enr-multifunctor-unit}.

\Fcatmultifunctorstep{Composition}

To show that $F$ preserves composition as in \cref{v-multifunctor-composition}, consider objects
\[y, \quad \angx = \ang{x_j}_{j=1}^n, \andspace \ang{w_j} = \ang{w_{ji}}_{i=1}^{\ell_j}\]
in $\M$ for $j \in \{1,\ldots,n\}$ with $\angw = \ang{\ang{w_j}}_{j=1}^n$.  We must show that the following diagram is commutative.
\begin{equation}\label{Fgammacommute}
\begin{tikzpicture}[xscale=5.5, yscale=1.5,vcenter]
\draw[0cell=.8]
(0,0) node (x11) {\M\scmap{\angx;y} \times \txprod_{j=1}^n \M\scmap{\ang{w_j};x_j}}
(x11)++(1,0) node (x12) {\N\scmap{F\angx;Fy} \times \txprod_{j=1}^n \N\scmap{F\ang{w_j};Fx_j}}
(x11)++(0,-1) node (x21) {\M\scmap{\angw;y}}
(x12)++(0,-1) node (x22) {\N\scmap{F\angw;Fy}}
;
\draw[1cell=.85]
(x11) edge node {F \times \smallprod_j F} (x12)
(x21) edge node {F} (x22)
(x11) edge node[swap] {\ga} (x21)
(x12) edge node {\ga} (x22)
;
\end{tikzpicture}
\end{equation}
Given
\begin{itemize}
\item an $n$-ary 1-cell $f \cn \angx \to y$ and
\item an $\ell_j$-ary 1-cell $f_j \cn \ang{w_j} \to x_j$ for each $j \in \{1,\ldots,n\}$
\end{itemize} 
in $\M$, we compute as follows.
\begin{equation}\label{GFgafangfj}
\left\{\begin{aligned}
&\phantom{==} GF\ga\scmap{f;\ang{f_j}} &&\\
&= \ga\Big( \eta_y \,;\, \ga\scmap{\ga\smscmap{f;\ang{f_j}}; \eta^\inv_{\angw}} \Big) & \phantom{=} & \text{by \cref{Ffdefinition}}\\
&= \scalebox{.75}{$\ga\Big( \ga\scmap{\eta_y; \ga\smscmap{f; \eta^\inv_{\angx}}} \,;\, \bang{\ga\scmap{\eta_{x_j}; \ga\smscmap{f_j; \eta^\inv_{\ang{w_j}}}}} \Big)$} & \phantom{=} & \text{by \cref{enr-multicategory-associativity,enr-multicategory-right-unity}}\\
&= \ga\scmap{GFf; GF\ang{f_j}} & \phantom{=} & \text{by \cref{Ffdefinition}}\\
&= G\ga\scmap{Ff; F\ang{f_j}} & \phantom{=} & \text{by \cref{v-multifunctor-composition} for $G$}
\end{aligned}\right.
\end{equation}
The computation \cref{GFgafangfj} and the injectivity of $G$ on $(\ell_1 + \cdots + \ell_n)$-ary 1-cells imply
\[F\ga\scmap{f;\ang{f_j}} = \ga\scmap{Ff; F\ang{f_j}} \inspace \N\scmap{F\angw;Fy}.\]
This proves the commutativity of \cref{Fgammacommute} on objects.

Suppose given
\begin{itemize}
\item an $n$-ary 2-cell $\theta \cn f \to f'$ in $\M\smscmap{\angx;y}$ and
\item an $\ell_j$-ary 2-cell $\theta_j \cn f_j \to f_j'$ in $\M\smscmap{\ang{w_j};x_j}$ for each $j \in \{1,\ldots,n\}$.
\end{itemize} 
Then
\begin{itemize}
\item the computation \cref{GFgafangfj} with $(f,f_j,\eta)$ and \cref{Ffdefinition} replaced by, respectively, $(\theta,\theta_j,1_\eta)$ and \cref{Fthetadefinition}, and
\item the injectivity of $G$ on $(\ell_1 + \cdots + \ell_n)$-ary 2-cells
\end{itemize} 
imply the commutativity of \cref{Fgammacommute} on morphisms.  This proves that $F$ preserves composition.  

The proof so far is still valid if $G$ is a non-symmetric $\Cat$-multifunctor.  In this case, we have shown that $F$ is a non-symmetric $\Cat$-multifunctor.

\Fcatmultifunctorstep{Symmetric Group Action}

This case only applies if $G$ is a $\Cat$-multifunctor.  With the notation above, $F$ preserves the symmetric group action in the sense of \cref{enr-multifunctor-equivariance} if, for each permutation $\sigma \in \Sigma_n$, the following diagram commutes.
\begin{equation}\label{Fsigmacommute}
\begin{tikzpicture}[xscale=3, yscale=1.3,vcenter]
\draw[0cell=.9]
(0,0) node (x11) {\M\scmap{\angx;y}}
(x11)++(1,0) node (x12) {\N\scmap{F\angx;Fy}}
(x11)++(0,-1) node (x21) {\M\scmap{\angx\sigma;y}}
(x12)++(0,-1) node (x22) {\N\scmap{F\angx\sigma;Fy}}
;
\draw[1cell=.9]
(x11) edge node {F} (x12)
(x21) edge node {F} (x22)
(x11) edge node[swap] {\sigma} (x21)
(x12) edge node {\sigma} (x22)
;
\end{tikzpicture}
\end{equation}
For an $n$-ary 1-cell $f \cn \angx \to y$ in $\M$, we compute as follows.
\begin{equation}\label{GFfsigma}
\left\{\begin{aligned}
&\phantom{==} G\big((Ff) \sigma\big) &&\\
&= (GFf) \sigma & \phantom{=} & \text{by \cref{enr-multifunctor-equivariance} for $G$}\\
&= \ga\scmap{\eta_y; \ga\smscmap{f; \eta^\inv_{\angx}}} \sigma & \phantom{=} & \text{by \cref{Ffdefinition}}\\
&= \ga\scmap{\eta_y; \ga\smscmap{f; \eta^\inv_{\angx}} \sigma} & \phantom{=} & \text{by bottom eq. \cref{enr-operadic-eq-2}}\\
&= \ga\scmap{\eta_y; \ga\smscmap{f\sigma; \eta^\inv_{\angx\sigma}}} & \phantom{=} & \text{by top eq. \cref{enr-operadic-eq-1}}\\
&= GF(f\sigma) & \phantom{=} & \text{by \cref{Ffdefinition}}
\end{aligned}\right.
\end{equation}
The injectivity of $G$ on $n$-ary 1-cells implies
\[(Ff)\sigma = F(f\sigma).\]
This proves the commutativity of \cref{Fsigmacommute} on objects.

The commutativity of \cref{Fsigmacommute} on morphisms is proved using
\begin{itemize}
\item the computation \cref{GFfsigma} with $f$, $\eta$, and \cref{Ffdefinition} replaced by, respectively, an $n$-ary 2-cell in $\M\smscmap{\angx;y}$, $1_{\eta}$, and \cref{Fthetadefinition}, and
\item the injectivity of $G$ on $n$-ary 2-cells. 
\end{itemize}
This finishes the proof that $F$ preserves the symmetric group action if $G$ is a $\Cat$-multifunctor.
\end{proof}

Recall from \cref{expl:catmultitransformation} the description of a (non-symmetric) $\Cat$-multinatural transformation or isomorphism.

\begin{lemma}\label{etacatmultinatural}
In the context of \cref{def:Fandeta},
\[\begin{tikzcd}[column sep=large]
1_{\M} \ar{r}{\eta} & GF
\end{tikzcd}\]
is a $\Cat$-multinatural isomorphism.  If $G$ is a non-symmetric $\Cat$-multifunctor, then $\eta$ is a non-symmetric $\Cat$-multinatural isomorphism.
\end{lemma}

\begin{proof}
In \cref{Fcatmultifunctor} we showed that $F \cn \M \to \N$ is a (non-symmetric) $\Cat$-multifunctor.  By definition \cref{etaxxgfx} each component 1-ary 1-cell $\eta_x \cn x \to GFx$ is an isomorphism in $\M$.  It remains to check the two $\Cat$-naturality conditions, \cref{catmultinaturality,catmultinaturalityiicell}, for $\eta$.

For an $n$-ary 1-cell $f \cn \angx \to y$ in $\M$ with $\angx = \ang{x_j}_{j=1}^n$, the following equalities prove the object $\Cat$-naturality condition \cref{catmultinaturality} for $\eta \cn 1_{\M} \to GF$.
\begin{equation}\label{etaobjcatnaturality}
\left\{\begin{aligned}
&\phantom{==} \ga\scmap{GFf; \eta_{\angx}} &&\\
&= \ga\Big(\ga\scmap{\eta_y; \ga\smscmap{f; \eta^\inv_{\angx}}} \sscs \eta_{\angx}\Big) & \phantom{=} & \text{by \cref{Ffdefinition}}\\
&= \ga\Big(\eta_y \sscs \ga\scmap{f; \bang{\ga\smscmap{\eta^\inv_{x_j}; \eta_{x_j}}}}\Big) & \phantom{=} & \text{by \cref{enr-multicategory-associativity}}\\
&= \ga\scmap{\eta_y;f} & \phantom{=} & \text{by \cref{enr-multicategory-right-unity}}
\end{aligned}\right.
\end{equation}
The morphism $\Cat$-naturality condition \cref{catmultinaturalityiicell} for $\eta$ is proved using the computation \cref{etaobjcatnaturality} by replacing $f$, $\eta$, and \cref{Ffdefinition} with, respectively, an $n$-ary 2-cell in $\M\smscmap{\angx;y}$, $1_{\eta}$, and \cref{Fthetadefinition}.
\end{proof}

\subsection*{Counit}

Next we define the counit for $(F,G)$.  The formula \cref{Gepzz} below comes from the right triangle identity \cref{triangleidentities} of an adjunction.

\begin{definition}\label{def:counitepz}
In the context of \cref{def:Fandeta}, define the data of a $\Cat$-multinatural isomorphism
\[\begin{tikzcd}[column sep=large]
FG \ar{r}{\epz}[swap]{\iso} & 1_{\N}
\end{tikzcd}\]
as follows.  For each object $z \in \N$, the $(Gz)$-component 1-ary 1-cell of $\eta$ in $\M$
\[\begin{tikzcd}[column sep=large]
Gz \ar{r}{\eta_{Gz}}[swap]{\iso} & GFGz
\end{tikzcd}\]
is an isomorphism \cref{etaxxgfx}, with inverse denoted by
\[\begin{tikzcd}[column sep=large]
GFGz \ar{r}{\eta^\inv_{Gz}}[swap]{\iso} & Gz. 
\end{tikzcd}\]
Since the multimorphism functor
\[\begin{tikzcd}[column sep=large]
\N\scmap{FGz;z} \ar{r}{G}[swap]{\iso} & \M\scmap{GFGz;Gz}
\end{tikzcd}\]
is an isomorphism, there exists a unique invertible 1-ary 1-cell in $\N$
\begin{equation}\label{epzzcomponent}
\begin{tikzcd}[column sep=large]
FGz \ar{r}{\epz_z}[swap]{\iso} & z
\end{tikzcd}
\end{equation}
such that
\begin{equation}\label{Gepzz}
G\epz_z = \eta^\inv_{Gz}.
\end{equation}
This defines the $z$-component of $\epz$.
\end{definition}

\begin{lemma}\label{epzcatmultiiso}
In the context of \cref{def:Fandeta,def:counitepz}, 
\[\begin{tikzcd}[column sep=large]
FG \ar{r}{\epz} & 1_{\N}
\end{tikzcd}\]
is a $\Cat$-multinatural isomorphism.   If $G$ is a non-symmetric $\Cat$-multifunctor, then $\epz$ is a non-symmetric $\Cat$-multinatural isomorphism.
\end{lemma}

\begin{proof}
Since each component of $\epz$ is an isomorphism \cref{epzzcomponent}, it suffices to check the two $\Cat$-naturality conditions, \cref{catmultinaturality,catmultinaturalityiicell}, for $\epz$.  Suppose $g \cn \angw \to z$ is an $n$-ary 1-cell in $\N$ with $\angw = \ang{w_j}_{j=1}^n$.  The injectivity of $G$ on $n$-ary 1-cells and the following equalities prove the object $\Cat$-naturality condition \cref{catmultinaturality} for $\epz \cn FG \to 1_{\N}$.
\begin{equation}\label{Ggafepzxj}
\left\{\begin{aligned}
&\phantom{==} G\ga\scmap{g; \epz_{\angw}} &&\\
&= \ga\scmap{Gg; G\epz_{\angw}} & \phantom{=} & \text{by \cref{v-multifunctor-composition} for $G$}\\
&= \ga\scmap{Gg; \eta^\inv_{G\angw}} & \phantom{=} & \text{by \cref{Gepzz}}\\
&= \ga\scmap{\eta^\inv_{Gz}; GFGg} & \phantom{=} & \text{by \cref{enr-multicategory-associativity,enr-multicategory-right-unity,enr-multicategory-left-unity,etaobjcatnaturality}}\\
&= \ga\scmap{G\epz_z; GFGg} & \phantom{=} & \text{by \cref{Gepzz}}\\
&= G\ga\scmap{\epz_z; FGg} & \phantom{=} & \text{by \cref{v-multifunctor-composition} for $G$}
\end{aligned}\right.
\end{equation}
The morphism  $\Cat$-naturality condition \cref{catmultinaturalityiicell} for $\epz$ is proved using
\begin{itemize}
\item the injectivity of $G$ on $n$-ary 2-cells and
\item the computation \cref{Ggafepzxj} with $g$, $\eta$, $\epz$, and \cref{etaobjcatnaturality} replaced by, respectively, an $n$-ary 2-cell in $\N$, $1_{\eta}$, $1_{\epz}$, and the morphism $\Cat$-naturality condition \cref{catmultinaturalityiicell} for $\eta$ (\cref{etacatmultinatural}).
\end{itemize}
This finishes the proof.
\end{proof}

\begin{lemma}\label{multieqtriangleids}
In the context of \cref{def:Fandeta,def:counitepz}, the triangle identities \cref{triangleidentities} hold for $(F,G,\eta,\epz)$.
\end{lemma}

\begin{proof}
The right triangle identity 
\[\ga\scmap{G\epz_z; \eta_{Gz}} = 1_{Gz}\]
follows from the definition \cref{Gepzz} of $\epz_z$ by applying $\ga\smscmap{-;\eta_{Gz}}$ to both sides. 

The left triangle identity follows from the following computation for each object $x \in \M$ and the injectivity of $G$ on 1-ary 1-cells.
\[\begin{aligned}
&\phantom{==} G\ga\scmap{\epz_{Fx}; F\eta_x} &&\\
&= \ga\scmap{G\epz_{Fx}; GF\eta_x} & \phantom{=} & \text{by \cref{v-multifunctor-composition} for $G$}\\
&= \ga\Big(\eta^\inv_{GFx} \sscs \ga\scmap{\eta_{GFx}; \ga\smscmap{\eta_x; \eta^\inv_x}} \Big) & \phantom{=} & \text{by \cref{Ffdefinition,Gepzz}}\\
&=1_{GFX} & \phantom{=} & \text{by \cref{enr-multicategory-right-unity}}\\
&= G1_{Fx} & \phantom{=} & \text{by \cref{enr-multifunctor-unit} for $G$}\\
\end{aligned}\]
This finishes the proof.
\end{proof}

\subsection*{$\Cat$-Multiequivalences}

We are now ready for the main result of this section, which is a practical characterization of a $\Cat$-multiequivalence (\cref{def:catmultiequivalence}).

\begin{theorem}[Multiequivalence]\label{thm:multiwhitehead}\index{Cat-multiequivalence@$\Cat$-multiequivalence!characterization}
Suppose
\[\begin{tikzcd}[column sep=large]
\N \ar{r}{G} & \M
\end{tikzcd}\]
is a $\Cat$-multifunctor between $\Cat$-multicategories.  Then the following three statements are equivalent.
\begin{enumerate}
\item\label{thm:multiwhitehead-i}
$G$ is the right adjoint of an adjoint $\Cat$-multiequivalence $(F,G, \eta,\epz)$.
\item\label{thm:multiwhitehead-ii}
$G$ is a $\Cat$-multiequivalence.
\item\label{thm:multiwhitehead-iii}
$G$ is
\begin{itemize}
\item essentially surjective on objects and
\item an isomorphism on each multimorphism category.
\end{itemize}
\end{enumerate}
Moreover, the non-symmetric variant is also true, where each of $\M$, $\N$, $G$, and $(F,G, \eta,\epz)$ is non-symmetric.
\end{theorem}

\begin{proof}
We consider each implication as follows.
\begin{itemize}
\item \eqref{thm:multiwhitehead-i} $\Rightarrow$ \eqref{thm:multiwhitehead-ii} holds by \cref{def:equivalences,def:catmultiequivalence}.
\item \eqref{thm:multiwhitehead-ii} $\Rightarrow$ \eqref{thm:multiwhitehead-iii} holds by \cref{catmultieqesff}.
\item \eqref{thm:multiwhitehead-iii} $\Rightarrow$ \eqref{thm:multiwhitehead-i} holds by \cref{Fcatmultifunctor,etacatmultinatural,epzcatmultiiso,multieqtriangleids}.
\end{itemize}
This completes the proof.
\end{proof}

Recall from \cref{def:iiequivalence} the notions of (adjoint) 2-equivalences.  Restricting the proofs of the lemmas referenced in the proof of \cref{thm:multiwhitehead} to the 1-ary case recovers the following characterization of a 2-equivalence \cite[7.5.8]{johnson-yau}.

\begin{theorem}[2-equivalence]\label{thm:iiequivalence}\index{2-equivalence!characterization}
Suppose
\[\begin{tikzcd}[column sep=large]
\C \ar{r}{G} & \B
\end{tikzcd}\]
is a 2-functor between 2-categories.  Then the following three statements are equivalent.
\begin{enumerate}
\item\label{thm:iiequivalence-i}
$G$ is the right adjoint of an adjoint 2-equivalence $(F,G, \eta,\epz)$.
\item\label{thm:iiequivalence-ii}
$G$ is a 2-equivalence.
\item\label{thm:iiequivalence-iii}
$G$ is
\begin{itemize}
\item essentially surjective on objects and
\item an isomorphism on each hom category.
\end{itemize}
\end{enumerate}
\end{theorem}

In \cite[7.5.8]{johnson-yau} this theorem is obtained by restricting the proof of an analogous characterization of a biequivalence between bicategories \cite[7.4.1]{johnson-yau}.  In contrast, here we obtain \cref{thm:iiequivalence} by restricting the proof of \cref{thm:multiwhitehead}, which is a characterization of a $\Cat$-multiequivalence.  See \cref{qu:multibicat} for a possible generalization involving a lax variant of a $\Cat$-multicategory.

\chapter{Pseudo Symmetry}
\label{ch:pseudosymmetry}
In this chapter we define \emph{pseudo symmetric $\Cat$-multifunctors} (\cref{def:pseudosmultifunctor}) between $\Cat$-multicategories.  They are weaker variants of $\Cat$-multifunctors that preserve the symmetric group action up to coherent isomorphisms.  The corresponding notion of a \emph{pseudo symmetric $\Cat$-multinatural transformation} is in \cref{def:psmftransformation}.  These notions are the 1-cells and 2-cells of a 2-category $\catmulticatps$ with small $\Cat$-multicategories as objects (\cref{thm:catmulticatps}).   \cref{thm:psmultiwhitehead} characterizes pseudo symmetric $\Cat$-multiequivalences in terms of essential surjectivity on objects and invertibility on multimorphism categories.  This is the pseudo symmetric analog of \cref{thm:multiwhitehead}.  Here is a summary table of some of the definitions and results in this chapter.
\smallskip
\begin{center}
\resizebox{0.85\width}{!}{
{\renewcommand{\arraystretch}{1.4}%
{\setlength{\tabcolsep}{1ex}
\begin{tabular}{|c|c|}\hline
pseudo symmetric $\Cat$-multifunctors & \cref{def:pseudosmultifunctor} \\ \hline
pseudo symmetric $\Cat$-multinatural transformations & \cref{def:psmftransformation} \\ \hline
2-category $\catmulticatps$ & \cref{thm:catmulticatps} \\ \hline
pseudo symmetric $\Cat$-multiequivalences & \cref{def:pscatmultiequivalence} \\ \hline
characterization of pseudo symmetric $\Cat$-multiequivalences & \cref{thm:psmultiwhitehead} \\ \hline 
\end{tabular}}}}
\end{center}
\smallskip
Pseudo symmetric $\Cat$-multifunctors are utilized in subsequent chapters as follows.

\subsection*{The Grothendieck Construction and Inverse $K$-Theory}

The prime example of a pseudo symmetric $\Cat$-multifunctor is the Grothendieck construction 
\[\begin{tikzpicture}[xscale=2.75,yscale=1.2,vcenter]
\draw[0cell=.9]
(0,0) node (dc) {\DCat}
(dc)++(1,0) node (pc) {\permcatsg};
\draw[1cell=.9]  
(dc) edge node[pos=.5] {\grod} (pc);
\end{tikzpicture}\]
on any small tight bipermutative category $(\cD,\oplus,\otimes)$ (\cref{thm:grocatmultifunctor}). 
\begin{itemize}
\item The domain $\Cat$-multicategory $\DCat$ (\cref{thm:dcatcatmulticat}) has additive symmetric monoidal functors $\Dplus \to \Cat$ as objects.  Its multimorphism categories (\cref{def:dcatnarycategory}) are defined using the full extent of the bipermutative structure of $\cD$ and Laplaza's Coherence \cref{thm:laplaza-coherence-1}.
\item The codomain $\Cat$-multicategory $\permcatsg$ has small permutative categories as objects (\cref{thm:permcatenrmulticat}).
\end{itemize}
An important consequence is that inverse $K$-theory $\groa \circ A$, which involves the Grothendieck construction $\groa$ on the small tight bipermutative category $\cA$ in \cref{ex:mandellcategory}, is a pseudo symmetric $\Cat$-multifunctor but \emph{not} a $\Cat$-multifunctor (\cref{thm:invKpseudosym}).

Taking the projection down to $\cD$ into account, the Grothendieck construction above lifts to a non-symmetric $\Cat$-multifunctor \cref{groddcatpfibd}
\[\begin{tikzpicture}[xscale=2.75,yscale=1.2,baseline={(x1.base)}]
\draw[0cell=.9]
(0,0) node (x1) {\DCat}
(x1)++(1,0) node (x2) {\pfibd.};
\draw[1cell=.9]  
(x1) edge node[pos=.5] {\grod} (x2);
\end{tikzpicture}\]
The codomain $\pfibd$ is a non-symmetric $\Cat$-multicategory with small permutative opfibrations over $\cD$ as objects (\cref{thm:pfibdmulticat}).  While these objects only use the additive structure $\Dplus$, the multimorphism categories of $\pfibd$ (\cref{def:pfibd}) make full use of the bipermutative structure of $\cD$.  The characterization of a $\Cat$-multiequivalence in \cref{thm:multiwhitehead} is used to show that this lifted Grothendieck construction $\grod$ is a non-symmetric $\Cat$-multiequivalence (\cref{thm:dcatpfibdeq}).

\subsection*{Organization}

In \cref{sec:pseudosmultifunctor} we define \emph{pseudo symmetric $\Cat$-multifunctors}.  These weaker variants of $\Cat$-multifunctors are only non-strict with respect to the symmetric group action.  They are required to \emph{strictly} preserve the colored units and the $\Cat$-multicategorical composition.  We define them this way because the Grothendieck construction $\grod$ has precisely these properties (\cref{thm:grocatmultifunctor}).  See \cref{expl:laxsmf,qu:multibicat} regarding further lax variants of $\Cat$-multicategories and $\Cat$-multifunctors.

In \cref{sec:pstransformation} we define \emph{pseudo symmetric $\Cat$-multinatural transformations} between  pseudo symmetric $\Cat$-multifunctors.  They have the same data as a $\Cat$-multinatural transformation, namely, component 1-ary 1-cells, and satisfy the $\Cat$-naturality conditions \cref{catmultinaturality,catmultinaturalityiicell}.  Moreover, each pseudo symmetric $\Cat$-multinatural  transformation is compatible with the pseudo symmetry isomorphisms of its domain and codomain.  There is a 2-category $\catmulticatps$ with pseudo symmetric $\Cat$-multifunctors as 1-cells and pseudo symmetric $\Cat$-multinatural transformations as 2-cells (\cref{thm:catmulticatps}).

In \cref{sec:pstcharacterization} we extend \cref{thm:multiwhitehead} to the pseudo symmetric context.  \cref{thm:psmultiwhitehead} proves that a pseudo symmetric $\Cat$-multifunctor is a pseudo symmetric $\Cat$-\emph{multiequivalence} if and only if it is essentially surjective on objects and an isomorphism on each multimorphism category.

We remind the reader of our left normalized bracketing \cref{expl:leftbracketing}.

\section{Pseudo Symmetric \texorpdfstring{$\Cat$}{Cat}-Multifunctors}
\label{sec:pseudosmultifunctor}

In this section we define a variant of a $\Cat$-multifunctor $F \cn \M \to \N$ between $\Cat$-multicategories, called a \emph{pseudo symmetric} $\Cat$-multifunctor, that preserves the symmetric group action \cref{enr-multifunctor-equivariance} up to coherent natural isomorphisms.  Pseudo symmetric $\Cat$-multifunctors are closed under composition (\cref{pseudosmfcomposite}), which is strictly associative and unital.  Transformations between them are discussed in \cref{sec:pstransformation}. 

Our main use of the notion of a pseudo symmetric $\Cat$-multifunctor is \cref{thm:grocatmultifunctor}.  This result  shows that, for each small tight bipermutative category $\cD$, the Grothendieck construction $\grod$ is a pseudo symmetric $\Cat$-multifunctor that is, in general, not a $\Cat$-multifunctor.  As a practical consequence, inverse $K$-theory is a pseudo symmetric $\Cat$-multifunctor but not a $\Cat$-multifunctor (\cref{thm:invKpseudosym}).

\begin{definition}\label{def:pseudosmultifunctor}\index{pseudo symmetric!Cat-multifunctor@$\Cat$-multifunctor}\index{multifunctor!pseudo symmetric}\index{Cat-multifunctor@$\Cat$-multifunctor!pseudo symmetric}
Suppose $\M$ and $\N$ are $\Cat$-multicategories.  A \emph{pseudo symmetric $\Cat$-multifunctor}
\[\begin{tikzcd}[column sep=2.5cm]
\M \ar{r}{\left(F \scs \{F_{\sigma,\angc,c'}\}\right)} & \N
\end{tikzcd}\]
consists of
\begin{itemize}
\item an object assignment $F \cn \ObM \to \ObN$;
\item for each $\smscmap{\angc;c'} \in \ProfMM$ with $\angc=\ang{c_j}_{j=1}^n$, a component functor
\[\begin{tikzcd}[column sep=large]
\M\mmap{c';\ang{c}} \ar{r}{F} & \N\mmap{Fc';F\ang{c}}
\end{tikzcd}\]
with $F\angc=\ang{Fc_j}_{j=1}^n$; and
\item for each permutation $\sigma \in \Sigma_n$ and object $\smscmap{\angc;c'}$ as above, a natural isomorphism
\begin{equation}\label{Fsigmaangccp}
\begin{tikzpicture}[xscale=3.5,yscale=1.3,vcenter]
\draw[0cell=.9]
(0,0) node (x11) {\M\scmap{\angc;c'}}
(x11)++(1,0) node (x12) {\N\scmap{F\angc;Fc'}}
(x11)++(0,-1) node (x21) {\M\scmap{\angc\sigma;c'}}
(x12)++(0,-1) node (x22) {\N\scmap{F\angc\sigma;Fc'}}
;
\draw[1cell=.9]
(x11) edge node (f1) {F} (x12)
(x21) edge node[swap] (f2) {F} (x22)
(x11) edge node[swap] {\sigma} (x21)
(x12) edge node {\sigma} (x22)
;
\draw[2cell=.9]
node[between=f1 and f2 at .6, 2label={above,F_{\sigma,\angc,c'}}, 2label={below,\iso}] {\Longrightarrow}
;
\end{tikzpicture}
\end{equation}
called the \index{pseudo symmetry isomorphism}\emph{pseudo symmetry isomorphism}.
\end{itemize}
We abbreviate $F_{\sigma,\angc,c'}$ to $F_{\sigma}$ when the domain object $\smscmap{\angc;c'}$ can be inferred from the context.  The above data are subject to
\begin{itemize}
\item the unit axiom \cref{enr-multifunctor-unit},
\item the composition axiom \cref{v-multifunctor-composition}, and
\item the four axioms \cref{pseudosmf-unity,pseudosmf-product,pseudosmf-topeq,pseudosmf-boteq} below with $\smscmap{\angc;c'}$ as above.
\end{itemize} 
\begin{description}
\item[Unit Permutation] 
For the identity permutation $\id_n \in \Sigma_n$, the following equality holds.\index{unit permutation axiom}
\begin{equation}\label{pseudosmf-unity}
F_{\id_n, \angc, c'} = 1_F
\end{equation}
\item[Product Permutation]
For $\sigma,\tau \in \Sigma_n$ the following pasting diagram equality holds.\index{product permutation axiom}
\begin{equation}\label{pseudosmf-product}
\begin{tikzpicture}[xscale=2.5,yscale=1.3,vcenter]
\def\s{.8} \def\t{.9}
\def\boundary{
\draw[0cell=\s]
(0,0) node (x11) {\M\scmap{\angc;c'}}
(x11)++(1,0) node (x12) {\N\scmap{F\angc;Fc'}}
(x11)++(0,-2) node (x31) {\M\scmap{\angc\sigma\tau;c'}}
(x12)++(0,-2) node (x32) {\N\scmap{F\angc\sigma\tau;Fc'}}
;
\draw[1cell=\s]
(x11) edge node (f1) {F} (x12)
(x31) edge node (f3) {F} (x32)
;}
\draw (0,0) node[font=\LARGE] (equal) {=};
\begin{scope}[shift={(-1.6,1)}]
\boundary
\draw[0cell=\s] 
(x11)++(0,-1) node (x21) {\M\scmap{\angc\sigma;c'}}
(x12)++(0,-1) node (x22) {\N\scmap{F\angc\sigma;Fc'}}
;
\draw[1cell=\s]
(x21) edge node (f2) {F} (x22)
(x11) edge node[swap] {\sigma} (x21)
(x12) edge node {\sigma} (x22)
(x21) edge node[swap] {\tau} (x31)
(x22) edge node {\tau} (x32)
;
\draw[2cell=\t]
node[between=f1 and f2 at .6, 2label={above,F_\sigma}] {\Longrightarrow}
node[between=f2 and f3 at .6, 2label={above,F_\tau}] {\Longrightarrow}
;
\end{scope}
\begin{scope}[shift={(.5,1)}]
\boundary
\draw[1cell=\s]
(x11) edge node[swap] {\sigma\tau} (x31)
(x12) edge node {\sigma\tau} (x32)
;
\draw[2cell=\t]
node[between=f1 and f3 at .6, 2label={above,F_{\sigma\tau}}] {\Longrightarrow}
;
\end{scope}
\end{tikzpicture}
\end{equation}
\item[Top Equivariance]
In the context of \cref{enr-operadic-eq-1} we use the following shorter notation.\label{not:topeq}
\[\left\{
\scalebox{.9}{$\begin{aligned}
\ang{\ang{c_j}}_{j=1}^n &= \big(\ang{c_1}, \ldots, \ang{c_n}\big) & \ang{\ang{c_{\sigma(j)}}}_{j=1}^n &= \big(\ang{c_{\sigma(1)}}, \ldots, \ang{c_{\sigma(n)}}\big)\\
\ang{F\ang{c_j}}_{j=1}^n &= \big(F\ang{c_1}, \ldots, F\ang{c_n}\big) & \ang{F\ang{c_{\sigma(j)}}}_{j=1}^n &= \big(F\ang{c_{\sigma(1)}}, \ldots, F\ang{c_{\sigma(n)}}\big)\\ 
\sigmabar &= \sigma\langle k_{\sigma(1)}, \ldots , k_{\sigma(n)} \rangle &&
\end{aligned}$}\right.\]
Then the following pasting diagram equality holds.\index{top equivariance axiom}
\begin{equation}\label{pseudosmf-topeq}
\begin{tikzpicture}[xscale=3.1,yscale=1.5,vcenter]
\def\s{.65} \def\t{.85}
\def\boundary{
\draw[0cell=\s]
(0,0) node[align=left] (x11) {$\M\scmap{\angcp;c''} \,\times$\\
$\txprod_{j=1}^n \M\scmap{\ang{c_j};c_j'}$}
(x11)++(1,0) node[align=left] (x12) {$\N\scmap{F\angcp;Fc''} \,\times$\\
$\txprod_{j=1}^n \N\scmap{F\ang{c_j};Fc_j'}$}
(x11)++(0,-2) node (x31) {\M\scmap{\ang{\ang{c_{\sigma(j)}}}_{j=1}^n;c''}}
(x12)++(0,-2) node (x32) {\N\scmap{\ang{F\ang{c_{\sigma(j)}}}_{j=1}^n;Fc''}}
;
\draw[1cell=\s]
(x11) edge node (f1) {F \times} node[swap] {\smallprod_j F} (x12)
(x31) edge node (f3) {F} (x32)
;}
\draw (0,0) node[font=\LARGE] (equal) {=};
\begin{scope}[shift={(-1.4,.7)}]
\boundary
\draw[0cell=\s] 
(x11)++(0,-1) node (x21) {\M\scmap{\ang{\ang{c_j}}_{j=1}^n;c''}}
(x12)++(0,-1) node (x22) {\N\scmap{\ang{F\ang{c_j}}_{j=1}^n;Fc''}}
;
\draw[1cell=\s]
(x21) edge node (f2) {F} (x22)
(x11) edge node[swap] {\ga} (x21)
(x12) edge node {\ga} (x22)
(x21) edge node[swap] {\sigmabar} (x31)
(x22) edge node {\sigmabar} (x32)
;
\draw[2cell=\t]
node[between=f1 and f2 at .6] {=}
node[between=f2 and f3 at .65, 2label={above,F_{\sigmabar}}] {\Longrightarrow}
;
\end{scope}
\begin{scope}[shift={(.4,.7)}]
\boundary
\draw[0cell=\s] 
(x11)++(0,-1.2) node[align=left] (x21) {$\M\scmap{\angcp\sigma;c''} \,\times$\\
$\txprod_{j=1}^n \M\scmap{\ang{c_{\sigma(j)}};c_{\sigma(j)}'}$}
(x12)++(0,-1.2) node[align=left] (x22) {$\N\scmap{F\angcp\sigma;Fc''} \,\times$\\
$\txprod_{j=1}^n \N\scmap{F\ang{c_{\sigma(j)}};Fc_{\sigma(j)}'}$}
;
\draw[1cell=\s]
(x21) edge node (f2) {F \times} node[swap] {\smallprod_j F} (x22)
(x11) edge node[swap] {(\sigma,\sigmainv)} (x21)
(x12) edge node {(\sigma,\sigmainv)} (x22)
(x21) edge node[swap] {\ga} (x31)
(x22) edge node {\ga} (x32)
;
\draw[2cell=\t]
node[between=f1 and f2 at .65, 2label={above,F_{\sigma} \times 1}] {\Longrightarrow}
node[between=f2 and f3 at .6] {=}
;
\end{scope}
\end{tikzpicture}
\end{equation}
\item[Bottom Equivariance]
In the context of \cref{enr-operadic-eq-2}, we use the following shorter notation.\label{not:boteq}
\[\left\{
\scalebox{.9}{$\begin{aligned}
\ang{\ang{c_j}\tau_j}_{j=1}^n &= \big(\ang{c_1}\tau_1, \ldots, \ang{c_n}\tau_n\big) & \tautimes &= \tau_1 \times \cdots \times \tau_n\\ 
\ang{F\ang{c_j}\tau_j}_{j=1}^n &= \big(F\ang{c_1}\tau_1, \ldots, F\ang{c_n}\tau_n\big) &&
\end{aligned}$}\right.\]
Then the following pasting diagram equality holds.\index{bottom equivariance axiom}
\begin{equation}\label{pseudosmf-boteq}
\begin{tikzpicture}[xscale=3.1,yscale=1.5,vcenter]
\def\s{.65} \def\t{.85}
\def\boundary{
\draw[0cell=\s]
(0,0) node[align=left] (x11) {$\M\scmap{\angcp;c''} \,\times$\\
$\txprod_{j=1}^n \M\scmap{\ang{c_j};c_j'}$}
(x11)++(1,0) node[align=left] (x12) {$\N\scmap{F\angcp;Fc''} \,\times$\\
$\txprod_{j=1}^n \N\scmap{F\ang{c_j};Fc_j'}$}
(x11)++(0,-2) node (x31) {\M\scmap{\ang{\ang{c_j}\tau_j}_{j=1}^n;c''}}
(x12)++(0,-2) node (x32) {\N\scmap{\ang{F\ang{c_j}\tau_j}_{j=1}^n;Fc''}}
;
\draw[1cell=\s]
(x11) edge node (f1) {F \times} node[swap] {\smallprod_j F} (x12)
(x31) edge node (f3) {F} (x32)
;}
\draw (0,0) node[font=\LARGE] (equal) {=};
\begin{scope}[shift={(-1.4,.7)}]
\boundary
\draw[0cell=\s] 
(x11)++(0,-1) node (x21) {\M\scmap{\ang{\ang{c_j}}_{j=1}^n;c''}}
(x12)++(0,-1) node (x22) {\N\scmap{\ang{F\ang{c_j}}_{j=1}^n;Fc''}}
;
\draw[1cell=\s]
(x21) edge node (f2) {F} (x22)
(x11) edge node[swap] {\ga} (x21)
(x12) edge node {\ga} (x22)
(x21) edge node[swap] {\tautimes} (x31)
(x22) edge node {\tautimes} (x32)
;
\draw[2cell=\t]
node[between=f1 and f2 at .6] {=}
node[between=f2 and f3 at .65, 2label={above,F_{\tautimes}}] {\Longrightarrow}
;
\end{scope}
\begin{scope}[shift={(.4,.7)}]
\boundary
\draw[0cell=\s] 
(x11)++(0,-1.2) node[align=left] (x21) {$\M\scmap{\angcp;c''} \,\times$\\
$\txprod_{j=1}^n \M\scmap{\ang{c_j}\tau_j;c_j'}$}
(x12)++(0,-1.2) node[align=left] (x22) {$\N\scmap{F\angcp;Fc''} \,\times$\\
$\txprod_{j=1}^n \N\scmap{F\ang{c_j}\tau_j;Fc_j'}$}
;
\draw[1cell=\s]
(x21) edge node (f2) {F \times} node[swap] {\smallprod_j F} (x22)
(x11) edge node[swap] {1 \times \smallprod_j \tau_j} (x21)
(x12) edge node {1 \times \smallprod_j \tau_j} (x22)
(x21) edge node[swap] {\ga} (x31)
(x22) edge node {\ga} (x32)
;
\draw[2cell=\t]
node[between=f1 and f2 at .65, 2label={above,1_F \times \smallprod_j F_{\tau_j}}] {\Longrightarrow}
node[between=f2 and f3 at .6] {=}
;
\end{scope}
\end{tikzpicture}
\end{equation}
\end{description}
This finishes the definition of a pseudo symmetric $\Cat$-multifunctor.

Moreover, the \emph{identity pseudo symmetric $\Cat$-multifunctor}\index{pseudo symmetric!identity - $\Cat$-multifunctor} of $\M$, denoted $1_{\M}$, consists of
\begin{itemize}
\item the identity object assignment on $\Ob\M$,
\item identity component functors, and
\item identity natural transformations as the pseudo symmetry isomorphisms.
\end{itemize}
This finishes the definition of $1_{\M}$.
\end{definition}

\begin{example}[$\Cat$-Multifunctors]\label{ex:idpsm}\index{multifunctor}
The domain and codomain of the pseudo symmetry isomorphism in \cref{Fsigmaangccp} form the boundary of the symmetric group action axiom \cref{enr-multifunctor-equivariance}.  Therefore, there is a canonical bijection between
\begin{itemize}
\item $\Cat$-multifunctors $F \cn \M \to \N$ and
\item pseudo symmetric $\Cat$-multifunctors $(F,\{F_\sigma\}) \cn \M \to \N$ with each pseudo symmetry isomorphism $F_\sigma$ given by the identity natural transformation.  
\end{itemize} 
In this work, when a $\Cat$-multifunctor $F$ is regarded as a pseudo symmetric $\Cat$-multifunctor, we always mean taking each $F_\sigma$ as the identity natural transformation.  See \cref{qu:pseudo} for a possible variation.
\end{example}

\begin{example}[Non-Symmetric $\Cat$-Multifunctors]\label{ex:nscatmultifunctors}\index{non-symmetric!multifunctor}
Each pseudo symmetric $\Cat$-multifunctor $(F,\{F_\sigma\}) \cn \M \to \N$ yields a non-symmetric $\Cat$-multifunctor $F \cn \M \to \N$ by forgetting the pseudo symmetry isomorphisms $\{F_\sigma\}$.
\end{example}

\begin{explanation}[Pseudo Symmetry Data]\label{expl:Fsigmadata}
The pseudo symmetry isomorphism $F_\sigma$ in \cref{Fsigmaangccp} consists of, for each $n$-ary 1-cell $p \in \M\smscmap{\angc;c'}$, an invertible $n$-ary 2-cell\label{not:Fsubsigmap}
\[\begin{tikzcd}[column sep=large]
F(p^\sigma) \ar{r}{F_{\sigma;p}}[swap]{\iso} & (Fp)^\sigma
\end{tikzcd} \inspace \N\scmap{F\angc\sigma; Fc'}\]
such that, for each $n$-ary 2-cell $f \cn p \to q$ in $\M\smscmap{\angc;c'}$, the naturality diagram
\[\begin{tikzpicture}[xscale=3,yscale=1.2,vcenter]
\draw[0cell=.9]
(0,0) node (x11) {F(p^\sigma)}
(x11)++(1,0) node (x12) {F(q^\sigma)}
(x11)++(0,-1) node (x21) {(Fp)^\sigma}
(x12)++(0,-1) node (x22) {(Fq)^\sigma}
;
\draw[1cell=.9]
(x11) edge node {F(f^\sigma)} (x12)
(x21) edge node (f2) {(Ff)^\sigma} (x22)
(x11) edge node[swap] {F_{\sigma;p}} (x21)
(x12) edge node {F_{\sigma;q}} (x22)
;
\end{tikzpicture}\]
is commutative.  Here and below we write the right $\sigma$-action on $p$ as\label{not:ptosigma}
\[p^\sigma = p \cdot \sigma \in \M\scmap{\angc\sigma;c'}\]
and the $p$-component of $F_\sigma$ as $F_{\sigma;p}$ or $F_{\sigma,p}$.
\end{explanation}

\begin{explanation}[Pseudo Symmetry Axioms]\label{expl:pseudosmf}
Consider \cref{def:pseudosmultifunctor} of a pseudo symmetric $\Cat$-multifunctor $(F,\{F_\sigma\})$.
\begin{enumerate}
\item\label{expl:pseudosmf-zero} The unit axiom \cref{enr-multifunctor-unit} and the composition axiom \cref{v-multifunctor-composition} are assumed to hold strictly.  These two axioms do not involve the symmetric group action and the pseudo symmetry isomorphisms.
\item\label{expl:pseudosmf-i} The unit permutation axiom \cref{pseudosmf-unity} and the product permutation axiom \cref{pseudosmf-product} are modeled after the symmetric group action axiom \cref{enr-multicategory-symmetry} of a $\Cat$-multicategory.  The product permutation axiom \cref{pseudosmf-product} means that, for each $n$-ary 1-cell $p \in \M\smscmap{\angc;c'}$, the diagram of $n$-ary 2-cells
\[\begin{tikzpicture}[xscale=1,yscale=1,vcenter]
\def\h{3}
\draw[0cell=.9]
(0,0) node (x11) {F(p^{\sigma\tau})}
(x11)++(\h,0) node (x12) {(Fp)^{\sigma\tau}}
(x11)++(\h/2,-1) node (x21) {(F(p^\sigma))^\tau}
;
\draw[1cell=.9]
(x11) edge node {F_{\sigma\tau;p}} (x12)
(x11) edge[shorten <=-.3ex] node[swap,pos=.4] {F_{\tau;\, p^\sigma}} (x21)
(x21) edge node[swap,pos=.6] {(F_{\sigma;p})^\tau} (x12)
;
\end{tikzpicture}\]
in $\N\scmap{F\angc\sigma\tau; Fc'}$ is commutative.
\item\label{expl:pseudosmf-iii} By the unit permutation axiom \cref{pseudosmf-unity}, the pseudo symmetry isomorphism $F_{\sigma,\angc,c'}$ is the identity natural transformation if $\angc$ is either the empty sequence or has length 1.  It follows that the restriction of a pseudo symmetric $\Cat$-multifunctor to the underlying 1-ary 2-categories is a 2-functor.  See \cref{locallysmalltwocat,def:twofunctor,ex:unarycategory}.
\item\label{expl:pseudosmf-ii} In the top and bottom equivariance axioms, \cref{pseudosmf-topeq,pseudosmf-boteq}, the four regions labeled $=$ are commutative by the composition axiom \cref{v-multifunctor-composition}, which is assumed to hold for a pseudo symmetric $\Cat$-multifunctor.  The top and bottom equivariance axioms of a pseudo symmetric $\Cat$-multifunctor are modeled after, respectively, the top and bottom equivariance axioms of a $\Cat$-multicategory, \cref{enr-operadic-eq-1,enr-operadic-eq-2}. 
\end{enumerate}
The top equivariance axiom \cref{pseudosmf-topeq} means that, for $n$-ary 1-cells 
\[p \in \M\scmap{\angcp;c''} \andspace p_j \in \M\scmap{\ang{c_j}; c_j'},\]
the following equality of $n$-ary 2-cells holds in $\N\scmap{\ang{F\ang{c_{\sigma(j)}}}_{j=1}^n; Fc''}$.  
\begin{equation}\label{Fsigmatopeq}
F_{\sigmabar;\, \ga\scmap{p; \ang{p_j}_{j=1}^n}} = 
\ga\lrscmap{F_{\sigma;p}; \bang{1_{Fp_{\sigma(j)}}}_{j=1}^n}
\end{equation}
The bottom equivariance axiom \cref{pseudosmf-boteq} means the equality of $n$-ary 2-cells
\begin{equation}\label{Fsigmaboteq}
F_{\tautimes;\, \ga\scmap{p; \ang{p_j}_{j=1}^n}} = 
\ga\lrscmap{1_{Fp}; \bang{F_{\tau_j;\, p_j}}_{j=1}^n}
\end{equation}
in $\N\scmap{\ang{F\ang{c_j}\tau_j}_{j=1}^n; Fc''}$.
\end{explanation}

\subsection*{Composition of Pseudo Symmetric $\Cat$-Multifunctors}

\begin{definition}\label{def:pseudosmfcomposition}\index{composition!pseudo symmetric $\Cat$-multifunctor}\index{pseudo symmetric!Cat-multifunctor@$\Cat$-multifunctor!composition}
Suppose given pseudo symmetric $\Cat$-multifunctors
\[\begin{tikzcd}[column sep=2.4cm]
\M \ar{r}{\left(F \scs \{F_{\sigma,\angc,c'}\}\right)} & \N \ar{r}{\left(G \scs \{G_{\sigma,\angd,d'}\}\right)} & \P. 
\end{tikzcd}\]
Their \emph{composite}\label{not:GFsigma}
\[\begin{tikzcd}[column sep=2.7cm]
\M \ar{r}{\left(GF \scs \{GF_{\sigma,\angc,c'}\}\right)} & \P 
\end{tikzcd}\]
consists of
\begin{itemize}
\item the composite function $GF \cn \Ob\M \to \Ob\P$;
\item for each $\smscmap{\angc;c'} \in \ProfMM$ the composite component functor
\[\begin{tikzcd}[column sep=large]
\M\scmap{\ang{c};c'} \ar{r}{F} & \N\scmap{F\ang{c};Fc'} \ar{r}{G} & \P\scmap{GF\angc;GFc'}
\end{tikzcd};\]
and
\item for each permutation $\sigma \in \Sigma_n$ and object $\smscmap{\angc;c'}$ as above, the natural isomorphism $GF_{\sigma,\angc,c'}$ given by the following pasting of pseudo symmetry isomorphisms.
\begin{equation}\label{pseudosymisopasting}
\begin{tikzpicture}[xscale=3,yscale=1.3,vcenter]
\draw[0cell=.85]
(0,0) node (x11) {\M\scmap{\angc;c'}}
(x11)++(1,0) node (x12) {\N\scmap{F\angc;Fc'}}
(x12)++(1.3,0) node (x13) {\P\scmap{GF\angc;GFc'}}
(x11)++(0,-1) node (x21) {\M\scmap{\angc\sigma;c'}}
(x12)++(0,-1) node (x22) {\N\scmap{F\angc\sigma;Fc'}}
(x13)++(0,-1) node (x23) {\P\scmap{GF\angc\sigma;GFc'}}
;
\draw[1cell=.9]
(x11) edge node (f1) {F} (x12)
(x12) edge node (g1) {G} (x13)
(x21) edge node[swap] (f2) {F} (x22)
(x22) edge node[swap] (g2) {G} (x23)
(x11) edge node[swap] {\sigma} (x21)
(x12) edge node {\sigma} (x22)
(x13) edge node {\sigma} (x23)
;
\draw[2cell=.9]
node[between=f1 and f2 at .6, 2label={above,F_{\sigma,\angc,c'}}, 2label={below,\iso}] {\Longrightarrow}
node[between=g1 and g2 at .6, 2label={above,G_{\sigma,F\angc,Fc'}}, 2label={below,\iso}] {\Longrightarrow}
;
\end{tikzpicture}
\end{equation} 
\end{itemize}
This finishes the definition of the composite $\big(GF, \{GF_{\sigma,\angc,c'}\}\big)$.
\end{definition}

\begin{explanation}\label{expl:GFsigma}
For each $n$-ary 1-cell $p \in \M\smscmap{\angc;c'}$, the $p$-component of $GF_\sigma$ in \cref{pseudosymisopasting} is given by the categorical composite
\[\begin{tikzpicture}[xscale=1,yscale=1,vcenter]
\def\h{3.5}
\draw[0cell=.9]
(0,0) node (x11) {GF(p^{\sigma})}
(x11)++(\h,0) node (x12) {(GFp)^{\sigma}}
(x11)++(\h/2,-1) node (x21) {G\left((Fp)^\sigma\right)}
;
\draw[1cell=.9]
(x11) edge node {(GF)_{\sigma;p}} (x12)
(x11) edge[shorten <=-.3ex] node[swap,pos=.4] {G(F_{\sigma;p})} (x21)
(x21) edge node[swap,pos=.6] {G_{\sigma;\, Fp}} (x12)
;
\end{tikzpicture}\]
in $\P\scmap{GF\angc\sigma;GFc'}$.
\end{explanation}

\begin{lemma}\label{pseudosmfcomposite}
In the context of \cref{def:pseudosmfcomposition}, the following statements hold.
\begin{enumerate}[label=(\roman*)]
\item\label{pseudosmfcomposite-i}
The composite $(GF,\{GF_{\sigma,\angc,c'}\})$ is a pseudo symmetric $\Cat$-multifunctor.
\item\label{pseudosmfcomposite-ii}
Composition of pseudo symmetric $\Cat$-multifunctors is associative and unital with respect to identity pseudo symmetric $\Cat$-multifunctors.
\end{enumerate}
\end{lemma}

\begin{proof}
For assertion \ref{pseudosmfcomposite-i}, first note that the composite $GF$ satisfies the unit axiom \cref{enr-multifunctor-unit} and the composition axiom \cref{v-multifunctor-composition} because each of $F$ and $G$ satisfies these axioms, which do not involve the symmetric group action.  The unit permutation axiom \cref{pseudosmf-unity} is true for $GF$ because the pasting \cref{pseudosymisopasting} of two identity natural transformations is an identity natural transformation.

Each of the other three axioms, \cref{pseudosmf-product,pseudosmf-topeq,pseudosmf-boteq}, for $GF$ is obtained by pasting together the respective pasting diagrams for $F$ and $G$.  For example, for the product permutation axiom \cref{pseudosmf-product} for $GF$, the left-hand side is the following pasting diagram.
\begin{equation}\label{prodpermleft}
\begin{tikzpicture}[xscale=3.2,yscale=1.3,vcenter]
\draw[0cell=.8]
(0,0) node (x11) {\M\scmap{\angc;c'}}
(x11)++(1,0) node (x12) {\N\scmap{F\angc;Fc'}}
(x12)++(1.1,0) node (x13) {\P\scmap{GF\angc;GFc'}}
(x11)++(0,-1) node (x21) {\M\scmap{\angc\sigma;c'}}
(x12)++(0,-1) node (x22) {\N\scmap{F\angc\sigma;Fc'}}
(x13)++(0,-1) node (x23) {\P\scmap{GF\angc\sigma;GFc'}}
(x21)++(0,-1) node (x31) {\M\scmap{\angc\sigma\tau;c'}}
(x22)++(0,-1) node (x32) {\N\scmap{F\angc\sigma\tau;Fc'}}
(x23)++(0,-1) node (x33) {\P\scmap{GF\angc\sigma\tau;GFc'}}
;
\draw[1cell=.85]
(x11) edge node (f1) {F} (x12)
(x12) edge node (g1) {G} (x13)
(x21) edge node (f2) {F} (x22)
(x22) edge node (g2) {G} (x23)
(x31) edge node (f3) {F} (x32)
(x32) edge node (g3) {G} (x33)
(x11) edge node[swap] {\sigma} (x21)
(x12) edge node {\sigma} (x22)
(x13) edge node {\sigma} (x23)
(x21) edge node[swap] {\tau} (x31)
(x22) edge node {\tau} (x32)
(x23) edge node {\tau} (x33)
;
\draw[2cell=.9]
node[between=f1 and f2 at .7, 2label={above,F_{\sigma}}] {\Longrightarrow}
node[between=g1 and g2 at .7, 2label={above,G_{\sigma}}] {\Longrightarrow}
node[between=f2 and f3 at .7, 2label={above,F_{\tau}}] {\Longrightarrow}
node[between=g2 and g3 at .7, 2label={above,G_{\tau}}] {\Longrightarrow}
;
\end{tikzpicture}
\end{equation}
By the product permutation axiom \cref{pseudosmf-product} for each of $F$ and $G$, the previous pasting diagram is equal to the following, which is $GF_{\sigma\tau}$.
\begin{equation}\label{prodpermright}
\begin{tikzpicture}[xscale=3.2,yscale=1.3,vcenter]
\draw[0cell=.8]
(0,0) node (x11) {\M\scmap{\angc;c'}}
(x11)++(1,0) node (x12) {\N\scmap{F\angc;Fc'}}
(x12)++(1.1,0) node (x13) {\P\scmap{GF\angc;GFc'}}
(x11)++(0,-1) node (x21) {\M\scmap{\angc\sigma\tau;c'}}
(x12)++(0,-1) node (x22) {\N\scmap{F\angc\sigma\tau;Fc'}}
(x13)++(0,-1) node (x23) {\P\scmap{GF\angc\sigma\tau;GFc'}}
;
\draw[1cell=.85]
(x11) edge node (f1) {F} (x12)
(x12) edge node (g1) {G} (x13)
(x21) edge node (f2) {F} (x22)
(x22) edge node (g2) {G} (x23)
(x11) edge node[swap] {\sigma\tau} (x21)
(x12) edge node {\sigma\tau} (x22)
(x13) edge node {\sigma\tau} (x23)
;
\draw[2cell=.9]
node[between=f1 and f2 at .7, 2label={above,F_{\sigma\tau}}] {\Longrightarrow}
node[between=g1 and g2 at .7, 2label={above,G_{\sigma\tau}}] {\Longrightarrow}
;
\end{tikzpicture}
\end{equation}
This proves the axiom \cref{pseudosmf-product} for $GF$.  The top and bottom equivariance axioms, \cref{pseudosmf-topeq,pseudosmf-boteq}, for $GF$ are verified in the same way.  This proves assertion \ref{pseudosmfcomposite-i}.

Assertion \cref{pseudosmfcomposite-ii} is true because
\begin{itemize}
\item composition of functions on objects,
\item composition of component functors, and
\item pasting of natural isomorphisms \cref{pseudosymisopasting}
\end{itemize}
are all associative and unital.
\end{proof}

\section{Pseudo Symmetric $\Cat$-Multinatural Transformations}
\label{sec:pstransformation}

In this section we discuss transformations between pseudo symmetric $\Cat$-multifunctors.
\begin{itemize}
\item These transformations are defined in \cref{def:psmftransformation}.
\item \cref{pstvcomp,psthcomp} show that identities, vertical composition, and horizontal composition of these transformations are well defined.
\item \cref{thm:catmulticatps} shows that there is a 2-category $\catmulticatps$ with small $\Cat$-multicategories as objects, pseudo symmetric $\Cat$-multifunctors as 1-cells, and pseudo symmetric $\Cat$-multinatural transformations as 2-cells.  
\end{itemize}
A pseudo symmetric $\Cat$-multifunctor (\cref{def:pseudosmultifunctor}) has not only an object assignment and multimorphism functors but also pseudo symmetry isomorphisms.  Thus a transformation between pseudo symmetric $\Cat$-multifunctors should also preserve the pseudo symmetry isomorphisms as in the next definition.

\begin{definition}\label{def:psmftransformation}
Suppose given pseudo symmetric $\Cat$-multifunctors as follows.
\[\begin{tikzcd}[column sep=2.4cm]
\M \ar[shift left]{r}{\left(F \scs \{F_{\sigma,\angc,c'}\}\right)} \ar[shift right]{r}[swap]{\left(G \scs \{G_{\sigma,\angc,c'}\}\right)} & \N
\end{tikzcd}\]
A \emph{pseudo symmetric $\Cat$-multinatural transformation}\index{pseudo symmetric!Cat-multinatural@$\Cat$-multinatural!transformation}\index{multinatural transformation!pseudo symmetric $\Cat$-} $\theta \cn F \to G$ consists of a component 1-ary 1-cell
\begin{equation}\label{pstcomponent}
\begin{tikzcd}[column sep=large]
Fc \ar{r}{\theta_c} & Gc
\end{tikzcd}
\qquad \text{in $\N$ for each $c \in \Ob\M$.}
\end{equation}
These data are required to satisfy
\begin{itemize}
\item the object $\Cat$-naturality condition \cref{catmultinaturality},
\item the morphism $\Cat$-naturality condition \cref{catmultinaturalityiicell}, and
\item the \emph{pseudo symmetry preservation axiom}\index{pseudo symmetry preservation axiom} \cref{pspreservation} below.  
\end{itemize}
For each $n$-ary 1-cell $p \in \M\smscmap{\angc;c'}$ and permutation $\sigma \in \Sigma_n$, the following equality of $n$-ary 2-cells holds in $\N\scmap{F\angc\sigma; Gc'}$.
\begin{equation}\label{pspreservation}
\ga\scmap{1_{\theta_{c'}}; F_{\sigma;p}} = 
\ga\scmap{G_{\sigma; p}; 1_{\theta_{\angc\sigma}}}
\end{equation}
Here $F_{\sigma;p}$ denotes the $p$-component of $F_\sigma$ and likewise for $G_{\sigma;p}$.  

A \emph{pseudo symmetric $\Cat$-multinatural isomorphism}\index{pseudo symmetric!Cat-multinatural@$\Cat$-multinatural!isomorphism}\index{multinatural isomorphism!pseudo symmetric $\Cat$-} is a pseudo symmetric $\Cat$-multinatural transformation with each component 1-ary 1-cell invertible with respect to $\ga$ in $\N$.

Moreover, for pseudo symmetric $\Cat$-multinatural transformations between pseudo symmetric $\Cat$-multifunctors, we define
\begin{itemize}
\item \emph{vertical composition}\index{vertical composition!pseudo symmetric $\Cat$-multinatural transformation} as in \cref{multinatvcomp},
\item \emph{horizontal composition}\index{horizontal composition!pseudo symmetric $\Cat$-multinatural transformation} as in \cref{multinathcomp}, and
\item the \emph{identity pseudo symmetric $\Cat$-multinatural transformation} $1_F$ as consisting of the component colored units \cref{thetaccomponent} 
\[(1_F)_c = 1_{Fc} \cn Fc \to Fc.\] 
\end{itemize}
Vertical composition and identities are well defined by \cref{pstvcomp} below.  Horizontal composition is well defined by \cref{psthcomp} below.  This finishes the definition.
\end{definition}

\begin{explanation}\label{expl:pstransformation}
Similar to \cref{catmultinatiicellpasting}, the pseudo symmetry preservation axiom \cref{pspreservation} says that the following two multicategorical composites in $\N$ are equal, with $(-)^\sigma$ denoting the right $\sigma$-action.
\begin{equation}\label{psprepasting}
\begin{tikzpicture}[xscale=1,yscale=1,vcenter]
\def\a{60} \def\s{.8} \def\v{-1.3}
\draw[0cell=\s]
(0,0) node (x11) {F\angc\sigma}
(x11)++(3.5,0) node (x12) {G\angc\sigma}
(x11)++(0,\v) node (x21) {Fc'}
(x12)++(0,\v) node (x22) {Gc'}
;
\draw[1cell=\s]  
(x11) edge node {\theta_{\angc\sigma}} (x12)
(x21) edge node {\theta_{c'}} (x22)
(x11) edge[bend right=\a] node[swap,pos=.42] {F(p^\sigma)} (x21)
(x11) edge[bend left=\a] node[pos=.42] {(Fp)^\sigma} (x21)
(x12) edge[bend right=\a] node[swap,pos=.42] {G(p^\sigma)} (x22)
(x12) edge[bend left=\a] node[pos=.42] {(Gp)^\sigma} (x22)
;
\draw[2cell=\s]
node[between=x11 and x21 at .65, 2label={above,F_{\sigma;p}}] {\Longrightarrow}
node[between=x12 and x22 at .65, 2label={above,G_{\sigma;p}}] {\Longrightarrow}
;
\end{tikzpicture}
\end{equation}
The domain $n$-ary 1-cells of the two $n$-ary 2-cell composites in \cref{psprepasting}, 
\begin{equation}\label{pstdomain}
\ga\scmap{G(p^\sigma); \theta_{\angc\sigma}} = \ga\scmap{\theta_{c'}; F(p^\sigma)},
\end{equation}
are equal by the object $\Cat$-naturality condition \cref{catmultinaturality} applied to $p^\sigma$.  The codomain $n$-ary 1-cells of the two $n$-ary 2-cell composites in \cref{psprepasting} are equal by the following equalities.
\begin{equation}\label{pstcodomain}
\left\{\begin{aligned}
\ga\scmap{(Gp)^\sigma; \theta_{\angc\sigma}} 
&= \ga\scmap{Gp; \theta_{\angc}} \cdot \sigma && \text{by \cref{enr-operadic-eq-1}}\\
&= \ga\scmap{\theta_{c'}; Fp} \cdot \sigma && \text{by \cref{catmultinaturality}}\\
&= \ga\scmap{\theta_{c'}; (Fp)^\sigma} && \text{by \cref{enr-operadic-eq-2}}
\end{aligned}\right.
\end{equation}
A pseudo symmetric $\Cat$-multinatural transformation has the same definition as a $\Cat$-multinatural transformation (\cref{expl:catmultitransformation}) with one additional axiom \cref{pspreservation} but \emph{no} extra data.
\end{explanation}

\begin{example}[$\Cat$-Multinatural Transformations]\label{ex:strictpst}\index{multinatural transformation}
Suppose $\theta \cn F \to G$ is a $\Cat$-multinatural transformation (\cref{expl:catmultitransformation}) between $\Cat$-multifunctors $F,G \cn \M \to \N$.  We regard each of $F$ and $G$ as a pseudo symmetric $\Cat$-multifunctor with identity pseudo symmetry isomorphisms (\cref{ex:idpsm}).  Then $\theta$ is a pseudo symmetric $\Cat$-multinatural transformation.  Its pseudo symmetry preservation axiom \cref{pspreservation} follows from the object $\Cat$-naturality condition \cref{catmultinaturality} as in \cref{pstdomain} and the functoriality of $\ga$.
\end{example}

\begin{example}[Non-Symmetric $\Cat$-Multinatural Transformations]\label{ex:nscattransformation}\index{multinatural transformation!non-symmetric}\index{non-symmetric!multinatural transformation}
Suppose $\theta \cn F \to G$ is a pseudo symmetric $\Cat$-multinatural transformation between pseudo symmetric $\Cat$-multifunctors.  Then $\theta$ yields a non-symmetric $\Cat$-multinatural transformation between the underlying non-symmetric $\Cat$-multifunctors of $F$ and $G$ (\cref{ex:nscatmultifunctors}), since $\theta$ satisfies the $\Cat$-naturality conditions \cref{catmultinaturality,catmultinaturalityiicell}.
\end{example}

\begin{lemma}\label{pstvcomp}
In \cref{def:psmftransformation},
\begin{itemize}
\item identity pseudo symmetric $\Cat$-multinatural transformations and
\item vertical composites of pseudo symmetric $\Cat$-multinatural transformations
\end{itemize}
are well defined.
\end{lemma}

\begin{proof}
The $\Cat$-naturality conditions \cref{catmultinaturality,catmultinaturalityiicell} do not involve the symmetric group action.  They are satisfied by the identities and vertical composites as part of the 2-category structure on $\catmulticatns$ (\cref{v-multicat-2cat}).

It remains to verify the pseudo symmetry preservation axiom \cref{pspreservation}.  For the identity $1_F$, the axiom \cref{pspreservation} holds by the right and left unity axioms, \cref{enr-multicategory-right-unity,enr-multicategory-left-unity}, since each component of $1_F$ is a colored unit.

Consider pseudo symmetric $\Cat$-multinatural transformations
\[\begin{tikzcd}[column sep=large]
F \ar{r}{\theta} & G \ar{r}{\psi} & H.
\end{tikzcd}\]
The axiom \cref{pspreservation} for the vertical composite $\psi\theta$ \cref{multinatvcomp} is obtained by pasting \cref{psprepasting} and the one for $\psi$ together as follows.
\[\begin{tikzpicture}[xscale=1,yscale=1,vcenter]
\def\a{60} \def\s{.8} \def\v{-1.3}
\draw[0cell=\s]
(0,0) node (x11) {F\angc\sigma}
(x11)++(3,0) node (x12) {G\angc\sigma}
(x12)++(3,0) node (x13) {H\angc\sigma}
(x11)++(0,\v) node (x21) {Fc'}
(x12)++(0,\v) node (x22) {Gc'}
(x13)++(0,\v) node (x23) {Hc'}
;
\draw[1cell=\s]  
(x11) edge node {\theta_{\angc\sigma}} (x12)
(x12) edge node {\psi_{\angc\sigma}} (x13)
(x21) edge node {\theta_{c'}} (x22)
(x22) edge node {\psi_{c'}} (x23)
;
\draw[1cell=.7]
(x11) edge[bend right=\a] node[swap,pos=.42] {F(p^\sigma)} (x21)
(x11) edge[bend left=\a] node[pos=.42] {(Fp)^\sigma} (x21)
(x12) edge[bend right=\a] node[swap,pos=.42] {G(p^\sigma)} (x22)
(x12) edge[bend left=\a] node[pos=.42] {(Gp)^\sigma} (x22)
(x13) edge[bend right=\a] node[swap,pos=.42] {H(p^\sigma)} (x23)
(x13) edge[bend left=\a] node[pos=.42] {(Hp)^\sigma} (x23)
;
\draw[2cell=\s]
node[between=x11 and x21 at .65, 2label={above,F_{\sigma;p}}] {\Longrightarrow}
node[between=x12 and x22 at .65, 2label={above,G_{\sigma;p}}] {\Longrightarrow}
node[between=x13 and x23 at .65, 2label={above,H_{\sigma;p}}] {\Longrightarrow}
;
\end{tikzpicture}\]
This is valid by the associativity axiom \cref{enr-multicategory-associativity} in $\N$ and the pseudo symmetry preservation axiom \cref{pspreservation} for $\theta$ and $\psi$.
\end{proof}

\begin{lemma}\label{psthcomp}
In \cref{def:psmftransformation}, horizontal composites of pseudo symmetric $\Cat$-multinatural transformations are well defined.
\end{lemma}

\begin{proof}
Consider the horizontal composite $\theta' * \theta$ \cref{multinathcomp} of the following pseudo symmetric $\Cat$-multinatural transformations.
\[\begin{tikzpicture}[xscale=1,yscale=1,vcenter]
\def\a{25}
\draw[0cell=.9]
(0,0) node (x11) {\M}
(x11)++(2.5,0) node (x12) {\N}
(x12)++(2.5,0) node (x13) {\P}
;
\draw[1cell=.8] 
(x11) edge[bend left=\a] node {(F, \{F_\sigma\})} (x12)
(x11) edge[bend right=\a] node[swap] {(G, \{G_\sigma\})} (x12)
(x12) edge[bend left=\a] node {(F', \{F'_\sigma\})} (x13)
(x12) edge[bend right=\a] node[swap] {(G', \{G'_\sigma\})} (x13)
;
\draw[2cell=1]
node[between=x11 and x12 at .45, rotate=-90, 2label={above,\theta}] {\Rightarrow}
node[between=x12 and x13 at .45, rotate=-90, 2label={above,\theta'}] {\Rightarrow}
;
\end{tikzpicture}\]
The $\Cat$-naturality conditions \cref{catmultinaturality,catmultinaturalityiicell} do not involve the symmetric group action and are satisfied by the horizontal composite as part of the 2-category structure on $\catmulticatns$ (\cref{v-multicat-2cat}).  It remains to verify the pseudo symmetry preservation axiom \cref{pspreservation}, which is the following equality in $\P\scmap{F'F\angc\sigma; G'Gc'}$ for each $n$-ary 1-cell $p \in \M\smscmap{\angc;c'}$ and permutation $\sigma \in \Sigma_n$.
\begin{equation}\label{pspreaxiomhcomp}
\ga\lrscmap{1_{(\theta' * \theta)_{c'}}; (F'F)_{\sigma;p}} 
= \ga\lrscmap{(G'G)_{\sigma;p}; 1_{(\theta' * \theta)_{\angc\sigma}}}
\end{equation}

To prove \cref{pspreaxiomhcomp}, we use
\begin{itemize}
\item the second expression in \cref{hcompexpressions} for the horizontal composite
\[(\theta' * \theta)_{c'} = \ga\scmap{G'\theta_{c'}; \theta'_{Fc'}}\]
and likewise for $(\theta' * \theta)_{c_{\sigma(j)}}$ for $j \in \{1,\ldots,n\}$;
\item the definition \cref{pseudosymisopasting} of the pseudo symmetry isomorphism
\[(F'F)_{\sigma;p} = F'_{\sigma;Fp} \circ F'(F_{\sigma;p})\] 
and likewise for $(G'G)_{\sigma;p}$; and
\item the functoriality and associativity \cref{enr-multicategory-associativity} of $\gamma$ in $\P$.
\end{itemize} 
Then the desired equality \cref{pspreaxiomhcomp} becomes the following equality.
\begin{equation}\label{pspreaxiomhcomp-ii}
\begin{split}
&\scalebox{.8}{$\ga\lrscmap{1_{G'\theta_{c'}}; \ga\lrscmap{1_{\theta'_{Fc'}}; F'_{\sigma; Fp}}} 
\circ \ga\lrscmap{1_{G'\theta_{c'}}; \ga\lrscmap{1_{\theta'_{Fc'}}; F'(F_{\sigma;p})}}$}\\
&= \scalebox{.8}{$\ga\lrscmap{\ga\lrscmap{G'_{\sigma;Gp}; 1_{G'\theta_{\angc\sigma}}} ; 1_{\theta'_{F\angc\sigma}}}
\circ \ga\lrscmap{\ga\lrscmap{G'(G_{\sigma;p}); 1_{G'\theta_{\angc\sigma}}} ; 1_{\theta'_{F\angc\sigma}}}$}
\end{split}
\end{equation}
We prove \cref{pspreaxiomhcomp-ii} by showing that the left, respectively, right, factors on the two sides are equal using
\begin{itemize}
\item the $\Cat$-naturality conditions \cref{catmultinaturality,catmultinaturalityiicell} for $\theta$ and $\theta'$,
\item the associativity \cref{enr-multicategory-associativity} of $\ga$, and
\item the fact that $G$', as a pseudo symmetric $\Cat$-multifunctor, strictly preserves the multicategorical composition \cref{v-multifunctor-composition}.
\end{itemize} 

The equality of the right factors in \cref{pspreaxiomhcomp-ii} are proved as follows.
\[\begin{aligned}
&\ga\lrscmap{1_{G'\theta_{c'}}; \ga\lrscmap{1_{\theta'_{Fc'}}; F'(F_{\sigma;p})}} &&\\
&= \ga\lrscmap{1_{G'\theta_{c'}}; \ga\lrscmap{G'(F_{\sigma;p}); 1_{\theta'_{F\angc\sigma}}}} && \text{by \cref{catmultinaturalityiicell} for $\theta'$}\\
&= \ga\lrscmap{\ga\lrscmap{1_{G'\theta_{c'}}; G'(F_{\sigma;p})}; 1_{\theta'_{F\angc\sigma}}} && \text{by \cref{enr-multicategory-associativity}}\\
&= \ga\lrscmap{G' \ga\lrscmap{1_{\theta_{c'}}; F_{\sigma;p}}; 1_{\theta'_{F\angc\sigma}}} && \text{by \cref{v-multifunctor-composition}}\\
&= \ga\lrscmap{G' \ga\lrscmap{G_{\sigma;p}; 1_{\theta_{\angc\sigma}}} ; 1_{\theta'_{F\angc\sigma}}} &&  \text{by \cref{pspreservation}}\\
&= \ga\lrscmap{\ga\lrscmap{G'(G_{\sigma;p}); 1_{G'\theta_{\angc\sigma}}} ; 1_{\theta'_{F\angc\sigma}}} &&  \text{by \cref{v-multifunctor-composition}}
\end{aligned}\]
The equality of the left factors in \cref{pspreaxiomhcomp-ii} are proved as follows.
\[\begin{aligned}
&\ga\lrscmap{1_{G'\theta_{c'}}; \ga\lrscmap{1_{\theta'_{Fc'}}; F'_{\sigma; Fp}}} &&\\
&= \ga\lrscmap{1_{G'\theta_{c'}}; \ga\lrscmap{G'_{\sigma; Fp}; 1_{\theta'_{F\angc\sigma}}}} && \text{by \cref{pspreservation} for $\theta'$}\\
&= \ga\lrscmap{\ga\lrscmap{1_{G'\theta_{c'}}; G'_{\sigma; Fp}} ; 1_{\theta'_{F\angc\sigma}}} && \text{by \cref{enr-multicategory-associativity}}\\
&= \ga\lrscmap{\ga\lrscmap{G'_{\sigma; Gp}; 1_{G'\theta_{\angc\sigma}}} ; 1_{\theta'_{F\angc\sigma}}} && \text{by \cref{topboteq} below}
\end{aligned}\]
The last equality follows from the following computation, which uses the top and bottom equivariance axioms, \cref{Fsigmatopeq,Fsigmaboteq}, for $G'$.
\begin{equation}\label{topboteq}
\left\{\begin{aligned}
&\ga\lrscmap{1_{G'\theta_{c'}} ; G'_{\sigma;\, Fp}} &&\\
&= G'_{\sigma;\, \ga\smscmap{\theta_{c'}; Fp}} && \text{by \cref{Fsigmaboteq} for $(\theta_{c'}, Fp)$}\\
&= G'_{\sigma;\, \ga\smscmap{Gp; \theta_{\angc}}} && \text{by \cref{catmultinaturality}}\\
&= \ga\lrscmap{G'_{\sigma;\, Gp}; 1_{G'\theta_{\angc\sigma}}} && \text{by \cref{Fsigmatopeq} for $(Gp, \theta_{\angc})$}
\end{aligned}\right.
\end{equation}
This finishes the proof of the desired equality \cref{pspreaxiomhcomp-ii}.
\end{proof}

\begin{explanation}\label{expl:topboteq}
The computation \cref{topboteq} asserts that the following two composite $n$-ary 2-cells in $\P$ are equal.
\[\begin{tikzpicture}[xscale=1,yscale=1,vcenter]
\def\a{60} \def\s{.8} \def\v{-1.3}
\draw[0cell=\s]
(0,0) node (x11) {G'F\angc\sigma}
(x11)++(4.5,0) node (x12) {G'G\angc\sigma}
(x11)++(0,\v) node (x21) {G'Fc'}
(x12)++(0,\v) node (x22) {G'Gc'}
;
\draw[1cell=\s]  
(x11) edge node {G'\theta_{\angc\sigma}} (x12)
(x21) edge node {G'\theta_{c'}} (x22)
(x11) edge[bend right=\a] node[swap,pos=.45] {G'(Fp)^\sigma} (x21)
(x11) edge[bend left=\a] node[pos=.45] {(G'Fp)^\sigma} (x21)
(x12) edge[bend right=\a] node[swap,pos=.45] {G'(Gp)^\sigma} (x22)
(x12) edge[bend left=\a] node[pos=.45] {(G'Gp)^\sigma} (x22)
;
\draw[2cell=\s]
node[between=x11 and x21 at .65, 2label={above,G'_{\sigma;\, Fp}}] {\Longrightarrow}
node[between=x12 and x22 at .65, 2label={above,G'_{\sigma;\, Gp}}] {\Longrightarrow}
;
\end{tikzpicture}\]
The top and bottom equivariance axioms, \cref{Fsigmatopeq,Fsigmaboteq}, of a pseudo symmetric $\Cat$-multifunctor are necessary to prove this equality.
\end{explanation}

\begin{theorem}\label{thm:catmulticatps}
There is a 2-category\index{2-category!of pseudo symmetric $\Cat$-multifunctors}\index{pseudo symmetric!Cat-multifunctor@$\Cat$-multifunctor!2-category} $\catmulticatps$ consisting of the following data.
\begin{itemize}
\item Its objects are small $\Cat$-multicategories.
\item For small $\Cat$-multicategories $\M$ and $\N$, the hom category 
\[\catmulticatps(\M,\N)\]
has
\begin{itemize}
\item pseudo symmetric $\Cat$-multifunctors $\M\to\N$ (\cref{def:pseudosmultifunctor}) as 1-cells and
\item pseudo symmetric $\Cat$-multinatural transformations (\cref{def:psmftransformation}) as 2-cells.
\end{itemize}
\item For pseudo symmetric $\Cat$-multinatural transformations, vertical composition, horizontal composition, and identity 2-cells are as in \cref{def:psmftransformation}.
\item The identity 1-cell $1_{\M}$ is the identity pseudo symmetric $\Cat$-multifunctor $1_{\M}$ (\cref{def:pseudosmultifunctor}).
\item Horizontal composition of 1-cells is the composition of pseudo symmetric $\Cat$-multifunctors (\cref{def:pseudosmfcomposition}).
\end{itemize}
\end{theorem}

\begin{proof}
Consider axioms \cref{twocat-i,twocat-ii,twocat-iii,twocat-iv} in \cref{def:twocategory} of a 2-category for $\catmulticatps$.
\begin{itemize}
\item By \cref{pseudosmfcomposite} the horizontal composition of 1-cells is well defined, and the 1-cell axiom \ref{twocat-iii} holds.
\item Vertical and horizontal compositions of 2-cells and identity 2-cells are well defined \cref{pstvcomp,psthcomp}.
\item The 2-cell axioms \cref{twocat-i,twocat-ii,twocat-iv} hold as they do in the 2-category $\catmulticatns$ (\cref{v-multicat-2cat}) because the symmetric group action is not involved.  
\end{itemize}  
This finishes the proof.
\end{proof}

\begin{explanation}[Lax Symmetric $\Cat$-Multifunctors]\label{expl:laxsmf}
The discussion in \cref{sec:pseudosmultifunctor} and this section does \emph{not} depend on the invertibility of the pseudo symmetry isomorphisms $F_\sigma$ \cref{Fsigmaangccp}.
\begin{itemize}
\item \cref{def:pseudosmultifunctor} of a pseudo symmetric $\Cat$-multifunctor still makes sense if each pseudo symmetry isomorphism $F_\sigma$ in \cref{Fsigmaangccp} is only required to be a natural \emph{transformation}, instead of a natural isomorphism.  We call this lax variant a \emph{lax symmetric $\Cat$-multifunctor}\index{lax symmetric $\Cat$-multifunctor} and $F_\sigma$ a \index{lax symmetry}\emph{lax symmetry}.
\item Composites of lax symmetric $\Cat$-multifunctors are well defined, with the composite lax symmetry given by the pasting \cref{pseudosymisopasting}.
\item The lax variants of \cref{pseudosmfcomposite,def:psmftransformation,pstvcomp,psthcomp,thm:catmulticatps} also hold.
\end{itemize}  
Therefore, there is a 2-category\label{not:catmulticatlax}\index{2-category!of lax symmetric $\Cat$-multifunctors}\index{lax symmetric $\Cat$-multifunctor!2-category}
\[\catmulticatlax\]
with
\begin{itemize}
\item small $\Cat$-multicategories as objects,
\item lax symmetric $\Cat$-multifunctors as 1-cells,
\item lax symmetric $\Cat$-multinatural transformations as 2-cells (as in \cref{def:psmftransformation}), and
\item identity 1- and 2-cells, vertical composition, and horizontal composition as in \cref{pseudosmfcomposite,pstvcomp,psthcomp}.
\end{itemize}
It is related to the 2-categories in \cref{v-multicat-2cat,thm:catmulticatps} via the following 2-functors. 
\begin{equation}\label{catmulticatpslaxns}
\begin{tikzpicture}[xscale=1,yscale=1,vcenter]
\def\s{.9} \def\v{-1} \def\h{.3}
\draw[0cell=\s]
(0,0) node (x11) {\catmulticat}
(x11)++(4,0) node (x12) {\catmulticatns}
(x11)++(\h,\v) node (x21) {\catmulticatps}
(x12)++(-\h,\v) node (x22) {\catmulticatlax}
;
\draw[1cell=\s]  
(x11) edge[right hook->] node[pos=.9] {i_1} (x21)
(x21) edge[right hook->] node {i_2} (x22)
(x22) edge node[pos=.1] {j} (x12)
;
\end{tikzpicture}
\end{equation}
Each of these 2-functors is the identity on
\begin{itemize}
\item objects, which are small $\Cat$-multicategories, and
\item 2-cells (\cref{ex:strictpst,ex:nscattransformation}).
\end{itemize} 
On 1-cells they are defined as follows.
\begin{itemize}
\item The assignment $i_1$ sends a $\Cat$-multifunctor $F$ to the pseudo symmetric $\Cat$-multifunctor $(F,\{\Id\})$ (\cref{ex:idpsm}).
\item The assignment $i_2$ is the inclusion of pseudo symmetric $\Cat$-multifunctors into lax symmetric $\Cat$-multifunctors.
\item The assignment $j$ sends a lax symmetric $\Cat$-multifunctor $(F,\{F_\sigma\})$ to its underlying non-symmetric $\Cat$-multifunctor $F$ (\cref{ex:nscatmultifunctors}).
\end{itemize}
See \cref{qu:multibicat,qu:pseudo} for related discussion.
\end{explanation}

\section{Pseudo Symmetric $\Cat$-Multiequivalences}
\label{sec:pstcharacterization}

In this section we show that the characterization of $\Cat$-multiequivalences in \cref{thm:multiwhitehead} also holds for the 2-category $\catmulticatps$ in \cref{thm:catmulticatps}.  See \cref{thm:psmultiwhitehead}.  In other words, a pseudo symmetric $\Cat$-multifunctor is a pseudo symmetric $\Cat$-multiequivalence if and only if it is essentially surjective on objects and an isomorphism on each multimorphism category.  For the \emph{if} direction, most of the work involves defining the pseudo symmetry isomorphisms for the adjoint and showing that the (co)unit preserves them.

\subsection*{Properties of Pseudo Symmetric $\Cat$-Multiequivalences}

The following notions are restrictions of those in \cref{def:equivalences} to the 2-category $\catmulticatps$ and analogs of those in \cref{def:catmultiequivalence}.

\begin{definition}\label{def:pscatmultiequivalence}
We define the following.
\begin{itemize}
\item A \emph{pseudo symmetric $\Cat$-multiequivalence}\index{pseudo symmetric!Cat-multiequivalence@$\Cat$-multiequivalence}\index{multiequivalence!pseudo symmetric $\Cat$-} is an equivalence in the 2-category $\catmulticatps$.
\item An \emph{adjoint pseudo symmetric $\Cat$-multiequivalence}\index{adjoint pseudo symmetric $\Cat$-multiequivalence}\index{pseudo symmetric!adjoint - $\Cat$-multiequivalence} is an adjoint equivalence in the 2-category $\catmulticatps$.\defmark
\end{itemize}
\end{definition}

\begin{explanation}\label{expl:pscatmultiequivalence}
In \cref{def:pscatmultiequivalence} the use of the 2-category $\catmulticatps$ is simply a matter of conceptual convenience.  Indeed, a pseudo symmetric $\Cat$-multifunctor (\cref{def:pseudosmultifunctor})
\[\begin{tikzcd}[column sep=huge]
\M \ar{r}{(F,\{F_\sigma\})} & \N
\end{tikzcd}\]
is a \emph{pseudo symmetric $\Cat$-multiequivalence} if there exist
\begin{itemize}
\item a pseudo symmetric $\Cat$-multifunctor $(G,\{G_\sigma\}) \cn \N \to \M$ and
\item pseudo symmetric $\Cat$-multinatural isomorphisms (\cref{def:psmftransformation})
\begin{equation}\label{pscatmultiequnit}
\begin{tikzcd}[column sep=large]
1_{\M} \ar{r}{\eta}[swap]{\iso} & GF
\end{tikzcd} \andspace 
\begin{tikzcd}[column sep=large]
FG \ar{r}{\epz}[swap]{\iso} & 1_{\N}.
\end{tikzcd}
\end{equation}
\end{itemize}
The quadruple $(F,G,\eta,\epz)$ is an \emph{adjoint pseudo symmetric $\Cat$-multiequivalence} if, moreover, the two triangle identities \cref{triangleidentities} are satisfied.
\end{explanation}

\begin{example}[$\Cat$-Multiequivalences]\label{ex:pscatmeq}
By the 2-functor 
\[\begin{tikzcd}[column sep=large]
\catmulticat \ar[hookrightarrow]{r}{i_1} & \catmulticatps
\end{tikzcd}\]
in \cref{catmulticatpslaxns}, a $\Cat$-multiequivalence $F \cn \M \to \N$ yields a pseudo symmetric $\Cat$-multiequivalence $(F,\{\Id\})$ when it is equipped with identity pseudo symmetry isomorphisms.  The adjoint variant is also true by the 2-functoriality of $i_1$.
\end{example}

\begin{example}[Non-Symmetric $\Cat$-Multiequivalences]\label{ex:nscatmeq}
By the 2-functor
\[\begin{tikzcd}[column sep=large]
\catmulticatps \ar{r}{j \circ i_2} & \catmulticatns
\end{tikzcd}\]
in \cref{catmulticatpslaxns}, a pseudo symmetric $\Cat$-multiequivalence $(F,\{F_\sigma\})$ yields a non-symmetric $\Cat$-multiequivalence $F$ by forgetting the pseudo symmetry isomorphisms.  The adjoint variant is also true by the 2-functoriality of $j \circ i_2$.
\end{example}

\begin{lemma}\label{pscatmultieqesff}\index{pseudo symmetric!Cat-multiequivalence@$\Cat$-multiequivalence!properties}
Each pseudo symmetric $\Cat$-multiequivalence is 
\begin{itemize}
\item essentially surjective on objects and
\item an isomorphism on each multimorphism category.
\end{itemize}
\end{lemma}

\begin{proof}
By \cref{ex:nscatmeq} each pseudo symmetric $\Cat$-multiequivalence $(F,\{F_\sigma\})$ yields a non-symmetric $\Cat$-multiequivalence $F$, which has the two stated properties by \cref{catmultieqesff}.
\end{proof}

\subsection*{Adjoint, Unit, and Counit}

The main part of \cref{thm:psmultiwhitehead} is the converse of \cref{pscatmultieqesff}.  We proceed as in \cref{sec:catmultieqtheorem}.  Given a pseudo symmetric $\Cat$-multifunctor that satisfies the two properties in \cref{pscatmultieqesff}, we define its adjoint, unit, and counit as in the non-symmetric case as follows.

\begin{definition}\label{def:psFetaepz}
Suppose given a pseudo symmetric $\Cat$-multifunctor (\cref{def:pseudosmultifunctor})
\[\begin{tikzcd}[column sep=huge]
\N \ar{r}{(G,\{G_\sigma\})} & \M
\end{tikzcd}\]
that is
\begin{itemize}
\item essentially surjective on objects and
\item an isomorphism on each multimorphism category.
\end{itemize}
Define 
\begin{itemize}
\item a non-symmetric $\Cat$-multifunctor $F \cn \M \to \N$ as in \cref{def:Fandeta};
\item a non-symmetric $\Cat$-multinatural isomorphism $\eta \cn 1_{\M} \fto{\iso} GF$ as in \cref{etaxxgfx}; and
\item a non-symmetric $\Cat$-multinatural isomorphism $\epz \cn FG \fto{\iso} 1_{\N}$ as in \cref{def:counitepz}.
\end{itemize}
By \cref{Fcatmultifunctor,etacatmultinatural,epzcatmultiiso}, respectively, $F$, $\eta$, and $\epz$ are well defined.  Moreover, by \cref{multieqtriangleids} the quadruple $(F,G,\eta,\epz)$ satisfies the triangle identities \cref{triangleidentities}.
\end{definition}

In the context of \cref{def:psFetaepz}, to finish the proof of \cref{thm:psmultiwhitehead}, we proceed as follows.
\begin{enumerate}[label=(\roman*)]
\item Define the pseudo symmetry isomorphisms $F_\sigma$ for $F$ (\cref{def:Fsigmapswhitehead}).
\item Show that $(F,\{F_\sigma\}) \cn \M \to \N$ is a pseudo symmetric $\Cat$-multifunctor (\cref{FFsigmaps}).
\item Show that each of $\eta$ and $\epz$ in \cref{def:psFetaepz} is a pseudo symmetric $\Cat$-multinatural isomorphism (\cref{etapscatmt,epzpscatmt}).
\item Conclude that the quadruple $(F,G,\eta,\epz)$ is an adjoint pseudo symmetric $\Cat$-multiequivalence (\cref{thm:psmultiwhitehead}).
\end{enumerate}
The following preliminary observation is needed to define $F_\sigma$.  Recall that $(-)^\sigma$ denotes the right $\sigma$-action.

\begin{lemma}\label{GFpsigma}
In the context of \cref{def:psFetaepz}, suppose $\sigma \in \Sigma_n$ is a permutation and 
\[p \in \M\smscmap{\angc;c'} \andspace q \in \N\smscmap{\angd;d'}\]
are $n$-ary 1-cells.  Then the following equalities hold.
\begin{equation}\label{GFpFGq}
\left\{\begin{split}
GF(p^\sigma) &= (GFp)^\sigma \in \M\scmap{GF\angc\sigma; GFc'}\\
FG(q^\sigma) &= (FGq)^\sigma \in \N\scmap{FG\angd\sigma; FGd'}
\end{split}\right.
\end{equation}
\end{lemma}

\begin{proof}
To prove the first equality in \cref{GFpFGq}, we compute as follows using the object $\Cat$-naturality of $\eta$ \cref{etaobjcatnaturality} and the top and bottom equivariance axioms, \cref{enr-operadic-eq-1,enr-operadic-eq-2}, in $\M$.
\begin{equation}\label{gaFGpsigma}
\left\{\begin{aligned}
& \ga\scmap{GF(p^\sigma); \eta_{\angc\sigma}} &&\\
&= \ga\scmap{\eta_{c'}; p^\sigma} && \text{by \cref{etaobjcatnaturality}}\\
&= \ga\scmap{\eta_{c'}; p} \cdot \sigma && \text{by \cref{enr-operadic-eq-2} in $\M$}\\
&= \ga\scmap{GFp; \eta_{\angc}} \cdot \sigma && \text{by \cref{etaobjcatnaturality}}\\
&= \ga\scmap{(GFp)^\sigma; \eta_{\angc\sigma}} && \text{by \cref{enr-operadic-eq-1} in $\M$}
\end{aligned}\right.
\end{equation}
Since each component 1-ary 1-cell $\eta_{c_j}$ is invertible by \cref{etaxxgfx}, we obtain the first equality in \cref{GFpFGq} by
\begin{itemize}
\item applying $\ga\scmap{-; \eta^\inv_{\angc\sigma}}$ to \cref{gaFGpsigma} and
\item using associativity \cref{enr-multicategory-associativity} and right unity \cref{enr-multicategory-right-unity} in $\M$.
\end{itemize}

The second equality in \cref{GFpFGq} is proved similarly by first obtaining the following equalities using the object $\Cat$-naturality of $\epz$ \cref{Ggafepzxj} and the top and bottom equivariance axioms in $\N$.
\begin{equation}\label{gaetadpFGqsigma}
\left\{\begin{aligned}
& \ga\scmap{\epz_{d'}; FG(q^\sigma)} &&\\
&= \ga\scmap{q^\sigma; \epz_{\angd\sigma}} && \text{by \cref{Ggafepzxj}}\\
&= \ga\scmap{q; \epz_{\angd}} \cdot \sigma && \text{by \cref{enr-operadic-eq-1} in $\N$}\\
&= \ga\scmap{\epz_{d'}; FGq} \cdot \sigma && \text{by \cref{Ggafepzxj}}\\
&= \ga\scmap{\epz_{d'}; (FGq)^\sigma} && \text{by \cref{enr-operadic-eq-2} in $\N$}
\end{aligned}\right.
\end{equation}
Since the component 1-ary 1-cell $\epz_{d'}$ is invertible by \cref{epzzcomponent},  we obtain the second equality in \cref{GFpFGq} by
\begin{itemize}
\item applying $\ga\scmap{\epz^\inv_{d'}; -}$ to \cref{gaetadpFGqsigma} and
\item using associativity \cref{enr-multicategory-associativity} and left unity \cref{enr-multicategory-left-unity} in $\N$.
\end{itemize}
This finishes the proof.
\end{proof}

\begin{definition}\label{def:Fsigmapswhitehead}
In the context of \cref{def:psFetaepz,GFpsigma}, define the data of a natural isomorphism 
\[\begin{tikzpicture}[xscale=3.5,yscale=1.3,vcenter]
\draw[0cell=.9]
(0,0) node (x11) {\M\scmap{\angc;c'}}
(x11)++(1,0) node (x12) {\N\scmap{F\angc;Fc'}}
(x11)++(0,-1) node (x21) {\M\scmap{\angc\sigma;c'}}
(x12)++(0,-1) node (x22) {\N\scmap{F\angc\sigma;Fc'}}
;
\draw[1cell=.9]
(x11) edge node (f1) {F} (x12)
(x21) edge node[swap] (f2) {F} (x22)
(x11) edge node[swap] {\sigma} (x21)
(x12) edge node {\sigma} (x22)
;
\draw[2cell=.9]
node[between=f1 and f2 at .55, 2label={above,F_{\sigma}}, 2label={below,\iso}] {\Longrightarrow}
;
\end{tikzpicture}\]
as follows.  For each $n$-ary 1-cell $p \in \M\smscmap{\angc;c'}$, the $(Fp)$-component of the pseudo symmetry isomorphism $G_\sigma$ has an inverse $n$-ary 2-cell
\[\begin{tikzcd}[column sep=large]
GF(p^\sigma) = (GFp)^\sigma \ar{r}{G^\inv_{\sigma; Fp}}[swap]{\iso} & G(Fp)^\sigma
\end{tikzcd} \inspace \M\scmap{GF\angc\sigma; GFc'}\]
with the equality from \cref{GFpFGq}.  Since the multimorphism functor
\begin{equation}\label{NFcsigmaFcp}
\begin{tikzcd}[column sep=large]
\N\scmap{F\angc\sigma; Fc'} \ar{r}{G}[swap]{\iso} & \M\scmap{GF\angc\sigma; GFc'}
\end{tikzcd}
\end{equation}
is an isomorphism, there exists a unique $n$-ary 2-cell
\begin{equation}\label{Fsigmappswhitehead}
\begin{tikzcd}[column sep=large]
F(p^\sigma) \ar{r}{F_{\sigma;p}}[swap]{\iso} & (Fp)^\sigma
\end{tikzcd} \inspace \N\scmap{F\angc\sigma; Fc'}
\end{equation}
such that
\begin{equation}\label{GFsigmapGinvsigmaFp}
G(F_{\sigma;p}) = G^\inv_{\sigma; Fp}.
\end{equation}
We define $F_{\sigma;p}$ as the $p$-component of $F_{\sigma}$.  Moreover, $F_{\sigma;p}$ is an invertible $n$-ary 2-cell by the fact that $G$ in \cref{NFcsigmaFcp} is an isomorphism.  The naturality of $F_\sigma$ follows from the naturality of $G_\sigma$ and the isomorphism $G$ \cref{NFcsigmaFcp}.
\end{definition}

\begin{explanation}\label{expl:Fsigmappswhitehead}
The definition \cref{GFsigmapGinvsigmaFp} of $F_\sigma$ means that the pasting
\begin{equation}\label{FsigmaGsigma}
\begin{tikzpicture}[xscale=3,yscale=1.3,vcenter]
\draw[0cell=.85]
(0,0) node (x11) {\M\scmap{\angc;c'}}
(x11)++(1,0) node (x12) {\N\scmap{F\angc;Fc'}}
(x12)++(1.2,0) node (x13) {\M\scmap{GF\angc;GFc'}}
(x11)++(0,-1) node (x21) {\M\scmap{\angc\sigma;c'}}
(x12)++(0,-1) node (x22) {\N\scmap{F\angc\sigma;Fc'}}
(x13)++(0,-1) node (x23) {\M\scmap{GF\angc\sigma;GFc'}}
;
\draw[1cell=.9]
(x11) edge node (f1) {F} (x12)
(x12) edge node (g1) {G} (x13)
(x21) edge node[swap] (f2) {F} (x22)
(x22) edge node[swap] (g2) {G} (x23)
(x11) edge node[swap] {\sigma} (x21)
(x12) edge node {\sigma} (x22)
(x13) edge node {\sigma} (x23)
;
\draw[2cell=.9]
node[between=f1 and f2 at .55, 2label={above,F_{\sigma}}, 2label={below,\iso}] {\Longrightarrow}
node[between=g1 and g2 at .55, 2label={above,G_{\sigma}}, 2label={below,\iso}] {\Longrightarrow}
;
\end{tikzpicture}
\end{equation}
is equal to the identity natural transformation
\[1_{GF(-)^\sigma} = 1_{(GF(-))^\sigma}\]
with the equality from \cref{GFpFGq}.  By the invertibility of the functor $G$ \cref{NFcsigmaFcp} and the pseudo symmetry isomorphism $G_\sigma$, the equality
\begin{equation}\label{GsigmaFpGFsigmap}
G_{\sigma; F(-)} \circ G(F_{\sigma;-}) = 1_{GF(-)^\sigma}
\end{equation}
uniquely determines $F_{\sigma}$.  Moreover, the pasting \cref{FsigmaGsigma} is the pseudo symmetry isomorphism $(GF)_\sigma$ in \cref{pseudosymisopasting} \emph{if} we already know that $(F,\{F_\sigma\})$ is a pseudo symmetric $\Cat$-multifunctor.
\end{explanation}

\begin{lemma}\label{FFsigmaps}
In the context of \cref{def:psFetaepz,def:Fsigmapswhitehead}, 
\[\begin{tikzcd}[column sep=huge]
\M \ar{r}{(F,\{F_\sigma\})} & \N
\end{tikzcd}\]
is a pseudo symmetric $\Cat$-multifunctor.
\end{lemma}

\begin{proof}
By \cref{Fcatmultifunctor} $F$ preserves the colored units \cref{enr-multifunctor-unit} and composition \cref{v-multifunctor-composition}.  It remains to check the axioms \cref{pseudosmf-unity,pseudosmf-product,pseudosmf-topeq,pseudosmf-boteq} for $(F,\{F_\sigma\})$.

\emph{Unit Permutation} \cref{pseudosmf-unity}.  This axiom says that 
\[F_{\id_n} = 1_F.\]
To prove this equality, consider an $n$-ary 1-cell $p \in \M\smscmap{\angc;c'}$.  The unit permutation axiom \cref{pseudosmf-unity} for $(G,\{G_\sigma\})$ yields the equalities of $n$-ary 2-cells
\[G^\inv_{\id_n; Fp} = 1_{GFp}^\inv = 1_{GFp}.\]
With $G$ the isomorphism in \cref{NFcsigmaFcp}, we apply $G^\inv$ to the definition \cref{GFsigmapGinvsigmaFp} of $F_\sigma$ to obtain the following equalities.
\[\begin{split}
F_{\id_n;p} &= G^\inv\left(G^\inv_{\id_n; Fp}\right)\\
&= G^\inv \left(1_{GFp}\right)\\
&= 1_{Fp}
\end{split}\]
This proves the unit permutation axiom for $F$.

\emph{Product Permutation} \cref{pseudosmf-product}. By
\begin{itemize}
\item the invertibility of the functor $G$ \cref{NFcsigmaFcp},
\item the invertibility of each pseudo symmetry isomorphism $G_\sigma$, and
\item the product permutation axiom \cref{pseudosmf-product} for $(G,\{G_\sigma\})$,
\end{itemize} 
it suffices to show that the two sides of \cref{pseudosmf-product} pasted with the corresponding sides for $(G,\{G_\sigma\})$ are equal.  The two resulting pastings are \cref{prodpermleft,prodpermright}.  Each of these pastings is equal to the identity natural transformation by \cref{GsigmaFpGFsigmap} for $\sigma$, $\tau$, and $\sigma\tau$.

\emph{Top and Bottom Equivariance}.  The axioms \cref{pseudosmf-topeq,pseudosmf-boteq} for $(F,\{F_\sigma\})$ are proved in the same way.  In each case, we
\begin{itemize}
\item paste each side of the desired equality with the corresponding pasting for $(G,\{G_\sigma\})$ and
\item observe that both sides are equal to the identity by \cref{GsigmaFpGFsigmap}.
\end{itemize}
This proves that $(F,\{F_\sigma\})$ is a pseudo symmetric $\Cat$-multifunctor.
\end{proof}

\begin{lemma}\label{etapscatmt}
In the context of \cref{def:psFetaepz,def:Fsigmapswhitehead}, 
\[\begin{tikzcd}[column sep=large]
1_{\M} \ar{r}{\eta} & GF
\end{tikzcd}\]
is a pseudo symmetric $\Cat$-multinatural isomorphism.
\end{lemma}

\begin{proof}
By \cref{etacatmultinatural} $\eta$ is a non-symmetric $\Cat$-multinatural isomorphism.  It remains to check the pseudo symmetry preservation axiom \cref{pspreservation}.  The axiom \cref{pspreservation} for $\eta$ says that, for each $n$-ary 1-cell $p \in \M\smscmap{\angc;c'}$, the following two composite $n$-ary 2-cells in $\M$ are equal.
\begin{equation}\label{etapspresaxiom}
\begin{tikzpicture}[xscale=1,yscale=1,vcenter]
\def\a{60} \def\s{.8} \def\v{-1.3}
\draw[0cell=\s]
(0,0) node (x11) {\angc\sigma}
(x11)++(3.5,0) node (x12) {GF\angc\sigma}
(x11)++(0,\v) node (x21) {c'}
(x12)++(0,\v) node (x22) {GFc'}
;
\draw[1cell=\s]  
(x11) edge node[pos=.35] {\eta_{\angc\sigma}} (x12)
(x21) edge node[pos=.35] {\eta_{c'}} (x22)
(x11) edge node[swap] {p^\sigma} (x21)
(x12) edge[bend right=\a] node[swap,pos=.45] {GF(p^\sigma)} (x22)
(x12) edge[bend left=\a] node[pos=.45] {(GFp)^\sigma = GF(p^\sigma)} (x22)
;
\draw[2cell=\s]
node[between=x12 and x22 at .65, 2label={above,(GF)_{\sigma;p}}] {\Longrightarrow}
;
\end{tikzpicture}
\end{equation}
By \cref{expl:Fsigmappswhitehead,FFsigmaps} there are equalities
\[(GF)_\sigma = 1_{GF(-)^\sigma} = 1_{(GF(-))^\sigma}.\]
Thus, the two composites in \cref{etapspresaxiom} are equal by the object $\Cat$-naturality condition for $\eta$, which is proved in \cref{etaobjcatnaturality}.
\end{proof}

\begin{lemma}\label{epzpscatmt}
In the context of \cref{def:psFetaepz,def:Fsigmapswhitehead}, 
\[\begin{tikzcd}[column sep=large]
FG \ar{r}{\epz} & 1_{\N}
\end{tikzcd}\]
is a pseudo symmetric $\Cat$-multinatural isomorphism.
\end{lemma}

\begin{proof}
By \cref{epzcatmultiiso} $\epz$ is a non-symmetric $\Cat$-multinatural isomorphism.  It remains to check the pseudo symmetry preservation axiom \cref{pspreservation}.  By the definition \cref{pseudosymisopasting} of $(FG)_\sigma$, the axiom \cref{pspreservation} for $\epz$ says that, for each $n$-ary 1-cell $q \in \N\smscmap{\angd;d'}$, the following two composite $n$-ary 2-cells in $\N$ are equal, with the unlabeled vertical arrow given by $F(Gq)^\sigma$.
\begin{equation}\label{epzpspresaxiom}
\begin{tikzpicture}[xscale=1,yscale=1,vcenter]
\def\a{82} \def\w{2} \def\s{.8} \def\v{-1.3} \def\u{-.2}
\draw[0cell=\s]
(0,0) node (x11) {FG\angd\sigma}
(x11)++(4,0) node (x12) {\angd\sigma}
(x11)++(0,\v) node (x21) {FGd'}
(x12)++(0,\v) node (x22) {d'}
;
\draw[1cell=\s]  
(x11) edge node[pos=.75] {\epz_{\angd\sigma}} (x12)
(x21) edge node[pos=.75] {\epz_{d'}} (x22)
(x12) edge node {q^\sigma} (x22)
(x11) edge[bend right=\a, looseness=\w] node[swap,pos=.45] (l) {FG(q^\sigma)} (x21)
(x11) edge (x21)
(x11) edge[bend left=\a, looseness=\w] node[pos=.45] (r) {(FGq)^\sigma} (x21)
;
\draw[2cell=.75]
node[between=l and r at .35, shift={(0,\u)}, 2label={above,F(G_{\sigma;q})}] {\Rightarrow}
node[between=l and r at .65, shift={(0,\u)}, 2label={above,F_{\sigma;\, Gq}}] {\Rightarrow}
;
\end{tikzpicture}
\end{equation}
It suffices to show that the two composites in \cref{epzpspresaxiom} have the same image under the isomorphism
\[\begin{tikzcd}[column sep=large]
\N\scmap{FG\angd\sigma; d'} \ar{r}{G}[swap]{\iso} & \M\scmap{GFG\angd\sigma; Gd'}.
\end{tikzcd}\]
By the functoriality of $\ga$, this means the following equality of $n$-ary 2-cells in $\M$.
\begin{equation}\label{epzpspresaxiom-ii}
G\left[\ga\lrscmap{1_{\epz_{d'}}; F_{\sigma;\, Gq}} \circ 
\ga\lrscmap{1_{\epz_{d'}}; F(G_{\sigma;q})} \right] = G1_{\ga\smscmap{q^\sigma; \epz_{\angd\sigma}}}
\end{equation}
By the functoriality and the composition axiom \cref{v-multifunctor-composition} for $G$, the equality \cref{epzpspresaxiom-ii} becomes the following equality.
\begin{equation}\label{epzpspresaxiom-iii}
\ga\lrscmap{1_{G\epz_{d'}}; G(F_{\sigma;\, Gq})} \circ 
\ga\lrscmap{1_{G\epz_{d'}}; GF(G_{\sigma;q})} 
= 1_{\ga\smscmap{G(q^\sigma); G\epz_{\angd\sigma}}}
\end{equation}
We prove \cref{epzpspresaxiom-iii} by computing the two terms on the left-hand side and observing that they are inverses of each other.

For the leftmost term in \cref{epzpspresaxiom-iii}, we compute as follows using, in particular, the top and bottom equivariance axioms, \cref{Fsigmaboteq,Fsigmatopeq}, for $G$.
\begin{equation}\label{epzpspresaxiom-left}
\left\{\begin{aligned}
& \ga\lrscmap{1_{G\epz_{d'}}; G(F_{\sigma;\, Gq})} &&\\
&= \ga\lrscmap{1_{G\epz_{d'}}; G_{\sigma;\, FGq}^\inv} && \text{by \cref{GFsigmapGinvsigmaFp}}\\
&= \ga\lrscmap{1_{G\epz_{d'}}; G_{\sigma;\, FGq}}^\inv && \text{by functoriality of $\ga$}\\
&= \big(G_{\sigma;\, \ga\smscmap{\epz_{d'}; FGq}}\big)^\inv && \text{by \cref{Fsigmaboteq} for $(\epz_{d'}, FGq)$}\\
&= \big(G_{\sigma;\, \ga\smscmap{q; \epz_{\angd}}}\big)^\inv && \text{by \cref{Ggafepzxj}}\\
&= \ga\lrscmap{G_{\sigma;q}; 1_{G\epz_{\angd\sigma}}}^\inv && \text{by \cref{Fsigmatopeq} for $(q, \epz_{\angd})$}
\end{aligned}\right.
\end{equation}
For the second term in \cref{epzpspresaxiom-iii}, we compute as follows.
\begin{equation}\label{epzpspresaxiom-right}
\left\{\begin{aligned}
& \ga\lrscmap{1_{G\epz_{d'}}; GF(G_{\sigma;q})} &&\\
&= \scalebox{.8}{$\ga\lrscmap{1_{\eta^\inv_{Gd'}}; \ga\lrscmap{1_{\eta_{Gd'}}; \ga\lrscmap{G_{\sigma;q}; 1_{\eta^\inv_{G\angd\sigma}}}}}$} && \text{by \cref{Fthetadefinition,Gepzz}}\\
&= \scalebox{.85}{$\ga\lrscmap{1_{\ga\smscmap{\eta^\inv_{Gd'}; \eta_{Gd'}}}; \ga\lrscmap{G_{\sigma;q}; 1_{\eta^\inv_{G\angd\sigma}}}}$} && \text{by \cref{enr-multicategory-associativity} and functoriality}\\
&= \ga\lrscmap{G_{\sigma;q}; 1_{\eta^\inv_{G\angd\sigma}}} && \text{by \cref{enr-multicategory-left-unity}}\\
&= \ga\lrscmap{G_{\sigma;q}; 1_{G\epz_{\angd\sigma}}} && \text{by \cref{Gepzz}}
\end{aligned}\right.
\end{equation}
Since the last terms in \cref{epzpspresaxiom-left,epzpspresaxiom-right} are inverses of each other, we have proved the desired equality \cref{epzpspresaxiom-iii}.
\end{proof}

\subsection*{Characterization of Pseudo Symmetric $\Cat$-Multiequivalences}

We are now ready for the main result of this section, which is a practical characterization of a pseudo symmetric $\Cat$-multiequivalence (\cref{def:pscatmultiequivalence}).  The next result is the pseudo symmetric variant of \cref{thm:multiwhitehead}.

\begin{theorem}\label{thm:psmultiwhitehead}\index{pseudo symmetric!Cat-multiequivalence@$\Cat$-multiequivalence!characterization}
Suppose
\[\begin{tikzcd}[column sep=huge]
\N \ar{r}{(G,\{G_\sigma\})} & \M
\end{tikzcd}\]
is a pseudo symmetric $\Cat$-multifunctor between $\Cat$-multicategories.  Then the following three statements are equivalent.
\begin{enumerate}
\item\label{thm:psmultiwhitehead-i}
$G$ is the right adjoint of an adjoint pseudo symmetric $\Cat$-multiequivalence $(F,G, \eta,\epz)$.
\item\label{thm:psmultiwhitehead-ii}
$G$ is a pseudo symmetric $\Cat$-multiequivalence.
\item\label{thm:psmultiwhitehead-iii}
$G$ is
\begin{itemize}
\item essentially surjective on objects and
\item an isomorphism on each multimorphism category.
\end{itemize}
\end{enumerate}
\end{theorem}

\begin{proof}
We consider each implication as follows.
\begin{itemize}
\item \eqref{thm:psmultiwhitehead-i} $\Rightarrow$ \eqref{thm:psmultiwhitehead-ii} holds by \cref{def:equivalences,def:pscatmultiequivalence}.
\item \eqref{thm:psmultiwhitehead-ii} $\Rightarrow$ \eqref{thm:psmultiwhitehead-iii} holds by \cref{pscatmultieqesff}.
\item \eqref{thm:psmultiwhitehead-iii} $\Rightarrow$ \eqref{thm:psmultiwhitehead-i} holds by \cref{FFsigmaps,etapscatmt,epzpscatmt,multieqtriangleids}.
\end{itemize}
This completes the proof.
\end{proof}

\part{Grothendieck Multiequivalence from Bipermutative-Indexed Categories to Permutative Opfibrations}
\label{part:multiequivalence}

\chapter{Enriched Multicategories of Indexed Categories}
\label{ch:diagram}
Recall that $(\Cat, \times, \boldone, [,])$ is the complete and cocomplete symmetric monoidal closed category of small categories and functors with
\begin{itemize}
\item the Cartesian product as the monoidal product and
\item diagram category as the closed structure (\cref{ex:cat}).
\end{itemize}   
Suppose $\cD$ is a small tight bipermutative category  (\cref{def:embipermutativecat}).  The main purpose of this chapter is to construct a $\Cat$-multicategory $\DCat$, with \emph{additive} symmetric monoidal functors as objects, that fully incorporates the bipermutative structure of $\cD$ (\cref{thm:dcatcatmulticat}).  We emphasize that the $\Cat$-multicategory structure on $\DCat$ is, in general, \emph{not} induced by a symmetric monoidal structure.

\subsection*{The Grothendieck Construction and Inverse $K$-Theory}

In later chapters we utilize the $\Cat$-multicategory $\DCat$ as follows.
\begin{enumerate}
\item In \cref{ch:multigro} $\DCat$ is the domain of the Grothendieck construction 
\[\begin{tikzcd}[column sep=large,every label/.append style={scale=.85}]
\DCat \ar{r}{\grod} & \permcatsg.
\end{tikzcd}\] 
The codomain $\permcatsg$ is a $\Cat$-multicategory with small permutative categories as objects.  Using the $\Cat$-multicategory structures on $\DCat$ and $\permcatsg$, we show that $\grod$ is a \emph{pseudo symmetric $\Cat$-multifunctor} (\cref{def:pseudosmultifunctor}).  
\item In \cref{ch:permfib} we lift $\grod$ to a non-symmetric $\Cat$-multifunctor
\[\begin{tikzcd}[column sep=large,every label/.append style={scale=.85}]
\DCat \ar{r}{\grod} & \pfibd.
\end{tikzcd}\]
The codomain $\pfibd$ is a non-symmetric $\Cat$-multicategory with small permutative opfibrations over $\cD$ as objects.
\item In \cref{ch:gromultequiv} we prove that the lifted $\grod$ in \cref{ch:permfib} is a non-symmetric $\Cat$-\emph{multiequivalence}.  This means that the lifted Grothendieck construction $\grod$ is
\begin{itemize}
\item essentially surjective on objects and
\item an isomorphism on each multimorphism category.
\end{itemize} 
\item In \cref{ch:multifunctorA,ch:invK} we apply these results to the indexing category $\cA$ (\cref{ex:mandellcategory}) for inverse $K$-theory $\groa A(-)$.  The pseudo symmetric $\Cat$-multifunctoriality of the Grothendieck construction $\groa$ implies that inverse $K$-theory is a pseudo symmetric $\Cat$-multifunctor but \emph{not} a $\Cat$-multifunctor.  
\end{enumerate}

\subsection*{Organization}

In \cref{sec:diagrampermutative} we first consider the diagram category $\Dcat$ when $(\cD,\otimes,\tu)$ is merely a permutative category (\cref{def:symmoncat}).  We recall the important fact that the symmetric monoidal closed structure on the diagram category $\DV$ induces a $\V$-multicategory structure on $\DV$.  Specializing to the case $\V = \Cat$, we obtain a $\Cat$-multicategory structure on $\Dcat$.  The rest of this section contains a detailed description of the $\Cat$-multicategory $\Dcat$.  To emphasize the monoidal structure in use, we also denote $\Dcat$ by $\Dtecat$ in subsequent sections.

With $(\cD,\oplus,\otimes)$ a small tight bipermutative category, in \cref{sec:additivenatmod} we define some of the data of $\DCat$. 
\begin{itemize}
\item Its objects are \emph{additive} symmetric monoidal functors, which are symmetric monoidal functors relative to the additive structure $\Dplus = (\cD, \oplus, \zero)$ of $\cD$.  This is in contrast to $\Dtecat$, whose objects are functors $\cD \to \Cat$ with no extra structure assumed. 
\item In each multimorphism category in $\DCat$, the objects are \emph{additive} natural transformations.  These are natural transformations that satisfy a unity axiom and an additivity axiom that generalize those of monoidal natural transformations; see \cref{expl:additivenattr} \eqref{expl:additivenattr-iv}. 
\item The morphisms between them are \emph{additive} modifications.  These are modifications that satisfy a unity axiom and an additivity axiom.
\end{itemize}
Additive natural transformations and additive modifications make use of both the additive structure $\Dplus$ and the multiplicative structure $\Dte$ of $\cD$.  Moreover, their additivity axioms involve Laplaza coherence isomorphisms in $\cD$; see \cref{laplazaapp,additivemodadditivity}.  This section ends with the verification that each positive arity multimorphism category in $\DCat$ is a subcategory of the corresponding multimorphism category of $\Dtecat$; see \cref{dcatnarycategory}.

In \cref{sec:diagrambipermutative} we define the $\Cat$-multicategory structure on $\DCat$ as the one inherited from $\Dtecat$ in an appropriate sense.  Most of the work involves checking that the $\Cat$-multicategory structure on $\Dtecat$ restricts to additive natural transformations and additive modifications.  In other words, we need to check that their unity and additivity axioms are preserved by the symmetric group action and composition in $\Dtecat$.

\subsection*{Summary}

There are two important differences between $\DCat$ and $\Dtecat$. 
\begin{enumerate}
\item The $\Cat$-multicategory structure on $\Dtecat$ is induced by its symmetric monoidal closed structure given by the Day convolution with respect to the multiplicative structure $\otimes$ in $\cD$.  On the other hand, the $\Cat$-multicategory structure on $\DCat$ does \emph{not} arise from a symmetric monoidal structure.  See \cref{expl:dcatbipermzero}.  
\item The $\Cat$-multicategory structure on $\DCat$ is defined using the one on $\Dtecat$.  However, $\DCat$ is \emph{not} a sub-$\Cat$-multicategory of $\Dtecat$.  See \cref{expl:dcatforgetdtecat}.
\end{enumerate} 
The following table summaries $\Dtecat$ and $\DCat$, with \emph{symmetric monoidal} and \emph{multimorphism category} abbreviated to, respectively, \emph{sm} and \emph{mc}.
\begin{center}
\resizebox{\columnwidth}{!}{%
{\renewcommand{\arraystretch}{1.4}%
{\setlength{\tabcolsep}{1ex}
\begin{tabular}{|c|cc|cc|}\hline
$\Cat$-multicategories & $\Dtecat$ &(\ref{expl:permindexedcat}) 
& $\DCat$ &(\ref{thm:dcatcatmulticat}) \\ \hline
objects 
& functors $\cD \to \Cat$ &(\ref{def:diagramcat}) 
& additive sm functors $\Dplus \to \Cat$ &(\ref{def:additivesmf}) \\ \hline
mc objects 
& natural transformations &(\ref{phitensorxz})
& additive natural transformations &(\ref{def:additivenattr}) \\ \hline
mc morphisms 
& modifications &(\ref{Phiphivarphi})
& additive modifications &(\ref{def:additivemodification})\\ \hline
induced by sm structure
& yes & (\ref{thm:permindexedcat})
& no & (\ref{expl:dcatbipermzero})
\\ \hline
\end{tabular}}}}
\end{center}
\smallskip
We remind the reader of our left normalized bracketing \cref{expl:leftbracketing} for iterated monoidal product.

\section{Enriched Multicategories of Permutative-Indexed Categories}
\label{sec:diagrampermutative}

Recall the diagram category $\DC$ (\cref{def:diagramcat}), which has $\cD$-diagrams in $\C$ as objects and natural transformations as morphisms, and the notion of a $\V$-multicategory (\cref{def:enr-multicategory}).  For a small permutative category $\cD$, in \cref{expl:dayconvolution,expl:dcatsmclosed} we described the symmetric monoidal structure on the diagram category $\Dcat$, which has $\cD$-indexed categories as objects, and its enriched, tensored, and cotensored structures over $\Cat$.  In this section we describe in detail the $\Cat$-enriched multicategory structure on $\Dcat$.  The following theorem is an instance of \cite[6.3.6]{cerberusIII} and guarantees the existence of the $\Cat$-multicategory $\Dcat$.

\begin{theorem}\label{thm:permindexedcat}\index{diagram category!multicategory structure}\index{enriched multicategory!diagram category}\index{Cat-multicategory@$\Cat$-multicategory!of diagrams}
Suppose $\V$ is a symmetric monoidal closed category that is complete and cocomplete.  Suppose $\cD$ is a small symmetric monoidal category.  Then the symmetric monoidal closed category $\DV$ in \cref{thm:Day} has the structure of a $\V$-multicategory with the same underlying $\V$-category as that of $\DV$.
\end{theorem}

Now we specialize to the case $\V = \Cat$, which is the main case of interest in this work.

\begin{explanation}[$\Cat$-Multicategory of Permutative-Indexed Categories]\label{expl:permindexedcat}
In the context of \cref{thm:permindexedcat}, suppose 
\begin{itemize}
\item $\V = (\Cat, \times, \boldone, [,])$ is the symmetric monoidal closed category of small categories and functors with the Cartesian product as the monoidal product and diagram category as the closed structure (\cref{ex:cat}) and
\item $(\cD,\otimes, \tu, \beta)$ is a small permutative category (\cref{def:symmoncat}).
\end{itemize}
By \cref{thm:permindexedcat} \label{not:Dcat}$\Dcat$ has the structure of a $\Cat$-multicategory, which we also denote by $\Dcat$, with the same underlying 2-category as that of $\Dcat$.  Here we describe this $\Cat$-multicategory structure in detail.
\begin{description}
\item[Objects]
The objects of the $\Cat$-multicategory $\Dcat$ are $\cD$-indexed categories (\cref{def:diagramcat}), that is, functors $\cD \to \Cat$.
\item[Multimorphism Categories in Arity 0]
Suppose $Z \cn \cD \to \Cat$ is a $\cD$-indexed category.  The \emph{0-ary multimorphism category} is given by
\[(\Dcat)\scmap{\ang{};Z} = (\Dcat)\srb{J,Z}\]
with $\ang{}$ the empty sequence and
\[J = \cD(\tu,-) \cn \cD \to \Cat\]
the unit diagram in \cref{dcatunit}.  As discussed in \cref{expl:dcatsmclosed}, the category $(\Dcat)\srb{J,Z}$ has
\begin{itemize}
\item natural transformations $J \to Z$ as objects,
\item modifications as morphisms,
\item identity modifications as identities, and
\item composition given by vertical composition of modifications.
\end{itemize}
By the Bicategorical Yoneda Lemma\index{Bicategorical!Yoneda Lemma} \cite[8.3.28]{johnson-yau} applied to the current 2-categorical context, there is a canonical isomorphism of categories
\begin{equation}\label{dcatzeroary}
(\Dcat)\scmap{\ang{};Z} = (\Dcat)\srb{J,Z} \iso Z\tu
\end{equation}
given by the assignments
\[\begin{tikzcd}[row sep=0ex,
/tikz/column 1/.append style={anchor=base east},
/tikz/column 2/.append style={anchor=base west}]
\phi \ar[mapsto]{r} & \phi_{\tu}(1_{\tu})\\
\Phi \ar[mapsto]{r} & (\Phi_{\tu})_{1_{\tu}}
\end{tikzcd}\]
for each
\begin{itemize}
\item natural transformation $\phi \cn J \to Z$ and
\item modification $\Phi \cn \phi \to \phi'$ with $\phi, \phi' \cn J \to Z$ natural transformations.
\end{itemize} 
\item[Multimorphism Categories in Positive Arity]
Suppose given $\cD$-indexed categories
\[Z, X_1, \ldots, X_n \cn \cD \to \Cat\]
with $n > 0$ and $\angX = (X_1, \ldots, X_n)$.  The \emph{$n$-ary multimorphism category} is given by
\begin{equation}\label{dcatnary}
(\Dcat)\scmap{\angX;Z} = (\Dcat)\brb{\textstyle \bigotimes_{j=1}^n X_j,Z}.
\end{equation}
For $n \geq 2$ the $\cD$-indexed category $\bigotimes_{j=1}^n X_j$ is the iterated Day convolution in \cref{expl:nday} with $\V = \Cat$.  As discussed in \cref{expl:dcatsmclosed}, the category $(\Dcat)\brb{\bigotimes_{j=1}^n X_j,Z}$ has
\begin{itemize}
\item natural transformations $\bigotimes_{j=1}^n X_j \to Z$ as objects,
\item modifications as morphisms,
\item identity modifications as identities, and
\item composition given by vertical composition of modifications.
\end{itemize}
Unraveling the iterated Day convolution \cref{iteratedday}, we explicitly describe these natural transformations and modifications below.

A \emph{natural transformation}\index{natural transformation}
\begin{equation}\label{phitensorxz}
\begin{tikzcd}[column sep=large]
\bigotimes_{j=1}^n X_j \ar{r}{\phi} & Z
\end{tikzcd}
\end{equation}
consists of, for each $n$-tuple of objects 
\begin{equation}\label{angaaonean}
\anga = \ang{a_j}_{j=1}^n \in \cD^n \andspace a = \txotimes_{j=1}^n a_j \in \cD,
\end{equation}
a component functor 
\begin{equation}\label{dcatnaryobjcomponent}
\begin{tikzcd}[column sep=huge]
\prod_{j=1}^n X_j a_j \ar{r}{\phi_{\anga}} & Za.
\end{tikzcd}
\end{equation}
For each $n$-tuple of morphisms $\ang{f_j \cn a_j \to a_j'}_{j=1}^n$ in $\cD$, denote
\begin{equation}\label{aprimef}
\angap = \ang{a_j'}_{j=1}^n \scs \quad 
a' = \txotimes_{j=1}^n a_j' \scs \andspace f = \txotimes_{j=1}^n f_j \cn a \to a'.
\end{equation}
The following \emph{naturality diagram}\index{naturality diagram} of functors is required to commute.
\begin{equation}\label{dcatnaryobjnaturality}
\begin{tikzpicture}[xscale=3.5,yscale=1.5,vcenter]
\draw[0cell=.9]
(0,0) node (x11) {\textstyle\prod_{j=1}^n X_j a_j}
(x11)++(1,0) node (x12) {Za}
(x11)++(0,-1) node (x21) {\textstyle\prod_{j=1}^n X_j a_j'}
(x12)++(0,-1) node (x22) {Za'}
;
\draw[1cell=.9]  
(x11) edge node {\phi_{\anga}} (x12)
(x12) edge node {Zf} (x22)
(x11) edge node[swap] {\textstyle\prod_{j=1}^n X_j f_j} (x21)
(x21) edge node {\phi_{\angap}} (x22)
;
\end{tikzpicture}
\end{equation}

A \emph{modification}\index{modification} $\Phi \cn \phi \to \varphi$ in $(\Dcat)\brb{\bigotimes_{j=1}^n X_j,Z}$ as in
\begin{equation}\label{Phiphivarphi}
\begin{tikzpicture}[xscale=2,yscale=1.7,baseline={(x1.base)}]
\draw[0cell=.9]
(0,0) node (x1) {\textstyle \bigotimes_{j=1}^n X_j}
(x1)++(.14,.03) node (x2) {\phantom{\textstyle \bigotimes Z}}
(x2)++(1,0) node (x3) {Z}
;
\draw[1cell=.9]  
(x2) edge[bend left, shorten <=-.5ex] node[pos=.43] {\phi} (x3)
(x2) edge[bend right, shorten <=-.5ex] node[swap,pos=.43] {\varphi} (x3)
;
\draw[2cell]
node[between=x2 and x3 at .45, rotate=-90, 2label={above,\Phi}] {\Rightarrow}
;
\end{tikzpicture}
\end{equation}
consists of a component natural transformation
\begin{equation}\label{dcatnarymodcomponent}
\begin{tikzpicture}[xscale=2.5,yscale=2,baseline={(x1.base)}]
\draw[0cell=.9]
(0,0) node (x1) {\textstyle \prod_{j=1}^n X_j a_j}
(x1)++(.17,.03) node (x2) {\phantom{\textstyle \bigotimes Z}}
(x2)++(1,0) node (x3) {Za}
;
\draw[1cell=.9]  
(x2) edge[bend left, shorten <=-.5ex] node[pos=.5] {\phi_{\anga}} (x3)
(x2) edge[bend right, shorten <=-.5ex] node[swap,pos=.5] {\varphi_{\anga}} (x3)
;
\draw[2cell]
node[between=x2 and x3 at .4, rotate=-90, 2label={above,\Phi_{\anga}}] {\Rightarrow}
;
\end{tikzpicture}
\end{equation}
for each $n$-tuple of objects $\anga = \{a_j\}_{j=1}^n$ in $\cD$ as in \cref{angaaonean}.  The \emph{modification axiom}\index{modification!axiom} in this case is equivalent to, in the context of \cref{angaaonean,aprimef}, the equality of the following two whiskered natural transformations.
\begin{equation}\label{dcatnarymodaxiom}
\begin{tikzpicture}[xscale=1,yscale=1,vcenter]
\def\d{25} \def\v{-1.5}
\draw[0cell=.85]
(0,0) node (x11) {\txprod_{j=1}^n X_j a_j}
(x11)++(.5,0) node (a) {\phantom{Z}}
(a)++(2.5,0) node (b) {\phantom{Z}}
(b)++(.1,0) node (x12) {Za}
(x11)++(0,\v) node (x21) {\txprod_{j=1}^n X_j a_j'}
(a)++(0,\v) node (a2) {\phantom{Z}}
(b)++(0,\v) node (b2) {\phantom{Z}}
(x12)++(0,\v) node (x22) {Za'}
;
\draw[1cell=.8] 
(a) edge[bend left=\d] node[pos=.4] {\phi_{\anga}} (b)
(a) edge[bend right=\d] node[swap,pos=.6] {\varphi_{\anga}} (b)
(a2) edge[bend left=\d] node[pos=.4] {\phi_{\angap}} (b2)
(a2) edge[bend right=\d] node[swap,pos=.6] {\varphi_{\angap}} (b2)
(x11) edge[transform canvas={xshift={1em}}] node[swap] {\txprod_{j=1}^n X_j f_j} (x21)
(x12) edge node {Zf} (x22)
;
\draw[2cell=.9]
node[between=a and b at .4, rotate=-90, 2label={above,\Phi_{\anga}}] {\Rightarrow}
node[between=a2 and b2 at .4, rotate=-90, 2label={above,\Phi_{\angap}}] {\Rightarrow}
;
\end{tikzpicture}
\end{equation}
This finishes the description of the category $(\Dcat)\smscmap{\angX;Z}$ in \cref{dcatnary}.
\item[Colored Units]
For a $\cD$-indexed category $Z \cn \cD \to \Cat$, the \emph{$Z$-colored unit} 
\begin{equation}\label{zcoloredunit}
\begin{tikzcd}[column sep=large]
Z \ar{r}{1_Z} & Z
\end{tikzcd}
\end{equation}
is the identity natural transformation.
\item[Symmetric Group Action]
For a permutation $\sigma \in \Sigma_n$, the \emph{symmetric group action}
\begin{equation}\label{dcatsigmaaction}
\begin{tikzcd}[column sep=large]
(\Dcat)\scmap{\angX;Z} \ar{r}{\sigma}[swap]{\iso} & (\Dcat)\scmap{\angX\sigma;Z}
\end{tikzcd}
\end{equation}
is the functor that sends a natural transformation $\phi \cn \bigotimes_{j=1}^n X_j \to Z$ as in \cref{phitensorxz} to the natural transformation
\[\begin{tikzcd}[column sep=large]
\bigotimes_{j=1}^n X_{\sigma(j)} \ar{r}{\phi^\sigma} & Z.
\end{tikzcd}\]
For an $n$-tuple of objects $\anga = \ang{a_j}_{j=1}^n$ in $\cD$ as in \cref{angaaonean}, $\phi^\sigma$ has the following component composite functor.
\begin{equation}\label{phisigmaanga}
\begin{tikzpicture}[xscale=3.5,yscale=1.5,vcenter]
\draw[0cell=.9]
(0,0) node (x11) {\textstyle\prod_{j=1}^n X_{\sigma(j)} a_j}
(x11)++(1,0) node (x12) {Za}
(x11)++(0,-1) node (x21) {\textstyle\prod_{j=1}^n X_j a_{\sigmainv(j)}}
(x12)++(0,-1) node (x22) {Za_{\sigmainv}}
;
\draw[1cell=.9]  
(x11) edge node {\phi^\sigma_{\anga}} (x12)
(x11) edge node {\iso} node[swap] {\sigma} (x21)
(x21) edge node {\phi_{\sigma\anga}} (x22)
(x22) edge node {\iso} node[swap] {Z(\sigmainv)} (x12)
;
\end{tikzpicture}
\end{equation}
The arrows in \cref{phisigmaanga} are given as follows.
\begin{itemize}
\item The left vertical arrow permutes the $n$ factors according to $\sigma \in \Sigma_n$.
\item The bottom horizontal arrow $\phi_{\sigma\anga}$ is the component functor of $\phi$ as in \cref{dcatnaryobjcomponent} for the $n$-tuple of objects
\begin{equation}\label{sigmaangaasigmainv}
\sigma\anga = \ang{a_{\sigmainv(j)}}_{j=1}^n \in \cD^n \withspace 
a_{\sigmainv} = \txotimes_{j=1}^n a_{\sigmainv(j)} \in \cD.
\end{equation}
\item The right vertical arrow $Z(\sigmainv)$ is the image under $Z$ of the unique coherence isomorphism \cite[XI.1]{maclane}
\begin{equation}\label{asigmainvtoa}
\begin{tikzcd}[column sep=large]
a_{\sigmainv} = \bigotimes_{j=1}^n a_{\sigmainv(j)} \ar{r}{\sigmainv}[swap]{\iso} & 
\bigotimes_{j=1}^n a_j = a
\end{tikzcd}
\inspace (\cD,\otimes,\beta)
\end{equation}
that permutes the $n$ factors according to $\sigmainv \in \Sigma_n$.
\end{itemize}
This finishes the description of the symmetric group action \cref{dcatsigmaaction} on natural transformations.

For a modification $\Phi \cn \phi \to \varphi$ as in \cref{Phiphivarphi}, its image under the symmetric group action $\sigma$ in \cref{dcatsigmaaction} is the following modification.
\[\begin{tikzpicture}[xscale=2,yscale=1.7,baseline={(x1.base)}]
\draw[0cell=.9]
(0,0) node (x1) {\textstyle \bigotimes_{j=1}^n X_{\sigma(j)}}
(x1)++(.24,.04) node (x2) {\phantom{\bigotimes Z}}
(x2)++(1,0) node (x3) {Z}
;
\draw[1cell=.9]  
(x2) edge[bend left, shorten <=-1ex] node[pos=.47] {\phi^\sigma} (x3)
(x2) edge[bend right, shorten <=-1ex] node[swap,pos=.46] {\varphi^\sigma} (x3)
;
\draw[2cell]
node[between=x2 and x3 at .42, rotate=-90, 2label={above,\Phi^\sigma}] {\Rightarrow}
;
\end{tikzpicture}\]
For an $n$-tuple of objects $\anga = \ang{a_j}_{j=1}^n$ in $\cD$ as in \cref{angaaonean}, $\Phi^\sigma$ has the component natural transformation $\Phi^\sigma_{\anga}$ given by the following whiskering.
\begin{equation}\label{Phisigmaanga}
\begin{tikzpicture}[xscale=1,yscale=1,vcenter]
\def\d{25} \def\v{-1.5}
\draw[0cell=.85]
(0,0) node (x11) {\txprod_{j=1}^n X_{\sigma(j)} a_j}
(x11)++(.75,0) node (a) {\phantom{Z}}
(a)++(2.5,0) node (b) {\phantom{Z}}
(b)++(.2,0) node (x12) {Za}
(x11)++(0,\v) node (x21) {\txprod_{j=1}^n X_j a_{\sigmainv(j)}}
(a)++(0,\v) node (a2) {\phantom{Z}}
(b)++(0,\v) node (b2) {\phantom{Z}}
(x12)++(0,\v) node (x22) {Za_{\sigmainv}}
;
\draw[1cell=.8] 
(a) edge[bend left=\d] node[pos=.4] {\phi^\sigma_{\anga}} (b)
(a) edge[bend right=\d] node[swap,pos=.6] {\varphi^\sigma_{\anga}} (b)
(a2) edge[bend left=\d] node[pos=.4] {\phi_{\sigma\anga}} (b2)
(a2) edge[bend right=\d] node[swap,pos=.6] {\varphi_{\sigma\anga}} (b2)
(x11) edge[transform canvas={xshift={1.4em}}] node[swap] {\sigma} node {\iso} (x21)
(x22) edge node {\iso} node[swap] {Z(\sigmainv)} (x12)
;
\draw[2cell=.9]
node[between=a and b at .4, rotate=-90, 2label={above,\Phi^\sigma_{\anga}}] {\Rightarrow}
node[between=a2 and b2 at .4, rotate=-90, 2label={above,\Phi_{\sigma\anga}}] {\Rightarrow}
;
\end{tikzpicture}
\end{equation}
In \cref{Phisigmaanga} $\Phi_{\sigma\anga}$ is the component natural transformation of $\Phi$ as in \cref{dcatnarymodcomponent} for the $n$-tuple of objects $\sigma\anga$ in \cref{sigmaangaasigmainv}.  This finishes the description of the symmetric group action \cref{dcatsigmaaction} on the $\Cat$-multicategory $\Dcat$.
\item[Composition]
Suppose given $\cD$-indexed categories
\[Z \scs \ang{X_j}_{j=1}^n \scs \ang{\ang{W_{ji}}_{i=1}^{\ell_j}}_{j=1}^n \cn \cD \to \Cat\]
with $n > 0$, $\angX = (X_1, \ldots, X_n)$,
\[\angWj = \big(W_{j1},\ldots,W_{j\ell_j}\big), \andspace \angW = \big(\ang{W_1}, \ldots, \ang{W_n}\big).\]
The \emph{multicategorical composition}
\begin{equation}\label{dcatgamma}
\begin{tikzpicture}[xscale=5,yscale=1.2,baseline={(x11.base)}]
\draw[0cell=.85]
(0,0) node (x11) {(\Dcat)\scmap{\angX;Z} \times \txprod_{j=1}^n (\Dcat)\scmap{\angWj;X_j}}
(x11)++(1,0) node (x12) {(\Dcat)\scmap{\angW;Z}}
;
\draw[1cell=.9]  
(x11) edge node {\gamma} (x12)
;
\end{tikzpicture}
\end{equation}
is the following composite functor.
\begin{equation}\label{dcatgammacomp}
\begin{tikzpicture}[xscale=3.5,yscale=1.2,vcenter]
\def\a{15}
\draw[0cell=.85]
(0,0) node (x11) {\textstyle \Dcat\brb{\bigotimes_{j=1}^n X_j,Z} \times \prod_{j=1}^n \Dcat\brb{\bigotimes_{i=1}^{\ell_j} W_{ji},X_j}}
(x11)++(1,-1) node (x12) {\textstyle \Dcat\brb{\bigotimes_{j=1}^n \bigotimes_{i=1}^{\ell_j} W_{ji},Z}}
(x11)++(0,-2) node (x2) {\textstyle \Dcat\brb{\bigotimes_{j=1}^n X_j,Z} \times \Dcat\brb{\bigotimes_{j=1}^n \bigotimes_{i=1}^{\ell_j} W_{ji}, \bigotimes_{j=1}^n X_j}}
;
\draw[1cell=.85]  
(x11) edge[bend left=\a, transform canvas={xshift={4ex}}] node {\gamma} (x12)
(x11) edge node[swap] {1 \times \otimes} (x2)
(x2) edge[bend right=\a, transform canvas={xshift={4ex}}] node[swap,pos=.6] {m} (x12)
;
\end{tikzpicture}
\end{equation}
\begin{itemize}
\item In \cref{dcatgammacomp} the arrow labeled $m$ is the horizontal composition in the 2-category $\Dcat$ (\cref{expl:dcatsmclosed}).
\item In the left vertical arrow in \cref{dcatgammacomp}, $\otimes$ is the Day convolution in
\begin{itemize}
\item \cref{eq:dayconvolution} for $\cD$-indexed categories,
\item \cref{daynattr} for natural transformations, and
\item \cref{daymodifications} for modifications.
\end{itemize} 
\end{itemize}
The diagram \cref{dcatgammacomp} assumes that each $\ell_j > 0$.

If $\ell_j = 0$ for some but not all $j \in \{1,\ldots,n\}$, then one makes the following adjustments in \cref{dcatgammacomp}.
\begin{itemize}
\item Since $\angWj = \ang{}$ for each index $j$ with $\ell_j=0$, along the left vertical arrow in \cref{dcatgammacomp} the $\cD$-indexed category $\bigotimes_{i=1}^{\ell_j} W_{ji}$ should be interpreted as the unit diagram $J = \cD(\tu,-)$ in \cref{dcatunit}, which is the monoidal unit in $\Dcat$ with respect to the Day convolution.
\item These copies of $J$ for $\ell_j=0$ are removed using the left and right unit isomorphisms in $\Dcat$ before the horizontal composition $m$.
\item In the codomain of $\ga$, the domain $\cD$-indexed category
\[\textstyle\bigotimes_{j=1}^n \bigotimes_{i=1}^{\ell_j} W_{ji}\]
should omit those indices $j$ with $\ell_j=0$.
\end{itemize}
Finally, if $\ell_j=0$ for all $j \in \{1,\ldots,n\}$, then there is one further adjustment.  Since $\angW$ is also the empty sequence $\ang{}$,  in the codomain of $\ga$, the domain $\cD$-indexed category should be $J$.
\end{description}
This completes the description of the $\Cat$-multicategory $\Dcat$ when $(\cD, \otimes, \tu, \beta)$ is a small permutative category.
\end{explanation}

\begin{explanation}[Composition in $\Dcat$]\label{expl:dcatgamma}
The composition $\ga$ in \cref{dcatgammacomp} can be described more explicitly as follows.  Suppose given natural transformations
\begin{equation}\label{phiphijwxz}
\begin{tikzcd}
\bigotimes_{j=1}^n X_j \ar{r}{\phi} & Z
\end{tikzcd}\andspace
\begin{tikzcd}
\bigotimes_{i=1}^{\ell_j} W_{ji} \ar{r}{\phi_j} & X_j
\end{tikzcd}
\end{equation}
for $1 \leq j \leq n$ with each $\ell_j > 0$.  Their composite is the natural transformation
\begin{equation}\label{gaphiphij}
\begin{tikzcd}[column sep=2cm]
\bigotimes_{j=1}^n \bigotimes_{i=1}^{\ell_j} W_{ji} \ar{r}{\ga\smscmap{\phi;\angphij}} & Z
\end{tikzcd}
\end{equation}
defined as follows.  Suppose given objects $a_{ji} \in \cD$ for $1 \leq j \leq n$ and $1 \leq i \leq \ell_j$, along with the following notation.
\begin{equation}\label{angajangaaja}
\left\{\begin{aligned}
\ang{a_j} &= (a_{j1},\ldots,a_{j\ell_j}) \in \cD^{\ell_j} & a_j &= \textstyle\bigotimes_{i=1}^{\ell_j} a_{ji} \in \cD\\
\anga &= (a_1,\ldots,a_n) \in \cD^n & a &= \textstyle\bigotimes_{j=1}^n a_j \in \cD
\end{aligned}\right.
\end{equation}
The corresponding component functor of $\ga\smscmap{\phi;\angphij}$ in \cref{gaphiphij} is the following composite.
\begin{equation}\label{dcatganattr}
\begin{tikzpicture}[xscale=4.5,yscale=1.2,vcenter]
\def\a{15}
\draw[0cell=.9]
(0,0) node (x11) {\textstyle \prod_{j=1}^n \prod_{i=1}^{\ell_j} W_{ji} a_{ji}}
(x11)++(1,0) node (x12) {Za}
(x11)++(.55,-1) node (x2) {\textstyle\prod_{j=1}^n X_j a_j}
;
\draw[1cell=.85]  
(x11) edge node {\ga\smscmap{\phi;\angphij}_{\ang{\ang{a_{ji}}_i}_j}} (x12)
(x11) edge node[swap,pos=.3] {\textstyle \prod_j (\phi_j)_{\angaj}} (x2)
(x2) edge node[swap,pos=.6] {\phi_{\anga}} (x12)
;
\end{tikzpicture}
\end{equation}
Here $\phi_{\anga}$ and $(\phi_j)_{\angaj}$ are the indicated component functors of, respective, $\phi$ and $\phi_j$ as in \cref{dcatnaryobjcomponent}.

Next, suppose given modifications $\Phi$ and $\Phi_j$ for $1 \leq j \leq n$ as follows.
\begin{equation}\label{PhiPhijwxz}
\begin{tikzpicture}[xscale=2,yscale=1.7,baseline={(x1.base)}]
\draw[0cell=.9]
(0,0) node (x1) {\textstyle \bigotimes_{j=1}^n X_j}
(x1)++(.13,.04) node (x2) {\phantom{\textstyle \bigotimes Z}}
(x2)++(1,0) node (x3) {Z}
;
\draw[1cell=.9]  
(x2) edge[bend left, shorten <=-.5ex] node[pos=.45] {\phi} (x3)
(x2) edge[bend right, shorten <=-.5ex] node[swap,pos=.45] {\varphi} (x3)
;
\draw[2cell]
node[between=x2 and x3 at .45, rotate=-90, 2label={above,\Phi}] {\Rightarrow}
;
\begin{scope}[shift={(2,.05)}]
\draw[0cell=.9]
(0,0) node (y1) {\textstyle \bigotimes_{i=1}^{\ell_j} W_{ji}}
(y1)++(.16,-.05) node (y2) {\phantom{\textstyle \bigotimes Z}}
(y2)++(1,0) node (y3) {X_j}
;
\draw[1cell=.9]  
(y2) edge[bend left, shorten <=-.5ex] node[pos=.45] {\phi_j} (y3)
(y2) edge[bend right, shorten <=-.5ex] node[swap,pos=.45] {\varphi_j} (y3)
;
\draw[2cell]
node[between=y2 and y3 at .42, rotate=-90, 2label={above,\Phi_j}] {\Rightarrow}
;
\end{scope}
\end{tikzpicture}
\end{equation}
Their composite modification
\[\begin{tikzpicture}[xscale=3.6,yscale=2.5,baseline={(x1.base)}]
\draw[0cell=.9]
(0,0) node (x1) {\textstyle \bigotimes_{j=1}^n \bigotimes_{i=1}^{\ell_j} W_{ji}}
(x1)++(.2,0) node (x2) {\phantom{\bigotimes Z}}
(x2)++(1,0) node (x3) {Z}
;
\draw[1cell=.9]  
(x2) edge[bend left, shorten <=-.75ex] node[pos=.45] {\ga\smscmap{\phi;\angphij}} (x3)
(x2) edge[bend right, shorten <=-.75ex] node[swap,pos=.45] {\ga\smscmap{\varphi;\ang{\varphi_j}}} (x3)
;
\draw[2cell=.9]
node[between=x2 and x3 at .3, rotate=-90, 2label={above,\ga\smscmap{\Phi;\ang{\Phi_j}}}] {\Rightarrow}
;
\end{tikzpicture}\]
has the following component horizontal composite natural transformation.
\begin{equation}\label{dcatgamod}
\begin{tikzpicture}[xscale=4,yscale=2,vcenter]
\def\a{30} \def\b{47}
\draw[0cell=.85]
(0,0) node (x11) {\textstyle \prod_{j=1}^n \prod_{i=1}^{\ell_j} W_{ji} a_{ji}}
(x11)++(1,0) node (x12) {\textstyle \prod_{j=1}^n X_j a_j}
(x12)++(.65,0) node (x13) {Za}
;
\draw[1cell=.85]  
(x11) edge[bend left=\a, shorten <=-1ex] node[pos=.4] {\textstyle \prod_j (\phi_j)_{\angaj}} (x12)
(x11) edge[bend right=\a, shorten <=-1ex] node[swap,pos=.4] {\textstyle \prod_j (\varphi_j)_{\angaj}} (x12)
(x12) edge[bend left=\b] node[pos=.45] {\phi_{\anga}} (x13)
(x12) edge[bend right=\b] node[swap,pos=.45] {\varphi_{\anga}} (x13)
;
\draw[2cell=.9]
node[between=x11 and x12 at .35, rotate=-90, 2label={above,\textstyle \prod_j (\Phi_j)_{\angaj}}] {\Rightarrow}
node[between=x12 and x13 at .45, rotate=-90, 2label={above,\Phi_{\anga}}] {\Rightarrow}
;
\end{tikzpicture}
\end{equation}
Here $\Phi_{\anga}$ and $(\Phi_j)_{\angaj}$ are the indicated component natural transformations of, respectively, $\Phi$ and $\Phi_j$ as in \cref{dcatnarymodcomponent}.

The discussion above assumes that each $\ell_j > 0$.  Suppose there is at least one index $j \in \{1,\ldots,n\}$ with $\ell_j = 0$.  Then one makes the appropriate adjustments as discussed after \cref{dcatgammacomp}, using the canonical isomorphism of categories in \cref{dcatzeroary}
\[(\Dcat)\scmap{\ang{};X_j} = (\Dcat)\srb{J,X_j} \iso X_j\tu.\]
For example, suppose $\ell_j = 0$ for all $j \in \{1,\ldots,n\}$.  Then the composition in \cref{dcatgammacomp}
\[\begin{tikzcd}[column sep=large]
\Dcat\brb{\bigotimes_{j=1}^n X_j, Z} \times \prod_{j=1}^n X_j\tu \ar{r}{\ga} & Z\tu
\end{tikzcd}\]
is given on objects by
\[\ga\scmap{\phi;\angx} = \phi_{\ang{\tu}} \angx \in Z\tu\]
for
\begin{itemize}
\item each natural transformation $\phi \cn \bigotimes_{j=1}^n X_j \to Z$ and
\item objects $x_j \in X_j\tu$ with $\angx = (x_1,\ldots,x_n)$.
\end{itemize}
Here $\ang{\tu}$ is the $n$-tuple with copies of $\tu \in \cD$, and
\[\begin{tikzcd}[column sep=large]
\prod_{j=1}^n X_j\tu \ar{r}{\phi_{\ang{\tu}}} & Z\tu
\end{tikzcd}\]
is the corresponding component functor of $\phi$.  This is well defined, since $\tu$ is the strict monoidal unit in $\cD$ and $\bigotimes_{j=1}^n \tu = \tu$.
\end{explanation}

\section{Additive Natural Transformations and Modifications}
\label{sec:additivenatmod}

Throughout this section we assume that
\begin{equation}\label{Dbipermutative}
\big(\cD, (\oplus, \zero, \betaplus), (\otimes, \tu, \betate), (\fal, \far)\big)
\end{equation}
is a small \emph{tight} bipermutative category (\cref{def:embipermutativecat}).  The tightness condition means that the factorization morphisms $\fal$ and $\far$ are natural isomorphisms.  For example, $\Finsk$, $\Fset$, $\Fskel$, and $\cA$ in \cref{sec:examples} are small tight bipermutative categories.  As in \cref{sec:bipermutativecat}, denote the additive structure of $\cD$ by
\begin{equation}\label{Dadditive}
\Dplus = (\cD, \oplus, \zero, \betaplus)
\end{equation}
and the multiplicative structure of $\cD$ by
\begin{equation}\label{Dmultiplicate}
\Dte = (\cD, \otimes, \tu, \betate).
\end{equation}
Each of $\Dplus$ and $\Dte$ is a permutative category (\cref{def:symmoncat}).  

Our next objective is to construct a variant of the $\Cat$-enriched multicategory $\Dtecat$ in \cref{sec:diagrampermutative}, denoted $\DCat$, that incorporates the bipermutative structure of $\cD$.  In this section, we define its
\begin{itemize}
\item objects, which are called additive symmetric monoidal functors (\cref{def:additivesmf}), and 
\item multimorphism categories in positive arity, with
\begin{itemize}
\item additive natural transformations (\cref{def:additivenattr}) as objects and
\item additive modifications (\cref{def:additivemodification}) as morphisms.
\end{itemize} 
\end{itemize}   
In \cref{sec:diagrambipermutative} we prove that $\DCat$ is a $\Cat$-multicategory.

\subsection*{Objects: Additive Symmetric Monoidal Functors}

We first define the objects in $\DCat$, which are symmetric monoidal functors with respect to the additive structure of $\cD$.

\begin{definition}\label{def:additivesmf}\index{additive!symmetric monoidal functor}\index{symmetric monoidal functor!additive}
An \emph{additive symmetric monoidal functor} is a symmetric monoidal functor (\cref{def:monoidalfunctor})
\[(X,X^2,X^0) \cn \Dplus \to \Cat.\]
We call $(X^2,X^0)$ the \emph{additive structure} of $X$ and suppress it from the notation if it is clear from the context.
\end{definition}

\begin{explanation}\label{ex:additivelysmf}
An additive symmetric monoidal functor $(X,X^2,X^0)$ is defined using the additive structure $\Dplus = (\cD, \oplus, \zero, \betaplus)$ of $\cD$ in \cref{Dadditive}.  Note that we do \emph{not} require $X^0$ or $X^2$ to be invertible in general. 
\begin{itemize}
\item The functor $X \cn \cD \to \Cat$ is a $\cD$-indexed category.
\item The \emph{unit constraint}\index{unit constraint} $X^0$ is a functor
\[\begin{tikzcd}[column sep=large]
\boldone \ar{r}{X^0} & X\zero.
\end{tikzcd}\]
So $X^0$ is determined by the \emph{unit object}\index{unit object}
\begin{equation}\label{additivesmfunitobj}
X^0* \in X\zero
\end{equation}
with $* \in \boldone$ the unique object and $\zero \in \cD$ the additive zero.
\item The \emph{monoidal constraint}\index{monoidal constraint} $X^2$ is a natural transformation with a component functor
\begin{equation}\label{additivesmfconstraint}
\begin{tikzcd}[column sep=large]
Xa \times Xb \ar{r}{X^2_{a,b}} & X(a \oplus b)
\end{tikzcd}
\end{equation}
for each pair of objects $a,b \in \cD$.
\end{itemize} 
The data $(X,X^2,X^0)$ satisfy the unity axiom \cref{monoidalfunctorunity}, the associativity axiom \cref{monoidalfunctorassoc}, and the braiding axiom \cref{monoidalfunctorbraiding} for a symmetric monoidal functor.  For future reference we state these axioms explicitly below.
\begin{description}
\item[Naturality] The naturality of $X^2$ states that, for morphisms 
\[f \cn a \to b \andspace f' \cn a' \to b'\]
in $\cD$, the following diagram commutes, with $(-)_* = X(-)$.
\begin{equation}\label{Xtwonaturality}
\begin{tikzcd}[column sep=large]
Xa \times Xa' \ar{d}[swap]{f_* \times f'_*} \ar{r}{X^2_{a,a'}} & X(a \oplus a') \ar{d}{(f \oplus f')_*}\\
Xb \times Xb' \ar{r}{X^2_{b,b'}} & X(b \oplus b')
\end{tikzcd}
\end{equation}
\item[Unity] The unity axiom \cref{monoidalfunctorunity} for $(X,X^2,X^0)$ states that, for objects $a \in \cD$ and $w \in Xa$ and morphisms $f \in Xa$, the following equalities hold in $Xa$.
\begin{equation}\label{asmfunity}
\left\{\begin{split}
X^2_{\zero,a} \left(X^0*, w\right) &= w = X^2_{a,\zero} \left(w, X^0*\right)\\
X^2_{\zero,a} \left(1_{X^0*}, f\right) &= f = X^2_{a,\zero} \left(f, 1_{X^0*}\right)
\end{split}\right.
\end{equation}
\item[Associativity] The associativity axiom \cref{monoidalfunctorassoc} for $(X,X^2,X^0)$ states that, for objects $a,b,c \in \cD$ and 
\[(x,y,z) \in Xa \times Xb \times Xc\]
either all objects or all morphisms, the equality 
\begin{equation}\label{asmfassociativity}
X^2_{a \oplus b, c} \left( X^2_{a,b}(x,y), z\right) = X^2_{a,b \oplus c} \left(x, X^2_{b,c}(y,z)\right)
\end{equation}
holds in $X(a \oplus b \oplus c)$.
\item[Braiding] The braiding axiom \cref{monoidalfunctorbraiding} for $(X,X^2,X^0)$ states the equality 
\begin{equation}\label{asmfbraiding}
\betaplus_* X^2_{a,b}(x,y) = X^2_{b,a}(y,x)
\end{equation}
in $X(b \oplus a)$ with
\begin{itemize}
\item $\betaplus \cn a \oplus b \to b \oplus a$ the braiding in $\cD$ and
\item $\betaplus_* = X\betaplus \cn X(a \oplus b) \to X(b \oplus a)$ its image under $X$.\defmark
\end{itemize} 
\end{description}
\end{explanation}

\subsection*{Multimorphism Categories: Additive Natural Transformations}
Next we define the objects in the multimorphism categories in $\DCat$.  These are subcategories of the corresponding multimorphism categories in $\Dtecat$ (\cref{expl:permindexedcat}), with $\Dte = (\cD, \otimes, \tu, \betate)$ the \emph{multiplicative} structure of $\cD$ in \cref{Dmultiplicate}.  Recall from \cref{dcatnary} that $(\Dtecat)\smscmap{\angX;Z}$ has
\begin{itemize}
\item natural transformations $\bigotimes_{j=1}^n X_j \to Z$ as objects \cref{phitensorxz},
\item modifications as morphisms \cref{Phiphivarphi}, and
\item identities and composition given by those of modifications under vertical composition.
\end{itemize} 
In each multimorphism category in $\DCat$, the objects are natural transformations in $\Dtecat$ satisfying two compatibility conditions, which are \cref{additivenattrunity,additivenattradditivity} below, with respect to the \emph{additive} structures of the additive symmetric monoidal functors in the input and the output.

\begin{definition}\label{def:additivenattr}
Suppose given additive symmetric monoidal functors 
\[(Z,Z^2,Z^0) \andspace \bang{(X_j, X_j^2, X_j^0)}_{j=1}^n \cn \Dplus \to \Cat\]
with $n > 0$ as in \cref{def:additivesmf}.  An \emph{additive natural transformation}\index{additive!natural transformation}\index{natural transformation!additive}
\begin{equation}\label{additivenattransformation}
\begin{tikzcd}[column sep=large]
\bang{(X_j, X_j^2, X_j^0)}_{j=1}^n \ar{r}{\phi} & (Z,Z^2,Z^0)
\end{tikzcd}
\end{equation}
is a natural transformation in $\Dtecat$ as in \cref{phitensorxz}
\begin{equation}\label{phixjz}
\begin{tikzcd}[column sep=large]
\bigotimes_{j=1}^n X_j \ar{r}{\phi} & Z
\end{tikzcd}
\end{equation}
that satisfies the following unity axiom \cref{additivenattrunity} and additivity axiom \cref{additivenattradditivity}. 
\begin{description}
\item[Unity] For each $i \in \{1,\ldots,n\}$ and objects $a_j \in \cD$ for $i \neq j \in \{1,\ldots,n\}$, consider the component functor of $\phi$ as in \cref{dcatnaryobjcomponent}
\begin{equation}\label{phiaizero}
\begin{tikzcd}[column sep=2cm]
X_1 a_1 \times \cdots \times X_i\zero \times \cdots \times X_n a_n \ar{r}{\phi_{(a_1,\ldots,\zero,\ldots,a_n)}} & Z\zero
\end{tikzcd}
\end{equation}
with
\begin{itemize}
\item $i$-th domain factor $X_i\zero$ and
\item $j$-th domain factor $X_j a_j$ for $j \neq i$.
\end{itemize}   
For each $j \neq i$, suppose $x_j \in X_j a_j$ is an object, and $f_j \in X_j a_j$ is a morphism.  Then the following object, respectively, morphism, equalities are required to hold in $Z\zero$, where $\phi = \phi_{(a_1,\ldots,\zero,\ldots,a_n)}$.\index{unity axiom!additive natural transformation}
\begin{equation}\label{additivenattrunity}
\left\{\begin{split}
\phi\big(x_1, \ldots, X_i^0*, \ldots, x_n\big) &= Z^0*\\
\phi\big(f_1, \ldots, 1_{X_i^0*}, \ldots, f_n\big) &= 1_{Z^0*}
\end{split}\right.
\end{equation}
Here $X_i^0* \in X_i\zero$ and $Z^0* \in Z\zero$ are the unit objects \cref{additivesmfunitobj} of, respectively, $X_i$ and $Z$.
\item[Additivity]
Suppose given objects $\anga = \ang{a_j}_{j=1}^n \in \cD^n$ and $a_i' \in \cD$ for some $i \in \{1,\ldots,n\}$, along with the following notation.
\begin{equation}\label{additivenattrnotation}
\left\{\scalebox{.9}{$
\begin{aligned}
a &= a_1 \otimes \cdots \otimes a_i \otimes \cdots \otimes a_n & a_i'' &= a_i \oplus a_i'\\
a' &= a_1 \otimes \cdots \otimes a_i' \otimes \cdots \otimes a_n & \angap &= (a_1, \ldots, a_i',\ldots,a_n)\\
a'' &= a_1 \otimes \cdots \otimes a_i'' \otimes \cdots \otimes a_n & \angapp &= (a_1, \ldots, a_i'', \ldots, a_n)\\
X\anga &= X_1 a_1 \times \cdots \times X_i a_i \times \cdots \times X_n a_n &&\\
X\angap &= X_1 a_1 \times \cdots \times X_i a_i' \times \cdots \times X_n a_n &&\\
X\angapp &= X_1 a_1 \times \cdots \times X_i a_i'' \times \cdots \times X_n a_n &&
\end{aligned}$}\right.
\end{equation}
\begin{itemize}
\item In each of $a$, $a'$, and $a''$, the $j$-th factor is $a_j$ for $j \neq i$, and the $i$-th factor is, respectively, $a_i$, $a_i'$, and $a_i''$.
\item In each of $\angap$ and $\angapp$, the $j$-th entry is $a_j$ for $j \neq i$, and the $i$-th entry is, respectively, $a_i'$ and $a_i''$.
\item In each of $X\anga$, $X\angap$, and $X\angapp$, the $j$-th factor is $X_j a_j$ for $j \neq i$, and the $i$-th factor is, respectively, $X_i a_i$, $X_i a_i'$, and $X_i a_i''$.
\end{itemize}
Then the following diagram of functors is required to commute.\index{additivity axiom!additive natural transformation}
\begin{equation}\label{additivenattradditivity}
\begin{tikzpicture}[xscale=2.5,yscale=1,vcenter]
\def\v{-1} \def\u{-1.5}
\draw[0cell=.85]
(0,0) node (x1) {X_1 a_1 \times \cdots \times (X_i a_i \times X_i a_i') \times \cdots \times X_n a_n}
(x1)++(-1,\v) node (x21) {X\angapp}
(x1)++(1,\v) node (x22) {X\anga \times X\angap}
(x21)++(0,\u) node (x31) {Za''}
(x22)++(0,\u) node (x32) {Za \times Za'} 
(x31)++(1,\v) node (x4) {Z(a \oplus a')}
;
\draw[1cell=.85]  
(x1) edge node[swap,pos=.9] {1 \times \cdots \times (X_i^2)_{a_i, a_i'} \times \cdots \times 1} (x21)
(x21) edge node[swap] {\phi_{\angapp}} (x31)
(x31) edge node[swap,pos=.5] {Z(\lap)} node {\iso} (x4)
(x1) edge node[pos=.7] {\varrho} (x22)
(x22) edge node {\phi_{\anga} \times \phi_{\angap}} (x32)
(x32) edge node[pos=.3] {Z^2_{a,a'}} (x4)
;
\end{tikzpicture}
\end{equation}
The detail of \cref{additivenattradditivity} is as follows.
\begin{itemize}
\item At the top node of \cref{additivenattradditivity}, the $i$-th factor is $X_i a_i \times X_i a_i'$, and the $j$-th factor is $X_j a_j$ for $j \neq i$.
\item $Z^2_{a,a'}$ and
\[\begin{tikzcd}[column sep=huge]
X_i a_i \times X_i a_i' \ar{r}{(X_i^2)_{a_i, a_i'}} & X_i(a_i \oplus a_i') = X_i a_i''
\end{tikzcd}\]
are the indicated component functors of the monoidal constraints \cref{additivesmfconstraint} of, respectively, $Z$ and $X_i$.
\item Each of $\phi_{\anga}$, $\phi_{\angap}$, and $\phi_{\angapp}$ is the indicated component functor \cref{dcatnaryobjcomponent} of $\phi$.
\item Omitting the $\otimes$ symbol to save space, the morphism
\begin{equation}\label{laplazaapp}
\begin{tikzcd}[column sep=tiny,row sep=tiny]
a'' \ar[equal, shorten >=-.5ex]{d} \ar{r}{\lap}[swap]{\iso} & 
a \oplus a' \ar[equal, shorten >=-.5ex]{d}\\
a_1 \cdots (a_i \oplus a_i') \cdots a_n & (a_1 \cdots a_i \cdots a_n) \oplus (a_1 \cdots a_i' \cdots a_n)
\end{tikzcd}
\end{equation}
in \cref{additivenattradditivity} is the unique Laplaza coherence isomorphism in the tight bipermutative category $\cD$ (\cref{ex:laplazacoherence}) that distributes over the sum.
\item $\varrho$ in \cref{additivenattradditivity} is the following composite functor.
\begin{equation}\label{varrhofactor}
\begin{tikzpicture}[xscale=2.5,yscale=1,vcenter]
\def\v{-1} \def\u{-1.5}
\draw[0cell=.85]
(0,0) node (x1) {X_1 a_1 \times \cdots \times (X_i a_i \times X_i a_i') \times \cdots \times X_n a_n}
(x1)++(0,-2) node (x21) {(X_1 a_1)^2 \times \cdots \times (X_i a_i \times X_i a_i') \times \cdots \times (X_n a_n)^2}
(x1)++(1,\v) node (x22) {X\anga \times X\angap}
;
\draw[1cell=.85]  
(x1) edge node[swap] {\varrho_1} (x21)
(x21) edge[transform canvas={xshift=1ex}] node[swap,pos=.7] {\iso} node[pos=.4] {\varrho_2} (x22)
(x1) edge[transform canvas={xshift=1ex}] node[pos=.7] {\varrho} (x22)
;
\end{tikzpicture}
\end{equation}
\begin{itemize}
\item $\varrho_1$ is the product of
\begin{enumerate}[label=(\roman*)]
\item the diagonal functors $X_j a_j \to (X_j a_j)^2$ for $j\neq i$ and
\item the identity functor of $X_i a_i \times X_i a_i'$.
\end{enumerate} 
\item $\varrho_2$ permutes the $2n$ factors according to the transpose permutation $\xitimes_{2,n} \in \Sigma_{2n}$ in \cref{fset-xitimes}.
\end{itemize}
\end{itemize}
\end{description}
This finishes the definition of an additive natural transformation.  An explicit form of the additivity axiom \cref{additivenattradditivity} is in \cref{additivityobjects}.
\end{definition}

\begin{explanation}[Additive Natural Transformations]\label{expl:additivenattr}
Consider \cref{def:additivenattr}.
\begin{enumerate}
\item\label{expl:additivenattr-i} An additive natural transformation is a natural transformation in $\Dtecat$ with two additional properties---the unity axiom \cref{additivenattrunity} and the additivity axiom \cref{additivenattradditivity}---but \emph{no} extra structure.
\item\label{expl:additivenattr-ii} The domain of the natural transformation $\phi \cn \bigotimes_{j=1}^n X_j \to Z$ in \cref{{phixjz}} is the iterated Day convolution $\bigotimes_{j=1}^n X_j$ in \cref{iteratedday}.  It is constructed using the \emph{multiplicative} structure $\Dte = (\cD, \otimes, \tu, \betate)$ and the underlying functors $X_j \cn \cD \to \Cat$.
\item\label{expl:additivenattr-iii} The component functor $\phi$ in \cref{phiaizero} is well defined because the iterated product in $\cD$
\[\overbracket[.5pt]{a_1 \otimes \cdots \otimes a_{i-1} \,\otimes}^{\text{empty if $i=1$}} \zero \overbracket[.5pt]{\otimes\, a_{i+1} \otimes \cdots \otimes a_n}^{\text{empty if $i=n$}} = \zero\]
by the multiplicative zero axiom \cref{ringcataxiommultzero} of a ring, hence also bipermutative, category.
\item\label{expl:additivenattr-iv} In the case $n=1$, an additive natural transformation 
\[\begin{tikzcd}[column sep=large]
X_1 \ar{r}{\phi} & Z
\end{tikzcd}\]
is precisely a monoidal natural transformation (\cref{def:monoidalnattr}).  The unity axiom \cref{additivenattrunity} and the additivity axiom \cref{additivenattradditivity} reduce to the diagrams in \cref{monnattr} when $n=1$.
\item\label{expl:additivenattr-v} In \cref{varrhofactor} $\varrho_2$ is the unique coherence isomorphism in the symmetric monoidal category $\Cat$ that permutes the factors to $X\anga \times X\angap$ in such a way that the relative order of the two factors in
\begin{itemize}
\item each $(X_j a_j)^2$ for $j \neq i$ and
\item $X_i a_i \times X_i a_i'$
\end{itemize}  
remain unchanged.
\item\label{expl:additivenattr-vi} The additivity axiom \cref{additivenattradditivity} is given explicitly as follows.  Suppose given 
\begin{equation}\label{angxxiprime}
\angx = \ang{x_j}_{j=1}^n \in \txprod_{j=1}^n X_j a_j \andspace x_i' \in X_i a_i'
\end{equation}
that are either all objects or all morphisms, along with the following notation.
\begin{equation}\label{angxprime}
\left\{\begin{aligned}
x_i'' &= (X_i^2)_{a_i, a_i'}(x_i, x_i') \in X_i a_i'' = X_i(a_i \oplus a_i')\\
\angxp &= (x_1, \ldots, x_i', \ldots, x_n) \in X\angap\\
\angxpp &= (x_1, \ldots, x_i'', \ldots, x_n) \in X\angapp
\end{aligned}\right.
\end{equation}
Then the additivity axiom \cref{additivenattradditivity} is the equality  
\begin{equation}\label{additivityobjects}
\lap_* \phi_{\angapp} \angxpp = Z^2_{a,a'} \brb{\phi_{\anga}\angx, \phi_{\angap}\angxp}
\end{equation}
in $Z(a \oplus a')$, with $\lap_* = Z(\lap)$.\defmark
\end{enumerate}
\end{explanation}

\subsection*{Multimorphism Categories: Additive Modifications}
Next we define the morphisms in the multimorphism categories in $\DCat$.  

\begin{definition}\label{def:additivemodification}
In the context of \cref{def:additivenattr}, suppose $\phi$ and $\psi$ as in 
\[\begin{tikzcd}[column sep=large]
\bang{(X_j, X_j^2, X_j^0)}_{j=1}^n \ar[shift left]{r}{\phi} \ar[shift right]{r}[swap]{\psi} & (Z,Z^2,Z^0)
\end{tikzcd}\]
are additive natural transformations.  An \emph{additive modification}\index{additive!modification}\index{modification!additive} $\Phi \cn \phi \to \psi$, which is also denoted 
\begin{equation}\label{additivemodtwocell}
\begin{tikzpicture}[xscale=1,yscale=1,baseline={(x1.base)}]
\draw[0cell=.9]
(0,0) node (x1) {\bang{(X_j, X_j^2, X_j^0)}_{j=1}^n}
(x1)++(1,0) node (x2) {\phantom{X}}
(x2)++(2,0) node (x3) {\phantom{X}}
(x3)++(.65,0) node (x4) {(Z,Z^2,Z^0),}
;
\draw[1cell=.9]  
(x2) edge[bend left] node[pos=.5] {\phi} (x3)
(x2) edge[bend right] node[swap,pos=.5] {\psi} (x3)
;
\draw[2cell]
node[between=x2 and x3 at .45, rotate=-90, 2label={above,\Phi}] {\Rightarrow}
;
\end{tikzpicture}
\end{equation}
is a modification
\[\begin{tikzpicture}[xscale=2,yscale=1.7,baseline={(x1.base)}]
\draw[0cell=.9]
(0,0) node (x1) {\textstyle \bigotimes_{j=1}^n X_j}
(x1)++(.15,.03) node (x2) {\phantom{\textstyle \bigotimes Z}}
(x2)++(1,0) node (x3) {Z}
;
\draw[1cell=.9]  
(x2) edge[bend left, shorten <=-.5ex] node[pos=.45] {\phi} (x3)
(x2) edge[bend right, shorten <=-.5ex] node[swap,pos=.45] {\psi} (x3)
;
\draw[2cell]
node[between=x2 and x3 at .45, rotate=-90, 2label={above,\Phi}] {\Rightarrow}
;
\end{tikzpicture}\]
as in \cref{Phiphivarphi} that satisfies the following unity and additivity axioms.
\begin{description}
\item[Unity] In the context of \cref{phiaizero,additivenattrunity}, the component morphism of $\Phi$ in $Z\zero$\index{unity axiom!additive modification}
\begin{equation}\label{additivemodunity}
\begin{tikzpicture}[xscale=2.2,yscale=1.4,vcenter]
\draw[0cell=.9]
(0,0) node (x1) {\phi(x_1, \ldots, X_i^0*, \ldots, x_n) = Z^0*}
(x1)++(0,-1) node (x2) {\psi(x_1, \ldots, X_i^0*, \ldots, x_n) = Z^0*}
;
\draw[1cell=.9]  
(x1) edge[transform canvas={xshift=-3ex}] node {\Phi_{(x_1, \ldots, X_i^0*, \ldots, x_n)}} (x2)
;
\end{tikzpicture}
\end{equation}
is the identity morphism $1_{Z^0*}$.
\item[Additivity] In the context of \cref{additivenattrnotation,additivenattradditivity}, the following two whiskered natural transformations are equal, where each $\Phi_?$ is the indicated component natural transformation of $\Phi$ as in \cref{dcatnarymodcomponent}.\index{additivity axiom!additive modification}
\begin{equation}\label{additivemodadditivity}
\begin{tikzpicture}[xscale=2.5,yscale=1,vcenter]
\def\a{20} \def\v{-1} \def\u{-1.5}
\draw[0cell=.85]
(0,0) node (x1) {X_1 a_1 \times \cdots \times (X_i a_i \times X_i a_i') \times \cdots \times X_n a_n}
(x1)++(-1,\v) node (x21) {X\angapp}
(x1)++(1,\v) node (x22) {X\anga \times X\angap}
(x21)++(0,\u) node (x31) {Za''}
(x22)++(0,\u) node (x32) {Za \times Za'} 
(x31)++(1,\v) node (x4) {Z(a \oplus a')}
;
\draw[1cell=.85]  
(x1) edge node[swap,pos=.9] {1 \times \cdots \times (X_i^2)_{a_i, a_i'} \times \cdots \times 1} (x21)
(x21) edge[bend right=\a, transform canvas={xshift=-1.5ex}] node[swap,pos=.4] {\phi_{\angapp}} (x31)
(x21) edge[bend left=\a, transform canvas={xshift=1.5ex}] node[pos=.6] {\psi_{\angapp}} (x31)
(x31) edge node[swap,pos=.5] {Z(\lap)} node {\iso} (x4)
(x1) edge node[pos=.7] {\varrho} (x22)
(x22) edge[bend right=\a, transform canvas={xshift=-3ex}] node[swap,pos=.4] {\phi_{\anga} \times \phi_{\angap}} (x32)
(x22) edge[bend left=\a, transform canvas={xshift=3ex}] node[pos=.6] {\psi_{\anga} \times \psi_{\angap}} (x32)
(x32) edge node[pos=.3] {Z^2_{a,a'}} (x4)
;
\draw[2cell=.85]
node[between=x21 and x31 at .65, 2label={above,\Phi_{\angapp}}] {\Longrightarrow}
node[between=x22 and x32 at .65, 2label={above,\Phi_{\anga} \times \Phi_{\angap}}] {\Longrightarrow}
;
\end{tikzpicture}
\end{equation}
\end{description}
This finishes the definition of an additive modification.
\end{definition}

\begin{explanation}[Additive Modifications]\label{expl:additivemodification}
Consider \cref{def:additivemodification}.
\begin{enumerate}
\item\label{expl:additivemodification-i} The morphism in the unity axiom \cref{additivemodunity} for $\Phi$ is the component of the component natural transformation of $\Phi$
\[\begin{tikzpicture}[xscale=3.5,yscale=2.3,baseline={(x1.base)}]
\draw[0cell=.9]
(0,0) node (x1) {X_1 a_1 \times \cdots \times X_i\zero \times \cdots \times X_n a_n}
(x1)++(.47,.01) node (x2) {\phantom{Z\zero}}
(x2)++(1,0) node (x3) {Z\zero}
;
\draw[1cell=.9]  
(x2) edge[bend left] node[pos=.54] {\phi_{(a_1,\ldots,\zero,\ldots,a_n)}} (x3)
(x2) edge[bend right] node[swap,pos=.54] {\psi_{(a_1,\ldots,\zero,\ldots,a_n)}} (x3)
;
\draw[2cell=.9]
node[between=x2 and x3 at .3, rotate=-90, 2label={above,\Phi_{(a_1,\ldots,\zero,\ldots,a_n)}}] {\Rightarrow}
;
\end{tikzpicture}\]
at the object $(x_1, \ldots, X_i^0*, \ldots, x_n)$.
\item\label{expl:additivemodification-ii}
In the context of \cref{angxprime,additivityobjects} with objects $x_j \in X_j a_j$ and $x_i' \in X_i a_i'$, the additivity axiom \cref{additivemodadditivity} for $\Phi$ is the morphism equality
\begin{equation}\label{additivitymod}
\lap_* \big(\Phi_{\angapp}\big)_{\angxpp} = Z^2_{a,a'} \brb{(\Phi_{\anga})_{\angx}, (\Phi_{\angap})_{\angxp}}
\end{equation}
in $Z(a \oplus a')$.\defmark
\end{enumerate}
\end{explanation}

\subsection*{Multimorphism Categories}
Next we define the multimorphism categories in $\DCat$ in positive arity and check that they are well defined.  Recall the multimorphism category $\Dtecat\smscmap{\angX;Z}$ in \cref{dcatnary}.

\begin{definition}\label{def:dcatnarycategory}
In the context of \cref{def:additivenattr} with $\angX = \ang{X_j}_{j=1}^n$, define the data of a subcategory
\begin{equation}\label{dcatangxz}
\DCat\scmap{\angX;Z}
\end{equation}
of $\Dtecat\smscmap{\angX;Z}$ as follows.
\begin{description}
\item[Objects] An object is an \emph{additive} natural transformation as in \cref{additivenattransformation}.
\item[Morphisms] A morphism is an \emph{additive} modification as in \cref{additivemodtwocell}.
\item[Identities] Identity morphisms are identity modifications.
\item[Composition] It is the vertical composition of modifications \cref{modificationvcomp}.
\end{description} 
This finishes the definition of $\DCat\smscmap{\angX;Z}$.
\end{definition}

\begin{lemma}\label{dcatnarycategory}
In the context of \cref{def:dcatnarycategory}, $\DCat\smscmap{\angX;Z}$ defines a subcategory of $\Dtecat\smscmap{\angX;Z}$ in \cref{dcatnary}.
\end{lemma}

\begin{proof}
Modifications are associative and unital with respect to vertical composition \cite[4.4.8 (2)]{johnson-yau}.  For an identity modification, the unity axiom \cref{additivemodunity} and the additivity axiom \cref{additivemodadditivity} reduce to those for an additive natural transformation, namely, \cref{additivenattrunity} and \cref{additivenattradditivity}.  Thus it remains to check that additive modifications are closed under vertical composition.

Suppose given additive modifications $\Phi$ and $\Psi$ as follows.
\[\begin{tikzpicture}[xscale=1,yscale=1,baseline={(x1.base)}]
\def\a{50}
\draw[0cell=.9]
(0,0) node (x1) {\bang{(X_j, X_j^2, X_j^0)}_{j=1}^n}
(x1)++(1.1,.03) node (x2) {\phantom{X}}
(x2)++(2,0) node (x3) {\phantom{X}}
(x3)++(.65,0) node (x4) {(Z,Z^2,Z^0)}
;
\draw[1cell=.9]  
(x2) edge[bend left=\a] node[pos=.5] {\phi} (x3)
(x2) edge node[pos=.25] {\psi} (x3)
(x2) edge[bend right=\a] node[swap,pos=.5] {\omega} (x3)
;
\draw[2cell=.9]
node[between=x2 and x3 at .48, shift={(0,.3)}, rotate=-90, 2label={above,\Phi}] {\Rightarrow}
node[between=x2 and x3 at .48, shift={(0,-.3)}, rotate=-90, 2label={above,\Psi}] {\Rightarrow}
;
\end{tikzpicture}\]
We must show that the vertical composite modification $\Psi\Phi \cn \phi \to \omega$ is additive.  The vertical composite $\Psi\Phi$ satisfies the unity axiom \cref{additivemodunity} by the morphism equality
\[1_{Z^0*} 1_{Z^0*} = 1_{Z^0*} \inspace Z\zero.\]
The following computation shows that the vertical composite $\Psi\Phi$ satisfies the additivity axiom \cref{additivitymod}.
\[\begin{split}
&\lap_*\big((\Psi\Phi)_{\angapp}\big)_{\angxpp}\\
&= \lap_* \Big(\big(\Psi_{\angapp}\big)_{\angxpp} \big(\Phi_{\angapp}\big)_{\angxpp}\Big)\\
&= \Big(\lap_* \big(\Psi_{\angapp}\big)_{\angxpp}\Big) \Big(\lap_* \big(\Phi_{\angapp}\big)_{\angxpp}\Big)\\
&= \Big(Z^2_{a,a'} \brb{(\Psi_{\anga})_{\angx}, (\Psi_{\angap})_{\angxp}}\Big) 
\Big(Z^2_{a,a'} \brb{(\Phi_{\anga})_{\angx}, (\Phi_{\angap})_{\angxp}}\Big)\\
&= Z^2_{a,a'} \brb{(\Psi_{\anga})_{\angx} (\Phi_{\anga})_{\angx}, (\Psi_{\angap})_{\angxp}(\Phi_{\angap})_{\angxp}}\\
&= Z^2_{a,a'} \brb{((\Psi\Phi)_{\anga})_{\angx}, ((\Psi\Phi)_{\angap})_{\angxp}}
\end{split}\]
The above equalities hold for the following reasons:
\begin{itemize}
\item The first and the last equalities follow from the definitions of the vertical composite modification $\Psi\Phi$.
\item The second equality follows from the functoriality of $\lap_* = Z(\lap)$.
\item The third equality follows from the additivity axiom \cref{additivitymod} for $\Phi$ and $\Psi$.
\item The fourth equality follows from the functoriality of $Z^2_{a,a'}$.
\end{itemize}
This proves that the vertical composite $\Psi\Phi$ is an additive modification.
\end{proof}

\section{Enriched Multicategories of Bipermutative-Indexed Categories}
\label{sec:diagrambipermutative}

We continue to assume that $\cD$ is a small tight bipermutative category as in \cref{Dbipermutative}.  In this section we construct the $\Cat$-enriched multicategory $\DCat$ that incorporates the bipermutative structure of $\cD$ (\cref{thm:dcatcatmulticat}).  While the $\Cat$-multicategory structure on $\DCat$ is inherited from that on $\Dtecat$ in \cref{expl:permindexedcat}, there are several subtleties that one should be aware of:
\begin{enumerate}
\item The $\Cat$-multicategory structure on $\Dtecat$ is induced by a symmetric monoidal closed structure (\cref{thm:permindexedcat}).  This is \emph{not} true for the $\Cat$-multicategory $\DCat$.  We discuss this point in \cref{expl:dcatbipermzero}.
\item Additive natural transformations (\cref{def:additivenattr}) are natural transformations with extra properties.  So we have to check that the $\Cat$-multicategory structure on $\Dtecat$ preserves additive natural transformations.  The same remark also applies to additive modifications (\cref{def:additivemodification}), which are modifications with extra properties.
\item $\DCat$ is, in general, \emph{not} a sub-$\Cat$-multicategory of $\Dtecat$ because the $\Cat$-multifunctor 
\[\DCat \to \Dtecat\]
is not injective on objects; see \cref{expl:dcatforgetdtecat}.
\end{enumerate}

\subsection*{Definitions}

Now we define the $\Cat$-multicategory structure on $\DCat$.

\begin{definition}\label{def:dcatcatmulticat}
For a small tight bipermutative category $\cD$ as in \cref{Dbipermutative}, define the data of a $\Cat$-multicategory $\DCat$ as follows.
\begin{description}
\item[Objects] 
These are additive symmetric monoidal functors $\Dplus \to \Cat$ (\cref{def:additivesmf}).
\item[Multimorphism Categories in Arity 0]
For each additive symmetric monoidal functor
\[(Z,Z^2,Z^0) \cn \Dplus \to \Cat,\]
define the 0-ary multimorphism category
\begin{equation}\label{dcatbipermzeroary}
\DCat\scmap{\ang{};Z} = Z\tu.
\end{equation}
\item[Multimorphism Categories in Positive Arity]
Suppose given additive symmetric monoidal functors 
\begin{equation}\label{zxjadditivesmf}
(Z,Z^2,Z^0) \andspace \bang{(X_j, X_j^2, X_j^0)}_{j=1}^n \cn \Dplus \to \Cat
\end{equation}
with $n > 0$ and $\angX = \ang{X_j}_{j=1}^n$.  Define the $n$-ary multimorphism category
\[\DCat\scmap{\angX;Z}\]
to be the category in \cref{dcatangxz}.
\item[Colored Units]
For an additive symmetric monoidal functor $(Z,Z^2,Z^0)$, the $Z$-colored unit is the identity natural transformation $1_Z \cn Z \to Z$.
\item[Symmetric Group Action and Composition]
These are restrictions of those in the $\Cat$-multicategory $\Dtecat$---namely, \cref{dcatsigmaaction,dcatgamma}---using
\begin{itemize}
\item \cref{dcatzeroary,dcatbipermzeroary} to identify multimorphism categories with arity 0 input and
\item \cref{dcatnarycategory} for multimorphism categories with positive arity input.
\end{itemize} 
\end{description}
This finishes the definition of $\DCat$.
\end{definition}

\begin{explanation}[$\DCat$ is Not Symmetric Monoidal]\label{expl:dcatbipermzero}\index{diagram category!not symmetric monoidal}
Here we discuss a subtle difference between $\Dtecat$ and $\DCat$.  On the one hand, the $\Cat$-multicategory structure on $\Dtecat$ is induced by the symmetric monoidal closed structure on $\Dtecat$ (\cref{expl:permindexedcat}).  Its monoidal unit is the unit diagram \cref{dcatunit}
\[J = \cD(\tu,-) \cn \cD \to \Cat.\]
By the definition of the induced $\Cat$-multicategory, for each $\cD$-indexed category $Z \cn \cD \to \Cat$, the 0-ary multimorphism category in $\Dtecat$ is 
\[(\Dtecat)\scmap{\ang{};Z} = (\Dcat)\srb{J,Z} \iso Z\tu,\]
with the canonical isomorphism given by the Bicategorical Yoneda Lemma.  

On the other hand, in \cref{dcatbipermzeroary} we directly define the 0-ary multimorphism category $\DCat\smscmap{\ang{};Z}$ to be $Z\tu$, without using $J$.  The reason is that, in general, $J = \cD(\tu,-)$ does \emph{not} admit the structure of a symmetric monoidal functor $\Dplus \to \Cat$.  Therefore, the $\Cat$-multicategory $\DCat$, which we will verify shortly, is, in general, \emph{not} induced by a restriction of the symmetric monoidal structure on $\Dtecat$.

In more detail, by definition a monoidal functor structure on $J$ has a unit constraint, which is a choice of a unit object \cref{additivesmfunitobj}
\[J^0* \in J\zero = \cD(\tu,\zero).\]
For a general bipermutative category $\cD$, the morphism set $\cD(\tu,\zero)$ may be empty, so such a unit object $J^0*$ does not exist.  Here are some concrete examples of small tight bipermutative categories where this morphism set is empty: 
\begin{enumerate}
\item In $\Finsk$ (\cref{ex:finsk}) there is an empty morphism set
\[\Finsk(\ufs{1},\ufs{0}) = \Set\brb{\{1\},\emptyset} = \emptyset.\]
\item In $\Fset$ (\cref{ex:Fset}) there is an empty morphism set
\[\Fset(1,0) = \emptyset.\]
\item In $\cA$ (\cref{ex:mandellcategory}) there is an empty morphism set
\[\cA\brb{(1),\ang{}} = \emptyset\]
because there are no functions $\ufs{1} = \{1\} \to \emptyset$.
\end{enumerate}
In each of these examples, $J$ cannot have a unit object and, therefore, does not admit the structure of a monoidal functor $\Dplus \to \Cat$.
\end{explanation}

\subsection*{Symmetric Group Action}

To show that $\DCat$ as in \cref{def:dcatcatmulticat} is a $\Cat$-multicategory, we first check that the symmetric group action, which is a restriction of the symmetric group action in $\Dtecat$, is well defined.

\begin{lemma}\label{dcatsigmaclosed}
Additive natural transformations and additive modifications are closed under the symmetric group action \cref{dcatsigmaaction} in $\Dtecat$.
\end{lemma}

\begin{proof}
Suppose given additive symmetric monoidal functors $Z$ and $X_j$ as in \cref{zxjadditivesmf} with $n>0$ and $\angX = \ang{X_j}_{j=1}^n$.  Suppose $\sigma \in \Sigma_n$ is a permutation.  Given an additive natural transformation $\phi \cn \bigotimes_{j=1}^n X_j \to Z$ as in \cref{phixjz}, we must show that the natural transformation
\[\begin{tikzcd}[column sep=large]
\bigotimes_{j=1}^n X_{\sigma(j)} \ar{r}{\phi^\sigma} & Z
\end{tikzcd}\]
satisfies the unity axiom \cref{additivenattrunity} and the additivity axiom \cref{additivenattradditivity} of an additive natural transformation.  By definition \cref{phisigmaanga}, for objects $\anga = \ang{a_j}_{j=1}^n$ in $\cD$, the component functor $\phi^\sigma_{\anga}$ is the following composite.
\[\begin{tikzpicture}[xscale=3.5,yscale=1.5,vcenter]
\draw[0cell=.9]
(0,0) node (x11) {\textstyle\prod_{j=1}^n X_{\sigma(j)} a_j}
(x11)++(1,0) node (x12) {Za}
(x11)++(0,-1) node (x21) {\textstyle\prod_{j=1}^n X_j a_{\sigmainv(j)}}
(x12)++(0,-1) node (x22) {Za_{\sigmainv}}
;
\draw[1cell=.9]  
(x11) edge node {\phi^\sigma_{\anga}} (x12)
(x11) edge node {\iso} node[swap] {\sigma} (x21)
(x21) edge node {\phi_{\sigma\anga}} (x22)
(x22) edge node {\iso} node[swap] {Z(\sigmainv)} (x12)
;
\end{tikzpicture}\]

For the unity axiom \cref{additivenattrunity} for $\phi^\sigma$, suppose $a_i = \zero$ for some index $i \in \{1,\ldots,n\}$.  The two equalities in the unity axiom hold for $\phi^\sigma_{\anga}$ by
\begin{itemize}
\item the unity axiom \cref{additivenattrunity} for $\phi_{\sigma\anga}$;
\item the uniqueness of Laplaza coherence isomorphisms (\cref{thm:laplaza-coherence-1}), which implies that the coherence isomorphism
\[\begin{tikzcd}[column sep=large]
a_{\sigmainv} = \bigotimes_{j=1}^n a_{\sigmainv(j)} = \zero \ar{r}{\sigmainv}[swap]{\iso} & 
\zero = \bigotimes_{j=1}^n a_j = a
\end{tikzcd}\]
is the identity morphism $1_\zero$ in $\cD$; and
\item the functoriality of $Z$, which implies that $Z(1_\zero) = 1_{Z\zero}$.
\end{itemize}

For the additivity axiom \cref{additivenattradditivity} for $\phi^\sigma$, in addition to the notation in \cref{sigmaangaasigmainv,additivenattrnotation}, we use the following notation. 
\begin{equation}\label{sigmaangaprime}
\left\{\begin{aligned}
\sigma\angap &= (a_{\sigmainv(1)}, \ldots, a_i', \ldots, a_{\sigmainv(n)}) \in \cD^n\\
a'_{\sigmainv} &= a_{\sigmainv(1)} \otimes \cdots \otimes a_i' \otimes \cdots \otimes a_{\sigmainv(n)} \in \cD\\
\sigma\angapp &= (a_{\sigmainv(1)}, \ldots, a_i'', \ldots, a_{\sigmainv(n)}) \in \cD^n\\
a''_{\sigmainv} &= a_{\sigmainv(1)} \otimes \cdots \otimes a_i'' \otimes \cdots \otimes a_{\sigmainv(n)} \in \cD
\end{aligned}\right.
\end{equation}
\[\left\{\begin{aligned}
X\sigma\anga &= X_1 a_{\sigmainv(1)} \times \cdots \times X_{\sigma(i)} a_i \times \cdots \times X_n a_{\sigmainv(n)}\\
X\sigma\angap &= X_1 a_{\sigmainv(1)} \times \cdots \times X_{\sigma(i)} a_i' \times \cdots \times X_n a_{\sigmainv(n)}\\
X\sigma\angapp &= X_1 a_{\sigmainv(1)} \times \cdots \times X_{\sigma(i)} a_i'' \times \cdots \times X_n a_{\sigmainv(n)}
\end{aligned}\right.\]
\[\left\{\begin{aligned}
X^\sigma\anga &= X_{\sigma(1)} a_1 \times \cdots \times X_{\sigma(i)} a_i \times \cdots \times X_{\sigma(n)} a_n\\
X^\sigma\angap &= X_{\sigma(1)} a_1 \times \cdots \times X_{\sigma(i)} a_i' \times \cdots \times X_{\sigma(n)} a_n\\
X^\sigma\angapp &= X_{\sigma(1)} a_1 \times \cdots \times X_{\sigma(i)} a_i'' \times \cdots \times X_{\sigma(n)} a_n
\end{aligned}\right.\]
\[\left\{\begin{aligned}
Y &= X_{\sigma(1)} a_1 \times \cdots \times (X_{\sigma(i)} a_i \times X_{\sigma(i)} a_i') \times \cdots \times X_{\sigma(n)} a_n\\
W &= X_1 a_{\sigmainv(1)} \times \cdots \times (X_{\sigma(i)} a_i \times X_{\sigma(i)} a_i') \times \cdots \times X_n a_{\sigmainv(n)}\\
\end{aligned}\right.\]
With these notation, the additivity axiom \cref{additivenattradditivity} for $\phi^\sigma$ asserts the commutativity of the boundary of the following diagram of functors.
\begin{equation}\label{sigmanatadditivity}
\begin{tikzpicture}[xscale=3,yscale=1,vcenter]
\def\v{-1} \def\u{-1.3}
\draw[0cell=.8]
(0,0) node (x1) {Y}
(x1)++(0,\u) node (w) {W}
(x1)++(-1,\v) node (x21) {X^\sigma\angapp}
(x1)++(1,\v) node (x22) {X^\sigma\anga \times X^\sigma\angap}
(x21)++(0,\u) node (x31) {X\sigma\angapp}
(x22)++(0,\u) node (x32) {X\sigma\anga \times X\sigma\angap}
(x31)++(0,\u) node (x41) {Za''_{\sigmainv}}
(x32)++(0,\u) node (x42) {Za_{\sigmainv} \times Za'_{\sigmainv}} 
(x41)++(0,\u) node (x51) {Za''}
(x42)++(0,\u) node (x52) {Za \times Za'}
(x51)++(1,\v) node (x6) {Z(a \oplus a')}
(x41)++(1,\v) node (z) {Z(a_{\sigmainv} \oplus a'_{\sigmainv})}
node[between=x21 and w at .5] {(\spadesuit)}
node[between=x22 and w at .5, shift={(-.3,0)}] {(\heartsuit)}
node[between=x51 and z at .5, shift={(-.4,.2)}] {(\clubsuit)}
node[between=x52 and z at .5, shift={(.4,.2)}] {(\diamondsuit)}
node[between=w and z at .6] {(\medstar)} 
;
\draw[1cell=.8]  
(x1) edge node[swap,pos=.4] {1 \times \cdots \times X_{\sigma(i)}^2 \times \cdots \times 1} (x21)
(x21) edge node[swap] {\sigma} (x31)
(x31) edge node[swap] {\phi_{\sigma\angapp}} (x41)
(x41) edge node[swap] {Z(\sigmainv)} (x51)
(x51) edge node[swap] {Z(\lap)} (x6)
(x1) edge node[pos=.5] {\varrho} (x22)
(x22) edge node {\sigma \times \sigma} (x32)
(x32) edge node {\phi_{\sigma\anga} \times \phi_{\sigma\angap}} (x42)
(x42) edge node {Z(\sigmainv) \times Z(\sigmainv)} (x52)
(x52) edge node {Z^2} (x6)
(x1) edge node {\barof{\sigma}} (w)
(w) edge node[pos=.8] {1 \times \cdots \times X_{\sigma(i)}^2 \times \cdots \times 1} (x31)
(w) edge node[swap,pos=.6] {\varrho} (x32)
(x41) edge node[pos=.3] {Z(\lap)} (z)
(x42) edge node[swap,pos=.3] {Z^2} (z)
(z) edge node[pos=.3] {Z(\sigmainv \oplus \sigmainv)} (x6)
;
\end{tikzpicture}
\end{equation}
The detail of the diagram \cref{sigmanatadditivity} is as follows.
\begin{itemize}
\item The isomorphism $\barof{\sigma} \cn Y \to W$ is the block permutation induced by $\sigma \in \Sigma_n$ that regards the product $X_{\sigma(i)} a_i \times X_{\sigma(i)} a_i'$ in $Y$ as a single block.
\item The region $(\spadesuit)$ commutes by the naturality of the braiding in the symmetric monoidal category $(\Cat,\times)$.
\item The region $(\heartsuit)$ commutes because the two composites become equal after projecting onto any factor of the common codomain, $X\sigma\anga \times X\sigma\angap$.
\item The region $(\clubsuit)$ commutes by the
\begin{itemize}
\item functoriality of $Z$ and
\item uniqueness of Laplaza coherence isomorphisms (\cref{thm:laplaza-coherence-1}).
\end{itemize} 
\item The region $(\diamondsuit)$ commutes by the naturality of the monoidal constraint $Z^2$ of the additive symmetric monoidal functor $(Z,Z^2,Z^0)$.
\item The hexagon $(\medstar)$ commutes by the additivity axiom \cref{additivenattradditivity} for $\phi$ 
\end{itemize}
Thus $\phi^\sigma$ is an additivity natural transformation.  This proves that additive natural transformations are closed under the symmetric group action \cref{dcatsigmaaction} in $\Dtecat$.

The proof that additive modifications are closed under the symmetric group action in $\Dtecat$ is obtained from the argument above by replacing $\phi$ with an additive modification $\Phi$.  This is valid for the following two reasons:
\begin{itemize}
\item $\Phi^\sigma$ \cref{Phisigmaanga} is defined in the same way as $\phi^\sigma$, with $\Phi_{\sigma\anga}$ replacing $\phi_{\sigma\anga}$.
\item The unity axiom \cref{additivemodunity} and the additivity axiom \cref{additivemodadditivity} for an additive modification are obtained from those for an additive natural transformation---namely, \cref{additivenattrunity,additivenattradditivity}---by replacing $\phi$ with $\Phi$.
\end{itemize} 
This finishes the proof.
\end{proof}

\subsection*{Composition}

Next we check that the composition in $\DCat$, which is a restriction of the composition in $\Dtecat$, is well defined.

\begin{lemma}\label{dcatgammaclosed}
Additive natural transformations and additive modifications are closed under the composition \cref{dcatgamma} in $\Dtecat$.
\end{lemma}

\begin{proof}
Suppose given additive symmetric monoidal functors (\cref{def:additivesmf})
\[Z \scs \ang{X_j}_{j=1}^n \scs \ang{\ang{W_{ji}}_{i=1}^{\ell_j}}_{j=1}^n \cn \Dplus \to \Cat\]
with $n > 0$, $\angX = (X_1, \ldots, X_n)$,
\[\angWj = \big(W_{j1},\ldots,W_{j\ell_j}\big), \andspace \angW = \big(\ang{W_1}, \ldots, \ang{W_n}\big).\]
Suppose given additive natural transformations (\cref{def:additivenattr})
\[\begin{tikzcd}
\bigotimes_{j=1}^n X_j \ar{r}{\phi} & Z
\end{tikzcd}\andspace
\begin{tikzcd}
\bigotimes_{i=1}^{\ell_j} W_{ji} \ar{r}{\phi_j} & X_j
\end{tikzcd}\]
for $1 \leq j \leq n$.  If some $\ell_j = 0$, then we make some appropriate adjustments in the argument below, as discussed in the last paragraph in \cref{expl:dcatgamma}.  In the rest of this proof, we assume that each $\ell_j > 0$.  

We must show that the composite natural transformation \cref{gaphiphij}
\[\begin{tikzcd}[column sep=2cm]
\bigotimes_{j=1}^n \bigotimes_{i=1}^{\ell_j} W_{ji} \ar{r}{\ga\smscmap{\phi;\angphij}} & Z
\end{tikzcd}\]
satisfies the unity axiom \cref{additivenattrunity} and the additivity axiom \cref{additivenattradditivity}.  The unity axiom for $\ga\smscmap{\phi;\angphij}$ follows from
\begin{itemize}
\item the unity axiom for $\phi$ and $\phi_j$ and
\item the definition \cref{dcatganattr} of each component functor of $\ga\smscmap{\phi;\angphij}$ as the composite
\[\begin{tikzcd}[column sep=huge]
\prod_{j=1}^n \prod_{i=1}^{\ell_j} W_{ji} a_{ji} \ar{r}{\prod_j (\phi_j)_{\angaj}} & \prod_{j=1}^n X_j a_j \ar{r}{\phi_{\anga}} & [-20pt] Za.
\end{tikzcd}\]
\end{itemize} 

To check the additivity axiom \cref{additivenattradditivity} for $\ga\smscmap{\phi;\angphij}$, suppose given objects
\begin{itemize}
\item $a_{ji} \in \cD$ for $j \in \{1,\ldots,n\}$ and $i \in \{1,\ldots,\ell_j\}$ and
\item $a_{hk}' \in \cD$ for one index $h \in \{1,\ldots,n\}$ and one index $k \in \{1,\ldots,\ell_h\}$,
\end{itemize}
along with the notation in \cref{angajangaaja} and below. 
\[\left\{\begin{aligned}
a_{hk}'' &= a_{hk} \oplus a_{hk}' & \abar &= a_1 \otimes \cdots \otimes (a_h \oplus a_h') \otimes \cdots \otimes a_n\\
a_h' &= a_{h1} \otimes \cdots \otimes a_{hk}' \otimes \cdots \otimes a_{h\ell_h} & a' &= a_1 \otimes \cdots \otimes a_h' \otimes \cdots \otimes a_n\\
a_h'' &= a_{h1} \otimes \cdots \otimes a_{hk}'' \otimes \cdots \otimes a_{h\ell_h} & a'' &= a_1 \otimes \cdots \otimes a_h'' \otimes \cdots \otimes a_n
\end{aligned}\right.\]
\[\left\{\begin{aligned}
A &= \textstyle \prod_{j=1}^{h-1} \prod_{i=1}^{\ell_j} W_{ji} a_{ji} & B_1 &= \textstyle \prod_{t=1}^{k-1} W_{ht} a_{ht} &&\\
C &= \textstyle \prod_{j=h+1}^{n} \prod_{i=1}^{\ell_j} W_{ji} a_{ji} & B_2 &= \textstyle \prod_{t=k+1}^{\ell_h} W_{ht} a_{ht} &&\\
B &= W_{hk} a_{hk} & B' &= W_{hk} a_{hk}' & B'' &= W_{hk} a_{hk}''
\end{aligned}\right.\]
\[\left\{\begin{aligned}
\phi^1 &= \textstyle \prod_{j=1}^{h-1} \phi_j & \phi^2 &= \textstyle \prod_{j=h+1}^{n} \phi_j & \phitil &= \textstyle \prod_{j=1}^{n} \phi_j = \phi^1 \times \phi_h \times \phi^2\\ 
T &= \textstyle \prod_{j=1}^{h-1} X_j a_j & V &= \textstyle \prod_{j=h+1}^{n} X_j a_j & \Ubar &= X_h(a_h \oplus a_h') \\
U &= X_h a_h & U' &= X_h a_h' & U'' &= X_h a_h''
\end{aligned}\right.\]
In $\phi^1$, $\phi^2$, and $\phitil$ above, the subscript in each component functor is omitted.  To further simplify the diagram below, we omit the symbol $\times$ and use juxtaposition to denote products.  With these notation and convention, the additivity axiom \cref{additivenattradditivity} for $\ga\smscmap{\phi;\angphij}$ asserts the commutativity of the boundary of the following diagram.
\begin{equation}\label{gammanatadditivity}
\begin{tikzpicture}[xscale=3.8,yscale=1.4,vcenter]
\def\s{.7} \def\t{-1.5} \def\u{-1.3} \def\v{-1} \def\w{-.6} 
\draw[0cell=\s]
(0,0) node (x1) {A (B_1 B B' B_2) C}
(x1)++(-1,\w) node (x21) {A (B_1 B'' B_2) C}
(x21)++(0,\t) node (x31) {T U'' V}
(x31)++(0,\t) node (x41) {Za''}
(x41)++(1,\w) node (x6) {Z(a \oplus a')}
(x1)++(1,\w) node (x22) {A (B_1 B B_2)(B_1 B' B_2) C}
(x22)++(0,\v) node (x32) {(AB_1 B B_2 C)(AB_1 B' B_2 C)}
(x32)++(0,\v) node (x42) {(TUV)(TU'V)}
(x42)++(0,\v) node (x52) {(Za)(Za')}
(x31)++(.67,-.25) node (z3) {T\Ubar V}
node[between=z3 and x22 at .5] (z1) {A(UU')C}
(z1)++(0,-.8) node (z2) {T(UU') V}
(z3)++(0,\v) node (z4) {Z\abar}
node[between=x21 and z1 at .5] {(\spadesuit)}
node[between=z2 and x32 at .4, shift={(.4,0)}] {(\heartsuit)}
node[between=z3 and z2 at .65, shift={(0,.5)}] {(\clubsuit)}
node[between=x31 and z4 at .5] {(\diamondsuit)}
node[between=x41 and x6 at .6, shift={(0,.6)}] {(\filledstar)}
node[between=z4 and x42 at .5] {(\filledlozenge)}
;
\draw[1cell=\s]  
(x1) edge node[swap,pos=.4] {1 \cdots W_{hk}^2  \cdots 1} (x21)
(x21) edge node[swap] {\phitil} (x31)
(x31) edge node[swap] {\phi} (x41)
(x41) edge node[swap] {Z(\lap)} (x6)
(x1) edge node[pos=.4] {1 \varrho 1} (x22)
(x22) edge node {\varrho} (x32)
(x32) edge node {\phitil \phitil} (x42)
(x42) edge node {\phi \phi} (x52)
(x52) edge node {Z^2} (x6)
(x22) edge node[swap] {1 \phi_h \phi_h 1} (z1)
(z1) edge node[swap,pos=.3] {\phi^1 (X_h^2) \phi^2} (z3)
(z1) edge node {\phi^1 1 \phi^2} (z2)
(z2) edge node {1 (X_h^2) 1} (z3)
(z2) edge node[swap] {\varrho} (x42)
(x31) edge node[pos=.2] {1 X_h(\lap) 1} (z3)
(x41) edge node[pos=.7] {Z(\lap)} (z4)
(z3) edge node {\phi} (z4)
(z4) edge node {Z(\lap)} (x6)
;
\end{tikzpicture}
\end{equation}
The detail of the diagram \cref{gammanatadditivity} is as follows.
\begin{itemize}
\item Each $\lap$ is a Laplaza coherence isomorphism (\cref{thm:laplaza-coherence-1}).
\item Each $\varrho$ is a composite as in \cref{varrhofactor} that involves a product of diagonal functors on some of the factors and identity functors, followed by a permutation.
\item The region $(\spadesuit)$ commutes by the additivity axiom \cref{additivenattradditivity} for $\phi_h$.
\item The regions $(\heartsuit)$ and $(\clubsuit)$ commute by construction.
\item The region $(\diamondsuit)$ commutes by the naturality \cref{dcatnaryobjnaturality} of $\phi$.
\item The region $(\filledstar)$ commutes by the
\begin{itemize}
\item functoriality of $Z$ and
\item uniqueness of Laplaza coherence isomorphisms (\cref{thm:laplaza-coherence-1}).
\end{itemize}
\item The region $(\filledlozenge)$ commutes by the additivity axiom \cref{additivenattradditivity} for $\phi$.
\end{itemize}
This proves that $\ga\smscmap{\phi;\angphij}$ satisfies the additivity axiom \cref{additivenattradditivity}.  We have shown that additive natural transformations are closed under the composition in $\Dtecat$.  

The proof that additive modifications are closed under the composition in $\Dtecat$ is obtained from the argument above by replacing $\phi$ and $\phi_j$ with, respectively, additive modifications $\Phi \cn \phi \to \varphi$ and $\Phi_j \cn \phi_j \to \varphi_j$.  This is valid for the following two reasons:
\begin{itemize}
\item Each component natural transformation of $\ga\smscmap{\Phi;\ang{\Phi_j}}$ is the horizontal composite \cref{dcatgamod}.
\item The unity axiom \cref{additivemodunity} and the additivity axiom \cref{additivemodadditivity} for an additive modification are obtained from those for an additive natural transformation---namely, \cref{additivenattrunity,additivenattradditivity}---by replacing $\phi$ with $\Phi$.
\end{itemize} 
This finishes the proof.
\end{proof}

\begin{explanation}[Unary Case]\label{expl:dcatgammaclosed}
As we mentioned in \cref{expl:additivenattr} \eqref{expl:additivenattr-iv}, additive natural transformations in arity 1 are precisely monoidal natural transformations.  With this in mind, the proof of \cref{dcatgammaclosed} generalizes the proof that monoidal natural transformations are closed under vertical composition.
\end{explanation}

We are now ready for the main observation of this chapter.

\begin{theorem}\label{thm:dcatcatmulticat}\index{diagram category!multicategory structure}\index{multicategory!of diagrams}\index{Cat-multicategory@$\Cat$-multicategory!of diagrams}
For each small tight bipermutative category $\cD$, $\DCat$ in \cref{def:dcatcatmulticat} is a $\Cat$-multicategory.
\end{theorem}

\begin{proof}
We note the following facts.
\begin{itemize}
\item By \cref{dcatnarycategory} each positive arity multimorphism category in $\DCat$ is a subcategory of the corresponding multimorphism category in $\Dtecat$. 
\item By \cref{dcatzeroary,dcatbipermzeroary} the arity 0 multimorphism categories in $\DCat$ and $\Dtecat$ are canonically isomorphic.
\item The colored units in $\Dtecat$ and $\DCat$ are both given by identity natural transformations.
\item By \cref{dcatsigmaclosed,dcatgammaclosed}, respectively, the symmetric group action and composition in $\Dtecat$ restrict to $\DCat$.
\end{itemize}
Thus the axioms of a $\Cat$-multicategory (\cref{def:enr-multicategory}) hold in $\DCat$ because they hold in $\Dtecat$ (\cref{thm:permindexedcat}).
\end{proof}

Since the $\Cat$-multicategory structure on $\DCat$ (\cref{thm:dcatcatmulticat})---the colored units, the symmetric group action, and the composition---is inherited from that on $\Dtecat$ (\cref{thm:permindexedcat}), we obtain the following result.  Recall from \cref{def:enr-multicategory-functor} the notion of a $\V$-multifunctor between $\V$-multicategories.

\begin{corollary}\label{dcatforgetdtecat}
For each small tight bipermutative category $\cD$, there is a forgetful $\Cat$-multifunctor
\[\begin{tikzcd}[column sep=large]
\DCat \ar{r}{U} & \Dtecat
\end{tikzcd}\]
that is defined as follows.
\begin{itemize}
\item $U$ sends each additive symmetric monoidal functor $(Z,Z^2,Z^0)$ to the underlying functor $Z \cn \cD \to \Cat$.
\item $U$ is the canonical isomorphism given by \cref{dcatzeroary,dcatbipermzeroary} for arity 0 multimorphism categories.
\item $U$ is the subcategory inclusion in \cref{dcatnarycategory} for positive arity multimorphism categories.
\end{itemize} 
\end{corollary}

\begin{explanation}[Non-injectivity on Objects]\label{expl:dcatforgetdtecat}
For the $\Cat$-multifunctor $U$ in \cref{dcatforgetdtecat}, each multimorphism functor is either an isomorphism or a subcategory inclusion.  However, in general, $U$ is \emph{not} injective on objects because a $\cD$-indexed category $Z \cn \cD \to \Cat$ may admit several different additive structures $(Z^2,Z^0)$.  Therefore, $\DCat$ is, in general, \emph{not} a sub-$\Cat$-multicategory of $\Dtecat$.
\end{explanation}

\begin{example}\label{ex:dcatexamples}
Suppose $\cD$ is any one of the small tight bipermutative categories
\[\Finsk,\quad \Fset, \quad \Fskel, \andspace \cA\]
in, respectively, \cref{ex:finsk,ex:Fset,ex:Fskel,ex:mandellcategory}.  Then \cref{thm:dcatcatmulticat,dcatforgetdtecat} apply to $\DCat$.  In particular, there is a $\Cat$-multicategory $\ACat$ whose objects are symmetric monoidal functors
\[\Aplus = \big(\cA,\oplus,\ang{},\betaplus\big) \to \big(\Cat,\times,\boldone,\xi\big).\]
As we will discuss in \cref{ch:multifunctorA,ch:invK}, $\ACat$ is an intermediate $\Cat$-multicategory in inverse $K$-theory.
\end{example}

\chapter{The Grothendieck Construction is a Pseudo Symmetric \texorpdfstring{$\Cat$}{Cat}-Multifunctor}
\label{ch:multigro}
For a small tight bipermutative category $(\cD,\oplus,\otimes)$, in \cref{thm:dcatcatmulticat} we observed that $\DCat$ is a $\Cat$-multicategory, with additive symmetric monoidal functors $\Dplus \to \Cat$ as objects (\cref{def:additivesmf}).  The main result of this chapter is \cref{thm:grocatmultifunctor}.  It shows that the Grothendieck construction of indexed categories extends to a \emph{pseudo symmetric} $\Cat$-multifunctor 
\[\begin{tikzcd}[column sep=large,every label/.append style={scale=.85}]
\DCat \ar{r}{\grod} & \permcatsg.
\end{tikzcd}\]
Furthermore, it is a $\Cat$-multifunctor if the multiplicative braiding in $\cD$ is the identity.  Conversely, if $\grod$ is a $\Cat$-multifunctor, then 
\[x \otimes y=y\otimes x\]
for all objects $x,y \in \cD$.  The codomain $\permcatsg$ (\cref{sec:permcatmulticat}) is a $\Cat$-multicategory with small permutative categories as objects.

\subsection*{Summary}

Extending part of the table in the introduction of \cref{ch:diagram}, the following table summaries the $\Cat$-multicategories $\DCat$ and $\permcatsg$ and the pseudo symmetric $\Cat$-multifunctor $\grod$.  We abbreviate \emph{symmetric monoidal}, \emph{multimorphism category}, and \emph{natural transformations} to, respectively, \emph{sm}, \emph{mc}, and \emph{nt}.
\begin{center}
\resizebox{\columnwidth}{!}{%
{\renewcommand{\arraystretch}{1.4}%
{\setlength{\tabcolsep}{1ex}
\begin{tabular}{|c|cc|cc|c|}\hline
& $\DCat$ &(\ref{thm:dcatcatmulticat}) & $\permcatsg$ &(\ref{thm:permcatenrmulticat}) & $\grod\,\,$ (\ref{thm:grocatmultifunctor}) \\ \hline
objects 
& additive sm functors &(\ref{def:additivesmf}) & small permutative categories &(\ref{def:symmoncat}) & \ref{grodxpermutative}\\ \hline
mc objects 
& additive nt &(\ref{def:additivenattr}) & strong multilinear functors &(\ref{def:nlinearfunctor}) & \ref{def:groconarityzero}, \ref{def:groconarityn}\\ \hline
mc morphisms 
& additive modifications &(\ref{def:additivemodification}) & multilinear nt &(\ref{def:nlineartransformation}) & \ref{def:groconmodification}\\ \hline
units & identity nt & (\ref{zcoloredunit}) & identity 1-linear functors & (\ref{ex:idonelinearfunctor}) & \ref{grodpreservesunits}\\ \hline
composition & \ref{dcatgamma}, \ref{expl:dcatgamma} && \ref{permcatgamma} && \ref{grodpreservescomposition}\\ \hline
symmetry & \ref{dcatsigmaaction} && \ref{permcatsymgroupaction} && \ref{gropseudosymmetry}\\ \hline
induced by sm structure
& no & (\ref{expl:dcatbipermzero}) & no & (\ref{rk:permcatnotmonoidal}) & no\\ \hline
\end{tabular}}}}
\end{center}
\smallskip

\subsection*{Pseudo Symmetry and Inverse $K$-Theory}

Recall from \cref{def:pseudosmultifunctor} that a pseudo symmetric $\Cat$-multifunctor \emph{strictly} preserves the colored units and multicategorical composition.  It preserves the symmetric group action up to the pseudo symmetry isomorphisms, which satisfy the four coherence axioms \cref{pseudosmf-unity,pseudosmf-product,pseudosmf-topeq,pseudosmf-boteq}.  There is a canonical bijective correspondence between (i) $\Cat$-multifunctors and (ii) pseudo symmetric $\Cat$-multifunctors with identity pseudo symmetry isomorphisms (\cref{ex:idpsm}).  

The components of the pseudo symmetry isomorphisms for the Grothendieck construction are defined in \cref{grodsigmacomponent}.  They make full use of the multiplicative braiding $\betate$ in $\cD$, so they are almost never identities in practice.  The Grothendieck construction $\grod$ is a $\Cat$-multifunctor if and only if the multiplicative braiding $\betate$ in $\cD$ happens to be the identity natural transformation.  For example, the pseudo symmetry isomorphisms for the Grothendieck construction are not identities for each of the small tight bipermutative categories $\Finsk$, $\Fset$, $\Fskel$, and $\cA$ in \cref{sec:examples}.  In each of these cases, the Grothendieck construction $\grod$ is genuinely a pseudo symmetric $\Cat$-multifunctor and \emph{not} a $\Cat$-multifunctor (\cref{ex:grodpseudosymmetric}).  

An important consequence of \cref{thm:grocatmultifunctor} is that the inverse $K$-theory functor, from $\Ga$-categories to small permutative categories, is a pseudo symmetric $\Cat$-multifunctor but \emph{not} a $\Cat$-multifunctor.  As we will discuss in \cref{ch:multifunctorA,ch:invK}, the inverse $K$-theory functor $\groa A(-)$ is a composite of two functors, the second of which is the Grothendieck construction $\groa$ with the indexing category $\cA$ (\cref{ex:mandellcategory}).  Since the multiplicative braiding in $\cA$ is not the identity, the pseudo symmetry isomorphisms of $\groa$ are not identities.

\subsection*{Organization}

In \cref{sec:permcatmulticat} we discuss the $\Cat$-multicategory $\permcatsg$.  It is the codomain of the Grothendieck construction $\grod$. 
\begin{itemize}
\item Its objects are small permutative categories.
\item Its $n$-ary 1-cells---that is, objects in the multimorphism categories---are \emph{strong} $n$-linear functors.  The strong condition, which is what the superscript $\mathsf{sg}$ refers to, means that each linearity constraint is a natural isomorphism. 
\item Its $n$-ary 2-cells are $n$-linear natural transformations.
\end{itemize} 

In \cref{sec:groadditivesmf} we first review the Grothendieck construction of an indexed category.  This construction associates to each functor $X \cn \cD \to \Cat$ a category $\grod X$, which is small if $\cD$ is small.  In the literature the Grothendieck construction is sometimes called the \emph{wreath product} and denoted $\cD \smallint X$.  Next we discuss the case with $\cD$ a permutative category and $X$ a symmetric monoidal functor.  In this case the Grothendieck construction has the canonical structure of a permutative category.  We emphasize that the symmetric monoidal functor $X$ is \emph{not} assumed to be strong.  In other words, its unit and monoidal constraints are not required to be invertible.  The construction $X \mapsto\! \grod X$ is the object assignment of the Grothendieck construction.

In \cref{sec:groadditivenattr} we extend the Grothendieck construction $\grod$ to $n$-ary 1-cells in $\DCat$, which are additive natural transformations (\cref{def:additivenattr}).  For such an additive natural transformation $\phi$, the Grothendieck construction $\grod\phi$ is a \emph{strong} $n$-linear functor.  The components of its linearity constraints $(\grod\phi)^2_i$ are defined in \cref{grodphilinearity}. 
\begin{itemize}
\item As a morphism in a Grothendieck construction $\grod Z$, the first entry in each component of $(\grod\phi)^2_i$ is a morphism in the indexing small tight bipermutative category $\cD$.  In fact, it is defined as the inverse of a Laplaza coherence isomorphism (\cref{ex:laplazacoherence}).  Proofs regarding these linearity constraints often involve the uniqueness in Laplaza's Coherence \cref{thm:laplaza-coherence-1}.  For example, in \cref{grodphifunctor} we prove each of the multilinearity constraint axioms for $\grod\phi$ by applying \cref{thm:laplaza-coherence-1}.
\item The second entry in each component of $(\grod\phi)^2_i$ is an identity morphism.
\end{itemize} 

In \cref{sec:groadditivemod} we extend the Grothendieck construction $\grod$ to $n$-ary 2-cells in $\DCat$, which are additive modifications (\cref{def:additivemodification}).  For such an additive modification $\Phi$, the Grothendieck construction $\grod\Phi$ is an $n$-linear natural transformation.  Each of its components is a morphism in a Grothendieck construction $\grod Z$; see \cref{grodPhiajxj}. 
\begin{itemize}
\item Its first entry is an identity morphism in $\cD$.
\item Its second entry is a component morphism of a component natural transformation of $\Phi$. 
\end{itemize}  
\cref{grodmultimorphismfunctor} shows that the Grothendieck constructions of
\begin{itemize}
\item additive symmetric monoidal functors $\Dplus \to \Cat$,
\item additive natural transformations, and
\item additive modifications
\end{itemize}
define functors from the multimorphism categories of $\DCat$ to those of $\permcatsg$.

In \cref{sec:gropreservesunitcomp} we show that the Grothendieck construction $\grod$ preserves the colored units and composition in the $\Cat$-multicategories $\DCat$ and $\permcatsg$.  For additive natural transformations, the proof that $\grod$ preserves composition involves showing that the relevant linearity constraints are equal.  This is once again an application of Laplaza's Coherence \cref{thm:laplaza-coherence-1}; see \cref{lapappaplusap,gagrophiphijconstraint}.

In \cref{sec:pseudosymgrothendieck} we first construct the pseudo symmetry isomorphisms $(\grod)_\sigma$ for the Grothendieck construction.  For an additive natural transformation $\phi$, each component of $(\grod)_{\sigma,\phi}$ is an isomorphism in a Grothendieck construction $\grod Z$; see \cref{grodsigmacomponent}. 
\begin{itemize}
\item Its first entry is a coherence isomorphism in the multiplicative structure $\Dte$ that permutes the tensor factors.
\item Its second entry is an identity morphism.
\end{itemize}   
In the rest of this section, we show that
\begin{itemize}
\item each component $(\grod)_{\sigma,\phi}$ is a multilinear natural isomorphism (\cref{grodsigmaphilinear,grodsigmaphinlinear}) and
\item it is a natural isomorphism in the $\phi$ variable (\cref{grodsigmanaturaliso}).
\end{itemize} 

In \cref{sec:enrmultigrothendieck} we finish the proof that the Grothendieck construction $\grod$ is a pseudo symmetric $\Cat$-multifunctor by checking the four axioms in \cref{def:pseudosmultifunctor}.  These proofs mostly involve unpacking the definitions of various constructions and do not involve \cref{thm:laplaza-coherence-1}.  We provide all the explicit steps because we want to leave no doubt that our many constructions and definitions are correct.  Moreover, $\grod$ is a $\Cat$-multifunctor if the multiplicative braiding in $\cD$ is the identity; a weaker converse is also true.

We remind the reader of our left normalized bracketing \cref{expl:leftbracketing} for iterated monoidal product.

\section{Enriched Multicategory of Permutative Categories}
\label{sec:permcatmulticat}

In this section we review the $\Cat$-enriched multicategory $\permcat$ with
\begin{itemize}
\item small permutative categories (\cref{def:symmoncat}) as objects,
\item $n$-linear functors (\cref{def:nlinearfunctor}) as $n$-ary 1-cells, and
\item $n$-linear natural transformations (\cref{def:nlineartransformation}) as $n$-ary 2-cells.
\end{itemize}   
The existence of the $\Cat$-multicategory $\permcat$ is stated in \cite[page 174]{elmendorf-mandell}.  A detailed proof is in \cite[Sections 6.5 and 6.6]{cerberusIII}, from which the presentation below is adapted.  Restricting to \emph{strong} multilinear functors yields the $\Cat$-enriched multicategory $\permcatsg$ that is the codomain of the Grothendieck construction.

\subsection*{Multilinear Functors}

The definitions of multilinear functors and multilinear natural transformations involve substitution of variables.  We use the following notation to simplify relevant presentation.

\begin{notation}\label{notation:compk}
Suppose 
\[\ang{x} = (x_1, \ldots, x_n)\]
is an $n$-tuple of symbols, and $x_i'$ is a symbol for some index $i \in \{1,\ldots,n\}$.  
\begin{itemize}
\item We denote by\label{not:compk}
\begin{equation}\label{compnotation}
\ang{x \compi x_i'} = \ang{x} \compi x_i' 
= \big(\overbracket[0.5pt]{x_1, \ldots, x_{i-1}}^{\text{empty if $i=1$}}, x_i', \overbracket[0.5pt]{x_{i+1}, \ldots, x_n}^{\text{empty if $i=n$}}\big)
\end{equation}
the $n$-tuple obtained from $\ang{x}$ by replacing its $i$-th entry by $x_i'$.
\item For $i \neq \ell \in \{1,\ldots,n\}$ and a symbol $x'_\ell$, we denote by
\begin{equation}\label{compcompnotation}
\ang{x \compi x_i' \comp_\ell x'_\ell} = \ang{x} \compi x_i' \comp_\ell x'_\ell
\end{equation}
the $n$-tuple obtained from $\ang{x \compi x_i'}$ by replacing its $\ell$-th entry by $x'_\ell$.\defmark
\end{itemize}
\end{notation}

The following definition of an $n$-linear functor is from \cite[3.2]{elmendorf-mandell}, as presented in \cite[6.5.4]{cerberusIII}.  For the rest of this section, suppose $\C_1, \ldots, \C_n$, and $\D$ are permutative categories (\cref{def:symmoncat}).  In each permutative category, the monoidal product, the monoidal unit, and the braiding are denoted, respectively, $\oplus$, $e$, and $\beta$.

\begin{definition}\label{def:nlinearfunctor}
An \emph{$n$-linear functor}\index{multilinear functor}\index{functor!multilinear}\index{n-linear functor@$n$-linear functor}
\begin{equation}\label{nlinearfunctorangcd}
\begin{tikzcd}[column sep=2.3cm]
\prod_{j=1}^n \C_j \ar{r}{\left(F,\, \ang{F^2_j}_{j=1}^n\right)} & \D
\end{tikzcd}
\end{equation}
consists of the following data.
\begin{itemize}
\item $F \cn \prod_{j=1}^n \C_j \to \D$ is a functor.
\item For each $j\in \{1,\ldots,n\}$, $F^2_j$ is a natural transformation, called the \index{linearity constraint}\emph{$j$-th linearity constraint}, with component morphisms
\begin{equation}\label{laxlinearityconstraints}
\begin{tikzcd}[column sep=large]
F\ang{X\compj X_j} \oplus F\ang{X\compj X_j'} \ar{r}{F^2_j}
& F\ang{X\compj (X_j \oplus X_j')} \in \D
\end{tikzcd}
\end{equation}
for objects $\ang{X} = \ang{X_j} \in \prod_{j=1}^n \C_j$ and $X_j' \in \C_j$.
\end{itemize}
The data above are subject to the following five axioms.
\begin{description}
\item[Unity]\index{unity!multilinear functor} 
$F$ satisfies the following unity axioms for $j \in \{1,\ldots,n\}$, objects $\angX$, and morphisms $\angf$ in $\prod_{j=1}^n \C_j$.
\begin{equation}\label{nlinearunity}
\left\{\begin{split}
F \ang{X \compj e} &= e\\
F \ang{f \compj 1_e} &= 1_e\end{split}\right.
\end{equation}
\item[Constraint Unity]\index{constraint unity!multilinear functor} 
\begin{equation}\label{constraintunity}
F^2_j = 1 \qquad \text{if any $X_i = e$ or if $X_j'=e$}.
\end{equation}
\item[Constraint Associativity]\index{constraint associativity!multilinear functor}\index{associativity!multilinear functor constraint} The following diagram commutes for each $i\in \{1,\ldots,n\}$ and objects $\ang{X} \in \prod_{j=1}^n \C_j$, with $X_i', X_i'' \in \C_i$.
    \begin{equation}\label{eq:ml-f2-assoc}
    \begin{tikzpicture}[x=55mm,y=12mm,vcenter]
    \tikzset{0cell/.append style={nodes={scale=.85}}}
    \tikzset{1cell/.append style={nodes={scale=.85}}}
      \draw[0cell] 
      (0,0) node (a) {
        F\ang{X\compi X_i}
        \oplus F\ang{X\compi X_i'}
        \oplus F\ang{X\compi X_i''}
      }
      (1,0) node (b) {
        F\ang{X\compi X_i}
        \oplus F\ang{X\compi (X_i' \oplus X_i'')}
      }
      (0,-1) node (c) {
        F\ang{X\compi (X_i \oplus X_i')}
        \oplus F\ang{X\compi X_i''}
      }
      (1,-1) node (d) {
        F\ang{X\compi (X_i \oplus X_i' \oplus X_i'')}
      }
      ;
      \draw[1cell] 
      (a) edge node {1 \oplus F^2_i} (b)
      (a) edge['] node {F^2_i \oplus 1} (c)
      (b) edge node {F^2_i} (d)
      (c) edge node {F^2_i} (d)
      ;
    \end{tikzpicture}
    \end{equation}
\item[Constraint Symmetry]\index{constraint symmetry!multilinear-functor}\index{symmetry!multilinear functor constraint} The following diagram commutes for each $i\in \{1,\ldots,n\}$ and objects $\ang{X} \in \prod_{j=1}^n \C_j$, with $X_i' \in \C_i$.
    \begin{equation}\label{eq:ml-f2-symm}
    \begin{tikzpicture}[x=40mm,y=12mm,vcenter]
\tikzset{0cell/.append style={nodes={scale=.85}}}
\tikzset{1cell/.append style={nodes={scale=.85}}}
      \draw[0cell] 
      (0,0) node (a) {
        F\ang{X\compi X_i}
        \oplus F\ang{X\compi X_i'}
      }
      (1,0) node (b) {
        F\ang{X\compi (X_i \oplus X_i')}
      }
      (0,-1) node (c) {
        F\ang{X\compi X_i'}
        \oplus F\ang{X\compi X_i}
      }
      (1,-1) node (d) {
        F\ang{X\compi (X_i' \oplus X_i)}
      }
      ;
      \draw[1cell] 
      (a) edge node {F^2_i} (b)
      (a) edge['] node {\beta} (c)
      (b) edge node {F\ang{1 \compi \beta}} (d)
      (c) edge node {F^2_i} (d)
      ;
    \end{tikzpicture}
    \end{equation}
  \item[Constraint 2-By-2]\index{constraint 2-by-2!multilinear functor}\index{2-by-2!multilinear functor constraint} The following diagram commutes for each
\[i \neq k\in \{1,\ldots,n\}, \quad \ang{X} \in \txprod_{j=1}^n \C_j, \quad
X_i' \in \C_i, \andspace X_k' \in \C_k.\]
\begin{equation}\label{eq:f2-2by2}
\begin{tikzpicture}[x=37mm,y=17mm,vcenter]
\tikzset{0cell-nomm/.style={commutative diagrams/every diagram,every cell,nodes={scale=.8}}}
\tikzset{1cell/.append style={nodes={scale=.8}}}
      \draw[0cell-nomm] 
      (0,0) node[align=left] (a) {$
        \phantom{\oplus}F\ang{X\compi X_i \compk X_k}
        \oplus F\ang{X\compi X_i' \compk X_k}$\\       
        $\oplus F\ang{X\compi X_i \compk X_k'}
        \oplus F\ang{X\compi X_i' \compk X_k'}
      $}
      (0,-1) node[align=left] (a') {$
        \phantom{\oplus}F\ang{X\compi X_i \compk X_k}
        \oplus F\ang{X\compi X_i \compk X_k'}$\\
        $\oplus F\ang{X\compi X_i' \compk X_k}
        \oplus F\ang{X\compi X_i' \compk X_k'}
      $}
      (.7,.75) node (b) {$
        F\ang{X\compi (X_i \oplus X_i') \compk X_k}
        \oplus F\ang{X\compi (X_i \oplus X_i') \compk X_k'} 
      $}
      (.7,-1.75) node (b') {$
        F\ang{X\compi X_i \compk (X_k \oplus X_k')}
        \oplus F\ang{X\compi X_i' \compk (X_k \oplus X_k')}
      $}
      (1.1,-.5) node (c) {$
        F\ang{X\compi (X_i \oplus X_i') \compk (X_k \oplus X_k')}
      $}
      ;
      \draw[1cell] 
      (a) edge[transform canvas={xshift=-4ex}] node[pos=.1] {F^2_i \oplus F^2_i} (b)
      (b) edge[transform canvas={shift={(.1,0)}}] node {F^2_k} (c)
      (a) edge[transform canvas={xshift=.5ex}] node[swap] {1 \oplus \beta \oplus 1} (a')
      (a') edge[transform canvas={xshift=-4ex}] node[pos=.1,swap] {F^2_k \oplus F^2_k} (b')
      (b') edge[',transform canvas={shift={(.1,0)}}] node {F^2_i} (c)
      ;
\end{tikzpicture}
\end{equation}
\end{description}
This finishes the definition of an $n$-linear functor.  

Moreover, we define the following.
\begin{itemize}
\item If $n=0$, then a \emph{0-linear functor}\index{0-linear functor} is a choice of an object in $\D$.
\item An $n$-linear functor $(F, \ang{F^2_j})$ is \index{strong multilinear functor}\emph{strong} if each linearity constraint $F^2_j$ is a natural isomorphism.
\item A \emph{multilinear functor} is an $n$-linear functor for some $n \geq 0$.\defmark
\end{itemize}  
\end{definition}

\begin{example}[Identities]\label{ex:idonelinearfunctor}
Each permutative category $\C$ has an \emph{identity 1-linear functor}\index{identity!1-linear functor}\index{1-linear functor!identity} 
\[\begin{tikzcd}[column sep=large]
\C \ar{r}{1_{\C}} & \C
\end{tikzcd}\] 
consisting of
\begin{itemize}
\item the identity functor on $\C$ and
\item the identity linearity constraint.
\end{itemize}
The latter is componentwise the identity morphism of $X \oplus X'$ for objects $X,X' \in \C$.  Thus $1_{\C}$ is a strong 1-linear functor.
\end{example}

\begin{example}[1-Linear Functors]\label{ex:onelinearfunctor}
A 1-linear functor\index{1-linear functor} between permutative categories
\[\begin{tikzcd}[column sep=huge]
\C \ar{r}{(F,F^2)} & \D
\end{tikzcd}\]
is precisely a strictly unital symmetric monoidal functor (\cref{def:monoidalfunctor}), since the constraint 2-by-2 axiom \cref{eq:f2-2by2} is vacuous when $n=1$.
\end{example}

\begin{example}[Bipermutative Categories]\label{ex:bipermnlinearfunctor}
Suppose given a \emph{tight} bipermutative category (\cref{def:embipermutativecat}) 
\[\big(\C, (\oplus,\zero,\betaplus), (\otimes,\tu,\betate),\fal,\far\big)\]
with additive structure $\Cplus = (\C,\oplus,\zero,\betaplus)$.  For each $n \geq 0$ there is a strong $n$-linear functor
\[\begin{tikzcd}[column sep=2.3cm]
\txprod_{j=1}^n \Cplus \ar{r}{\brb{\otimes,\ang{\otimes^2_j}_{j=1}^n}} & \Cplus
\end{tikzcd}\]
defined as follows.
\begin{itemize}
\item The functor 
\[\begin{tikzcd}[column sep=large]
(\Cplus)^n \ar{r}{\otimes} & \Cplus
\end{tikzcd}\]
is any one of the equal iterates of $\otimes \cn \C^2 \to \C$.  It is
\begin{itemize}
\item the multiplicative unit $\tu \in \C$ if $n=0$ and
\item the identity functor $1_{\C}$ if $n=1$.
\end{itemize} 
\item For $n \geq 1$ and $j \in \{1,\ldots,n\}$, the $j$-th linearity constraint
\begin{equation}\label{laplinearityconst}
\begin{tikzcd}
(X_1 \otimes \cdots \otimes X_j \otimes \cdots \otimes X_n) \oplus (X_1 \otimes \cdots \otimes X_j' \otimes \cdots \otimes X_n) \ar[shorten <=-1ex]{d}{\lap_j^\inv}[swap]{\otimes^2_j \,=}\\
X_1 \otimes \cdots \otimes (X_j \oplus X_j') \otimes \cdots \otimes X_n
\end{tikzcd}
\end{equation}
is the inverse of the unique Laplaza coherence isomorphism $\lap_j$ (\cref{ex:laplazacoherence}) that distributes over the sum $X_j \oplus X_j'$ in the $j$-th tensor factor.  The naturality of $\otimes^2_j = \lap_j^\inv$ follows from the naturality of Laplaza coherence isomorphisms.
\end{itemize}
The axioms of an $n$-linear functor hold as follows.
\begin{itemize}
\item The unity axiom \cref{nlinearunity} follows from the multiplicative zero axiom \cref{ringcataxiommultzero} in $\C$.
\item The four constraint axioms, \cref{constraintunity,eq:ml-f2-assoc,eq:ml-f2-symm,eq:f2-2by2}, follow from the uniqueness in Laplaza's Coherence \cref{thm:laplaza-coherence-1}.
\end{itemize}

For instance, if $n=2$ then
\begin{itemize}
\item $\otimes^2_1 = \lap_1^\inv$ is the left factorization morphism $\fal$, and
\item $\otimes^2_2 = \lap_2^\inv$ is the right factorization morphism $\far$ \cref{ringcatfactorization}. 
\end{itemize} 
The four constraint axioms, \cref{constraintunity,eq:ml-f2-assoc,eq:ml-f2-symm,eq:f2-2by2}, follow from, respectively, the ring category axioms \cref{ringcataxiomzerofact}, \cref{ringcataxiomifa}, \cref{ringcataxiomfactsyml}, and \cref{ringcataxiomtwotwo} in $\C$.
\end{example}

\subsection*{Multilinear Natural Transformations}

The following multilinear analogs of natural transformations are from \cite[page 174]{elmendorf-mandell}, as presented in \cite[6.5.11]{cerberusIII}.

\begin{definition}\label{def:nlineartransformation}
Suppose $F$ and $G$ are $n$-linear functors as displayed below.
\begin{equation}\label{nlineartransformation}
\begin{tikzpicture}[xscale=2,yscale=2,baseline={(a.base)}]
\draw[0cell=.9]
(0,0) node (a) {\textstyle\prod_{j=1}^n \C_j}
(a)++(.15,.03) node (b) {\phantom{D}}
(b)++(1,0) node (c) {\D}
;
\draw[1cell=.8]  
(b) edge[bend left] node[pos=.5] {\big( F, \ang{F^2_j}\big)} (c)
(b) edge[bend right] node[swap,pos=.5] {\big( G, \ang{G^2_j}\big)} (c)
;
\draw[2cell] 
node[between=b and c at .45, rotate=-90, 2label={above,\theta}] {\Rightarrow}
;
\end{tikzpicture}
\end{equation}
An \index{natural transformation!multilinear}\index{multilinear natural transformation}\index{n-linear natural transformation@$n$-linear natural transformation}\emph{$n$-linear natural transformation} $\theta \cn F \to G$ is a natural transformation of underlying functors that satisfies the following two axioms.
\begin{description}
\item[Unity]\index{unity axiom!multilinear natural transformation}
\begin{equation}\label{ntransformationunity}
\theta_{\ang{X}} = 1_e \qquad \text{if any $X_i = e \in \C_i$}.
\end{equation}
\item[Constraint Compatibility]\index{constraint compatibility condition} 
The diagram
\begin{equation}\label{eq:monoidal-in-each-variable}
\begin{tikzpicture}[x=45mm,y=12mm,vcenter]
\tikzset{0cell/.append style={nodes={scale=.9}}}
\tikzset{1cell/.append style={nodes={scale=.9}}}
      \draw[0cell] 
      (0,0) node (a) {
        F\ang{X\compi X_i}
        \oplus F\ang{X\compi X_i'}
      }
      (1,0) node (b) {
        F\ang{X\compi (X_i \oplus X_i')}
      }
      (0,-1) node (c) {
        G\ang{X\compi X_i}
        \oplus G\ang{X\compi X_i'}
      }
      (1,-1) node (d) {
        G\ang{X\compi (X_i \oplus X_i')}
      }
      ;
      \draw[1cell] 
      (a) edge node {F^2_i} (b)
      (a) edge['] node {\theta \oplus \theta} (c)
      (b) edge node {\theta} (d)
      (c) edge node {G^2_i} (d)
      ;
  \end{tikzpicture}
\end{equation}
commutes for each $i \in \{1,\ldots,n\}$, objects $\ang{X} \in \prod_{j=1}^n \C_j$, and $X_i' \in \C_i$.
\end{description}
This finishes the definition of an $n$-linear natural transformation.  

Moreover, we define the following.
\begin{itemize}
\item An \emph{$n$-linear natural isomorphism}\index{multilinear natural transformation!isomorphism}\index{n-linear natural transformation@$n$-linear natural transformation!isomorphism} is an invertible $n$-linear natural transformation.
\item A \emph{0-linear natural transformation}\index{0-linear natural transformation} is a morphism in $\D$.
\item A \emph{multilinear natural transformation} is an $n$-linear natural transformation for some $n \geq 0$.\defmark
\end{itemize}
\end{definition}

\begin{definition}\label{def:permiicatsus}\index{multimorphism category!permutative category}\index{permutative category!multimorphism category}
For small permutative categories $\angC = \ang{\C_j}_{j=1}^n$ and $\D$, define the category
\begin{equation}\label{permcatangcd}
\permcat\mmap{\D; \ang{\C}} = 
\permcat\mmap{\D; \C_1, \ldots, \C_n}
\end{equation}
with
\begin{itemize}
\item $n$-linear functors $\prod_{j=1}^n \C_j \to \D$ (\Cref{def:nlinearfunctor}) as objects,
\item $n$-linear natural transformations (\Cref{def:nlineartransformation}) as morphisms, 
\item identity natural transformations as identity morphisms, and
\item vertical composition of natural transformations as composition.
\end{itemize}
This finishes the definition of the category $\permcatsmangcd$.

Moreover, we define the following.
\begin{itemize}
\item The \emph{strong} variant\index{multimorphism category!permutative category!strong variant}\index{permutative category!multimorphism category!strong variant}
\begin{equation}\label{permcatsgangcd}
\permcatsg\mmap{\D; \ang{\C}} = 
\permcatsg\mmap{\D; \C_1, \ldots, \C_n}
\end{equation}
is the full subcategory of $\permcatsmangcd$ with strong $n$-linear functors as objects.
\item If $n=0$ then we define
\begin{equation}\label{permcatemptyd}
\permcat\scmap{\ang{};\D} = \permcatsg\scmap{\ang{};\D} = \D
\end{equation}
as categories.\defmark
\end{itemize}
\end{definition}

In \cite[Sections 6.5 and 6.6]{cerberusIII} $\permcat$ is denoted $\permcatsu$, with the superscript \textsf{su} referring to the strict unity conditions in \cref{def:nlinearfunctor,def:nlineartransformation}.  Next we describe the $\Cat$-enriched multicategory structure on $\permcat$, beginning with the symmetric group action.

\subsection*{Symmetric Group Action}

Suppose given
\begin{itemize}
\item $n$-linear functors $F,G \cn \prod_{j=1}^n \C_j \to \D$ and
\item an $n$-linear natural transformation $\theta \cn F \to G$ as in \cref{nlineartransformation}.
\end{itemize} 
For a permutation $\sigma \in \Sigma_n$, the \emph{symmetric group action functor}\index{symmetric group!action!permutative category}\index{permutative category!symmetric group action}
\begin{equation}\label{permcatsymgroupaction}
\begin{tikzcd}[column sep=large]
\permcat\mmap{\D; \ang{\C}} \ar{r}{\sigma}[swap]{\iso} & \permcat\mmap{\D; \ang{\C}\sigma}
\end{tikzcd}
\end{equation}
sends the $n$-linear functors and natural transformation in \cref{nlineartransformation} to the following composite functors and whiskered natural transformation.
\begin{equation}\label{permcatsigmaaction}
\begin{tikzpicture}[xscale=4,yscale=3,vcenter]
\draw[0cell=.9]
(0,0) node (a) {\prod_{j=1}^n \C_j}
(a)++(.06,.03) node (a') {\phantom{\D}}
(a')++(.6,0) node (b) {\D}
(a)++(-.6,0) node (c) {\prod_{j=1}^n \C_{\sigma(j)}}
;
\draw[1cell=.9]
(c) edge node {\sigma} (a)
;
\draw[1cell=.8]  
(a') edge[bend left] node[pos=.47] {(F, \ang{F^2_j})} (b)
(a') edge[bend right] node[swap,pos=.47] {(G, \ang{G^2_j})} (b)
;
\draw[2cell] 
node[between=a' and b at .45, rotate=-90, 2label={above,\theta}] {\Rightarrow}
;
\end{tikzpicture}
\end{equation}
For each $j \in \{1,\ldots,n\}$, the $j$-th linearity constraint of $F^\sigma = F \circ \sigma$ is the composite in $\D$
\begin{equation}\label{fsigmatwoj}
\begin{tikzcd}[column sep=large, row sep=small]
\Fsigma \ang{A} \oplus \Fsigma \ang{A \compj A_j'} \ar{r}{(\Fsigma)^2_j} \ar[equal,shorten <=-1ex]{d} 
& \Fsigma \ang{A \compj (A_j \oplus A_j')} \ar[equal,shorten <=-1ex]{d}\\
F\big( \sigma\ang{A} \big) \oplus F\big( \sigma\ang{A} \comp_{\sigma(j)} A_j' \big) \ar{r}{F^2_{\sigma(j)}} & F\big( \sigma\ang{A} \comp_{\sigma(j)} (A_j \oplus A_j') \big)
\end{tikzcd}
\end{equation}
for objects 
\[\ang{A} = \ang{A_j}_{j=1}^n \in \txprod_{j=1}^n \C_{\sigma(j)} \andspace A_j' \in \C_{\sigma(j)}.\]  
If $F$ is strong---that is, each $F^2_j$ is a natural isomorphism---then $F^\sigma$ is also strong.

\subsection*{Composition}

For the composition, suppose given, for each $j \in \{1,\ldots,n\}$,
\begin{itemize}
\item small permutative categories $\ang{\B_j} = \ang{\B_{ji}}_{i=1}^{k_j}$ with $k_j \geq 0$,
\item $k_j$-linear functors $F'_j, G'_j \cn \prod_{i=1}^{k_j} \B_{ji} \to \C_j$, and
\item a $k_j$-linear natural transformation $\theta_j \cn F_j' \to G_j'$ as follows.
\end{itemize} 
\begin{equation}\label{permcatbcdata}
\begin{tikzpicture}[xscale=2.2,yscale=2,vcenter]
\draw[0cell=.9]
(0,0) node (a) {\textstyle\prod_{i=1}^{k_j} \B_{ji}}
(a)++(.17,-.03) node (b) {\phantom{\C_j}}
(b)++(1,0) node (c) {\C_j}
;
\draw[1cell=.85]  
(b) edge[bend left] node[pos=.48] {F'_j} (c)
(b) edge[bend right] node[swap,pos=.48] {G'_j} (c)
;
\draw[2cell] 
node[between=b and c at .45, rotate=-90, 2label={above,\theta_j}] {\Rightarrow}
;
\end{tikzpicture}
\end{equation}
With $\ang{\B} = \ang{\ang{\B_j}}_{j=1}^n$, the \emph{multicategorical composition functor}\index{composition!permutative category}\index{permutative category!composition}
\begin{equation}\label{permcatgamma}
\begin{tikzcd}[cells={nodes={scale=.9}}]
\permcat\mmap{\D;\ang{\C}} \times \txprod_{j=1}^n \permcat\mmap{\C_j;\ang{\B_j}} \ar{r}{\gamma}
& \permcat\mmap{\D;\ang{\B}}
\end{tikzcd}
\end{equation}
sends the data \cref{nlineartransformation,permcatbcdata} to the following composite functors and horizontal composite natural transformations.
\begin{equation}\label{permcatcomposite}
\begin{tikzpicture}[xscale=1,yscale=1,vcenter]
\def\a{30} \def\b{30}
\draw[0cell=.9]
(0,0) node (x11) {\txprod_{j=1}^n \txprod_{i=1}^{k_j} \B_{ji}}
(x11)++(.7,0) node (s1) {\phantom{X}}
(s1)++(2,0) node (s2) {\phantom{X}}
(s2)++(.4,-.02) node (x12) {\txprod_{j=1}^n \C_j}
(x12)++(.4,0) node (s3) {\phantom{X}}
(s3)++(2,0) node (x13) {\D}
;
\draw[1cell=.85]  
(s1) edge[bend left=\a] node {\smallprod_j F_j'} (s2)
(s1) edge[bend right=\a] node[swap] {\smallprod_j G_j'} (s2)
(s3) edge[bend left=\b] node {F} (x13)
(s3) edge[bend right=\b] node[swap] {G} (x13)
;
\draw[2cell=.9]
node[between=s1 and s2 at .37, rotate=-90, 2label={above,\smallprod_j \theta_j}] {\Rightarrow}
node[between=s3 and x13 at .45, rotate=-90, 2label={above,\theta}] {\Rightarrow}
;
\end{tikzpicture}
\end{equation}

If any $k_j=0$, then the empty product $\prod_{i=1}^{k_j} \B_{ji}$ is, by definition, the terminal category $\boldone$.  In this case,
\begin{itemize}
\item $F_j'$ is determined by the object $F_j'* \in \C_j$, and
\item $\theta_j$ is determined by the component morphism 
\[(\theta_j)_* \cn F_j'* \to G_j'* \inspace \C_j.\]
\end{itemize}
In \cref{permcatcomposite}
\begin{itemize}
\item the leftmost product should omit those indices $j$ with $k_j=0$, and 
\item $F_j'*$, $G_j'*$, and $(\theta_j)_*$ are used in $\C_j$ in the middle product.
\end{itemize} 
If $k_j=0$ for all $j \in \{1,\ldots,n\}$, then the leftmost product in \cref{permcatcomposite} is $\boldone$.

\subsection*{Composite Linearity Constraints}

To describe the linearity constraints of the composite functor $F \circ \prod_j F'_j$ in \cref{permcatcomposite}, suppose given objects
\begin{itemize}
\item $\ang{W_j} = \ang{W_{ji}}_{i=1}^{k_j} \in \prod_{i=1}^{k_j} \B_{ji}$ for each $j \in \{1,\ldots,n\}$,
\item $\ang{W} = \ang{\ang{W_j}}_{j=1}^n \in \prod_{j=1}^n \prod_{i=1}^{k_j} \B_{ji}$,
\item $W_{ji}' \in \B_{ji}$ for some choice of $(j,i)$ with $\ell = k_1 + \cdots + k_{j-1} + i$, and
\item $\ang{F'W} = \ang{F'_j\ang{W_j}}_{j=1}^n \in \prod_{j=1}^n \C_j$.
\end{itemize} 
The following objects appear in \cref{ffjlinearity} below.
\[\begin{split}
\ang{W_j \compi W_{ji}'} &= \big(\overbracket[.5pt]{W_{j1},\ldots, W_{j,i-1}}^{\text{empty if $i=1$}},\, W_{ji}' \, , \overbracket[.5pt]{W_{j,i+1},\ldots, W_{jk_j}}^{\text{empty if $i=k_j$}}\big)\\
\bang{W_j \compi (W_{ji} \oplus W_{ji}')} &= \big(\underbracket[.5pt]{W_{j1},\ldots, W_{j,i-1}}_{\text{empty if $i=1$}},\, W_{ji} \oplus W_{ji}' \, , \underbracket[.5pt]{W_{j,i+1},\ldots, W_{jk_j}}_{\text{empty if $i=k_j$}}\big).
\end{split}\]
The \emph{$\ell$-th linearity constraint} $\big(F \circ \txprod_j F'_j\big)^2_\ell$ is defined as the following composite morphism in $\D$.
\begin{equation}\label{ffjlinearity}
\begin{tikzpicture}[xscale=3,yscale=1,vcenter]
\def\h{1} \def\v{-1}
\draw[0cell=.75] 
(0,0) node (x0) {F\ang{F'W} \oplus F\bang{F'W \compj F'_j\ang{W_j \compi W_{ji}'}}}
(x0)++(-\h,\v) node (x1) {(F \circ \txprod_j F'_j)\ang{W} \oplus (F \circ \txprod_j F'_j)\ang{W \comp_\ell W_{ji}'}}
(x1)++(1.3*\h,\v) node (x2) {F\bang{F'W \compj (F'_j\ang{W_j} \oplus F'_j\ang{W_j \compi W_{ji}'})}} 
(x1)++(0,2*\v) node (x3) {(F \circ \txprod_j F'_j)\bang{W \comp_\ell (W_{ji} \oplus W_{ji}')}}
(x3)++(\h,\v) node (x4) {F\bang{F'W \compj F'_j\ang{W_j \compi (W_{ji} \oplus W_{ji}')}}}
;
\draw[1cell=.75] 
(x1) edge[-,double equal sign distance, shorten <=1ex, transform canvas={xshift=-5ex}] (x0)
(x0) edge node[pos=.7] {\Ftwoj} (x2)
(x2) edge node[pos=.3] {F\ang{1 \compj (F'_j)^2_i}} (x4)
(x1) edge node[swap] {\big(F \circ \txprod_j F'_j\big)^2_\ell} (x3)
(x3) edge[-,double equal sign distance, shorten >=1ex, transform canvas={xshift=-4ex}] (x4)
; 
\end{tikzpicture}
\end{equation}
If $F$ and each $F'_j$ are strong---that is, each of their linearity constraints is a natural isomorphism---then $F \circ \txprod_j F'_j$ is also strong.

A detailed proof of the following theorem for $\permcat$ is in \cite[Section 6.6]{cerberusIII}.  That proof also applies to the strong variant $\permcatsg$, with component categories defined in \cref{permcatsgangcd,permcatemptyd}, by restricting to strong multilinear functors.

\begin{theorem}\label{thm:permcatenrmulticat}\index{permutative category!Cat-multicategory@$\Cat$-multicategory}\index{Cat-multicategory@$\Cat$-multicategory!of permutative categories}\index{multicategory!of permutative categories}
There is a $\Cat$-multicategory $\permcat$ with
\begin{itemize}
\item small permutative categories (\cref{def:symmoncat}) as objects,
\item multimorphism categories $\permcatsmangcd$ as in \cref{permcatangcd},
\item identity 1-linear functors (\cref{ex:idonelinearfunctor}) as colored units,
\item symmetric group action functor as in \cref{permcatsymgroupaction}, and
\item composition functor as in \cref{permcatgamma}.
\end{itemize}
Moreover, the strong variant $\permcatsg$\label{not:permcatsgcatm} is also a $\Cat$-multicategory.
\end{theorem}

\begin{remark}\label{rk:permcatnotmonoidal}
The $\Cat$-multicategory structure on $\permcat$ and $\permcatsg$ are \emph{not} induced by some symmetric monoidal structure on the $\Cat$-multicategory of small pointed multicategories.  This is explained in detail in \cite[5.7.23 and 10.2.17]{cerberusIII}.  The root cause for the failure of $\permcat$ and $\permcatsg$ to admit a monoidal structure is the non-existence of a suitable monoidal unit.  This is similar to the reason for the non-existence of a monoidal structure on $\DCat$ (\cref{expl:dcatbipermzero}).
\end{remark}

\section{The Grothendieck Construction of Additive Symmetric Monoidal Functors}
\label{sec:groadditivesmf}

Now we begin to prove that the Grothendieck construction is a pseudo symmetric $\Cat$-multifunctor when the indexing category is a small tight bipermutative category.  As the first step, in this section we review
\begin{itemize}
\item the Grothendieck construction $\grod X$ of an indexed category $X \cn \cD \to \Cat$ (\cref{def:diagramcat}) and 
\item its permutative structure when $(\cD, \oplus, \zero, \betaplus)$ is a small permutative category with $X$ a symmetric monoidal functor.
\end{itemize} 
The assignment $X \mapsto\! \grod X$ in the permutative case is the object assignment of the pseudo symmetric $\Cat$-multifunctor $\grod$.  The Grothendieck constructions of additive natural transformations and additive modifications, along with the pseudo symmetric $\Cat$-multifunctoriality, are discussed in subsequent sections.

\subsection*{The Grothendieck Construction as a Category}

A detailed reference for the Grothendieck construction is \cite[Ch.\! 10]{johnson-yau}, which discusses the more general situation when $X$ is a lax functor, or a pseudofunctor for some results.  That reference actually discusses the \emph{contravariant}\index{Grothendieck construction!contravariant} Grothendieck construction for $\cD^\op \to \Cat$.  All the results there apply here as well by considering the opposite category of $\cD^\op$, that is, $\cD$.  When we cite results in \cite[Ch.\! 10]{johnson-yau}, we are referring to the \emph{covariant} version.

\begin{definition}\label{def:grothendieckconst}
Suppose $\cD$ is a small category, and $X \cn \cD \to \Cat$ is a $\cD$-indexed category.  The \emph{Grothendieck construction}\index{Grothendieck construction} of $X$, denoted $\grod X$, is the category defined by the following data.
\begin{description}
\item[Objects] An \emph{object} in $\grod X$ is a pair of objects 
\[(a,x) \in \cD \times Xa.\]
\item[Morphisms] A \emph{morphism}
\[\begin{tikzcd}[column sep=large]
(a,x) \ar{r}{(f,p)} & (b,y) 
\end{tikzcd}\inspace \grod X\]
consists of
\begin{itemize}
\item a morphism $f \cn a \to b$ in $\cD$ and
\item a morphism $p \cn f_*x \to y$ in $Xb$ with the functor 
\[f_* = Xf \cn Xa \to Xb.\]
\end{itemize}
\item[Identities] For each object $(a,x)$ in $\grod X$, the \emph{identity morphism} consists of a pair of identity morphisms: 
\begin{equation}\label{groconidentity}
1_{(a,x)} = (1_a, 1_x).
\end{equation}
\item[Composition] For another morphism 
\[\begin{tikzcd}[column sep=large]
(b,y) \ar{r}{(g,q)} & (c,z) 
\end{tikzcd}\inspace \grod X\]
the \emph{composite} of 
\[\begin{tikzcd}[column sep=large]
(a,x) \ar{r}{(f,p)} & (b,y) \ar{r}{(g,q)} & (c,z) 
\end{tikzcd}\]
is given by
\begin{equation}\label{groconcomposition}
(g,q) \circ (f,p) = \big(gf, q \circ g_* p\big).
\end{equation}
\end{description}
This finishes the definition of $\grod X$.
\end{definition}

The reader should not confuse the notation $\grod X$ for the Grothendieck construction with the coend in \cref{eq:dayconvolution} or the end in \cref{eq:dayhom}.  The following observation is a special case of \cite[10.1.8]{johnson-yau}.

\begin{lemma}\label{groconstcategory}
In the context of \cref{def:grothendieckconst}, $\grod X$ is a small category.
\end{lemma}

\subsection*{The Grothendieck Construction as a Permutative Category}

Next we discuss the permutative variant of the Grothendieck construction.  A reference is \cite[4.29]{gjo1}, which discusses the more general situation when the codomain is $\iicat$ instead of $\Cat$.  Recall from \cref{def:monoidalfunctor} the notion of a symmetric monoidal functor between symmetric monoidal categories.

\begin{definition}\label{def:permutativegroconst}\index{Grothendieck construction!symmetric monoidal functor}\index{symmetric monoidal functor!Grothendieck construction}
Suppose $(\cD, \oplus, \zero, \betaplus)$ is a small permutative category, and 
\[(X,X^2,X^0) \cn \cD \to \Cat\]
is a symmetric monoidal functor.  Define the data of a permutative category on the Grothendieck construction $\grod X$ as follows.
\begin{description}
\item[Monoidal Unit] It is the object 
\begin{equation}\label{groconunit}
(\zero, X^0*) \in \grod X
\end{equation} 
where $X^0* \in X\zero$ is the unit object \cref{additivesmfunitobj}.
\item[Monoidal Product on Objects] For objects $(a,x)$ and $(a',x')$ in $\grod X$, their monoidal product is defined as the object
\begin{equation}\label{gboxobjects}
(a,x) \gbox (a',x') = \big(a \oplus a', X^2_{a,a'}(x,x')\big) \inspace \grod X
\end{equation}
with 
\[\begin{tikzcd}[column sep=large]
Xa \times Xa' \ar{r}{X^2_{a,a'}} & X(a \oplus a')
\end{tikzcd}\]
the $(a,a')$-component functor of the monoidal constraint $X^2$ \cref{additivesmfconstraint}. 
\item[Monoidal Product on Morphisms] For morphisms
\[\begin{tikzcd}
(a,x) \ar{r}{(f,p)} & (b,y) 
\end{tikzcd} \andspace
\begin{tikzcd}
(a',x') \ar{r}{(f',p')} & (b',y') 
\end{tikzcd} \inspace \grod X\]
their monoidal product is defined as the following morphism.
\begin{equation}\label{gboxmorphisms}
\begin{tikzcd}[column sep=2.8cm, row sep=small, cells={nodes={scale=.9}}]
(a,x) \gbox (a',x') \ar[equal, transform canvas={xshift=-.6ex}, shorten <=-.5ex]{d} \ar{r}{(f,p) \gbox (f',p')} 
& (b,y) \gbox (b',y') \ar[equal, transform canvas={xshift=-.6ex}, shorten <=-.5ex]{d}\\
\big(a \oplus a', X^2_{a,a'}(x,x')\big) \ar{r}{\left(f \oplus f' \scs X^2_{b,b'}(p,p')\right)} 
& \big(b \oplus b', X^2_{b,b'}(y,y')\big) 
\end{tikzcd}
\end{equation}
\item[Braiding] The component of the braiding $\betabox$ for objects $(a,x)$ and $(a',x')$ in $\grod X$ is defined as the following isomorphism.
\begin{equation}\label{groconbeta}
\begin{tikzcd}[column sep=large, row sep=tiny]
(a,x) \gbox (a',x') \ar[equal, transform canvas={xshift=-.6ex}, shorten <=-.5ex]{d} \ar{r}{\betabox} 
& (a',x') \gbox (a,x) \ar[equal, transform canvas={xshift=.7ex}, shorten <=-.5ex]{d}\\
\big(a \oplus a', X^2_{a,a'}(x,x')\big) \ar{r}{(\betaplus,1)}[swap]{\iso} 
& \big(a' \oplus a, X^2_{a',a}(x',x)\big)
\end{tikzcd}
\end{equation}
\end{description}
This finishes the definition.
\end{definition}

\begin{lemma}\label{grodxpermutative}\index{Grothendieck construction!permutative category}\index{permutative category!Grothendieck construction}
In the context of \cref{def:permutativegroconst}, the quadruple
\[\big(\grod X, \gbox, (\zero, X^0*), \betabox\big)\]
is a small permutative category.
\end{lemma}

\begin{proof}
By \cref{groconstcategory} $\grod X$ is a small category.  The functoriality of 
\[\begin{tikzcd}[column sep=large]
\grod X \times \grod X \ar{r}{\gbox} & \grod X
\end{tikzcd}\]
follows from
\begin{itemize}
\item the functoriality of the monoidal product $\oplus \cn \cD \times \cD \to \cD$,
\item the functoriality of each component functor of the monoidal constraint $X^2$ \cref{additivesmfconstraint}, and
\item the naturality of $X^2$ \cref{Xtwonaturality}.
\end{itemize}
The associativity of $\gbox$ follows from 
\begin{itemize}
\item the associativity of $\oplus$ in $\cD$ and
\item the associativity axiom \cref{asmfassociativity} for $(X,X^2)$.
\end{itemize}
The left and right unity of $\gbox$ follows from
\begin{itemize}
\item the unity of $\zero \in \cD$ and
\item the unity axiom \cref{asmfunity} for $(X,X^2,X^0)$.
\end{itemize}
Thus $\big(\grod X, \gbox, (\zero, X^0*)\big)$ is a strict monoidal category.

The braiding $\betabox$ \cref{groconbeta} is well defined by the braiding axiom \cref{asmfbraiding} for $(X,X^2)$.  The naturality of $\betabox$ follows from
\begin{itemize}
\item the naturality of $\betaplus$ in $\cD$ and
\item the braiding axiom \cref{asmfbraiding} for $(X,X^2)$.
\end{itemize}   
The symmetry and hexagon axioms \cref{symmoncatsymhexagon} for $(\grod X, \betabox)$ follow from
\begin{itemize}
\item the same axioms in the permutative category $\cD$ and
\item the functoriality of each component functor of $X^2$ \cref{additivesmfconstraint}.
\end{itemize}
Thus $(\grod X, \betabox)$ is a strict symmetric monoidal category.
\end{proof}

In what follows, we regard the Grothendieck construction $\grod X$ as a permutative category via  \cref{grodxpermutative}.  Note that, while $\cD$ is required to be a small permutative category, \emph{neither} $X^0$ nor $X^2$ is required to be an isomorphism.

\section{The Grothendieck Construction of Additive Natural Transformations}
\label{sec:groadditivenattr}

For the rest of this chapter, we assume that 
\[\big(\cD, (\oplus, \zero, \betaplus), (\otimes, \tu, \betate), (\fal, \far)\big)\]
is a small \emph{tight} bipermutative category (\cref{def:embipermutativecat}) with 
\begin{itemize}
\item additive structure $\Dplus = (\cD, \oplus, \zero, \betaplus)$ and
\item multiplicative structure $\Dte = (\cD, \otimes, \tu, \betate)$.
\end{itemize}
Recall the $\Cat$-multicategory $\DCat$ in \cref{thm:dcatcatmulticat} with
\begin{itemize}
\item additive symmetric monoidal functors $\Dplus \to \Cat$ as objects,
\item additive natural transformations as objects in the multimorphism categories, and
\item additive modifications as morphisms in the multimorphism categories.
\end{itemize}
Also recall the $\Cat$-multicategory $\permcatsg$ in \cref{thm:permcatenrmulticat} with
\begin{itemize}
\item small permutative categories as objects,
\item \emph{strong} multilinear functors as objects in the multimorphism categories, and
\item multilinear natural transformations as morphisms in the multimorphism categories.
\end{itemize}
This section contains the next step in proving that the Grothendieck construction 
\[\begin{tikzcd}[column sep=large,every label/.append style={scale=.85}]
\DCat \ar{r}{\grod} & \permcatsg
\end{tikzcd}\]
is a pseudo symmetric $\Cat$-multifunctor.  The object assignment $X \mapsto\! \grod X$ is given by \cref{grodxpermutative}.  In this section we show that the Grothendieck construction sends each additive natural transformation between additive symmetric monoidal functors to a strong multilinear functor; see \cref{grodphifunctor}.

\subsection*{Arity 0}

First we discuss the arity 0 case, which has its own definition and is easier than the positive arity case.

\begin{definition}\label{def:groconarityzero}
Suppose given an additive symmetric monoidal functor (\cref{def:additivesmf})
\[(Z,Z^2,Z^0) \cn \Dplus \to \Cat.\]
Using the definitions \cref{dcatbipermzeroary,permcatemptyd} of arity 0 multimorphism categories in, respectively, $\DCat$ and $\permcatsg$, we define the functor\index{Grothendieck construction!0-ary multimorphism functor}
\begin{equation}\label{groarityzero}
\begin{tikzcd}[column sep=large,every label/.append style={scale=.85}]
Z\tu = \DCat\scmap{\ang{};Z} \ar{r}{\grod} & \permcatsg\scmap{\ang{};\grod Z} = \grod Z
\end{tikzcd}
\end{equation}
by sending
\begin{itemize}
\item each object $x \in Z\tu$ to the object
\[(\tu,x) \in \grod Z\]
and
\item each morphism $f \cn x \to y$ in $Z\tu$ to the morphism
\[\begin{tikzcd}[column sep=large]
(\tu,x) \ar{r}{(1_{\tu},f)} & (\tu,y)
\end{tikzcd} \inspace \grod Z.\]
\end{itemize} 
These object and morphism assignments define a functor because $Z(1_{\tu}) = 1_{Z\tu}$ by the functoriality of $Z$.
\end{definition}

As discussed in
\begin{itemize}
\item the first paragraph in the proof of \cref{dcatgammaclosed} and
\item the paragraph after \cref{permcatcomposite},
\end{itemize}  
proofs involving the arity 0 case require some slight adjustments.  With this in mind, in the rest of this chapter we discuss exclusively the positive arity case, using \cref{def:groconarityzero} to make the necessary adjustments when some inputs have arity 0.

\subsection*{Positive Arity}

Next we define the Grothendieck construction on additive natural transformations (\cref{def:additivenattr}).  Recall that an \emph{additive natural transformation} is a natural transformation in $\Dtecat$ \cref{phitensorxz} that satisfies the unity axiom \cref{additivenattrunity} and the additivity axiom \cref{additivenattradditivity}.  Also recall from \cref{def:nlinearfunctor} that an \emph{$n$-linear functor} is a functor equipped with $n$ linearity constraints \cref{laxlinearityconstraints} that satisfy the axioms \cref{nlinearunity,constraintunity,eq:ml-f2-assoc,eq:ml-f2-symm,eq:f2-2by2}.

\begin{definition}\label{def:groconarityn}
Suppose 
\[(Z,Z^2,Z^0) \andspace \bang{(X_j,X_j^2,X_j^0)}_{j=1}^n \cn \Dplus \to \Cat\]
are additive symmetric monoidal functors with $n>0$.  Suppose 
\[\begin{tikzcd}[column sep=large]
\bigotimes_{j=1}^n X_j \ar{r}{\phi} & Z
\end{tikzcd}\]
is an additive natural transformation \cref{additivenattransformation}.  Define the data of an $n$-linear functor\index{additive!natural transformation!Grothendieck construction}\index{Grothendieck construction!additive natural transformation}
\begin{equation}\label{groadditivenattr}
\begin{tikzcd}[column sep=huge, every label/.append style={scale=.9}]
\prod_{j=1}^n \big(\grod X_j\big) \ar{r}{\grod \phi} & \grod Z
\end{tikzcd}
\end{equation}
with the following object assignment \cref{grodphiangajxj}, morphism assignment \cref{grodphiangfjpj}, and linearity constraints \cref{grodphilinearity}.
\begin{description}
\item[Object Assignment] Suppose given an $n$-tuple of objects
\begin{equation}\label{angajxj}
\ang{(a_j,x_j)} \in \txprod_{j=1}^n \big(\grod X_j\big)
\end{equation}
consisting of objects $a_j \in \cD$ and $x_j \in X_ja_j$ for $j \in \{1,\ldots,n\}$.  First define the following objects.
\begin{equation}\label{atensorangxj}
\left\{\begin{aligned}
\anga &= \ang{a_j}_{j=1}^n \in \cD^n & a &= \txotimes_{j=1}^n a_j \in \cD\\
\angx &= \ang{x_j}_{j=1}^n \in \txprod_{j=1}^n X_j a_j && 
\end{aligned}\right.
\end{equation}
Then $\grod \phi$ sends $\ang{(a_j,x_j)}$ to the object
\begin{equation}\label{grodphiangajxj}
(\grod \phi) \ang{(a_j,x_j)} = \brb{a, \phi_{\anga} \angx} \inspace \grod Z
\end{equation}
with $\phi_{\anga}$ the $\anga$-component functor of $\phi$ as in \cref{dcatnaryobjcomponent}:
\[\begin{tikzcd}[column sep=large]
\prod_{j=1}^n X_j a_j \ar{r}{\phi_{\anga}} & Za.
\end{tikzcd}\]
\item[Morphism Assignment]
Suppose given an $n$-tuple of morphisms
\begin{equation}\label{angphijfj}
\bang{(f_j, p_j) \cn (a_j,x_j) \to (b_j,y_j)} \in \txprod_{j=1}^n \big(\grod X_j\big)
\end{equation}
consisting of, for each $j \in \{1,\ldots,n\}$, morphisms
\begin{itemize}
\item $f_j \cn a_j \to b_j$ in $\cD$ and
\item $p_j \cn (f_j)_* x_j \to y_j$ in $X_j b_j$ with 
\[(f_j)_* = Xf_j \cn Xa_j \to Xb_j.\]
\end{itemize} 
In addition to the objects in \cref{atensorangxj}, define the following objects and morphisms.
\begin{equation}\label{angbyp}
\left\{\begin{aligned}
\angb &= \ang{b_j}_{j=1}^n \in \cD^n & b &= \txotimes_{j=1}^n b_j \in \cD\\
\angy &= \ang{y_j}_{j=1}^n \in \txprod_{j=1}^n X_j b_j & f &= \txotimes_{j=1}^n f_j \in \cD\\
\angp &= \ang{p_j}_{j=1}^n \in \txprod_{j=1}^n X_j b_j &&
\end{aligned}\right.
\end{equation}
Then $\grod \phi$ sends $\ang{(f_j,p_j)}$ to the following morphism in $\grod Z$.  
\begin{equation}\label{grodphiangfjpj}
\begin{tikzpicture}[xscale=5,yscale=1,vcenter]
\draw[0cell=.9]
(0,0) node (x11) {(\grod \phi)\ang{(a_j,x_j)}}
(x11)++(1,0) node (x12) {(\grod \phi)\ang{(b_j,y_j)}} 
(x11)++(0,-1) node (x21) {\brb{a,\phi_{\anga}\angx}}
(x12)++(0,-1) node (x22) {\brb{b,\phi_{\angb}\angy}}
;
\draw[1cell=.9]  
(x11) edge node {(\grod \phi)\ang{(f_j,p_j)}} (x12)
(x21) edge node {\srb{f,\phi_{\angb}\angp}} (x22)
(x11) edge[equal] (x21)
(x12) edge[equal] (x22)
;
\draw[2cell] 
;
\end{tikzpicture}
\end{equation}
This is well defined by the naturality diagram \cref{dcatnaryobjnaturality} for $\phi$.
\item[Linearity Constraints]
To define the linearity constraints of $\grod \phi$, suppose given an index $i \in \{1,\ldots,n\}$ and an object 
\[(a_i', x_i') \in \grod X_i\]
with $a_i' \in \cD$ and $x_i' \in X_i a_i'$.  In addition to the objects in \cref{angajxj,atensorangxj}, define the following objects using the $\compi$ notation in \cref{compnotation}.
\begin{equation}\label{grodphiconstraintnot}
\left\{\begin{aligned}
a_i'' &= a_i \oplus a_i' \in \cD\\
a' &= a_1 \otimes \cdots \otimes a_i' \otimes \cdots \otimes a_n \in \cD \qquad \angap = \anga \compi a_i' \in \cD^n\\
a'' &= a_1 \otimes \cdots \otimes a_i'' \otimes \cdots \otimes a_n \in \cD \qquad \angapp = \anga \compi a_i'' \in \cD^n\\
x_i'' &= (X_i^2)_{a_i,a_i'} (x_i,x_i') \in X_i a_i''\\
\angxp &= \angx \compi x_i' \in X_1 a_1 \times \cdots \times X_i a_i' \times \cdots \times X_n a_n\\
\angxpp &= \angx \compi x_i'' \in X_1 a_1 \times \cdots \times X_i a_i'' \times \cdots \times X_n a_n
\end{aligned}\right.
\end{equation}
The above definition of $x_i''$ uses the monoidal constraint \cref{additivesmfconstraint} of $X_i$:
\[\begin{tikzcd}[column sep=huge]
X_i a_i \times X_i a_i' \ar{r}{(X_i^2)_{a_i,a_i'}} & X_i (a_i \oplus a_i') = X_i a_i''.
\end{tikzcd}\]
For the objects $\ang{(a_j,x_j)}$ and $(a_i',x_i')$, the $i$-th linearity constraint \cref{laxlinearityconstraints} of $\grod \phi$ has the following domain object in $\grod Z$.
\begin{equation}\label{grodphiconstraintdom}
\left\{\begin{aligned}
& (\grod \phi) \ang{(a_j,x_j)} \gbox (\grod \phi) \big( \ang{(a_j,x_j)} \compi (a_i',x_i')\big) &&\\
&= \brb{a,\phi_{\anga}\angx} \gbox \brb{a',\phi_{\angap}\angxp} & \phantom{=} & \text{by \cref{grodphiangajxj}}\\
&= \left(a \oplus a' \scs Z^2_{a,a'}(\phi_{\anga}\angx \scs \phi_{\angap}\angxp)\right) & \phantom{=} & \text{by \cref{gboxobjects}}\\
&= \left(a \oplus a' \scs \lap_* \phi_{\angapp} \angxpp\right) & \phantom{=} & \text{by \cref{additivityobjects}}
\end{aligned}\right.
\end{equation}
In the previous line, $\lap$ is the unique Laplaza coherence isomorphism\index{Laplaza coherence isomorphism} in $\cD$
\[\begin{tikzcd}[column sep=large]
a'' \ar{r}{\lap}[swap]{\iso} & a \oplus a'
\end{tikzcd}\]
in \cref{laplazaapp} that distributes over the sum $a_i \oplus a_i'$ in $a''$.  

The codomain object in $\grod Z$ of the same $i$-th linearity constraint of $\grod \phi$ is the following.
\begin{equation}\label{grodphiconstraintcod}
\left\{\begin{aligned}
& \scalebox{.9}{$(\grod \phi) \Big( \ang{(a_j,x_j)} \compi \big((a_i,x_i) \gbox (a_i',x_i') \big) \Big)$} &&\\
&= (\grod \phi) \Big( \ang{(a_j,x_j)} \compi \srb{a_i'', x_i''}\Big) & & \text{by \cref{gboxobjects,grodphiconstraintnot}}\\
&= \brb{a'',\phi_{\angapp} \angxpp} & & \text{by \cref{grodphiangajxj}}
\end{aligned}\right.
\end{equation}
Using the last expression in each of \cref{grodphiconstraintdom,grodphiconstraintcod}, the \emph{$i$-th linearity constraint} $(\grod \phi)^2_i$ of $\grod \phi$ \cref{laxlinearityconstraints} is the isomorphism in $\grod Z$
\begin{equation}\label{grodphilinearity}
\begin{tikzcd}[column sep=huge]
\left(a \oplus a' \scs \lap_* \phi_{\angapp} \angxpp\right) \ar{r}{(\lap^\inv, 1)}[swap]{\iso} 
& \brb{a'',\phi_{\angapp} \angxpp}
\end{tikzcd}
\end{equation}
with
\[\begin{tikzcd}
(\lap^\inv)_* \lap_* \phi_{\angapp} \angxpp = \phi_{\angapp} \angxpp \ar{r}{1} & \phi_{\angapp} \angxpp
\end{tikzcd}\]
the identity morphism of $\phi_{\angapp} \angxpp$ in $Za''$.
\end{description}
This finishes the definition of $\grod \phi$ in \cref{groadditivenattr}.
\end{definition}

\begin{lemma}\label{grodphifunctor}\index{Grothendieck construction!strong multilinear functor}\index{multilinear functor!Grothendieck construction}
In the context of \cref{def:groconarityn}, the data
\[\begin{tikzcd}[column sep=3.4cm, every label/.append style={scale=.9}]
\prod_{j=1}^n \big(\grod X_j\big) \ar{r}{\brb{\grod \phi, \ang{(\grod \phi)_i^2}_{i=1}^n}} & \grod Z
\end{tikzcd}\]
form a strong $n$-linear functor.
\end{lemma}

\begin{proof}
By \cref{def:nlinearfunctor} we must verify the following statements.
\begin{enumerate}[label=(\roman*)]
\item\label{grodphifunctor-i} $\grod \phi$, as defined in \cref{grodphiangajxj,grodphiangfjpj}, is a functor.
\item\label{grodphifunctor-ii} For each $i \in \{1,\ldots,n\}$, the $i$-th linearity constraint $(\grod \phi)_i^2$ in \cref{grodphilinearity} is a natural isomorphism.
\item\label{grodphifunctor-iii} The axioms \cref{nlinearunity,constraintunity,eq:ml-f2-assoc,eq:ml-f2-symm,eq:f2-2by2} hold for $\brb{\grod \phi, \ang{(\grod \phi)_i^2}_{i=1}^n}$.
\end{enumerate}
We begin with statement \ref{grodphifunctor-i}.

\medskip
\ref{grodphifunctor-i} \emph{Functoriality}.  First we observe that $\grod \phi$ preserves identity morphisms.  By definition \cref{groconidentity} an identity morphism in the Grothendieck construction consists of a pair of identity morphisms.  In \cref{angphijfj} if for each $j \in \{1,\ldots,n\}$ 
\[\begin{split}
f_j &= 1 \cn a_j \to a_j \in \cD \andspace\\
p_j &= 1 \cn x_j \to x_j \in X_j a_j
\end{split}\]
then 
\[\begin{split}
f &= \txotimes_{j=1}^n f_j = 1_a \andspace\\
\phi_{\angb}\angp &= \phi_{\anga}\ang{1} = 1
\end{split}\]
by the functoriality of $\otimes$ in $\cD$ and the component functor $\phi_{\angb} = \phi_{\anga}$.  Thus $(\grod \phi)\ang{(f_j,p_j)}$ is the identity morphism by definition \cref{grodphiangfjpj}.  

Similarly, $\grod \phi$ preserves composition of morphisms \cref{groconcomposition} by
\begin{itemize}
\item the functoriality of $\otimes$ in $\cD$ and each component functor of $\phi$ and
\item the naturality diagram \cref{dcatnaryobjnaturality} for $\phi$ with respect to morphisms in $\cD$.
\end{itemize}
Thus $\grod \phi$ is a functor.

\medskip
\ref{grodphifunctor-ii} \emph{Naturality}.  Next we show that each linearity constraint $(\grod\phi)_i^2$ in \cref{grodphilinearity}, which is componentwise an isomorphism, is a natural isomorphism.  In addition to the morphisms $\ang{(f_j,p_j)}$ in \cref{angphijfj}, suppose given a morphism
\[\begin{tikzcd}[column sep=large]
(a_i',x_i') \ar{r}{(f_i',p_i')} & (b_i',y_i')
\end{tikzcd} \inspace \grod X_i\]
with
\begin{itemize}
\item $f_i' \cn a_i' \to b_i'$ a morphism in $\cD$ and
\item $p_i' \cn (f_i')_* x_i' \to y_i'$ a morphism in $Xb_i'$. 
\end{itemize} 
We extend the notation in \cref{grodphiconstraintnot} to
\begin{itemize}
\item the objects $\ang{(b_j,y_j)}$ and $(b_i',y_i')$ and
\item the morphisms $\ang{(f_j,p_j)}$ and $(f_i',p_i')$.
\end{itemize} 
Then the naturality of $(\grod\phi)_i^2$ means the commutativity of the following diagram in $\grod Z$.
\begin{equation}\label{grodphinatdiagram}
\begin{tikzpicture}[xscale=4,yscale=1.4,vcenter]
\def\h{1} \def\v{-1}
\draw[0cell=.85] 
(0,0) node (x11) {\brb{a \oplus a', \lap_*\phi_{\angapp}\angxpp}}
(x11)++(1,0) node (x12) {\brb{a'', \phi_{\angapp}\angxpp}}
(x11)++(0,-1) node (x21) {\brb{b \oplus b', \lap_*\phi_{\angbpp}\angypp}}
(x12)++(0,-1) node (x22) {\brb{b'', \phi_{\angbpp}\angypp}}
;
\draw[1cell=.85] 
(x11) edge node {(\lap^\inv,1)} (x12)
(x21) edge node {(\lap^\inv,1)} (x22)
(x11) edge[transform canvas={xshift=2.5em}] node[swap] {\srb{f\oplus f', \lap_*\phi_{\angbpp}\angppp}} (x21)
(x12) edge[transform canvas={xshift=-2em}] node {\srb{f'', \phi_{\angbpp}\angppp}} (x22)
; 
\end{tikzpicture}
\end{equation}
To check the commutativity of the diagram \cref{grodphinatdiagram}, we consider the two arguments separately.
\begin{itemize}
\item In the first argument, the morphism equality
\[\lap^\inv (f \oplus f') = f'' \lap^\inv \cn a \oplus a' \to b'' \inspace \cD\]
follows from the naturality of the Laplaza coherence isomorphism $\lap$ and the definitions of $f$, $f'$, and $f''$.
\item The morphism equality in $Zb''$ in the second argument follows from the following equalities.
\[\begin{split}
1 \circ \lap_*^\inv \lap_* \phi_{\angbpp}\angppp &= \phi_{\angbpp}\angppp\\
&= \phi_{\angbpp}\angppp \circ 1\\
&= \phi_{\angbpp}\angppp \circ  (f_*'' 1)
\end{split}\]
\end{itemize}
Thus each linearity constraint $(\grod\phi)_i^2$ is a natural isomorphism.

\medskip
\ref{grodphifunctor-iii} \emph{Multilinearity}.  Next we check the axioms \cref{nlinearunity,constraintunity,eq:ml-f2-assoc,eq:ml-f2-symm,eq:f2-2by2} for $\grod \phi$.  

\emph{Unity} \cref{nlinearunity}. 
By the definitions \cref{groconidentity,groconunit} of, respectively, an identity morphism and the monoidal unit in the Grothendieck construction, the unity axiom for $\grod \phi$ states the following equalities in $\grod Z$.
\begin{equation}\label{grodphiunityaxiom}
\left\{\begin{split}
(\grod\phi)\Big(\ang{(a_j,x_j)} \compi (\zero, X_i^0*) \Big) &= (\zero, Z^0*)\\
(\grod\phi)\Big(\ang{(f_j,p_j)} \compi \big(1_{\zero}, 1_{X_i^0*}\big) \Big) &= \big(1_{\zero}, 1_{Z^0*}\big)
\end{split}\right.
\end{equation}
The equalities \cref{grodphiunityaxiom} hold by
\begin{itemize}
\item the definitions \cref{grodphiangajxj,grodphiangfjpj} of $\grod\phi$ on, respectively, objects and morphisms,
\item the fact that $\zero \in \cD$ is a strict multiplicative zero \cref{ringcataxiommultzero}, and
\item the unity axiom \cref{additivenattrunity} for the additive natural transformation $\phi$.
\end{itemize} 

\emph{Constraint Unity} \cref{constraintunity}.
This axiom for $\grod\phi$ states that, if either
\[\begin{split}
(a_i',x_i') &= (\zero,X_i^0*) \in \grod X_i \orspace\\
(a_j,x_j) &= (\zero,X_j^0*) \in \grod X_j \foranyspace j \in \{1,\ldots,n\},
\end{split}\]
then the $i$-th linearity constraint $(\lap^\inv,1)$ in \cref{grodphilinearity} is the identity morphism.  This axiom holds because, if either
\[a_i' = \zero \orspace a_j = \zero \foranyspace j \in \{1,\ldots,n\},\]
then the Laplaza coherence isomorphism $\lap \cn a'' \to a \oplus a'$ in \cref{laplazaapp} is the identity morphism by the uniqueness in \cref{thm:laplaza-coherence-1}.

\emph{Other Constraint Axioms}. 
Similar to the constraint unity axiom, each of
\begin{itemize}
\item the constraint associativity axiom \cref{eq:ml-f2-assoc},
\item the constraint symmetry axiom \cref{eq:ml-f2-symm}, and
\item the constraint 2-by-2 axiom \cref{eq:f2-2by2}
\end{itemize}  
for $\grod\phi$ holds by the uniqueness in Laplaza's Coherence \cref{thm:laplaza-coherence-1}.  Indeed, recall that each component of the braiding $\betabox$ in the Grothendieck construction has the form $(\betaplus, 1)$, with
\begin{itemize}
\item $\betaplus$ the braiding in $\Dplus = (\cD,\oplus,\zero,\betaplus)$ and
\item the second argument an identity morphism.
\end{itemize}   
Thus in the diagram for each of the axioms \cref{eq:ml-f2-assoc,eq:ml-f2-symm,eq:f2-2by2} for $\grod\phi$, each composite has second argument given by an identity morphism in some component category of $Z$.  In each case \cref{thm:laplaza-coherence-1} gives the desired morphism equality in the first argument in $\cD$.  For example, in the constraint symmetry axiom \cref{eq:ml-f2-symm} for $\grod\phi$, the first argument is the following diagram in $\cD$, with $\otimes$ abbreviated to juxtaposition to save space.
\[\begin{tikzpicture}[xscale=3,yscale=1.4,vcenter]
\def\h{1} \def\v{-1}
\draw[0cell=.9] 
(0,0) node (x11) {a \oplus a'}
(x11)++(1,0) node (x12) {a_1 \cdots (a_i \oplus a_i') \cdots a_n}
(x11)++(0,-1) node (x21) {a' \oplus a}
(x12)++(0,-1) node (x22) {a_1 \cdots (a_i' \oplus a_i) \cdots a_n}
;
\draw[1cell=.9] 
(x11) edge node {\lap^\inv} (x12)
(x21) edge node {\lap^\inv} (x22)
(x11) edge node[swap] {\betaplus} (x21)
(x12) edge node {1 \cdots \betaplus \cdots 1} (x22)
; 
\end{tikzpicture}\]
This diagram commutes by \cref{thm:laplaza-coherence-1}.
\end{proof}

\begin{remark}[Unary Case]\label{rk:grodphifunctorunary}
In \cref{grodphifunctor} with $n=1$, the pair
\[\big(\grod\phi, (\grod\phi)^2\big) \cn \grod X \to \grod Z\]
is a \emph{strict} symmetric monoidal functor for the following reasons.
\begin{itemize}
\item It is a 1-linear functor, which is the same thing as a strictly unital symmetric monoidal functor (\cref{ex:onelinearfunctor}).
\item The linearity constraint $(\grod\phi)^2$ is componentwise $(\lap^\inv,1)$ in \cref{grodphilinearity}.  This is an identity morphism because the Laplaza coherence isomorphism 
\[\begin{tikzcd}[column sep=large]
a \oplus a' \ar{r}{\lap}[swap]{\iso} & a \oplus a'
\end{tikzcd}\]
is the identity morphism by uniqueness (\cref{thm:laplaza-coherence-1}).
\end{itemize}
For $n > 1$ the linearity constraints $(\grod\phi)^2_i$ are natural isomorphisms but not identities in general.
\end{remark}

\section{The Grothendieck Construction of Additive Modifications}
\label{sec:groadditivemod}

As in \cref{sec:groadditivenattr} we continue to assume that $\cD$ is a small tight bipermutative category with additive structure $\Dplus$ and multiplicative structure $\Dte$.  Recall from \cref{def:additivemodification} that an \emph{additive modification} is a modification \cref{Phiphivarphi} that satisfies the unity axiom \cref{additivemodunity} and the additivity axiom \cref{additivemodadditivity}.  Also recall from \cref{def:nlineartransformation} that an \emph{$n$-linear natural transformation} is a natural transformation that satisfies the unity axiom \cref{ntransformationunity} and the constraint compatibility axiom \cref{eq:monoidal-in-each-variable}.

This section contains the next step in the construction of the pseudo symmetric $\Cat$-multifunctorial Grothendieck construction 
\[\begin{tikzcd}[column sep=large,every label/.append style={scale=.85}]
\DCat \ar{r}{\grod} & \permcatsg.
\end{tikzcd}\]
The object assignment $X \mapsto\! \grod X$ is given by \cref{grodxpermutative}.  The object assignment of the  multimorphism functors is in 
\begin{itemize}
\item \cref{def:groconarityzero} for arity 0 inputs and
\item \cref{def:groconarityn} for positive arity inputs.
\end{itemize}
In this section we show that the Grothendieck construction
\begin{itemize}
\item sends additive modifications in $\DCat$ to multilinear natural transformations in $\permcatsg$ and
\item preserves vertical composition.
\end{itemize}   
Since the multimorphism functors in arity 0 are constructed in \cref{def:groconarityzero}, we only need to consider the positive arity case.  
\begin{itemize}
\item The Grothendieck construction of an additive modification is in \cref{def:groconmodification}.
\item Its naturality and $n$-linearity are proved in \cref{grodPhinatural,grodPhinlinear}.
\item \cref{grodmultimorphismfunctor} shows that $\grod$ defines a functor between multimorphism categories.
\end{itemize}

\begin{definition}\label{def:groconmodification}
Suppose given
\begin{itemize}
\item additive symmetric monoidal functors (\cref{def:additivesmf})
\[(Z,Z^2,Z^0) \andspace \bang{(X_j,X_j^2,X_j^0)}_{j=1}^n \cn \Dplus \to \Cat\]
with $n>0$ and 
\item an additive modification 
\[\begin{tikzpicture}[xscale=1,yscale=1,baseline={(x1.base)}]
\draw[0cell=.9]
(0,0) node (x1) {\bang{(X_j, X_j^2, X_j^0)}_{j=1}^n}
(x1)++(1,0) node (x2) {\phantom{X}}
(x2)++(2,0) node (x3) {\phantom{X}}
(x3)++(.65,0) node (x4) {(Z,Z^2,Z^0)}
;
\draw[1cell=.9]  
(x2) edge[bend left] node[pos=.5] {\phi} (x3)
(x2) edge[bend right] node[swap,pos=.5] {\psi} (x3)
;
\draw[2cell]
node[between=x2 and x3 at .45, rotate=-90, 2label={above,\Phi}] {\Rightarrow}
;
\end{tikzpicture}\]
as in \cref{additivemodtwocell} with $\phi$ and $\psi$ additive natural transformations \cref{additivenattransformation}.
\end{itemize}
Define the data of an $n$-linear natural transformation\index{Grothendieck construction!additive modification}\index{additive!modification!Grothendieck construction}
\begin{equation}\label{grodPhi}
\begin{tikzpicture}[xscale=2.4,yscale=2,baseline={(a.base)}]
\draw[0cell=.9]
(0,0) node (a) {\textstyle\prod_{j=1}^n (\grod X_j)}
(a)++(.28,.01) node (b) {\phantom{D}}
(b)++(1,0) node (c) {\grod Z}
;
\draw[1cell=.85]  
(b) edge[bend left] node[pos=.53] {\grod \phi} (c)
(b) edge[bend right] node[swap,pos=.53] {\grod \psi} (c)
;
\draw[2cell] 
node[between=b and c at .35, rotate=-90, 2label={above,\grod\Phi}] {\Rightarrow}
;
\end{tikzpicture}
\end{equation}
as follows.  Using the notation in \cref{angajxj,atensorangxj,grodphiangajxj}, for each $n$-tuple of objects
\[\ang{(a_j,x_j)} \in \txprod_{j=1}^n \big(\grod X_j\big)\]
the component morphism of $\grod\Phi$ is defined as the following morphism in $\grod Z$.
\begin{equation}\label{grodPhiajxj}
\begin{tikzpicture}[xscale=5,yscale=1,vcenter]
\draw[0cell=.9]
(0,0) node (x11) {(\grod \phi)\ang{(a_j,x_j)}}
(x11)++(1,0) node (x12) {(\grod \psi)\ang{(a_j,x_j)}}
(x11)++(0,-1) node (x21) {\brb{a,\phi_{\anga}\angx}}
(x12)++(0,-1) node (x22) {\brb{a,\psi_{\anga}\angx}}
;
\draw[1cell=.9]  
(x11) edge node {(\grod \Phi)_{\ang{(a_j,x_j)}}} (x12)
(x21) edge node {\brb{1_a,(\Phi_{\anga})_{\angx}}} (x22)
(x11) edge[equal] (x21)
(x12) edge[equal] (x22)
;
\draw[2cell] 
;
\end{tikzpicture}
\end{equation}
Here
\[\begin{tikzpicture}[xscale=2.4,yscale=2,baseline={(x1.base)}]
\draw[0cell=.9]
(0,0) node (x1) {\textstyle \prod_{j=1}^n X_j a_j}
(x1)++(.17,.03) node (x2) {\phantom{\textstyle \bigotimes Z}}
(x2)++(1,0) node (x3) {Za}
;
\draw[1cell=.9]  
(x2) edge[bend left, shorten <=-.5ex] node[pos=.5] {\phi_{\anga}} (x3)
(x2) edge[bend right, shorten <=-.5ex] node[swap,pos=.5] {\psi_{\anga}} (x3)
;
\draw[2cell]
node[between=x2 and x3 at .4, rotate=-90, 2label={above,\Phi_{\anga}}] {\Rightarrow}
;
\end{tikzpicture}\]
is the $\anga$-component natural transformation of $\Phi$ as in \cref{dcatnarymodcomponent}.  This finishes the definition of $\grod\Phi$
\end{definition}

\begin{lemma}\label{grodPhinatural}
In the context of \cref{def:groconmodification}, $\grod\Phi$ is a natural transformation.
\end{lemma}

\begin{proof}
Using the notation in \cref{angphijfj,angbyp,grodphiangfjpj}, suppose given any $n$-tuple of morphisms
\[\bang{(f_j, p_j) \cn (a_j,x_j) \to (b_j,y_j)} \in \txprod_{j=1}^n \big(\grod X_j\big).\]
The naturality of $\grod\Phi$ means the commutativity of the following diagram in $\grod Z$.
\begin{equation}\label{grodPhinatdiag}
\begin{tikzpicture}[xscale=4.5,yscale=1.4,vcenter]
\def\h{1} \def\v{-1}
\draw[0cell=.9] 
(0,0) node (x11) {\srb{a,\phi_{\anga}\angx}}
(x11)++(1,0) node (x12) {\srb{a,\psi_{\anga}\angx}}
(x11)++(0,-1) node (x21) {\srb{b,\phi_{\angb}\angy}}
(x12)++(0,-1) node (x22) {\srb{b,\psi_{\angb}\angy}}
;
\draw[1cell=.9] 
(x11) edge node {\brb{1_a,(\Phi_{\anga})_{\angx}}} (x12)
(x21) edge node {\brb{1_b,(\Phi_{\angb})_{\angy}}} (x22)
(x11) edge node[swap] {\srb{f,\phi_{\angb}\angp}} (x21)
(x12) edge node {\srb{f,\psi_{\angb}\angp}} (x22)
; 
\end{tikzpicture}
\end{equation}
In the first argument in \cref{grodPhinatdiag}, each composite is equal to $f \cn a \to b$ in $\cD$.  

The second argument in \cref{grodPhinatdiag} forms the boundary in the following diagram in $Zb$.
\begin{equation}\label{grodPhinatzb}
\begin{tikzpicture}[xscale=4.3,yscale=1.2,vcenter]
\def\u{-1.2} \def\v{-.9}
\draw[0cell=.9] 
(0,0) node (x11) {f_* \phi_{\anga}\angx}
(x11)++(1,0) node (x12) {f_* \psi_{\anga}\angx}
(x11)++(0,\v) node (x21) {\phi_{\angb}\ang{(f_j)_*x_j}}
(x12)++(0,\v) node (x22) {\psi_{\angb}\ang{(f_j)_*x_j}}
(x21)++(0,\u) node (x31) {\phi_{\angb}\angy}
(x22)++(0,\u) node (x32) {\psi_{\angb}\angy} 
;
\draw[1cell=.9] 
(x11) edge node {f_*(\Phi_{\anga})_{\angx}} (x12)
(x21) edge node {(\Phi_{\angb})_{\ang{(f_j)_*x_j}}} (x22)
(x31) edge node {(\Phi_{\angb})_{\angy}} (x32)
(x11) edge[equal] (x21)
(x12) edge[equal] (x22)
(x21) edge node[swap] {\phi_{\angb}\angp} (x31)
(x22) edge node {\psi_{\angb}\angp} (x32)
; 
\end{tikzpicture}
\end{equation}
The detail of the diagram \cref{grodPhinatzb} is as follows.
\begin{itemize}
\item $f_* = Zf \cn Za \to Zb$.
\item $(f_j)_* = X_j f_j \cn X_j a_j \to X_j b_j$.
\item The top rectangle commutes by the modification axiom for $\Phi$ \cref{dcatnarymodaxiom}.
\item The bottom rectangle commutes by the naturality of the component natural transformation $\Phi_{\angb}$
\end{itemize}
This proves that $\grod\Phi$ is a natural transformation.
\end{proof}

\begin{lemma}\label{grodPhinlinear}\index{Grothendieck construction!multilinear natural transformation}\index{multilinear natural transformation!Grothendieck construction}
In the context of \cref{def:groconmodification}, $\grod\Phi$ is an $n$-linear natural transformation.
\end{lemma}

\begin{proof}
\cref{grodPhinatural} shows that $\grod\Phi$ is a natural transformation.  It remains to check the $n$-linearity axioms \cref{ntransformationunity,eq:monoidal-in-each-variable}.

\emph{Unity} \cref{ntransformationunity}.  This axiom for $\grod\Phi$ says that, if
\[(a_i,x_i) = (\zero, X_i^0*) \in \grod X_i \forsomespace i \in \{1,\ldots,n\},\]
then 
\[(\grod\Phi)_{\ang{(a_j,x_j)}} = \brb{1_a,(\Phi_{\anga})_{\angx}} = \brb{1_\zero,1_{Z^0*}} \in \grod Z.\]
This axiom holds by
\begin{itemize}
\item the fact that $\zero \in \cD$ is the strict multiplicative zero \cref{ringcataxiommultzero}, which implies 
\[a = \txotimes_{j=1}^n a_j = \zero \in \cD,\] 
and
\item the unity axiom \cref{additivemodunity} for the additive modification $\Phi$.
\end{itemize}

\emph{Constraint Compatibility} \cref{eq:monoidal-in-each-variable}.
Using the notation in \cref{def:groconarityn}, the constraint compatibility axiom for $\grod\Phi$ states the commutativity of the following diagram in $\grod Z$.
\begin{equation}\label{grodPhiconstraintcomp}
\begin{tikzpicture}[xscale=4.2,yscale=1.4,vcenter]
\def\h{1} \def\v{-1}
\draw[0cell=.9] 
(0,0) node (x11) {\left(a \oplus a' \scs \lap_* \phi_{\angapp} \angxpp\right)}
(x11)++(1,0) node (x12) {\brb{a'',\phi_{\angapp} \angxpp}}
(x11)++(0,-1) node (x21) {\left(a \oplus a' \scs \lap_* \psi_{\angapp} \angxpp\right)}
(x12)++(0,-1) node (x22) {\brb{a'',\psi_{\angapp} \angxpp}}
;
\draw[1cell=.85]
(x11) edge node {(\lap^\inv,1)} (x12)
(x21) edge node {(\lap^\inv,1)} (x22)
(x11) edge[transform canvas={xshift=2em}] node[swap] {\brb{1_{a \oplus a'},\lap_* (\Phi_{\angapp})_{\angxpp}}} (x21)
(x12) edge[transform canvas={xshift=-2em}] node {\brb{1_{a''},(\Phi_{\angapp})_{\angxpp}}} (x22)
; 
\end{tikzpicture}
\end{equation}
The detail of the diagram \cref{grodPhiconstraintcomp} is as follows.
\begin{itemize}
\item The top and bottom $(\lap^\inv,1)$ are the $i$-th linearity constraints \cref{grodphilinearity} of, respectively, $\grod\phi$ and $\grod\psi$.  
\item The left vertical arrow comes from the following computation.
\[\left\{\begin{aligned}
&\brb{1_a,(\Phi_{\anga})_{\angx}} \gbox \brb{1_{a'},(\Phi_{\angap})_{\angxp}} &&\\
&= \left(1_a \oplus 1_{a'} \scs Z^2_{a,a'} \left((\Phi_{\anga})_{\angx} \scs (\Phi_{\angap})_{\angxp}\right)\right) & \phantom{=} & \text{by \cref{gboxmorphisms}}\\
&= \left(1_{a \oplus a'} \scs \lap_* (\Phi_{\angapp})_{\angxpp}\right) & \phantom{=} & \text{by \cref{additivitymod}}
\end{aligned}\right.\]
\item In the first argument in $\cD$, each composite is equal to 
\[\lap^\inv \cn a \oplus a' \to a''.\]
\item In the second argument in $Za''$, each composite is equal to the morphism
\[\lap^\inv_* \lap_* (\Phi_{\angapp})_{\angxpp} = (\Phi_{\angapp})_{\angxpp}.\]
\end{itemize}
This proves that \cref{grodPhiconstraintcomp} is commutative.
\end{proof}

\begin{lemma}\label{grodmultimorphismfunctor}\index{Grothendieck construction!multimorphism functor}\index{multimorphism functor!Grothendieck construction}
In the context of \cref{def:groconarityn,def:groconmodification}, the object and morphism assignments,
\[\phi \mapsto \grod\phi \andspace \Phi \mapsto \grod\Phi,\]
define a functor
\begin{equation}\label{grodfunctors}
\begin{tikzcd}[column sep=large,every label/.append style={scale=.85}]
\DCat\scmap{\angX;Z} \ar{r}{\grod} & \permcatsg\scmap{\ang{\grod X};\grod Z}
\end{tikzcd}
\end{equation}
with $\angX = \ang{X_j}_{j=1}^n$ and $\ang{\grod X} = \ang{\grod X_j}_{j=1}^n$.
\end{lemma}

\begin{proof}
The object and morphism assignments are well defined by, respectively, \cref{grodphifunctor,grodPhinlinear}.  It remains to check that $\grod$ preserves identity morphisms and composition.

\emph{Preservation of Identities}.  Consider the identity modification 
\[1_\phi \cn \phi \to \phi\]
of an additive natural transformation $\phi$ in $\DCat\scmap{\angX;Z}$.  Each component natural transformation of $1_\phi$ is an identity.  By definitions \cref{groconidentity,grodPhiajxj} each component morphism of $\grod 1_{\phi}$ is an identity morphism.

\emph{Preservation of Composition}.  Suppose given vertically composable additive modifications $\Phi$ and $\Psi$ in $\DCat\scmap{\angX;Z}$ as follows.
\[\begin{tikzpicture}[xscale=1,yscale=1,baseline={(x1.base)}]
\def\a{50}
\draw[0cell=.9]
(0,0) node (x1) {\bang{(X_j, X_j^2, X_j^0)}_{j=1}^n}
(x1)++(1.1,.03) node (x2) {\phantom{X}}
(x2)++(2,0) node (x3) {\phantom{X}}
(x3)++(.65,0) node (x4) {(Z,Z^2,Z^0)}
;
\draw[1cell=.9]  
(x2) edge[bend left=\a] node[pos=.5] {\phi} (x3)
(x2) edge node[pos=.25] {\psi} (x3)
(x2) edge[bend right=\a] node[swap,pos=.5] {\omega} (x3)
;
\draw[2cell=.9]
node[between=x2 and x3 at .48, shift={(0,.3)}, rotate=-90, 2label={above,\Phi}] {\Rightarrow}
node[between=x2 and x3 at .48, shift={(0,-.3)}, rotate=-90, 2label={above,\Psi}] {\Rightarrow}
;
\end{tikzpicture}\]
Using the notation in \cref{angajxj,atensorangxj,grodphiangajxj}, suppose given an $n$-tuple of objects
\[\ang{(a_j,x_j)} \in \txprod_{j=1}^n \big(\grod X_j\big).\]
By definition \cref{grodPhiajxj} the component morphism of $\grod (\Psi\Phi)$ at $\ang{(a_j,x_j)}$ is the following morphism in $\grod Z$.
\begin{equation}\label{grodPsiPhi}
\begin{tikzpicture}[xscale=1,yscale=1,baseline={(x1.base)}]
\draw[0cell=.9]
(0,0) node (x1) {\brb{a,\phi_{\anga}\angx}}
(x1)++(5,0) node (x2) {\brb{a,\omega_{\anga}\angx}}
;
\draw[1cell=.9]  
(x1) edge node {\brb{1_a,((\Psi\Phi)_{\anga})_{\angx}}} (x2)
;
\end{tikzpicture}
\end{equation}
By the definitions \cref{modificationvcomp,groconcomposition} of, respectively, the vertical composite $\Psi\Phi$ and the composition in $\grod Z$, the morphism \cref{grodPsiPhi} is equal to the composite
\begin{equation}\label{grodPsiPhib}
\begin{tikzpicture}[xscale=.9,yscale=1,baseline={(x1.base)}]
\draw[0cell=.8]
(0,0) node (x1) {\brb{a,\phi_{\anga}\angx}}
(x1)++(4,0) node (x2) {\brb{a,\psi_{\anga}\angx}}
(x2)++(4,0) node (x3) {\brb{a,\omega_{\anga}\angx}.}
;
\draw[1cell=.8]  
(x1) edge node {\brb{1_a,(\Phi_{\anga})_{\angx}}} (x2)
(x2) edge node {\brb{1_a,(\Psi_{\anga})_{\angx}}} (x3)
;
\end{tikzpicture}
\end{equation}
By definition \cref{grodPhiajxj} the composite \cref{grodPsiPhib} is equal to the component morphism at $\ang{(a_j,x_j)}$ of the vertical composite natural transformation
\[\begin{tikzpicture}[xscale=1,yscale=1,baseline={(x1.base)}]
\draw[0cell=.9]
(0,0) node (x1) {\grod\phi}
(x1)++(2.5,0) node (x2) {\grod\psi}
(x2)++(2.5,0) node (x3) {\grod\omega.}
;
\draw[1cell=.9]  
(x1) edge node {\grod\Phi} (x2)
(x2) edge node {\grod\Psi} (x3)
;
\end{tikzpicture}\]
Therefore, $\grod (\Psi\Phi)$ is equal to $(\grod \Psi) (\grod \Phi)$.
\end{proof}

\section{The Grothendieck Construction Preserves Units and Composition}
\label{sec:gropreservesunitcomp}

As in \cref{sec:groadditivenattr,sec:groadditivemod} we continue to assume that $\cD$ is a small tight bipermutative category with
\begin{itemize}
\item additive structure $\Dplus = (\cD,\oplus,\zero,\betaplus)$ and
\item multiplicative structure $\Dte = (\cD,\otimes,\tu,\betate)$ (\cref{def:embipermutativecat}).
\end{itemize}   
In this section we show that the Grothendieck construction $\grod$ preserves the colored units and multicategorical composition.  Recall the following constructions from previous sections.
\begin{itemize}
\item The object assignment $X \mapsto\! \grod X$ is given by \cref{grodxpermutative}.
\item The multimorphism functors in arity 0 are in \cref{def:groconarityzero}.
\item For the multimorphism functors in positive arity, the object and morphism assignments are in, respectively,  \cref{def:groconarityn,def:groconmodification}.
\item The functoriality of the multimorphism functors is proved in \cref{grodmultimorphismfunctor}.
\end{itemize}
Also recall the $\Cat$-multicategories $\DCat$ and $\permcatsg$ in, respectively, \cref{thm:dcatcatmulticat,thm:permcatenrmulticat}.

\begin{lemma}\label{grodpreservesunits}
The Grothendieck construction
\[\begin{tikzcd}[column sep=large,every label/.append style={scale=.85}]
\DCat \ar{r}{\grod} & \permcatsg
\end{tikzcd}\]
preserves the colored units in the sense of \cref{enr-multifunctor-unit}.
\end{lemma}

\begin{proof}
Suppose given an additive symmetric monoidal functor
\[(Z,Z^2,Z^0) \cn \Dplus \to \Cat.\]
Consider \cref{def:groconarityn} for the identity natural transformation $\phi = 1_Z \cn Z \to Z$.
\begin{itemize}
\item By the definitions \cref{grodphiangajxj,grodphiangfjpj} of, respectively, object and morphism assignments, 
\[\grod 1_Z \cn \grod Z \to \grod Z\]
is the identity functor.
\item The only linearity constraint \cref{grodphilinearity} of $\grod 1_Z$ is the identity natural transformation because the unique Laplaza coherence isomorphism
\[\begin{tikzcd}[column sep=large]
a_1 \oplus a_1' \ar{r}{\lap}[swap]{\iso} & a_1 \oplus a_1' 
\end{tikzcd}
\forspace a_1,a_1' \in \cD\]
is the identity morphism by the uniqueness in \cref{thm:laplaza-coherence-1}.
\end{itemize} 
Thus $\grod 1_Z$ is the identity 1-linear functor of $\grod Z$.
\end{proof}

\begin{lemma}\label{grodpreservescomposition}
The Grothendieck construction
\[\begin{tikzcd}[column sep=large,every label/.append style={scale=.85}]
\DCat \ar{r}{\grod} & \permcatsg
\end{tikzcd}\]
preserves composition in the sense of \cref{v-multifunctor-composition}.
\end{lemma}

\begin{proof}
Suppose given additive symmetric monoidal functors
\begin{equation}\label{zxwadditivesmf}
\scalebox{.9}{$(Z,Z^2,Z^0) \scs \bang{(X_j,X_j^2,X_j^0)}_{j=1}^n \scs 
\bang{\bang{(W_{ji},W_{ji}^2,W_{ji}^0)}_{i=1}^{\ell_j}}_{j=1}^n \cn \Dplus \to \Cat$}
\end{equation}
with $n >0$, along with the following notation. 
\begin{equation}\label{xwgrothendieck}
\left\{\begin{aligned}
\angX &= \ang{X_j}_{j=1}^n & \ang{\grod X} &= \ang{\grod X_j}_{j=1}^n\\
\angWj &= \ang{W_{ji}}_{i=1}^{\ell_j} & \ang{\grod W_j} &= \ang{\grod W_{ji}}_{i=1}^{\ell_j}\\
\angW &= \ang{\angWj}_{j=1}^n & \ang{\grod W} &= \ang{\ang{\grod W_j}}_{j=1}^n 
\end{aligned}\right.
\end{equation}
We must show that the following diagram of functors is commutative, with the left and right $\ga$ as in, respectively, \cref{expl:dcatgamma,permcatgamma}.
\begin{equation}\label{grodpreservesgamma}
\begin{tikzpicture}[xscale=5.2,yscale=1.5,vcenter]
\def\h{1} \def\v{-1}
\draw[0cell=.8] 
(0,0) node[align=left] (x11) {$\DCat\scmap{\angX;Z} \,\times$\\ $\txprod_{j=1}^n \DCat\scmap{\angWj;X_j}$}
(x11)++(1,0) node[align=left] (x12) {$\permcatsg\scmap{\ang{\grod X};\grod Z} \,\times$\\ $\txprod_{j=1}^n \permcatsg\scmap{\ang{\grod W_j}; \grod X_j}$}
(x11)++(0,-1) node (x21) {\DCat\scmap{\angW;Z}}
(x12)++(0,-1) node (x22) {\permcatsg\scmap{\ang{\grod W};\grod Z}}
;
\draw[1cell=.85] 
(x11) edge node {\grod \times \smallprod_j\! \grod} (x12)
(x21) edge node {\grod} (x22)
(x11) edge node[swap] {\ga} (x21)
(x12) edge node {\ga} (x22)
; 
\end{tikzpicture}
\end{equation}

\emph{Additive Natural Transformations and Objects}.  To show that \cref{grodpreservesgamma} is commutative, suppose given additive natural transformations
\begin{equation}\label{phixjzphijwjixj}
\begin{tikzcd}
\bigotimes_{j=1}^n X_j \ar{r}{\phi} & Z
\end{tikzcd}\andspace
\begin{tikzcd}
\bigotimes_{i=1}^{\ell_j} W_{ji} \ar{r}{\phi_j} & X_j
\end{tikzcd}
\end{equation}
for $j \in \{1,\ldots,n\}$ and objects
\begin{equation}\label{ajiwjiprodgrodwji}
\ang{(a_{ji},w_{ji})} \in \txprod_{j=1}^n \txprod_{i=1}^{\ell_j} (\grod W_{ji})
\end{equation}
along with the following notation.
\begin{equation}\label{angwji}
\left\{\scalebox{.9}{$\begin{aligned}
a_j &= \txotimes_{i=1}^{\ell_j} a_{ji} \in \cD & a &= \txotimes_{j=1}^n a_j \in \cD\\
\ang{a_j} &= \ang{a_{ji}}_{i=1}^{\ell_j} \in \cD^{\ell_j} & \anga &= \ang{a_j}_{j=1}^n \in \cD^n\\
\ang{w_j} &= \ang{w_{ji}}_{i=1}^{\ell_j} \in \txprod_{i=1}^{\ell_j} W_{ji}a_{ji} & \angw &= \ang{\ang{w_j}}_{j=1}^n \in \txprod_{j=1}^n \txprod_{i=1}^{\ell_j} W_{ji}a_{ji}\\
\end{aligned}$}\right.
\end{equation}
Then the following computation shows that the two composites in the diagram \cref{grodpreservesgamma}, when applied to \cref{phixjzphijwjixj} and the objects \cref{ajiwjiprodgrodwji}, are equal.
\begin{equation}\label{grodgaphiphijaw}
\left\{\begin{aligned}
& \big( \grod \ga\smscmap{\phi;\ang{\phi_j}} \big) \ang{(a_{ji},w_{ji})} &&\\
&= \brb{a, \ga\smscmap{\phi;\ang{\phi_j}}_{\ang{\ang{a_{ji}}_i}_j} \angw} & \phantom{=} & \text{by \cref{grodphiangajxj}}\\
&= \brb{a, \phi_{\anga} \bang{(\phi_j)_{\ang{a_j}} \ang{w_j}}_{j=1}^n} & \phantom{=} & \text{by \cref{dcatganattr}}\\
&= (\grod \phi) \bang{\brb{a_j, (\phi_j)_{\ang{a_j}} \ang{w_j}}}_{j=1}^n & \phantom{=} & \text{by \cref{grodphiangajxj}}\\
&= (\grod \phi) \bang{(\grod \phi_j) \ang{(a_{ji},w_{ji})}_{i=1}^{\ell_j}}_{j=1}^n & \phantom{=} & \text{by \cref{grodphiangajxj}}\\
&= \ga\scmap{\grod \phi; \ang{\grod \phi_j}} \ang{(a_{ji},w_{ji})} & \phantom{=} & \text{by \cref{permcatcomposite}}
\end{aligned}\right.
\end{equation}

\emph{Additive Natural Transformations and Morphisms}.  Suppose given morphisms
\begin{equation}\label{fjipjiprodgrodwji}
\ang{(f_{ji},p_{ji})} \in \txprod_{j=1}^n \txprod_{i=1}^{\ell_j} (\grod W_{ji}).
\end{equation}
The computation \cref{grodgaphiphijaw}, with
\begin{itemize}
\item $\ang{(f_{ji},p_{ji})}$ in place of $(a_{ji},w_{ji})$ and
\item \cref{grodphiangfjpj} in place of \cref{grodphiangajxj},
\end{itemize}  
shows that the two composites in the diagram \cref{grodpreservesgamma}, when applied to \cref{phixjzphijwjixj} and the morphisms \cref{fjipjiprodgrodwji}, are equal.  So far we have shown that the two composites in the diagram \cref{grodpreservesgamma}, when applied to \cref{phixjzphijwjixj}, are equal as functors.  We will show that they also have the same linearity constraints after the following paragraph.

\emph{Additive Modifications}.  Suppose given additive modifications
\[\brb{\Phi,\ang{\Phi_j}} \in \DCat\scmap{\angX;Z} \times \txprod_{j=1}^n \DCat\scmap{\angWj;X_j}.\]
The computation \cref{grodgaphiphijaw}, with
\begin{itemize}
\item $\srb{\Phi,\ang{\Phi_j}}$ in place of $\srb{\phi,\ang{\phi_j}}$,
\item \cref{grodPhiajxj} in place of \cref{grodphiangajxj}, and
\item \cref{dcatgamod} in place of \cref{dcatganattr},
\end{itemize}  
shows that the two composites in the diagram \cref{grodpreservesgamma} are equal on additive modifications.  

\emph{Linearity Constraints}.  To show that the diagram \cref{grodpreservesgamma} is commutative, it remains to show that the functors
\[\grod \ga\smscmap{\phi;\ang{\phi_j}} \andspace \ga\scmap{\grod \phi; \ang{\grod \phi_j}} \cn \txprod_{j=1}^n \txprod_{i=1}^{\ell_j} (\grod W_{ji}) \to \grod Z\]
have the same linearity constraints.  In addition to the objects in \cref{ajiwjiprodgrodwji}, suppose given an object
\[(a_{ji}',w_{ji}') \in \grod W_{ji}\]
for some $j \in \{1,\ldots,n\}$ and some $i \in \{1,\ldots,\ell_j\}$, along with the notation in  \cref{angwji} and below.
\[\left\{\begin{aligned}
a_{ji}'' &= a_{ji} \oplus a_{ji}' & k &= \ell_1 + \cdots + \ell_{j-1} + i\\
a_j' &= a_{j1} \otimes \cdots \otimes a_{ji}' \otimes \cdots \otimes a_{j\ell_j} & a' &= a_1 \otimes \cdots \otimes a_j' \otimes \cdots \otimes a_n\\
a_j'' &= a_{j1} \otimes \cdots \otimes a_{ji}'' \otimes \cdots \otimes a_{j\ell_j} & a'' &= a_1 \otimes \cdots \otimes a_j'' \otimes \cdots \otimes a_n 
\end{aligned}\right.\]
With these objects and notation, by definition \cref{grodphilinearity} the $k$-th linearity constraint of $\grod \ga\smscmap{\phi;\ang{\phi_j}}$ has the form
\[(\lap^\inv,1) \in \grod Z.\]
Here $\lap$ is the unique Laplaza coherence isomorphism\index{Laplaza coherence isomorphism} in $\cD$ (\cref{thm:laplaza-coherence-1})
\begin{equation}\label{lapappaplusap}
\begin{tikzcd}[column sep=large]
a'' \ar{r}{\lap}[swap]{\iso} & a \oplus a'
\end{tikzcd}
\end{equation}
that distributes over the sum $a_{ji} \oplus a_{ji}'$ in $a''$. 

On the other hand, by definition \cref{ffjlinearity} the $k$-th linearity constraint of $\ga\scmap{\grod \phi; \ang{\grod \phi_j}}$ is the following composite in $\grod Z$.
\begin{equation}\label{gagrophiphijconstraint}
\left\{\scalebox{.9}{$\begin{aligned}
&\big((\grod \phi)\bang{1 \compj (\grod\phi_j)_i^2}\big) \circ (\grod\phi)_j^2 &&\\
&= \big((\grod \phi)\bang{1 \compj (\lap^\inv_j,1)} \big) \circ (\lap^\inv_0,1) & \phantom{=} & \text{by \cref{grodphilinearity}}\\
&= \big(1 \otimes \cdots \otimes \lap^\inv_j \otimes \cdots \otimes 1, \phi\ang{1}\big) \circ (\lap^\inv_0,1) & \phantom{=} & \text{by \cref{grodphiangfjpj}}\\
&= \big(1 \otimes \cdots \otimes \lap^\inv_j \otimes \cdots \otimes 1, 1\big) \circ (\lap^\inv_0,1) & \phantom{=} & \text{by functoriality}\\
&= \big((1 \otimes \cdots \otimes \lap^\inv_j \otimes \cdots \otimes 1) \lap^\inv_0, 1 \big) & \phantom{=} & \text{by \cref{groconcomposition}}\\
&= (\lap^\inv,1) & \phantom{=} & \text{by \cref{thm:laplaza-coherence-1}}
\end{aligned}$}\right.
\end{equation}
The detail of \cref{gagrophiphijconstraint} is as follows.
\begin{itemize}
\item The first, second, and fourth equalities follow from, respectively, the definitions \cref{grodphilinearity} , \cref{grodphiangfjpj}, and \cref{groconcomposition} as stated.
\item $\lap_0$ is the unique Laplaza coherence isomorphism in $\cD$
\[\begin{tikzcd}[column sep=large]
a_1 \otimes \cdots \otimes (a_j \oplus a_j') \otimes \cdots \otimes a_n \ar{r}{\lap_0}[swap]{\iso} & a \oplus a'
\end{tikzcd}\]
that distributes over the sum $a_j \oplus a_j'$ in the domain. 
\item $\lap_j$ is the unique Laplaza coherence isomorphism in $\cD$
\[\begin{tikzcd}[column sep=large]
a_j'' \ar{r}{\lap_j}[swap]{\iso} & a_j \oplus a_j'
\end{tikzcd}\]
that distributes over the sum $a_{ji} \oplus a_{ji}'$ in $a_j''$.  Omitting the $\otimes$ symbols to save space, the composite in $\cD$
\[\begin{tikzcd}[column sep=huge]
a'' \ar{r}{1 \cdots \lap_j \cdots 1} & 
a_1 \cdots (a_j \oplus a_j') \cdots a_n \ar{r}{\lap_0} & [-20pt] a \oplus a'
\end{tikzcd}\]
is equal to $\lap$ in \cref{lapappaplusap} by the uniqueness in \cref{thm:laplaza-coherence-1}.  This yields the last equality in \cref{gagrophiphijconstraint}.
\item After the second equality, $\ang{1}$ is an $n$-tuple of identity morphisms.  By the functoriality of each component functor of $\phi$, $\phi\ang{1}$ is also an identity morphism.  This yields the third equality in \cref{gagrophiphijconstraint}.
\end{itemize}
This proves that 
\[\grod \ga\smscmap{\phi;\ang{\phi_j}} \andspace \ga\scmap{\grod \phi; \ang{\grod \phi_j}}\]
have the same linearity constraints.
\end{proof}

\section{Pseudo Symmetry Isomorphisms for the Grothendieck Construction}
\label{sec:pseudosymgrothendieck}

For each small tight bipermutative category $\cD$, in this section we define the pseudo symmetry isomorphisms for the Grothendieck construction $\grod$ and prove that they are well-defined natural isomorphisms.  The four coherence axioms, \cref{pseudosmf-unity,pseudosmf-product,pseudosmf-topeq,pseudosmf-boteq}, of a pseudo symmetric $\Cat$-multifunctor are proved in \cref{grodpseudosymaxioms}.  Here is an outline of this section.
\begin{itemize}
\item The pseudo symmetry isomorphisms of $\grod$ are in \cref{def:gropseudosymmetry}.
\item \cref{grodsigmaphilinear,grodsigmaphinlinear} show that each component of each pseudo symmetry isomorphism is an $n$-linear natural isomorphism.
\item \cref{grodsigmanaturaliso} establishes the naturality of the pseudo symmetry isomorphisms.
\end{itemize}  

First we define the natural isomorphisms that will be shown to be a part of the pseudo symmetric $\Cat$-multifunctor structure on the Grothendieck construction.

\begin{definition}\label{def:gropseudosymmetry}
Suppose given additive symmetric monoidal functors $Z$ and $\angX$ as in \cref{zxwadditivesmf,xwgrothendieck} and a permutation $\sigma \in \Sigma_n$.  Define the data of a natural isomorphism
\begin{equation}\label{gropseudosymmetry}
\begin{tikzpicture}[xscale=4,yscale=1.3,vcenter]
\def\h{1} \def\v{-1}
\draw[0cell=.8] 
(0,0) node (x11) {\DCat\scmap{\angX;Z}}
(x11)++(1,0) node (x12) {\permcatsg\scmap{\ang{\grod X};\grod Z}}
(x11)++(0,-1) node (x21) {\DCat\scmap{\angX\sigma;Z}}
(x12)++(0,-1) node (x22) {\permcatsg\scmap{\ang{\grod X}\sigma;\grod Z}}
;
\draw[1cell=.9] 
(x11) edge node (f1) {\grod} (x12)
(x21) edge node[swap] (f2){\grod} (x22)
(x11) edge node[swap] {\sigma} (x21)
(x12) edge node {\sigma} (x22)
;
\draw[2cell=.9]
node[between=f1 and f2 at .63, 2label={above,(\grod)_{\sigma,\angX,Z}}] {\Longrightarrow}
;
\end{tikzpicture}
\end{equation}
called the \index{Grothendieck construction!pseudo symmetry isomorphism}\index{pseudo symmetry isomorphism!Grothendieck construction}\emph{pseudo symmetry isomorphism}, as follows.  Suppose given an additive natural transformation \cref{additivenattransformation}
\[\begin{tikzcd}[column sep=large]
\bigotimes_{j=1}^n X_j \ar{r}{\phi} & Z
\end{tikzcd}\]
with
\begin{itemize}
\item image under the $\sigma$-action \cref{phisigmaanga}
\[\begin{tikzcd}[column sep=large]
\bigotimes_{j=1}^n X_{\sigma(j)} \ar{r}{\phi^\sigma} & Z
\end{tikzcd}\]
and
\item associated $n$-linear functor \cref{groadditivenattr} 
\[\begin{tikzcd}[column sep=huge, every label/.append style={scale=.8}]
\prod_{j=1}^n (\grod X_j) \ar{r}{\grod\phi} & \grod Z.
\end{tikzcd}\]
\end{itemize}
Suppose given an $n$-tuple of objects
\begin{equation}\label{ajxjgroxsigmaj}
\ang{(a_j,x_j)} \in \txprod_{j=1}^n (\grod X_{\sigma(j)})
\end{equation}
with each $a_j \in \cD$ and $x_j \in X_{\sigma(j)} a_j$.  To define the component of $(\grod)_{\sigma,\angX,Z}$ at $\phi$ and $\ang{(a_j,x_j)}$, first note that its domain and codomain objects in $\grod Z$ are as follows.
\begin{equation}\label{grodsigmadomain}
\left\{\begin{aligned}
& (\grod \phi^\sigma) \ang{(a_j,x_j)} &&\\
&\phantom{==}= \brb{a, \phi^\sigma_{\anga}\angx} & \phantom{=} & \text{by \cref{grodphiangajxj}}\\
&\phantom{==}= \brb{a, \sigmainv_* \phi_{\sigma\anga} \sigma\angx} & \phantom{=} & \text{by \cref{phisigmaanga}}\\
& (\grod\phi)^\sigma \ang{(a_j,x_j)} &&\\
&\phantom{==}= (\grod\phi) \bang{\brb{a_{\sigmainv(j)}, x_{\sigmainv(j)}}} & \phantom{=} & \text{by \cref{permcatsigmaaction}}\\
&\phantom{==}= \brb{a_{\sigmainv}, \phi_{\sigma\anga} \sigma\angx} & \phantom{=} & \text{by \cref{sigmaangaasigmainv,grodphiangajxj}}
\end{aligned}\right.
\end{equation}
Using \cref{grodsigmadomain} we define the component of $(\grod)_{\sigma,\angX,Z}$ at $\phi$ and $\ang{(a_j,x_j)}$ as the isomorphism in $\grod Z$
\begin{equation}\label{grodsigmacomponent}
\begin{tikzcd}[column sep=large]
\brb{a, \sigmainv_* \phi_{\sigma\anga} \sigma\angx} \ar{r}{(\sigma,1)}[swap]{\iso} &
\brb{a_{\sigmainv}, \phi_{\sigma\anga} \sigma\angx}.
\end{tikzcd}
\end{equation}
This finishes the definition of $(\grod)_{\sigma,\angX,Z}$.

As in \cref{sec:pseudosmultifunctor} we usually abbreviate $(\grod)_{\sigma,\angX,Z}$ to $(\grod)_{\sigma}$.  Its components at $\phi$ and $\brb{\phi,\ang{(a_j,x_j)}}$ are denoted by\label{not:grodsigmaphi}
\[(\grod)_{\sigma,\phi} \andspace (\grod)_{\sigma,\phi,\ang{(a_j,x_j)}} \scs\]
respectively.
\end{definition}

\begin{explanation}\label{expl:grodsigma}
The two arguments in the isomorphism
\[(\grod)_{\sigma,\phi,\ang{(a_j,x_j)}} = (\sigma,1)\]
in \cref{grodsigmacomponent} are as follows.
\begin{itemize}
\item $\sigma$ is the unique coherence isomorphism in $\Dte = (\cD,\otimes,\tu,\betate)$
\begin{equation}\label{atoasigmainv}
\begin{tikzcd}[column sep=large]
a = \bigotimes_{j=1}^n a_j \ar{r}{\sigma}[swap]{\iso} & \bigotimes_{j=1}^n a_{\sigmainv(j)} = a_{\sigmainv} 
\end{tikzcd}
\end{equation}
that permutes the $n$ factors according to $\sigma \in \Sigma_n$.  This is the inverse of the coherence isomorphism in \cref{asigmainvtoa}. 
\item The second argument $1$ is the identity morphism
\[\begin{tikzcd}[column sep=large]
\sigma_* \sigmainv_* \phi_{\sigma\anga} \sigma\angx \ar{r}{=} & \phi_{\sigma\anga} \sigma\angx
\end{tikzcd}\]
in $Za_{\sigmainv}$.\defmark
\end{itemize} 
\end{explanation}

\subsection*{Componentwise Multilinearity}

Recall that the $n$-ary 2-cells in the $\Cat$-multicategory $\permcatsg$ are $n$-linear natural transformations between $n$-linear functors (\cref{def:nlineartransformation}).  An $n$-linear natural transformation is a natural transformation that satisfies the unity axiom \cref{ntransformationunity} and the constraint compatibility axiom \cref{eq:monoidal-in-each-variable}.  Next we show that each component of the pseudo symmetry isomorphism $(\grod)_{\sigma}$ is an $n$-linear natural isomorphism, starting with its naturality.

\begin{lemma}\label{grodsigmaphilinear}
In the context of \cref{def:gropseudosymmetry}, the component of $(\grod)_{\sigma}$ at $\phi$
\begin{equation}\label{grodsigmaphicomponent}
\begin{tikzpicture}[xscale=2.7,yscale=2.2,baseline={(a.base)}]
\draw[0cell=.9]
(0,0) node (a) {\txprod_{j=1}^n (\grod X_{\sigma(j)})}
(a)++(.32,.01) node (b) {\phantom{D}}
(b)++(1,0) node (c) {\grod Z}
;
\draw[1cell=.8]  
(b) edge[bend left] node[pos=.53] {\grod \phi^\sigma} (c)
(b) edge[bend right] node[swap,pos=.53] {(\grod\phi)^\sigma} (c)
;
\draw[2cell=.9] 
node[between=b and c at .3, rotate=-90, 2label={above,(\grod)_{\sigma,\phi}}] {\Rightarrow}
;
\end{tikzpicture}
\end{equation}
is a natural isomorphism. 
\end{lemma}

\begin{proof}  
By definition \cref{grodsigmacomponent} each component of $(\grod)_{\sigma,\phi}$ is an isomorphism in $\grod Z$.

To show that $(\grod)_{\sigma,\phi}$ is a natural transformation, suppose given an $n$-tuple of morphisms
\[\bang{(f_j,p_j) \cn (a_j,x_j) \to (b_j,y_j)} \in \txprod_{j=1}^n (\grod X_{\sigma(j)})\]
with 
\begin{itemize}
\item each $f_j \cn a_j \to b_j$ a morphism in $\cD$ and
\item each $p_j \cn (f_j)_* x_j \to y_j$ a morphism in $X_{\sigma(j)} b_j$, where 
\[(f_j)_* = X_{\sigma(j)} f_j \cn X_{\sigma(j)} a_j \to X_{\sigma(j)} b_j.\] 
\end{itemize}
Adapting the notation for $\ang{(a_j,x_j)}$ to $\ang{(b_j,y_j)}$ and $\ang{(f_j,p_j)}$ as in \cref{angbyp}, we must show that the following naturality diagram for $(\grod)_{\sigma,\phi}$ is commutative in $\grod Z$.  
\begin{equation}\label{grodsigmaphinatdiagram}
\begin{tikzpicture}[xscale=4,yscale=1.4,vcenter]
\draw[0cell=.9] 
(0,0) node (x11) {\brb{a, \sigmainv_* \phi_{\sigma\anga} \sigma\angx}}
(x11)++(1,0) node (x12) {\brb{a_{\sigmainv}, \phi_{\sigma\anga} \sigma\angx}}
(x11)++(0,-1) node (x21) {\brb{b, \sigmainv_* \phi_{\sigma\angb} \sigma\angy}}
(x12)++(0,-1) node (x22) {\brb{b_{\sigmainv}, \phi_{\sigma\angb} \sigma\angy}}
;
\draw[1cell=.9] 
(x11) edge node {(\sigma,1)} (x12)
(x21) edge node {(\sigma,1)} (x22)
(x11) edge[transform canvas={xshift=2em}] node[swap] {\brb{f, \sigmainv_* \phi_{\sigma\angb} \sigma\angp}} (x21)
(x12) edge[transform canvas={xshift=-2em}] node {\brb{f_{\sigmainv}, \phi_{\sigma\angb} \sigma\angp}} (x22)
;
\end{tikzpicture}
\end{equation}
The arrows in the diagram \cref{grodsigmaphinatdiagram} are as follows.
\begin{itemize}
\item The top and bottom horizontal arrows $(\sigma,1)$ are the components of $(\grod)_{\sigma,\phi}$ in \cref{grodsigmacomponent} at, respectively, $\ang{(a_j,x_j)}$ and $\ang{(b_j,y_j)}$.
\item The left and right vertical arrows are, respectively, 
\[(\grod\phi^\sigma)\ang{(f_j,p_j)} \andspace (\grod\phi)^\sigma\ang{(f_j,p_j)}.\]
\end{itemize}
The first argument in the diagram \cref{grodsigmaphinatdiagram} is the following diagram in $\Dte = (\cD, \otimes, \betate)$.  It commutes by the naturality of the multiplicative braiding $\betate$.
\[\begin{tikzpicture}[xscale=2,yscale=1.2,vcenter]
\draw[0cell] 
(0,0) node (x11) {a}
(x11)++(1,0) node (x12) {a_{\sigmainv}}
(x11)++(0,-1) node (x21) {b}
(x12)++(0,-1) node (x22) {b_{\sigmainv}}
;
\draw[1cell] 
(x11) edge node {\sigma} (x12)
(x21) edge node {\sigma} (x22)
(x11) edge node[swap] {f} (x21)
(x12) edge node {f_{\sigmainv}} (x22)
;
\end{tikzpicture}\]
The second argument in the diagram \cref{grodsigmaphinatdiagram} yields two composites in $Zb_{\sigmainv}$.  They are equal by the following computation.
\[\begin{split}
1 \circ \sigma_* \sigmainv_* \phi_{\sigma\angb} \sigma\angp 
&= \phi_{\sigma\angb} \sigma\angp\\
&= \phi_{\sigma\angb} \sigma\angp \circ 1\\
&= \phi_{\sigma\angb} \sigma\angp \circ (f_{\sigmainv})_* 1
\end{split}\]
This shows that the diagram \cref{grodsigmaphinatdiagram} is commutative.
\end{proof}

Next we check the multilinearity of $(\grod)_{\sigma,\phi}$.

\begin{lemma}\label{grodsigmaphinlinear}
In the context of \cref{def:gropseudosymmetry}, the component $(\grod)_{\sigma,\phi}$ is an $n$-linear natural isomorphism.
\end{lemma}

\begin{proof}
\cref{grodsigmaphilinear} shows that $(\grod)_{\sigma,\phi}$ is a natural isomorphism.  It remains to check the two axioms in \cref{def:nlineartransformation}.

\emph{Unity} \cref{ntransformationunity}.
The unity axiom for $(\grod)_{\sigma,\phi}$ states that, if for some $i \in \{1,\ldots,n\}$
\[(a_i,x_i) = \brb{\zero,X_{\sigma(i)}^0 \!*} \in \grod X_{\sigma(i)}\]
is the monoidal unit \cref{groconunit} in $\grod X_{\sigma(i)}$, then the component \cref{grodsigmacomponent}
\[(\grod)_{\sigma,\phi,\ang{(a_j,x_j)}} = (\sigma,1)\]
is the identity morphism.  This axiom holds because, if $a_i = \zero$, which is the multiplicative zero in $\cD$, then the uniqueness in \cref{thm:laplaza-coherence-1} implies that the coherence isomorphism in $\cD$ \cref{atoasigmainv}
\[\begin{tikzcd}[column sep=large]
a = \zero \ar{r}{\sigma} & \zero = a_{\sigmainv}
\end{tikzcd}\]
is the identity morphism $1_\zero$.  Therefore, 
\[(\sigma,1) = (1_\zero,1) \cn (\zero,Z^0*) \to (\zero,Z^0*)\]
is the identity morphism of the monoidal unit in $\grod Z$ by the unity axiom \cref{additivenattrunity} for $\phi$.

\emph{Constraint Compatibility} \cref{eq:monoidal-in-each-variable}. 
In addition to the objects $\ang{(a_j,x_j)}$ in \cref{ajxjgroxsigmaj}, suppose given an object
\[(a_i',x_i') \in \grod X_{\sigma(i)}\]
for some $i \in \{1,\ldots,n\}$ and objects $a_i' \in \cD$ and $x_i' \in X_{\sigma(i)} a_i'$.  We use the notation in 
\begin{itemize}
\item \cref{sigmaangaprime} for $a_{\sigmainv}'$ and $a_{\sigmainv}''$ and
\item \cref{atensorangxj,grodphiconstraintnot} for $a$, $a'$, $a''$, $\anga$, $\angap$, $\angapp$, $\angx$, $\angxp$, and $\angxpp$ with $X_{\sigma(j)}$ in place of $X_j$ for each $j \in \{1,\ldots,n\}$.
\end{itemize}
With these objects and notation, the constraint compatibility axiom for $(\grod)_{\sigma,\phi}$ states that the following diagram in $\grod Z$ is commutative.
\begin{equation}\label{grodsigmaphiconstraint}
\begin{tikzpicture}[xscale=5.5,yscale=1.4,vcenter]
\def\s{.8}
\draw[0cell=\s] 
(0,0) node (x11) {\brb{a, \sigmainv_* \phi_{\sigma\anga} \sigma\angx} \gbox \brb{a', \sigmainv_* \phi_{\sigma\angap} \sigma\angxp}}
(x11)++(1,0) node (x12) {\brb{a'', \sigmainv_* \phi_{\sigma\angapp} \sigma\angxpp}}
(x11)++(0,-1) node (x21) {\brb{a_{\sigmainv}, \phi_{\sigma\anga} \sigma\angx} \gbox \brb{a'_{\sigmainv}, \phi_{\sigma\angap} \sigma\angxp}}
(x12)++(0,-1) node (x22) {\brb{a''_{\sigmainv}, \phi_{\sigma\angapp} \sigma\angxpp}}
;
\draw[1cell=\s] 
(x11) edge node {(\grod \phi^\sigma)^2_i} (x12)
(x21) edge node {\left((\grod\phi)^\sigma\right)^2_i} (x22)
(x11) edge[transform canvas={xshift=-.5ex}] node[swap] {(\sigma,1) \gbox (\sigma,1)} (x21)
(x12) edge node {(\sigma,1)} (x22)
;
\end{tikzpicture}
\end{equation}
The detail of the diagram \cref{grodsigmaphiconstraint} is as follows.
\begin{itemize}
\item By 
\begin{itemize}
\item the definitions \cref{gboxobjects,gboxmorphisms} of the monoidal product $\gbox$ in $\grod Z$ and
\item the functoriality of each component functor of the monoidal constraint $Z^2$ \cref{additivesmfconstraint},
\end{itemize}  
the left vertical morphism in \cref{grodsigmaphiconstraint}, $(\sigma,1) \gbox (\sigma,1)$, is the following morphism.
\begin{equation}\label{sigmaplussigmaone}
\begin{tikzpicture}[xscale=5.5,yscale=1.4,vcenter]
\def\s{.9}
\draw[0cell=\s] 
(0,0) node (x11) {\Big(a \oplus a' \scs Z^2_{a,a'}\brb{\sigmainv_* \phi_{\sigma\anga} \sigma\angx, \sigmainv_* \phi_{\sigma\angap} \sigma\angxp} \Big)}
(x11)++(0,-1) node (x21) {\Big(a_{\sigmainv} \oplus a'_{\sigmainv} \scs Z^2_{a_{\sigmainv}, a'_{\sigmainv}} \brb{\phi_{\sigma\anga} \sigma\angx, \phi_{\sigma\angap} \sigma\angxp} \Big)}
;
\draw[1cell=\s] 
(x11) edge[shorten >=-1ex, transform canvas={xshift=-1em}] node[swap,pos=.6] {(\sigma \oplus \sigma,1)} (x21)
;
\end{tikzpicture}
\end{equation}
\item By definition \cref{grodphilinearity}, the top horizontal morphism in \cref{grodsigmaphiconstraint} is
\begin{equation}\label{grodphisigmatwoi}
(\grod \phi^\sigma)^2_i = (\lap^\inv,1).
\end{equation}
Here $\lap$ is the unique Laplaza coherence isomorphism in $\cD$
\[\begin{tikzcd}[column sep=large]
a'' \ar{r}{\lap}[swap]{\iso} & a \oplus a'
\end{tikzcd}\]
that distributes over the sum $a_i \oplus a_i'$ in the $i$-th tensor factor in $a''$. 
\item By \cref{fsigmatwoj,grodphilinearity}, the bottom horizontal morphism in \cref{grodsigmaphiconstraint} is
\begin{equation}\label{sigmaofgrodphitwoi}
\left((\grod\phi)^\sigma\right)^2_i = \big(\lap_{\sigmainv}^\inv \scs 1\big).
\end{equation}
Here $\lap_{\sigmainv}$ is the unique Laplaza coherence isomorphism in $\cD$
\[\begin{tikzcd}[column sep=large]
a''_{\sigmainv} \ar{r}{\lap_{\sigmainv}}[swap]{\iso} & a_{\sigmainv} \oplus a'_{\sigmainv}
\end{tikzcd}\]
that distributes over the sum $a_i \oplus a_i'$ in the $\sigma(i)$-th tensor factor in $a''_{\sigmainv}$. 
\end{itemize}

Using \cref{sigmaplussigmaone,grodphisigmatwoi,sigmaofgrodphitwoi} the first argument in the diagram \cref{grodsigmaphiconstraint} is the following diagram in $\cD$.
\[\begin{tikzpicture}[xscale=3,yscale=1.2,vcenter]
\draw[0cell] 
(0,0) node (x11) {a \oplus a'}
(x11)++(1,0) node (x12) {a''} 
(x11)++(0,-1) node (x21) {a_{\sigmainv} \oplus a'_{\sigmainv}}
(x12)++(0,-1) node (x22) {a''_{\sigmainv}}
;
\draw[1cell=.9] 
(x11) edge node {\lapinv} (x12)
(x21) edge node {\lapinv_{\sigmainv}} (x22)
(x11) edge node[swap] {\sigma \oplus \sigma} (x21)
(x12) edge node {\sigma} (x22)
;
\end{tikzpicture}\]
This diagram commutes by the uniqueness in \cref{thm:laplaza-coherence-1}.  

Next, each of the two composites in the diagram \cref{grodsigmaphiconstraint} has the second argument given by the identity morphism of $\phi_{\sigma\angapp}\sigma\angxpp$.  This is true because, in each of the four morphisms in the diagram \cref{grodsigmaphiconstraint}, the second argument is an identity morphism.  This shows that $(\grod)_{\sigma,\phi}$ satisfies the constraint compatibility axiom.
\end{proof}

\subsection*{Naturality of the Pseudo Symmetry Isomorphisms}

Next we show that the pseudo symmetry isomorphisms for the Grothendieck construction are natural with respect to the morphisms in the multimorphism category $\DCat\scmap{\angX;Z}$.

\begin{lemma}\label{grodsigmanaturaliso}
In the context of \cref{def:gropseudosymmetry}, the pseudo symmetry isomorphism
\[\begin{tikzpicture}[xscale=1,yscale=1,vcenter]
\def\s{.9} \def\a{30}
\draw[0cell=\s] 
(0,0) node (x1) {\DCat\scmap{\angX;Z}}
(x1)++(.9,0) node (x2) {\phantom{X}}
(x2)++(2.5,0) node (x3) {\phantom{X}}
(x3)++(1.7,0) node (x4) {\permcatsg\scmap{\ang{\grod X}\sigma;\grod Z}} 
;
\draw[1cell=\s] 
(x2) edge[bend left=\a] node[pos=.5] {\grod (-)^\sigma} (x3)
(x2) edge[bend right=\a] node[swap,pos=.5] {(\grod (-))^\sigma} (x3)
;
\draw[2cell]
node[between=x2 and x3 at .35, rotate=-90, 2label={above,(\grod)_\sigma}] {\Rightarrow}
;
\end{tikzpicture}\]
is a natural isomorphism.
\end{lemma}

\begin{proof}
\cref{grodsigmaphinlinear} shows that each component of $(\grod)_\sigma$ is an $n$-linear natural isomorphism.  To show its naturality, suppose given an additive modification $\Phi$ as in \cref{additivemodtwocell}
\[\begin{tikzpicture}[xscale=1,yscale=1,baseline={(x1.base)}]
\draw[0cell=.9]
(0,0) node (x1) {\bang{(X_j, X_j^2, X_j^0)}_{j=1}^n}
(x1)++(1,0) node (x2) {\phantom{X}}
(x2)++(2,0) node (x3) {\phantom{X}}
(x3)++(.65,0) node (x4) {(Z,Z^2,Z^0).}
;
\draw[1cell=.9]  
(x2) edge[bend left] node[pos=.5] {\phi} (x3)
(x2) edge[bend right] node[swap,pos=.5] {\psi} (x3)
;
\draw[2cell]
node[between=x2 and x3 at .45, rotate=-90, 2label={above,\Phi}] {\Rightarrow}
;
\end{tikzpicture}\]
Naturality with respect to $\Phi$ means the commutativity of the following diagram of $n$-linear natural transformations in $\permcatsg\scmap{\ang{\grod X}\sigma;\grod Z}$.
\begin{equation}\label{grodsigmanatPhi}
\begin{tikzpicture}[xscale=3,yscale=1.3,vcenter]
\def\s{.9}
\draw[0cell=\s] 
(0,0) node (x11) {\grod \phi^\sigma}
(x11)++(1,0) node (x12) {(\grod\phi)^\sigma}
(x11)++(0,-1) node (x21) {\grod \psi^\sigma}
(x12)++(0,-1) node (x22) {(\grod\psi)^\sigma}
;
\draw[1cell=\s] 
(x11) edge node {(\grod)_{\sigma,\phi}} (x12)
(x21) edge node {(\grod)_{\sigma,\psi}} (x22)
(x11) edge node[swap] {\grod \Phi^\sigma} (x21)
(x12) edge node {(\grod \Phi)^\sigma} (x22)
;
\end{tikzpicture}
\end{equation}
It suffices to check the commutativity of \cref{grodsigmanatPhi} when it is evaluated at a general $n$-tuple of objects as in \cref{ajxjgroxsigmaj}
\[\ang{(a_j,x_j)} \in \txprod_{j=1}^n (\grod X_{\sigma(j)}).\]
By definitions \cref{Phisigmaanga,permcatsigmaaction,grodPhiajxj,grodsigmacomponent}, the diagram \cref{grodsigmanatPhi} evaluated at $\ang{(a_j,x_j)}$ is the following diagram in $\grod Z$.
\begin{equation}\label{grodsigmanatPhiobj}
\begin{tikzpicture}[xscale=4,yscale=1.4,vcenter]
\draw[0cell=.9] 
(0,0) node (x11) {\brb{a, \sigmainv_* \phi_{\sigma\anga} \sigma\angx}}
(x11)++(1,0) node (x12) {\brb{a_{\sigmainv}, \phi_{\sigma\anga} \sigma\angx}}
(x11)++(0,-1) node (x21) {\brb{a, \sigmainv_* \psi_{\sigma\anga} \sigma\angx}}
(x12)++(0,-1) node (x22) {\brb{a_{\sigmainv}, \psi_{\sigma\anga} \sigma\angx}}
;
\draw[1cell=.9] 
(x11) edge node {(\sigma,1)} (x12)
(x21) edge node {(\sigma,1)} (x22)
(x11) edge[transform canvas={xshift=1.5em}] node[swap] {\brb{1_a, \sigmainv_* (\Phi_{\sigma\anga})_{\sigma\angx}}} (x21)
(x12) edge[transform canvas={xshift=-1.5em}] node {\brb{1_{a_{\sigmainv}}, (\Phi_{\sigma\anga})_{\sigma\angx}}} (x22)
;
\end{tikzpicture}
\end{equation}
The detail of the diagram \cref{grodsigmanatPhiobj} is as follows.
\begin{itemize}
\item The first argument of each of the two composites is the coherence isomorphism
\[\begin{tikzcd}[column sep=large]
a \ar{r}{\sigma}[swap]{\iso} & a_{\sigmainv}
\end{tikzcd} \inspace \Dte.\]
\item The second argument in the diagram \cref{grodsigmanatPhiobj} yields two composites in $Za_{\sigmainv}$.  They are equal by the following computation.
\[\begin{split}
1 \circ \sigma_* \sigmainv_* (\Phi_{\sigma\anga})_{\sigma\angx} 
&= (\Phi_{\sigma\anga})_{\sigma\angx}\\
&= (\Phi_{\sigma\anga})_{\sigma\angx} \circ 1\\
&= (\Phi_{\sigma\anga})_{\sigma\angx} \circ (1_{a_{\sigmainv}})_* 1
\end{split}\]
\end{itemize}
This shows that the diagram \cref{grodsigmanatPhiobj} is commutative.
\end{proof}

\section{Pseudo Symmetric \texorpdfstring{$\Cat$}{Cat}-Multifunctoriality}
\label{sec:enrmultigrothendieck}

As in \cref{sec:groadditivenattr,sec:groadditivemod,sec:gropreservesunitcomp,sec:pseudosymgrothendieck} we continue to assume that $\cD$ is a small tight bipermutative category with
\begin{itemize}
\item additive structure $\Dplus = (\cD,\oplus,\zero,\betaplus)$ and
\item multiplicative structure $\Dte = (\cD,\otimes,\tu,\betate)$ (\cref{def:embipermutativecat}).
\end{itemize}
In \cref{grodpreservesunits,grodpreservescomposition} we showed that the Grothendieck construction $\grod$ strictly preserves the colored units and multicategorical composition.  In \cref{sec:pseudosymgrothendieck} we constructed the pseudo symmetry isomorphisms for the Grothendieck construction and proved that they are natural isomorphisms.  To show that $\grod$ is a pseudo symmetric $\Cat$-multifunctor, next we check that its pseudo symmetry isomorphisms satisfy the four axioms in \cref{def:pseudosmultifunctor}.

\begin{lemma}\label{grodpseudosymaxioms}
The pseudo symmetry isomorphisms $(\grod)_\sigma$ in \cref{gropseudosymmetry} satisfy the axioms \cref{pseudosmf-unity,pseudosmf-product,pseudosmf-topeq,pseudosmf-boteq}.
\end{lemma}

\begin{proof}
Suppose given additive symmetric monoidal functors 
\[Z \scs \ang{X_j}_{j=1}^n \scs \bang{\ang{W_{ji}}_{i=1}^{\ell_j}}_{j=1}^n \cn \Dplus \to \Cat\] 
as in \cref{zxwadditivesmf,xwgrothendieck}.

\medskip
\emph{The Unit Permutation Axiom} \cref{pseudosmf-unity}.

Suppose $\id_n \in \Sigma_n$ is the identity permutation.  By definition \cref{grodsigmacomponent} there is an equality of morphisms in $\grod Z$
\[(\grod)_{\id_n,\phi,\ang{(a_j,x_j)}} = (1_a,1) 
\cn \brb{a,\phi_{\anga}\angx} \to \brb{a,\phi_{\anga}\angx}\]
for 
\begin{itemize}
\item each additive natural transformation $\phi \in \DCat\scmap{\angX;Z}$ as in \cref{additivenattransformation} and
\item each $n$-tuple of objects as in \cref{ajxjgroxsigmaj}
\[\ang{(a_j,x_j)} \in \txprod_{j=1}^n (\grod X_j).\]
\end{itemize}
Since $(1_a,1)$ is an identity morphism \cref{groconidentity}, this shows that $(\grod)_{\id_n}$ is the identity natural transformation of $\grod$.

\medskip
\emph{The Product Permutation Axiom} \cref{pseudosmf-product}.

Suppose given
\begin{itemize}
\item permutations $\sigma,\tau \in \Sigma_n$,
\item an additive natural transformation $\phi \in \DCat\scmap{\angX;Z}$, and
\item an $n$-tuple of objects
\[\ang{(a_j,x_j)} \in \txprod_{j=1}^n (\grod X_{\sigma\tau(j)})\]
with each $a_j \in \cD$ and $x_j \in X_{\sigma\tau(j)} a_j$.
\end{itemize}
By definition \cref{permcatsigmaaction,grodsigmacomponent}, for the above data the left pasting diagram in the product permutation axiom \cref{pseudosmf-product} is the composite in $\grod Z$
\[(\sigma,1) \circ (\tau,1) = (\sigma\tau,1).\]
The right-hand side above is the component of $(\grod)_{\sigma\tau}$ at $\phi$ and $\ang{(a_j,x_j)}$.  This proves the product permutation axiom.

\medskip
\emph{The Top Equivariance Axiom} \cref{pseudosmf-topeq}.

Suppose given
\begin{itemize}
\item a permutation $\sigma \in \Sigma_n$,
\item additive natural transformations
\begin{equation}\label{phiphijtensorxjz}
\begin{tikzcd}
\txotimes_{j=1}^n X_j \ar{r}{\phi} & Z
\end{tikzcd}\andspace
\begin{tikzcd}
\txotimes_{i=1}^{\ell_j} W_{ji} \ar{r}{\phi_j} & X_j
\end{tikzcd}
\end{equation}
for $j \in \{1,\ldots,n\}$ and $i \in \{1,\ldots,\ell_j\}$, and
\item an $(\ell_{\sigma(1)} + \cdots + \ell_{\sigma(n)})$-tuple of objects
\[\ang{a,w} = \ang{\ang{(a_{ji}, w_{ji})}_i}_j \in \txprod_{j=1}^n \txprod_{i=1}^{\ell_{\sigma(j)}} (\grod W_{\sigma(j),i})\]
with each $a_{ji} \in \cD$ and $w_{ji} \in W_{\sigma(j),i}\, a_{ji}$.
\end{itemize}
We use the notation below, where $\sigmabar$ is the block permutation induced by $\sigma$ that permutes $n$ blocks of lengths $\{\ell_{\sigma(j)}\}_{j=1}^n$.
\begin{equation}\label{phijgphij}
\left\{
\scalebox{.9}{$\begin{aligned}
\ang{\phi_j} &= (\phi_1,\ldots,\phi_n) & g &= \gamma\scmap{\phi;\ang{\phi_j}} \qquad\quad \phi^{\sigma(j)} = \grod \phi_{\sigma(j)}\\
\ell &= \ell_{\sigma(1)} + \cdots + \ell_{\sigma(n)} & \sigmabar &= \sigma\ang{\ell_{\sigma(1)}, \ldots, \ell_{\sigma(n)}} \in \Sigma_\ell\\
\ang{a_j} &= \ang{a_{ji}}_{i=1}^{\ell_{\sigma(j)}} \in \cD^{\ell_{\sigma(j)}} & \ang{w_j} &= \ang{w_{ji}}_{i=1}^{\ell_{\sigma(j)}} \in \txprod_{i=1}^{\ell_{\sigma(j)}} W_{\sigma(j),i}\, a_{ji} \\
\anga &= \ang{\ang{a_j}}_{j=1}^n \in \cD^\ell & \angw &= \ang{\ang{w_j}}_{j=1}^n \in \txprod_{j=1}^n \txprod_{i=1}^{\ell_{\sigma(j)}} W_{\sigma(j),i}\, a_{ji}\\
a_j &= \txotimes_{i=1}^{\ell_{\sigma(j)}} a_{ji} \in \cD & a &= \txotimes_{j=1}^n a_j \in \cD \qquad\quad a_{\sigmabarinv} = \txotimes_{j=1}^n a_{\sigmainv(j)} \in \cD 
\end{aligned}$}\right.
\end{equation}
By definition \cref{grodsigmacomponent}, for the above data the left pasting diagram in \cref{pseudosmf-topeq} yields the following morphism in $\grod Z$.
\begin{equation}\label{grodsigmabargaw}
\begin{tikzpicture}[xscale=5,yscale=1.2,vcenter]
\draw[0cell=.9] 
(0,0) node (x11) {\brb{a, \sigmabarinv_* g_{\sigmabar\anga} \sigmabar\angw}}
(x11)++(1,0) node (x12) {\brb{a_{\sigmabarinv}, g_{\sigmabar\anga} \sigmabar\angw}}
;
\draw[1cell=.9] 
(x11) edge node {(\grod)_{\sigmabar,g,\ang{a,w}}} node[swap] {= (\sigmabar,1)} (x12)
;
\end{tikzpicture}
\end{equation}
The right pasting diagram in \cref{pseudosmf-topeq} yields the following morphism in $\grod Z$, with the same domain, respectively, codomain, as \cref{grodsigmabargaw}.
\begin{equation}\label{gagrodsigmaphione}
\left\{\begin{aligned}
&\ga \scmap{(\grod)_{\sigma,\phi}; \ang{1_{\phi^{\sigma(j)}}}_{j=1}^n}_{\ang{a,w}} &&\\
&= (\grod)_{\sigma,\phi, \bang{\phi^{\sigma(j)} \ang{(a_{ji},w_{ji})}_{i=1}^{\ell_{\sigma(j)}}}_{j=1}^n} & \phantom{M} & \text{by \cref{permcatcomposite}}\\
&= (\grod)_{\sigma,\phi, \bang{(a_j, (\phi_{\sigma(j)})_{\ang{a_j}} \ang{w_j})}_{j=1}^n} & \phantom{M} & \text{by \cref{grodphiangajxj}}\\
&= (\sigma,1) & \phantom{M} & \text{by \cref{grodsigmacomponent}}
\end{aligned}\right.
\end{equation}
The morphisms $(\sigmabar,1)$ and $(\sigma,1)$ in, respectively, \cref{grodsigmabargaw,gagrodsigmaphione}, are equal because the following diagram in $\cD$ is commutative by the definition of the block permutation $\sigmabar$.
\[\begin{tikzpicture}[xscale=3.5,yscale=1,vcenter]
\draw[0cell=.9] 
(0,0) node (x11) {\txotimes_{j=1}^n \txotimes_{i=1}^{\ell_{\sigma(j)}} a_{ji}}
(x11)++(1,0) node (x12) {\txotimes_{j=1}^n \txotimes_{i=1}^{\ell_j} a_{\sigmainv(j),i}}
(x11)++(0,-1) node (x21) {\txotimes_{j=1}^n a_j}
(x12)++(0,-1) node (x22) {\txotimes_{j=1}^n a_{\sigmainv(j)}}
;
\draw[1cell=.9] 
(x11) edge node {\sigmabar} (x12)
(x21) edge node {\sigma} (x22)
(x11) edge[equal] (x21)
(x12) edge[equal] (x22)
;
\end{tikzpicture}\]
This proves the top equivariance axiom.

\medskip
\emph{The Bottom Equivariance Axiom} \cref{pseudosmf-boteq}.

Suppose given
\begin{itemize}
\item permutations $\tau_j \in \Sigma_{\ell_j}$ for $j \in \{1,\ldots,n\}$,
\item additive natural transformations $\phi$ and $\phi_j$ as in \cref{phiphijtensorxjz}, and
\item an $(\ell_1 + \cdots + \ell_n)$-tuple of objects
\[\ang{a,w} = \ang{\ang{(a_{ji}, w_{ji})}_i}_j \in \txprod_{j=1}^n \txprod_{i=1}^{\ell_j} (\grod W_{j,\tau_j(i)})\]
with each $a_{ji} \in \cD$ and $w_{ji} \in W_{j,\tau_j(i)}\, a_{ji}$.
\end{itemize}
We use the following notation, along with $\ang{\phi_j}$ and $g$ in \cref{phijgphij}.
\[\left\{
\scalebox{.9}{$\begin{aligned}
\phibar &= \grod \phi & \ell &= \ell_1 + \cdots + \ell_n\\
\tautimes &= \tau_1 \times \cdots \times \tau_n \in \Sigma_\ell & \taumtimes &= (\tautimes)^\inv\\
\ang{a_j} &= \ang{a_{ji}}_{i=1}^{\ell_j} \in \cD^{\ell_j} & \ang{w_j} &= \ang{w_{ji}}_{i=1}^{\ell_j} \in \txprod_{i=1}^{\ell_j} W_{j,\tau_j(i)}\, a_{ji}\\
\anga &= \ang{\ang{a_j}}_{j=1}^n \in \cD^\ell & \angw &= \ang{\ang{w_j}}_{j=1}^n \in \txprod_{j=1}^n \txprod_{i=1}^{\ell_j} W_{j,\tau_j(i)}\, a_{ji}\\
a_j &= \txotimes_{i=1}^{\ell_j} a_{ji} \in \cD & a &= \txotimes_{j=1}^n a_j \in \cD\\
(a_j)_{\tau_j^\inv} &= \txotimes_{i=1}^{\ell_j} a_{j,\tau_j^\inv(i)} & a_{\taumtimes} &= \txotimes_{j=1}^n (a_j)_{\tau_j^\inv} 
\end{aligned}$}\right.\]
By definition \cref{grodsigmacomponent}, for the above data the left pasting diagram in \cref{pseudosmf-boteq} yields the following morphism in $\grod Z$.
\begin{equation}\label{grodtautimesgaw}
\begin{tikzpicture}[xscale=5.5,yscale=1.2,vcenter]
\draw[0cell=.9] 
(0,0) node (x11) {\brb{a, \taumtimes_* g_{\tautimes\anga} \tautimes\angw}}
(x11)++(1,0) node (x12) {\brb{a_{\taumtimes}, g_{\tautimes\anga} \tautimes\angw}}
;
\draw[1cell=.9] 
(x11) edge node {(\grod)_{\tautimes,g,\ang{a,w}}} node[swap] {= (\tautimes,1)} (x12)
;
\end{tikzpicture}
\end{equation}

The right pasting diagram in \cref{pseudosmf-boteq} yields the following morphism in $\grod Z$, with the same domain, respectively, codomain, as \cref{grodtautimesgaw}.
\begin{equation}\label{gaonegrodphi}
\left\{\begin{split}
&\ga \scmap{1_{\phibar}; \ang{(\grod)_{\tau_j,\phi_j}}_{j=1}^n}_{\ang{a,w}}&&\\
&= \phibar \bang{(\grod)_{\tau_j, \phi_j, \ang{(a_{ji},w_{ji})}_{i=1}^{\ell_j}}}_{j=1}^n & \phantom{M} & \text{by \cref{permcatcomposite}}\\
&= \phibar \bang{(\tau_j, 1)}_{j=1}^n & \phantom{M} & \text{by \cref{grodsigmacomponent}}\\
&= \brb{\txotimes_{j=1}^n \tau_j, \phi\ang{1}} & \phantom{M} & \text{by \cref{grodphiangfjpj}}\\
&= \brb{\txotimes_{j=1}^n \tau_j, 1} & \phantom{M} & \text{by functoriality}
\end{split}\right.
\end{equation}
In the middle line of \cref{gaonegrodphi}, $(\tau_j,1)$ is the following morphism in $\grod X_j$.
\[\begin{tikzpicture}[xscale=5.3,yscale=1.2,vcenter]
\draw[0cell=.9] 
(0,0) node (x11) {\brb{a_j, (\tau_j^\inv)_* (\phi_j)_{\tau_j\ang{a_j}} \tau_j\ang{w_j}}}
(x11)++(1,0) node (x12) {\brb{(a_j)_{\tau_j^\inv}, (\phi_j)_{\tau_j\ang{a_j}} \tau_j\ang{w_j}}}
;
\draw[1cell=.9] 
(x11) edge node {(\tau_j,1)} (x12)
;
\end{tikzpicture}\]
In the fourth line of \cref{gaonegrodphi}, $\phi\ang{1}$ is the morphism in $Za_{\taumtimes}$
\[\phi\ang{1} = \phi_{\ang{(a_j)_{\tau_j^\inv}}_{j=1}^n} (\overbracket[.5pt]{1,\ldots,1}^{n}).\]
This is the identity morphism by the functoriality of each component functor of $\phi$.  This yields the last equality in \cref{gaonegrodphi}.

The morphisms $(\tautimes,1)$ and $\brb{\txotimes_{j=1}^n \tau_j, 1}$ in, respectively, \cref{grodtautimesgaw,gaonegrodphi}, are equal because the following diagram in $\cD$ is commutative by the definition of the permutation $\tautimes$.
\[\begin{tikzpicture}[xscale=3.5,yscale=1,vcenter]
\draw[0cell=.9] 
(0,0) node (x11) {\txotimes_{j=1}^n \txotimes_{i=1}^{\ell_j} a_{ji}}
(x11)++(1,0) node (x12) {\txotimes_{j=1}^n \txotimes_{i=1}^{\ell_j} a_{j,\tau_j^\inv(i)}}
(x11)++(0,-1) node (x21) {\txotimes_{j=1}^n a_j}
(x12)++(0,-1) node (x22) {\txotimes_{j=1}^n (a_j)_{\tau_j^\inv}}
;
\draw[1cell=.9] 
(x11) edge node {\tautimes} (x12)
(x21) edge node {\txotimes_j \tau_j} (x22)
(x11) edge[equal] (x21)
(x12) edge[equal, shorten <=-.5ex] (x22)
;
\end{tikzpicture}\]
This proves the bottom equivariance axiom.
\end{proof}

Next is the main result of this chapter.  It relates the $\Cat$-multicategories $\DCat$ and $\permcatsg$ in, respectively, \cref{thm:dcatcatmulticat,thm:permcatenrmulticat}, via the Grothendieck construction.  

\begin{theorem}\label{thm:grocatmultifunctor}\index{bipermutative category!Grothendieck construction}\index{Grothendieck construction!bipermutative category}
Suppose $\cD$ is a small tight bipermutative category.  
\begin{enumerate}[label=(\roman*)]
\item\label{thm:grocatmultifunctor-i} The Grothendieck construction
\begin{equation}\label{grodpscatmultifunctor}
\begin{tikzcd}[column sep=large,every label/.append style={scale=.85}]
\DCat \ar{r}{\grod} & \permcatsg
\end{tikzcd}
\end{equation}
is a \index{pseudo symmetric!Cat-multifunctor@$\Cat$-multifunctor!Grothendieck construction}\index{Cat-multifunctor@$\Cat$-multifunctor!pseudo symmetric!Grothendieck construction}\index{Grothendieck construction!pseudo symmetric Cat-multifunctor@pseudo symmetric $\Cat$-multifunctor}pseudo symmetric $\Cat$-multifunctor.
\item\label{thm:grocatmultifunctor-ii} If the multiplicative braiding in $\cD$ is the identity natural transformation, then $\grod$ is a $\Cat$-multifunctor.
\item\label{thm:grocatmultifunctor-iii} If $\grod$ is a $\Cat$-multifunctor, then 
\begin{equation}\label{xyisyx}
x \otimes y = y \otimes x
\end{equation}
for all objects $x,y \in \cD$.
\end{enumerate}
\end{theorem}

\begin{proof}
For the Grothendieck construction $\grod$ we established the following.
\begin{itemize}
\item The object assignment $X \mapsto\! \grod X$ is in \cref{grodxpermutative}.
\item The multimorphism functors 
\[\begin{tikzcd}[column sep=large,every label/.append style={scale=.85}]
\DCat\scmap{\angX;Z} \ar{r}{\grod} & \permcatsg\scmap{\ang{\grod X};\grod Z}
\end{tikzcd}\]
are in \cref{def:groconarityzero,grodmultimorphismfunctor} for, respectively, arity 0 and positive arity.
\item It preserves the colored units by \cref{grodpreservesunits}.
\item It preserves the composition by \cref{grodpreservescomposition}.
\item Its pseudo symmetry isomorphisms $(\grod)_\sigma$ are in \cref{def:gropseudosymmetry}.
\item These pseudo symmetry isomorphisms are well-defined natural isomorphisms by \cref{grodsigmaphilinear,grodsigmaphinlinear,grodsigmanaturaliso}.
\item The four axioms in \cref{def:pseudosmultifunctor} for the pseudo symmetry isomorphisms are proved in \cref{grodpseudosymaxioms}.
\end{itemize}
This shows that $\grod$ is a pseudo symmetric $\Cat$-multifunctor, proving assertion \cref{thm:grocatmultifunctor-i}.  

For \cref{thm:grocatmultifunctor-ii}, suppose that the multiplicative braiding $\betate$ in $\cD$ is the identity. Since $\grod$ preserves the colored units and composition (\cref{grodpreservesunits,grodpreservescomposition}), we must show that $\grod$ preserves the symmetric group action.  In other words, we want to show that the boundary of the diagram \cref{gropseudosymmetry} is commutative.  Since $\betate = 1$, the two objects in \cref{grodsigmadomain} are equal.  So 
\[\grod \phi^\sigma = (\grod\phi)^\sigma\]
for each additive natural transformation $\phi$ in $\DCat\smscmap{\angX;Z}$.  This means that the boundary diagram \cref{gropseudosymmetry} is commutative on objects.  

To prove the commutativity of the boundary diagram \cref{gropseudosymmetry} on morphisms, suppose $\Phi \cn \phi \to \psi$ is an additive modification in $\DCat\smscmap{\angX;Z}$ and $\ang{(a_j,x_j)}$ is an $n$-tuple of objects as in \cref{ajxjgroxsigmaj}.  Then there are equalities as follows.
\begin{equation}\label{grodPhisigmaangajxj}
\left\{\begin{aligned}
& (\grod \Phi)^\sigma_{\ang{(a_j,x_j)}} && \\
&= (\grod \Phi)_{\ang{a_{\sigmainv(j)}, x_{\sigmainv(j)}}} && \text{by \cref{permcatsigmaaction}}\\
&= \left(1_{a_{\sigmainv}} \scs \left(\Phi_{\sigma\anga}\right)_{\sigma\angx}\right) && \text{by \cref{grodPhiajxj}}\\
&= \left(1_a \scs Z(\id_n) \left(\Phi_{\sigma\anga}\right)_{\sigma\angx}\right) && \scalebox{.8}{$\text{by $\betate = 1$ and functoriality of $Z$}$}\\
&= \left(1_a \scs Z(\sigmainv) \left(\Phi_{\sigma\anga}\right)_{\sigma\angx}\right) && \text{by $\betate = 1$}\\
&= \left(1_a \scs \big(\Phi^\sigma_{\anga}\big)_{\angx}\right) && \text{by \cref{Phisigmaanga}}\\
&= \left(\grod \Phi^\sigma\right)_{\ang{(a_j,x_j)}} && \text{by \cref{grodPhiajxj}}
\end{aligned}\right.
\end{equation}
This shows that the boundary diagram \cref{gropseudosymmetry} is commutative on morphisms.  Thus $\grod$ is a $\Cat$-multifunctor, proving assertion \cref{thm:grocatmultifunctor-ii}.

For \cref{thm:grocatmultifunctor-iii}, suppose that $\grod$ is a $\Cat$-multifunctor.  So $\grod$ preserves the symmetric group action, which means that the boundary of the diagram \cref{gropseudosymmetry} is commutative.  The first component of the computation \cref{grodsigmadomain}, which uses the notation in \cref{asigmainvtoa}, shows that
\[\txotimes_{j=1}^n a_j = a = a_{\sigmainv} = \txotimes_{j=1}^n a_{\sigmainv(j)}.\]
This implies the condition \cref{xyisyx}.
\end{proof}

\begin{example}[Obstruction to $\Cat$-Multifunctoriality]\label{ex:grodpseudosymmetric}
For the small tight bipermutative category $\cA$ in \cref{ex:mandellcategory}, the condition \cref{xyisyx} is \emph{not} satisfied by the definition \cref{Amtensorn} of $\otimes$ in $\cA$.  By \cref{thm:grocatmultifunctor} \cref{thm:grocatmultifunctor-i,thm:grocatmultifunctor-iii} the Grothendieck construction $\groa$ is a pseudo symmetric $\Cat$-multifunctor but \emph{not} a $\Cat$-multifunctor.
\end{example}

\begin{explanation}[Further Criteria for $\Cat$-Multifunctoriality]\label{expl:grodcatmultifunctor} 
Under the hypothesis of \cref{thm:grocatmultifunctor} \cref{thm:grocatmultifunctor-iii}---that $\grod$ is a $\Cat$-multifunctor---something beyond the condition \cref{xyisyx} is true.  Since 
\[\grod \phi^\sigma = (\grod \phi)^\sigma,\]
the second component of the computation \cref{grodsigmadomain} shows that the diagram 
\[\begin{tikzpicture}[xscale=3.5,yscale=1.3,vcenter]
\draw[0cell=.85]
(0,0) node (x11) {\txprod_{j=1}^n X_j a_j}
(x11)++(1,0) node (x12) {Za}
(x11)++(0,-1) node (x21) {Za}
(x12)++(0,-1) node (x22) {Za_{\sigmainv}}
;
\draw[1cell=.9]  
(x11) edge node {\phi_{\anga}} (x12)
(x11) edge node[swap] {\phi_{\anga}} (x21)
(x21) edge node {Z(\sigma)} (x22)
(x22) edge[equal] (x12)
;
\end{tikzpicture}\]
is commutative on objects for any
\begin{itemize}
\item permutation $\sigma \in \Sigma_n$, 
\item object $\anga = \ang{a_j}_{j=1}^n \in \cD^n$, 
\item additive symmetric monoidal functors $Z$ and $\angX = \ang{X_j}_{j=1}^n$ in $\DCat$, and
\item additive natural transformation $\phi \in \DCat\smscmap{\angX;Z}$.
\end{itemize} 
Moreover, for any additive modification $\Phi \cn \phi \to \psi$ in $\DCat\smscmap{\angX;Z}$, the second component of the computation \cref{grodPhisigmaangajxj} yields the morphism equality
\[\big(\Phi_{\anga}\big)_{\angx} = Z(\sigma) \big(\Phi_{\anga}\big)_{\angx}\]
in $Za$ for each object $\angx \in \txprod_{j=1}^n X_j a_j$.

In the conclusion of \cref{thm:grocatmultifunctor} \cref{thm:grocatmultifunctor-iii}, we only include the condition \cref{xyisyx} because it is easy to state and sufficient for our main application, which is the case $\cA$ in \cref{ex:grodpseudosymmetric}.
\end{explanation}

\chapter{Permutative Opfibrations from Bipermutative-Indexed Categories}
\label{ch:permfib}
\subsection*{Context}

Suppose $(\cD,\oplus,\otimes)$ is a small tight bipermutative category (\cref{def:embipermutativecat}).  We proved in \cref{thm:grocatmultifunctor} that there is a pseudo symmetric $\Cat$-multifunctor (\cref{def:pseudosmultifunctor})
\begin{equation}\label{dcatgrodpermcatsg}
\begin{tikzcd}[column sep=large,every label/.append style={scale=.85}]
\DCat \ar{r}{\grod} & \permcatsg
\end{tikzcd}
\end{equation}
with object assignment given by the Grothendieck construction 
\[\begin{tikzcd}
(X,X^2,X^0) \ar[mapsto]{r} & \big(\grod X, \gbox, (\zero, X^0*), \betabox\big).
\end{tikzcd}\]
Moreover, $\grod$ is a $\Cat$-multifunctor if the multiplicative braiding $\betate$ in $\cD$ is the identity natural transformation.  The weaker converse says that, if $\grod$ is a $\Cat$-multifunctor, then $x \otimes y = y \otimes x$ for all objects $x,y \in \cD$.  

The Grothendieck construction $\grod X$ has a canonical first-factor projection functor 
\[\begin{tikzcd}[column sep=large]
\grod X \ar{r}{U_X} & \cD.
\end{tikzcd}\]
This functor $U_X$ has some nice lifting properties that make it into a permutative opfibration over $\cD$.  As the natural codomain of the Grothendieck construction, we would like to consider not just small permutative categories but small permutative opfibrations over $\cD$.  

\subsection*{Purpose}

The main purpose of this chapter is to lift the Grothendieck construction \cref{dcatgrodpermcatsg} to a non-symmetric $\Cat$-multifunctor
\begin{equation}\label{dcatgrodpfibd}
\begin{tikzcd}[column sep=large,every label/.append style={scale=.85}]
\DCat \ar{r}{\grod} & \pfibd.
\end{tikzcd}
\end{equation}
\begin{itemize}
\item The codomain $\pfibd$ is a non-symmetric $\Cat$-multicategory with small permutative opfibrations over $\cD$ as objects.
\item The objects in the multimorphism categories in $\pfibd$ are called \emph{opcartesian strong $n$-linear functors}.  These are strong $n$-linear functors with extra properties that generalize Cartesian functors between Grothendieck fibrations. 
\item The morphisms in the multimorphism categories in $\pfibd$ are called \emph{opcartesian $n$-linear transformations}.  These are $n$-linear natural transformations with an extra property that generalize vertical natural transformations between Cartesian functors. 
\item The non-symmetric $\Cat$-multicategory structure on $\pfibd$ is inherited from $\permcatsg$.
\end{itemize}
Permutative opfibrations over $\cD$ are defined using only the additive structure $\Dplus$.  On the other hand, opcartesian $n$-linear functors require the additive structure $\Dplus$, the multiplicative structure $\Dte$, and Laplaza's Coherence \cref{thm:laplaza-coherence-1}.

The Grothendieck constructions \cref{dcatgrodpermcatsg,dcatgrodpfibd} fit together in the following commutative diagram of non-symmetric $\Cat$-multifunctors.
\[\begin{tikzpicture}[xscale=2.75,yscale=1.2,vcenter]
\draw[0cell=.9]
(0,0) node (dc) {\DCat}
(dc)++(1,0) node (pc) {\permcatsg}
(pc)++(0,1) node (pf) {\pfibd}
;
\draw[1cell=.9]  
(dc) edge node[pos=.6] {\grod} (pc)
(pf) edge node {\U} (pc)
(dc) edge[bend left=20] node[pos=.7] {\grod} (pf)
;
\end{tikzpicture}\]
The non-symmetric $\Cat$-multifunctor $\U$ sends each small permutative opfibration over $\cD$ to its domain small permutative category.  Between multimorphism categories, $\U$ is given by subcategory inclusions.  Moreover, if the multiplicative braiding $\betate$ in $\cD$ is the identity natural transformation, then the following two statements hold:
\begin{itemize}
\item The symmetric group action in $\permcatsg$ restricts to $\pfibd$ to make it into a $\Cat$-multicategory.
\item Both $\U$ and $\grod \cn \DCat \to \pfibd$ are $\Cat$-multifunctors.
\end{itemize}

A weaker converse is also true: if $\pfibd$ inherits a well-defined symmetric group action from $\permcatsg$, then $\otimes = \tensorop$ in $\cD$.  The very restrictive condition $\otimes = \tensorop$ fails in most bipermutative categories of interest, including $\Finsk$, $\Fset$, $\Fskel$, and $\cA$ in \cref{sec:examples}.  Therefore, in each of those cases, $\pfibd$ is genuinely a non-symmetric $\Cat$-multicategory and \emph{not} a $\Cat$-multicategory.   

In \cref{ch:gromultequiv} we show that the non-symmetric $\Cat$-multifunctor $\grod$ in \cref{dcatgrodpfibd} is actually a non-symmetric \emph{$\Cat$-multiequivalence}.  This means that $\grod$ is
\begin{itemize}
\item essentially surjective on objects between the underlying 1-ary categories and
\item an isomorphism between multimorphism categories.
\end{itemize} 
Furthermore, $\grod$ is a $\Cat$-multiequivalence if $\betate$ in $\cD$ is the identity.

\subsection*{Summary}

Extending the tables in the introductions of \cref{ch:diagram,ch:multigro}, the following table summaries $\DCat$, $\pfibd$, and $\grod$.  We abbreviate \emph{symmetric monoidal}, \emph{multimorphism category},  \emph{natural transformations}, and \emph{opcartesian} to, respectively, \emph{sm}, \emph{mc}, \emph{nt}, and \emph{op}.
\begin{center}
\resizebox{\columnwidth}{!}{%
{\renewcommand{\arraystretch}{1.4}%
{\setlength{\tabcolsep}{1ex}
\begin{tabular}{|c|cc|cc|c|}\hline
& $\DCat$ &(\ref{thm:dcatcatmulticat}) & $\pfibd$ &(\ref{thm:pfibdmulticat}, \ref{cor:pfibdtopermcatsg}) & $\grod\,\,$ (\ref{thm:dcatpfibd}) \\ \hline
objects 
& additive sm functors &(\ref{def:additivesmf}) & small permutative opfibrations &(\ref{def:permutativefibration}) & \ref{grodxpermopfib}\\ \hline
mc objects 
& additive nt &(\ref{def:additivenattr}) & op strong multilinear functors &(\ref{def:opcartnlinearfunctor}) & \ref{grodfactorarityzero}, \ref{grodfactorarityn}\\ \hline
mc morphisms 
& additive modifications &(\ref{def:additivemodification}) & op multilinear transformations &(\ref{def:opcartnlineartr}) & \ref{grodfactorarityzero}, \ref{grodfactorarityn}\\ \hline
colored units & identity nt & (\ref{zcoloredunit}) & identity op 1-linear functors & (\ref{ex:idomf}) & \ref{grodpreservesunits}, \ref{def:pfibd}\\ \hline
composition & \ref{dcatgamma}, \ref{expl:dcatgamma} && \ref{pfibdcompobjwelldef}, \ref{pfibdcompmorwelldef} && \ref{grodpreservescomposition}, \ref{grodfactorarityzero}, \ref{grodfactorarityn}\\ \hline
symmetry & \ref{dcatsigmaaction} && \ref{pfibdsymwelldefobj}, \ref{pfibdsymwelldefmor} && \ref{thm:grocatmultifunctor}, \ref{grodfactorarityn}\\ \hline
sm
& no & (\ref{expl:dcatbipermzero}) & no & (\ref{pfibdsymwelldefobj} \eqref{pfibdsymwelldefobj-ii}) & no\\ \hline
\end{tabular}}}}
\end{center}
\smallskip

\subsection*{Organization}

In \cref{sec:permopfib} we first review the definition of a split opfibration.  A functor $P \cn \E \to \B$ is an opfibration if and only if its opposite functor $P^\op$ is a Grothendieck fibration.  A split opfibration is an opfibration with a cleavage that satisfies the unitarity and multiplicativity conditions.  For a permutative category $\E$, a \emph{permutative opfibration} $P \cn \E \to \Dplus$ is a split opfibration that satisfies the following three conditions:
\begin{itemize}
\item $P$ is a strict symmetric monoidal functor.
\item Chosen opcartesian lifts are closed under the monoidal product in $\E$.
\item Each component of the braiding in $\E$ is the chosen opcartesian lift of the corresponding component of the braiding in $\cD$.
\end{itemize} 
Small permutative opfibrations over $\cD$ are the objects in $\pfibd$.

In \cref{sec:opmultilinear} we define opcartesian strong $n$-linear functors and opcartesian $n$-linear transformations.  They are, respectively, the $n$-ary 1-cells and $n$-ary 2-cells in $\pfibd$.  By definition an opcartesian $n$-linear functor is an $n$-linear functor with extra properties but \emph{no} extra structure, and the same remark applies to opcartesian $n$-linear transformations.  Each multimorphism category in $\pfibd$ is defined as a subcategory of the corresponding multimorphism category in $\permcatsg$.

In \cref{sec:popfibmulticat} we define the (non-symmetric) $\Cat$-multicategory structure on $\pfibd$ as the one inherited from $\permcatsg$; see \cref{thm:pfibdmulticat,cor:pfibdtopermcatsg}.  First we check that the composition in $\pfibd$ is well defined.  Then we check that, if $\betate = 1$ in $\cD$, then the symmetric group action in $\pfibd$ is well defined.  The weaker converse says that if the symmetric group action in $\pfibd$ is well defined, then $\otimes = \tensorop$ in $\cD$.

In \cref{sec:popfibgrothendieck} we construct the lifted Grothendieck construction in \cref{dcatgrodpfibd}; see \cref{thm:dcatpfibd}.  We already know that
\begin{itemize}
\item $\grod$ in \cref{dcatgrodpermcatsg} is a pseudo symmetric $\Cat$-multifunctor and
\item each multimorphism category in $\pfibd$ is a subcategory of the corresponding multimorphism category in $\permcatsg$.
\end{itemize} 
Thus most of the work in lifting $\grod$ to $\pfibd$ is to show that $\grod$ lands in the multimorphism categories in $\pfibd$.

We remind the reader of our left normalized bracketing \cref{expl:leftbracketing} for iterated monoidal product.

\section{Permutative Opfibrations}
\label{sec:permopfib}

To define permutative opfibrations, we first define opfibrations.  A general reference for opfibrations is \cite[Ch.\! 9]{johnson-yau}, which discusses fibrations.  Results there apply here as well by considering the opposite functor between the opposite categories.  Permutative opfibrations are in \cref{def:permutativefibration}.  Small permutative opfibrations are the objects in a non-symmetric $\Cat$-multicategory $\pfibd$ that we will discuss in subsequent sections.

\begin{definition}\label{def:opfibration}
Suppose $P \cn \E \to \B$ is a functor between categories.  Define the following.
\begin{enumerate}[label=(\roman*)]
\item\label{def:opfibration-i} A \emph{fore-lift}\index{fore-lift} is a pair $\fl{y}{f}$ consisting of
\begin{itemize}
\item an object $y$ in $\E$ and
\item a morphism $f \cn Py \to a$ in $\B$ for some object $a$.
\end{itemize} 
\item\label{def:opfibration-ii} A \emph{lift}\index{lift} of a fore-lift $\fl{y}{f}$ is a morphism $\fbar \cn y \to y_f$ in $\E$ such that
\[Py_f = a \andspace P\fbar = f \inspace B,\]
which we depict as follows.
\[\begin{tikzpicture}[xscale=2.5, yscale=1.2]
\draw[0cell=.9]
(0,0) node (Y) {y}
(Y)++(0,-1) node (Yf) {y_f}
(Y)++(.3,-.4) node (s) {} 
(s)++(.5,0) node (t) {} 
(t)++(.3,.4) node (Yp) {Py}
(Yp)++(0,-1) node (A) {a}
;
\draw[1cell=.9]
(Y) edge node[swap] {\fbar} (Yf)
(s) edge[|->] node {P} (t)
(Yp) edge node {f} (A)
;
\end{tikzpicture}\]
If $y$ is clear from the context, then we also say that $\fbar$ is a \emph{lift of $f$}.
\item\label{def:opfibration-iii} A \emph{fore-raise}\index{fore-raise} is a triple $\fr{g}{h}{f}$ consisting of
\begin{itemize}
\item two morphisms $g \cn y \to x$ and $h \cn y \to z$ in $\E$ with the same domain and
\item a morphism $f \cn Px \to Pz$ in $\B$ such that
\end{itemize}
\[f (Pg) = Ph.\]
\item\label{def:opfibration-iv} A \emph{raise}\index{raise} of a fore-raise $\fr{g}{h}{f}$ is a morphism $\ftil \cn x \to z$ in $\E$ such that
\[\ftil g = h \in \E \andspace P\ftil = f \in \B,\]
which we depict as follows.
\[\begin{tikzpicture}[xscale=2.5, yscale=1.3]
\draw[0cell=.9]
(0,0) node (X) {x}
(X)++(1,0) node (Y) {y}
(X)++(.5,.85) node (Z) {z}
(Y)++(.3,.4) node (s) {}
(s)++(.5,0) node (t) {}
(t)++(.3,-.4) node (Xp) {Px}
(Xp)++(1,0) node (Yp) {Py}
(Xp)++(.5,.85) node (Zp) {Pz}
;
\draw[1cell=.9]
(Y) edge node[swap] {g} (X)
(Y) edge node[swap] {h} (Z)
(X) edge node {\ftil} (Z)
(s) edge[|->] node {P} (t)
(Yp) edge node[swap] {Pg} (Xp)
(Yp) edge node[swap] {Ph} (Zp)
(Xp) edge node {f} (Zp)
;
\end{tikzpicture}\]
\item\label{def:opfibration-v} A morphism $g \cn y \to x$ in $\E$ is called an \emph{opcartesian morphism}\index{opcartesian morphism} if each fore-raise of the form $\fr{g}{h}{f}$ has a \emph{unique} raise.
\item\label{def:opfibration-vi} A lift of a fore-lift $\fl{y}{f}$ is called an \emph{opcartesian lift}\index{opcartesian lift}\index{lift!opcartesian} if it is also an opcartesian morphism.
\item\label{def:opfibration-vii} $P \cn \E \to \B$ is called an \emph{opfibration}\index{opfibration} over $\B$ if each fore-lift admits an opcartesian lift.  In this case, $P$, $\E$, and $\B$ are called, respectively, the \index{projection functor}\emph{projection functor}, the \index{total category}\index{category!total}\emph{total category}, and the \index{base category}\index{category!base}\emph{base category}.  A \index{cleavage}\emph{cleavage} is a choice of an opcartesian lift $\fbar$ for each fore-lift $\fl{y}{f}$, in which case $\fbar$ is called the \emph{chosen opcartesian lift}\index{chosen opcartesian lift}\index{opcartesian lift!chosen} of $\fl{y}{f}$.  An opfibration equipped with a cleavage is called a \index{opfibration!cloven}\index{cloven opfibration}\emph{cloven opfibration}. 
\item\label{def:opfibration-viii} A \emph{split opfibration}\index{split opfibration}\index{opfibration!split} is a cloven opfibration $P \cn \E \to \B$ such that the following two axioms hold.
\begin{description}
\item[Unitarity]\index{unitarity} For each object $y \in \E$, the identity morphism $1_y$ is the chosen opcartesian lift of $\fl{y}{1_{Py}}$.
\item[Multiplicativity]\index{multiplicativity} Suppose given
\begin{itemize}
\item an object $y$ in $\E$ and
\item composable morphisms $f \cn Py \to a$ and $g \cn a \to a'$ in $\B$.
\end{itemize} 
Suppose
\begin{itemize}
\item $\fbar \cn y \to y_f$ is the chosen opcartesian lift of $\fl{y}{f}$ and
\item $\gbar \cn y_f \to (y_f)_g$ is the chosen opcartesian lift of $\fl{y_f}{g}$, as depicted below.
\end{itemize} 
\[\begin{tikzpicture}[xscale=2.5, yscale=1.2]
\draw[0cell=.9]
(0,0) node (Y) {y}
(Y)++(0,-1) node (Yf) {y_f}
(Yf)++(0,-1) node (Yfg) {(y_f)_g}
(Yf)++(.3,0) node (s) {} 
(s)++(.5,0) node (t) {} 
(t)++(.3,1) node (Yp) {Py}
(Yp)++(0,-1) node (A) {a}
(A)++(0,-1) node (Ap) {a'}
;
\draw[1cell=.9]
(Y) edge node[swap] {\fbar} (Yf)
(Yf) edge node[swap] {\gbar} (Yfg)
(s) edge[|->] node {P} (t)
(Yp) edge node {f} (A)
(A) edge node {g} (Ap)
;
\end{tikzpicture}\]
Then the composite 
\[\begin{tikzcd}[column sep=large]
y \ar{r}{\gbar \fbar} & (y_f)_g
\end{tikzcd} \inspace \E\]
is the chosen opcartesian lift of the fore-lift $\fl{y}{gf}$.
\end{description}
\item\label{def:opfibration-ix} An opfibration $P \cn \E \to \B$ is \emph{small}\index{opfibration!small} if both $\E$ and $\B$ are small categories. 
\end{enumerate}
This finishes the definition.
\end{definition}

\begin{explanation}\label{expl:opfibterminology}
Consider \cref{def:opfibration}.
\begin{enumerate}
\item\label{expl:opfibterminology-i} A functor $P \cn \E \to \B$ is an opfibration if and only if its opposite functor $P^\op \cn \E^\op \to \B^\op$ is a fibration in the sense of \cite[9.1.1]{johnson-yau}.  In the literature an opfibration is sometimes called a \index{cofibred category}\index{category!cofibred}\emph{cofibred category}.
\item\label{expl:opfibterminology-ii} The data $\fl{y}{f}$ and $\fr{g}{h}{f}$ for, respective, a fore-lift and a fore-raise are often not given names in the literature.  In the setting of fibrations in \cite[Ch.\! 9]{johnson-yau}, the corresponding concepts are called, respectively, a \emph{pre-lift} and a \emph{pre-raise}.  A lift is sometimes called a \emph{transport morphism}\index{transport morphism} in the literature.
\item\label{expl:opfibterminology-iii} A cloven opfibration can always be arranged to satisfy the unitarity axiom.  However, the multiplicativity axiom is not automatic.  By the \index{Grothendieck Fibration Theorem}Grothendieck Fibration Theorem \cite[9.5.6]{johnson-yau}, for each small category $\B$, there is an explicitly-constructed 2-monad\index{2-monad} on $\Cat/\B$ whose pseudo, respectively, strict, algebras are in canonical bijections with cloven, respectively, split, opfibrations over $\B$.
\item An opfibration is a functor with additional properties, namely, the existence of an opcartesian lift for each fore-lift, but no extra structure.  On the other hand, a \emph{split} opfibration is an opfibration with extra structure, namely, the \emph{chosen} opcartesian lift of each fore-lift.  As a consequence, notions that are built upon split opfibrations usually have appropriate axioms regarding the chosen opcartesian lifts.  Examples include  \cref{def:permutativefibration,def:opcartnlinearfunctor}.\defmark
\end{enumerate}
\end{explanation}

\begin{example}\label{ex:idopfibration}
Consider \cref{def:opfibration}.
\begin{enumerate}
\item\label{ex:idopfibration-i}
For a functor $P \cn \E \to \B$, each isomorphism in $\E$ is an opcartesian morphism.  In particular, each identity morphism in $\E$ is an opcartesian morphism.
\item\label{ex:idopfibration-ii}
Each identity functor $1_{\B} \cn \B \to \B$ is a split opfibration.  For each fore-lift $\fl{y}{f}$, the chosen opcartesian lift is $f$.\defmark
\end{enumerate}
\end{example}

Next we consider the situation when each of $P$, $\E$, and $\B$ is equipped with a permutative structure (\cref{def:symmoncat,def:monoidalfunctor}).  Eventually we will assume that the base category $\B$ is a small tight bipermutative category, but for the next definition it only needs a permutative structure.

\begin{definition}\label{def:permutativefibration}
Suppose $(\E,\oplus,e,\betaplus)$ and $(\cD,\oplus,e,\betaplus)$ are permutative categories.  A \emph{permutative opfibration}\index{permutative opfibration}\index{opfibration!permutative} over $\cD$
\[\begin{tikzcd}[column sep=large]
\E \ar{r}{P} & \cD
\end{tikzcd}\]
is a split opfibration (\cref{def:opfibration} \cref{def:opfibration-viii}) that satisfies the following three axioms.
\begin{enumerate}[label=(\roman*)]
\item\label{def:permutativefibration-i} 
$P$ is a strict symmetric monoidal functor.  This is called the \index{strict symmetry axiom}\emph{strict symmetry axiom}.
\item\label{def:permutativefibration-ii}
Suppose given fore-lifts
\begin{itemize}
\item $\fl{x}{f}$ with chosen opcartesian lift $\fbar$ and
\item $\fl{y}{g}$ with chosen opcartesian lift $\gbar$.
\end{itemize}
Then $\fbar \oplus \gbar$ is the chosen opcartesian lift of the fore-lift $\fl{x \oplus y}{f \oplus g}$.  This is called the \index{monoidal product axiom}\emph{monoidal product axiom}.
\item\label{def:permutativefibration-iii}
For any two objects $x,y \in \E$, the braiding
\[\begin{tikzcd}[column sep=large]
x \oplus y \ar{r}{\betaplus_{x,y}} & y \oplus x
\end{tikzcd} \inspace \E\]
is the chosen opcartesian lift of the fore-lift
\[\bfl{x \oplus y}{\betaplus_{Px,Py} \cn Px \oplus Py \to Py \oplus Px}.\]
This is called the \index{braiding axiom}\emph{braiding axiom}.
\end{enumerate}
This finishes the definition of a permutative opfibration.  Moreover, a permutative opfibration as above is \emph{small} if both $\E$ and $\cD$ are small categories. 
\end{definition}

\begin{explanation}\label{expl:permutativefibration}
Consider \cref{def:permutativefibration}.
\begin{enumerate}
\item\label{expl:permutativefibration-0} A permutative opfibration is a split opfibration with additional properties, namely, the axioms \cref{def:permutativefibration-i,def:permutativefibration-ii,def:permutativefibration-iii}, but \emph{no} extra structure.
\item\label{expl:permutativefibration-i} The strict symmetry axiom \cref{def:permutativefibration-i} means that the unit constraint $P^0$ and the monoidal constraint $P^2$ are identities (\cref{def:monoidalfunctor}).  They imply the following equalities for all objects $x,y \in \E$ and morphisms $p,q \in \E$.
\begin{equation}\label{strictsmfunctor}
\left\{\begin{split}
Pe &= e\\
P(x \oplus y) &= Px \oplus Py\\
P(p \oplus q) &= Pp \oplus Pq\\
P\betaplus_{x,y} &= \betaplus_{Px,Py}
\end{split}\right.
\end{equation}
The last three equalities in \cref{strictsmfunctor} ensure that axioms \cref{def:permutativefibration-ii,def:permutativefibration-iii} are well defined.  Conversely, if the functor $P$ satisfies the equalities in \cref{strictsmfunctor}, then we can define the unit constraint $P^0$ and the monoidal constraint $P^2$ to be identities, and $(P,1,1)$ is a strict symmetric monoidal functor.
\item\label{expl:permutativefibration-ii} The braiding axiom \cref{def:permutativefibration-iii} subsumes the unitarity condition of a split opfibration.  Indeed, suppose $x$ is the monoidal unit $e$ in $\E$.  Then 
\[\begin{split}
\betaplus_{e,y} &= 1_y \in \E \andspace\\ 
\betaplus_{Pe,Py} &= \betaplus_{e,Py} = 1_{Py} \in \cD
\end{split}\]
by the first equality in \cref{strictsmfunctor} and the right unit diagram \cref{braidedunity} in $\E$ and $\cD$.  Thus in this case the braiding axiom says that $1_y$ is the chosen opcartesian lift of $\fl{y}{1_{Py}}$, which is the unitarity condition.
\end{enumerate}
See \cref{qu:permopfibration} for a potential 2-monadic description of permutative opfibrations that is analogous to the Grothendieck Fibration Theorem mentioned in \cref{expl:opfibterminology} \eqref{expl:opfibterminology-iii}
\end{explanation}

\begin{example}[Identity Permutative Opfibrations]\label{ex:idpermopfibration}
For each permutative category $\cD$, the identity functor $1_{\cD} \cn \cD \to \cD$ is a permutative opfibration.  In fact, the identity functor is a split opfibration (\cref{ex:idopfibration} \eqref{ex:idopfibration-ii}) and satisfies the three axioms in \cref{def:permutativefibration}.
\end{example}

\section{Opcartesian Multilinear Functors and Transformations}
\label{sec:opmultilinear}

Throughout the rest of this chapter, we assume that 
\begin{equation}\label{Dstbipermutative}
\big(\cD, (\oplus, \zero, \betaplus), (\otimes, \tu, \betate), (\fal, \far)\big)
\end{equation}
is a small tight bipermutative category (\cref{def:embipermutativecat}) with 
\begin{itemize}
\item additive structure $\Dplus = (\cD, \oplus, \zero, \betaplus)$ and
\item multiplicative structure $\Dte = (\cD, \otimes, \tu, \betate)$.
\end{itemize}
We use $\cD$ as our base category.  As the next step in constructing the non-symmetric $\Cat$-multicategory $\pfibd$ with small permutative opfibrations over $\cD$ (\cref{def:permutativefibration}) as objects, in this section we define its multimorphism categories. 
\begin{itemize}
\item The $n$-ary 1-cells (\cref{def:opcartnlinearfunctor}) are \emph{opcartesian strong $n$-linear functors}.  They are strong $n$-linear functors that are also compatible with the given permutative opfibrations.
\item The $n$-ary 2-cells (\cref{def:opcartnlineartr}) are \emph{opcartesian $n$-linear transformations}.  They are $n$-linear natural transformations that are also compatible with the given permutative opfibrations.
\item The multimorphism categories of $\pfibd$ are in \cref{def:pfibd}.  They are subcategories of corresponding multimorphism categories of the $\Cat$-multicategory $\permcatsg$ (\cref{thm:permcatenrmulticat}).  We check that the multimorphism categories of $\pfibd$ are well defined in \cref{pfibdcomponentwelldef}.
\end{itemize}
The rest of the non-symmetric $\Cat$-multicategory structure of $\pfibd$ is discussed in \cref{sec:popfibmulticat}.

\subsection*{Opcartesian Multilinear Functors}

For fibrations 
\[P \cn \E \to \B \andspace P' \cn \E' \to \B,\] 
a \emph{Cartesian functor}\index{Cartesian functor}\index{functor!Cartesian} $F \cn P \to P'$ \cite[9.1.14 (1)]{johnson-yau} is a functor $F \cn \E \to \E'$ that
\begin{itemize}
\item satisfies $P = P'F$ and
\item sends Cartesian morphisms in $\E$ to Cartesian morphisms in $\E'$.
\end{itemize}
Recall that an \emph{$n$-linear functor} is a functor $F$ equipped with $n$ linearity constraints $\ang{F^2_j}_{j=1}^n$ and is subject to five axioms (\cref{def:nlinearfunctor}).  The multilinear analogs of Cartesian functors for permutative opfibrations with base category $\cD$ are multilinear functors that are compatible with the split opfibration structure as in the next definition.  Recall the notation $\compi$ in \cref{compnotation}.  Unless otherwise specified, in a permutative category we denote its monoidal product, monoidal unit, and braiding by, respectively, $\oplus$, $e$, and $\betaplus$.

\begin{definition}\label{def:opcartnlinearfunctor}
Suppose given permutative opfibrations (\cref{def:permutativefibration})
\[\begin{tikzcd}[column sep=large]
\C_k \ar{r}{P_k} & \cD
\end{tikzcd} \forspace k \in \{0,1,\ldots,n\}\]
with $\ang{P} = \ang{P_j}_{j=1}^n$ and $n \geq 0$.  An \emph{opcartesian $n$-linear functor}\index{opcartesian n-linear functor@opcartesian $n$-linear functor}\index{multilinear functor!opcartesian}
\begin{equation}\label{Fangcczero}
\begin{tikzcd}[column sep=2.3cm]
\ang{P} \ar{r}{\brb{F,\ang{F^2_j}_{j=1}^n}} & P_0
\end{tikzcd}
\end{equation}
is an $n$-linear functor \cref{nlinearfunctorangcd}
\[\begin{tikzcd}[column sep=2.3cm]
\prod_{j=1}^n \C_j \ar{r}{\brb{F,\ang{F^2_j}_{j=1}^n}} & \C_0
\end{tikzcd}\]
that satisfies the axioms \cref{def:opcartnlinearfunctor-i,def:opcartnlinearfunctor-ii,def:opcartnlinearfunctor-iii,def:opcartnlinearfunctor-iv} below.
\begin{enumerate}[label=(\roman*)]
\item\label{def:opcartnlinearfunctor-i} If $\ang{f_j} \in \txprod_{j=1}^n \C_j$ is an $n$-tuple of opcartesian morphisms, then $F\ang{f_j} \in \C_0$ is an opcartesian morphism.
\item\label{def:opcartnlinearfunctor-ii} The following diagram of functors is commutative, where $\otimes$ denotes any one of the equal iterates of $\otimes \cn \cD^2 \to \cD$.  The empty tensor product in $\cD$ is interpreted as the multiplicative unit $\tu \in \cD$.
\begin{equation}\label{def:opcartnlinearfunctor-diag}
\begin{tikzpicture}[xscale=2.5,yscale=1.3,vcenter]
\draw[0cell=.9]
(0,0) node (x11) {\txprod_{j=1}^n \C_j}
(x11)++(1,0) node (x12) {\C_0} 
(x11)++(0,-1) node (x21) {\cD^n}
(x12)++(0,-1) node (x22) {\cD}
;
\draw[1cell=.9]  
(x11) edge node {F} (x12)
(x21) edge node {\otimes} (x22)
(x11) edge node[swap] {\smallprod_j P_j} (x21)
(x12) edge node {P_0} (x22)
;
\end{tikzpicture}
\end{equation}
We call \cref{def:opcartnlinearfunctor-diag} the \index{projection axiom}\emph{projection axiom}.
\item\label{def:opcartnlinearfunctor-iii} Suppose that $\fbar_j \in \C_j$ is the chosen opcartesian lift of a fore-lift $\fl{x_j}{f_j}$ for $j \in \{1,\ldots,n\}$.  Then $F\ang{\fbar_j} \in \C_0$ is the chosen opcartesian lift of the fore-lift 
\[\bfl{F\ang{x_j}}{\txotimes_{j=1}^n f_j}.\]
\item\label{def:opcartnlinearfunctor-iv} 
Suppose given objects $\ang{x_j} \in \prod_{j=1}^n \C_j$ and $x_i' \in \C_i$ for some $i \in \{1,\ldots,n\}$, along with the following notation.
\begin{equation}\label{atensorpjxj}
\left\{\begin{split}
a &= P_1 x_1 \otimes \cdots \otimes P_i x_i \otimes \cdots \otimes P_n x_n\\
a' &= P_1 x_1 \otimes \cdots \otimes P_i x_i' \otimes \cdots \otimes P_n x_n\\
a'' &= P_1 x_1 \otimes \cdots \otimes (P_i x_i \oplus P_i x_i') \otimes \cdots \otimes P_n x_n
\end{split}\right.
\end{equation}
Denote by
\begin{equation}\label{laplazapx}
\begin{tikzcd}[column sep=large]
a'' \ar{r}{\lap}[swap]{\iso} & a \oplus a'
\end{tikzcd}\inspace \cD
\end{equation}
the unique Laplaza coherence isomorphism\index{Laplaza coherence isomorphism} (\cref{thm:laplaza-coherence-1}) that distributes over the sum $P_i x_i \oplus P_i x_i'$ in the $i$-th tensor factor in $a''$.  Then the component of the $i$-th linearity constraint
\begin{equation}\label{Ftwoixjxiprime}
\begin{tikzcd}[column sep=large]
F\ang{x_j} \oplus F\ang{x_j \compi x_i'} \ar{r}{F^2_i} & F\ang{x_j \compi (x_i \oplus x_i')}
\end{tikzcd}
\end{equation}
is the chosen opcartesian lift of the fore-lift
\begin{equation}\label{foreliftlapinverse}
\bfl{F\ang{x_j} \oplus F\ang{x_j \compi x_i'}}{\lap^\inv} \inspace \C_0.
\end{equation}
We call this the \index{constraint lift axiom}\emph{constraint lift axiom}.
\end{enumerate}
This finishes the definition of an opcartesian $n$-linear functor.  

Moreover, we define the following.
\begin{itemize}
\item An \emph{opcartesian strong $n$-linear functor}\index{opcartesian n-linear functor@opcartesian $n$-linear functor!strong} is an opcartesian $n$-linear functor $(F,\ang{F^2_j})$ in which each linearity constraint $F^2_j$ is a natural isomorphism.
\item An \emph{opcartesian (strong) multilinear functor}\index{opcartesian multilinear!functor} is an opcartesian (strong) $n$-linear functor for some $n \geq 0$.
\item If the projection functors $\ang{P_k}_{k=0}^n$ are unambiguously determined by the context, then we sometimes denote an opcartesian $n$-linear functor $\ang{P} \to P_0$ by $\angC \to \C_0$, where $\angC = \ang{\C_j}_{j=1}^n$.\defmark
\end{itemize}
\end{definition}

\begin{explanation}\label{expl:opcartnlinearfunctor}
Consider \cref{def:opcartnlinearfunctor}.
\begin{enumerate}
\item\label{expl:opcartnlinearfunctor-one}
An opcartesian $n$-linear functor is an $n$-linear functor with the extra properties \cref{def:opcartnlinearfunctor-i,def:opcartnlinearfunctor-ii,def:opcartnlinearfunctor-iii,def:opcartnlinearfunctor-iv} but \emph{no} extra structure.
\item\label{expl:opcartnlinearfunctor-0}
By definition a 0-linear functor is a functor $F \cn \boldone \to \C_0$, which is determined by the object $F* \in \C_0$.  In this case, axioms \cref{def:opcartnlinearfunctor-i,def:opcartnlinearfunctor-iii,def:opcartnlinearfunctor-iv} are vacuous, while the projection axiom \cref{def:opcartnlinearfunctor-ii} is the equality
\begin{equation}\label{pzerofstarone}
P_0 (F*) = \tu \inspace \cD.
\end{equation}
In other words, an \emph{opcartesian 0-linear functor}\index{opcartesian 0-linear!functor}\index{0-linear functor!opcartesian} is a choice of an object in the $P_0$-preimage of $\tu \in \cD$.
\item\label{expl:opcartnlinearfunctor-i} Axioms \cref{def:opcartnlinearfunctor-i,def:opcartnlinearfunctor-ii} are $n$-variable generalizations of the two axioms of a Cartesian functor, which we recalled before \cref{def:opcartnlinearfunctor}.  Axiom \cref{def:opcartnlinearfunctor-i} actually follows from axiom \cref{def:opcartnlinearfunctor-iii}; see \cref{preserveopcmorphisms}.
\item\label{expl:opcartnlinearfunctor-ii} Axiom \cref{def:opcartnlinearfunctor-iii} says that the functor $F$ preserves chosen opcartesian lifts.  It is well defined because
\[P_0 F\ang{\fbar_j} = \txotimes_{j=1}^n P_j \fbar_j = \txotimes_{j=1}^n f_j.\]
The first equality follows from the projection axiom \cref{def:opcartnlinearfunctor-diag}.  The second equality follows from the assumption that $\fbar_j$ is the chosen opcartesian lift of $\fl{x_j}{f_j}$.
\item\label{expl:opcartnlinearfunctor-iii} A part of the constraint lift axiom \cref{def:opcartnlinearfunctor-iv} is the morphism equality
\begin{equation}\label{pzeroftwoilapinv}
P_0 F^2_i = \lap^\inv \cn a \oplus a' \to a'' \inspace \cD
\end{equation} 
with $\lap$ in \cref{laplazapx} and $F^2_i$ in \cref{Ftwoixjxiprime}.  For the domain, respectively, codomain, in \cref{pzeroftwoilapinv} we use the following equalities.
\[\left\{\begin{aligned}
P_0 & \big(F\ang{x_j} \oplus F\ang{x_j \compi x_i'}\big) &&\\
&= P_0 F\ang{x_j} \oplus P_0 F\ang{x_j \compi x_i'} & \phantom{=} & \text{by \cref{strictsmfunctor}}\\
&= a \oplus a' & \phantom{=} & \text{by \cref{def:opcartnlinearfunctor-diag}}\\
P_0 & F \ang{x_j \compi (x_i \oplus x_i')} &&\\
&= P_1 x_1 \otimes \cdots \otimes P_i (x_i \oplus x_i') \otimes \cdots \otimes P_n x_n & \phantom{=} & \text{by \cref{def:opcartnlinearfunctor-diag}}\\
&= P_1 x_1 \otimes \cdots \otimes (P_i x_i \oplus P_i x_i') \otimes \cdots \otimes P_n x_n & \phantom{=} & \text{by \cref{strictsmfunctor}}\\
&= a'' &&
\end{aligned}\right.\]
The Laplaza coherence isomorphism $\lap$ in \cref{laplazapx} has appeared several times before, namely
\begin{itemize}
\item in \cref{laplazaapp} in the additivity axiom of an additive natural transformation,
\item in \cref{additivemodadditivity} in the additivity axiom of an additive modification, and
\item in \cref{grodphiconstraintdom,grodphilinearity} in the linearity constraints of $\grod\phi$ for an additive natural transformation $\phi$.\defmark
\end{itemize}
\end{enumerate}
\end{explanation}

While \cref{def:opcartnlinearfunctor} \cref{def:opcartnlinearfunctor-i} is the natural $n$-variable variant of the preservation of Cartesian morphisms by Cartesian functors, the following observation says that it is actually redundant.

\begin{lemma}\label{preserveopcmorphisms}
In the context of \cref{def:opcartnlinearfunctor}, axiom \cref{def:opcartnlinearfunctor-i} follows from axiom \cref{def:opcartnlinearfunctor-iii}.
\end{lemma}

\begin{proof}
Suppose $f_j \cn x_j \to y_j$ in $\C_j$ is an opcartesian morphism for each $j \in \{1,\ldots,n\}$.  We want to show that $F\ang{f_j}$ is an opcartesian morphism as a consequence of axiom \cref{def:opcartnlinearfunctor-iii}.  We use the following facts from \cite[9.1.4 (1), (4), and (5)]{johnson-yau} for a functor $P \cn \E \to \B$:
\begin{itemize}
\item Each isomorphism in $\E$ is an opcartesian morphism.
\item Opcartesian morphisms are closed under composition.
\item If both $\fbar \cn y \to z$ and $\fbar' \cn y \to z'$ are opcartesian lifts of a fore-lift $\fl{y}{f}$, then there exists a unique isomorphism $g \cn z \fto{\iso} z'$ in $\E$ such that 
\[g\fbar = \fbar' \andspace Pg = 1_{Pz}.\]
\end{itemize}

Consider the fore-lift $\fl{x_j}{P_j f_j}$ with respect to $P_j \cn \C_j \to \cD$.  
\begin{itemize}
\item Since $f_j$ is an opcartesian morphism, it is an opcartesian lift of $\fl{x_j}{P_j f_j}$.
\item Since $P_j$ is a split opfibration, $\fl{x_j}{P_j f_j}$ has a \emph{chosen} opcartesian lift $\fbar_j \cn x_j \to w_j$.  
\end{itemize}
Therefore, there is a unique isomorphism $g_j \cn w_j \fto{\iso} y_j$ in $\C_j$ such that
\[g_j \fbar_j = f_j \andspace P_j g_j = 1_{P_j y_j}.\]
The functoriality of $F$ implies
\begin{equation}\label{Fangfjfactor}
F\ang{f_j} = F\ang{g_j \fbar_j} = F\ang{g_j} \circ F\ang{\fbar_j} \inspace \C_0.
\end{equation}
Now we observe the following.
\begin{itemize}
\item By axiom \cref{def:opcartnlinearfunctor-iii} $F\ang{\fbar_j}$ is the chosen opcartesian lift of $\bfl{F\ang{x_j}}{\otimes_j P_j f_j}$. 
\item Since each $g_j$ is an isomorphism, so is $F\ang{g_j}$ by the functoriality of $F$.  Thus $F\ang{g_j}$ is an opcartesian morphism.
\end{itemize}
The factorization \cref{Fangfjfactor} shows that $F\ang{f_j}$ is the composite of two opcartesian morphisms, so it is also an opcartesian morphism.
\end{proof}

In view of \cref{preserveopcmorphisms}, to check that an $n$-linear functor is opcartesian, it is unnecessary to check axiom \cref{def:opcartnlinearfunctor-i}.  We include axiom \cref{def:opcartnlinearfunctor-i} in \cref{def:opcartnlinearfunctor} solely to preserve the symmetry with the definition of a Cartesian functor.

\begin{example}[Identity Opcartesian 1-Linear Functors]\label{ex:idomf}
Suppose $P \cn \C \to \cD$ is a permutative opfibration (\cref{def:permutativefibration}).  Then the identity 1-linear functor 
\[\begin{tikzcd}[column sep=large]
\C \ar{r}{1_{\C}} & \C
\end{tikzcd}\]
in \cref{ex:idonelinearfunctor} is an opcartesian 1-linear functor.  Indeed, axioms \cref{def:opcartnlinearfunctor-i,def:opcartnlinearfunctor-ii,def:opcartnlinearfunctor-iii} in \cref{def:opcartnlinearfunctor} follow from the fact that $1_{\C}$ is the identity functor on $\C$.  

For the constraint lift axiom \cref{def:opcartnlinearfunctor-iv}, the Laplaza coherence isomorphism \cref{laplazapx}
\[\begin{tikzcd}[column sep=large]
Px \oplus Px' \ar{r}{\lap}[swap]{\iso} & Px \oplus Px' \forspace x,x' \in \C
\end{tikzcd}\]
is the identity morphism by uniqueness.  So $\lap^\inv$ is also the identity morphism.  The only linearity constraint of $1_{\C}$ is the identity natural transformation by definition.  The chosen opcartesian lift of the fore-lift \cref{foreliftlapinverse}
\[\bfl{x \oplus x'}{1_{Px \oplus Px'}} = \bfl{x \oplus x'}{1_{P(x \oplus x')}}\]
is the identity morphism $1_{x \oplus x'}$ by the unitarity condition of a split opfibration (\cref{def:opfibration} \cref{def:opfibration-viii}), including $P$.
\end{example}

\begin{example}[Opcartesian 1-Linear Functors]\label{ex:idoponelinear}\index{opcartesian 1-linear functor}\index{1-linear functor!opcartesian}
Generalizing \cref{ex:idomf}, for permutative opfibrations
\[\begin{tikzcd}[column sep=large]
\C_1 \ar{r}{P_1} & \cD
\end{tikzcd} \andspace
\begin{tikzcd}[column sep=large]
\C_0 \ar{r}{P_0} & \cD,
\end{tikzcd}\]
an opcartesian 1-linear functor
\[\begin{tikzcd}[column sep=large]
P_1 \ar{r}{F} & P_0
\end{tikzcd}\]
is precisely a \emph{strict} symmetric monoidal functor
\[\begin{tikzcd}[column sep=large]
\C_1 \ar{r}{F} & \C_0
\end{tikzcd}\]
that
\begin{itemize}
\item satisfies $P_1 = P_0 F$ and
\item preserves chosen opcartesian lifts.
\end{itemize}
Indeed, the conditions in the previous two items are axioms \cref{def:opcartnlinearfunctor-ii,def:opcartnlinearfunctor-iii} in \cref{def:opcartnlinearfunctor}.  Axiom \cref{def:opcartnlinearfunctor-i} holds by \cref{preserveopcmorphisms}.  A 1-linear functor (\cref{def:nlinearfunctor})
\[\begin{tikzcd}[column sep=large]
\C_1 \ar{r}{(F,F^2)} & \C_0
\end{tikzcd}\]
is precisely a strictly unital symmetric monoidal functor.  Moreover, since $n=1$, the Laplaza coherence isomorphism $\lap$ in \cref{laplazapx} is the identity morphism.  Then
\begin{itemize}
\item axiom \cref{def:opcartnlinearfunctor-iv} in \cref{def:opcartnlinearfunctor} and
\item the unitarity condition (\cref{def:opfibration} \cref{def:opfibration-viii}) of the split opfibration $P_0$
\end{itemize}  
imply that the only linearity constraint $F^2$ is the identity.  Thus $F$ is, in fact, a strict symmetric monoidal functor.
\end{example}

\begin{example}[Bipermutative Categories]\label{ex:bipermopnlinear}
Since $(\cD,\oplus,\otimes)$ is a small tight bipermutative category (\cref{def:embipermutativecat}), the identity functor
\[\begin{tikzcd}[column sep=large]
\Dplus \ar{r}{1_{\Dplus}} & \Dplus
\end{tikzcd}\]
is a permutative opfibration (\cref{ex:idpermopfibration}).  For each $n \geq 0$ there is an opcartesian strong $n$-linear functor
\begin{equation}\label{onedntooned}
\begin{tikzcd}[column sep=2.3cm]
\ang{1_{\Dplus}}_{j=1}^n \ar{r}{\brb{\otimes,\ang{\lap_j^\inv}_{j=1}^n}} & 1_{\Dplus}.
\end{tikzcd}
\end{equation}
Indeed, by \cref{ex:bipermnlinearfunctor} there is a strong $n$-linear functor
\[\begin{tikzcd}[column sep=2.3cm]
\txprod_{j=1}^n \Dplus \ar{r}{\brb{\otimes,\ang{\lap_j^\inv}_{j=1}^n}} & \Dplus.
\end{tikzcd}\]
Its $j$-th linearity constraint is \cref{laplinearityconst}
\[\otimes^2_j = \lap_j^\inv\] 
where $\lap_j$ is the unique Laplaza coherence isomorphism (\cref{ex:laplazacoherence}) that distributes the sum in the $j$-th tensor factor. 
\begin{itemize}
\item Axioms \cref{def:opcartnlinearfunctor-i,def:opcartnlinearfunctor-ii,def:opcartnlinearfunctor-iii} in \cref{def:opcartnlinearfunctor} hold because each $1_{\Dplus}$ is the identity functor.
\item The constraint lift axiom \cref{def:opcartnlinearfunctor-iv} holds because $\otimes^2_i = \lap_i^\inv$.\defmark
\end{itemize}
\end{example}

\subsection*{Opcartesian Multilinear Transformations}

For
\begin{itemize}
\item fibrations $P \cn \E \to \B$ and $P' \cn \E' \to \B$ and
\item Cartesian functors $F,G \cn P \to P'$,
\end{itemize}  
a \emph{vertical natural transformation}\index{vertical natural transformation}\index{natural transformation!vertical} $\theta \cn F \to G$ \cite[9.1.14 (2)]{johnson-yau} is a natural transformation that satisfies
\[1_{P'} * \theta = 1_P.\]
Recall that an \emph{$n$-linear natural transformation} is a natural transformation that satisfies the unity axiom and the constraint compatibility axiom (\cref{def:nlineartransformation}).  The multilinear analogs of vertical natural transformations are compatible with the functors $\ang{P_k}_{k=0}^n$ in the following sense.

\begin{definition}\label{def:opcartnlineartr}
Suppose given permutative opfibrations (\cref{def:permutativefibration})
\[\begin{tikzcd}[column sep=large]
\C_k \ar{r}{P_k} & \cD
\end{tikzcd} \forspace k \in \{0,1,\ldots,n\}\]
with $\ang{P} = \ang{P_j}_{j=1}^n$ and $n \geq 0$.  Suppose given two opcartesian $n$-linear functors (\cref{def:opcartnlinearfunctor})
\[\begin{tikzcd}[column sep=2.3cm]
\ang{P} \ar[shift left]{r}{\brb{F,\ang{F^2_j}_{j=1}^n}} 
\ar[shift right]{r}[swap]{\brb{G,\ang{G^2_j}_{j=1}^n}} & P_0.
\end{tikzcd}\]
An \emph{opcartesian $n$-linear transformation}\index{opcartesian n-linear transformation@opcartesian $n$-linear transformation}\index{n-linear natural transformation@$n$-linear natural transformation!opcartesian}
\begin{equation}\label{thetafftwoggtwo}
\begin{tikzcd}[column sep=large]
\brb{F,\ang{F^2_j}_{j=1}^n} \ar{r}{\theta} & \brb{G,\ang{G^2_j}_{j=1}^n}
\end{tikzcd}
\end{equation}
is an $n$-linear natural transformation \cref{nlineartransformation}
\[\begin{tikzpicture}[xscale=2,yscale=2,baseline={(a.base)}]
\draw[0cell=.9]
(0,0) node (a) {\textstyle\prod_{j=1}^n \C_j}
(a)++(.16,0) node (b) {\phantom{D}}
(b)++(1,0) node (c) {\C_0}
;
\draw[1cell=.8]  
(b) edge[bend left] node[pos=.5] {\big( F, \ang{F^2_j}\big)} (c)
(b) edge[bend right] node[swap,pos=.5] {\big( G, \ang{G^2_j}\big)} (c)
;
\draw[2cell] 
node[between=b and c at .45, rotate=-90, 2label={above,\theta}] {\Rightarrow}
;
\end{tikzpicture}\]
such that the following two whiskered natural transformations are equal.
\begin{equation}\label{thetapzeropjtensor}
\begin{tikzpicture}[xscale=2,yscale=2,vcenter]
\def\v{-.65}
\draw[0cell=.9]
(0,0) node (a) {\txprod_{j=1}^n \C_j}
(a)++(.16,0) node (b) {\phantom{D}}
(b)++(1,0) node (c) {\C_0}
(a)++(0,\v) node (dn) {\cD^n}
(c)++(0,\v) node (d) {\cD}
;
\draw[1cell=.8]  
(b) edge[bend left] node[pos=.5] {F} (c)
(b) edge[bend right] node[swap,pos=.5] {G} (c)
(a) edge node[swap] {\smallprod_j P_j} (dn)
(c) edge node {P_0} (d)
(dn) edge node {\otimes} (d)
;
\draw[2cell] 
node[between=b and c at .45, rotate=-90, 2label={above,\theta}] {\Rightarrow}
;
\end{tikzpicture}
\end{equation}
We call \cref{thetapzeropjtensor} the \index{transformation projection axiom}\emph{transformation projection axiom}.  An \emph{opcartesian multilinear transformation}\index{opcartesian multilinear!transformation}\index{multilinear natural transformation!opcartesian} is an opcartesian $n$-linear transformation for some $n \geq 0$.
\end{definition}

\begin{explanation}\label{expl:opcartnlineartr}
Consider \cref{def:opcartnlineartr}.
\begin{enumerate}
\item\label{expl:opcartnlineartr-i}
An opcartesian $n$-linear transformation is an $n$-linear natural transformation with the extra property \cref{thetapzeropjtensor} but \emph{no} extra structure.
\item\label{expl:opcartnlineartr-ii} 
By definition a 0-linear natural transformation $\theta \cn F \to G$ is a morphism 
\begin{equation}\label{thetafgstar}
\begin{tikzcd}[column sep=large]
F* \ar{r}{\theta_*} & G*
\end{tikzcd} \inspace \C_0.
\end{equation}
Using \cref{pzerofstarone}, in this case the transformation projection axiom \cref{thetapzeropjtensor} is the morphism equality
\begin{equation}\label{pzerothetastarone}
P_0 \theta_* = 1_{\tu} \cn \tu \to \tu \inspace \cD.
\end{equation}
In other words, an \emph{opcartesian 0-linear transformation}\index{opcartesian 0-linear!transformation} is a morphism as in \cref{thetafgstar} that is in the $P_0$-preimage of $\tu \in \cD$.
\item\label{expl:opcartnlineartr-iii} 
For $n \geq 1$ the three composite functors in the transformation projection axiom \cref{thetapzeropjtensor} are equal by the projection axiom \cref{def:opcartnlinearfunctor-diag} for $F$ and $G$.  The axiom \cref{thetapzeropjtensor} means that, for each $n$-tuple of objects $\ang{x_j} \in \prod_{j=1}^n \C_j$, there is a morphism equality
\begin{equation}\label{pzerothetaangxjone}
P_0 \theta_{\ang{x_j}} = 1_{\otimes_{j=1}^n P_j x_j} \cn P_0 F\ang{x_j} \to P_0 G\ang{x_j}
\end{equation}
in $\cD$.\defmark
\end{enumerate}
\end{explanation}

\subsection*{Multimorphism Categories}

For small permutative categories $\C_0$ and $\angC = \ang{\C_j}_{j=1}^n$, recall from \cref{permcatsgangcd,permcatemptyd} the category
\[\permcatsg\scmap{\angC;\C_0} = \permcatsg\scmap{\C_1,\ldots,\C_n;\C_0}\]
with
\begin{itemize}
\item \emph{strong} $n$-linear functors $\txprod_{j=1}^n \C_j \to \C_0$ as objects,
\item $n$-linear natural transformations as morphisms,
\item identity natural transformations as identity morphisms, and
\item vertical composition of natural transformations as composition.
\end{itemize} 
In the arity 0 case, $\permcatsg\smscmap{\ang{};\C_0}$ is the category $\C_0$.

\begin{definition}\label{def:pfibd}
Suppose given small permutative opfibrations (\cref{def:permutativefibration})
\[\begin{tikzcd}[column sep=large]
\C_k \ar{r}{P_k} & \cD
\end{tikzcd} \forspace k \in \{0,1,\ldots,n\}\]
with $\angC = \ang{\C_j}_{j=1}^n$, $\ang{P} = \ang{P_j}_{j=1}^n$, and $n \geq 0$.  Define\label{not:pfibdangppzero}\index{permutative opfibration!multimorphism category}\index{multimorphism category!permutative opfibration}
\[\pfibd\scmap{\ang{P};P_0} = \pfibd\scmap{P_1,\ldots,P_n;P_0}\]
as the subcategory of $\permcatsg\scmap{\angC;\C_0}$ with
\begin{itemize}
\item opcartesian \emph{strong} $n$-linear functors (\cref{def:opcartnlinearfunctor}) as objects and
\item opcartesian $n$-linear transformations (\cref{def:opcartnlineartr}) as morphisms.
\end{itemize}
This finishes the definition of $\pfibd\scmap{\ang{P};P_0}$. 

Moreover, if the projection functors $\ang{P_k}_{k=0}^n$ are unambiguously determined by the context, then we also denote $\pfibd\scmap{\ang{P};P_0}$ by $\pfibd\scmap{\angC;\C_0}$.
\end{definition}

First we check that the category $\pfibd\scmap{\ang{P};P_0}$ is well defined.

\begin{lemma}\label{pfibdcomponentwelldef}
In the context of \cref{def:pfibd}, 
\[\pfibd\scmap{\ang{P};P_0}\]
is a well-defined subcategory of $\permcatsg\scmap{\angC;\C_0}$.
\end{lemma}

\begin{proof}
In the arity 0 case, the assertion is that $\pfibd\smscmap{\ang{};P_0}$ is a subcategory of $\C_0$.  This is true because $\pfibd\smscmap{\ang{};P_0}$ is the $P_0$-preimage of $\tu \in \cD$ by
\begin{itemize}
\item \cref{expl:opcartnlinearfunctor} \eqref{expl:opcartnlinearfunctor-0} for objects and
\item \cref{expl:opcartnlineartr} \eqref{expl:opcartnlineartr-ii} for morphisms.
\end{itemize} 

Suppose $n > 0$.  The identity natural transformation $1_F$ of each opcartesian $n$-linear functor $F$ is an opcartesian $n$-linear transformation.  Indeed, in this case the transformation projection axiom \cref{thetapzeropjtensor} reduces to the projection axiom \cref{def:opcartnlinearfunctor-diag} for $F$.  

It remains to check that opcartesian $n$-linear transformations are closed under vertical composition.  Suppose
\[\begin{tikzcd}[column sep=large]
\brb{F,\ang{F^2_j}_{j=1}^n} \ar{r}{\theta} 
& \brb{G,\ang{G^2_j}_{j=1}^n} \ar{r}{\psi} & \brb{H,\ang{H^2_j}_{j=1}^n}
\end{tikzcd}\]
are opcartesian $n$-linear transformations in $\pfibd\scmap{\ang{P};P_0}$.  The vertical composite $n$-linear natural transformation
\[\begin{tikzcd}[column sep=large]
\brb{F,\ang{F^2_j}_{j=1}^n} \ar{r}{\psi\theta} & \brb{H,\ang{H^2_j}_{j=1}^n}
\end{tikzcd}\]
satisfies the transformation projection axiom \cref{pzerothetaangxjone} by the following computation for each $n$-tuple of objects $\ang{x_j} \in \prod_{j=1}^n \C_j$.
\[\begin{aligned}
P_0 (\psi\theta)_{\ang{x_j}} 
&= P_0 \big( \psi_{\ang{x_j}} \theta_{\ang{x_j}} \big) &&\\
&= \big( P_0 \psi_{\ang{x_j}} \big) \big( P_0 \theta_{\ang{x_j}} \big) & \phantom{=} & \text{by functoriality of $P_0$}\\
&= \big(1_{\otimes_{j=1}^n P_j x_j}\big) \big(1_{\otimes_{j=1}^n P_j x_j}\big) & \phantom{=} & \text{by \cref{pzerothetaangxjone} for $\theta$ and $\psi$}\\
&= 1_{\otimes_{j=1}^n P_j x_j} &&
\end{aligned}\]
This finishes the proof.
\end{proof}

\section{The Multicategory of Permutative Opfibrations}
\label{sec:popfibmulticat}

As in \cref{Dstbipermutative} we continue to assume that $(\cD,\oplus,\otimes)$ is a small tight bipermutative category.  In \cref{def:pfibd} we constructed the multimorphism categories of $\pfibd$.  Recall from \cref{def:enr-multicategory} that a non-symmetric $\V$-multicategory is a $\V$-multicategory without the symmetric group action and the corresponding axioms.  Beyond the objects and the multimorphism categories, a non-symmetric $\V$-multicategory is determined by its colored units and composition, subject to
\begin{itemize}
\item the associativity axiom \cref{enr-multicategory-associativity} and
\item the unity axioms \cref{enr-multicategory-right-unity,enr-multicategory-left-unity}.
\end{itemize} 

In this section we construct the rest of the non-symmetric $\Cat$-multicategory structure on $\pfibd$.  In the rare situation that the multiplicative braiding in $\cD$ is the identity, $\pfibd$ is a $\Cat$-multicategory.  This (non-symmetric) $\Cat$-multicategory comes with a (non-symmetric) $\Cat$-multifunctor
\[\begin{tikzcd}[column sep=large]
\pfibd \ar{r}{\U} & \permcatsg
\end{tikzcd}\]
that sends each small permutative opfibration over $\cD$ to its total category.  Here is an outline of this section.
\begin{itemize}
\item \cref{def:pfibdmulticat} specifies the colored units, composition, and symmetric group action in $\pfibd$ by restricting those in the $\Cat$-multicategory $\permcatsg$ (\cref{thm:permcatenrmulticat}).
\item \cref{pfibdcompobjwelldef,pfibdcompmorwelldef} show that the composition in $\permcatsg$ restricts to $\pfibd$.
\item \cref{pfibdsymwelldefobj,pfibdsymwelldefmor} show that the symmetric group action in $\permcatsg$ restricts to $\pfibd$ \emph{if} the multiplicative braiding $\betate$ in $\cD$ is the identity.  The last assumption is used in the diagrams \cref{pfibdnonsymmetrydiag,pfibdnonsymmetrymordiag}.  A weaker converse holds: if the symmetric group action in $\permcatsg$ restricts to $\pfibd$, then $\otimes = \tensorop$ in $\cD$.
\item The main \cref{thm:pfibdmulticat} proves the existence of the non-symmetric $\Cat$-multicategory $\pfibd$, which is, furthermore, a $\Cat$-multicategory if $\betate = 1$ in $\cD$.
\item \cref{cor:pfibdtopermcatsg} states that there exists a non-symmetric $\Cat$-multifunctor $\U$ from $\pfibd$ to $\permcatsg$ that sends each small permutative opfibration over $\cD$ to the total category.  Moreover, $\U$ is a $\Cat$-multifunctor if $\betate = 1$ in $\cD$.
\end{itemize}

\begin{definition}\label{def:pfibdmulticat}
For a small tight bipermutative category $\cD$ as in \cref{Dstbipermutative}, define the data of a $\Cat$-multicategory $\pfibd$ as follows.
\begin{description}
\item[Objects] An object in $\pfibd$ is a small permutative opfibration over $\cD$ (\cref{def:permutativefibration}).
\item[Multimorphisms] The multimorphism categories in $\pfibd$ are the ones in \cref{def:pfibd}:
\begin{itemize}
\item The objects are opcartesian strong multilinear functors (\cref{def:opcartnlinearfunctor}).
\item The morphisms are opcartesian multilinear transformations (\cref{def:opcartnlineartr}).
\item Identity morphisms are identity natural transformations.
\item Composition is vertical composition of natural transformations.
\end{itemize}
\item[Units] For each small permutative opfibration $P \cn \C \to \cD$, the $P$-colored unit is the identity opcartesian 1-linear functor $1_{\C} \cn \C \to \C$ (\cref{ex:idomf}).
\item[Composition] It is the restriction of the composition in $\permcatsg$ \cref{permcatgamma} via the subcategory inclusions in \cref{pfibdcomponentwelldef}.
\item[Symmetric Group Action] It is the restriction of the symmetric group action in $\permcatsg$ \cref{permcatsymgroupaction} via the subcategory inclusions in \cref{pfibdcomponentwelldef}.
\end{description}
This finishes the definition of $\pfibd$.
\end{definition}

\subsection*{Composition}

To show that the composition in $\pfibd$ is well defined, we begin with $n$-ary 1-cells.

\begin{lemma}\label{pfibdcompobjwelldef}
In the context of \cref{def:pfibdmulticat}, the composition in $\pfibd$ is well defined on opcartesian strong multilinear functors.
\end{lemma}

\begin{proof}
Suppose given
\begin{itemize}
\item small permutative opfibrations over $\cD$ (\cref{def:permutativefibration})
\begin{itemize}
\item $P_k \cn \C_k \to \cD$ for $k \in \{0,1,\ldots,n\}$ and
\item $Q_{ji} \cn \B_{ji} \to \cD$ for $j \in \{1,\ldots,n\}$ and $i \in \{1,\ldots,\ell_j\}$,
\end{itemize}
\item an opcartesian $n$-linear functor (\cref{def:opcartnlinearfunctor})
\begin{equation}\label{angpfftwojpzero}
\begin{tikzcd}[column sep=2.3cm]
\ang{P} \ar{r}{\brb{F,\ang{F^2_j}_{j=1}^n}} & P_0
\end{tikzcd}
\end{equation}
with $\ang{P} = \ang{P_j}_{j=1}^n$, and
\item for each $j \in \{1,\ldots,n\}$ an opcartesian $\ell_j$-linear functor 
\begin{equation}\label{angqjgjpj}
\begin{tikzcd}[column sep=2.5cm]
\ang{Q_j} \ar{r}{\brb{G_j,\bang{(G_j)^2_i}_{i=1}^{\ell_j}}} & P_j
\end{tikzcd}
\end{equation}
with $\ang{Q_j} = \ang{Q_{ji}}_{i=1}^{\ell_j}$.
\end{itemize}
With $\ell = \ell_1 + \cdots + \ell_n$, the composite $\ell$-linear functor
\begin{equation}\label{prodgjf}
\begin{tikzcd}[column sep=large]
\txprod_{j=1}^n \txprod_{i=1}^{\ell_j} \B_{ji} \ar{r}{\smallprod_j G_j} 
& \txprod_{j=1}^n \C_j \ar{r}{F} & [-5pt] \C_0
\end{tikzcd}
\end{equation}
has linearity constraints in \cref{ffjlinearity}.  It is strong if $F$ and each $G_j$ are strong.  We must show that the composite $\ell$-linear functor in \cref{prodgjf} satisfies axioms \cref{def:opcartnlinearfunctor-i,def:opcartnlinearfunctor-ii,def:opcartnlinearfunctor-iii,def:opcartnlinearfunctor-iv} in \cref{def:opcartnlinearfunctor}.  

If some $\ell_j = 0$, then the argument below needs a slight adjustment using
\begin{itemize}
\item \cref{expl:opcartnlinearfunctor} \eqref{expl:opcartnlinearfunctor-0}, which identifies opcartesian 0-linear functors with objects in the preimage of $\tu \in \cD$;
\item the fact that $\tu$ is the strict multiplicative unit in $\cD$; and
\item the unitarity condition of a split opfibration (\cref{def:opfibration} \cref{def:opfibration-viii}).
\end{itemize} 
With this in mind, in the rest of this proof we assume that each $\ell_j > 0$.

\emph{Axiom \cref{def:opcartnlinearfunctor-i}}.  This holds because $F$ and each $G_j$ preserve opcartesian morphisms.

\emph{Axiom \cref{def:opcartnlinearfunctor-ii}}.  The diagram \cref{def:opcartnlinearfunctor-diag} for the composite $F \circ (\prod_j G_j)$ in \cref{prodgjf} is the boundary of the following diagram.
\[\begin{tikzpicture}[xscale=2.2,yscale=1.3,vcenter]
\draw[0cell=.9]
(0,0) node (x11) {\txprod_{j=1}^n \C_j}
(x11)++(1,0) node (x12) {\C_0} 
(x11)++(-1.4,0) node (x10) {\txprod_{j=1}^n \txprod_{i=1}^{\ell_j} \B_{ji}}
(x10)++(0,-1) node (x20) {\cD^{\ell_1+\cdots+\ell_n}}
(x11)++(0,-1) node (x21) {\cD^n}
(x12)++(0,-1) node (x22) {\cD}
;
\draw[1cell=.9]  
(x10) edge node {\smallprod_j G_j} (x11)
(x11) edge node {F} (x12)
(x20) edge node {\smallprod_j \otimes} (x21)
(x21) edge node {\otimes} (x22)
(x10) edge node[swap] {\smallprod_{j,i}\, Q_{ji}} (x20)
(x11) edge node[swap] {\smallprod_j P_j} (x21)
(x12) edge node {P_0} (x22)
;
\end{tikzpicture}\]
The left and right squares commute by the projection axiom \cref{def:opcartnlinearfunctor-diag} for, respectively, $\ang{G_j}_{j=1}^n$ and $F$.

\emph{Axiom \cref{def:opcartnlinearfunctor-iii}}.  This holds because $F$ and each $G_j$ preserve chosen opcartesian lifts.

\emph{Axiom \cref{def:opcartnlinearfunctor-iv}}.  By definition \cref{ffjlinearity} for each $q \in \{1,\ldots,n\}$ and $p \in \{1,\ldots,\ell_q\}$ with 
\[r = \ell_1 + \cdots + \ell_{q-1} + p,\]
the $r$-th linearity constraint of $(F \circ \prod_j G_j)$ in \cref{prodgjf} is the composite
\begin{equation}\label{fprodjgjtwor}
(F \circ \txprod_j G_j)^2_r = F \bang{1 \circ_q (G_q)^2_p} \circ F^2_q.
\end{equation}
Here
\begin{itemize}
\item $F^2_q$ is the $q$-th linearity constraint of $F$, and
\item $(G_q)^2_p$ is the $p$-th linearity constraint of $G_q$.
\end{itemize}
Axiom \cref{def:opcartnlinearfunctor-iv} requires each component of $(F \circ \prod_j G_j)^2_r$ to be the chosen opcartesian lift of $\lap^\inv_r$, where $\lap_r$ is the unique Laplaza coherence isomorphism in $\cD$ that distributes over the sum in the $r$-th tensor factor in its domain.  This axiom holds for the following reasons.
\begin{enumerate}
\item\label{fprodjgj-i} By axiom \cref{def:opcartnlinearfunctor-iv} for $F$, each component of $F^2_q$ is the chosen opcartesian lift of $\lap_q^\inv$, where $\lap_q$ is the unique Laplaza coherence isomorphism in $\cD$ that distributes over the sum in the $q$-th tensor factor in its domain.
\item\label{fprodjgj-ii} By axiom \cref{def:opcartnlinearfunctor-iv} for $G_q$, each component of $(G_q)^2_p$ is the chosen opcartesian lift of $\lap_p^\inv$, where $\lap_p$ is the unique Laplaza coherence isomorphism in $\cD$ that distributes over the sum in the $p$-th tensor factor in its domain.
\item\label{fprodjgj-iii} By the unitarity axiom of a split opfibration (\cref{def:opfibration} \cref{def:opfibration-viii}), each component of each entry in the $n$-tuple
\[\bang{1 \circ_q (G_q)^2_p} = \bang{1, \ldots, (G_q)^2_p, \ldots, 1}\]
is a chosen opcartesian lift.  By
\begin{itemize}
\item \eqref{fprodjgj-ii} above,
\item the functoriality of each $P_j$, and
\item axiom \cref{def:opcartnlinearfunctor-iii} for $F$,
\end{itemize}  
each component of
\[F \bang{1 \circ_q (G_q)^2_p}\]
is the chosen opcartesian lift of 
\[1 \otimes \cdots \otimes \lap_p^\inv \otimes \cdots \otimes 1.\]
\item\label{fprodjgj-iv} By the uniqueness in \cref{thm:laplaza-coherence-1}, there is a morphism equality in $\cD$
\[\lap_q \circ \big(1 \otimes \cdots \otimes \lap_p \otimes \cdots \otimes 1) = \lap_r.\]
\item\label{fprodjgj-v} By the multiplicativity axiom of a split opfibration (\cref{def:opfibration} \cref{def:opfibration-viii}) and \eqref{fprodjgj-i}, \eqref{fprodjgj-iii}, and \eqref{fprodjgj-iv} above, each component of the right-hand side of \cref{fprodjgjtwor} is the chosen opcartesian lift of 
\[\big(1 \otimes \cdots \otimes \lap_p^\inv \otimes \cdots \otimes 1\big) \circ \lap_q^\inv = \lap_r^\inv.\]
\end{enumerate}
This finishes the proof.
\end{proof}

Next we consider composition of $n$-ary 2-cells.

\begin{lemma}\label{pfibdcompmorwelldef}
In the context of \cref{def:pfibdmulticat}, the composition in $\pfibd$ is well defined on opcartesian multilinear transformations.
\end{lemma}

\begin{proof}
Reusing the notation in \cref{angpfftwojpzero,angqjgjpj,prodgjf}, suppose given
\begin{itemize}
\item an opcartesian $n$-linear transformation 
\[\begin{tikzcd}[column sep=large]
F \ar{r}{\theta} & F'
\end{tikzcd}\]
in $\pfibd\scmap{\ang{P};P_0}$ and
\item for each $j \in \{1,\ldots,n\}$ an opcartesian $\ell_j$-linear transformation 
\[\begin{tikzcd}[column sep=large]
G_j \ar{r}{\theta_j} & G_j'
\end{tikzcd}\]
in $\pfibd\scmap{\ang{Q_j};P_j}$.
\end{itemize} 
By the definition \cref{permcatcomposite} of the composition, the transformation projection axiom \cref{thetapzeropjtensor} for the composite $\ell$-linear natural transformation $\ga\scmap{\theta;\ang{\theta_j}}$ states that the two boundary whiskered natural transformations in the following diagram are equal.
\[\begin{tikzpicture}[xscale=2,yscale=2,baseline={(a.base)}]
\def\v{-.7} \def\d{25}
\draw[0cell=.9]
(0,0) node (a) {\textstyle\prod_{j=1}^n \C_j}
(a)++(.16,0) node (b) {\phantom{D}}
(b)++(1,0) node (c) {\C_0}
(a)++(-.16,0) node (z) {\phantom{D}}
(z)++(-1.3,0) node (y) {\phantom{D}}
(y)++(-.33,0) node (x) {\txprod_{j=1}^n \txprod_{i=1}^{\ell_j} \B_{ji}}
(x)++(0,\v) node (w1) {\cD^{\ell_1+\cdots+\ell_n}}
(a)++(0,\v) node (w2) {\cD^n}
(c)++(0,\v) node (w3) {\cD} 
;
\draw[1cell=.8]  
(b) edge[bend left=\d] node[pos=.5] {F} (c)
(b) edge[bend right=\d] node[swap,pos=.5] {F'} (c)
(y) edge[bend left=\d] node[pos=.5] {\smallprod_j G_j} (z)
(y) edge[bend right=\d] node[swap,pos=.5] {\smallprod_j G_j'} (z)
(x) edge node[swap] {\smallprod_{j,i}\, Q_{ji}} (w1)
(a) edge node[swap] {\smallprod_j P_j} (w2)
(c) edge node {P_0} (w3)
(w1) edge node[pos=.3] {\smallprod_j\, \otimes} (w2)
(w2) edge node {\otimes} (w3)
;
\draw[2cell] 
node[between=b and c at .45, rotate=-90, 2label={above,\theta}] {\Rightarrow}
node[between=y and z at .38, rotate=-90, 2label={above,\smallprod_j \theta_j}] {\Rightarrow}
;
\end{tikzpicture}\]
The two boundary whiskered natural transformations above are equal by the transformation projection axiom \cref{thetapzeropjtensor} for $\ang{\theta_j}$ and $\theta$.
\end{proof}

\subsection*{Symmetry}

Next we show that the symmetric group action in $\pfibd$, which is a restriction of the one in $\permcatsg$, is well defined when the multiplicative braiding in $\cD$ is the identity, together with a weaker converse statement.  For a monoidal category $(\C,\otimes)$, denote by $\tensorop$, which is called the \emph{opposite monoidal product}, the composite functor
\begin{equation}\label{oppositetensor}
\begin{tikzcd}[column sep=large]
\C^2 \ar{r}[pos=.6]{\xi} \ar[bend left]{rr}{\tensorop} & \C^2 \ar{r}[pos=.4]{\otimes} & \C
\end{tikzcd}
\end{equation}
where $\xi$ swaps the two arguments.

\begin{lemma}\label{pfibdsymwelldefobj}
In the context of \cref{def:pfibdmulticat}, the following statements hold.
\begin{enumerate}
\item\label{pfibdsymwelldefobj-i} If the multiplicative braiding $\betate$ in $\cD$ is the identity natural transformation, then the symmetric group action in $\pfibd$ is well defined on opcartesian strong multilinear functors.
\item\label{pfibdsymwelldefobj-ii} If the symmetric group action in $\pfibd$ is well defined on opcartesian strong multilinear functors, then $\otimes = \tensorop$ in $\cD$.
\end{enumerate}
\end{lemma}

\begin{proof}
First suppose $\betate = 1$ in $\cD$.  We want to show that the symmetric group action in $\permcatsg$ restricts to $\pfibd$.  Suppose given small permutative opfibrations over $\cD$ (\cref{def:permutativefibration})
\[\begin{tikzcd}[column sep=large]
\C_k \ar{r}{P_k} & \cD
\end{tikzcd}
\forspace k \in \{0,1,\ldots,n\}\]
with $\angC = \ang{\C_j}_{j=1}^n$ and $\ang{P} = \ang{P_j}_{j=1}^n$.  For a permutation $\sigma \in \Sigma_n$, to show that the symmetric group action
\begin{equation}
\begin{tikzcd}[column sep=large]
\pfibd\scmap{\ang{P};P_0} \ar{r}{\sigma} & \pfibd\scmap{\ang{P}\sigma;P_0}
\end{tikzcd}
\end{equation}
is well defined on objects, suppose given an opcartesian $n$-linear functor (\cref{def:opcartnlinearfunctor})
\[\begin{tikzcd}[column sep=2.3cm]
\ang{P} \ar{r}{\brb{F,\ang{F^2_j}_{j=1}^n}} & P_0.
\end{tikzcd}\]
We must show that the $n$-linear functor 
\begin{equation}\label{sigmathenf}
\begin{tikzcd}[column sep=large]
\txprod_{j=1}^n \C_{\sigma(j)} \ar{r}{\sigma} & \txprod_{j=1}^n \C_j \ar{r}{F} & \C_0
\end{tikzcd}
\end{equation}
with the linearity constraints in \cref{fsigmatwoj} satisfies axioms \cref{def:opcartnlinearfunctor-i,def:opcartnlinearfunctor-ii,def:opcartnlinearfunctor-iii,def:opcartnlinearfunctor-iv} in \cref{def:opcartnlinearfunctor}.  

\emph{Axiom \cref{def:opcartnlinearfunctor-i}}.  This axiom holds for $F^\sigma = F \circ \sigma$ because
\begin{itemize}
\item the permutation $\sigma$ sends an $n$-tuple of opcartesian morphisms to an $n$-tuple of opcartesian morphisms and
\item $F$ preserves opcartesian morphisms by axiom \cref{def:opcartnlinearfunctor-i}.
\end{itemize}
Note that axiom \cref{def:opcartnlinearfunctor-i} is valid even without the assumption $\betate = 1$.

\emph{Axiom \cref{def:opcartnlinearfunctor-ii}}.  The diagram \cref{def:opcartnlinearfunctor-diag} for $F^\sigma$ in \cref{sigmathenf} is the boundary of the following diagram.
\begin{equation}\label{pfibdnonsymmetrydiag}
\begin{tikzpicture}[xscale=2.2,yscale=1.3,vcenter]
\draw[0cell=.9]
(0,0) node (x11) {\txprod_{j=1}^n \C_{\sigma(j)}}
(x11)++(1.2,0) node (x12) {\txprod_{j=1}^n \C_j} 
(x12)++(1,0) node (x13) {\C_0}
(x11)++(0,-1) node (x21) {\cD^n}
(x12)++(0,-1) node (x22) {\cD^n}
(x13)++(0,-1) node (x23) {\cD}
;
\draw[1cell=.9]  
(x11) edge node {\sigma} (x12)
(x12) edge node {F} (x13)
(x21) edge node {\sigma} (x22)
(x22) edge node {\otimes} (x23)
(x11) edge node[swap] {\smallprod_j P_{\sigma(j)}} (x21)
(x12) edge node[swap] {\smallprod_j P_j} (x22)
(x13) edge node {P_0} (x23)
(x21) edge[bend right=50] node[swap] (te) {\otimes} (x23)
;
\draw[2cell=.9]
node[between=te and x22 at .5, rotate=90, 2label={above,\betate}, 2label={below,\iso}] {\Rightarrow};
\end{tikzpicture}
\end{equation}
The detail of the diagram \cref{pfibdnonsymmetrydiag} is as follows.
\begin{itemize}
\item The left square commutes by the naturality of the braiding in $(\Cat,\times)$.
\item The right square commutes by axiom \cref{def:opcartnlinearfunctor-ii} for $F$.
\end{itemize}
In the bottom region, the natural isomorphism
\[\begin{tikzcd}[column sep=large]
\otimes \ar{r}{\betate}[swap]{\iso} & \otimes \circ \sigma
\end{tikzcd}\]
is an iterate of the multiplicative braiding $\betate$ in $\cD$.  Under the current assumption that $\betate = 1$, the bottom region in \cref{pfibdnonsymmetrydiag} strictly commutes, proving axiom \cref{def:opcartnlinearfunctor-ii} for $F^\sigma$.

\emph{Axiom \cref{def:opcartnlinearfunctor-iii}}.  This axiom for $F^\sigma$ follows from
\begin{itemize}
\item axiom \cref{def:opcartnlinearfunctor-ii} for $F^\sigma$ and
\item axiom \cref{def:opcartnlinearfunctor-iii} for $F$.
\end{itemize} 

\emph{Axiom \cref{def:opcartnlinearfunctor-iv}}.  This axiom for $F^\sigma$ follows from
\begin{itemize}
\item the definition \cref{fsigmatwoj} of the $i$-th linearity constraint $(F^\sigma)^2_i = F^2_{\sigma(i)}$,
\item the assumption that $\betate = 1$ in $\cD$,
\item axiom \cref{def:opcartnlinearfunctor-ii} for $F^\sigma$, and
\item axiom \cref{def:opcartnlinearfunctor-iv} for $F$.
\end{itemize} 
This finishes the proof of assertion \eqref{pfibdsymwelldefobj-i}.

For assertion \eqref{pfibdsymwelldefobj-ii}, suppose that the symmetric group action in $\pfibd$ is well defined on opcartesian strong multilinear functors.  For the opcartesian strong $2$-linear functor \cref{onedntooned}
\[\begin{tikzcd}[column sep=2.3cm]
\ang{1_{\Dplus}}_{j=1}^2 \ar{r}{\brb{\otimes,\ang{\lap_j^\inv}_{j=1}^2}} & 1_{\Dplus}
\end{tikzcd}\]
and $\sigma = (1,2) \in \Sigma_2$, the boundary of the commutative diagram \cref{pfibdnonsymmetrydiag} yields the equality
\[\otimes = \otimes \circ \xi = \tensorop.\]
This proves assertion \eqref{pfibdsymwelldefobj-ii}.
\end{proof}

Next we consider opcartesian multilinear transformations.

\begin{lemma}\label{pfibdsymwelldefmor}
In the context of \cref{def:pfibdmulticat}, the following statements hold.
\begin{enumerate}
\item\label{pfibdsymwelldefmor-i} If the multiplicative braiding $\betate$ in $\cD$ is the identity natural transformation, then the symmetric group action in $\pfibd$ is well defined on opcartesian multilinear transformations.
\item\label{pfibdsymwelldefmor-ii} If the symmetric group action in $\pfibd$ is well defined on opcartesian multilinear transformations, then $\otimes = \tensorop$ in $\cD$.
\end{enumerate}
\end{lemma}

\begin{proof}
First suppose $\betate = 1$ in $\cD$.  Reusing the notation in the proof of \cref{pfibdsymwelldefobj}, suppose given an opcartesian $n$-linear transformation $\theta$ in $\pfibd\scmap{\ang{P};P_0}$ as follows.
\[\begin{tikzpicture}[xscale=2,yscale=2,baseline={(a.base)}]
\draw[0cell=.9]
(0,0) node (a) {\textstyle\prod_{j=1}^n \C_j}
(a)++(.16,0) node (b) {\phantom{D}}
(b)++(1,0) node (c) {\C_0}
;
\draw[1cell=.8]  
(b) edge[bend left] node[pos=.5] {\big( F, \ang{F^2_j}\big)} (c)
(b) edge[bend right] node[swap,pos=.5] {\big( G, \ang{G^2_j}\big)} (c)
;
\draw[2cell] 
node[between=b and c at .45, rotate=-90, 2label={above,\theta}] {\Rightarrow}
;
\end{tikzpicture}\]
By definition \cref{permcatsigmaaction}, the axiom \cref{thetapzeropjtensor} for $\theta^\sigma = \theta * 1_{\sigma}$ is the equality
\begin{equation}\label{pzerothetasigma}
1_{P_0} * \theta * 1_\sigma = 1_{\otimes \circ (\smallprod_j P_{\sigma(j)})}
\end{equation}
as in the boundary of the following diagram.
\begin{equation}\label{pfibdnonsymmetrymordiag}
\begin{tikzpicture}[xscale=2.2,yscale=1.3,vcenter]
\draw[0cell=.9]
(0,0) node (x11) {\txprod_{j=1}^n \C_{\sigma(j)}}
(x11)++(1.2,0) node (x12) {\txprod_{j=1}^n \C_j} 
(x12)++(.15,0) node (t) {\phantom{\C_0}}
(x12)++(1,0) node (x13) {\C_0}
(x11)++(0,-1) node (x21) {\cD^n}
(x12)++(0,-1) node (x22) {\cD^n}
(x13)++(0,-1) node (x23) {\cD}
;
\draw[1cell=.85]  
(x11) edge node {\sigma} (x12)
(t) edge[bend left=40] node {F} (x13)
(t) edge[bend right=40] node[swap] {G} (x13)
(x21) edge node {\sigma} (x22)
(x22) edge node {\otimes} (x23)
(x11) edge node[swap] {\smallprod_j P_{\sigma(j)}} (x21)
(x12) edge node[swap] {\smallprod_j P_j} (x22)
(x13) edge node {P_0} (x23)
(x21) edge[bend right=50] node[swap] (te) {\otimes} (x23)
;
\draw[2cell=.9]
node[between=x12 and x13 at .53, rotate=-90, 2label={above,\theta}] {\Rightarrow}
node[between=te and x22 at .5, rotate=90, 2label={above,\betate}, 2label={below,\iso}] {\Rightarrow}
;
\end{tikzpicture}
\end{equation}
The detail of the diagram \cref{pfibdnonsymmetrymordiag} is as follows.
\begin{itemize}
\item The left square commutes by the naturality of the braiding in $(\Cat,\times)$.
\item To its right, the equality 
\[1_{P_0} * \theta = 1_{\otimes \circ (\smallprod_j P_j)}\]
holds by the axiom \cref{thetapzeropjtensor} for $\theta$.
\item Under the current assumption that $\betate = 1$ in $\cD$, the bottom region strictly commutes.
\end{itemize}  
This proves the equality \cref{pzerothetasigma}.

For assertion \eqref{pfibdsymwelldefmor-ii}, suppose that the symmetric group action in $\pfibd$ is well defined on opcartesian multilinear transformations.  Similar to the proof of \cref{pfibdsymwelldefobj} \eqref{pfibdsymwelldefobj-ii}, the desired equality, $\otimes = \tensorop$ in $\cD$, follows from considering the diagram \cref{pfibdnonsymmetrymordiag} in the case $n=2$ with
\begin{itemize}
\item $\sigma \in \Sigma_2$ the nonidentity permutation,
\item $F = G$ the opcartesian strong 2-linear functor \cref{onedntooned}
\[\begin{tikzcd}[column sep=2.3cm]
\ang{1_{\Dplus}}_{j=1}^2 \ar{r}{\brb{\otimes,\ang{\lap_j^\inv}_{j=1}^2}} & 1_{\Dplus},
\end{tikzcd}\]
and
\item $\theta = 1_F$ the identity opcartesian 2-linear transformation of $F$.
\end{itemize} 
This finishes the proof.
\end{proof}

\subsection*{Multicategory Structure}

Next is the main result of this section.  Recall the $\Cat$-multicategory $\permcatsg$ (\cref{thm:permcatenrmulticat}).

\begin{theorem}\label{thm:pfibdmulticat}\index{permutative opfibration!non-symmetric multicategory}\index{non-symmetric!multicategory!permutative opfibration}
Suppose $\cD$ is a small tight bipermutative category.
\begin{enumerate}[label=(\roman*)]
\item\label{thm:pfibdmulticat-i}
$\pfibd$ in \cref{def:pfibdmulticat} is a non-symmetric $\Cat$-multicategory.
\item\label{thm:pfibdmulticat-ii}
If the multiplicative braiding $\betate$ in $\cD$ is the identity natural transformation, then $\pfibd$ is a $\Cat$-multicategory.
\item\label{thm:pfibdmulticat-iii} If the non-symmetric $\Cat$-multicategory $\pfibd$ is a $\Cat$-multicategory, then $\otimes = \tensorop$ in $\cD$.
\end{enumerate}
\end{theorem}

\begin{proof}
By \cref{pfibdcomponentwelldef} each multimorphism category $\pfibd\scmap{\ang{P};P_0}$ is a subcategory of the corresponding multimorphism category $\permcatsg\scmap{\angC;\C_0}$. 
\begin{itemize}
\item In both $\pfibd$ and $\permcatsg$, the colored units are given by identity 1-linear functors. 
\item The composition in $\pfibd$, which is a restriction of the composition in $\permcatsg$, is well defined by \cref{pfibdcompobjwelldef,pfibdcompmorwelldef}.
\end{itemize}
Therefore, the associativity and unity axioms, \cref{enr-multicategory-associativity,enr-multicategory-right-unity,enr-multicategory-left-unity}, hold in $\pfibd$ as they do in the $\Cat$-multicategory $\permcatsg$.  This proves assertion \cref{thm:pfibdmulticat-i}.

For assertion \cref{thm:pfibdmulticat-ii}, under the assumption that $\betate = 1$ in $\cD$, the symmetric group action in $\pfibd$ is well defined by \cref{pfibdsymwelldefobj} \eqref{pfibdsymwelldefobj-i} and \cref{pfibdsymwelldefmor} \eqref{pfibdsymwelldefmor-i}.  The symmetry group action axiom, \cref{enr-multicategory-symmetry}, and the equivariance axioms, \cref{enr-operadic-eq-1,enr-operadic-eq-2}, hold in $\pfibd$ as they do in the $\Cat$-multicategory $\permcatsg$.

Assertion \cref{thm:pfibdmulticat-iii} follows from \cref{pfibdsymwelldefobj} \eqref{pfibdsymwelldefobj-ii} and \cref{pfibdsymwelldefmor} \eqref{pfibdsymwelldefmor-ii}.
\end{proof}

Recall from \cref{def:enr-multicategory-functor} that a non-symmetric $\V$-multifunctor is a $\V$-multifunctor without the symmetric group action axiom \cref{enr-multifunctor-equivariance}.  In other words, a non-symmetric $\V$-multifunctor is only required to preserve the colored units \cref{enr-multifunctor-unit} and the composition \cref{v-multifunctor-composition}.  The following observation follows from \cref{thm:pfibdmulticat}.

\begin{corollary}\label{cor:pfibdtopermcatsg}\index{non-symmetric!multifunctor}\index{multifunctor!non-symmetric}
Suppose $\cD$ is a small tight bipermutative category.
\begin{enumerate}[label=(\roman*)]
\item\label{cor:pfibdtopermcatsg-i}
There exists a non-symmetric $\Cat$-multifunctor
\begin{equation}\label{pfibdtopermcatsg}
\begin{tikzcd}[column sep=large]
\pfibd \ar{r}{\U} & \permcatsg
\end{tikzcd}
\end{equation}
defined by
\begin{itemize}
\item the object assignment
\[\U(P \cn \C \to \cD) = \C\]
for each small permutative opfibration $P$ over $\cD$ and
\item the subcategory inclusions in \cref{pfibdcomponentwelldef} for multimorphism categories.
\end{itemize}
\item\label{cor:pfibdtopermcatsg-ii}
If the multiplicative braiding $\betate$ in $\cD$ is the identity natural transformation, then $\U$ is a $\Cat$-multifunctor.
\end{enumerate}
\end{corollary}

\begin{explanation}\label{expl:pfibdtopermcatsg}
Although $\U$ in \cref{pfibdtopermcatsg} consists of subcategory inclusions between multimorphism categories, $\pfibd$ is, in general, \emph{not} a sub-$\Cat$-multicategory of $\permcatsg$.  The reason is that $\U$ is, in general, not injective on objects because a small permutative category $\C$ may admit several different permutative opfibrations $\C \to \cD$.
\end{explanation}

\begin{example}\label{ex:pfibdcatmulticat}
Consider the small tight bipermutative categories
\[\Finsk,\quad \Fset, \quad \Fskel, \andspace \cA\]
in, respectively, \cref{ex:finsk,ex:Fset,ex:Fskel,ex:mandellcategory}.  With $\cD$ denoting any one of them, 
\begin{itemize}
\item $\pfibd$ is a non-symmetric $\Cat$-multicategory by \cref{thm:pfibdmulticat} \cref{thm:pfibdmulticat-i}, and
\item there is a non-symmetric $\Cat$-multifunctor $\U$ as in \cref{pfibdtopermcatsg}.
\end{itemize}
Since $\otimes \neq \tensorop$ in $\cD$, $\pfibd$ does \emph{not} have a well-defined symmetric group action by \cref{pfibdsymwelldefobj} \eqref{pfibdsymwelldefobj-ii}.  Therefore, $\pfibd$ is \emph{not} a $\Cat$-multicategory.
\end{example}

\section{Permutative Opfibrations via the Grothendieck Construction}
\label{sec:popfibgrothendieck}

As in \cref{Dstbipermutative} we continue to assume that $(\cD,\oplus,\otimes)$ is a small tight bipermutative category.  For an additive symmetric monoidal functor (\cref{def:additivesmf})
\[(X,X^2,X^0) \cn \Dplus = (\cD,\oplus,\zero,\betaplus) \to \Cat,\]
its Grothendieck construction $\grod X$ (\cref{def:grothendieckconst}) is a small permutative category (\cref{def:permutativegroconst,grodxpermutative})
\[\big(\grod X, \gbox, (\zero, X^0*), \betabox\big).\]
By \cref{thm:grocatmultifunctor} $X \mapsto\! \grod X$ is the object assignment of a pseudo symmetric $\Cat$-multifunctor (\cref{def:pseudosmultifunctor})
\[\begin{tikzcd}[column sep=large,every label/.append style={scale=.85}]
\DCat \ar{r}{\grod} & \permcatsg.
\end{tikzcd}\]
It is, furthermore, a $\Cat$-multifunctor if the multiplicative braiding $\betate$ in $\cD$ is the identity natural transformation.  

In this section we show that $\grod$ lifts to a non-symmetric $\Cat$-multifunctor against the non-symmetric $\Cat$-multifunctor $\U$ in \cref{pfibdtopermcatsg} as follows.
\[\begin{tikzpicture}[xscale=2.75,yscale=1.2,vcenter]
\draw[0cell=.9]
(0,0) node (dc) {\DCat}
(dc)++(1,0) node (pc) {\permcatsg}
(pc)++(0,1) node (pf) {\pfibd}
;
\draw[1cell=.9]  
(dc) edge node[pos=.6] {\grod} (pc)
(pf) edge node {\U} (pc)
(dc) edge[densely dashed,bend left=20] node[pos=.7] {\grod} (pf)
;
\end{tikzpicture}\]
Moreover, if $\betate = 1$ in $\cD$, then this lift is a $\Cat$-multifunctor.  In other words, it strictly preserves the symmetric group action.

\subsection*{Extension of Object Assignment}

First we extend the Grothendieck construction $\grod X$ to a permutative opfibration $U_X$ over $\cD$ (\cref{def:permutativefibration}).

\begin{lemma}\label{grodxpermopfib}
For each additive symmetric monoidal functor
\[(X,X^2,X^0) \cn \Dplus \to \Cat,\]
the first-factor projection\index{first-factor projection}
\begin{equation}\label{Usubx}
\begin{tikzcd}[column sep=large]
\grod X \ar{r}{U_X} & \cD
\end{tikzcd}
\end{equation}
defined by
\begin{itemize}
\item $U_X(a,x) = a$ for each object $(a,x)$ in $\grod X$ and
\item $U_X(f,p) = f$ for each morphism $(f,p)$ in $\grod X$
\end{itemize}
is a small permutative opfibration over $\cD$.
\end{lemma}

\begin{proof}
Given
\begin{itemize}
\item an object $(a,x) \in \grod X$ and 
\item a morphism $f \cn a \to b \in \cD$,
\end{itemize}  
the chosen opcartesian lift (\cref{def:opfibration} \cref{def:opfibration-vii}) of the fore-lift $\fl{(a,x)}{f}$ is defined as the morphism in $\grod X$
\begin{equation}\label{fonechosenoplift}
\begin{tikzcd}[column sep=huge]
(a,x) \ar{r}{(f,1_{f_*x})} & (b,f_*x)
\end{tikzcd}
\end{equation}
with 
\[f_* = Xf \cn Xa \to Xb.\]  
This definition makes $U_X$ into a split opfibration over $\cD$ by \cite[10.1.11]{johnson-yau}.  Axioms \cref{def:permutativefibration-i,def:permutativefibration-ii,def:permutativefibration-iii} in \cref{def:permutativefibration} follow from
\begin{itemize}
\item the definition of the permutative structure $\big(\gbox, (\zero, X^0*), \betabox\big)$ in \cref{groconunit,gboxobjects,gboxmorphisms,groconbeta} and
\item the functoriality of each component functor of the monoidal constraint $X^2$ \cref{additivesmfconstraint}.
\end{itemize} 
This finishes the proof.
\end{proof}

The assignment $X \mapsto U_X$, with $U_X \cn \grod X \to \cD$ in \cref{Usubx}, is the object assignment of the lifted Grothendieck construction.

\subsection*{Lifting Multimorphism Functors}

Next we observe that each multimorphism functor of $\grod$ lands in $\pfibd$.  Recall that the multimorphism functors of $\grod$ are in
\begin{itemize}
\item \cref{def:groconarityzero} for arity 0 inputs and
\item \cref{def:groconarityn,def:groconmodification,grodmultimorphismfunctor} for positive arity inputs.
\end{itemize} 
For the next several results, suppose given additive symmetric monoidal functors
\begin{equation}\label{zxdplustocat}
(Z,Z^2,Z^0) \andspace \bang{(X_j,X_j^2,X_j^0)}_{j=1}^n \cn \Dplus \to \Cat
\end{equation}
with $n \geq 0$.  We consider the solid-arrow functors
\begin{equation}\label{grodcomponentfactors}
\begin{tikzpicture}[xscale=4,yscale=1.4,vcenter]
\draw[0cell=.85]
(0,0) node (dc) {\DCat\scmap{\angX;Z}}
(dc)++(1,0) node (pc) {\permcatsg\scmap{\ang{\grod X};\grod Z}}
(pc)++(0,1) node (pf) {\pfibd\scmap{\ang{U_X};U_Z}}
;
\draw[1cell=.85]  
(dc) edge node[pos=.5] {\grod} (pc)
(pf) edge[right hook->] node {\U} (pc)
(dc) edge[densely dashed,bend left=20] (pf)
;
\end{tikzpicture}
\end{equation}
with
\begin{itemize}
\item $\angX = \ang{X_j}_{j=1}^n$, $\ang{\grod X} = \ang{\grod X_j}_{j=1}^n$, and $\ang{U_X} = \ang{U_{X_j}}_{j=1}^n$,\label{not:anggrodx}
\item each 
\[\begin{tikzcd}[column sep=large]
\grod X_j \ar{r}{U_{X_j}} & \cD
\end{tikzcd} \andspace 
\begin{tikzcd}[column sep=large]
\grod Z \ar{r}{U_Z} & \cD
\end{tikzcd}\]
as in \cref{Usubx},
\item $\grod$ the functor in \cref{groarityzero} if $n=0$ and \cref{grodfunctors} if $n>0$, and
\item $\U$ the subcategory inclusion in \cref{pfibdcomponentwelldef}.
\end{itemize}
\cref{grodfactorarityzero,grodfactorarityn} show that $\grod$ lifts against $\U$ as the dotted arrow.  First we consider the arity 0 case.

\begin{lemma}\label{grodfactorarityzero}
In the context of \cref{grodcomponentfactors} with $n=0$, the functor
\[\begin{tikzcd}[column sep=large,every label/.append style={scale=.85}]
\DCat\scmap{\ang{};Z} \ar{r}{\grod} & \permcatsg\scmap{\ang{};\grod Z}
\end{tikzcd}\]
factors through the subcategory $\pfibd\scmap{\ang{};U_Z}$ and yields an isomorphism
\[\begin{tikzcd}[column sep=large,every label/.append style={scale=.85}]
\DCat\scmap{\ang{};Z} \ar{r}{\grod}[swap]{\iso} & \pfibd\scmap{\ang{};U_Z}.
\end{tikzcd}\]
\end{lemma}

\begin{proof}
The subcategory
\[\pfibd\scmap{\ang{};U_Z} \bigsubset \grod Z\]
is the $U_Z$-preimage of $\tu \in \cD$ by
\begin{itemize}
\item \cref{expl:opcartnlinearfunctor} \eqref{expl:opcartnlinearfunctor-0} for objects and
\item \cref{expl:opcartnlineartr} \eqref{expl:opcartnlineartr-ii} for morphisms.
\end{itemize} 
In other words, in $\pfibd\scmap{\ang{};U_Z}$
\begin{itemize}
\item each object has the form $(\tu,x)$ for some object $x \in Z\tu$ and
\item each morphism has the form 
\[\begin{tikzcd}[column sep=large]
(\tu,x) \ar{r}{(1_{\tu},p)} & (\tu,y)
\end{tikzcd}\]
for some morphism $p \cn x \to y$ in $Z\tu$.
\end{itemize}
Therefore, the multimorphism functor in \cref{groarityzero} factors through the subcategory 
\[\pfibd\scmap{\ang{};U_Z} \bigsubset \grod Z = \permcatsg\scmap{\ang{};\grod Z}\]
and yields an isomorphism
\[\begin{tikzcd}[column sep=large,every label/.append style={scale=.85}]
Z\tu = \DCat\scmap{\ang{};Z} \ar{r}{\grod}[swap]{\iso} & \pfibd\scmap{\ang{};U_Z}
\end{tikzcd}\]
of categories.
\end{proof}

Next we consider the positive arity case.  In some diagrams below, the symbols $\exists$ and $!$ mean, respectively, \emph{there exists} and \emph{unique}.

\begin{lemma}\label{grodfactorarityn}
In the context of \cref{grodcomponentfactors} with $n>0$, the image of the functor
\[\begin{tikzcd}[column sep=large,every label/.append style={scale=.85}]
\DCat\scmap{\angX;Z} \ar{r}{\grod} & \permcatsg\scmap{\ang{\grod X};\grod Z}
\end{tikzcd}\]
lies in the subcategory $\pfibd\scmap{\ang{U_X};U_Z}$ of the codomain. 
\end{lemma}

\begin{proof}
First we consider objects, which are additive natural transformations, followed by morphisms, which are additive modifications.

\emph{Objects}.  Given an additive natural transformation \cref{additivenattransformation}
\[\begin{tikzcd}[column sep=large]
\bang{(X_j, X_j^2, X_j^0)}_{j=1}^n \ar{r}{\phi} & (Z,Z^2,Z^0)
\end{tikzcd} \inspace \DCat\scmap{\angX;Z},\]
we must show that the strong $n$-linear functor (\cref{grodphifunctor})
\begin{equation}\label{grodphistrongnlinear}
\begin{tikzcd}[column sep=3.4cm, every label/.append style={scale=.8}]
\prod_{j=1}^n \big(\grod X_j\big) \ar{r}{\brb{\grod \phi, \ang{(\grod \phi)_i^2}_{i=1}^n}} & \grod Z
\end{tikzcd}
\end{equation}
is opcartesian.  Using \cref{preserveopcmorphisms} we need to verify axioms \cref{def:opcartnlinearfunctor-ii,def:opcartnlinearfunctor-iii,def:opcartnlinearfunctor-iv} in \cref{def:opcartnlinearfunctor}.  We use the notation in \cref{def:groconarityn}.

\emph{Axiom \cref{def:opcartnlinearfunctor-ii}}.
The diagram \cref{def:opcartnlinearfunctor-diag} for $\grod\phi$ is as follows.
\[\begin{tikzpicture}[xscale=3,yscale=1.3,vcenter]
\draw[0cell=.9]
(0,0) node (x11) {\txprod_{j=1}^n (\grod X_j)}
(x11)++(1,0) node (x12) {\grod Z} 
(x11)++(0,-1) node (x21) {\cD^n}
(x12)++(0,-1) node (x22) {\cD}
;
\draw[1cell=.9]  
(x11) edge node {\grod\phi} (x12)
(x21) edge node {\otimes} (x22)
(x11) edge node[swap] {\smallprod_j U_{X_j}} (x21)
(x12) edge node {U_Z} (x22)
;
\end{tikzpicture}\]
This diagram is commutative by
\begin{itemize}
\item the definition \cref{Usubx} of $U_?$ as the first-factor projection and
\item the definitions \cref{grodphiangajxj,grodphiangfjpj} of $\grod\phi$ on, respectively, objects and morphisms.
\end{itemize}
This proves axiom \cref{def:opcartnlinearfunctor-ii}.

\emph{Axiom \cref{def:opcartnlinearfunctor-iii}}.
For each $j\in \{1,\ldots,n\}$, suppose given a fore-lift $\fl{(a_j,x_j)}{f_j}$ with respect to $U_{X_j} \cn \grod X_j \to \cD$.  Its chosen opcartesian lift is the morphism \cref{fonechosenoplift}
\[\begin{tikzcd}[column sep=2.2cm]
(a_j,x_j) \ar{r}{\brb{f_j,1_{(f_j)_*x_j}}} & \brb{b_j,(f_j)_* x_j}
\end{tikzcd} \inspace \grod X_j.\]
By
\begin{itemize}
\item the definition \cref{grodphiangfjpj} of $\grod\phi$ on morphisms and
\item the functoriality of $\phi_{\angb}$,
\end{itemize}
there are morphism equalities in $\grod Z$ as follows.
\[\begin{split}
&(\grod\phi) \bang{(f_j \scs 1_{(f_j)_*x_j})}_{j=1}^n \cn (\grod\phi)\ang{(a_j,x_j)} \to (\grod\phi)\bang{b_j,(f_j)_*x_j}\\
&= \brb{f, \phi_{\angb}\bang{1_{(f_j)_*x_j}}_{j=1}^n} \cn \brb{a, \phi_{\anga}\angx} \to \brb{b,\phi_{\angb} \ang{(f_j)_* x_j}}\\
&= \brb{f, 1_{f_*\phi_{\anga}\angx}} \cn \brb{a, \phi_{\anga}\angx} \to \brb{b, f_*\phi_{\anga}\angx} 
\end{split}\]
By definition \cref{fonechosenoplift} the last morphism is the chosen opcartesian lift of the fore-lift
\[\bfl{(\grod\phi)\ang{(a_j,x_j)}}{\txotimes_{j=1}^n f_j} 
= \bfl{\brb{a, \phi_{\anga}\angx}}{f}\]
with respect to $U_Z \cn \grod Z \to \cD$.  Thus $\grod\phi$ preserves chosen opcartesian lifts, proving axiom \cref{def:opcartnlinearfunctor-iii}.

\emph{Axiom \cref{def:opcartnlinearfunctor-iv}}.
This axiom follows from
\begin{itemize}
\item the definition \cref{grodphilinearity} of each linearity constraint $(\grod\phi)^2_i$ and
\item the definition \cref{fonechosenoplift} of the chosen opcartesian lift of each fore-lift with respect to $U_Z \cn \grod Z \to \cD$.
\end{itemize}
We have shown that $\grod\phi$ in \cref{grodphistrongnlinear} is opcartesian, so it is an object in the subcategory $\pfibd\scmap{\ang{U_X};U_Z}$.

\emph{Morphisms}.
Given an additive modification \cref{additivemodtwocell}
\[\begin{tikzpicture}[xscale=1,yscale=1,baseline={(x1.base)}]
\draw[0cell=.9]
(0,0) node (x1) {\bang{(X_j, X_j^2, X_j^0)}_{j=1}^n}
(x1)++(1,0) node (x2) {\phantom{X}}
(x2)++(2,0) node (x3) {\phantom{X}}
(x3)++(.65,0) node (x4) {(Z,Z^2,Z^0)}
;
\draw[1cell=.9]  
(x2) edge[bend left] node[pos=.5] {\phi} (x3)
(x2) edge[bend right] node[swap,pos=.5] {\psi} (x3)
;
\draw[2cell]
node[between=x2 and x3 at .45, rotate=-90, 2label={above,\Phi}] {\Rightarrow}
;
\end{tikzpicture}
\inspace \DCat\scmap{\angX;Z},\]
we must prove the axiom \cref{thetapzeropjtensor} for $\grod\Phi$ (\cref{grodPhinlinear}), which asserts that the following two whiskered natural transformations are equal.
\[\begin{tikzpicture}[xscale=2.4,yscale=2,baseline={(a.base)}]
\def\v{-.65}
\draw[0cell=.9]
(0,0) node (a) {\textstyle\prod_{j=1}^n (\grod X_j)}
(a)++(.28,0) node (b) {\phantom{D}}
(b)++(1,0) node (c) {\grod Z}
(a)++(0,\v) node (x21) {\cD^n}
(c)++(0,\v) node (x22) {\cD}
;
\draw[1cell=.85]  
(b) edge[bend left] node[pos=.53] {\grod \phi} (c)
(b) edge[bend right] node[swap,pos=.53] {\grod \psi} (c)
(a) edge node[swap] {\smallprod_j U_{X_j}} (x21)
(c) edge node {U_Z} (x22)
(x21) edge node[pos=.35] {\otimes} (x22)
;
\draw[2cell] 
node[between=b and c at .35, rotate=-90, 2label={above,\grod\Phi}] {\Rightarrow}
;
\end{tikzpicture}\]
These two natural transformations are equal by
\begin{itemize}
\item the definition \cref{grodPhiajxj} of each component 
\[(\grod\Phi)_{\ang{(a_j,x_j)}} = \brb{1_a, (\Phi_{\anga})_{\angx}}\]
and
\item the definition \cref{Usubx} of $U_?$ as the first-factor projection.
\end{itemize}
This shows that $\grod\Phi$ is a morphism in the subcategory $\pfibd\scmap{\ang{U_X};U_Z}$.
\end{proof}

We are now ready for the main result of this chapter.

\begin{theorem}\label{thm:dcatpfibd}\index{bipermutative category!Grothendieck construction}\index{Grothendieck construction!bipermutative category}\index{non-symmetric!multifunctor!Grothendieck construction}\index{diagram category!non-symmetric multifunctor}\index{permutative opfibration!non-symmetric multifunctor}
Suppose $\cD$ is a small tight bipermutative category.
\begin{enumerate}[label=(\roman*)]
\item\label{thm:dcatpfibd-i}
There is a non-symmetric $\Cat$-multifunctor
\begin{equation}\label{groddcatpfibd}
\begin{tikzpicture}[xscale=2.75,yscale=1.2,baseline={(dc.base)}]
\draw[0cell=.9]
(0,0) node (dc) {\DCat}
(dc)++(1,0) node (pf) {\pfibd}
;
\draw[1cell=.9]  
(dc) edge node {\grod} (pf)
;
\end{tikzpicture}
\end{equation}
defined by
\begin{itemize}
\item the object assignment 
\[X \mapsto \big(U_X \cn \grod X \to \cD\big)\] 
in \cref{Usubx} and
\item the multimorphism functors from \cref{grodfactorarityzero,grodfactorarityn}.
\end{itemize}
\item\label{thm:dcatpfibd-ii} 
The diagram of non-symmetric $\Cat$-multifunctors
\[\begin{tikzpicture}[xscale=2.75,yscale=1.2,vcenter]
\draw[0cell=.9]
(0,0) node (dc) {\DCat}
(dc)++(1,0) node (pc) {\permcatsg}
(pc)++(0,1) node (pf) {\pfibd}
;
\draw[1cell=.9]  
(dc) edge node[pos=.6] {\grod} (pc)
(pf) edge node {\U} (pc)
(dc) edge[bend left=20] node[pos=.7] {\grod} (pf)
;
\end{tikzpicture}\]
is commutative, with the bottom $\grod$ and $\U$ from, respectively, \cref{grodpscatmultifunctor,pfibdtopermcatsg}.
\item\label{thm:dcatpfibd-iii} 
If the multiplicative braiding $\betate$ in $\cD$ is the identity natural transformation, then $\grod$ in \cref{groddcatpfibd} is a $\Cat$-multifunctor.
\end{enumerate} 
\end{theorem}

\begin{proof}
By \cref{thm:pfibdmulticat} $\pfibd$ is a non-symmetric $\Cat$-multicategory whose colored units and composition (\cref{def:pfibdmulticat}) are both inherited from $\permcatsg$ (\cref{thm:permcatenrmulticat}) via the subcategory inclusions in \cref{pfibdcomponentwelldef}.  The latter are the multimorphism functors of $\U$ in \cref{pfibdtopermcatsg}.  If $\betate = 1$ in $\cD$, then $\pfibd$ is a $\Cat$-multicategory whose symmetric group action is inherited from $\permcatsg$.

Assertions \cref{thm:dcatpfibd-i,thm:dcatpfibd-ii} follow from \cref{thm:grocatmultifunctor} \cref{thm:grocatmultifunctor-i}, \cref{thm:pfibdmulticat} \cref{thm:pfibdmulticat-i}, and \cref{grodxpermopfib,grodfactorarityzero,grodfactorarityn}.

Assertion \cref{thm:dcatpfibd-iii} follows from \cref{thm:grocatmultifunctor} \cref{thm:grocatmultifunctor-ii} and \cref{thm:pfibdmulticat} \cref{thm:pfibdmulticat-ii}.
\end{proof}

\begin{example}\label{ex:dcatgrodpfibd}
Suppose $\cD$ is any one of the small tight bipermutative categories
\[\Finsk,\quad \Fset, \quad \Fskel, \andspace \cA\]
in, respectively, \cref{ex:finsk,ex:Fset,ex:Fskel,ex:mandellcategory}.  Then the following statements hold.
\begin{itemize}
\item $\pfibd$ is a non-symmetric $\Cat$-multicategory and \emph{not} a $\Cat$-multicategory by \cref{ex:pfibdcatmulticat}.  
\item \cref{thm:dcatpfibd} \cref{thm:dcatpfibd-i,thm:dcatpfibd-ii} apply to $\cD$.
\item $\grod$ in \cref{groddcatpfibd} is \emph{not} a $\Cat$-multifunctor because $\pfibd$ is not a $\Cat$-multicategory.\defmark
\end{itemize} 
\end{example}

\chapter{The Grothendieck Construction is a \texorpdfstring{$\Cat$}{Cat}-Multiequivalence}
\label{ch:gromultequiv}
\subsection*{Context}

Throughout this chapter, we assume that 
\[\big(\cD, (\oplus, \zero, \betaplus), (\otimes, \tu, \betate), (\fal, \far)\big)\]
is a small tight bipermutative category (\cref{def:embipermutativecat}) with 
\begin{itemize}
\item additive structure $\Dplus = (\cD, \oplus, \zero, \betaplus)$ and
\item multiplicative structure $\Dte = (\cD, \otimes, \tu, \betate)$.
\end{itemize}
In \cref{thm:dcatpfibd} \cref{thm:dcatpfibd-ii} we observe that the pseudo symmetric $\Cat$-multifunctor
\[\begin{tikzpicture}[xscale=2.75,yscale=1.2,vcenter]
\draw[0cell=.9]
(0,0) node (dc) {\DCat}
(dc)++(1,0) node (pc) {\permcatsg};
\draw[1cell=.9]  
(dc) edge node[pos=.5] {\grod} (pc);
\end{tikzpicture}\]
in \cref{thm:grocatmultifunctor} factors through the non-symmetric $\Cat$-multicategory $\pfibd$ (\cref{thm:pfibdmulticat}) into two non-symmetric $\Cat$-multifunctors, $\grod$ \cref{groddcatpfibd} and $\U$ \cref{pfibdtopermcatsg}.
\[\begin{tikzpicture}[xscale=3, yscale=1]
\draw[0cell=.9]
(0,0) node (x11) {\DCat}
(x11)++(1,0) node (x12) {\permcatsg}
(x11)++(.5,-1) node (x21) {\pfibd}
;
\draw[1cell=.9]
(x11) edge node {\grod} (x12)
(x11) edge node[swap,pos=.3] {\grod} (x21)
(x21) edge node[swap,pos=.7] {\U} (x12)
;
\end{tikzpicture}\]

\subsection*{Purpose}

Recall from \cref{def:catmultiequivalence} that a \emph{$\Cat$-multiequivalence} is an equivalence in the 2-category $\catmulticat$ of small $\Cat$-multicategories.  The non-symmetric variant is defined similarly.  By \cref{thm:multiwhitehead} a $\Cat$-multifunctor is a $\Cat$-multiequivalence if and only if it is essentially surjective on objects and an isomorphism on each multimorphism category.

The main purpose of this chapter is to strengthen the factorization above by showing that the non-symmetric $\Cat$-multifunctor
\begin{equation}\label{grodnscatmultieq}
\begin{tikzpicture}[xscale=2.75,yscale=1.2,baseline={(x1.base)}]
\draw[0cell=.9]
(0,0) node (x1) {\DCat}
(x1)++(1,0) node (x2) {\pfibd};
\draw[1cell=.9]  
(x1) edge node[pos=.5] {\grod} (x2);
\end{tikzpicture}
\end{equation}
is a non-symmetric $\Cat$-\emph{multiequivalence}.  Moreover, if the multiplicative braiding in $\cD$ is the identity natural transformation, then this $\grod$ is a $\Cat$-multiequivalence.  See \cref{thm:dcatpfibdeq}.  These assertions are proved using \cref{thm:multiwhitehead}.

The following table summaries the main properties of the two $\grod$ and $\U$.  We abbreviate \emph{pseudo symmetric} and \emph{non-symmetric} to, respectively, \emph{ps} and \emph{ns}.
\begin{center}
\resizebox{.85\width}{!}{%
{\renewcommand{\arraystretch}{1.4}%
{\setlength{\tabcolsep}{1ex}
\begin{tabular}{|cc|c|c|}\hline
\multicolumn{3}{|c|}{$\Cat$-multifunctor/multiequivalence} & if $\betate = 1$ in $\cD$\\ \hline
$\grod \cn \DCat \to \permcatsg$ & (\ref{thm:grocatmultifunctor}) & ps $\Cat$-multifunctor & $\Cat$-multifunctor \\ \hline
$\U \cn \pfibd \to \permcatsg$ & (\ref{cor:pfibdtopermcatsg}) & ns $\Cat$-multifunctor & $\Cat$-multifunctor\\ \hline
$\grod \cn \DCat \fto{\sim} \pfibd$ & (\ref{thm:dcatpfibdeq}) & ns $\Cat$-multiequivalence & $\Cat$-multiequivalence\\ \hline
\end{tabular}}}}
\end{center}
\smallskip

\subsection*{Implication for Inverse $K$-Theory}

The inverse $K$-theory functor (\cref{ch:multifunctorA,ch:invK}) is the composite $\groa \circ A$ of two functors, with $\cA$ the small tight bipermutative category in \cref{ex:mandellcategory}.  Applied to $\cA$, \cref{thm:dcatpfibdeq} implies that the pseudo symmetric $\Cat$-multifunctor $\groa$ in inverse $K$-theory factors through $\pfiba$ as displayed below.
\[\begin{tikzpicture}[xscale=3, yscale=1]
\draw[0cell=.9]
(0,0) node (x11) {\ACat}
(x11)++(1,0) node (x12) {\permcatsg}
(x11)++(.5,-1) node (x21) {\pfiba}
;
\draw[1cell=.9]
(x11) edge node {\groa} (x12)
(x11) edge node[swap,pos=.3] {\groa} node[pos=.7] {\sim} (x21)
(x21) edge node[swap,pos=.7] {\U} (x12)
;
\end{tikzpicture}\]
Its first factor is the non-symmetric $\Cat$-multiequivalence 
\[\begin{tikzpicture}[xscale=2.75,yscale=1.2,vcenter]
\draw[0cell=.9]
(0,0) node (x1) {\ACat}
(x1)++(1,0) node (x2) {\pfiba};
\draw[1cell=.9]  
(x1) edge node[pos=.5] {\groa} node[swap] {\sim} (x2);
\end{tikzpicture}\]
given by the Grothendieck construction.

\subsection*{Organization}

As we mentioned above, we prove that $\grod$ in \cref{grodnscatmultieq} is a non-symmetric $\Cat$-multiequivalence using \cref{thm:multiwhitehead}.  In other words, we must show that $\grod$ is essentially surjective on objects and an isomorphism on each multimorphism category.  In \cref{sec:essentiallysurjective} we first show that $\grod$ is essentially surjective on objects.  This means that each small permutative opfibration over $\cD$ is, up to an opcartesian 1-linear isomorphism, the Grothendieck construction of some additive symmetric monoidal functor. 

In \cref{sec:catfaithful} we show that each multimorphism functor of $\grod$ is injective on both objects and morphism sets.  Injectivity on objects means that each additive natural transformation is uniquely determined by its Grothendieck construction as a strong multilinear functor (\cref{grodarityninjobj}).  Injectivity on morphism sets means that each additive modification is uniquely determined by its Grothendieck construction as a multilinear natural transformation (\cref{grodarityninjmor}).

In \cref{sec:objectfullness} we show that each multimorphism functor of $\grod$ is surjective on objects.  This means that each opcartesian strong $n$-linear functor in $\pfibd$ is equal to the Grothendieck construction of some additive natural transformation in $\DCat$.

In \cref{sec:morphismfullness} we show that each multimorphism functor of $\grod$ is surjective on morphism sets, finishing the proof of \cref{thm:dcatpfibdeq}.  This means that each opcartesian $n$-linear transformation is equal to the Grothendieck construction of some additive modification.  The proof that $\grod$ is, furthermore, a $\Cat$-multiequivalence if $\betate = 1$ in $\cD$ requires no extra effort because the lemmas in this chapter do not depend on any special property of $\betate$.  Thus the symmetric case of \cref{thm:multiwhitehead} applies when $\betate = 1$ in $\cD$.

We remind the reader of our left normalized bracketing \cref{expl:leftbracketing} for iterated monoidal product.

\section{Essential Surjectivity on Objects}
\label{sec:essentiallysurjective}

This section contains the first main step toward showing that the Grothendieck construction $\grod$ in \cref{thm:dcatpfibd} \cref{thm:dcatpfibd-i} is a non-symmetric $\Cat$-multiequivalence.  We show that $\grod$ is essentially surjective on objects between the underlying 1-ary categories.  To clarify the presentation, we separate the proof of essential surjectivity into the following two parts:
\begin{enumerate}
\item \cref{splitopfgrothendieck} shows that any split opfibration $P \cn \E \to \B$ can be realized by a Grothendieck construction as $U_X \cn \grob X \to \B$, up to an isomorphism $\varphi \cn \grob X \to \E$ that preserves both opcartesian morphisms and chosen opcartesian lifts.  This step is only about split opfibrations and does not involve any permutative structure.
\item The essential surjectivity of $\grod$ is proved in \cref{grodesssurjective}.  For a small permutative opfibration $P \cn \E \to \cD$, we first equip $X \cn \cD \to \Cat$ in \cref{splitopfgrothendieck} with the structure of an additive symmetric monoidal functor $\Dplus \to \Cat$.  Then we show that $\varphi$ is an opcartesian 1-linear isomorphism in $\pfibd$. 
\end{enumerate}
As stated in the first paragraph of \cref{sec:permopfib}, results in \cite[Ch.\! 9 and 10]{johnson-yau} about fibrations and the Grothendieck construction there are applicable in our current context of opfibrations by considering the opposite functor between the opposite categories.

\begin{lemma}\label{splitopfgrothendieck}\index{split opfibration!Grothendieck construction}\index{Grothendieck construction!split opfibration}
Suppose $P \cn \E \to \B$ is a split opfibration.  Then there exist
\begin{itemize}
\item a functor $X \cn \B \to \Cat$ and
\item an isomorphism of categories
\begin{equation}\label{varphigrobxtoe}
\begin{tikzcd}[column sep=large]
\grob X \ar{r}{\varphi}[swap]{\iso} & \E
\end{tikzcd}
\end{equation}
that
\begin{enumerate}[label=(\roman*)]
\item makes the diagram
\[\begin{tikzpicture}[xscale=2.5, yscale=1]
\draw[0cell=.9]
(0,0) node (x11) {\grob X}
(x11)++(1,0) node (x12) {\E}
(x11)++(.5,-1) node (x21) {\B}
;
\draw[1cell=.9]
(x11) edge node {\varphi} node[swap] {\iso} (x12)
(x12) edge node[pos=.4] {P} (x21)
(x11) edge[shorten <=-.6ex] node[swap,pos=.3] {U_X} (x21)
;
\end{tikzpicture}\]
commutative and
\item preserves opcartesian morphisms and chosen opcartesian lifts.
\end{enumerate}
\end{itemize}
\end{lemma}

\begin{proof}
First we construct the functor $X \cn \B \to \Cat$.  

\emph{Objects}.  Given an object $a \in \B$, define the subcategory\index{fiber subcategory}
\begin{equation}\label{Xapinva}
Xa = P^\inv(a) \bigsubset \E.
\end{equation}
In other words, 
\begin{itemize}
\item an object $x \in Xa$ is an object in $\E$ such that $Px = a$, and
\item a morphism $p \cn x \to y$ in $Xa$ is a morphism in $\E$ such that $Pp = 1_a$ in $\B$.
\end{itemize}

\emph{Morphisms}.  Given a morphism $f \cn a \to b$ in $\B$, the functor
\[\begin{tikzcd}[column sep=large]
P^\inv(a) = Xa \ar{r}{Xf} & Xb = P^\inv(b)
\end{tikzcd}\]
is defined as follows.
\begin{itemize}
\item Each object $y \in P^\inv(a)$ is sent to the codomain
\begin{equation}\label{xfyyfpinvb}
(Xf)y = y_f \in P^\inv(b)
\end{equation}
of the \emph{chosen} opcartesian lift $\fbar \cn y \to y_f$ of the fore-lift $\fl{y}{f}$ with respect to $P$ (\cref{def:opfibration} \cref{def:opfibration-vii}).
\item Each morphism $p \cn y \to y'$ in $P^\inv(a)$ is sent to the \emph{unique} raise 
\begin{equation}\label{xfppfyfyfprime}
(Xf)p = p_f \cn y_f \to y'_f \inspace P^\inv(b)
\end{equation} 
of the fore-raise $\fr{\fbar}{\fbar' p}{1_b}$ (\cref{def:opfibration} \cref{def:opfibration-iv}), as depicted below.
\[\begin{tikzpicture}[xscale=2, yscale=1.2]
\def\d{.5}
\draw[0cell=.9]
(0,0) node (x11) {y'_f}
(x11)++(1,0) node (x12) {y'}
(x11)++(0,-1) node (x21) {y_f}
(x12)++(0,-1) node (x22) {y}
(x12)++(\d,-.5) node (s) {}
(s)++(.6,0) node (t) {}
(t)++(\d,.5) node (d11) {b}
(d11)++(1,0) node (d12) {a}
(d11)++(0,-1) node (d21) {b}
(d12)++(0,-1) node (d22) {a}
;
\draw[1cell=.9]
(x22) edge node[swap] {\fbar} (x21)
(x21) edge[densely dashed] node {p_f} (x11)
(x22) edge node[swap] {p} (x12)
(x12) edge node[swap] {\fbar'} (x11)
(s) edge[|->] node {P} (t)
(d22) edge node[swap] {f} (d21)
(d21) edge node {1_b} (d11)
(d22) edge node[swap] {1_a} (d12)
(d12) edge node[swap] {f} (d11)
;
\end{tikzpicture}\]
Here $\fbar$ and $\fbar'$ are the chosen opcartesian lifts of, respectively, $\fl{y}{f}$ and $\fl{y'}{f}$.  The unique existence of $p_f$ follows from the fact that $\fbar$ is an opcartesian morphism.
\end{itemize}

\emph{Functoriality}.  The functoriality of $X \cn \B \to \Cat$ follows from
\begin{itemize}
\item \cite[10.4.7]{johnson-yau}, which shows that $X$ is a strictly unitary pseudofunctor, and
\item the assumption that $P$ is a \emph{split} opfibration, which forces the lax functoriality constraint of $X$ to be the identity.
\end{itemize}
Moreover, by \cite[10.1.11]{johnson-yau} the first-factor projection functor \cref{Usubx}
\[\begin{tikzcd}[column sep=large]
\grob X \ar{r}{U_X} & \B
\end{tikzcd}\]
is a split opfibration.

\emph{The Functor $\varphi \cn \grob X \to \E$}.  For each object $(a,x) \in \grob X$ with $a \in \B$ and $x \in P^\inv(a)$, define the object
\begin{equation}\label{varphiaxx}
\varphi(a,x) = x \in \E.
\end{equation}
Suppose given a morphism 
\[\begin{tikzcd}[column sep=large]
(a,x) \ar{r}{(f,p)} & (b,y) \inspace \grob X
\end{tikzcd}\]
with
\begin{itemize}
\item $f \cn a \to b$ in $\B$,
\item $p \cn (Xf)x = x_f \to y$ in $P^\inv(b)$, and
\item $\fbar \cn x \to x_f$ the chosen opcartesian lift of $\fl{x}{f}$ with respect to $P$.
\end{itemize} 
Define $\varphi(f,p)$ as the following composite in $\E$.
\begin{equation}\label{varphifpfbarp}
\begin{tikzpicture}[xscale=4, yscale=1,vcenter]
\draw[0cell]
(0,0) node (x11) {\varphi(a,x)}
(x11)++(1,0) node (x12) {\varphi(b,y)}
(x11)++(0,-1) node (x21) {x}
(x21)++(.5,0) node (x22) {x_f}
(x12)++(0,-1) node (x23) {y}
;
\draw[1cell=.9]
(x11) edge node {\varphi(f,p)} (x12)
(x11) edge[equal] (x21)
(x12) edge[equal] (x23)
(x21) edge node {\fbar} (x22)
(x22) edge node {p} (x23)
;
\end{tikzpicture}
\end{equation}

\emph{Properties of $\varphi$}.  The assignments \cref{varphiaxx,varphifpfbarp} have the following properties.
\begin{itemize}
\item By \cite[10.4.12]{johnson-yau} $\varphi \cn \grob X \to \E$ is a functor such that 
\[U_X = P\varphi \cn \grob X \to \B.\]
\item By \cite[10.4.23]{johnson-yau} $\varphi$ is an isomorphism and preserves opcartesian morphisms.
\end{itemize}
To see that $\varphi$ preserves chosen opcartesian lifts, consider a fore-lift $\fl{(a,x)}{f}$ with respect to $U_X$.  Its chosen opcartesian lift is \cref{fonechosenoplift}
\[\begin{tikzcd}[column sep=huge]
(a,x) \ar{r}{(f,1_{x_f})} & (b,f_*x = x_f)
\end{tikzcd} \inspace \grob X.\]
By \cref{varphifpfbarp} its image under $\varphi$ is
\[\varphi(f,1_{x_f}) = 1_{x_f} \circ \fbar = \fbar \inspace \E.\]
This is the chosen opcartesian lift of $\bfl{\varphi(a,x) = x}{f}$ with respect to $P$. 
\end{proof}

Restricting to 1-ary multimorphism categories, each non-symmetric $\Cat$-multicategory has an underlying 2-category (\cref{ex:unarycategory,locallysmalltwocat}).  By the same restriction, a non-symmetric $\Cat$-multifunctor has an underlying 2-functor.  Next we observe that the underlying 2-functor of the Grothendieck construction $\grod$ is essentially surjective on objects.

\begin{lemma}\label{grodesssurjective}\index{Grothendieck construction!essential surjectivity}
Suppose $\cD$ is a small tight bipermutative category.  Then the non-symmetric $\Cat$-multifunctor in \cref{groddcatpfibd}
\[\begin{tikzpicture}[xscale=2.75,yscale=1.2,vcenter]
\draw[0cell=.9]
(0,0) node (dc) {\DCat}
(dc)++(1,0) node (pf) {\pfibd};
\draw[1cell=.9]  
(dc) edge node {\grod} (pf);
\end{tikzpicture}\]
is essentially surjective on objects.
\end{lemma} 

\begin{proof}
Suppose given an object in $\pfibd$, that is, a small permutative opfibration (\cref{def:permutativefibration})
\[\begin{tikzcd}[column sep=large]
(\E,\oplus,e,\beta) \ar{r}{P} & \Dplus = (\cD,\oplus,\zero,\betaplus).
\end{tikzcd}\]
We must show that there exist
\begin{itemize}
\item an object in $\DCat$, that is, an additive symmetric monoidal functor (\cref{def:additivesmf})
\[(X,X^2,X^0) \cn \Dplus \to \Cat\]
and
\item an invertible opcartesian 1-linear functor in $\pfibd$ (\cref{def:opcartnlinearfunctor})
\[\big(\grod X \fto{U_X} \cD\big) \fto[\iso]{\varphi} \big(\E \fto{P} \cD\big).\] 
\end{itemize}
By \cref{ex:idoponelinear} $\varphi$ is precisely a strict symmetric monoidal isomorphism \cref{strictsmfunctor}
\[\begin{tikzcd}[column sep=large]
\grod X \ar{r}{\varphi}[swap]{\iso} & \E
\end{tikzcd}\]
that 
\begin{itemize}
\item satisfies $U_X = P\varphi$ and
\item preserves chosen opcartesian lifts.
\end{itemize}
We use \cref{grodxpermopfib} and the functors in \cref{splitopfgrothendieck} with $\B = \cD$, namely
\[\begin{tikzcd}[column sep=large]
\cD \ar{r}{X} & \Cat
\end{tikzcd} \andspace
\begin{tikzcd}[column sep=large]
\grod X \ar{r}{\varphi}[swap]{\iso} & \E.
\end{tikzcd}\]
It remains to
\begin{itemize}
\item equip the functor $X$, defined in \cref{Xapinva,xfyyfpinvb,xfppfyfyfprime}, with the structure $(X^2,X^0)$ of a symmetric monoidal functor and
\item show that the functor $\varphi$, defined in \cref{varphiaxx,varphifpfbarp}, is strict symmetric monoidal.
\end{itemize}   
To clarify the presentation, we divide the rest of this proof into several steps, starting with the constructions of $X^0$ and $X^2$.

\newcounter{esssur-step}
\newcommand{\esssurstep}[1]{\medskip\textbf{Step \stepcounter{esssur-step}\arabic{esssur-step}}: #1\medskip}

\esssurstep{The Unit Constraint of $X$}

Define the unit constraint 
\[\begin{tikzcd}[column sep=large]
\boldone \ar{r}{X^0} & X\zero = P^\inv(\zero) \bigsubset \E
\end{tikzcd}\]
by the unit object \cref{additivesmfunitobj}
\begin{equation}\label{Xzerostare}
X^0* = e \in \E.
\end{equation}
This is well defined because $P$ is strict symmetry monoidal \cref{strictsmfunctor}, which includes the property 
\[Pe = \zero \inspace \cD.\]

\esssurstep{The Monoidal Constraint of $X$}

For objects $a,b \in \cD$ the $(a,b)$-component functor of the monoidal constraint $X^2$ \cref{additivesmfconstraint},
\begin{equation}\label{Xtwoabaplusb}
\begin{tikzcd}[column sep=large]
P^\inv(a) \times P^\inv(b) = Xa \times Xb \ar{r}{X^2_{a,b}} & X(a \oplus b) = P^\inv(a \oplus b),
\end{tikzcd}
\end{equation}
is defined by
\begin{equation}\label{Xtwoabxplusy}
X^2_{a,b}(x,y) = x \oplus y \in \E
\end{equation}
for both objects and morphisms $(x,y) \in P^\inv(a) \times P^\inv(b)$.
\begin{itemize}
\item The assignment $X^2_{a,b}$ is well defined because $P$ strictly preserves the monoidal product \cref{strictsmfunctor}.
\item The functoriality of $X^2_{a,b}$ follows from the functoriality of $\oplus$ in $\E$.
\end{itemize}

\esssurstep{The Naturality of the Monoidal Constraint}

For morphisms 
\[\begin{tikzcd}
a \ar{r}{f} & b
\end{tikzcd}\andspace 
\begin{tikzcd}
a' \ar{r}{f'} & b'
\end{tikzcd}\inspace \cD,\] 
the naturality diagram for $X^2$ \cref{Xtwoabaplusb} is as follows, with $(-)_* = X(-)$.
\begin{equation}\label{Xtwonatdiagram}
\begin{tikzpicture}[xscale=3.5, yscale=1.3,vcenter]
\draw[0cell=.9]
(0,0) node (x11) {P^\inv(a) \times P^\inv(a')}
(x11)++(1,0) node (x12) {P^\inv(a \oplus a')}
(x11)++(0,-1) node (x21) {P^\inv(b) \times P^\inv(b')}
(x12)++(0,-1) node (x22) {P^\inv(b \oplus b')}
;
\draw[1cell=.9]
(x11) edge node {\oplus} (x12)
(x11) edge node[swap] {f_* \times f'_*} (x21)
(x12) edge[transform canvas={xshift=-1em}] node {(f \oplus f')_*} (x22)
(x21) edge node {\oplus} (x22)
;
\end{tikzpicture}
\end{equation}

\emph{Objects}.  The following equalities in $\E$ show that the diagram \cref{Xtwonatdiagram} commutes for any pair of objects $(x,x') \in P^\inv(a) \times P^\inv(a')$.
\[\begin{aligned}
(f \oplus f')_* (x \oplus x') &= (x \oplus x')_{f \oplus f'} & \phantom{M} & \text{by \cref{xfyyfpinvb}}\\
&= x_f \oplus x'_{f'} & \phantom{M} & \text{by \cref{def:permutativefibration} \cref{def:permutativefibration-ii}}\\
&= f_* x \oplus f'_* x' & \phantom{M} & \text{by \cref{xfyyfpinvb}}
\end{aligned}\]
In the computation above, we remind the reader of our notation:
\begin{itemize}
\item $\fbar \cn x \to x_f$ is the chosen opcartesian lift of $\fl{x}{f}$ with respect to $P$.
\item $\fbar' \cn x' \to x'_{f'}$ is the chosen opcartesian lift of $\fl{x'}{f'}$.
\item $\barof{f \oplus f'} = \fbar \oplus \fbar' \cn x \oplus x' \to (x \oplus x')_{f \oplus f'}$ is the chosen opcartesian lift of $\fl{x \oplus x'}{f \oplus f'}$, in which the equality follows from axiom \cref{def:permutativefibration-ii} in \cref{def:permutativefibration}.
\end{itemize}

\emph{Morphisms}.  To check the commutativity of \cref{Xtwonatdiagram} on morphisms, suppose given a pair of morphisms
\begin{equation}\label{morphismsppprime}
\begin{tikzcd}
x \ar{r}{p} & y \in \Pinv(a)
\end{tikzcd} \andspace 
\begin{tikzcd}
x' \ar{r}{p'} & y' \in \Pinv(a').
\end{tikzcd}
\end{equation}
By \cref{def:permutativefibration} \cref{def:permutativefibration-ii} and \cref{xfppfyfyfprime}, the morphism in $\Pinv(b \oplus b')$ 
\[(f \oplus f')_* (p \oplus p') = (p \oplus p')_{f \oplus f'} \cn (x \oplus x')_{f \oplus f'} \to (y \oplus y')_{f \oplus f'}\]
is the \emph{unique} raise of the fore-raise
\begin{equation}\label{foreraiseffprime}
\bfr{\fbar \oplus \fbar'}{(\fbar \oplus \fbar')(p \oplus p')}{1_{b \oplus b'}}
\end{equation}
as depicted below.
\[\begin{tikzpicture}[xscale=2, yscale=1.2]
\def\d{.5} \def\h{.75}
\draw[0cell=.9]
(0,0) node (x11) {(y \oplus y')_{f \oplus f'} = y_f \oplus y'_{f'}}
(x11)++(1.7,0) node (x12) {y \oplus y'}
(x11)++(0,-1) node (x21) {(x \oplus x')_{f \oplus f'} = x_f \oplus x'_{f'}}
(x12)++(0,-1) node (x22) {x \oplus x'}
(x12)++(\d,-\h) node (s) {}
(s)++(.5,0) node (t) {}
(t)++(\d,\h) node (d11) {b \oplus b'}
(d11)++(1,0) node (d12) {a \oplus a'}
(d11)++(0,-1) node (d21) {b \oplus b'}
(d12)++(0,-1) node (d22) {a \oplus a'}
;
\draw[1cell=.85]
(x22) edge node[swap] {\fbar \oplus \fbar'} (x21)
(x21) edge[densely dashed, shorten >=-.5ex, transform canvas={xshift=1.7ex}] node {(p \oplus p')_{f \oplus f'}} (x11)
(x22) edge node[swap] {p \oplus p'} (x12)
(x12) edge node[swap] {\fbar \oplus \fbar'} (x11)
(s) edge[|->] node {P} (t)
(d22) edge node[swap] {f \oplus f'} (d21)
(d21) edge node {1_{b \oplus b'}} (d11)
(d22) edge node[swap] {1_{a \oplus a'}} (d12)
(d12) edge node[swap] {f \oplus f'} (d11)
;
\end{tikzpicture}\]
By \cref{xfppfyfyfprime} for $p_f$ and $p'_{f'}$ the sum 
\[p_f \oplus p'_{f'} = f_*p \oplus f'_*p'\]
is also a raise of the fore-raise \cref{foreraiseffprime}.  So uniqueness implies that
\[f_*p \oplus f'_*p' = (f \oplus f')_* (p \oplus p').\]
This proves that the diagram \cref{Xtwonatdiagram} is commutative.

\esssurstep{Symmetric Monoidal Functor Axioms}

Next we check that $(X,X^2,X^0)$ is a symmetric monoidal functor.  

\emph{Unity}.  Since $\Dplus$ is a permutative category, by definitions \cref{Xzerostare,Xtwoabxplusy}, the unity diagrams \cref{monoidalfunctorunity} are as follows.
\[\begin{tikzpicture}[xscale=2.5, yscale=1.2,vcenter]
\def\s{.85}
\draw[0cell=\s]
(0,0) node (x11) {\boldone \times P^\inv(b)}
(x11)++(1,0) node (x12) {P^\inv(b)}
(x11)++(0,-1) node (x21) {P^\inv(\zero) \times P^\inv(b)}
(x12)++(0,-1) node (x22) {P^\inv(\zero \oplus b)}
;
\draw[1cell=\s]
(x11) edge node {\iso} (x12)
(x11) edge node[swap] {e \times 1} (x21)
(x21) edge node {\oplus} (x22)
(x22) edge[equal] (x12)
;
\begin{scope}[shift={(2,0)}]
\draw[0cell=\s]
(0,0) node (x11) {\Pinv(a) \times \boldone}
(x11)++(1,0) node (x12) {P^\inv(a)}
(x11)++(0,-1) node (x21) {P^\inv(a) \times P^\inv(\zero)}
(x12)++(0,-1) node (x22) {P^\inv(a \oplus \zero)}
;
\draw[1cell=\s]
(x11) edge node {\iso} (x12)
(x11) edge node[swap] {1 \times e} (x21)
(x21) edge node {\oplus} (x22)
(x22) edge[equal] (x12)
;
\end{scope}
\end{tikzpicture}\]
These two diagrams commute because $e$ is the monoidal unit in the permutative category $\E$.

\emph{Associativity}.  The associativity diagram \cref{monoidalfunctorassoc} of $(X,X^2)$ is as follows for objects $a,b,c \in \cD$.
\[\begin{tikzpicture}[xscale=4.5, yscale=1.2,vcenter]
\def\s{.85}
\draw[0cell=\s]
(0,0) node (x11) {(\Pinv(a) \times P^\inv(b)) \times \Pinv(c)}
(x11)++(1,0) node (x12) {\Pinv(a) \times (P^\inv(b) \times \Pinv(c))}
(x11)++(0,-1) node (x21) {P^\inv(a \oplus b) \times P^\inv(c)}
(x12)++(0,-1) node (x22) {\Pinv(a) \oplus P^\inv(b \oplus c)}
(x21)++(0,-1) node (x31) {P^\inv(a \oplus b \oplus c)}
(x22)++(0,-1) node (x32) {\Pinv(a \oplus b \oplus c)}
;
\draw[1cell=\s]
(x11) edge node {\iso} (x12)
(x31) edge[equal] (x32)
(x11) edge node[swap] {\oplus \times 1} (x21)
(x21) edge node[swap] {\oplus} (x31)
(x12) edge node {1 \times \oplus} (x22)
(x22) edge node {\oplus} (x32)
;
\end{tikzpicture}\]
This diagram commutes because $\oplus$ in $\E$ is strictly associative.

\emph{Braiding}.  Denoting by
\begin{itemize}
\item $\betaplus \cn a \oplus a' \to a' \oplus a$ the braiding in $\cD$,
\item $\betaplus_* = X\betaplus \cn \Pinv(a \oplus a') \to \Pinv(a' \oplus a)$ as in \cref{Xapinva,xfyyfpinvb,xfppfyfyfprime}, and
\item $\xi$ the braiding in $(\Cat,\times)$,
\end{itemize}
the compatibility diagram \cref{monoidalfunctorbraiding} of $(X,X^2)$ with the braidings is as follows.
\begin{equation}\label{xioplusbetaplusstar}
\begin{tikzpicture}[xscale=3, yscale=1.2,vcenter]
\def\s{.9}
\draw[0cell=\s]
(0,0) node (x11) {\Pinv(a) \times P^\inv(a')}
(x11)++(1,0) node (x12) {P^\inv(a \oplus a')}
(x11)++(0,-1) node (x21) {P^\inv(a') \times P^\inv(a)}
(x12)++(0,-1) node (x22) {P^\inv(a' \oplus a)}
;
\draw[1cell=\s]
(x11) edge node {\oplus} (x12)
(x12) edge node {\betaplus_*} (x22)
(x11) edge node[swap] {\xi} (x21)
(x21) edge node {\oplus} (x22)
;
\end{tikzpicture}
\end{equation}
The following equalities in $\E$ show that the diagram \cref{xioplusbetaplusstar} commutes for any pair of objects $(x,x') \in P^\inv(a) \times P^\inv(a')$.
\[\begin{aligned}
\betaplus_* (x \oplus x') &= (x \oplus x')_{\betaplus} & \phantom{M} & \text{by \cref{xfyyfpinvb}}\\
&= x' \oplus x & \phantom{M} & \text{by \cref{def:permutativefibration} \cref{def:permutativefibration-iii}}\\
&= (\oplus \circ \xi)(x,x') &&
\end{aligned}\]

To check the commutativity of \cref{xioplusbetaplusstar} on morphisms, suppose given a pair of morphisms
$(p,p') \in \Pinv(a) \times \Pinv(a')$ as in \cref{morphismsppprime}.  By \cref{def:permutativefibration} \cref{def:permutativefibration-iii} and \cref{xfppfyfyfprime}, the morphism in $\Pinv(a' \oplus a)$ 
\[\betaplus_* (p \oplus p') = (p \oplus p')_{\betaplus} \cn (x \oplus x')_{\betaplus} \to (y \oplus y')_{\betaplus}\]
is the \emph{unique} raise of the fore-raise
\begin{equation}\label{foreraisebetaplusppprime}
\bfr{\beta_{x,x'}}{\beta_{y,y'}(p \oplus p')}{1_{a' \oplus a}}
\end{equation}
as depicted below.
\[\begin{tikzpicture}[xscale=2, yscale=1.2]
\def\d{.5} \def\h{.75}
\draw[0cell=.9]
(0,0) node (x11) {(y \oplus y')_{\betaplus} = y' \oplus y}
(x11)++(1.5,0) node (x12) {y \oplus y'}
(x11)++(0,-1) node (x21) {(x \oplus x')_{\betaplus} = x' \oplus x}
(x12)++(0,-1) node (x22) {x \oplus x'}
(x12)++(\d,-\h) node (s) {}
(s)++(.5,0) node (t) {}
(t)++(\d,\h) node (d11) {a' \oplus a}
(d11)++(1,0) node (d12) {a \oplus a'}
(d11)++(0,-1) node (d21) {a' \oplus a}
(d12)++(0,-1) node (d22) {a \oplus a'}
;
\draw[1cell=.85]
(x22) edge node[swap] {\beta_{x,x'}} (x21)
(x21) edge[densely dashed, shorten >=-.3ex, transform canvas={xshift=2ex}] node {(p \oplus p')_{\betaplus}} (x11)
(x22) edge node[swap] {p \oplus p'} (x12)
(x12) edge node[swap] {\beta_{y,y'}} (x11)
(s) edge[|->] node {P} (t)
(d22) edge node[swap] {\betaplus_{a,a'}} (d21)
(d21) edge node {1_{a' \oplus a}} (d11)
(d22) edge node[swap] {1_{a \oplus a'}} (d12)
(d12) edge node[swap] {\betaplus_{a,a'}} (d11)
;
\end{tikzpicture}\]
By the naturality of $\beta$, the sum $p' \oplus p$ is also a raise of the fore-raise \cref{foreraisebetaplusppprime}.  So uniqueness implies that
\[\betaplus_* (p \oplus p') = (p \oplus p')_{\betaplus} = p' \oplus p.\]
This proves that the diagram \cref{xioplusbetaplusstar} is commutative.  We have shown that $(X,X^2,X^0)$ is a symmetric monoidal functor.

\esssurstep{Strict Symmetric Monoidality of $\varphi$}

Next we check that the categorical isomorphism
\[\begin{tikzcd}[column sep=large]
\big(\grod X, \gbox, (\zero,e), \betabox\big) \ar{r}{\varphi}[swap]{\iso} & (\E,\oplus,e,\beta),
\end{tikzcd}\]
defined in \cref{varphiaxx,varphifpfbarp}, is a strict symmetric monoidal functor as in \cref{strictsmfunctor}.  

\emph{Monoidal Unit}.  The functor $\varphi$ preserves the monoidal unit because, by definition \cref{varphiaxx},
\[\varphi(\zero,e) = e \in \E.\]

\emph{Objects}.  The functor $\varphi$ preserves the monoidal product on objects by the following equalities in $\E$ for objects $(a,x)$ and $(b,y) \in \grod X$.
\[\begin{aligned}
\varphi\big((a,x) \gbox (b,y)\big)
&= \varphi(a \oplus b, x \oplus y) & \phantom{M} & \text{by \cref{gboxobjects,Xtwoabxplusy}}\\
&= x \oplus y & \phantom{M} & \text{by \cref{varphiaxx}}\\
&= \varphi(a,x) \oplus \varphi(b,y) & \phantom{M} & \text{by \cref{varphiaxx}}
\end{aligned}\]

\emph{Morphisms}.  To show that $\varphi$ preserves the monoidal product on morphisms, suppose given two morphisms in $\grod X$,
\[\begin{tikzcd}[column sep=large]
(a,x) \ar{r}{(f,p)} & (a',x')
\end{tikzcd} \andspace 
\begin{tikzcd}[column sep=large]
(b,y) \ar{r}{(g,q)} & (b',y'),
\end{tikzcd}\]
with
\begin{itemize}
\item $f \cn a \to a'$ and $g \cn b\to b'$ in $\cD$,
\item $p \cn f_*x = x_f \to x'$ in $Xa' = \Pinv(a')$ by \cref{Xapinva,xfyyfpinvb}, and
\item $q \cn g_*y = y_g \to y'$ in $Xb' = \Pinv(b')$.
\end{itemize}
Here 
\[\begin{tikzcd}
x \ar{r}{\fbar} & x_f
\end{tikzcd} \andspace 
\begin{tikzcd}
y \ar{r}{\gbar} & y_g
\end{tikzcd} \inspace \E\]
are the chosen opcartesian lifts of, respectively, the fore-lifts $\fl{x}{f}$ and $\fl{y}{g}$ with respect to $P$.  The functor $\varphi$ preserves the monoidal product on morphisms by the following morphism equalities in $\E$.
\[\begin{aligned}
\varphi\big((f,p) \gbox (g,q)\big)
&= \varphi(f \oplus g, p \oplus q) & \phantom{M} & \text{by \cref{gboxmorphisms,Xtwoabxplusy}}\\
&= (p \oplus q) \barof{(f \oplus g)} & \phantom{M} & \text{by \cref{varphifpfbarp}}\\
&= (p \oplus q) (\fbar \oplus \gbar) & \phantom{M} & \text{by \cref{def:permutativefibration} \cref{def:permutativefibration-ii}}\\
&= p\fbar \oplus q\gbar & \phantom{M} & \text{by functoriality of $\oplus$}\\
&= \varphi(f,p) \oplus \varphi(g,q) & \phantom{M} & \text{by \cref{varphifpfbarp}}
\end{aligned}\]

\emph{Braiding}.
The functor $\varphi$ preserves the braiding by the following morphism equalities in $\E$.
\[\begin{aligned}
\varphi\betabox_{(a,x),(b,y)}
&= \varphi(\betaplus_{a,b}, 1_{y \oplus x}) & \phantom{M} & \text{by \cref{groconbeta,Xtwoabxplusy}}\\
&= 1_{y \oplus x} \barof{\betaplus_{a,b}} & \phantom{M} & \text{by \cref{varphifpfbarp}}\\
&= \beta_{x,y} & \phantom{M} & \text{by \cref{def:permutativefibration} \cref{def:permutativefibration-iii}}\\
&= \beta_{\varphi(a,x),\varphi(b,y)} & \phantom{M} & \text{by \cref{varphiaxx}}
\end{aligned}\]
Therefore, we may define the unit constraint $\varphi^0$ and the monoidal constraint $\varphi^2$ of $\varphi$ to be the identities, and $(\varphi, \varphi^2, \varphi^0)$ is a strict symmetric monoidal functor.
\end{proof}

\section{Faithfulness on Multimorphism Categories}
\label{sec:catfaithful}

For the rest of this chapter, suppose 
\begin{equation}\label{zangxjdpluscat}
(Z,Z^2,Z^0) \andspace \bang{(X_j,X_j^2,X_j^0)}_{j=1}^n \cn \Dplus \to \Cat
\end{equation}
are additive symmetric monoidal functors (\cref{def:additivesmf}) with $n \geq 0$.  We use the notation in \cref{grodcomponentfactors}, so 
\begin{equation}\label{angxgrodxux}
\angX = \ang{X_j}_{j=1}^n, \quad \ang{\grod X} = \ang{\grod X_j}_{j=1}^n, \andspace
\ang{U_X} = \ang{U_{X_j}}_{j=1}^n
\end{equation}
with each $U_{X_j} \cn \grod X_j \to \cD$ the small permutative opfibration in \cref{grodxpermopfib}.

By \cref{grodfactorarityzero} the arity 0 multimorphism functor is an isomorphism
\[\begin{tikzpicture}[xscale=4.5,yscale=1.2,vcenter]
\draw[0cell=.9]
(0,0) node (dc) {Z\tu = \DCat\scmap{\ang{};Z}}
(dc)++(1,0) node (pf) {\pfibd\scmap{\ang{};U_Z} = U_Z^\inv(\tu)};
\draw[1cell=.9]  
(dc) edge node {\grod} node[swap] {\iso} (pf);
\end{tikzpicture}\]
with $U_Z \cn \grod Z \to \cD$ (\cref{grodxpermopfib}).  Our next objective is to show that the multimorphism functors of $\grod$ are also isomorphisms for positive arity inputs.  In this section, we show that they are injective on objects and morphisms; see \cref{grodarityninjobj,grodarityninjmor}.  By \cref{def:dcatnarycategory} the objects and morphisms in $\DCat\scmap{\angX;Z}$ are, respectively, additive natural transformations \cref{additivenattransformation} and additive modifications \cref{additivemodtwocell}.

\begin{lemma}\label{grodarityninjobj}\index{Grothendieck construction!object injectivity}
For $n>0$ the functor from \cref{grodfactorarityn}
\[\begin{tikzpicture}[xscale=3.8,yscale=1.2,vcenter]
\draw[0cell=.9]
(0,0) node (dc) {\DCat\scmap{\angX;Z}}
(dc)++(1,0) node (pf) {\pfibd\scmap{\ang{U_X};U_Z}};
\draw[1cell=.9]  
(dc) edge node {\grod} (pf);
\end{tikzpicture}\]
is injective on objects.
\end{lemma}

\begin{proof}
Suppose given additive natural transformations \cref{additivenattransformation}
\begin{equation}\label{phiandpsi}
\begin{tikzcd}[column sep=large]
\bang{(X_j, X_j^2, X_j^0)}_{j=1}^n \ar[shift left]{r}{\phi} \ar[shift right]{r}[swap]{\psi} & (Z,Z^2,Z^0)
\end{tikzcd}
\end{equation}
such that the following two strong $n$-linear functors (\cref{grodphifunctor}) are equal.  
\begin{equation}\label{grodphigrodpsi}
\begin{tikzcd}[column sep=3.4cm, every label/.append style={scale=.8}]
\prod_{j=1}^n (\grod X_j) \ar[shift left]{r}{\brb{\grod \phi, \ang{(\grod \phi)_i^2}_{i=1}^n}} 
\ar[shift right]{r}[swap]{\brb{\grod \psi, \ang{(\grod \psi)_i^2}_{i=1}^n}} & \grod Z
\end{tikzcd}
\end{equation}
We must show that $\phi = \psi$.  These additive natural transformations consist of component functors \cref{dcatnaryobjcomponent}
\begin{equation}\label{phiapsia}
\begin{tikzcd}[column sep=huge]
\prod_{j=1}^n X_j a_j \ar[shift left]{r}{\phi_{\anga}} \ar[shift right]{r}[swap]{\psi_{\anga}} & Za
\end{tikzcd}
\end{equation}
for each $n$-tuple of objects 
\begin{equation}\label{angajtensoraj}
\anga = \ang{a_j}_{j=1}^n \in \cD^n \withspace a = \txotimes_{j=1}^n a_j \in \cD.
\end{equation}
We must show that $\phi_{\anga} = \psi_{\anga}$ as functors.  

\emph{Objects}.  For each $n$-tuple of objects
\begin{equation}\label{angxjobject}
\angx = \ang{x_j}_{j=1}^n \in \txprod_{j=1}^n X_j a_j
\end{equation}
the desired equality,
\[\phi_{\anga}\angx = \psi_{\anga}\angx \inspace Za,\]
follows from the following equalities in $\grod Z$.
\[\begin{aligned}
\brb{a,\phi_{\anga}\angx} 
&= (\grod\phi) \ang{(a_j,x_j)} & \phantom{M} & \text{by \cref{grodphiangajxj}}\\
&= (\grod\psi) \ang{(a_j,x_j)} & \phantom{M} & \text{by \cref{grodphigrodpsi}}\\
&= \brb{a,\psi_{\anga}\angx} & \phantom{M} & \text{by \cref{grodphiangajxj}}
\end{aligned}\]

\emph{Morphisms}.  For each $n$-tuple of morphisms
\begin{equation}\label{angpjxjyj}
\angp = \bang{p_j \cn x_j \to y_j}_{j=1}^n \in \txprod_{j=1}^n X_j a_j
\end{equation} 
consider the $n$-tuple of morphisms
\[\bang{(1_{a_j},p_j) \cn (a_j,x_j) \to (a_j,y_j)}_{j=1}^n \in \txprod_{j=1}^n (\grod X_j).\]
The desired equality,
\[\phi_{\anga}\angp = \psi_{\anga}\angp \inspace Za,\]
follows from the following morphism equalities in $\grod Z$.
\[\begin{aligned}
\brb{1_a,\phi_{\anga}\angp} 
&= (\grod\phi) \ang{(1_{a_j},p_j)} & \phantom{M} & \text{by \cref{grodphiangfjpj}}\\
&= (\grod\psi) \ang{(1_{a_j},p_j)} & \phantom{M} & \text{by \cref{grodphigrodpsi}}\\
&= \brb{1_a,\psi_{\anga}\angp} & \phantom{M} & \text{by \cref{grodphiangfjpj}}
\end{aligned}\]
This proves that $\phi_{\anga} = \psi_{\anga}$ in \cref{phiapsia}.
\end{proof}

\begin{remark}\label{rk:grodarityninjobj}
The proof of \cref{grodarityninjobj} uses only the assignments of $\grod\phi$ and $\grod\psi$ on objects and morphisms of the form $\ang{(1_{a_j},p_j)}$.  Their linearity constraints, $\ang{(\grod\phi)^2_i}_{i=1}^n$ and $\ang{(\grod\psi)^2_i}_{i=1}^n$, are not used in the proof.
\end{remark}

\begin{lemma}\label{grodarityninjmor}\index{Grothendieck construction!morphism injectivity}
For $n>0$ the functor from \cref{grodfactorarityn}
\[\begin{tikzpicture}[xscale=3.8,yscale=1.2,vcenter]
\draw[0cell=.9]
(0,0) node (dc) {\DCat\scmap{\angX;Z}}
(dc)++(1,0) node (pf) {\pfibd\scmap{\ang{U_X};U_Z}};
\draw[1cell=.9]  
(dc) edge node {\grod} (pf);
\end{tikzpicture}\]
is injective on morphism sets.
\end{lemma}

\begin{proof}
Suppose given additive natural transformations $\phi$ and $\psi$ as in \cref{phiandpsi} and additive modifications \cref{additivemodtwocell}
\begin{equation}\label{PhiandTheta}
\begin{tikzpicture}[xscale=1,yscale=1,baseline={(x1.base)}]
\def\a{25}
\draw[0cell=.9]
(0,0) node (x1) {\bang{(X_j, X_j^2, X_j^0)}_{j=1}^n}
(x1)++(1,0) node (x2) {\phantom{X}}
(x2)++(2.5,0) node (x3) {\phantom{X}}
(x3)++(.65,0) node (x4) {(Z,Z^2,Z^0)}
;
\draw[1cell=.9]  
(x2) edge[bend left=\a] node[pos=.5] {\phi} (x3)
(x2) edge[bend right=\a] node[swap,pos=.5] {\psi} (x3)
;
\draw[2cell]
node[between=x2 and x3 at .4, rotate=-90, 2label={below,\Phi}] {\Rightarrow}
node[between=x2 and x3 at .6, rotate=-90, 2label={above,\Theta}] {\Rightarrow}
;
\end{tikzpicture}
\end{equation}
such that the following two $n$-linear natural transformations \cref{grodPhi} are equal.
\begin{equation}\label{grodPhiTheta}
\begin{tikzpicture}[xscale=1,yscale=1,baseline={(x1.base)}]
\def\a{25}
\draw[0cell=.9]
(0,0) node (x1) {\txprod_{j=1}^n (\grod X_j)}
(x1)++(.67,0) node (x2) {\phantom{X}}
(x2)++(3,0) node (x3) {\grod Z}
;
\draw[1cell=.85]  
(x2) edge[bend left=\a] node[pos=.53] {\grod \phi} (x3)
(x2) edge[bend right=\a] node[swap,pos=.53] {\grod \psi} (x3)
;
\draw[2cell=.9] 
node[between=x2 and x3 at .42, rotate=-90, 2label={below,\grod\Phi}] {\Rightarrow}
node[between=x2 and x3 at .55, rotate=-90, 2label={above,\grod\Theta}] {\Rightarrow}
;
\end{tikzpicture}
\end{equation}
We must show that $\Phi = \Theta$.

With the notation in \cref{phiapsia,angajtensoraj}, $\Phi$ and $\Theta$ consist of the following component natural transformations \cref{dcatnarymodcomponent}.
\begin{equation}\label{PhiaThetaa}
\begin{tikzpicture}[xscale=1,yscale=1,baseline={(x1.base)}]
\def\a{25}
\draw[0cell=.9]
(0,0) node (x1) {\txprod_{j=1}^n X_j a_j}
(x1)++(.5,0) node (x2) {\phantom{X}}
(x2)++(3,0) node (x3) {Za}
;
\draw[1cell=.9]  
(x2) edge[bend left=\a] node[pos=.53] {\phi_{\anga}} (x3)
(x2) edge[bend right=\a] node[swap,pos=.53] {\psi_{\anga}} (x3)
;
\draw[2cell]
node[between=x2 and x3 at .44, rotate=-90, 2label={below,\Phi_{\anga}}] {\Rightarrow}
node[between=x2 and x3 at .57, rotate=-90, 2label={above,\Theta_{\anga}}] {\Rightarrow}
;
\end{tikzpicture}
\end{equation}
For each $n$-tuple of objects $\angx \in \txprod_{j=1}^n X_j a_j$ as in \cref{angxjobject}, the desired equality,
\[(\Phi_{\anga})_{\angx} = (\Theta_{\anga})_{\angx} \cn \phi_{\anga}\angx \to \psi_{\anga}\angx \inspace Za,\]
follows from the following morphism equalities in $\grod Z$.
\[\begin{aligned}
\brb{1_a,(\Phi_{\anga})_{\angx}} 
&= (\grod\Phi)_{\ang{(a_j,x_j)}} & \phantom{M} & \text{by \cref{grodPhiajxj}}\\
&= (\grod\Theta)_{\ang{(a_j,x_j)}} & \phantom{M} & \text{by \cref{grodPhiTheta}}\\
&= \brb{1_a,(\Theta_{\anga})_{\angx}} & \phantom{M} & \text{by \cref{grodPhiajxj}}
\end{aligned}\]
This proves that $\Phi_{\anga} = \Theta_{\anga}$ in \cref{PhiaThetaa}.
\end{proof}

\section{Object Fullness on Multimorphism Categories}
\label{sec:objectfullness}

We continue to use the notation in \cref{zangxjdpluscat,angxgrodxux,angajtensoraj,angxjobject,angpjxjyj} with $n>0$.  In \cref{grodarityninjobj,,grodarityninjmor} we observe that the component functors (\cref{grodfactorarityn}) 
\[\begin{tikzpicture}[xscale=3.8,yscale=1.2,vcenter]
\draw[0cell=.9]
(0,0) node (dc) {\DCat\scmap{\angX;Z}}
(dc)++(1,0) node (pf) {\pfibd\scmap{\ang{U_X};U_Z}};
\draw[1cell=.9]  
(dc) edge node {\grod} (pf);
\end{tikzpicture}\] 
are injective on objects and morphism sets.  In this section we show that these component functors are surjective on objects.  Recall from \cref{def:pfibd} that 
\[\pfibd\scmap{\ang{U_X};U_Z} \bigsubset \permcatsg\scmap{\ang{\grod X};\grod Z}\]
is the subcategory with
\begin{itemize}
\item opcartesian strong $n$-linear functors (\cref{def:opcartnlinearfunctor}) as objects and
\item opcartesian $n$-linear transformations (\cref{def:opcartnlineartr}) as morphisms.
\end{itemize}

Here is an outline of this section.
\begin{itemize}
\item In \cref{def:grodobjectfull} we associate to each opcartesian strong $n$-linear functor $F$ in $\pfibd\scmap{\ang{U_X};U_Z}$ the data of an additive natural transformation $\phi$ in $\DCat\scmap{\angX;Z}$.  The projection axiom \cref{def:opcartnlinearfunctor-diag} plays a crucial role in this definition, as we discuss in \cref{expl:grodobjectfull}.
\item \cref{phiangafunctor,phinaturaltr} show that $\phi$ is a natural transformation \cref{phitensorxz}.
\item \cref{phiunityaxiom,phiadditivityaxiom} show that $\phi$ is an \emph{additive} natural transformation (\cref{def:additivenattr}).
\item \cref{grodobjectfullfphi} shows that $\grod\phi = F$ as strong $n$-linear functors.  This is sufficient to show the object surjectivity of the component functors of $\grod$ because an opcartesian $n$-linear functor is an $n$-linear functor with extra properties but no extra structure (\cref{expl:opcartnlinearfunctor} \eqref{expl:opcartnlinearfunctor-one}).
\end{itemize}

\begin{definition}\label{def:grodobjectfull}
Suppose given an opcartesian strong $n$-linear functor
\begin{equation}\label{anguxfuz}
\begin{tikzpicture}[xscale=1,yscale=1,baseline={(x1.base)}]
\draw[0cell=1]
(0,0) node (x1) {\ang{U_X}}
(x1)++(3,0) node (x2) {U_Z.};
\draw[1cell=.9]  
(x1) edge node {\brb{F,\ang{F^2_j}_{j=1}^n}} (x2);
\end{tikzpicture}
\end{equation}
This means a strong $n$-linear functor (\cref{def:nlinearfunctor})
\begin{equation}\label{anguxfuzperm}
\begin{tikzcd}[column sep=2.2cm]
\prod_{j=1}^n (\grod X_j) \ar{r}{\brb{F,\ang{F^2_j}_{j=1}^n}} & \grod Z
\end{tikzcd}
\end{equation}
that satisfies axioms \cref{def:opcartnlinearfunctor-i,def:opcartnlinearfunctor-ii,def:opcartnlinearfunctor-iii,def:opcartnlinearfunctor-iv} in \cref{def:opcartnlinearfunctor}.  Define the data of an additive natural transformation \cref{additivenattransformation}
\begin{equation}\label{Fpreimagephi}
\begin{tikzcd}[column sep=large]
\bang{(X_j, X_j^2, X_j^0)}_{j=1}^n \ar{r}{\phi} & (Z,Z^2,Z^0)
\end{tikzcd}
\end{equation}
with component functors
\begin{equation}\label{Fpreimagecomponent}
\begin{tikzcd}[column sep=large]
\prod_{j=1}^n X_j a_j \ar{r}{\phi_{\anga}} & Za
\end{tikzcd}
\end{equation}
as follows, where we use the notation in \cref{angajtensoraj,angxjobject,angpjxjyj}.
\begin{description}
\item[Objects] For each $n$-tuple of objects $\angx \in \txprod_{j=1}^n X_j a_j$, define the object
\begin{equation}\label{Fpreimageobject}
\phi_{\anga}\angx = F_2\angx \inspace Za
\end{equation}
where $F_2\angx$ is the second entry in the object
\begin{equation}\label{Fangajxjaftwox}
F\ang{(a_j,x_j)} = \big(a,F_2\angx\big) \inspace \grod Z.
\end{equation}
\item[Morphisms] For each $n$-tuple of morphisms $\angp \in \txprod_{j=1}^n X_j a_j$, define the morphism
\begin{equation}\label{Fpreimagemor}
\phi_{\anga}\angp = F_2\angp \cn F_2\angx \to F_2\angy \inspace Za
\end{equation}
where $F_2\angp$ is the second entry in the morphism
\begin{equation}\label{Fangoneajpj}
F\ang{(1_{a_j},p_j)} = \big(1_a, F_2\angp\big) \inspace \grod Z.
\end{equation}
\end{description}
This finishes the definition of $\phi$ in \cref{Fpreimagephi}.
\end{definition}

\begin{explanation}\label{expl:grodobjectfull}
The projection axiom \cref{def:opcartnlinearfunctor-diag} for $F$ in \cref{anguxfuzperm} is the following commutative diagram of functors.
\begin{equation}\label{projectionaxiomf}
\begin{tikzpicture}[xscale=2.8,yscale=1.3,vcenter]
\draw[0cell=.9]
(0,0) node (x11) {\txprod_{j=1}^n (\grod X_j)}
(x11)++(1,0) node (x12) {\grod Z} 
(x11)++(0,-1) node (x21) {\cD^n}
(x12)++(0,-1) node (x22) {\cD}
;
\draw[1cell=.9]  
(x11) edge node {F} (x12)
(x21) edge node {\otimes} (x22)
(x11) edge node[swap] {\smallprod_j U_{X_j}} (x21)
(x12) edge node {U_Z} (x22)
;
\end{tikzpicture}
\end{equation}
Since each $U_?$ \cref{Usubx} is the first-factor projection, evaluating the commutative diagram \cref{projectionaxiomf} at the object 
\[\ang{(a_j,x_j)} \in \txprod_{j=1}^n (\grod X_j)\] 
yields the equalities
\[U_Z F\ang{(a_j,x_j)} = \txotimes_{j=1}^n U_{X_j}(a_j,x_j) = \txotimes_{j=1}^n a_j = a.\]
So the first entry of the object 
\[F\ang{(a_j,x_j)} \in \grod Z\]
is $a \in \cD$.  This is the meaning of \cref{Fangajxjaftwox}.

Likewise, evaluating the commutative diagram \cref{projectionaxiomf} at the morphism
\[\bang{(1_{a_j},p_j) \cn (a_j,x_j) \to (a_j,y_j)} \in \txprod_{j=1}^n (\grod X_j)\]
yields the morphism equalities
\[U_Z F\ang{(1_{a_j},p_j)} = \txotimes_{j=1}^n U_{X_j} (1_{a_j},p_j) = \txotimes_{j=1}^n 1_{a_j} = 1_a.\]
So the first entry of the morphism 
\[F\ang{(1_{a_j},p_j)} \in \grod Z\]
is $1_a$.  This is the meaning of \cref{Fangoneajpj}.
\end{explanation}

In the rest of this section we show that $\phi$ in \cref{Fpreimagephi} is an additive natural transformation (\cref{def:additivenattr}) such that 
\[\brb{\grod\phi, \ang{(\grod\phi)^2_j}_{j=1}^n} = \brb{F, \ang{F^2_j}_{j=1}^n}\]
as strong $n$-linear functors, with the left-hand side from \cref{grodphifunctor}.

\subsection*{Naturality}

We begin by showing that each component of $\phi$ is well defined.

\begin{lemma}\label{phiangafunctor}
Each component $\phi_{\anga}$ in \cref{Fpreimagecomponent} is a functor.
\end{lemma}

\begin{proof}
There are morphism equalities in $\grod Z$ as follows.
\[\begin{aligned}
\big(1_a, F_2\ang{1_{x_j}}\big) 
&=F\ang{(1_{a_j}, 1_{x_j})} & \phantom{M} & \text{by \cref{Fangoneajpj}}\\
&= F\ang{1_{(a_j,x_j)}} & \phantom{M} & \text{by \cref{groconidentity}}\\
&= 1_{F\ang{(a_j,x_j)}} & \phantom{M} & \text{by functoriality of $F$}\\ 
&=1_{(a, F_2\angx)} & \phantom{M} & \text{by \cref{Fangajxjaftwox}}\\
&= \big(1_a, 1_{F_2\angx}\big) & \phantom{M} & \text{by \cref{groconidentity}}
\end{aligned}\]
So $\phi_{\anga}$ preserves identity morphisms by definitions \cref{Fpreimageobject,Fpreimagemor}.

For morphisms
\[\begin{tikzcd}[column sep=large]
x_j \ar{r}{p_j} & y_j \ar{r}{q_j} & z_j
\end{tikzcd} \inspace X_ja_j\]
there are morphism equalities in $\grod Z$ as follows.
\[\begin{aligned}
\big(1_a, F_2\ang{q_j p_j}\big) 
&= F\ang{(1_{a_j}, q_j p_j)} & \phantom{M} & \text{by \cref{Fangoneajpj}}\\
&= F\bang{(1_{a_j}, q_j) \circ (1_{a_j}, p_j)} & \phantom{M} & \text{by \cref{groconcomposition}}\\
&= F\ang{(1_{a_j}, q_j)} \circ F\ang{(1_{a_j}, p_j)} & \phantom{M} & \text{by functoriality of $F$}\\
&= \big(1_a, F_2\angq\big) \circ \big(1_a, F_2\angp\big) & \phantom{M} & \text{by \cref{Fangoneajpj}}\\
&= \big(1_a, F_2\angq \circ F_2\angp\big) & \phantom{M} & \text{by \cref{groconcomposition}}
\end{aligned}\]
So $\phi_{\anga}$ preserves composition of morphisms by definition \cref{Fpreimagemor}.
\end{proof}

Next we show that
\[\begin{tikzcd}[column sep=large]
\bigotimes_{j=1}^n X_j \ar{r}{\phi} & Z
\end{tikzcd}\]
is a natural transformation \cref{phitensorxz}.  Below we write $Zf$ as $f_*$ and $X_j f_j$ as $(f_j)_*$.

\begin{lemma}\label{phinaturaltr}
With the notation in \cref{aprimef,dcatnaryobjnaturality} and the component functors $\phi_{\anga}$ in \cref{Fpreimagecomponent}, the following naturality diagram commutes. 
\begin{equation}\label{phinaturalitydiagram}
\begin{tikzpicture}[xscale=3,yscale=1.3,vcenter]
\draw[0cell=.9]
(0,0) node (x11) {\txprod_{j=1}^n X_j a_j}
(x11)++(1,0) node (x12) {Za}
(x11)++(0,-1) node (x21) {\txprod_{j=1}^n X_j a_j'}
(x12)++(0,-1) node (x22) {Za'}
;
\draw[1cell=.9]  
(x11) edge node {\phi_{\anga}} (x12)
(x12) edge node {f_*} (x22)
(x11) edge node[swap] {\txprod_{j=1}^n (f_j)_*} (x21)
(x21) edge node {\phi_{\angap}} (x22)
;
\end{tikzpicture}
\end{equation}
\end{lemma}

\begin{proof}
We check the commutativity of \cref{phinaturalitydiagram} first on objects and then on morphisms.

\emph{Objects}.  Suppose given an $n$-tuple of objects $\angx \in \txprod_{j=1}^n X_j a_j$.  
Consider the fore-lift with respect to $U_Z \cn \grod Z \to \cD$
\begin{equation}\label{Fajxjtensorfj}
\bfl{F\ang{(a_j,x_j)}}{\txotimes_{j=1}^n f_j} = \bfl{(a, \phi_{\anga}\angx)}{f}
\end{equation}
with the equality from \cref{Fpreimageobject,Fangajxjaftwox}.  There are two ways to obtain its chosen opcartesian lift as follows.
\begin{itemize}
\item By definition \cref{fonechosenoplift} the chosen opcartesian lift of \cref{Fajxjtensorfj} is the morphism
\begin{equation}\label{Fajxjtensorfj-i}
\begin{tikzpicture}[xscale=1,yscale=1,baseline={(x11.base)}]
\draw[0cell=.9]
(0,0) node (x11) {\brb{a, \phi_{\anga}\angx}}
(x11)++(4.5,0) node (x12) {\brb{a', f_* \phi_{\anga}\angx}}
;
\draw[1cell=.9]  
(x11) edge node {\brb{f,1_{f_* \phi_{\anga}\angx}}} (x12)
;
\end{tikzpicture}
\inspace \grod Z.
\end{equation}
\item On the other hand, for each $j \in \{1,\ldots,n\}$, the chosen opcartesian lift of the fore-lift $\fl{(a_j,x_j)}{f_j}$ with respect to $U_{X_j} \cn \grod X_j \to \cD$ is, by definition \cref{fonechosenoplift}, the morphism
\begin{equation}\label{ajxjfjoneajpfjxj}
\begin{tikzpicture}[xscale=1,yscale=1,baseline={(x11.base)}]
\draw[0cell=.9]
(0,0) node (x11) {(a_j,x_j)}
(x11)++(3.5,0) node (x12) {\brb{a_j', (f_j)_* x_j}}
;
\draw[1cell=.9]  
(x11) edge node {\brb{f_j,1_{(f_j)_* x_j}}} (x12)
;
\end{tikzpicture}
\inspace \grod X_j.
\end{equation}
Since $F$ preserves chosen opcartesian lifts (\cref{def:opcartnlinearfunctor} \cref{def:opcartnlinearfunctor-iii}), the chosen opcartesian lift of \cref{Fajxjtensorfj} is also given by
\begin{equation}\label{Fajxjtensorfj-ii}
\begin{tikzpicture}[xscale=1,yscale=.8,vcenter]
\draw[0cell=.9]
(0,0) node (x11) {F\ang{(a_j,x_j)}}
(x11)++(4,0) node (x12) {F\bang{\brb{a_j', (f_j)_* x_j}}}
(x11)++(0,-1) node (x21) {\brb{a, \phi_{\anga}\angx}}
(x12)++(0,-1) node (x22) {\brb{a', \phi_{\angap}\ang{(f_j)_* x_j}}}
;
\draw[1cell=.9]  
(x11) edge node {F\ang{(f_j,1)}} (x12)
(x11) edge[equal, shorten >=-.5ex] (x21)
(x12) edge[equal, shorten <=-.5ex, shorten >=-.3ex] (x22)
;
\end{tikzpicture}
\end{equation}
with each equality from \cref{Fpreimageobject,Fangajxjaftwox}.
\end{itemize}
The codomains of \cref{Fajxjtensorfj-i,Fajxjtensorfj-ii} yield the object equality
\[f_* \phi_{\anga}\angx = \phi_{\angap}\ang{(f_j)_* x_j} \inspace Za'.\]
This proves the commutativity of \cref{phinaturalitydiagram} on objects.

\emph{Morphisms}.  Suppose given an $n$-tuple of morphisms $\angp \in \prod_{j=1}^n X_ja_j$ as in \cref{angpjxjyj}.  For each $j \in \{1,\ldots,n\}$, the following diagram in $\grod X_j$ is commutative by \cref{groconcomposition}.
\begin{equation}\label{fjonefjstarpj}
\begin{tikzpicture}[xscale=4,yscale=1.3,vcenter]
\draw[0cell=.9]
(0,0) node (x11) {(a_j,x_j)}
(x11)++(1,0) node (x12) {\brb{a_j', (f_j)_* x_j}}
(x11)++(0,-1) node (x21) {(a_j,y_j)}
(x12)++(0,-1) node (x22) {\brb{a_j', (f_j)_* y_j}}
;
\draw[1cell=.9]  
(x11) edge node {\brb{f_j, 1_{(f_j)_* x_j}}} (x12)
(x12) edge[transform canvas={xshift=-1em}] node {\brb{1_{a_j'}, (f_j)_* p_j}} (x22)
(x11) edge node[swap] {(1_{a_j}, p_j)} (x21)
(x21) edge node {\brb{f_j, 1_{(f_j)_* y_j}}} (x22)
;
\end{tikzpicture}
\end{equation}
Indeed, in each of the two composites in \cref{fjonefjstarpj},
\begin{itemize}
\item the first component is $f_j \cn a_j \to a_j'$ in $\cD$, and
\item the second component is
\[\begin{tikzcd}[column sep=huge]
(f_j)_* x_j \ar{r}{(f_j)_* p_j} & (f_j)_* y_j
\end{tikzcd} \inspace X_j a_j'.\] 
\end{itemize}
In the commutative diagram \cref{fjonefjstarpj}, as in \cref{ajxjfjoneajpfjxj},
\begin{itemize}
\item the top $\brb{f_j, 1_{(f_j)_* x_j}}$ is the chosen opcartesian lift of $\fl{(a_j,x_j)}{f_j}$ and
\item the bottom $\brb{f_j, 1_{(f_j)_* y_j}}$ is the chosen opcartesian lift of $\fl{(a_j,y_j)}{f_j}$ 
\end{itemize}
with respect to $U_{X_j} \cn \grod X_j \to \cD$.

Applying $F$ to the commutative diagrams \cref{fjonefjstarpj} for $j \in \{1,\ldots,n\}$ and using \cref{Fpreimagemor,Fangoneajpj,Fajxjtensorfj-i,Fajxjtensorfj-ii} yield the following commutative diagram in $\grod Z$.
\[\begin{tikzpicture}[xscale=5,yscale=1.3,vcenter]
\draw[0cell=.9]
(0,0) node (x11) {\brb{a,\phi_{\anga}\angx}}
(x11)++(1,0) node (x12) {\brb{a', \phi_{\angap} \ang{(f_j)_* x_j}}}
(x11)++(0,-1) node (x21) {\brb{a,\phi_{\anga}\angy}}
(x12)++(0,-1) node (x22) {\brb{a', \phi_{\angap} \ang{(f_j)_* y_j}}}
;
\draw[1cell=.9]  
(x11) edge node {\brb{f, 1_{f_* \phi_{\anga}\angx}}} (x12)
(x12) edge[transform canvas={xshift=-2em}] node {\brb{1_{a'}, \phi_{\angap} \ang{(f_j)_* p_j}}} (x22)
(x11) edge[transform canvas={xshift=1em}] node[swap] {\brb{1_a, \phi_{\anga}\angp}} (x21)
(x21) edge node {\brb{f, 1_{f_* \phi_{\anga}\angy}}} (x22)
;
\end{tikzpicture}\]
The second components of these two equal composites yield the morphism equality
\[f_* \phi_{\anga}\angp = \phi_{\angap} \ang{(f_j)_* p_j} \inspace Za'.\]
This proves the commutativity of \cref{phinaturalitydiagram} on morphisms.
\end{proof}

\subsection*{Unity and Additivity}

\cref{phiangafunctor,phinaturaltr} show that 
\[\begin{tikzcd}[column sep=large]
\bigotimes_{j=1}^n X_j \ar{r}{\phi} & Z,
\end{tikzcd}\]
with component functors $\phi_{\anga}$ in \cref{Fpreimagecomponent}, is a natural transformation as in \cref{phitensorxz}.  To show that $\phi$ is an additive natural transformation, next we check the two axioms in \cref{def:additivenattr}, starting with the unity axiom.

\begin{lemma}\label{phiunityaxiom}
The natural transformation \cref{Fpreimagecomponent}
\[\begin{tikzcd}[column sep=large]
\bigotimes_{j=1}^n X_j \ar{r}{\phi} & Z
\end{tikzcd}\]
satisfies the unity axiom \cref{additivenattrunity}.
\end{lemma}

\begin{proof}
With $\compi$ as in \cref{compnotation}, the unity axiom asserts the equalities in $Z\zero$
\[\left\{\begin{split}
\phi_{\anga \compi \zero} \big( \angx \compi (X_i^0*)\big) &= Z^0*\\
\phi_{\anga \compi \zero} \big( \angp \compi 1_{X_i^0*}\big) &= 1_{Z^0*}
\end{split}\right.\]
for objects $\angx$ and morphisms $\angp$ in $\prod_{j=1}^n X_j a_j$.  The unity axiom for objects is proved by the following equalities in $\grod Z$.
\begin{equation}\label{phiunitycomputation}
\left\{\begin{aligned}
&\phantom{==} (\zero, Z^0*) &&\\
&= F\bang{(a_j,x_j) \compi (\zero, X_i^0*)} & \phantom{M} & \text{by \cref{nlinearunity,groconunit}}\\
&= \scalebox{.8}{$\brb{a_1 \otimes \cdots \otimes \zero \otimes \cdots \otimes a_n, \phi_{\anga \compi \zero} (\angx \compi X_i^0*)}$} & \phantom{M} & \text{by \cref{Fpreimageobject,Fangajxjaftwox}}\\
&= \brb{\zero, \phi_{\anga \compi \zero} (\angx \compi X_i^0*)} & \phantom{M} & \text{by \cref{ringcataxiommultzero}}
\end{aligned}\right.
\end{equation}
The unity axiom for morphisms is proved by the computation \cref{phiunitycomputation} by  
\begin{itemize}
\item replacing $(\zero, Z^0*)$, $\ang{(a_j,x_j)}$, and $(\zero,X_i^0*)$ with, respectively, 
\[\brb{1_{\zero}, 1_{Z^0*}}, \quad \ang{(1_{a_j},p_j)}, \andspace \brb{1_{\zero}, 1_{X_i^0*}}\]
and
\item using \cref{groconidentity}
\[\brb{1_{\zero}, 1_{Z^0*}} = 1_{(\zero, Z^0*)} \inspace \grod Z\]
and the fact that $F$ preserves identity morphisms.
\end{itemize}
This finishes the proof of the unity axiom for $\phi$.
\end{proof}

\begin{lemma}\label{phiadditivityaxiom}
The natural transformation \cref{Fpreimagecomponent}
\[\begin{tikzcd}[column sep=large]
\bigotimes_{j=1}^n X_j \ar{r}{\phi} & Z
\end{tikzcd}\]
satisfies the additivity axiom \cref{additivityobjects}.
\end{lemma}

\begin{proof}
We first prove the additivity axiom on objects and then on morphisms.

\emph{Objects}.  Using the notation in \cref{additivenattrnotation,laplazaapp,angxxiprime,angxprime} for objects, we compute the (co)domain objects in $\grod Z$ of the $i$-th linearity constraint, $F^2_i$, as follows.
\begin{equation}\label{Ftwoicodomain}
\left\{\scalebox{.9}{$
\begin{aligned}
&F\ang{(a_j,x_j)} \gbox F\bang{(a_j,x_j) \compi (a_i',x_i')} &&\\
&\phantom{M} = \brb{a, \phi_{\anga}\angx} \gbox \brb{a', \phi_{\angap}\angxp} & \phantom{M} & \text{by \cref{Fpreimageobject,Fangajxjaftwox}}\\
&\phantom{M} = \left(a \oplus a' \scs Z^2_{a,a'} \brb{\phi_{\anga}\angx, \phi_{\angap}\angxp}\right) & \phantom{M} & \text{by \cref{gboxobjects}}\\
&F\lrang{(a_j,x_j) \compi \big((a_i,x_i) \gbox (a_i',x_i') \big)} &&\\
&\phantom{M} = F\lrang{(a_j,x_j) \compi \brb{a_i'',x_i''}} & \phantom{M} & \text{by \cref{gboxobjects}}\\
&\phantom{M} = \brb{a'', \phi_{\angapp} \angxpp} & \phantom{M} & \text{by \cref{Fpreimageobject,Fangajxjaftwox}}
\end{aligned}$}\right.
\end{equation}
By the constraint lift axiom \cref{def:opcartnlinearfunctor-iv} in \cref{def:opcartnlinearfunctor}, this component of $F^2_i$ is the chosen opcartesian lift of the fore-lift
\[\bfl{F\ang{(a_j,x_j)} \gbox F\bang{(a_j,x_j) \compi (a_i',x_i')} }{\lapinv}\]
with respect to $U_Z \cn \grod Z \to \cD$.  By \cref{fonechosenoplift,Ftwoicodomain} this chosen opcartesian lift, which is $F^2_i$, is also given by the following morphism in $\grod Z$.
\begin{equation}\label{Ftwoilapinvone}
\begin{tikzpicture}[xscale=5,yscale=1,baseline={(x11.base)}]
\draw[0cell=.85]
(0,0) node (x11) {\left(a \oplus a' \scs Z^2_{a,a'} \brb{\phi_{\anga}\angx, \phi_{\angap}\angxp}\right)}
(x11)++(1,0) node (x12) {\brb{a'', \phi_{\angapp} \angxpp}};
\draw[1cell=.9]  
(x11) edge node {(\lapinv,1)} (x12);
\end{tikzpicture}
\end{equation}
Its second entry is the identity morphism in $Za''$
\[\begin{tikzpicture}[xscale=4,yscale=1,vcenter]
\draw[0cell=.9]
(0,0) node (x11) {\lapinv_* Z^2_{a,a'} \brb{\phi_{\anga}\angx, \phi_{\angap}\angxp}}
(x11)++(1,0) node (x12) {\phi_{\angapp} \angxpp.};
\draw[1cell=.9]  
(x11) edge node {1} (x12);
\end{tikzpicture}\]
Applying $\lap_*$ yields the equality in $Z(a \oplus a')$
\begin{equation}\label{lapstarphiappxppztwo}
\lap_* \phi_{\angapp} \angxpp = Z^2_{a,a'} \brb{\phi_{\anga}\angx, \phi_{\angap}\angxp}.
\end{equation}
This proves the additivity axiom \cref{additivityobjects} for objects.

\emph{Morphisms}.  For the additivity axiom \cref{additivityobjects} for morphisms, suppose given morphisms
\[\begin{split}
\angp = \bang{p_j \cn x_j \to y_j}_{j=1}^n & \in \txprod_{j=1}^n X_j a_j \andspace\\
p_i' \cn x_i' \to y_i' & \in X_ia_i'.
\end{split}\]
We extend the notation in \cref{angxprime,additivityobjects} to morphisms by replacing each occurrence of $x$ with $p$.  By \cref{Ftwoicodomain,Ftwoilapinvone}, the naturality diagram of $F^2_i$ for these morphisms is the following commutative diagram in $\grod Z$.
\[\begin{tikzpicture}[xscale=5,yscale=1.4,vcenter]
\draw[0cell=.85]
(0,0) node (x11) {\left(a \oplus a' \scs Z^2_{a,a'} \brb{\phi_{\anga}\angx, \phi_{\angap}\angxp}\right)}
(x11)++(1,0) node (x12) {\brb{a'', \phi_{\angapp} \angxpp}}
(x11)++(0,-1) node (x21) {\left(a \oplus a' \scs Z^2_{a,a'} \brb{\phi_{\anga}\angy, \phi_{\angap}\angyp}\right)}
(x12)++(0,-1) node (x22) {\brb{a'', \phi_{\angapp} \angypp}}
;
\draw[1cell=.8]  
(x11) edge node {(\lapinv,1)} (x12)
(x12) edge[transform canvas={xshift=-2em}] node {\brb{1, \phi_{\angapp} \angppp}} (x22)
(x11) edge[transform canvas={xshift=4em}] node[swap] {\brb{1, Z^2_{a,a'} \brb{\phi_{\anga}\angp, \phi_{\angap}\angpp}}} (x21)
(x21) edge node {(\lapinv,1)} (x22)
;
\end{tikzpicture}\]
The second component of this commutative diagram is the morphism equality in $Za''$
\[\lapinv_* Z^2_{a,a'} \brb{\phi_{\anga}\angp, \phi_{\angap}\angpp} = \phi_{\angapp} \angppp.\]
Applying $\lap_*$ yields the morphism equality in $Z(a \oplus a')$
\[\lap_* \phi_{\angapp} \angppp = Z^2_{a,a'} \brb{\phi_{\anga}\angp, \phi_{\angap}\angpp}.\]
This proves the additivity axiom \cref{additivityobjects} for morphisms.
\end{proof}

\subsection*{Preimage of $F$}

By \cref{phiangafunctor,phinaturaltr,phiunityaxiom,phiadditivityaxiom}, 
\[\begin{tikzcd}[column sep=large]
\bang{(X_j, X_j^2, X_j^0)}_{j=1}^n \ar{r}{\phi} & (Z,Z^2,Z^0)
\end{tikzcd}\]
in \cref{Fpreimagephi} is an additive natural transformation (\cref{def:additivenattr}).  Next we show that its Grothendieck construction, in the sense of \cref{grodphifunctor}, is equal to $F$.

\begin{lemma}\label{grodobjectfullfphi}\index{Grothendieck construction!object surjectivity}
In the context of \cref{def:grodobjectfull}, there is an equality 
\[\brb{\grod\phi, \ang{(\grod\phi)^2_j}_{j=1}^n} = \brb{F, \ang{F^2_j}_{j=1}^n}\]
of strong $n$-linear functors.
\end{lemma}

\begin{proof}
We must prove the following two statements.
\begin{enumerate}
\item\label{grodobjectfullfphi-i} There is an equality of functors
\[\grod\phi = F \cn \txprod_{j=1}^n (\grod X_j) \to \grod Z.\]
\item\label{grodobjectfullfphi-ii} 
In the context of \cref{Ftwoicodomain,Ftwoilapinvone,lapstarphiappxppztwo}, the following two morphisms in $\grod Z$ are equal.
\[\begin{tikzpicture}[xscale=5,yscale=1,baseline={(x11.base)}]
\draw[0cell=1]
(0,0) node (x11) {\brb{a \oplus a', \lap_* \phi_{\angapp}\angxpp}}
(x11)++(1,0) node (x12) {\brb{a'', \phi_{\angapp} \angxpp}};
\draw[1cell=.8]  
(x11) edge[transform canvas={yshift=.5ex}] node {F^2_i = (\lapinv,1)} (x12)
(x11) edge[transform canvas={yshift=-.5ex}] node[swap] {(\grod\phi)^2_i} (x12);
\end{tikzpicture}\]
\end{enumerate}
Once \eqref{grodobjectfullfphi-i} is proved, \eqref{grodobjectfullfphi-ii} follows from the definition \cref{grodphilinearity} of $(\grod\phi)^2_i$.  We prove \eqref{grodobjectfullfphi-i} first on objects and then on morphisms.

\emph{Objects}.  For each $n$-tuple of objects 
\[\ang{(a_j,x_j)} \in \txprod_{j=1}^n (\grod X_j),\]
the following equalities in $\grod Z$ prove \eqref{grodobjectfullfphi-i} on objects.
\[\begin{aligned}
(\grod\phi)\ang{(a_j,x_j)}
&= \brb{a,\phi_{\anga}\angx} & \phantom{M} & \text{by \cref{grodphiangajxj}}\\
&= F\ang{(a_j,x_j)} & \phantom{M} & \text{by \cref{Fpreimageobject,Fangajxjaftwox}}
\end{aligned}\]

\emph{Morphisms}.  Suppose given an $n$-tuple of morphisms
\[\bang{(f_j,p_j) \cn (a_j,x_j) \to (b_j,y_j)} \in \txprod_{j=1}^n (\grod X_j).\]
We must show that its images under $\grod\phi$ and $F$ are equal.
By definition \cref{groconcomposition} each morphism $(f_j,p_j)$ in $\grod X_j$ factors as follows.
\begin{equation}\label{phijfjfactors}
\begin{tikzpicture}[xscale=3.5,yscale=1,vcenter]
\draw[0cell=.9]
(0,0) node (x11) {(a_j,x_j)}
(x11)++(1,0) node (x12) {(b_j,y_j)}
(x11)++(.5,-1) node (x21) {\brb{b_j, (f_j)_* x_j}}
;
\draw[1cell=.9]  
(x11) edge node {(f_j,p_j)} (x12)
(x11) edge[transform canvas={xshift=-1em}] node[swap,pos=.3] {(f_j \scs 1_{(f_j)_* x_j})} (x21)
(x21) edge[transform canvas={xshift=1em}] node[swap,pos=.7] {(1_{b_j} \scs p_j)} (x12)
;
\end{tikzpicture}
\end{equation}
Now we observe the following for the first morphism $\brb{f_j,1_{(f_j)_* x_j}}$ in \cref{phijfjfactors}.
\begin{itemize}
\item By definition \cref{fonechosenoplift} $(f_j,1)$ is the chosen opcartesian lift of the fore-lift $\fl{(a_j,x_j)}{f_j}$ with respect to $U_{X_j} \cn \grod X_j \to \cD$.
\item Since $F$ preserves chosen opcartesian lifts (\cref{def:opcartnlinearfunctor} \cref{def:opcartnlinearfunctor-iii}), the image
\[\begin{tikzpicture}[xscale=3.8,yscale=1,baseline={(x11.base)}]
\draw[0cell=1]
(0,0) node (x11) {F\ang{(a_j,x_j)}}
(x11)++(1,0) node (x12) {F\bang{\brb{b_j, (f_j)_* x_j}}}
;
\draw[1cell=.9]  
(x11) edge node {F\ang{(f_j,1)}} (x12);
\end{tikzpicture} \inspace \grod Z\]
is the chosen opcartesian lift of the fore-lift
\[\bfl{F\ang{(a_j,x_j)}}{\txotimes_{j=1}^n f_j} = \bfl{\brb{a,\phi_{\anga}\angx}}{f}\]
with respect to $U_Z \cn \grod Z \to \cD$.  By definition \cref{fonechosenoplift} this chosen opcartesian lift, which is $F\ang{(f_j,1)}$, is also given by 
\begin{equation}\label{fonefstarphi}
\begin{tikzpicture}[xscale=4.5,yscale=1,baseline={(x11.base)}]
\draw[0cell=.9]
(0,0) node (x11) {\brb{a,\phi_{\anga}\angx}}
(x11)++(1,0) node (x12) {\brb{b,f_*\phi_{\anga}\angx}.}
;
\draw[1cell=.9]  
(x11) edge node {\brb{f, 1_{f_*\phi_{\anga}\angx}}} (x12);
\end{tikzpicture}
\end{equation}
\end{itemize} 
The following morphism equalities in $\grod Z$ prove \eqref{grodobjectfullfphi-i} on morphisms. 
\[\begin{aligned}
&\phantom{==} (\grod\phi)\ang{(f_j,p_j)} &&\\
&= \brb{f,\phi_{\angb}\angp} & \phantom{M} & \text{by \cref{grodphiangfjpj}}\\
&= \brb{1_b,\phi_{\angb}\angp} \circ \brb{f,1_{f_*\phi_{\anga}\angx}} & \phantom{M} & \text{by \cref{groconcomposition}}\\
&= F\ang{(1_{b_j},p_j)} \circ F\bang{\brb{f_j,1_{(f_j)_* x_j}}} & \phantom{M} & \text{by \cref{Fpreimagemor,Fangoneajpj,fonefstarphi}}\\
&= F\ang{(f_j,p_j)} & \phantom{M} & \text{by \cref{phijfjfactors}}
\end{aligned}\]
This finishes the proof of \eqref{grodobjectfullfphi-i}.
\end{proof}

\section{Morphism Fullness and the Main Theorem}
\label{sec:morphismfullness}

We continue to use the notation in \cref{zangxjdpluscat,angxgrodxux,angajtensoraj,angxjobject,angpjxjyj} with $n>0$.  In \cref{grodarityninjobj,grodarityninjmor,grodobjectfullfphi} we observe that the component functors (\cref{grodfactorarityn}) 
\[\begin{tikzpicture}[xscale=3.8,yscale=1.2,vcenter]
\draw[0cell=.9]
(0,0) node (dc) {\DCat\scmap{\angX;Z}}
(dc)++(1,0) node (pf) {\pfibd\scmap{\ang{U_X};U_Z}};
\draw[1cell=.9]  
(dc) edge node {\grod} (pf);
\end{tikzpicture}\] 
are (i) injective on objects and morphism sets and (ii) surjective on objects.  In this section we show that these component functors are also surjective on morphism sets.  

Recall the following.
\begin{itemize}
\item $\DCat\scmap{\angX;Z}$ has additive natural transformations as objects and additive modifications as morphisms (\cref{def:additivenattr,def:additivemodification}).
\item $\pfibd\scmap{\ang{U_X};U_Z}$ has opcartesian strong $n$-linear functors as objects and opcartesian $n$-linear transformations as morphisms (\cref{def:opcartnlinearfunctor,def:opcartnlineartr}).
\end{itemize}
Here is an outline of this section.
\begin{itemize}
\item In \cref{def:grodmorphismfull} we associate to each opcartesian $n$-linear transformation $\omega \cn \grod\phi \to \grod\psi$ between the Grothendieck constructions of two additive natural transformations the data of an additive modification $\Phi \cn \phi \to \psi$.
\item \cref{Phianganattransformation,Phiphipsimodification} show that $\Phi$ is a modification.
\item \cref{Phiunityaxiom,Phiadditivityaxiom} show that $\Phi$ is an additive modification.
\item \cref{grodPhiomega} shows that $\grod\Phi$ and $\omega$ are equal as opcartesian $n$-linear transformations.  This proves that the multimorphism functors of $\grod$ are surjective on morphism sets.
\item Using the Multiequivalence \cref{thm:multiwhitehead} and the lemmas in this chapter, \cref{thm:dcatpfibdeq} proves that
\[\begin{tikzpicture}[xscale=2.75,yscale=1.2,baseline={(dc.base)}]
\draw[0cell=.9]
(0,0) node (dc) {\DCat}
(dc)++(1,0) node (pf) {\pfibd};
\draw[1cell=.9]  
(dc) edge node {\grod} (pf);
\end{tikzpicture}\]
is a non-symmetric $\Cat$-multiequivalence.  If, moreover, $\betate = 1$ in $\cD$, then $\grod$ is a $\Cat$-multiequivalence.
\end{itemize}

\begin{definition}\label{def:grodmorphismfull}
Suppose given additive natural transformations
\[\begin{tikzcd}[column sep=large]
\bang{(X_j, X_j^2, X_j^0)}_{j=1}^n \ar[shift left]{r}{\phi} \ar[shift right]{r}[swap]{\psi} & (Z,Z^2,Z^0)
\end{tikzcd}\]
and an opcartesian $n$-linear transformation
\begin{equation}\label{Omegaopcartesian}
\begin{tikzcd}[column sep=large]
\brb{\grod\phi, \ang{(\grod\phi)^2_j}_{j=1}^n} \ar{r}{\omega} 
& \brb{\grod\psi, \ang{(\grod\psi)^2_j}_{j=1}^n}
\end{tikzcd}
\end{equation}
with the (co)domain given by \cref{grodphifunctor}.  The structure $\omega$ is given by an $n$-linear natural transformation (\cref{def:nlineartransformation})
\begin{equation}\label{Omegagrodphipsi}
\begin{tikzpicture}[xscale=1,yscale=1,baseline={(x1.base)}]
\def\a{25}
\draw[0cell=.9]
(0,0) node (x1) {\txprod_{j=1}^n (\grod X_j)}
(x1)++(.75,0) node (x2) {\phantom{D}}
(x2)++(2.5,0) node (x3) {\phantom{D}}
(x3)++(.23,0.02) node (x4) {\grod Z}
;
\draw[1cell=.8]  
(x2) edge[bend left=\a] node[pos=.5] {\brb{\grod\phi, \ang{(\grod\phi)^2_j}_{j=1}^n}} (x3)
(x2) edge[bend right=\a] node[swap,pos=.5] {\brb{\grod\psi, \ang{(\grod\psi)^2_j}_{j=1}^n}} (x3)
;
\draw[2cell] 
node[between=x2 and x3 at .45, rotate=-90, 2label={above,\omega}] {\Rightarrow}
;
\end{tikzpicture}
\end{equation}
that satisfies the transformation projection axiom \cref{thetapzeropjtensor}.  Define the data of an additive modification
\begin{equation}\label{Omegapreimage}
\begin{tikzpicture}[xscale=1,yscale=1,baseline={(x1.base)}]
\def\a{25}
\draw[0cell=.9]
(0,0) node (x1) {\bang{(X_j, X_j^2, X_j^0)}_{j=1}^n}
(x1)++(1,0) node (x2) {\phantom{X}}
(x2)++(2,0) node (x3) {\phantom{X}}
(x3)++(.65,0) node (x4) {(Z,Z^2,Z^0)}
;
\draw[1cell=.9]  
(x2) edge[bend left=\a] node[pos=.5] {\phi} (x3)
(x2) edge[bend right=\a] node[swap,pos=.5] {\psi} (x3)
;
\draw[2cell]
node[between=x2 and x3 at .45, rotate=-90, 2label={above,\Phi}] {\Rightarrow}
;
\end{tikzpicture}
\end{equation}
with $\anga$-component natural transformations \cref{dcatnarymodcomponent}
\begin{equation}\label{Omegapreimagecomponent}
\begin{tikzpicture}[xscale=2.5,yscale=2,baseline={(x1.base)}]
\draw[0cell=.9]
(0,0) node (x1) {\txprod_{j=1}^n X_j a_j}
(x1)++(.17,.03) node (x2) {\phantom{\textstyle \bigotimes Z}}
(x2)++(1,0) node (x3) {Za}
;
\draw[1cell=.9]  
(x2) edge[bend left, shorten <=-.5ex] node[pos=.5] {\phi_{\anga}} (x3)
(x2) edge[bend right, shorten <=-.5ex] node[swap,pos=.5] {\psi_{\anga}} (x3)
;
\draw[2cell]
node[between=x2 and x3 at .4, rotate=-90, 2label={above,\Phi_{\anga}}] {\Rightarrow}
;
\end{tikzpicture}
\end{equation}
as follows.  For each $n$-tuple of objects $\angx \in \txprod_{j=1}^n X_j a_j$, define the morphism
\begin{equation}\label{Phiangaangx}
(\Phi_{\anga})_{\angx} = \omega_{\angx} \cn \phi_{\anga}\angx \to \psi_{\anga}\angx \inspace Za
\end{equation}
where $\omega_{\angx}$ is the second entry in the $\ang{(a_j,x_j)}$-component morphism
\begin{equation}\label{Omegaangajxj}
\omega_{\ang{(a_j,x_j)}} = \brb{1_a,\omega_{\angx}} \cn (\grod\phi)\ang{(a_j,x_j)} \to (\grod\psi)\ang{(a_j,x_j)}
\end{equation}
in $\grod Z$.  This finishes the definition of $\Phi$ in \cref{Omegapreimage}.
\end{definition}

\begin{explanation}\label{expl:grodmorphismfull}
In \cref{Omegaangajxj} the (co)domain objects are
\[\begin{split}
(\grod\phi)\ang{(a_j,x_j)} &= \brb{a, \phi_{\anga}\angx} \andspace\\
(\grod\psi)\ang{(a_j,x_j)} &= \brb{a, \psi_{\anga}\angx}
\end{split}\]
in $\grod Z$ by definition \cref{grodphiangajxj}.  The transformation projection axiom \cref{thetapzeropjtensor} for $\omega$ in \eqref{Omegaopcartesian} means that
\[U_Z \omega_{\ang{(a_j,x_j)}} = 1_{\otimes_{j=1}^n U_{X_j} (a_j,x_j)} = 1_a \inspace \cD.\]
So the first entry in the morphism
\[\begin{tikzcd}[column sep=huge]
\brb{a, \phi_{\anga}\angx} \ar{r}{\omega_{\ang{(a_j,x_j)}}} & \brb{a, \psi_{\anga}\angx}
\end{tikzcd} \inspace \grod Z\]
is the identity morphism $1_a \in \cD$.  This is the meaning of \cref{Omegaangajxj}.
\end{explanation}

In the rest of this section we show that $\Phi \cn \phi \to \psi$ in \cref{Omegapreimage} is an additive modification such that
\[\grod\Phi = \omega\]
as $n$-linear natural transformations, with $\grod\Phi$ given by \cref{grodPhinlinear}.

\subsection*{Naturality and the Modification Axiom}

To show that $\Phi$ is an additive modification, first we check that its components are natural transformations.

\begin{lemma}\label{Phianganattransformation}
In \cref{Omegapreimagecomponent} each component 
\[\begin{tikzcd}[column sep=large]
\phi_{\anga} \ar{r}{\Phi_{\anga}} & \psi_{\anga}
\end{tikzcd}\]
is a natural transformation between functors $\prod_{j=1}^n X_ja_j \to Za$.
\end{lemma}

\begin{proof}
For each $n$-tuple of morphisms
\[\angp = \bang{p_j \cn x_j \to y_j}_{j=1}^n \in \txprod_{j=1}^n X_j a_j\]
the desired naturality diagram for $\Phi_{\anga}$ is the following diagram in $Za$.
\begin{equation}\label{Phianganaturalitydiag}
\begin{tikzpicture}[xscale=3,yscale=1.2,vcenter]
\draw[0cell=.9]
(0,0) node (x11) {\phi_{\anga}\angx}
(x11)++(1,0) node (x12) {\psi_{\anga}\angx}
(x11)++(0,-1) node (x21) {\phi_{\anga}\angy}
(x12)++(0,-1) node (x22) {\psi_{\anga}\angy}
;
\draw[1cell=.9]  
(x11) edge node {(\Phi_{\anga})_{\angx}} (x12)
(x21) edge node {(\Phi_{\anga})_{\angy}} (x22)
(x11) edge node[swap] {\phi_{\anga}\angp} (x21)
(x12) edge node {\psi_{\anga}\angp} (x22)
;
\end{tikzpicture}
\end{equation}
To show that \cref{Phianganaturalitydiag} is commutative, consider the $n$-tuple of morphisms
\[\bang{(1_{a_j},p_j) \cn (a_j,x_j) \to (a_j,y_j)}_{j=1}^n \in \txprod_{j=1}^n (\grod X_j).\]
The naturality of $\omega \cn \grod\phi \to \grod\psi$ in \cref{Omegagrodphipsi} yields the following commutative diagram in $\grod Z$.
\[\begin{tikzpicture}[xscale=4.3,yscale=1.3,vcenter]
\draw[0cell=.9]
(0,0) node (x11) {(\grod\phi)\ang{(a_j,x_j)}}
(x11)++(1,0) node (x12) {(\grod\psi)\ang{(a_j,x_j)}}
(x11)++(0,-1) node (x21) {(\grod\phi)\ang{(a_j,y_j)}}
(x12)++(0,-1) node (x22) {(\grod\psi)\ang{(a_j,y_j)}}
;
\draw[1cell=.85]  
(x11) edge node {\omega_{\ang{(a_j,x_j)}}} (x12)
(x21) edge node {\omega_{\ang{(a_j,y_j)}}} (x22)
(x11) edge[transform canvas={xshift=1em}] node[swap] {(\grod\phi)\ang{(1_{a_j},p_j)}} (x21)
(x12) edge[transform canvas={xshift=-1em}] node {(\grod\psi)\ang{(1_{a_j},p_j)}} (x22)
;
\end{tikzpicture}\]
By the definitions 
\begin{itemize}
\item \cref{grodphiangajxj,grodphiangfjpj} of $\grod\phi$ on, respectively, objects and morphisms and
\item \cref{Phiangaangx,Omegaangajxj} of $\Phi_{\anga}$, 
\end{itemize}
the previous commutative diagram is equal to the following diagram.
\[\begin{tikzpicture}[xscale=4.3,yscale=1.3,vcenter]
\draw[0cell=.9]
(0,0) node (x11) {\brb{a,\phi_{\anga}\angx}}
(x11)++(1,0) node (x12) {\brb{a,\psi_{\anga}\angx}}
(x11)++(0,-1) node (x21) {\brb{a,\phi_{\anga}\angy}}
(x12)++(0,-1) node (x22) {\brb{a,\psi_{\anga}\angy}}
;
\draw[1cell=.8]  
(x11) edge node {\brb{1_a,(\Phi_{\anga})_{\angx}}} (x12)
(x21) edge node {\brb{1_a,(\Phi_{\anga})_{\angy}}} (x22)
(x11) edge[transform canvas={xshift=1em}] node[swap] {\brb{1_a, \phi_{\anga}\angp}} (x21)
(x12) edge[transform canvas={xshift=-1em}] node {\brb{1_a, \psi_{\anga}\angp}} (x22)
;
\end{tikzpicture}\]
The second component of this commutative diagram yields the desired diagram \cref{Phianganaturalitydiag}.
\end{proof}

Next we show that $\Phi$ is a modification.

\begin{lemma}\label{Phiphipsimodification}
In \cref{Omegapreimage} $\Phi \cn \phi \to \psi$ is a modification.
\end{lemma}

\begin{proof}
The desired modification axiom \cref{dcatnarymodaxiom} states that, for each $n$-tuple of objects $\angx \in \prod_{j=1}^n X_ja_j$, the following morphism equality holds in $Za'$.
\begin{equation}\label{fstaromegaangx}
f_* (\Phi_{\anga})_{\angx} = (\Phi_{\angap})_{\ang{(f_j)_* x_j}} \cn f_* \phi_{\anga}\angx \to \psi_{\angap}\ang{(f_j)_* x_j}
\end{equation}
To prove the equality in \cref{fstaromegaangx}, consider the $n$-tuple of morphisms
\[\bang{\brb{f_j,1_{(f_j)_* x_j}} \cn (a_j,x_j) \to \brb{a_j', (f_j)_* x_j}} \in \txprod_{j=1}^n (\grod X_j).\]
The naturality of $\omega \cn \grod\phi \to \grod\psi$ in \cref{Omegagrodphipsi} yields the following commutative diagram in $\grod Z$.
\[\begin{tikzpicture}[xscale=5,yscale=1.4,vcenter]
\draw[0cell=.85]
(0,0) node (x11) {(\grod\phi)\ang{(a_j,x_j)}}
(x11)++(1,0) node (x12) {(\grod\phi)\bang{\brb{a_j',(f_j)_* x_j}}}
(x11)++(0,-1) node (x21) {(\grod\psi)\ang{(a_j,x_j)}}
(x12)++(0,-1) node (x22) {(\grod\psi)\bang{\brb{a_j',(f_j)_* x_j}}}
;
\draw[1cell=.85]  
(x11) edge node {(\grod\phi)\ang{(f_j,1)}} (x12)
(x21) edge node {(\grod\psi)\ang{(f_j,1)}} (x22)
(x11) edge[transform canvas={xshift=1em}] node[swap] {\omega_{\ang{(a_j,x_j)}}} (x21)
(x12) edge[transform canvas={xshift=-1em}] node {\omega_{\brb{a_j',(f_j)_* x_j}}} (x22)
;
\end{tikzpicture}\]
By \cref{grodphiangajxj,grodphiangfjpj,Phiangaangx,Omegaangajxj} the previous commutative diagram is equal to the following diagram.
\begin{equation}\label{omeganatfjone}
\begin{tikzpicture}[xscale=4.5,yscale=1.3,vcenter]
\draw[0cell=.9]
(0,0) node (x11) {\brb{a,\phi_{\anga}\angx}}
(x11)++(1,0) node (x12) {\brb{a',\phi_{\angap}\ang{(f_j)_* x_j}}}
(x11)++(0,-1) node (x21) {\brb{a,\psi_{\anga}\angx}}
(x12)++(0,-1) node (x22) {\brb{a',\psi_{\angap}\ang{(f_j)_* x_j}}}
;
\draw[1cell=.8]
(x11) edge node {\brb{f,\phi_{\angap} \ang{1}}} (x12)
(x21) edge node {\brb{f,\psi_{\angap} \ang{1}}} (x22)
(x11) edge[transform canvas={xshift=1em}] node[swap] {\brb{1_a, (\Phi_{\anga})_{\angx}}} (x21)
(x12) edge[transform canvas={xshift=-2em}] node {\brb{1_{a'}, (\Phi_{\angap})_{\ang{(f_j)_* x_j}}}} (x22)
;
\end{tikzpicture}
\end{equation}
By the functoriality of $\phi_{\angap}$ and $\psi_{\angap}$, there are identity morphisms
\[\begin{split}
\phi_{\angap} \ang{1} &= 1_{\phi_{\angap}\ang{(f_j)_* x_j}} \andspace\\ 
\psi_{\angap} \ang{1} &= 1_{\psi_{\angap}\ang{(f_j)_* x_j}}
\end{split}\]
in $Za'$.  Thus the second component of the commutative diagram \cref{omeganatfjone} yields the desired equality \cref{fstaromegaangx}.
\end{proof}

\subsection*{Unity and Additivity}

Next we show that $\Phi$ is an \emph{additive} modification (\cref{def:additivemodification}).

\begin{lemma}\label{Phiunityaxiom}
$\Phi \cn \phi \to \psi$ in \cref{Omegapreimage} satisfies the unity axiom \cref{additivemodunity}.
\end{lemma}

\begin{proof}
With $\compi$ as in \cref{compnotation}, there are morphism equalities in $\grod Z$ as follows.
\[\begin{aligned}
&\phantom{==} \brb{1_{\zero}, (\Phi_{\anga \compi \zero})_{\angx \compi (X_i^0*)}} &&\\
&= \omega_{\ang{(a_j,x_j)} \compi (\zero,X_i^0*)} & \phantom{M} & \text{by \cref{ringcataxiommultzero,Phiangaangx,Omegaangajxj}}\\
&= 1_{(\zero,Z^0*)} & \phantom{M} & \text{by \cref{ntransformationunity,groconunit}}\\
&= \brb{1_{\zero}, 1_{Z^0*}} & \phantom{M} & \text{by \cref{groconidentity}}
\end{aligned}\]
The morphism equality
\[(\Phi_{\anga \compi \zero})_{\angx \compi (X_i^0*)} = 1_{Z^0*} \inspace Z\zero\]
is the unity axiom for $\Phi$.
\end{proof}

\begin{lemma}\label{Phiadditivityaxiom}
$\Phi \cn \phi \to \psi$ in \cref{Omegapreimage} satisfies the additivity axiom \cref{additivitymod}.
\end{lemma}

\begin{proof}
In the context of \cref{grodphiconstraintnot,grodphiconstraintdom,grodphiconstraintcod,grodphilinearity,Phiangaangx,Omegaangajxj}, the constraint compatibility axiom \cref{eq:monoidal-in-each-variable} for $\omega \cn \grod\phi \to \grod\psi$ in \cref{Omegagrodphipsi} is the following commutative diagram in $\grod Z$.
\[\begin{tikzpicture}[xscale=5.2,yscale=1.4,vcenter]
\draw[0cell=.9]
(0,0) node (x11) {\big(a \oplus a' \scs Z^2_{a,a'} \brb{\phi_{\anga}\angx, \phi_{\angap}\angxp} \big)}
(x11)++(1,0) node (x12) {\brb{a'', \phi_{\angapp}\angxpp}}
(x11)++(0,-1) node (x21) {\big(a \oplus a' \scs Z^2_{a,a'} \brb{\psi_{\anga}\angx, \psi_{\angap}\angxp} \big)}
(x12)++(0,-1) node (x22) {\brb{a'', \psi_{\angapp}\angxpp}}
;
\draw[1cell=.8]
(x11) edge node {(\lapinv,1)} (x12)
(x21) edge node {(\lapinv,1)} (x22)
(x11) edge[transform canvas={xshift=5em}] node[swap] {\left(1_{a \oplus a'} \scs Z^2_{a,a'} \brb{(\Phi_{\anga})_{\angx}, (\Phi_{\angap})_{\angxp}} \right)} (x21)
(x12) edge[transform canvas={xshift=-2.3em}] node {\brb{1_{a''}, (\Phi_{\angapp})_{\angxpp}}} (x22)
;
\end{tikzpicture}\]
The second entry of this commutative diagram yields the morphism equality
\[(\Phi_{\angapp})_{\angxpp} = \lapinv_* Z^2_{a,a'} \brb{(\Phi_{\anga})_{\angx}, (\Phi_{\angap})_{\angxp}} \inspace Za''.\]
Applying $\lap_*$ yields the equality
\[\lap_* (\Phi_{\angapp})_{\angxpp} = Z^2_{a,a'} \brb{(\Phi_{\anga})_{\angx}, (\Phi_{\angap})_{\angxp}} \inspace Z(a \oplus a').\]
This proves the additivity axiom for $\Phi$.
\end{proof}

\subsection*{Preimage of $\omega$}

\cref{Phiphipsimodification,Phiunityaxiom,Phiadditivityaxiom} show that $\Phi \cn \phi \to \psi$ in \cref{Omegapreimage} is an additive modification (\cref{def:additivemodification}).  By \cref{grodPhinlinear,grodfactorarityn} its Grothendieck construction $\grod\Phi$ is an opcartesian $n$-linear transformation.

\begin{lemma}\label{grodPhiomega}\index{Grothendieck construction!morphism surjectivity}
In the context of \cref{def:grodmorphismfull}, there is an equality 
\[\grod\Phi = \omega \cn \grod\phi \to \grod\psi\]
of opcartesian $n$-linear transformations.
\end{lemma}

\begin{proof}
Opcartesian $n$-linear transformations are $n$-linear natural transformations with the extra property \cref{thetapzeropjtensor} but no extra structure (\cref{expl:opcartnlineartr} \eqref{expl:opcartnlineartr-i}).  Thus it suffices to show that $\grod\Phi$ and $\omega$ are equal as natural transformations.  For each $n$-tuple of objects
\[\ang{(a_j,x_j)} \in \txprod_{j=1}^n (\grod X_j)\]
there are equalities in $\grod Z$ as follows.
\[\begin{aligned}
(\grod\Phi)_{\ang{(a_j,x_j)}} 
&= \brb{1_a, (\Phi_{\anga})_{\angx}} & \phantom{M} & \text{by \cref{grodPhiajxj}}\\
&= \omega_{\ang{(a_j,x_j)}} & \phantom{M} & \text{by \cref{Phiangaangx,Omegaangajxj}}
\end{aligned}\]
This proves $\grod\Phi = \omega$ as natural transformations.
\end{proof}

\subsection*{Multiequivalence}

We are now ready for the main result of this chapter.

\begin{theorem}\label{thm:dcatpfibdeq}\index{Grothendieck construction!non-symmetric Cat-multiequivalence@non-symmetric $\Cat$-multiequivalence}\index{multiequivalence!Cat-@$\Cat$-!Grothendieck construction}\index{Cat-multiequivalence@$\Cat$-multiequivalence!Grothendieck construction}
Suppose $\cD$ is a small tight bipermutative category.
\begin{enumerate}[label=(\roman*)]
\item\label{thm:dcatpfibdeq-i}
The non-symmetric $\Cat$-multifunctor in \cref{groddcatpfibd}
\[\begin{tikzpicture}[xscale=2.75,yscale=1.2,baseline={(dc.base)}]
\draw[0cell=.9]
(0,0) node (dc) {\DCat}
(dc)++(1,0) node (pf) {\pfibd};
\draw[1cell=.9]  
(dc) edge node {\grod} (pf);
\end{tikzpicture}\]
is a non-symmetric $\Cat$-multiequivalence.
\item\label{thm:dcatpfibdeq-ii} 
If the multiplicative braiding $\betate$ in $\cD$ is the identity natural transformation, then $\grod$ is a $\Cat$-multiequivalence.
\end{enumerate} 
\end{theorem}

\begin{proof}
We have proved the following regarding $\grod$.
\begin{itemize}
\item It is essentially surjective on objects (\cref{grodesssurjective}).
\item It is an isomorphism on arity 0 multimorphism categories (\cref{grodfactorarityzero}).
\item Between positive arity multimorphism categories, $\grod$ is
\begin{itemize}
\item injective on objects (\cref{grodarityninjobj}),
\item injective on morphism sets (\cref{grodarityninjmor}), 
\item surjective on objects (\cref{grodobjectfullfphi}), and
\item surjective on morphism sets (\cref{grodPhiomega}). 
\end{itemize}
\end{itemize}
Assertion \cref{thm:dcatpfibdeq-i} now follows from the Multiequivalence \cref{thm:multiwhitehead}.

Moreover, under the assumption that $\betate = 1$, $\grod$ is a $\Cat$-multifunctor between $\Cat$-multicategories by \cref{thm:dcatpfibd} \cref{thm:dcatpfibd-iii}.  Thus assertion \cref{thm:dcatpfibdeq-ii} also follows from \cref{thm:multiwhitehead} and the lemmas stated in the previous paragraph.
\end{proof}

\begin{example}\label{ex:dcatgrodpfibdeq}
Suppose $\cD$ is any one of the small tight bipermutative categories
\[\Finsk,\quad \Fset, \quad \Fskel, \andspace \cA\]
in, respectively, \cref{ex:finsk,ex:Fset,ex:Fskel,ex:mandellcategory}.  Then $\pfibd$ is a non-symmetric $\Cat$-multicategory and \emph{not} a $\Cat$-multicategory by \cref{ex:pfibdcatmulticat}.  The commutative diagram of non-symmetric $\Cat$-multifunctors
\[\begin{tikzpicture}[xscale=2.75,yscale=1.2,vcenter]
\draw[0cell=.9]
(0,0) node (dc) {\DCat}
(dc)++(1,0) node (pc) {\permcatsg}
(pc)++(0,1) node (pf) {\pfibd}
;
\draw[1cell=.9]  
(dc) edge node[pos=.6] {\grod} (pc)
(pf) edge node {\U} (pc)
(dc) edge[bend left=20] node[pos=.7] {\grod} node[swap,pos=.7] {\sim} (pf)
;
\end{tikzpicture}\]
in \cref{thm:dcatpfibd} \cref{thm:dcatpfibd-ii} shows that the pseudo symmetric $\Cat$-multifunctor
\[\begin{tikzpicture}[xscale=2.75,yscale=1.2,vcenter]
\draw[0cell=.9]
(0,0) node (dc) {\DCat}
(dc)++(1,0) node (pc) {\permcatsg};
\draw[1cell=.9]  
(dc) edge node[pos=.5] {\grod} (pc);
\end{tikzpicture}\]
in \cref{grodpscatmultifunctor}, regarded as a non-symmetric $\Cat$-multifunctor, factors as
\begin{itemize}
\item the non-symmetric $\Cat$-multiequivalence $\grod$ in \cref{thm:dcatpfibdeq}
\item followed by the non-symmetric $\Cat$-multifunctor $\U$ in \cref{pfibdtopermcatsg}.
\end{itemize} 
The case $\cD = \cA$ is used in \cref{thm:invKpseudosym} \cref{thm:invKpseudosym-iii} in the context of inverse $K$-theory.
\end{example}

\part{Pseudo Symmetric Enriched Multifunctorial Inverse \texorpdfstring{$K$}{K}-Theory}
\label{part:invK}

\chapter{The \texorpdfstring{$\Cat$}{Cat}-Multifunctor $A$}
\label{ch:multifunctorA}
The main topic of \cref{part:invK} is the pseudo symmetric $\Cat$-multifunctorial inverse $K$-theory 
\[\begin{tikzpicture}[xscale=1,yscale=1,vcenter]
\def\v{.6}
\draw[0cell=1]
(0,0) node (x1) {\Gacat}
(x1)++(2.8,0) node (x2) {\ACat}
(x2)++(2.8,0) node (x3) {\permcatsg.}
(x1)++(0,\v) node[inner sep=0pt] (s) {}
(x3)++(0,\v) node[inner sep=0pt] (t) {}
;
\draw[1cell=.85] 
(x1) edge node {A} (x2)
(x2) edge node {\groa} (x3)
(x1) edge[-,shorten >=-1pt] (s)
(s) edge[-,shorten <=-1pt,shorten >=-1pt] node {\cP} (t)
(t) edge[shorten <=-1pt] (x3)
;
\end{tikzpicture}\]
The Grothendieck construction 
\[\begin{tikzpicture}[xscale=2.8,yscale=1,vcenter]
\draw[0cell=1]
(0,0) node (x1) {\ACat}
(x1)++(1,0) node (x2) {\permcatsg}
;
\draw[1cell=.85] 
(x1) edge[transform canvas={yshift=-.4ex}] node {\groa} (x2);
\end{tikzpicture}\]
is the pseudo symmetric $\Cat$-multifunctor in \cref{grodpscatmultifunctor} with the following data:
\begin{itemize}
\item $\cA$ is the small tight bipermutative category in \cref{ex:mandellcategory}.
\item $\ACat$ is the $\Cat$-multicategory in \cref{thm:dcatcatmulticat}.
\item $\permcatsg$ is the $\Cat$-multicategory in \cref{thm:permcatenrmulticat}.
\end{itemize} 
In this chapter we construct the $\Cat$-multifunctor (\cref{thm:Acatmultifunctor})
\[\begin{tikzcd}[column sep=large]
\Gacat \ar{r}{A} & \ACat.
\end{tikzcd}\]
We emphasize that---unlike the pseudo symmetric $\Cat$-multifunctor $\groa$ that is \emph{not} a $\Cat$-multifunctor---$A$ is a $\Cat$-multifunctor in the symmetric sense (\cref{def:enr-multicategory-functor}).  In other words, $A$ strictly preserves the symmetric group action, colored units, and composition in $\Gacat$ and $\ACat$.

\subsection*{Summary}

Reusing part of the table in the introduction of \cref{ch:permfib} for $\cD = \cA$ in \cref{ex:mandellcategory}, the following table summaries $\Gacat$, $\ACat$, and $A$.  We abbreviate \emph{symmetric monoidal}, \emph{multimorphism category}, and \emph{natural transformations} to, respectively, \emph{sm}, \emph{mc}, and \emph{nt}.
\begin{center}
\resizebox{\columnwidth}{!}{%
{\renewcommand{\arraystretch}{1.4}%
{\setlength{\tabcolsep}{1ex}
\begin{tabular}{|c|cc|cc|c|}\hline
& $\Gacat$ & (\ref{thm:Gacatsmc}) & $\ACat$ & (\ref{thm:dcatcatmulticat}) & $A \,\,$ (\ref{thm:Acatmultifunctor}) \\ \hline
objects 
& $\Ga$-categories & (\ref{def:gammacategory}) & additive sm functors & (\ref{def:additivesmf}) & \ref{def:AX}, \ref{AXsmf} \\ \hline
mc objects 
& pointed nt & (\ref{expl:Fpointedmap}, \ref{Gacatemptyz}, \ref{Fxonexntoz}) & additive nt & (\ref{def:additivenattr}) & \ref{Aarityzeroobject}, \ref{angAXAFAZ} \\ \hline
mc morphisms 
& pointed modifications & (\ref{expl:Fpointedmap}, \ref{Gacatemptyz}, \ref{thetaordmcomponent}) & additive modifications & (\ref{def:additivemodification}) & \ref{Aarityzeromorphism}, \ref{AFAangXAZ} \\ \hline
colored units & identity nt & (\ref{Gacatcoloredunit}) & identity nt & (\ref{zcoloredunit}) & \ref{thm:Acatmultifunctor} \\ \hline
composition & \ref{Gacatgammacomp} && \ref{dcatgamma}, \ref{expl:dcatgamma} && \ref{Apreservescomp} \\ \hline
symmetry & \ref{rightactionfunctor} && \ref{dcatsigmaaction} && \ref{Apreservessymmetry} \\ \hline
sm & yes & (\ref{thm:Gacatsmc}) & no & (\ref{expl:dcatbipermzero}) & no\\ \hline
\end{tabular}}}}
\end{center}
\smallskip

\subsection*{Organization}

To prepare for the discussion of $A$ and inverse $K$-theory, in \cref{sec:Gammacategories} we review the symmetric monoidal closed structure on the category of $\Ga$-categories.  While $\Ga$-categories are relatively simple to define, the main subtlety of the symmetric monoidal closed structure on $\Gacat$ is that it uses the \emph{pointed} Day convolution \cref{pointedday}, along with the pointed version of the monoidal unit \cref{Gacatunit} and pointed hom diagram \cref{Fpointedhom}.  This symmetric monoidal closed structure induces a $\Cat$-multicategory structure on $\Gacat$, which is explained in detail in \cref{sec:Gacatmulticat}.

To make the discussion of the $\Cat$-multifunctoriality of $A$ easier to understand, in \cref{sec:functorA} we first construct the underlying functor of $A$.  In particular, this section contains the object assignment of $A$, which sends a $\Ga$-category to an object in $\ACat$.  The $\Cat$-multifunctor $A$ is constructed in two stages. 
\begin{itemize}
\item In \cref{sec:multimorphismA} we construct the multimorphism functors of $A$ and observe the following.
\begin{itemize}
\item $A$ is an isomorphism of categories in arity 0 (\cref{Aarityzeroiso}). 
\item In positive arity $A$ is injective on both objects and morphism sets (\cref{Amultimorfunctor}).
\end{itemize} 
\item In \cref{sec:Acatmultifunctor} we prove that the object assignment in \cref{sec:functorA} and the multimorphism functors in \cref{sec:multimorphismA} define a $\Cat$-multifunctor $A$ from $\Gacat$ to $\ACat$.
\end{itemize} 

We remind the reader of our left normalized bracketing \cref{expl:leftbracketing} for iterated monoidal product.

\section{The Symmetric Monoidal Category of \texorpdfstring{$\Ga$}{Gamma}-Categories}
\label{sec:Gammacategories}

The purpose of this section is to discuss the symmetric monoidal closed category of $\Ga$-categories.  In \cref{sec:Gacatmulticat} we discuss its induced $\Cat$-multicategory structure.  Here is an outline of this section.  
\begin{itemize}
\item We begin with a discussion of pointed categories and the symmetric monoidal closed structure on the category $\pCat$ of small pointed categories (\cref{theorem:pC-monoidal}).
\item Then we define $\Ga$-categories (\cref{def:gammacategory}) and discuss its symmetric monoidal closed structure (\cref{thm:Gacatsmc}).
\item We describe explicitly the monoidal unit, the pointed Day convolution, and the pointed mapping object for $\Ga$-categories in, respectively, \cref{expl:Ftu}, \ref{expl:pointedday}, and \ref{expl:Fpointedmap}.
\end{itemize}
We emphasize that the symmetric monoidal structure on $\Gacat$ uses the \emph{pointed} Day convolution in \cref{pointedday}, which is different from the unpointed version in \cref{eq:dayconvolution}.  As a result, the associated monoidal unit \cref{Gacatunit} and internal hom \cref{Fpointedhom} are also different from the unpointed versions in \cref{eq:dayunit,eq:dayhom}.  Much more detailed and general discussion of pointed categories and pointed Day convolution is in \cite[Ch.\! 4]{cerberusIII}.

\subsection*{Pointed Categories}

\begin{definition}[Pointed Diagrams]\label{def:pointedcategory}
A \emph{pointed category}\index{pointed!category}\index{category!pointed} is a pair $(\C,*)$ consisting of a category $\C$ and a distinguished object $* \in \C$, which is called the \index{basepoint}\emph{basepoint}.  The basepoint is also denoted by the functor $\iota^\C \cn \boldone \to \C$ with $\iota^\C(*) = *$.
\begin{itemize}
\item For pointed categories $(\C,*)$ and $(\D,*)$, a \emph{pointed functor}\index{pointed!functor}\index{functor!pointed}
\[\begin{tikzcd}[column sep=large]
(\C,*) \ar{r}{F} & (\D,*)
\end{tikzcd}\]
is a functor $F \cn \C \to \D$ such that $F(*) = *$.
\item For pointed functors $F,G \cn (\C,*) \to (\D,*)$, a \emph{pointed natural transformation}\index{pointed!natural transformation}\index{natural transformation!pointed}
\[\begin{tikzcd}[column sep=large]
F \ar{r}{\theta} & G
\end{tikzcd}\]
is a natural transformation $\theta \cn F \to G$ such that
\[\theta_* = 1_*.\]
\end{itemize}
Denote by $\pCat$\label{not:pCat} the category with
\begin{itemize}
\item small pointed categories as objects and
\item pointed functors as morphisms. 
\end{itemize}
This finishes the definition.
\end{definition}

For an elementary discussion of (co)limits in $\Cat$, the reader is referred to \cite[Section 1.4]{yau-involutive}.  For categories $\C$ and $\D$ with $\C$ small, recall that the \emph{diagram category} $[\C,\D]$ has functors $\C \to \D$ as objects and natural transformations as morphisms (\cref{def:diagramcat}).  

\begin{definition}\label{def:smashproduct}
Suppose $(\C,\iota^\C)$ and $(\D,\iota^\D)$ are small pointed categories.  Define the following small pointed categories.
\begin{itemize}
\item The \index{smash unit}\index{unit!smash}\index{pointed!smash unit}\emph{smash unit}\index{smash unit} is defined as the coproduct in $\Cat$\label{not:smashunit}
\begin{equation}\label{smashunit}
E = \boldone \bincoprod \boldone
\end{equation}
with the basepoint given by the first copy of $\boldone$.
\item The \index{wedge product}\index{product!wedge}\index{pointed!wedge}\emph{wedge product} $\C \wed \D$ is defined as the pushout in $\Cat$ of the span\label{not:wedgeproduct}
\begin{equation}\label{wedgeproduct}
\C \xleftarrow{\iota^\C} \boldone \fto{\iota^\D} \D
\end{equation}
with the basepoint given by the composite 
\[\boldone \to \C \to \C \wed \D.\]
\item The \index{smash product}\index{product!smash}\index{pointed!smash product}\emph{smash product} $\C \sma \D$ is defined as the following pushout in $\Cat$.
\begin{equation}\label{eq:smash-pushout}
\begin{tikzpicture}[x=50mm,y=13mm,vcenter]
\draw[0cell=.85]
(0,0) node (a) {(\C \times \boldone) \bincoprod (\boldone \times \D)}
(1,0) node (b) {\C \times \D}
(0,-1) node (c) {\boldone}
(1,-1) node (d) {\C \sma \D}
;
\draw[1cell=.8]
(a) edge node {(1_\C \times \iota^\D) \bincoprod (\iota^\C \times 1_\D)} (b)
(c) edge node {\iota^{\C \sma \D}} (d)
(a) edge node {} (c)
(b) edge node {\pr} (d)
;
\end{tikzpicture}
\end{equation}
The functor $\iota^{\C \sma \D}$ makes $\C \sma \D$ into a pointed category.
\item The \index{diagram category!pointed}\index{pointed!diagram category}\emph{pointed diagram category} $[\C,\D]_*$ is defined as the following pullback in $\Cat$.
\begin{equation}\label{eq:pHom}
\begin{tikzpicture}[x=30mm,y=13mm,vcenter]
\draw[0cell=.85]
(0,0) node (a) {[\C,\D]_*}
(1,0) node (b) {\boldone}
(0,-1) node (c) {[\C,\D]}
(1,-1) node (d) {[\boldone,\D]}
;
\draw[1cell=.9]
(a) edge node {} (b)
(c) edge node {[\iota^\C,\D]} (d)
(a) edge node {} (c)
(b) edge node {\iota^\D} (d)
;
\end{tikzpicture}
\end{equation}
The composite functor
\[\boldone \fto{\iso} [\C,\boldone] \fto{[\C,\iota^\D]} [\C,\D]
\]
induces a canonical functor 
\[\boldone \to [\C,\D]_*\]
making $[\C,\D]_*$ into a pointed category.
\end{itemize}
This finishes the definition.
\end{definition}

The following result is a special case of \cite[4.2.3]{cerberusIII}.

\begin{theorem}\label{theorem:pC-monoidal}\index{smash product!symmetric monoidal}\index{symmetric monoidal!smash product}
In the context of \cref{def:smashproduct}, 
\[\big(\pCat,\sma,E,[,]_*\big)\]
is a complete and cocomplete symmetric monoidal closed category.
\end{theorem}

\begin{explanation}[Pointed Diagram Category]\label{ex:pointeddiagramcat}
Unraveling the pullback \cref{eq:pHom}, the pointed diagram category $[\C,\D]_*$ is the category with
\begin{itemize}
\item pointed functors $(\C,*) \to (\D,*)$ as objects,
\item pointed natural transformations as morphisms,
\item identity natural transformations as identity morphisms, and
\item vertical composition of natural transformations as composition.
\end{itemize}
The basepoint in $[\C,\D]_*$ is given by the constant functor $\C \to \D$ at the basepoint in $\D$.  Using this explicit description, the pointed diagram category $[\C,\D]_*$ is well defined as long as $\C$ is small.
\end{explanation}

\subsection*{$\Ga$-Categories}

To define $\Ga$-categories, first recall from \cref{ex:Fskel} the small tight bipermutative category 
\[\big(\Fskel, (\vee,\ord{0},\xiwedge), (\sma,\ord{1},\xisma), \fal = (\deltar)^\inv, \far = 1\big).\]
It has
\begin{itemize}
\item pointed finite sets $\ord{n} = \{0,\ldots,n\}$ for $n \geq 0$ as objects and
\item pointed functions as morphisms.
\end{itemize} 
The basepoint of $\ord{n}$ is 0.  
\begin{itemize}
\item The additive structure $\vee$ is given by the wedge of pointed finite sets.
\item The multiplicative structure $\sma$ is given by the smash product of pointed finite sets.
\end{itemize}
The following two statements hold by \cref{thm:smcbipermuative}.
\begin{itemize}
\item The right factorization $\far = (\deltal)^\inv$ is the identity, since $\deltal = 1$ by \cref{Fskeldeltal}.
\item The left factorization $\fal = (\deltar)^\inv$, with $\deltar$ in \cref{Fskeldeltar}, is the pointed bijection given by the inverse of the tetris permutation \cref{fset-deltarformula} for nonzero elements.
\end{itemize}

\begin{definition}\label{def:Fpunc}
For objects $\ord{m}, \ord{n} \in \Fskel$, the \emph{zero morphism}\index{zero morphism}\index{morphism!zero} $\ord{m} \to \ord{n}$ is the composite
\[\ord{m} \to \ord{0} \to \ord{n}.\]
A \emph{nonzero morphism}\index{nonzero morphism}\index{morphism!nonzero} in $\Fskel$ is a morphism that does not factor through $\ord{0}$.  The subset of nonzero morphisms $\ord{m} \to \ord{n}$ in $\Fskel$ is denoted
\begin{equation}\label{nonzeromorphism}
\Fpunc\scmap{\ord{m};\ord{n}}.
\end{equation}
Note that $\Fpunc\smscmap{\ord{m};\ord{n}}$ is empty if either $\ord{m}$ or $\ord{n}$ is $\ord{0}$.
\end{definition}

\begin{definition}[$\Ga$-Categories]\label{def:gammacategory}
For the pointed categories $(\Cat,\boldone)$ and $(\Fskel,\ord{0})$, define the pointed diagram category (\cref{ex:pointeddiagramcat}) 
\[\Gacat = [\Fskel,\Cat]_*.\]
\begin{itemize}
\item An object in $\Gacat$ is a pointed functor 
\[(\Fskel,\ord{0}) \to (\Cat,\boldone)\]
and is called a \index{category!Ga-@$\Ga$-}\index{Ga-category@$\Ga$-category}\emph{$\Ga$-category}.
\item A morphism in $\Gacat$ is a pointed natural transformation and is called a \index{morphism!Ga-category@$\Ga$-category}\index{Ga-category@$\Ga$-category!morphism}\emph{$\Ga$-category morphism}.\defmark
\end{itemize}
\end{definition}

\begin{explanation}[Canonical Basepoints]\label{expl:gammacategory}
Suppose $X \cn (\Fskel,\ord{0}) \to (\Cat,\boldone)$ is a $\Ga$-category. 
\begin{itemize}
\item For each $n \geq 0$, the unique pointed function $\iota_n \cn \ord{0} \to \ord{n}$\label{not:iotan} yields a functor
\[\begin{tikzcd}[column sep=large]
\boldone = X\ord{0} \ar{r}{X\iota_n} & X\ord{n}.
\end{tikzcd}\]
The image of the object $* \in \boldone$ in $X\ord{n}$ under $X\iota_n$ is also denoted by $*$.  Each $X\ord{n} \in \Cat$ is regarded as a pointed category with basepoint $*$.
\item For each morphism $f \cn \ord{m} \to \ord{n}$ in $\Fskel$, the functor
\[\begin{tikzcd}[column sep=large]
X\ord{m} \ar{r}{Xf} & X\ord{n}
\end{tikzcd}\]
is a pointed functor, since $f\iota_m = \iota_n$ in $\Fskel$.
\end{itemize} 
Therefore, $X$ canonically yields a pointed functor\index{Ga-category@$\Ga$-category!canonical basepoint}
\begin{equation}\label{pointedgammacat}
\begin{tikzcd}[column sep=large]
(\Fskel,\ord{0}) \ar{r}{X} & (\pCat,\boldone)
\end{tikzcd}
\end{equation}
with $\pCat$ the category of small pointed categories and pointed functors (\cref{def:pointedcategory}).

Furthermore, suppose $F \cn X \to Y$ is a morphism of $\Ga$-categories, that is, a pointed natural transformation.  It consists of a component functor
\[\begin{tikzcd}[column sep=large]
X\ord{n} \ar{r}{F_{\ord{n}}} & Y\ord{n}
\end{tikzcd}
\foreachspace \ord{n} \in \Fskel\]
that is natural in $\ord{n}$.  By naturality each $F_{\ord{n}}$ is a pointed functor.  Therefore, $F \cn X \to Y$ is a pointed natural transformation when $X$ and $Y$ are regarded as pointed functors with codomain $(\pCat,\boldone)$ as in \cref{pointedgammacat}.  This discussion implies that there is a canonical isomorphism of categories
\[\Gacat \iso [\Fskel,\pCat]_*.\]
In other words, in \cref{def:gammacategory} we can equivalently use $(\pCat,\boldone)$ instead of $(\Cat,\boldone)$ to define the category of $\Ga$-categories.
\end{explanation}

\subsection*{Symmetric Monoidal Closed Structure on $\Gacat$}

We use the notion of a $\V$-category (\cref{def:enriched-category}) for $(\pCat, \sma, E, [,]_*)$ in \cref{theorem:pC-monoidal}.  In \cref{Fhathom,pointedday,Gacatunit} below, the empty wedge is defined as the chosen terminal category $\boldone$.

\begin{definition}\label{def:Fhat}
Define $\Fhat$ as the $\pCat$-category with the same objects as $\Fskel$ (\cref{ex:Fskel}) and hom objects
\begin{equation}\label{Fhathom}
\Fhat \smscmap{\ord{m};\ord{n}} =
\bigvee_{\Fpunc\smscmap{\ord{m};\ord{n}}} E
\end{equation}
in $\pCat$ for $\ord{m}, \ord{n} \in \Fskel$.
\end{definition}

\begin{explanation}\label{expl:Fhat}
The wedge in \cref{Fhathom} is defined in \cref{wedgeproduct}, with
\begin{itemize}
\item $E = \boldone \bincoprod \boldone$ the smash unit \cref{smashunit} and
\item $\Fpunc\smscmap{\ord{m};\ord{n}}$ the set \cref{nonzeromorphism} of nonzero morphisms $\ord{m} \to \ord{n}$.
\end{itemize} 
Since $\Fpunc\smscmap{\ord{0};\ord{n}}$ and $\Fpunc\smscmap{\ord{m};\ord{0}}$ are empty sets, there are pointed categories
\[\Fhat\smscmap{\ord{0};\ord{n}} = \boldone = \Fhat\smscmap{\ord{m};\ord{0}}\]
for $\ord{m}, \ord{n} \in \Fskel$.
\end{explanation}

\begin{definition}\label{def:Gacatsymmetricmonoidal}
Suppose $X$ and $Y$ are $\Ga$-categories.  Define the following $\Ga$-categories.
\begin{itemize}
\item The \emph{pointed Day convolution}\index{Ga-category@$\Ga$-category!pointed Day convolution}\index{pointed!Day convolution}\index{Day convolution!pointed} of $X$ and $Y$ is the $\Ga$-category given by the $\Cat$-coend\index{coend}
\begin{equation}\label{pointedday}
X \Gotimes Y =
\ecint^{(\ord{m}, \ord{n}) \,\in\, \Fhat\sma\Fhat}
\bigvee_{\Fpunc\smscmap{\ord{m} \sma \ord{n}; -}}
\big(X\ord{m} \sma Y\ord{n}\big).
\end{equation}
In \cref{pointedday}:
\begin{itemize}
\item $\Fhat \sma \Fhat$ is the tensor product (\cref{definition:vtensor-0}) of two copies of the $\pCat$-category $\Fhat$ in \cref{def:Fhat}.
\item The smash product $X\ord{m} \sma Y\ord{n}$ is defined in \cref{eq:smash-pushout}, with each of $X$ and $Y$ regarded as a pointed functor $(\Fskel,\ord{0}) \to (\pCat,\boldone)$ as in \cref{pointedgammacat}.
\item $\Fpunc\smscmap{-;-}$ is the set of nonzero morphisms \cref{nonzeromorphism}.
\end{itemize} 
\item The \emph{monoidal unit diagram} is the $\Ga$-category\index{monoidal unit diagram!Ga-category@$\Ga$-category}\index{Ga-category@$\Ga$-category!monoidal unit diagram} given by the wedge
\begin{equation}\label{Gacatunit}
\ftu = \Fhat \smscmap{\ord{1};-} = \bigvee_{\Fpunc\smscmap{\ord{1};-}} E.
\end{equation}
\item The \emph{pointed hom diagram} for $X$ and $Y$ is the $\Ga$-category\index{hom diagram!pointed}\index{pointed!hom diagram}\index{Ga-category@$\Ga$-category!pointed hom diagram} given by the $\Cat$-end
\begin{equation}\label{Fpointedhom}
\GHom(X,Y) = \ecint_{\ord{n} \,\in\, \Fhat}\, \big[X\ord{n} \scs Y(- \sma \ord{n}) \big]_*
\end{equation}
with $[,]_*$ the pointed diagram category in \cref{eq:pHom}.
\end{itemize}
Both $\Gotimes$ and $\GHom$ are extended to $\Ga$-category morphisms componentwise.

Moreover, define the following.
\begin{itemize}
\item The \emph{braiding}\index{braiding!pointed Day convolution}\index{pointed!Day convolution!braiding} for the pointed Day convolution
\begin{equation}\label{eq:Fxi}
\begin{tikzcd}[column sep=large]
X \sma Y \ar{r}{\xiday} & Y \sma X
\end{tikzcd}
\end{equation}
is the pointed natural transformation induced by the braiding for the smash product in $\pCat$
\[\begin{tikzcd}[column sep=large]
X\ord{m} \sma Y\ord{n} \ar{r}{\xi} & Y\ord{n} \sma X\ord{m},
\end{tikzcd}\]
the multiplicative braiding $\xisma$ in $\Fskel$ \cref{xisma}, and the universal properties of coends.
\item Evaluation at the object $\ord{1} \in \Fskel$ defines a symmetric monoidal functor \index{evaluation!at $\ord{1}$}$\ev_{\ord{1}}$ that has a strong symmetric monoidal left adjoint \index{evaluation!at $\ord{1}$!left adjoint}$L_{\ord{1}}$:\label{not:Lordone}
\[\begin{tikzpicture}[xscale=2.5,yscale=1,vcenter]
\draw[0cell]
(0,0) node (x) {\pCat}
(1,0) node (y) {\Gacat.}
;
\draw[1cell=.85] 
(x) edge[transform canvas={yshift=.5ex}] node {L_{\ord{1}}} (y)
(y) edge[transform canvas={yshift=-.4ex}] node {\ev_{\ord{1}}} (x)
;
\end{tikzpicture}\]
\item The \emph{pointed mapping object} is the pointed category\index{pointed!mapping object}\index{Ga-category@$\Ga$-category!pointed mapping object}\index{mapping object!pointed}
\begin{equation}\label{Fpointedmap}
\begin{split}
\GMap(X,Y) &= \big(\GHom(X,Y)\big)(\ord{1})\\
&= \ecint_{\ord{n} \,\in\, \Fhat}\, \left[X\ord{n} \scs Y\ord{n} \right]_*
\end{split}
\end{equation}
given by evaluating $\GHom(X,Y)$ in \cref{Fpointedhom} at $\ord{1} \in \Fskel$.
\end{itemize}
This finishes the definition.
\end{definition}

Recall the notions of \emph{(co)tensored} and a \emph{$\V$-multicategory} in, respectively, \cref{def:tensored,def:enr-multicategory}.  The following result is the special case of \cite[4.3.37 and 6.3.6]{cerberusIII} for $\Gacat$.

\begin{theorem}\label{thm:Gacatsmc}
In the context of \cref{def:Gacatsymmetricmonoidal}, the following statements hold.
\begin{enumerate}
\item\label{thm:Gacatsmc-i} $\big( \Gacat, \Gotimes, \ftu, \GHom\! \big)$ is a complete and cocomplete \index{Ga-category@$\Ga$-category!symmetric monoidal closed}symmetric monoidal closed category.  
\item\label{thm:Gacatsmc-ii} The adjunction $(L_{\ord{1}},\ev_{\ord{1}})$ makes $\Gacat$ enriched, \index{tensored}\index{Ga-category@$\Ga$-category!tensored and cotensored}tensored, and \index{cotensored}cotensored over $\pCat$ with hom object $\GMap$.
\item\label{thm:Gacatsmc-iii} The structure in \eqref{thm:Gacatsmc-i} and \eqref{thm:Gacatsmc-ii} induce a $\pCat$-multicategory structure on $\Gacat$.
\end{enumerate}
\end{theorem}

We also regard $\Gacat$ as enriched over $\Cat$ by forgetting the basepoint in each pointed category.  This is valid because the basepoint-forgetting functor
\[\pCat \to \Cat\]
is symmetric monoidal.  \cref{sec:Gacatmulticat} contains a detailed description of the $\Cat$-multicategory structure on $\Gacat$.  We end this section with an explicit description of each of the monoidal unit diagram $\ftu$, the pointed Day convolution $\Gotimes$, and the pointed mapping object $\GMap$.

\begin{explanation}[Monoidal Unit Diagram]\label{expl:Ftu}\index{monoidal unit diagram!Ga-category@$\Ga$-category}\index{Ga-category@$\Ga$-category!monoidal unit diagram}
Nonzero morphisms $\ord{1} \to \ord{n}$ bijectively correspond to elements in the unpointed finite set with $n$ elements
\[\{1,\ldots,n\} = \ufs{n}.\]  
So for the monoidal unit diagram $\ftu$ in \cref{Gacatunit}, there is a discrete pointed category
\begin{equation}\label{Fhatonen}
\ftu\ord{n} = \Fhat\smscmap{\ord{1};\ord{n}} 
= \bigvee_{i \,\in\, \ufs{n}} \, (\boldone \bincoprod \boldone) \iso \ord{n}
\end{equation}
for each $\ord{n} \in \Fskel$, with $\vee$ the wedge product in \cref{wedgeproduct}.
\end{explanation}

\begin{explanation}[Pointed Day Convolution]\label{expl:pointedday}\index{Ga-category@$\Ga$-category!pointed Day convolution}\index{pointed!Day convolution}\index{Day convolution!pointed}
The pointed Day convolution \cref{pointedday} is characterized by the following universal property.  For $\Ga$-categories $X$, $Y$, and $Z$, a $\Ga$-category morphism
\[\begin{tikzcd}[column sep=large]
X \Gotimes Y \ar{r}{\psi} & Z
\end{tikzcd}\]
consists of precisely a collection of pointed functors
\[\begin{tikzcd}[column sep=huge]
X\ord{m} \sma Y\ord{n} \ar{r}{\psi_{\ord{m},\ord{n}}} & Z\ord{mn}
\end{tikzcd} \forspace \ord{m}, \ord{n} \in \Fskel\]
such that, for each pair of nonzero morphisms
\[\begin{tikzcd}
\ord{m} \ar{r}{f} & \ord{p}
\end{tikzcd} \andspace 
\begin{tikzcd}
\ord{n} \ar{r}{g} & \ord{q}
\end{tikzcd}\]
in $\Fskel$, the diagram
\begin{equation}\label{ptdaynaturality}
\begin{tikzpicture}[xscale=3,yscale=1.2,vcenter]
\draw[0cell=.9]
(0,0) node (x11) {X\ord{m} \sma Y\ord{n}}
(x11)++(1,0) node (x12) {Z\ord{mn}}
(x11)++(0,-1) node (x21) {X\ord{p} \sma Y\ord{q}}
(x12)++(0,-1) node (x22) {Z\ord{pq}}
;
\draw[1cell=.9] 
(x11) edge node {\psi_{\ord{m},\ord{n}}} (x12)
(x21) edge node {\psi_{\ord{p},\ord{q}}} (x22)
(x11) edge node[swap] {Xf \sma Yg} (x21)
(x12) edge node {Z(f \sma g)} (x22)
;
\end{tikzpicture}
\end{equation}
commutes.  The smash product $X\ord{m} \sma Y\ord{n}$ is the pushout in \cref{eq:smash-pushout} using the canonical basepoints in \cref{expl:gammacategory}.
\end{explanation}

\begin{explanation}[Pointed Mapping Object]\label{expl:Fpointedmap}\index{pointed!mapping object}\index{Ga-category@$\Ga$-category!pointed mapping object}\index{mapping object!pointed}
Unraveling the $\Cat$-end in \cref{Fpointedmap} and using \cref{ex:pointeddiagramcat} to interpret each pointed diagram category $[X\ord{n} \scs Y\ord{n}]_*$, the pointed category $\GMap(X,Y)$ can be described explicitly as follows.  

An \emph{object} $F$ in $\GMap(X,Y)$ is a $\Ga$-category morphism $F \cn X \to Y$, that is, a pointed natural transformation.  More explicitly, it consists of pointed functors
\[\begin{tikzcd}[column sep=large]
X\ord{n} \ar{r}{F_{\ord{n}}} & Y\ord{n}
\end{tikzcd} \forspace \ord{n} \in \Fskel\]
such that, for each nonzero morphism $f \cn \ord{m} \to \ord{n}$ in $\Fskel$, the following diagram of functors commutes.
\begin{equation}\label{GMapnaturality}
\begin{tikzpicture}[xscale=2.5,yscale=1.2,vcenter]
\draw[0cell=.9]
(0,0) node (x11) {X\ord{m}}
(x11)++(1,0) node (x12) {Y\ord{m}}
(x11)++(0,-1) node (x21) {X\ord{n}}
(x12)++(0,-1) node (x22) {Y\ord{n}}
;
\draw[1cell=.9] 
(x11) edge node {F_{\ord{m}}} (x12)
(x21) edge node {F_{\ord{n}}} (x22)
(x11) edge node[swap] {Xf} (x21)
(x12) edge node {Yf} (x22)
;
\end{tikzpicture}
\end{equation}
The basepoint in $\GMap(X,Y)$ consists of the constant functors $X\ord{n} \to Y\ord{n}$ at the basepoint of $Y\ord{n}$ for $\ord{n} \in \Fskel$.

For objects $F,G \in \GMap(X,Y)$, a \emph{morphism} $\theta \cn F \to G$ is a pointed modification.  More explicitly, it consists of pointed natural transformations
\[\begin{tikzpicture}[xscale=2.4,yscale=1.5,baseline={(x11.base)}]
\def\a{35}
\draw[0cell=.9]
(0,0) node (x11) {X\ord{n}}
(x11)++(1,0) node (x12) {Y\ord{n}}
;
\draw[1cell=.9] 
(x11) edge[bend left=\a] node {F_{\ord{n}}} (x12)
(x11) edge[bend right=\a] node[swap] {G_{\ord{n}}} (x12)
;
\draw[2cell]
node[between=x11 and x12 at .44, rotate=-90, 2label={above,\theta_{\ord{n}}}] {\Rightarrow}
;
\end{tikzpicture} \forspace \ord{n} \in \Fskel\]
such that, for each nonzero morphism $f \cn \ord{m} \to \ord{n}$ in $\Fskel$, the following two whiskered natural transformations are equal.
\begin{equation}\label{GMapmodaxiom}
\begin{tikzpicture}[xscale=2.4,yscale=1.3,vcenter]
\def\a{35}
\draw[0cell=.8]
(0,0) node (x11) {X\ord{m}}
(x11)++(1,0) node (x12) {Y\ord{m}}
(x11)++(0,-1) node (x21) {X\ord{n}}
(x12)++(0,-1) node (x22) {Y\ord{n}}
;
\draw[1cell=.8] 
(x11) edge[bend left=\a] node[pos=.35] {F_{\ord{m}}} (x12)
(x21) edge[bend left=\a] node[pos=.35] {F_{\ord{n}}} (x22)
(x11) edge[bend right=\a] node[swap,pos=.65] {G_{\ord{m}}} (x12)
(x21) edge[bend right=\a] node[swap,pos=.65] {G_{\ord{n}}} (x22)
(x11) edge node[swap] {Xf} (x21)
(x12) edge node {Yf} (x22)
;
\draw[2cell=.9]
node[between=x11 and x12 at .44, rotate=-90, 2label={above,\theta_{\ord{m}}}] {\Rightarrow}
node[between=x21 and x22 at .44, rotate=-90, 2label={above,\theta_{\ord{n}}}] {\Rightarrow}
;
\end{tikzpicture}
\end{equation}
This axiom is called the \emph{modification axiom}.  

The \emph{identity morphism} $1_F$ of an object $F$ consists of identity natural transformations
\[(1_F)_{\ord{n}} = 1_{F_{\ord{n}}}.\]
Composition of morphisms is given componentwise by vertical composition of natural transformations.

There is a similar description for $\GHom(X,Y)$ in \cref{Fpointedhom}.
\end{explanation}

\section{The \texorpdfstring{$\Cat$}{Cat}-Multicategory of \texorpdfstring{$\Ga$}{Gamma}-Categories}
\label{sec:Gacatmulticat}

By \cref{thm:Gacatsmc} \eqref{thm:Gacatsmc-iii} $\Gacat$ has the structure of a $\pCat$-multicategory, hence also a $\Cat$-multicategory by forgetting the basepoint in each pointed category.  In this section we explain in detail the $\Cat$-multicategory structure on $\Gacat$.\index{Ga-category@$\Ga$-category!Cat-multicategory@$\Cat$-multicategory}\index{Cat-multicategory@$\Cat$-multicategory!of Ga-categories@of $\Ga$-categories}  This discussion is needed for the discussion of inverse $K$-theory, whose domain is $\Gacat$, in \cref{ch:invK}.  

The following discussion of $\Gacat$ is analogous to \cref{expl:permindexedcat,expl:dcatgammaclosed} for $\Dcat$.  The main difference is that, on the one hand, $\Dcat$ is entirely in the unpointed setting.  In particular, the Day convolution \cref{eq:dayconvolution}, the hom diagram \cref{eq:dayhom}, and the unit diagram \cref{eq:dayunit} in $\Dcat$ are all unpointed.  On the other hand, in the current setting of $\Gacat$ (\cref{def:gammacategory}), we use the \emph{pointed} Day convolution and associated monoidal unit, pointed hom diagram, and pointed mapping object in \cref{def:Gacatsymmetricmonoidal}.

\subsection*{Objects and Multimorphism Categories}

An object in $\Gacat$ is a $\Ga$-category (\cref{def:gammacategory}), that is, a pointed functor
\begin{equation}\label{GacategoryZ}
\begin{tikzcd}[column sep=large]
(\Fskel, \ord{0}) \ar{r}{Z} & (\Cat,\boldone).
\end{tikzcd}
\end{equation}
For $\Ga$-categories $Z$ and $\angX = \ang{X_j}_{j=1}^n$, the $n$-ary multimorphism category is
\begin{equation}\label{Gacatnmorphism}
\Gacat\scmap{\ang{X};Z} = \GMap \big( X_1 \Gotimes \cdots \Gotimes X_n \scs Z \big)
\end{equation}
with
\begin{itemize}
\item $\Gotimes$ the pointed Day convolution \cref{pointedday} and
\item $\GMap$ the pointed mapping object \cref{Fpointedmap}.
\end{itemize} 
To unravel $\Gacat\scmap{\ang{X};Z}$ further, we consider the following cases.

\subsection*{Multimorphism Categories in Arity 1}

The 1-ary multimorphism category
\[\Gacat\scmap{X_1;Z} = \GMap(X_1,Z)\]
is as described in \cref{expl:Fpointedmap}.

\subsection*{Multimorphism Categories in Arity 0}

The 0-ary multimorphism category is
\begin{equation}\label{Gacatemptyz}
\Gacat\scmap{\ang{};Z} = \GMap(\ftu,Z)
\end{equation}
with $\ftu$ the monoidal unit diagram \cref{Gacatunit}.  By \cref{expl:Ftu,expl:Fpointedmap}, an \emph{object} $F$ in $\GMap(\ftu,Z)$ consists of pointed functors
\[\begin{tikzcd}[column sep=large]
\ftu\ord{n} \iso \ord{n} \ar{r}{F_{\ord{n}}} & Z\ord{n}
\end{tikzcd} \forspace \ord{n} \in \Fskel\]
that are natural for nonzero morphisms in $\Fskel$ as in \cref{GMapnaturality}.  Moreover, since $\ord{n} = \{0,\ldots,n\}$ is a \emph{discrete} pointed category, $F_{\ord{n}}$ is determined by the objects 
\begin{equation}\label{Fnjobject}
F_{\ord{n}}(j) \in Z\ord{n} \qquad \text{for $j \in \ord{n}$ such that} \qquad F_{\ord{n}}(0) = *,
\end{equation}
the basepoint in $Z\ord{n}$.  

In fact, the objects $F_{\ord{n}}(j)$ in \cref{Fnjobject} are determined by the object $F_{\ord{1}}(1) \in Z\ord{1}$.  In more detail, for each $1 \leq j \leq n$, consider the morphism in $\Fskel$
\begin{equation}\label{iotajordoneordn}
\begin{tikzcd}[column sep=large]
\ord{1} \ar{r}{\nu_j} & \ord{n}
\end{tikzcd}
\stspace \nu_j(1) = j.
\end{equation}
The naturality diagram \cref{GMapnaturality} for $\nu_j$ is the following commutative diagram.
\begin{equation}\label{iotajnaturality}
\begin{tikzpicture}[xscale=1,yscale=1.2,vcenter]
\draw[0cell=.9]
(0,0) node (x11) {\ord{1}}
(x11)++(2.2,0) node (x12) {Z\ord{1}}
(x11)++(0,-1) node (x21) {\ord{n}}
(x12)++(0,-1) node (x22) {Z\ord{n}}
(x11)++(-.6,0) node (j1) {\ftu\ord{1}}
(x21)++(-.6,.03) node (jn) {\ftu\ord{n}} 
node[between=x11 and j1 at .45] {\iso}
node[between=x21 and jn at .45] {\iso}
;
\draw[1cell=.9] 
(x11) edge node {F_{\ord{1}}} (x12)
(x21) edge node {F_{\ord{n}}} (x22)
(x11) edge node[swap] {\nu_j} (x21)
(x12) edge node {Z\nu_j} (x22)
;
\end{tikzpicture}
\end{equation}
Evaluating at the object $1 \in \ord{1}$ yields the object equality
\begin{equation}\label{Fordnjziotaj}
F_{\ord{n}}(j) = (Z\nu_j) \big(F_{\ord{1}}(1)\big)
\inspace Z\ord{n}.
\end{equation}  
Thus an object $F$ in $\GMap(\ftu,Z)$ is determined by the object $F_{\ord{1}}(1) \in Z\ord{1}$.  Moreover, the naturality diagrams \cref{GMapnaturality} for nonzero morphisms $f \in \Fskel$ are determined by the naturality diagrams \cref{iotajnaturality} for $\nu_j$ with $1 \leq j \leq n$.  This is true by the computation \cref{Zfthetaordmi} with $\theta$ replaced by $F$.

For objects $F,G \in \GMap(\ftu,Z)$, a \emph{morphism} $\theta \cn F \to G$ is determined by morphisms in $Z\ord{n}$
\begin{equation}\label{thetaordnj}
\begin{tikzcd}[column sep=large]
F_{\ord{n}}(j) \ar{r}{\theta_{\ord{n},j}} & G_{\ord{n}}(j)
\end{tikzcd} \forspace 0 \leq j \leq n
\end{equation}
such that
\begin{itemize}
\item $\theta_{\ord{n},0} = 1_*$ for each $n \geq 0$ and
\item the modification axiom \cref{GMapmodaxiom} is satisfied.
\end{itemize} 
In fact, the morphisms $\theta_{\ord{n},j}$ in \cref{thetaordnj} are determined by the single morphism 
\begin{equation}\label{thetaordoneone}
\begin{tikzcd}[column sep=large]
F_{\ord{1}}(1) \ar{r}{\theta_{\ord{1},1}} & G_{\ord{1}}(1)
\end{tikzcd}
\inspace Z\ord{1}
\end{equation}
subject to the modification axiom \cref{GMapmodaxiom} for $f = \nu_j$ in \cref{iotajordoneordn} for $1 \leq j \leq n$.  In more detail, the modification axiom for $\nu_j$ states that the two whiskered natural transformations below are equal.
\begin{equation}\label{thetamodaxiomiotaj}
\begin{tikzpicture}[xscale=2.4,yscale=1.3,vcenter]
\def\a{35}
\draw[0cell=.8]
(0,0) node (x11) {\ftu\ord{1}}
(x11)++(1,0) node (x12) {Z\ord{1}}
(x11)++(0,-1) node (x21) {\ftu\ord{n}}
(x12)++(0,-1) node (x22) {Z\ord{n}}
;
\draw[1cell=.8] 
(x11) edge[bend left=\a] node[pos=.35] {F_{\ord{1}}} (x12)
(x21) edge[bend left=\a] node[pos=.35] {F_{\ord{n}}} (x22)
(x11) edge[bend right=\a] node[swap,pos=.65] {G_{\ord{1}}} (x12)
(x21) edge[bend right=\a] node[swap,pos=.65] {G_{\ord{n}}} (x22)
(x11) edge node[swap] {\nu_j} (x21)
(x12) edge node {Z\nu_j} (x22)
;
\draw[2cell=.9]
node[between=x11 and x12 at .44, rotate=-90, 2label={above,\theta_{\ord{1}}}] {\Rightarrow}
node[between=x21 and x22 at .44, rotate=-90, 2label={above,\theta_{\ord{n}}}] {\Rightarrow}
;
\end{tikzpicture}
\end{equation}
Evaluating at the object $1 \in \ord{1}$ yields the morphism equality
\begin{equation}\label{thetaordnjziota}
\theta_{\ord{n},j} = (Z\nu_j) \theta_{\ord{1},1} \inspace Z\ord{n}.
\end{equation}
Moreover, for any nonzero morphism $f \cn \ord{m} \to \ord{n}$ in $\Fskel$, the modification axiom \cref{GMapmodaxiom}---which asserts the equality of the two whiskered natural transformations below---is automatically satisfied.
\begin{equation}\label{modaxiomf}
\begin{tikzpicture}[xscale=2.4,yscale=1.3,vcenter]
\def\a{35}
\draw[0cell=.8]
(0,0) node (x11) {\ftu\ord{m}}
(x11)++(1,0) node (x12) {Z\ord{m}}
(x11)++(0,-1) node (x21) {\ftu\ord{n}}
(x12)++(0,-1) node (x22) {Z\ord{n}}
;
\draw[1cell=.8] 
(x11) edge[bend left=\a] node[pos=.35] {F_{\ord{m}}} (x12)
(x21) edge[bend left=\a] node[pos=.35] {F_{\ord{n}}} (x22)
(x11) edge[bend right=\a] node[swap,pos=.65] {G_{\ord{m}}} (x12)
(x21) edge[bend right=\a] node[swap,pos=.65] {G_{\ord{n}}} (x22)
(x11) edge node[swap] {f} (x21)
(x12) edge node {Zf} (x22)
;
\draw[2cell=.9]
node[between=x11 and x12 at .44, rotate=-90, 2label={above,\theta_{\ord{m}}}] {\Rightarrow}
node[between=x21 and x22 at .44, rotate=-90, 2label={above,\theta_{\ord{n}}}] {\Rightarrow}
;
\end{tikzpicture}
\end{equation}
Indeed, for each $1 \leq i \leq m$ the following diagram in $\Fskel$ is commutative.
\begin{equation}\label{fiotai}
\begin{tikzpicture}[xscale=1,yscale=1,vcenter]
\draw[0cell=.9]
(0,0) node (x11) {\ord{1}}
(x11)++(-1,-1) node (x21) {\ord{m}}
(x11)++(1,-1) node (x22) {\ord{n}}
;
\draw[1cell=.9] 
(x11) edge node[swap] {\nu_i} (x21)
(x21) edge node {f} (x22)
(x11) edge node {\nu_{f(i)}} (x22)
;
\end{tikzpicture}
\end{equation}
Thus there are morphism equalities in $Z\ord{n}$ as follows.
\begin{equation}\label{Zfthetaordmi}
\left\{\begin{aligned}
(Zf) \theta_{\ord{m},i} 
&= (Zf) \big( (Z\nu_i) \theta_{\ord{1},1} \big) & \phantom{=} & \text{by \cref{thetaordnjziota}}\\
&= (Z\nu_{f(i)}) \theta_{\ord{1},1} & \phantom{=} & \text{by \cref{fiotai}}\\
&= \theta_{\ord{n},f(i)} & \phantom{=} & \text{by \cref{thetaordnjziota}}
\end{aligned}\right.
\end{equation}
This proves the modification axiom \cref{modaxiomf} for $f \cn \ord{m} \to \ord{n}$.

\subsection*{Multimorphism Categories in Higher Arity}

For arity $n \geq 2$, we combine \cref{expl:pointedday,expl:Fpointedmap}.  An \emph{object} $F$ in the $n$-ary multimorphism category $\Gacat\scmap{\ang{X};Z}$ in \cref{Gacatnmorphism} is a $\Ga$-category morphism, that is, a pointed natural transformation
\begin{equation}\label{Fxonexntoz}
\begin{tikzcd}[column sep=large]
X_1 \Gotimes \cdots \Gotimes X_n \ar{r}{F} & Z.
\end{tikzcd}
\end{equation}
Such a $\Ga$-category morphism $F$ consists of, for each $n$-tuple of objects $\bang{\ord{m_j}}_{j=1}^n$ in $\Fskel$, a component pointed functor
\begin{equation}\label{Fordmcomponents}
\begin{tikzcd}[column sep=huge]
\txsma_{j=1}^n X_j\ord{m_j} \ar{r}{F_{\ang{\ord{m_j}}}} 
& Z\big(\Gotimes_{j=1}^n \ord{m_j}\big) = Z\big(\ord{m_1 \cdots\, m_n}\big)
\end{tikzcd}
\end{equation}
such that, for each $n$-tuple of nonzero morphisms 
\begin{equation}\label{anghjmorphisms}
\bang{h_j \cn \ord{m_j} \to \ord{p_j}}_{j=1}^n \inspace \Fskel,
\end{equation}
the following \emph{naturality diagram} commutes.
\begin{equation}\label{Fordmnaturality}
\begin{tikzpicture}[xscale=3.5,yscale=1.4,vcenter]
\draw[0cell=.85]
(0,0) node (x11) {\txsma_{j=1}^n X_j\ord{m_j}}
(x11)++(1,0) node (x12) {Z\big(\ord{m_1 \cdots m_n}\big)}
(x11)++(0,-1) node (x21) {\txsma_{j=1}^n X_j\ord{p_j}}
(x12)++(0,-1) node (x22) {Z\big(\ord{p_1 \cdots p_n}\big)}
;
\draw[1cell=.9] 
(x11) edge node {F_{\ang{\ord{m_j}}}} (x12)
(x21) edge node {F_{\ang{\ord{p_j}}}} (x22)
(x11) edge node[swap] {\sma_j\, X_j h_j} (x21)
(x12) edge node {Z(\sma_j\, h_j)} (x22)
;
\end{tikzpicture}
\end{equation}

For objects $F,G \in \Gacat\scmap{\ang{X};Z}$, a \emph{morphism} $\theta \cn F \to G$ is a pointed modification.  It consists of, for each $n$-tuple of objects $\bang{\ord{m_j}}_{j=1}^n$ in $\Fskel$, a component pointed natural transformation
\begin{equation}\label{thetaordmcomponent}
\begin{tikzpicture}[xscale=1,yscale=1,baseline={(x11.base)}]
\def\d{25} 
\draw[0cell=.9]
(0,0) node (x11) {\txsma_{j=1}^n X_j\ord{m_j}}
(x11)++(.6,0) node (a) {\phantom{Z}}
(a)++(2.5,0) node (b) {\phantom{Z}}
(b)++(.7,0) node (x12) {Z\big(\ord{m_1 \cdots\, m_n}\big)}
;
\draw[1cell=.8] 
(a) edge[bend left=\d] node {F_{\ang{\ord{m_j}}}} (b)
(a) edge[bend right=\d] node[swap] {G_{\ang{\ord{m_j}}}} (b)
;
\draw[2cell=.9]
node[between=a and b at .4, rotate=-90, 2label={above,\theta_{\ang{\ord{m_j}}}}] {\Rightarrow}
;
\end{tikzpicture}
\end{equation}
such that, for each $n$-tuple of nonzero morphisms $\ang{h_j}$ as in \cref{anghjmorphisms}, the following two whiskered natural transformations are equal.
\begin{equation}\label{thetaordmodification}
\begin{tikzpicture}[xscale=1,yscale=1,vcenter]
\def\d{25} \def\v{-1.5}
\draw[0cell=.85]
(0,0) node (x11) {\txsma_{j=1}^n X_j\ord{m_j}}
(x11)++(.6,0) node (a) {\phantom{Z}}
(a)++(2.5,0) node (b) {\phantom{Z}}
(b)++(.6,0) node (x12) {Z\big(\ord{m_1 \cdots m_n}\big)}
(x11)++(0,\v) node (x21) {\txsma_{j=1}^n X_j\ord{p_j}}
(a)++(0,\v) node (a2) {\phantom{Z}}
(b)++(0,\v) node (b2) {\phantom{Z}}
(x12)++(0,\v) node (x22) {Z\big(\ord{p_1 \cdots p_n}\big)}
;
\draw[1cell=.8] 
(a) edge[bend left=\d] node[pos=.4] {F_{\ang{\ord{m_j}}}} (b)
(a) edge[bend right=\d] node[swap,pos=.6] {G_{\ang{\ord{m_j}}}} (b)
(a2) edge[bend left=\d] node[pos=.4] {F_{\ang{\ord{p_j}}}} (b2)
(a2) edge[bend right=\d] node[swap,pos=.6] {G_{\ang{\ord{p_j}}}} (b2)
(x11) edge[transform canvas={xshift={1.3em}}] node[swap] {\sma_j\, X_j h_j} (x21)
(x12) edge[transform canvas={xshift={-1.3em}}] node {Z(\sma_j\, h_j)} (x22)
;
\draw[2cell=.9]
node[between=a and b at .4, rotate=-90, 2label={above,\theta_{\ang{\ord{m_j}}}}] {\Rightarrow}
node[between=a2 and b2 at .4, rotate=-90, 2label={above,\theta_{\ang{\ord{p_j}}}}] {\Rightarrow}
;
\end{tikzpicture}
\end{equation}
This is called the \emph{modification axiom}.

\subsection*{Colored Units}

For a $\Ga$-category $Z$ as in \cref{GacategoryZ}, the $Z$-colored unit \cref{ccoloredunit}
\begin{equation}\label{Gacatcoloredunit}
\begin{tikzcd}[column sep=large]
\boldone \ar{r}{1_Z} & \Gacat\smscmap{Z;Z} = \GMap(Z,Z)
\end{tikzcd}
\end{equation}
is the identity natural transformation of the functor $Z \cn \Fskel \to \Cat$.

\subsection*{Symmetric Group Action}

Suppose $\sigma \in \Sigma_n$ is a permutation.  The right action functor \cref{rightsigmaaction}
\begin{equation}\label{rightactionfunctor}
\begin{tikzcd}[column sep=large, row sep=tiny, cells={nodes={scale=.85}}]
\Gacat\scmap{\angX;Z} \ar[equal, shorten <=-.4ex]{d} \ar{r}{\sigma}[swap]{\iso}
& \Gacat\scmap{\angX\sigma;Z} \ar[equal, shorten <=-.4ex]{d}\\
\GMap\brb{\txsma_{j=1}^n X_j, Z} & \GMap\brb{\txsma_{j=1}^n X_{\sigma(j)}, Z} 
\end{tikzcd}
\end{equation}
sends an object $F$ as in \cref{Fxonexntoz} to the following vertical composite pointed natural transformation.
\begin{equation}\label{Gacatsigmaaction}
\begin{tikzpicture}[xscale=1,yscale=.6,vcenter]
\draw[0cell=.9]
(0,0) node (x11) {\txsma_{j=1}^n X_{\sigma(j)}}
(x11)++(2.5,-1) node (x12) {Z}
(x11)++(0,-2) node (x21) {\txsma_{j=1}^n X_j}
;
\draw[1cell=.9]  
(x11) edge node[pos=.3] {F^\sigma} (x12)
(x11) edge node {\iso} node[swap] {\sigma} (x21)
(x21) edge node[swap,pos=.4] {F} (x12)
;
\end{tikzpicture}
\end{equation}
The arrow $\sigma$ permutes the factors $\ang{X_{\sigma(j)}}_{j=1}^n$ using the braiding \cref{eq:Fxi} for the pointed Day convolution.  More explicitly, for each $n$-tuple of objects $\bang{\ord{m_j}}_{j=1}^n$ in $\Fskel$, the component pointed functor of $F^\sigma$ is the following composite.
\begin{equation}\label{Fsigmacomponent}
\begin{tikzpicture}[xscale=1,yscale=1.4,vcenter]
\draw[0cell=.9]
(0,0) node (x11) {\txsma_{j=1}^n X_{\sigma(j)} \ord{m_j}}
(x11)++(4.5,0) node (x12) {Z\big(\ord{m_1 \cdots m_n} \big)}
(x11)++(0,-1) node (x21) {\txsma_{j=1}^n X_j \ord{m_{\sigmainv(j)}}}
(x12)++(0,-1) node (x22) {Z\big(\ord{m_{\sigmainv(1)} \cdots m_{\sigmainv(n)}} \big)}
;
\draw[1cell=.9]  
(x11) edge node {F^\sigma_{\ang{\ord{m_j}}}} (x12)
(x11) edge node {\iso} node[swap] {\sigma} (x21)
(x21) edge node {F_{\ang{\ord{m_{\sigmainv(j)}}}}} (x22)
(x22) edge[shorten <=-.4ex] node {\iso} node [swap] {Z(\sigmainv)} (x12)
;
\end{tikzpicture}
\end{equation}
The detail of \cref{Fsigmacomponent} is as follows.
\begin{itemize}
\item The left vertical arrow $\sigma$ permutes the factors $\bang{X_{\sigma(j)} \ord{m_j}}_{j=1}^n$ according to $\sigma$ using the braiding for the smash product \cref{eq:smash-pushout} in $\pCat$.  The latter is induced by the braiding for the Cartesian product in $\Cat$.
\item The bottom horizontal arrow is the component pointed functor of $F$ as in \cref{Fordmcomponents} for the $n$-tuple of objects $\bang{\ord{m_{\sigmainv(j)}}}_{j=1}^n$.
\item The right vertical arrow is the image under $Z$ of the isomorphism in $\Fskel$
\begin{equation}\label{sigmainvFskel}
\begin{tikzcd}[column sep=large, row sep=small, cells={nodes={scale=.85}}]
\txsma_{j=1}^n \ord{m_{\sigmainv(j)}} \ar[equal]{d} \ar{r}{\sigmainv}[swap]{\iso} 
& \txsma_{j=1}^n \ord{m_j} \ar[equal]{d}\\
\ord{m_{\sigmainv(1)} \cdots\, m_{\sigmainv(n)}} & \ord{m_1 \cdots\, m_n} 
\end{tikzcd}
\end{equation}
that permutes the factors $\bang{\ord{m_{\sigmainv(j)}}}_{j=1}^n$ according to $\sigmainv$ using the multiplicative braiding $\xisma$ \cref{xisma} in $\Fskel$. 
\end{itemize}

Suppose $\theta \cn F \to G$ is a morphism in $\Gacat\scmap{\angX;Z}$ with component pointed natural transformations as in \cref{thetaordmcomponent}.  The right action functor $\sigma$ in \cref{rightactionfunctor} sends $\theta$ to the horizontal composite pointed modification
\[\theta^\sigma = \theta * 1_\sigma \cn F^\sigma \to G^\sigma\]
as in \cref{Gacatsigmaaction} with $F$ replaced by $\theta$.  Its component pointed natural transformations are given by the horizontal composites
\begin{equation}\label{thetasigmacomponent}
\theta^\sigma_{\ang{\ord{m_j}}} = 1_{Z(\sigmainv)} * \theta_{\ang{\ord{m_{\sigmainv(j)}}}} * 1_\sigma
\end{equation}
as in \cref{Fsigmacomponent} with $F$ replaced by $\theta$.

\subsection*{Multicategorical Composition}

Suppose given, for each $j \in \{1,\ldots,n\}$, an $\ell_j$-tuple of $\Ga$-categories $\ang{W_j} = \ang{W_{ji}}_{i=1}^{\ell_j}$ along with 
\begin{equation}\label{ellangW}
\ell = \ell_1 + \cdots + \ell_n \andspace \angW = \ang{\ang{W_j}}_{j=1}^n
\end{equation}
the concatenated $\ell$-tuple.  The multicategorical composition \cref{eq:enr-defn-gamma}
\begin{equation}\label{Gacatgamma}
\begin{tikzpicture}[xscale=5,yscale=1.2,baseline={(x11.base)}]
\draw[0cell=.85]
(0,0) node (x11) {(\Gacat)\scmap{\angX;Z} \times \txprod_{j=1}^n (\Gacat)\scmap{\angWj;X_j}}
(x11)++(1,0) node (x12) {(\Gacat)\scmap{\angW;Z}}
;
\draw[1cell=.9]  
(x11) edge node {\gamma} (x12)
;
\end{tikzpicture}
\end{equation}
is the following composite functor.
\begin{equation}\label{Gacatgammacomp}
\begin{tikzpicture}[xscale=3.5,yscale=1.2,vcenter]
\def\a{15}
\draw[0cell=.85]
(0,0) node (x11) {\textstyle \GMap\brb{\txsma_{j=1}^n X_j,Z} \times \prod_{j=1}^n \GMap\brb{\txsma_{i=1}^{\ell_j} W_{ji},X_j}}
(x11)++(1,-1) node (x12) {\textstyle \GMap\brb{\txsma_{j=1}^n \txsma_{i=1}^{\ell_j} W_{ji},Z}}
(x11)++(0,-2) node (x2) {\textstyle \GMap\brb{\txsma_{j=1}^n X_j,Z} \times \GMap\brb{\txsma_{j=1}^n \txsma_{i=1}^{\ell_j} W_{ji}, \txsma_{j=1}^n X_j}}
;
\draw[1cell=.85]  
(x11) edge[bend left=\a, transform canvas={xshift={4ex}}] node {\gamma} (x12)
(x11) edge[transform canvas={xshift={-1.5em}}] node[swap] {1 \times \sma} (x2)
(x2) edge[bend right=\a, transform canvas={xshift={4ex}}] node[swap,pos=.6] {m} (x12)
;
\end{tikzpicture}
\end{equation}
The detail of \cref{Gacatgammacomp} is as follows.
\begin{itemize}
\item The arrow labeled $m$ is the composition in $\Gacat$ as a $\Cat$-category.  It is given by
\begin{itemize}
\item vertical composition of pointed natural transformations on objects and
\item horizontal composition of pointed modifications on morphisms.
\end{itemize}
\item In the left vertical arrow, $\sma$ is the pointed Day convolution \cref{pointedday}, which extends to pointed natural transformations and pointed modifications componentwise.
\end{itemize}
The explicit formulas of $\gamma$ in terms of component functors and component natural transformations are as in  \cref{dcatganattr,dcatgamod} with the smash product instead of the Cartesian product.

If some $\ell_j = 0$, then \cref{Gacatgammacomp} needs to be slightly adjusted as in the discussion following \cref{dcatgammacomp} using \cref{Gacatemptyz} and the monoidal unit diagram $\ftu$ \cref{Gacatunit}.

\section{The Functor $A$}
\label{sec:functorA}

The purpose of this section is to define the functor
\[\begin{tikzcd}[column sep=large]
\Gacat \ar{r}{A} & \ACat
\end{tikzcd}\]
that is one of the two constituent functors of the inverse $K$-theory functor $\cP$ (\cref{sec:invKfunctor}).
\begin{itemize}
\item The object assignment of $A$ is in \cref{def:AX}.
\item The morphism assignment of $A$ is in \cref{def:AF}.
\item The functoriality of $A$ is proved in \cref{Afunctor}.
\end{itemize}

\subsection*{Object Assignment of $A$}

Recall the small tight bipermutative category
\[\big(\cA, (\oplus,\ang{},\betaplus), (\otimes, (1), \betate), \fal = (\deltar)^\inv, \far = (\deltal)^\inv = 1\big)\]
in \cref{ex:mandellcategory} with
\begin{itemize}
\item additive structure $\Aplus = (\cA,\oplus,\ang{},\betaplus)$ and
\item multiplicative structure $\Ate = (\cA,\otimes, (1), \betate)$. 
\end{itemize}
The inverse $K$-theory functor (\cref{def:mandellinvK}) only uses the additive structure $\Aplus$.  The multiplicative structure $\Ate$ is used when we consider the $\Cat$-multifunctoriality of $A$; see \cref{monen}.

For each integer $n \geq 0$, recall that
\begin{itemize}
\item $\ufs{n} = \{1,\ldots,n\}$ denotes the unpointed finite set with $\ufs{0} = \emptyset$, and
\item $\ord{n} = \{0,1,\ldots,n\}$ denotes the pointed finite set with basepoint 0.
\end{itemize} 
First we define the object assignment of $A$.

\begin{definition}\label{def:AX}
Suppose $X$ is a $\Gamma$-category (\cref{def:gammacategory}).  Define a functor\index{functor!A@$A$}
\begin{equation}\label{functorAX}
\begin{tikzcd}[column sep=large]
\cA \ar{r}{AX} & \Cat
\end{tikzcd}
\end{equation}
as follows.  
\begin{description}
\item[Objects] For an object $m = \ang{m_i}_{i=1}^p \in \cA$ with each $m_i > 0$, define the category
\begin{equation}\label{AXmdef}
(AX)(m) = \begin{cases}
\txprod_{i=1}^p X\ord{m_i} & \text{if $p>0$ and}\\
X\ord{0} = \boldone & \text{if $p=0$.}
\end{cases}
\end{equation}
\item[Morphisms]
For a morphism 
\[\begin{tikzcd}[column sep=large]
\ang{m_i}_{i=1}^p = m \ar{r}{\phi} & n = \ang{n_j}_{j=1}^q
\end{tikzcd} \inspace \cA\]
as in \cref{Amorphismphi}, the functor
\begin{equation}\label{phistarAphi}
\begin{tikzcd}[column sep=2.3cm]
\txprod_{i=1}^p X\ord{m_i} = (AX)(m) \ar{r}{\phi_* \,=\, (AX)(\phi)} 
& (AX)(n) = \txprod_{j=1}^q X\ord{n_j}
\end{tikzcd}
\end{equation}
is defined as follows, where $j \in \{1,\ldots,q\}$.
\begin{enumerate}[label=(\roman*)]
\item\label{phistarAphi-i} If $q = 0$ then $\phi_*$ is the unique functor to the terminal category $\boldone$. 
\item\label{phistarAphi-ii} If $p=0$ and $q>0$, then the composite of $\phi_*$ with the $j$-th coordinate projection is defined as the composite
\[\begin{tikzcd}[column sep=large]
(AX)\ang{} \ar[equal]{d} \ar{r}{\phi_*} & (AX)(n) = \txprod_{j=1}^q X\ord{n_j} \ar[shorten <=-1ex]{d}{\text{project}}\\
X\ord{0} \ar{r}{X\iota_{n_j}} & X\ord{n_j}
\end{tikzcd}\]
with $\iota_{n_j} \cn \ord{0} \to \ord{n_j}$ the unique morphism of pointed finite sets.
\item\label{phistarAphi-iii} Suppose $p,q > 0$ for this and the next cases.  If $\phi^{\inv}\big(\ufs{n_j}\big) = \emptyset$, then the composite of $\phi_*$ with the $j$-th coordinate projection is defined as the following composite.
\[\begin{tikzcd}[column sep=large]
(AX)(m) = \txprod_{i=1}^p X\ord{m_i} \ar[shorten <=-1ex]{d} \ar{r}{\phi_*} 
& (AX)(n) = \txprod_{j=1}^q X\ord{n_j} \ar[shorten <=-1ex]{d}{\text{project}}\\
\boldone = X\ord{0} \ar{r}{X\iota_{n_j}} & X\ord{n_j}
\end{tikzcd}\]
\item\label{phistarAphi-iv} If $\phi^{\inv}\big(\ufs{n_j}\big) \neq \emptyset$, then there is a unique index $i \in \{1,\ldots,p\}$ such that 
\[\emptyset \neq \phi^{\inv}\big(\ufs{n_j}\big) \subset \ufs{m_i}.\]
Define the morphism of pointed finite sets
\begin{equation}\label{phiijpointed}
\begin{tikzcd}[column sep=large]
\ord{m_i} \ar{r}{\phi_{i,j}} & \ord{n_j}
\end{tikzcd} \inspace \Fskel
\end{equation}
by
\[\phi_{i,j}(x) = \begin{cases}
\phi(x) & \text{if $x \in \phi^{\inv}\big(\ufs{n_j}\big) \subset \ufs{m_i}$ and}\\
0 & \text{if $x \in \ord{m_i} \setminus \phi^{\inv}\big(\ufs{n_j}\big)$}.
\end{cases}\]
The composite of $\phi_*$ with the $j$-th coordinate projection is defined as the following composite.
\begin{equation}\label{xofphiij}
\begin{tikzcd}[column sep=large]
(AX)(m) = \txprod_{i=1}^p X\ord{m_i} \ar[shorten <=-1ex]{d}[swap]{\text{project}} \ar{r}{\phi_*} 
& (AX)(n) = \txprod_{j=1}^q X\ord{n_j} \ar[shorten <=-1ex]{d}{\text{project}}\\
X\ord{m_i} \ar{r}{(\phi_{i,j})_* \,=\, X\phi_{i,j}} & X\ord{n_j}
\end{tikzcd}
\end{equation}
\end{enumerate}
\end{description}
The functoriality of $AX$ follows from the functoriality of $X \cn (\Fskel,\ord{0}) \to (\Cat,\boldone)$.  This finishes the definition of the functor $AX$.
\end{definition}

Recall that the additive structure $\Aplus = (\cA,\oplus,\ang{},\betaplus)$ is given by concatenation $\oplus$ of sequences on objects, with the empty sequence $\ang{}$ as the strict monoidal unit (\cref{ex:mandellcategory}).

\begin{lemma}\label{AXsmf}
In the context of \cref{def:AX}, the triple
\[\begin{tikzcd}[column sep=3cm]
\Aplus \ar{r}{\big(AX, (AX)^2, (AX)^0\big)} & \Cat
\end{tikzcd}\]
is a strictly unital strong symmetric monoidal functor with monoidal constraint $(AX)^2$ given by the associativity and unit isomorphisms for the Cartesian product.
\end{lemma}

\begin{proof}
The unit constraint 
\begin{equation}\label{AXzero}
\begin{tikzcd}[column sep=huge]
\boldone \ar{r}{(AX)^0} & (AX)\ang{} = X\ord{0} = \boldone
\end{tikzcd}
\end{equation}
is the identity functor by definition.  For objects $m = \ang{m_i}_{i=1}^p$ and $n = \ang{n_j}_{j=1}^p$ in $\cA$, the monoidal constraint
\begin{equation}\label{AXtwo}
\begin{tikzcd}[column sep=small, row sep=small, cells={nodes={scale=.85}}]
(AX)(m) \times (AX)(n) \ar[equal]{d} \ar{r}{(AX)^2}[swap]{\iso} & (AX)(m \oplus n) \ar[equal]{d}\\
\big(\txprod_{i=1}^p X\ord{m_i}\big) \times \big(\txprod_{j=1}^q X\ord{n_j}\big) & X\ord{m_1} \times \cdots \times X\ord{m_p} \times X\ord{n_1} \times \cdots \times X\ord{n_q}
\end{tikzcd}
\end{equation}
is given by
\begin{itemize}
\item the associativity isomorphism for the Cartesian product if $p,q > 0$,
\item the left unit isomorphism for the Cartesian product if $p=0$, and
\item the right unit isomorphism for the Cartesian product if $q=0$.
\end{itemize}
The naturality of $(AX)^2$ and the axioms of a symmetric monoidal functor (\cref{def:monoidalfunctor}) for $\big(AX,(AX)^2, (AX)^0\big)$ follow from the corresponding properties for the Cartesian product.
\end{proof}

For a $\Ga$-category $X$, we always regard $AX$ as a strictly unital strong symmetric monoidal functor as in \cref{AXsmf}.

\subsection*{Morphism Assignment of $A$}

Next we define the assignment of $A$ on morphisms.

\begin{definition}\label{def:AF}
For a $\Ga$-category morphism $F \cn X \to Y$ (\cref{def:gammacategory}), define a natural transformation
\begin{equation}\label{natAF}
\begin{tikzpicture}[xscale=2.2,yscale=1.5,baseline={(x11.base)}]
\def\a{35}
\draw[0cell=.9]
(0,0) node (x11) {\cA}
(x11)++(1,0) node (x12) {\Cat}
;
\draw[1cell=.8] 
(x11) edge[bend left=\a] node {AX} (x12)
(x11) edge[bend right=\a] node[swap] {AY} (x12)
;
\draw[2cell]
node[between=x11 and x12 at .37, rotate=-90, 2label={above,AF}] {\Rightarrow}
;
\end{tikzpicture}
\end{equation}
by, for each object $m = \ang{m_i}_{i=1}^p \in \cA$, the component functor
\begin{equation}\label{AFm}
(AF)_m = 
\begin{cases}
\txprod_{i=1}^p F_{\ord{m_i}} \cn (AX)(m) \to (AY)(m) & \text{if $p>0$ and}\\
1_\boldone \cn (AX)\ang{} \to (AY)\ang{} & \text{if $p=0$.}
\end{cases}
\end{equation}
The naturality of $AF \cn AX \to AY$ follows from the pointed naturality of $F$. 
\end{definition}

\begin{explanation}[Injectivity]\label{expl:AFinjective}
The assignment $F \mapsto AF$ in \cref{def:AF} is injective.  In other words, if $F,G \cn X \to Y$ are two $\Ga$-category morphisms such that $AF = AG$ as natural transformations, then $F = G$.  This is true by the $p=1$ case of \cref{AFm}.
\end{explanation}

Recall from \cref{def:monoidalnattr} that a \emph{monoidal natural transformation} between monoidal functors is a natural transformation that is compatible with the unit constraints and the monoidal constraints in the sense of \cref{monnattr}.

\begin{lemma}\label{AFmonoidal}
For a $\Ga$-category morphism $F \cn X \to Y$, the natural transformation in \cref{def:AF}
\[\begin{tikzcd}[column sep=large]
AX \ar{r}{AF} & AY
\end{tikzcd}\]
is a monoidal natural transformation, which is, furthermore, a natural isomorphism if and only if $F$ is.
\end{lemma}

\begin{proof}
For the first assertion, we consider the two constraint diagrams in \cref{monnattr}.
\begin{itemize}
\item The unit constraint diagram is commutative because each of $(AF)_{\ang{}}$, $(AX)^0$, and $(AY)^0$ is the identity functor $1_{\boldone}$. 
\item The commutativity of the monoidal constraint diagram for $AF$ follows from the naturality of the associativity and unit isomorphisms for the Cartesian product. 
\end{itemize} 
The second assertion about natural isomorphism follows from the definition \cref{AFm} of $(AF)_m$ as the product of the component functors $F_{\ord{m_i}}$.
\end{proof}

\subsection*{Functoriality of $A$}

To express the functoriality of $A$ in \cref{def:AX,def:AF} and to define the inverse $K$-theory functor, we use the following concept.

\begin{definition}\label{def:underlyingcategory}
For a (non-symmetric) $\Cat$-multicategory $\M$ (\cref{def:enr-multicategory}), its \emph{underlying category}\index{underlying category}\index{category!underlying} $\M_1$ is the category defined as follows.
\begin{itemize}
\item Its objects are those in $\M$.
\item Its morphism sets are 
\[\M_1(a,b) = \Ob\big(\M\smscmap{a;b}\big)\]
for objects $a,b \in \M$.
\item Identity morphisms are the colored units in $\M$.
\item Composition of morphisms is the multicategorical composition in $\M$ restricted to the objects in the 1-ary multimorphism categories.
\end{itemize}
\cref{def:twocategory,locallysmalltwocat,ex:unarycategory} imply that $\M_1$ is a well-defined category.  If the context is clear, then we omit the subscript and denote $\M_1$ by $\M$ to simplify the notation.
\end{definition}

\begin{example}[Permutative Opfibrations]\label{ex:pfibdi}\index{permutative opfibration!underlying category}\index{underlying category!permutative opfibration}
Suppose $(\cD,\oplus,\otimes)$ is a small tight bipermutative category (\cref{def:embipermutativecat}).  The underlying category of the non-symmetric $\Cat$-multicategory $\pfibd$ (\cref{thm:pfibdmulticat}) has
\begin{itemize}
\item permutative opfibrations over $\cD$ (\cref{def:permutativefibration}) as objects;
\item opcartesian 1-linear functors (\cref{def:opcartnlinearfunctor}) as morphisms; and
\item identities and composition (\cref{def:pfibdmulticat}) given by those of symmetric monoidal functors.
\end{itemize}
Permutative opfibrations over $\cD$ only use the additive structure $\oplus$.  Thus, by the description of opcartesian 1-linear functors in \cref{ex:idoponelinear}, only the additive structure $\oplus$, but \emph{not} the multiplicative structure $\otimes$, is used in the underlying category of $\pfibd$. 
\end{example}

\begin{example}[$\Ga$-Categories]\label{ex:Gacati}\index{Ga-category@$\Ga$-category!underlying category}\index{underlying category!Ga-category@$\Ga$-category}
The underlying category of the $\Cat$-multicategory $\Gacat$ (\cref{thm:Gacatsmc,sec:Gacatmulticat}) is the category $[\Fskel,\Cat]_*$ of $\Ga$-categories and their morphisms (\cref{def:gammacategory}).  In \cref{Afunctor,def:mandellinvK} we consider the underlying category of $\Gacat$, which is also denoted $\Gacat$.
\end{example}

\begin{example}[Bipermutative-Indexed Categories]\label{ex:DCati}\index{diagram category!underlying category}\index{underlying category!diagram category}
Suppose $(\cD,\oplus,\otimes)$ is a small tight bipermutative category (\cref{def:embipermutativecat}).  The underlying category of the $\Cat$-multicategory $\DCat$ (\cref{thm:dcatcatmulticat}) has
\begin{itemize}
\item symmetric monoidal functors $\Dplus \to \Cat$ as objects;
\item monoidal natural transformations (\cref{expl:additivenattr} \eqref{expl:additivenattr-iv}) as morphisms; 
\item identity natural transformations as identities; and
\item vertical composition of natural transformations \cref{dcatgammacomp} as composition.
\end{itemize}
Note that only the additive structure $\oplus$, but \emph{not} the multiplicative structure $\otimes$, is used in the underlying category of $\DCat$.  In \cref{Afunctor,def:mandellinvK} we consider the underlying category of $\ACat$, which is also denoted $\ACat$, for the small tight bipermutative category $\cA$ in \cref{ex:mandellcategory}.  
\end{example}

\begin{example}[Permutative Categories]\label{ex:permcatsgone}\index{permutative category!underlying category}\index{underlying category!permutative category}
The underlying category of the $\Cat$-multicategory $\permcatsg$ (\cref{thm:permcatenrmulticat}) has
\begin{itemize}
\item small permutative categories as objects;
\item strictly unital strong symmetric monoidal functors (\cref{ex:onelinearfunctor}) as morphisms; and
\item identities and composition given by those of monoidal functors. 
\end{itemize}
In \cref{def:mandellinvK} we consider the underlying category of $\permcatsg$, which is also denoted $\permcatsg$.
\end{example}

The next result uses the underlying categories of $\Gacat$ and $\ACat$ in, respectively, \cref{ex:Gacati,ex:DCati}.

\begin{lemma}\label{Afunctor}
The object and morphism assignments
\[X \mapsto AX \andspace F \mapsto AF\]
in, respectively, \cref{def:AX,def:AF} define a functor
\begin{equation}\label{AGacatACati}
\begin{tikzcd}[column sep=large]
\Gacat \ar{r}{A} & \ACat.
\end{tikzcd}
\end{equation}
\end{lemma}

\begin{proof}
The object and morphism assignments are well defined by, respectively, \cref{AXsmf,AFmonoidal}.  The functoriality of $A$ follows from the definition \cref{AFm} of each component functor $(AF)_m$ and the functoriality of the Cartesian product.
\end{proof}

The functor $A$ in \cref{AGacatACati} is extended to a $\Cat$-multifunctor in \cref{thm:Acatmultifunctor}.

\section{Multimorphism Functors of \texorpdfstring{$A$}{A}}
\label{sec:multimorphismA}

Our next objective is to extend the functor
\[\begin{tikzcd}[column sep=large]
\Gacat \ar{r}{A} & \ACat
\end{tikzcd}\] 
in \cref{AGacatACati} to a $\Cat$-multifunctor; see \cref{thm:Acatmultifunctor}.  In this section we first construct its multimorphism functors.\index{multimorphism functor!A@$A$}  For context recall the following.
\begin{itemize}
\item $\Gacat$ (\cref{def:gammacategory}) has the structure of a $\Cat$-multicategory (\cref{thm:Gacatsmc}) induced by a symmetric monoidal closed structure.  The monoidal product is given by the \emph{pointed} Day convolution \cref{pointedday}.
\item $\DCat$ is a $\Cat$-multicategory for each small tight bipermutative category $\cD$ (\cref{thm:dcatcatmulticat}).  In particular, $\ACat$ is a $\Cat$-multicategory for the small tight bipermutative category $\cA$ in \cref{ex:mandellcategory}.  The $\Cat$-multicategory structure on $\ACat$ is \emph{not} induced by a monoidal structure (\cref{expl:dcatbipermzero}).
\end{itemize} 
Here is an outline of this section.
\begin{itemize}
\item The multimorphism functors of $A$ in arity 0 are in \cref{def:Aarityzero}.  They are verified to be isomorphisms in \cref{Aarityzeroiso}.
\item The multimorphism functors of $A$ in positive arity are constructed in
\begin{itemize}
\item \cref{def:Aaritynfunctorobj,AFadditivenatural} for objects and 
\item \cref{def:Aaritynfunctormor,Athetaadditivemod} for morphisms.
\end{itemize} 
\item \cref{Amultimorfunctor} shows that these assignments define a functor.
\end{itemize}

\subsection*{Multimorphism Functors of $A$ in Arity 0}

The object assignment of $A$, $X \mapsto AX$ for $\Ga$-categories $X$, is in \cref{def:AX,AXsmf}.  To define the multimorphism functors of $A$, first we consider the case where the input has arity 0. 

\begin{definition}\label{def:Aarityzero}
Suppose $Z$ is a $\Ga$-category.  Define a functor
\begin{equation}\label{Aarityzero}
\begin{tikzpicture}[xscale=1,yscale=.9,vcenter]
\draw[0cell=.9]
(0,0) node (x11) {\Gacat\scmap{\ang{};Z}}
(x11)++(3.5,0) node (x12) {\ACat\scmap{\ang{};AZ}}
(x11)++(0,-1) node (x21) {\GMap(\ftu,Z)}
(x12)++(0,-1) node (x22) {(AZ)(1) = Z\ord{1}}
;
\draw[1cell=.9]  
(x11) edge node {A} (x12)
(x11) edge[equal] (x21)
(x12) edge[equal] (x22)
;
\end{tikzpicture}
\end{equation}
as follows.  In \cref{Aarityzero} the equality on the left is from \cref{Gacatemptyz}, and the equalities on the right are from \cref{dcatbipermzeroary,AXmdef}.
\begin{description}
\item[Objects] By \cref{Fordnjziotaj} an object $F$ in $\GMap(\ftu,Z)$ is determined by the object $F_{\ord{1}}(1)$ in $Z\ord{1}$.  We define the object
\begin{equation}\label{Aarityzeroobject}
AF = F_{\ord{1}}(1) \inspace Z\ord{1}.
\end{equation}
\item[Morphisms] By \cref{thetaordnjziota} a morphism $\theta \cn F \to G$ in $\GMap(\ftu,Z)$ is determined by the morphism $\theta_{\ord{1},1}$ in \cref{thetaordoneone}.  We define the morphism
\begin{equation}\label{Aarityzeromorphism}
A\theta = \theta_{\ord{1},1} \cn AF = F_{\ord{1}}(1) \to AG = G_{\ord{1}}(1) \inspace Z\ord{1}.
\end{equation}
\end{description}
This finishes the definition of $A$ in \cref{Aarityzero}.
\end{definition}

\begin{lemma}\label{Aarityzeroiso}
In the context of \cref{def:Aarityzero},  
\[\begin{tikzpicture}[xscale=1,yscale=.9,vcenter]
\draw[0cell=.9]
(0,0) node (x11) {\Gacat\scmap{\ang{};Z}}
(x11)++(3.5,0) node (x12) {\ACat\scmap{\ang{};AZ}};
\draw[1cell=.9]  
(x11) edge node {A} node[swap] {\iso} (x12);
\end{tikzpicture}\]
is an isomorphism of categories.
\end{lemma}

\begin{proof}
The functoriality of $A$ in \cref{Aarityzero} follows from the fact that identity morphisms and composition in $\GMap(\ftu,Z)$ are given componentwise by those of natural transformations (\cref{expl:Fpointedmap}).  The functor $A$ is an isomorphism for the following reasons.
\begin{itemize}
\item By \cref{Fordnjziotaj,thetaordnjziota} $A$ is injective on, respectively, objects and morphisms. 
\item Given an object $z \in Z\ord{1}$, we define an object $F$ in $\GMap(\ftu,Z)$ using \cref{Fnjobject,iotajordoneordn,Fordnjziotaj}:
\begin{equation}\label{Fordnjdefinitions}
F_{\ord{n}}(j) = \begin{cases}
* \in Z\ord{n} & \text{if $j=0$ and}\\
(Z\nu_j)(z) \in Z\ord{n} & \text{if $1 \leq j \leq n$}.
\end{cases}
\end{equation}
The naturality diagram \cref{iotajnaturality} for $\nu_j$ holds because
\begin{equation}\label{Fordoneonez}
\left\{\begin{aligned}
F_{\ord{1}}(1) &= (Z\nu_1)(z) & \phantom{=} & \text{by \cref{Fordnjdefinitions}}\\ 
&= z & \phantom{=} & \text{by functoriality of $Z$.}
\end{aligned}\right.
\end{equation}
By \cref{Aarityzeroobject,Fordoneonez}
\[AF = F_{\ord{1}}(1) = z,\] 
so $A$ is surjective on objects.
\item The surjectivity of $A$ on morphism sets is proved similarly using \cref{thetaordoneone,thetaordnjziota,Aarityzeromorphism}.
\end{itemize}
Therefore, $A$ is an isomorphism of categories.
\end{proof}

\subsection*{Multimorphism Functors of $A$ in Positive Arity: Object Assignment}

Recall the following.
\begin{itemize}
\item With $\cD = \cA$ in \cref{def:additivenattr}, the $n$-ary 1-cells in $\ACat$ are additive natural transformations.  The latter are natural transformations as in \cref{phitensorxz} that satisfy the unity axiom \cref{additivenattrunity} and the additivity axiom \cref{additivenattradditivity}.  In terms of objects and morphisms, the additivity axiom is \cref{additivityobjects}.
\item By \cref{AXsmf} $A$ sends each $\Ga$-category to a strictly unital strong symmetric monoidal functor $\Aplus \to \Cat$.
\end{itemize} 
Next we define the object assignment of the multimorphism functor
\[\begin{tikzpicture}[xscale=1,yscale=.9,vcenter]
\draw[0cell=.9]
(0,0) node (x11) {\Gacat\scmap{\angX;Z}}
(x11)++(4,0) node (x12) {\ACat\scmap{\ang{AX};AZ}}
;
\draw[1cell=.9]  
(x11) edge node {A} (x12);
\end{tikzpicture}\]
where the input $\angX$ has arity $n > 0$.  When the input has positive arity, we need to consider iterated products in $\cA$ using its multiplicative structure $\otimes$.  To simplify the presentation, we define the following.

\begin{notation}\label{not:cAtensorobjects}
For $1 \leq j \leq n$ suppose given an object
\begin{equation}\label{mjangmji}
m^j = \bang{m^j_{i_j}}_{i_j=1}^{r_j} \in \cA
\end{equation}
with length $r_j \geq 0$ and each $m^j_{i_j} > 0$.  For $1 \leq i_j \leq r_j$ for each $1 \leq j \leq n$, define the positive integer
\begin{equation}\label{monenionein}
m^{1,\ldots,n}_{i_1,\ldots,i_n} = \txprod_{j=1}^n m^j_{i_j}
\end{equation}
and the object
\begin{equation}\label{monen}
\begin{aligned}
m^{1,\ldots,n} &= \txotimes_{j=1}^n m^j \in \cA\\
&= \bang{ \cdots\, \bang{m^{1,\ldots,n}_{i_1,\ldots,i_n}}_{i_1=1}^{r_1} \,\cdots\,}_{i_n=1}^{r_n}
\end{aligned}
\end{equation}
with length $r_1 \cdots\, r_n$.  Here $\otimes$ is the multiplicative structure \cref{Amtensorn} in $\cA$.
\end{notation}

\begin{explanation}\label{expl:monen}
As we mentioned just after \cref{Amtensorn}, the monoidal product
\[\ang{m_i}_{i=1}^p \otimes \ang{n_j}_{j=1}^q = \ang{\ang{m_i n_j}_{i=1}^p}_{j=1}^q \inspace \cA\]
is visualized as the $q \times p$ matrix with $(j,i)$-entry given by the positive integer $m_i n_j$, which is also regarded as the unpointed finite set $\ufs{m_i n_j}$.  Similarly, the $n$-fold monoidal product $m^{1,\ldots,n}$ in \cref{monen} may be visualized as an $n$-dimensional matrix.  A typical entry has the form $m^{1,\ldots,n}_{i_1,\ldots,i_n}$ in \cref{monenionein}, which is also regarded as the unpointed finite set
\[\ufs{m^{1,\ldots,n}_{i_1,\ldots,i_n}} = \txprod_{j=1}^n \ufs{m^j_{i_j}}.\]
If any $r_j = 0$, then $m^j$ is the empty sequence $\ang{} \in \cA$, and so is $m^{1,\ldots,n}$.
\end{explanation}

\begin{definition}\label{def:Aaritynfunctorobj}
Suppose given
\begin{itemize}
\item $\Ga$-categories $Z,X_1, \ldots, X_n$ with $n > 0$, $\angX = \ang{X_j}_{j=1}^n$, and $\ang{AX} = \ang{AX_j}_{j=1}^n$, and
\item an object $F$ in the $n$-ary multimorphism category $\Gacat\scmap{\angX;Z}$, that is, a pointed natural transformation \cref{Fxonexntoz}
\[\begin{tikzcd}[column sep=large]
\txsma_{j=1}^n X_j \ar{r}{F} & Z.
\end{tikzcd}\]
\end{itemize}
Define the data of an additive natural transformation\index{additive!natural transformation!A@$A$} (\cref{def:additivenattr})
\begin{equation}\label{angAXAFAZ}
\begin{tikzcd}[column sep=large]
\ang{AX} \ar{r}{AF} & AZ
\end{tikzcd}
\end{equation}
as follows.  For each $n$-tuple of objects 
\begin{equation}\label{angmangmj}
\ang{m} = \ang{m^j}_{j=1}^n \in \cA^n
\end{equation}
with each $m^j$ as in \cref{mjangmji}, the component functor $(AF)_{\angm}$ in the sense of \cref{dcatnaryobjcomponent} is defined by the following commutative diagrams for $1 \leq i_j \leq r_j$ and $1 \leq j \leq n$.
\begin{equation}\label{AFangm}
\begin{tikzpicture}[xscale=6,yscale=1,vcenter]
\def\v{-1} \def\u{-1.3}
\draw[0cell=.8]
(0,0) node (x11) {\txprod_{j=1}^n (AX_j)m^j}
(x11)++(1,0) node (x12) {(AZ) m^{1,\ldots,n}}
(x11)++(0,\v) node (x21) {\txprod_{j=1}^n \txprod_{i_j=1}^{r_j} \, X_j \ord{m^j_{i_j}}}
(x12)++(0,\v) node (x22) {\txprod_{i_n=1}^{r_n} \,\cdots\, \txprod_{i_1=1}^{r_1} \, Z\ord{m^{1,\ldots,n}_{i_1,\ldots,i_n}}}
(x21)++(0,\u) node (x31) {\txprod_{j=1}^n X_j \ord{m^j_{i_j}}}
(x22)++(0,\u) node (x32) {Z\ord{m^{1,\ldots,n}_{i_1,\ldots,i_n}}}
(x31)++(.5,0) node (x4) {\txsma_{j=1}^n X_j \ord{m^j_{i_j}}}
;
\draw[1cell=.8]  
(x11) edge node {(AF)_{\angm}} (x12)
(x11) edge[equal, shorten <=-.5ex] (x21)
(x12) edge[equal] (x22)
(x21) edge[shorten <=-.5ex, shorten >=-.5ex] node[swap] {\smallprod_j \, \pr_{i_j}} (x31)
(x22) edge[shorten <=-.5ex] node {\pr_{i_1,\ldots,i_n}} (x32)
(x31) edge node[pos=.5] {\pr} (x4)
(x4) edge node[pos=.5] {F_{\ang{m^j_{i_j}}_{j=1}^n}} (x32)
;
\end{tikzpicture}
\end{equation}
This finishes the definition of $AF$ in \cref{angAXAFAZ}.
\end{definition}

\begin{explanation}[Components]\label{expl:AFangm}
The detail of the diagram \cref{AFangm} is as follows.
\begin{itemize}
\item Each of $AX_j$ and $AZ$ is a strictly unital strong symmetric monoidal functor $\Aplus \to \Cat$ (\cref{AXsmf}).  The categories $(AX_j)m^j$ and $(AZ)m^{1,\ldots,n}$ are defined in \cref{AXmdef}.
\item If any $r_j=0$, then $m^{1,\ldots,n}$ in \cref{monen} is the empty sequence, and $(AF)_{\angm}$ is the unique functor to the terminal category $(AZ)\ang{} = \boldone$.
\item Each arrow $\pr_?$ is the $?$-coordinate projection.
\item The arrow $\pr$ is the universal arrow from the Cartesian product to the smash product \cref{eq:smash-pushout}.  It is defined as the identity functor if $n=1$.
\item $F_{\ang{m^j_{i_j}}_{j=1}^n}$ is the $\ang{m^j_{i_j}}_{j=1}^n$-component functor, in the sense of \cref{Fordmcomponents}, of the given pointed natural transformation $F \cn \txsma_{j=1}^n X_j \to Z$.
\end{itemize}
If $n=1$ then $(AF)_{\angm}$ agrees with the earlier definition in \cref{AFm}.
\end{explanation}

\begin{explanation}[Injectivity]\label{expl:AFinj}
The assignment $F \mapsto AF$ in \cref{def:Aaritynfunctorobj} is injective.  This is true by
\begin{itemize}
\item the special case of \cref{AFangm} with each $r_j = 1$, which means that each $m^j \in \cA$ has length 1, and
\item the universal property of the smash product \cref{eq:smash-pushout}.
\end{itemize}
The special case with $n=1$ is in \cref{expl:AFinjective}.
\end{explanation}

\begin{explanation}[Formulas]\label{expl:AFangmformulas}
The component functor $(AF)_{\angm}$ in \cref{AFangm} can be described in terms of object and morphism assignments as follows.  Suppose given
\begin{equation}\label{angxjij}
x^j = \bang{x^j_{i_j}}_{i_j=1}^{r_j} \in \txprod_{i_j=1}^{r_j} X_j \ord{m^j_{i_j}} = (AX_j)m^j 
\forspace 1 \leq j \leq n
\end{equation}
that are either all objects or all morphisms.  The image of
\[\ang{x^j}_{j=1}^n \in \txprod_{j=1}^n (AX_j)m^j\]
under the component functor $(AF)_{\angm}$ is
\begin{equation}\label{AFangmangxj}
\begin{split}
(AF)_{\angm} \ang{x^j}_{j=1}^n 
&= \lrang{ \,\cdots\, \Bang{ \scalebox{.9}{$F_{\ang{m^j_{i_j}}_{j=1}^n} \bang{x^j_{i_j}}_{j=1}^n$} }_{i_1=1}^{r_1} \,\cdots\,}_{i_n=1}^{r_n}\\
& \in (AZ)m^{1,\ldots,n} = \txprod_{i_n=1}^{r_n} \,\cdots\, \txprod_{i_1=1}^{r_1} \, \Big(Z\ord{m^{1,\ldots,n}_{i_1,\ldots,i_n}}\Big).
\end{split}
\end{equation}
Here
\[\bang{x^j_{i_j}}_{j=1}^n \in \txprod_{j=1}^n X_j \ord{m^j_{i_j}}\]
also denotes its image in $\txsma_{j=1}^n X_j \ord{m^j_{i_j}}$ under the universal arrow $\pr$.
\end{explanation}

Next we show that $AF$ in \cref{angAXAFAZ} defines an additive natural transformation (\cref{def:additivenattr}).

\begin{lemma}\label{AFadditivenatural}
In the context of \cref{def:Aaritynfunctorobj}, $AF$ is an additive natural transformation.
\end{lemma}

\begin{proof}
For each $n$-tuple of morphisms in $\cA$, the naturality diagram \cref{dcatnaryobjnaturality} for $AF$ follows from
\begin{itemize}
\item the universal property of Cartesian products,
\item the naturality of each $\pr_?$ and $\pr$ in \cref{AFangm}, and
\item the naturality diagram \cref{Fordmnaturality} for the given pointed natural transformation $F$.
\end{itemize}
Thus $AF$ is a natural transformation in the sense of \cref{phitensorxz}.  It remains to check the two axioms in \cref{def:additivenattr}.

\medskip
\emph{Unity} \cref{additivenattrunity}.  The unit object \cref{additivesmfunitobj} of $AZ \cn \Aplus \to \Cat$ belongs to the terminal category $(AZ)\ang{} = \boldone$.  So the unity axiom holds.

\medskip
\emph{Additivity} \cref{additivenattradditivity}.  This axiom holds for the following reasons.
\begin{itemize}
\item Each component functor $(AF)_{\angm}$ in \cref{angAXAFAZ} is defined by its coordinate projections as a component functor of $F$.
\item By definition \cref{AXtwo} the monoidal constraints $(AX_j)^2$ and $(AZ)^2$ are given by the associativity and unit isomorphisms for the Cartesian product.
\end{itemize}
Thus, in the desired diagram \cref{additivenattradditivity} in the current context, each composite followed by a coordinate projection to a factor of $(AZ)(?)$ is induced by the same component functor of $F$.
\end{proof}

\subsection*{Multimorphism Functors of $A$ in Positive Arity: Morphism Assignment}

Next we define the morphism assignment of the multimorphism functor
\[\begin{tikzpicture}[xscale=1,yscale=.9,vcenter]
\draw[0cell=.9]
(0,0) node (x11) {\Gacat\scmap{\angX;Z}}
(x11)++(4,0) node (x12) {\ACat\scmap{\ang{AX};AZ}}
;
\draw[1cell=.9]  
(x11) edge node {A} (x12);
\end{tikzpicture}\]
where the input $\angX$ has arity $n > 0$.  Recall that a morphism in $\Gacat\scmap{\angX;Z}$ is a pointed modification with component pointed natural transformations as in \cref{thetaordmcomponent} that satisfy the modification axiom \cref{thetaordmodification}.

\begin{definition}\label{def:Aaritynfunctormor}
In the context of \cref{def:Aaritynfunctorobj}, suppose given a morphism 
\[\begin{tikzpicture}[xscale=2,yscale=1.7,baseline={(x1.base)}]
\draw[0cell=.9]
(0,0) node (x1) {\txsma_{j=1}^n X_j}
(x1)++(.18,.03) node (x2) {\phantom{Z}}
(x2)++(1,0) node (x3) {Z}
;
\draw[1cell=.8]  
(x2) edge[bend left] node[pos=.47] {F} (x3)
(x2) edge[bend right] node[swap,pos=.47] {G} (x3)
;
\draw[2cell]
node[between=x2 and x3 at .45, rotate=-90, 2label={above,\theta}] {\Rightarrow}
;
\end{tikzpicture}
\inspace \Gacat\scmap{\angX;Z}.\]  
Define the data of an additive modification\index{additive!modification!A@$A$} (\cref{def:additivemodification})
\begin{equation}\label{AFAangXAZ}
\begin{tikzpicture}[xscale=2.2,yscale=1.7,baseline={(x11.base)}]
\def\a{30}
\draw[0cell=.9]
(0,0) node (x11) {\ang{AX}}
(x11)++(.12,.03) node (x2) {\phantom{Z}}
(x2)++(1,0) node (x12) {AZ}
;
\draw[1cell=.8] 
(x2) edge[bend left=\a] node {AF} (x12)
(x2) edge[bend right=\a] node[swap] {AG} (x12)
;
\draw[2cell]
node[between=x2 and x12 at .38, rotate=-90, 2label={above,A\theta}] {\Rightarrow}
;
\end{tikzpicture}
\end{equation}
by requiring that, for each $n$-tuple of objects $\angm \in \cA^n$ as in \cref{angmangmj}, the two whiskered natural transformations in the following diagram be equal.
\begin{equation}\label{Athetaangm}
\begin{tikzpicture}[xscale=6,yscale=1,vcenter]
\def\v{-1} \def\u{-1.3}
\draw[0cell=.8]
(0,0) node (x11) {\txprod_{j=1}^n (AX_j)m^j}
(x11)++(1,0) node (x12) {(AZ) m^{1,\ldots,n}}
(x11)++(0,\v) node (x21) {\txprod_{j=1}^n \txprod_{i_j=1}^{r_j} \, X_j \ord{m^j_{i_j}}}
(x12)++(0,\v) node (x22) {\txprod_{i_n=1}^{r_n} \,\cdots\, \txprod_{i_1=1}^{r_1} \, Z\ord{m^{1,\ldots,n}_{i_1,\ldots,i_n}}}
(x21)++(0,\u) node (x31) {\txprod_{j=1}^n X_j \ord{m^j_{i_j}}}
(x22)++(0,\u) node (x32) {Z\ord{m^{1,\ldots,n}_{i_1,\ldots,i_n}}}
(x31)++(.5,0) node (x4) {\txsma_{j=1}^n X_j \ord{m^j_{i_j}}}
;
\draw[1cell=.8]  
(x11) edge[out=20,in=160] node[pos=.4] {(AF)_{\angm}} (x12)
(x11) edge[out=-20,in=200] node[pos=.6,swap] {(AG)_{\angm}} (x12)
(x11) edge[equal, shorten <=-.5ex] (x21)
(x12) edge[equal] (x22)
(x21) edge[shorten <=-.5ex, shorten >=-.5ex] node[swap] {\smallprod_j \, \pr_{i_j}} (x31)
(x22) edge[shorten <=-.5ex] node {\pr_{i_1,\ldots,i_n}} (x32)
(x31) edge node[pos=.5] {\pr} (x4)
(x4) edge[out=20,in=160] node[pos=.38] {F} (x32)
(x4) edge[out=-20,in=200] node[pos=.55,swap] {G} (x32)
;
\draw[2cell=.9]
node[between=x11 and x12 at .43, rotate=-90, 2label={above,(A\theta)_{\angm}}] {\Rightarrow}
node[between=x4 and x32 at .5, rotate=-90, 2label={above,\theta}] {\Rightarrow}
;
\end{tikzpicture}
\end{equation}
In this diagram, $\theta \cn F \to G$ denotes the $\ang{m^j_{i_j}}_{j=1}^n$-component pointed natural transformation of $\theta$.
\end{definition}

\begin{explanation}[Components and Injectivity]\label{expl:Athetaangm}
Consider \cref{def:Aaritynfunctormor}.
\begin{enumerate}
\item The diagram \cref{Athetaangm} defining the component $(A\theta)_{\angm}$ is obtained from the diagram \cref{AFangm}, which defines the component functor $(AF)_{\angm}$, by replacing $F$ with $\theta$. 
\item The assignment $\theta \mapsto A\theta$ is injective.  This is true by
\begin{itemize}
\item the special case of \cref{Athetaangm} with each $r_j = 1$ and
\item the fact that the universal arrow $\pr$ from the Cartesian product to the smash product \cref{eq:smash-pushout} is surjective on objects.
\end{itemize}
This extends the injectivity of the assignment $F \mapsto AF$ in \cref{expl:AFinj}.\defmark
\end{enumerate}
\end{explanation}

\begin{explanation}[Formulas]\label{expl:Athetaformulas}
Using the notation in \cref{expl:AFangmformulas}, the component natural transformation $(A\theta)_{\angm}$ in \cref{Athetaangm} is given by the morphism
\begin{equation}\label{Athetaangmangxj}
\begin{split}
\left((A\theta)_{\angm}\right)_{\ang{x^j}_{j=1}^n} 
&= \lrang{ \,\cdots\, \Bang{ \theta_{\ang{x^j_{i_j}}_{j=1}^n} }_{i_1=1}^{r_1} \,\cdots\,}_{i_n=1}^{r_n}\\
& \in (AZ)m^{1,\ldots,n} = \txprod_{i_n=1}^{r_n} \,\cdots\, \txprod_{i_1=1}^{r_1} \, \Big(Z\ord{m^{1,\ldots,n}_{i_1,\ldots,i_n}}\Big)
\end{split}
\end{equation}
for each object $\ang{x^j}_{j=1}^n$ in $\txprod_{j=1}^n (AX_j)m^j$.
\end{explanation}

\begin{lemma}\label{Athetaadditivemod}
In the context of \cref{def:Aaritynfunctormor}, $A\theta$ is an additive modification.
\end{lemma}

\begin{proof}
The naturality of each component $(A\theta)_{\angm}$ in \cref{Athetaangm} follows from
\begin{itemize}
\item the universal property of Cartesian products,
\item the naturality of each $\pr_?$ and $\pr$ in \cref{Athetaangm},
\item the coordinate-wise definition \cref{AFangm} of $(AF)_{\angm}$ and $(AG)_{\angm}$, and
\item the naturality of each component natural transformation of the given pointed modification $\theta$.
\end{itemize}
The modification axiom \cref{dcatnarymodaxiom} for $A\theta$ follows for the same reasons from the modification axiom \cref{thetaordmodification} for $\theta$.  Thus $A\theta$ is a modification in the sense of \cref{Phiphivarphi}.  The unity axiom \cref{additivemodunity} and the additivity axiom \cref{additivemodadditivity} for $A\theta$ are proved as in \cref{AFadditivenatural} with $F$ replaced by $\theta$.
\end{proof}

\begin{lemma}\label{Amultimorfunctor}
The object and morphism assignments
\[F \mapsto AF \andspace \theta \mapsto A\theta\]
in, respectively, \cref{def:Aaritynfunctorobj,def:Aaritynfunctormor} define a functor
\begin{equation}\label{Aaritynfunctor}
\begin{tikzpicture}[xscale=1,yscale=.9,baseline={(x11.base)}]
\draw[0cell=.9]
(0,0) node (x11) {\Gacat\scmap{\angX;Z}}
(x11)++(4,0) node (x12) {\ACat\scmap{\ang{AX};AZ}}
;
\draw[1cell=.9]  
(x11) edge node {A} (x12);
\end{tikzpicture}
\end{equation}
that is injective on objects and morphism sets.
\end{lemma}

\begin{proof}
The object and morphism assignments are well defined by, respective, \cref{AFadditivenatural,Athetaadditivemod}.  By \cref{dcatnary,expl:Fpointedmap,Gacatnmorphism}, in both the domain and the codomain of $A$, identity morphisms and composition are given componentwise by those of natural transformations.  Thus the functoriality of $A$ follows from the definition \cref{Athetaangm} of $A\theta$ as given componentwise by the corresponding component of $\theta$.  The injectivity of $A$ on objects and morphism sets is discussed in \cref{expl:AFinj,expl:Athetaangm}.
\end{proof}

\section{\texorpdfstring{$\Cat$}{Cat}-Multifunctoriality of \texorpdfstring{$A$}{A}}
\label{sec:Acatmultifunctor}

The purpose of this section is to show that $A$ with
\begin{itemize}
\item object assignment $X \mapsto AX$ in \cref{AXsmf} and 
\item multimorphism functors in \cref{Amultimorfunctor}
\end{itemize}
is a $\Cat$-multifunctor\index{Cat-multifunctor@$\Cat$-multifunctor!A@$A$} in the sense of \cref{def:enr-multicategory-functor}.  Here is an outline of this section.
\begin{itemize}
\item \cref{Apreservessymmetry} shows that $A$ preserves the symmetric group action of the $\Cat$-multicategories $\Gacat$ and $\ACat$.
\item \cref{Apreservescomp} shows that $A$ preserves the composition in $\Gacat$ and $\ACat$.
\item The main result on the $\Cat$-multifunctoriality of $A$ is \cref{thm:Acatmultifunctor}.
\end{itemize}

\subsection*{Preservation of Symmetric Group Action}

Recall that the symmetric group action is discussed in
\begin{itemize}
\item \cref{Fsigmacomponent,thetasigmacomponent} for the $\Cat$-multicategory $\Gacat$ and
\item \cref{phisigmaanga,Phisigmaanga} for the $\Cat$-multicategory $\ACat$.
\end{itemize} 

\begin{lemma}\label{Apreservessymmetry}
The functors $A$ in \cref{Aaritynfunctor} preserve the symmetric group action.
\end{lemma}

\begin{proof}
We must show that, for each permutation $\sigma \in \Sigma_n$, the following diagram of functors commutes.
\begin{equation}\label{Asigmaaction}
\begin{tikzpicture}[xscale=1,yscale=1.3,vcenter]
\draw[0cell=.9]
(0,0) node (x11) {\Gacat\scmap{\angX;Z}}
(x11)++(4,0) node (x12) {\ACat\scmap{\ang{AX};AZ}}
(x11)++(0,-1) node (x21) {\Gacat\scmap{\angX\sigma;Z}}
(x12)++(0,-1) node (x22) {\ACat\scmap{\ang{AX}\sigma;AZ}}
;
\draw[1cell=.9]  
(x11) edge node {A} (x12)
(x21) edge node {A} (x22)
(x11) edge node {\iso} node[swap] {\sigma} (x21)
(x12) edge node[swap] {\iso} node {\sigma} (x22)
;
\end{tikzpicture}
\end{equation}
Suppose $F$ is an object in $\Gacat\scmap{\angX;Z}$, which means a pointed natural transformation $F \cn \txsma_{j=1}^n X_j \to Z$ as in \cref{Fxonexntoz}.  With $(-)^\sigma$ denoting the right $\sigma$-action, to show that the additive natural transformations
\[\begin{tikzpicture}[xscale=1,yscale=1,vcenter]
\draw[0cell=.9]
(0,0) node (x11) {\ang{AX}\sigma}
(x11)++(3,0) node (x12) {\phantom{AZ}}
(x12)++(0,.03) node (az) {AZ} 
;
\draw[1cell=.9]  
(x11) edge[transform canvas={yshift={.5ex}}] node {A(F^\sigma)} (x12)
(x11) edge[transform canvas={yshift={-.5ex}}] node[swap] {(AF)^\sigma} (x12)
;
\end{tikzpicture}\]
are equal, it suffices to show that their corresponding component functors \cref{dcatnaryobjcomponent} are equal (\cref{expl:additivenattr} \eqref{expl:additivenattr-i}).  Moreover, by the definition \cref{AXmdef} of $AZ$, it suffices to show that the corresponding component functors of $A(F^\sigma)$ and $(AF)^\sigma$ followed by each coordinate projection are equal.  Thus, using \cref{not:cAtensorobjects} and the abbreviations
\begin{equation}\label{mtilde}
\tau = \sigmainv, \qquad \mtil = m^{1,\ldots,n}_{i_1,\ldots,i_n}, \andspace 
\mtil^{\tau} = m^{\sigmainv(1),\ldots,\sigmainv(n)}_{i_{\sigmainv(1)},\ldots,i_{\sigmainv(n)}}
\end{equation}
it suffices to show that the following diagram, which will be explained further below, is commutative.
\begin{equation}\label{Asigmaactiondiag}
\begin{tikzpicture}[xscale=1,yscale=1,vcenter]
\def\v{-1} \def\u{-1.3} \def\h{2.7} \def\g{1.3}
\draw[0cell=.75]
(0,0) node (x1) {\txprod_{j=1}^n \txprod_{i_j=1}^{r_j} X_{\sigma(j)} \ord{m^j_{i_j}} = \txprod_{j=1}^n (AX_{\sigma(j)}) m^j}
(x1)++(-\h,\v) node (x21) {\txprod_{j=1}^n X_{\sigma(j)} \ord{m^j_{i_j}}}
(x1)++(\h,\v) node (x22) {\phantom{\txprod_{j=1}^n \txprod_{i_{\tau(j)}=1}^{r_{\tau(j)}} \, X_j \ord{m^{\tau(j)}_{i_{\tau(j)}}}}}
(x22)++(1,0) node (x23) {\txprod_{j=1}^n \txprod_{i_{\tau(j)}=1}^{r_{\tau(j)}} \, X_j \ord{m^{\tau(j)}_{i_{\tau(j)}}} = \txprod_{j=1}^n (AX_j) m^{\tau(j)}}
(x21)++(\h,\v) node (x3) {\txprod_{j=1}^n \, X_j \ord{m^{\tau(j)}_{i_{\tau(j)}}}}
(x21)++(0,\u) node (x41) {\txsma_{j=1}^n X_{\sigma(j)} \ord{m^j_{i_j}}}
(x22)++(0,\u) node (x42) {\phantom{\txprod_{i_{\tau(n)}=1}^{r_{\tau(n)}} \,\cdots\, \txprod_{i_{\tau(1)}=1}^{r_{\tau(1)}} \, Z \ord{\mtil^\tau}}}
(x42)++(1,0) node (x43) {\txprod_{i_{\tau(n)}=1}^{r_{\tau(n)}} \,\cdots\, \txprod_{i_{\tau(1)}=1}^{r_{\tau(1)}} \, Z \ord{\mtil^\tau} = (AZ)m^{\tau(1),\ldots,\tau(n)}}
(x41)++(0,\u) node (x51) {\txsma_{j=1}^n X_j \ord{m^{\tau(j)}_{i_{\tau(j)}}}} 
(x42)++(0,\u) node (x52) {\phantom{\txprod_{i_n=1}^{r_n} \,\cdots\, \txprod_{i_1=1}^{r_1} \, Z \ord{\mtil}}}
(x52)++(.8,0) node (x53) {\txprod_{i_n=1}^{r_n} \,\cdots\, \txprod_{i_1=1}^{r_1} \, Z \ord{\mtil} = (AZ)m^{1,\ldots,n}}
(x51)++(\g,\v) node (x61) {Z \ord{\mtil^\tau}}
(x52)++(-\g,\v) node (x62) {Z \ord{\mtil}}
;
\draw[0cell]
node[between=x1 and x3 at .45] {\spadesuit}
node[between=x41 and x3 at .5] {\blacktriangle}
node[between=x3 and x61 at .5] {\clubsuit}
node[between=x61 and x52 at .6] {\filledmedlozenge}
;
\draw[1cell=.7]
(x1) edge[out=180,in=90] node[swap,pos=.7] {\smallprod_j \,\pr_{i_j}} (x21)
(x21) edge[shorten <=-1ex] node[swap] {\pr} (x41)
(x41) edge[shorten <=-1ex] node[swap] {\sigma} (x51)
(x51) edge[out=-90,in=180] node[swap,pos=.2] {F} (x61)
(x61) edge node[pos=.6] {Z(\sigmainv)} (x62)
(x1) edge[out=0,in=90] node[pos=.6] {\sigma} (x22)
(x22) edge node {(AF)_{\sigma\angm}} (x42)
(x42) edge node {(AZ)(\sigmainv)} (x52)
(x52) edge[out=-90,in=0] node[pos=.3] {\pr_{i_1,\ldots,i_n}} (x62)
(x21) edge node {\sigma} (x3)
(x22) edge[shorten <=-2ex] node[swap,pos=.3] {\smallprod_j \, \pr_{i_{\tau(j)}}} (x3)
(x3) edge[shorten <=-1ex] node[swap,pos=.6] {\pr} (x51)
(x42) edge node[swap,pos=.2] {\pr_{i_{\tau(1)}, \ldots, i_{\tau(n)}}} (x61)
;
\end{tikzpicture}
\end{equation}
The detail of the diagram \cref{Asigmaactiondiag} is as follows.
\begin{itemize}
\item By \cref{AFangm,Fsigmacomponent} the left-bottom boundary composite is the composite functor
\[\pr_{i_1,\ldots,i_n} \circ A(F^\sigma)_{\angm}.\]
\item By \cref{phisigmaanga} the right boundary composite is the composite functor
\[\pr_{i_1,\ldots,i_n} \circ (AF)^\sigma_{\angm}.\]
\item The region $\spadesuit$ is commutative by the naturality of the braiding for the Cartesian product.
\item The region $\blacktriangle$ is commutative by the universal property of the pushout \cref{eq:smash-pushout} that defines the smash product.
\item The region $\clubsuit$ is commutative by \cref{AFangm} applied to the $n$-tuple of objects 
\[\sigma\angm = \bang{m^{\sigmainv(j)}}_{j=1}^n \in \cA^n.\]
\item The region $\filledmedlozenge$ is commutative by \cref{Abetate,xofphiij}.
\end{itemize}
This proves that the diagram \cref{Asigmaaction} is commutative on objects.

To prove that the diagram \cref{Asigmaaction} is commutative on morphisms, we reuse the proof for objects above by replacing \cref{phisigmaanga,AFangm,Fsigmacomponent} with, respectively, \cref{Phisigmaanga,Athetaangm,thetasigmacomponent}.
\end{proof}

\subsection*{Preservation of Composition}

Next we show that $A$ preserves the composition in the $\Cat$-multicategories $\Gacat$ and $\ACat$, which are defined in, respectively, \cref{Gacatgammacomp} and \cref{dcatgammacomp}.

\begin{lemma}\label{Apreservescomp}
The functors $A$ in \cref{Aaritynfunctor} preserve the multicategorical composition.
\end{lemma}

\begin{proof}
Using the notation in \cref{ellangW}, we must show that the following diagram of functors is commutative.
\begin{equation}\label{Agamma}
\begin{tikzpicture}[xscale=1,yscale=1.3,vcenter]
\draw[0cell=.8]
(0,0) node (x11) {(\Gacat)\scmap{\angX;Z} \times \txprod_{j=1}^n (\Gacat)\scmap{\angWj;X_j}}
(x11)++(5.2,0) node (x12) {(\Gacat)\scmap{\angW;Z}}
(x11)++(0,-1) node (x21) {(\ACat)\scmap{\ang{AX};AZ} \times \txprod_{j=1}^n (\ACat)\scmap{\ang{AW_j};AX_j}}
(x12)++(0,-1) node (x22) {(\ACat)\scmap{\ang{AW};AZ}} 
;
\draw[1cell=.9]  
(x11) edge node {\gamma} (x12)
(x21) edge node {\gamma} (x22)
(x11) edge node[swap] {A \times \smallprod_j A} (x21)
(x12) edge node {A} (x22)
;
\end{tikzpicture}
\end{equation}
Suppose given an object $\scmap{F;\ang{F_j}_{j=1}^n}$ in the upper-left category in \cref{Agamma}, so 
\[\begin{tikzcd}
\txsma_{j=1}^n X_j \ar{r}{F} & Z
\end{tikzcd} \andspace 
\begin{tikzcd}
\txsma_{i=1}^{\ell_j} W_{ji} \ar{r}{F_j} & X_j
\end{tikzcd}\]
are pointed natural transformations as in \cref{Fxonexntoz}.  We must show that the two additive natural transformations below are equal.
\begin{equation}\label{AgagaA}
\begin{tikzpicture}[xscale=1,yscale=1,baseline={(x11.base)}]
\draw[0cell=.9]
(0,0) node (x11) {\ang{AW}}
(x11)++(3.5,0) node (x12) {\phantom{AZ}}
(x12)++(0,.03) node (az) {AZ} 
;
\draw[1cell=.8]  
(x11) edge[transform canvas={yshift={.5ex}}] node {A\ga\scmap{F;\ang{F_j}}} (x12)
(x11) edge[transform canvas={yshift={-.5ex}}] node[swap] {\ga\scmap{AF; \ang{AF_j}}} (x12)
;
\end{tikzpicture}
\end{equation}
As in the proof of \cref{Apreservessymmetry}, it suffices to show that in \cref{AgagaA} the corresponding component functors \cref{dcatnaryobjcomponent} followed by each coordinate projection are equal.  This is proved in the diagram \cref{Agammadiag} using the following notation.

Suppose given objects
\[p^{ji} = \bang{p^{ji}_{k_{ji}}}_{k_{ji}=1}^{r_{ji}} \in \cA 
\qquad \text{for $1 \leq j \leq n$ and $1 \leq i \leq \ell_j$}\]
with each $p^{ji}_{k_{ji}} > 0$ and $r_{ji} \geq 0$, along with the following notation.
\[\left\{\begin{aligned}
p^{j\bdot} &= \ang{p^{ji}}_{i=1}^{\ell_j} \in \cA^{\ell_j} 
& p^j &= \txotimes_{i=1}^{\ell_j} p^{ji} \in \cA \qquad \ang{p^j}_j = \ang{p^j}_{j=1}^n \in \cA^n\\
p^{\bdot\bdot} &= \bang{p^{j\bdot}}_{j=1}^n \in \cA^{\ell}
& p &= \txotimes_{j=1}^n p^j \in \cA\\
(k_j) &= \big(k_{j\ell_j}, \ldots, k_{j1}\big) & (k) &= \big((k_n), \ldots, (k_1)\big)\\
\end{aligned}\right.\]
Moreover, in (smash) products the indices $j$, $i$, and $k_{ji}$ run through, respectively,
\[1 \leq j \leq n, \qquad 1 \leq i \leq \ell_j, \andspace 1 \leq k_{ji} \leq r_{ji}.\]
Here are the key examples that appear in \cref{Agammadiag}.
\[\left\{\begin{split}
\txprod_{j,i} &= \txprod_{j=1}^n \txprod_{i=1}^{\ell_j} \qquad \txprod_j \txsma_i = \txprod_{j=1}^n \txsma_{i=1}^{\ell_j}\\
\txprod_{j,i,k_{ji}} &= \txprod_{j=1}^n \txprod_{i=1}^{\ell_j} \txprod_{k_{ji}=1}^{r_{ji}}\\
\txprod_{j,(k_j)} &= \txprod_{j=1}^n \txprod_{k_{j\ell_j}=1}^{r_{j\ell_j}} \,\cdots\, \txprod_{k_{j1}=1}^{r_{j1}}\\
\txprod_{(k)} &= \txprod_{k_{n\ell_n}=1}^{r_{n\ell_n}} \,\cdots\, \txprod_{k_{n1}=1}^{r_{n1}} \,\cdots\, \txprod_{k_{1\ell_1}=1}^{r_{1\ell_1}} \,\cdots\, \txprod_{k_{11}=1}^{r_{11}}
\end{split}\right.\]
With the notation above, we consider the following diagram of functors.
\begin{equation}\label{Agammadiag}
\begin{tikzpicture}[xscale=1,yscale=1,vcenter]
\def\t{-.7} \def\v{-1} \def\u{-1.3} \def\h{3.7} \def\g{.3}
\draw[0cell=.75]
(0,0) node (x1) {\txprod_{j,i} (AW_{ji})p^{ji} = \txprod_{j,i,k_{ji}} W_{ji} \ord{p^{ji}_{k_{ji}}}}
(x1)++(-\h,\t) node (x2) {\txprod_{j,(k_j)} X_j \Big(\ord{\smallprod_i p^{ji}_{k_{ji}}}\Big)}
(x2)++(0,\v) node (x3) {\txprod_j (AX_j)p^j}
(x3)++(0,\u) node (x4) {(AZ)p}
(x4)++(0,\v) node (x5) {\txprod_{(k)} Z\Big(\ord{\smallprod_{j,i} p^{ji}_{k_{ji}}}\Big)}
(x5)++(\h,-\g) node (z) {Z\Big(\ord{\smallprod_{j,i} p^{ji}_{k_{ji}}}\Big)}
(x1)++(\h,\t) node (y1) {\txprod_{j,i} W_{ji} \ord{p^{ji}_{k_{ji}}}}
(z)++(\h,\g) node (y3) {\txsma_j X_j \Big(\ord{\smallprod_i p^{ji}_{k_{ji}}}\Big)}
node[between=y1 and y3 at .5] (y2) {\txsma_{j,i} W_{ji} \ord{p^{ji}_{k_{ji}}}}
node[between=x2 and y3 at .4] (z2) {\txprod_j X_j \Big(\ord{\smallprod_i p^{ji}_{k_{ji}}}\Big)}
node[between=z2 and y1 at .5] (z1) {\txprod_j \txsma_i W_{ji} \ord{p^{ji}_{k_{ji}}}}
;
\draw[1cell=.7]
(x1) edge[out=180,in=30] node[swap,pos=.85] {\smallprod_j (AF_j)_{p^{j\bdot}}} (x2)
(x2) edge[equal, shorten <=-.5ex] (x3)
(x3) edge node[swap] {(AF)_{\ang{p^j}_j}} (x4)
(x4) edge[equal] (x5)
(x5) edge node {\pr} (z)
(x1) edge[out=0,in=150,shorten >=-1ex] node[pos=.85] {\smallprod_{j,i}\, \pr} (y1)
(y1) edge node {\pr} (y2)
(y2) edge node {\sma_j F_j} (y3)
(y3) edge node[swap] {F} (z)
(y1) edge[shorten <=-.5ex] node[swap,pos=.3] {\smallprod_j \, \pr} (z1)
(z1) edge[shorten >=-.5ex] node (p) {\pr} (y2)
(z1) edge[shorten <=-1ex] node[swap,pos=.3] {\smallprod_j F_j} (z2)
(z2) edge node {\pr} (y3)
(x2) edge node {\smallprod_j \, \pr} (z2)
;
\draw[0cell]
node[between=x1 and z2 at .4] {\spadesuit}
node[between=y1 and p at .6] {\blacktriangle}
node[between=z2 and y2 at .5] {\clubsuit}
node[between=x3 and z at .5] {\filledmedlozenge}
;
\end{tikzpicture}
\end{equation}
The detail of the diagram \cref{Agammadiag} is as follows.
\begin{itemize}
\item By \cref{dcatganattr} the top-left-bottom boundary composite is
\[\pr_{\ang{\ang{k_{ji}}_{i=1}^{\ell_j}}_{j=1}^n} \circ \ga\scmap{AF;\ang{AF_j}}_{p^{\bdot\bdot}}.\]
\item By \cref{Gacatgammacomp,AFangm} the top-right-bottom boundary composite is
\[\pr_{\ang{\ang{k_{ji}}_{i=1}^{\ell_j}}_{j=1}^n} \circ \big(A\ga\scmap{F;\ang{F_j}}\big)_{p^{\bdot\bdot}}.\]
\item The region $\spadesuit$ is commutative by \cref{AFangm} for $(AF_j)_{p^{j\bdot}}$.
\item The region $\blacktriangle$ is commutative by the universal property of the pushout \cref{eq:smash-pushout} that defines the smash product.
\item The region $\clubsuit$ is commutative by the naturality of the pushout \cref{eq:smash-pushout}.
\item The region $\filledmedlozenge$ is commutative by \cref{AFangm}.
\end{itemize}
This proves that the diagram \cref{Agamma} is commutative on objects.

To prove that the diagram \cref{Agamma} is commutative on morphisms, we reuse the proof for objects above by replacing \cref{dcatganattr,AFangm} with, respectively, \cref{dcatgamod,Athetaangm}.
\end{proof}

Next is the main result of this chapter.

\begin{theorem}\label{thm:Acatmultifunctor}\index{Ga-category@$\Ga$-category!Cat-multifunctor@$\Cat$-multifunctor}\index{Cat-multifunctor@$\Cat$-multifunctor!Ga-category@$\Ga$-category}
There is a $\Cat$-multifunctor
\begin{equation}\label{Acatmultifunctor}
\begin{tikzcd}[column sep=large]
\Gacat \ar{r}{A} & \ACat
\end{tikzcd}
\end{equation}
with
\begin{itemize}
\item $\Gacat$ the $\Cat$-multicategory in \cref{thm:Gacatsmc};
\item $\ACat$ the $\Cat$-multicategory in \cref{thm:dcatcatmulticat} for $\cA$ in \cref{ex:mandellcategory};
\item object assignment in \cref{AXsmf};
\item arity 0 multimorphism functors in \cref{Aarityzeroiso}; and
\item positive arity multimorphism functors in \cref{Amultimorfunctor}.
\end{itemize}
Moreover, the $\Cat$-multifunctor $A$ extends the functor in \cref{AGacatACati}.
\end{theorem}

\begin{proof}
Consider the axioms of a $\Cat$-multifunctor in \cref{def:enr-multicategory-functor}.
\begin{itemize}
\item $A$ preserves the colored units because, for each $\Ga$-category $X$, each of $A1_X$ and $1_{AX}$ is componentwise an identity functor by \cref{def:dcatcatmulticat,AFangm,AXmdef}.
\item $A$ preserves the symmetric group action by \cref{Apreservessymmetry}.
\item $A$ preserves the composition by \cref{Apreservescomp}.
\end{itemize}
Therefore, $A$ is a $\Cat$-multifunctor.  The second assertion holds for the following reasons.
\begin{itemize}
\item The object assignment $X \mapsto AX$ is the one in \cref{AXsmf} in both cases.
\item The $n=1$ case of \cref{AFangm} agrees with \cref{AFm}.
\end{itemize} 
Thus, $A$ extends the functor in \cref{AGacatACati}.
\end{proof}

\chapter{Inverse \texorpdfstring{$K$}{K}-Theory is a Pseudo Symmetric \texorpdfstring{$\Cat$}{Cat}-Multifunctor}
\label{ch:invK}
The purpose of this chapter is to observe that the inverse $K$-theory functor 
\[\begin{tikzpicture}[xscale=1,yscale=1,vcenter]
\def\v{.6}
\draw[0cell=1]
(0,0) node (x1) {\Gacat}
(x1)++(2.5,0) node (x2) {\ACat}
(x2)++(2.8,0) node (x3) {\permcatsg}
(x1)++(0,\v) node[inner sep=0pt] (s) {}
(x3)++(0,\v) node[inner sep=0pt] (t) {}
;
\draw[1cell=.85] 
(x1) edge node {A} (x2)
(x2) edge node {\groa} (x3)
(x1) edge[-,shorten >=-1pt] (s)
(s) edge[-,shorten >=-1pt,shorten <=-1pt] node {\cP} (t)
(t) edge[shorten <=-1pt] (x3)
;
\end{tikzpicture}\]
extends to a pseudo symmetric $\Cat$-multifunctor; see \cref{thm:invKpseudosym}.  In fact, $A$ is a $\Cat$-multifunctor, and $\groa$ is a pseudo symmetric $\Cat$-multifunctor.  Moreover, there is a refinement \cref{invKAgroaU} 
\[\begin{tikzpicture}[xscale=1,yscale=1,vcenter]
\def\v{.6}
\draw[0cell=1]
(0,0) node (x1) {\Gacat}
(x1)++(2.5,0) node (x2) {\ACat}
(x2)++(2.5,0) node (x3) {\pfiba}
(x3)++(2.8,0) node (x4) {\permcatsg}
(x1)++(0,\v) node[inner sep=0pt] (s) {}
(x4)++(0,\v) node[inner sep=0pt] (t) {}
(x2)++(0,-\v) node[inner sep=0pt] (p) {}
(x4)++(0,-\v) node[inner sep=0pt] (q) {}
;
\draw[1cell=.85] 
(x1) edge node {A} (x2)
(x2) edge node {\groa} node[swap] {\sim} (x3)
(x3) edge node {\U} (x4)
(x1) edge[-,shorten >=-1pt] (s)
(s) edge[-,shorten >=-1pt,shorten <=-1pt] node {\cP} (t)
(t) edge[shorten <=-1pt] (x4)
(x2) edge[-,shorten >=-1pt] (p)
(p) edge[-,shorten >=-1pt,shorten <=-1pt] node[pos=.6] {\groa} (q)
(q) edge[shorten <=-1pt] (x4)
;
\end{tikzpicture}\]
that factors $\cP$ through a non-symmetric $\Cat$-\emph{multiequivalence} $\groa$ that takes into account the opfibration structure of the Grothendieck construction.  Here is a summary table.
\smallskip
\begin{center}
\resizebox{0.85\width}{!}{
{\renewcommand{\arraystretch}{1.4}%
{\setlength{\tabcolsep}{1ex}
\begin{tabular}{|c|c|c|}\hline
& (-) $\Cat$-multifunctor & ref \\ \hline
$\groa \cn \ACat \to \permcatsg$ & pseudo symmetric & \ref{grodpscatmultifunctor}\\ \hline
$\U \cn \pfiba \to \permcatsg$ & non-symmetric & \ref{pfibdtopermcatsg}\\ \hline
$\groa \cn \ACat \fto{\sim} \pfiba$ & non-symmetric $\Cat$-multiequivalence & \ref{thm:dcatpfibdeq}\\ \hline
$A \cn \Gacat \to \ACat$ & symmetric & \ref{Acatmultifunctor}\\ \hline
$\cP \cn \Gacat \to \permcatsg$ & pseudo symmetric & \ref{thm:invKpseudosym}\\ \hline
\end{tabular}}}}
\end{center}
\medskip

\subsection*{Preservation of Algebraic Structures}

As a pseudo symmetric $\Cat$-multifunctor, $\cP$ preserves algebraic structures parametrized by
\begin{itemize}
\item non-symmetric $\Cat$-multifunctors (\cref{cor:Ppreservesnonsym}) and
\item pseudo symmetric $\Cat$-multifunctors from $\Cat$-multicategories (\cref{cor:Ppreservespseudosym}).  
\end{itemize}
In particular, inverse $K$-theory $\cP$ preserves pseudo symmetric $E_\infty$-algebras (\cref{cor:Pefinitypsalg}).  We emphasize that inverse $K$-theory $\cP$ is \emph{not} a $\Cat$-multifunctor in the symmetric sense (\cref{def:enr-multicategory-functor}) because it does not strictly preserve the symmetric group action (\cref{expl:Ppsisos}).  In general, $\cP$ does not preserve algebraic structures parametrized by $\Cat$-multicategories, such as $E_\infty$-algebras (\cref{def:barratteccles}).  Instead, $E_\infty$-algebras in $\Gacat$ are sent by $\cP$ to pseudo symmetric $E_\infty$-algebras in $\permcatsg$ (\cref{ex:Peinftypseudo}).

\subsection*{Some History}

Inverse $K$-theory $\cP$ is important because it is a homotopy inverse of Segal $K$-theory\index{Segal K-theory@Segal $K$-theory} \cite{segal}; see \cite[Chapter 8]{cerberusIII} for a detailed discussion of Segal $K$-theory.  This homotopy inverse property of $\cP$ is one of the main results of Mandell \cite{mandell_inverseK}.  It is a sharpening of the result of Thomason \cite{thomason} that says that there is an equivalence between the stable homotopy categories of small permutative categories and of $\Ga$-categories.

Although Segal $K$-theory provides an equivalence of stable homotopy categories, it is \emph{not} a multifunctor, even in the non-symmetric sense.  In particular, Segal $K$-theory does not, in general, preserve algebraic structures such as monoids.  This fact is noted in \cite{elmendorf-mandell} and is one of the main reasons for the development of Elmendorf-Mandell $K$-theory $\Kem$, which is $\Cat$-multifunctorial.  A detailed discussion of $\Kem$ is in \cite[Chapters 10--13]{cerberusIII}.  In view of the non-multifunctoriality of Segal $K$-theory, the pseudo symmetric $\Cat$-multifunctoriality of inverse $K$-theory $\cP$ (\cref{thm:invKpseudosym}) is somewhat unexpected.

\subsection*{Organization}

To prepare for the discussion of inverse $K$-theory as a pseudo symmetric $\Cat$-multifunctor, in \cref{sec:invKfunctor} we review the inverse $K$-theory functor $\cP = \groa \circ A$ from the underlying category of $\Gacat$ to the underlying category of $\permcatsg$.

In \cref{sec:invKpseudosym} we prove the main \cref{thm:invKpseudosym} on the pseudo symmetric $\Cat$-multifunctoriality of inverse $K$-theory $\cP$ by combining some of the results in earlier chapters.  The associative operad, which parametrizes monoids, is a non-symmetric operad.  Thus $\cP$ preserves monoids (\cref{ex:Preservesnonsym}).  This means that $\cP$ sends each monoid in $\Gacat$ to a small tight ring category.

In \cref{sec:unravelingP} we explain in detail the pseudo symmetric $\Cat$-multifunctorial structure of $\cP$ by describing its multimorphism functors and pseudo symmetry isomorphisms.  In particular, the formulas \cref{groasigmaafmx,monensigmam} for the pseudo symmetry isomorphism $\cP_\sigma$ show that $\cP_\sigma$ is the identity if and only if $\sigma$ is the identity permutation.  In \cref{expl:invKrelated,expl:invKrelatedelm} we discuss related work in the literature \cite{elmendorf-multi-inverse,johnson-yau-invK} about the multifunctoriality of $\cP$.

To prepare for the discussion of $\cP$ preserving pseudo symmetric $E_\infty$-algebras, in \cref{sec:barratteccles} we review the categorical Barratt-Eccles $E_\infty$-operad $\BE$ and its Coherence \cref{mapfrombarratteccles}.  Small permutative categories and small (tight) bipermutative categories are $E_\infty$-algebras in $\Cat$, $\permcat$, and $\permcatsg$ (\cref{ex:Einftyalgebras}).  In \cref{ex:Einftygacat} we explain the structure of an $E_\infty$-algebra in $\Gacat$.

In \cref{sec:apppseudo} we define pseudo symmetric $E_\infty$-algebras as pseudo symmetric $\Cat$-multifunctors from the Barratt-Eccles operad $\BE$.  As a special case of the main \cref{thm:invKpseudosym}, inverse $K$-theory $\cP$ preserves pseudo symmetric $E_\infty$-algebras.  To better understand the structure involved, in \cref{expl:psEinftyalgebra,expl:psEinftyGacat,expl:psEinftypermcat} we unpack the data and axioms of pseudo symmetric $E_\infty$-algebras in a general $\Cat$-multicategory, $\Gacat$, $\permcat$, and $\permcatsg$.

We remind the reader of our left normalized bracketing \cref{expl:leftbracketing} for iterated monoidal product.

\section{The Inverse $K$-Theory Functor}
\label{sec:invKfunctor}

In this section we review the inverse $K$-theory functor $\cP$ due to Mandell \cite{mandell_inverseK}; see also \cite[Section 6]{johnson-yau-invK}.  The inverse $K$-theory functor is a composite of two functors---$A$ and $\groa$ with $\cA$ in \cref{ex:mandellcategory}---from $\Ga$-categories to small permutative categories.  Here is an outline of this section.
\begin{itemize}
\item The inverse $K$-theory functor $\cP = \groa \circ A$ is in \cref{def:mandellinvK}.  
\item \cref{expl:PX} is a detailed description of the functor $\cP$.
\item \cref{ex:Pterminal} shows that $\cP X$ recovers the permutative category $\Aplus$ when $X$ is the terminal $\Ga$-category.
\end{itemize}

\subsection*{Inverse $K$-Theory}

For a small tight bipermutative category $(\cD,\oplus,\otimes)$ (\cref{def:embipermutativecat}), the Grothendieck construction $\grod$ has the following properties.
\begin{itemize}
\item It sends a symmetric monoidal functor $X \cn \Dplus \to \Cat$ to a small permutative category $\grod X$ (\cref{grodxpermutative}).
\item It sends a monoidal natural transformation (\cref{expl:additivenattr} \eqref{expl:additivenattr-iv}) 
\[\begin{tikzpicture}[xscale=1,yscale=1,baseline={(x11.base)}]
\def\d{25} 
\draw[0cell=.9]
(0,0) node (a) {\Dplus}
(a)++(2.2,0) node (b) {\phantom{Z}}
(b)++(.1,0) node (x12) {\Cat}
;
\draw[1cell=.8] 
(a) edge[bend left=\d] node {X} (b)
(a) edge[bend right=\d] node[swap] {Z} (b)
;
\draw[2cell=1]
node[between=a and b at .47, rotate=-90, 2label={above,\phi}] {\Rightarrow}
;
\end{tikzpicture}\]
between symmetric monoidal functors to a strict symmetric monoidal functor $\grod\phi \cn \grod X \to \grod Z$ (\cref{rk:grodphifunctorunary}).
\end{itemize} 
Thus the Grothendieck construction (\cref{thm:grocatmultifunctor}) restricts to a functor
\begin{equation}\label{grodone}
\begin{tikzpicture}[xscale=1,yscale=1,baseline={(x11.base)}]
\draw[0cell=.9]
(0,0) node (x11) {\DCat}
(x11)++(2.7,0) node (x12) {\permcatsg};
\draw[1cell=.9]  
(x11) edge node {\grod} (x12);
\end{tikzpicture}
\end{equation}
between the underlying categories in \cref{ex:DCati,ex:permcatsgone}.  The following definition is due to Mandell \cite[4.3 and 4.5]{mandell_inverseK}.  

\begin{definition}\label{def:mandellinvK}\index{inverse K-theory@inverse $K$-theory!functor}\index{functor!inverse K-theory@inverse $K$-theory}
The \emph{inverse $K$-theory functor} $\cP$ is defined as the composite functor
\begin{equation}\label{mandellinvK}
\begin{tikzpicture}[xscale=2.75,yscale=1,vcenter]
\draw[0cell=.9]
(0,0) node (x11) {\Gacat}
(x11)++(1,0) node (x12) {\permcatsg}
(x11)++(.42,-1) node (x21) {\ACat}
;
\draw[1cell=.9] 
(x11) edge node {\cP} (x12)
(x11) edge node[swap,pos=.3] {A} (x21)
(x21) edge node[swap,pos=.7] {\groa} (x12)
;
\end{tikzpicture}
\end{equation}
with the following data.
\begin{itemize}
\item $\Gacat$, $\ACat$, and $\permcatsg$ are the underlying categories of the respective $\Cat$-multicategories in \cref{ex:Gacati,ex:DCati,ex:permcatsgone}.
\item $A$ is the functor in \cref{AGacatACati}.
\item $\groa$ is the functor in \cref{grodone} with $\cD = \cA$ in \cref{ex:mandellcategory}.
\end{itemize}
This finishes the definition of the functor $\cP$.
\end{definition}

\begin{explanation}[Codomain of $\cP$]\label{expl:Pcodomain}
For the inverse $K$-theory functor \cref{mandellinvK}
\[\begin{tikzpicture}[xscale=1,yscale=1,baseline={(x11.base)}]
\draw[0cell=.9]
(0,0) node (x11) {\Gacat}
(x11)++(4,0) node (x12) {\permcatsg};
\draw[1cell=.9]  
(x11) edge node {\cP = \groa A(-)} (x12);
\end{tikzpicture}\]
we can restrict the codomain to the subcategory of $\permcatsg$ consisting of \emph{strict} symmetric monoidal functors.  The reason is that each monoidal natural transformation in $\ACat$ is sent by the Grothendieck construction
\[\begin{tikzpicture}[xscale=1,yscale=1,baseline={(x11.base)}]
\draw[0cell=.9]
(0,0) node (x11) {\ACat}
(x11)++(2.7,0) node (x12) {\permcatsg};
\draw[1cell=.9]  
(x11) edge node {\groa} (x12);
\end{tikzpicture}\]
to a strict symmetric monoidal functor (\cref{rk:grodphifunctorunary}).  However, when we consider the pseudo symmetric $\Cat$-multifunctoriality of $\cP$ in \cref{sec:invKpseudosym}, the linearity constraints are natural isomorphisms but not identities in general (\cref{expl:invKarityn}). 
\end{explanation}

\begin{explanation}[Inverse $K$-Theory of $\Ga$-Categories]\label{expl:PX}\index{inverse K-theory@inverse $K$-theory!Ga-category@$\Ga$-category}\index{Ga-category@$\Ga$-category!inverse K-theory@inverse $K$-theory}
We unpack the inverse $K$-theory functor $\cP = \groa A(-)$ in \cref{mandellinvK}.  Suppose $X$ is a $\Ga$-category.
\begin{description}
\item[Underlying Category]
The inverse $K$-theory
\begin{equation}\label{mandellPX}
\cP X = \groa AX
\end{equation}
is the Grothendieck construction of the functor $AX \cn \cA \to \Cat$ in \cref{functorAX}.  Interpreting \cref{def:grothendieckconst} in this context, an \emph{object} in the category $\cP X$ is a pair $(m,x)$ with
\begin{itemize}
\item $m = \ang{m_i}_{i=1}^p$ an object in $\cA$ (\cref{ex:mandellcategory}) and
\item $x = \ang{x_i}_{i=1}^p$ an object in $(AX)(m) = \prod_{i=1}^p X\ord{m_i}$ in \cref{AXmdef}.
\end{itemize}
A \emph{morphism} 
\begin{equation}\label{pofxmorphism}
\begin{tikzcd}[column sep=large]
(m,x) \ar{r}{(\phi, f)} & (n,y)
\end{tikzcd} \inspace \cP X
\end{equation}
consists of
\begin{itemize}
\item a morphism 
\[\begin{tikzcd}[column sep=large]
m = \ang{m_i}_{i=1}^p \ar{r}{\phi} & n = \ang{n_j}_{j=1}^q
\end{tikzcd} \inspace \cA\]
as in \cref{Amorphismphi} and 
\item a morphism 
\[f = \ang{f_j}_{j=1}^q \cn \phi_*(x) \to y \inspace (AX)(n) = \txprod_{j=1}^q X\ord{n_j}\]
with 
\[\phi_* = (AX)(\phi) \cn (AX)(m) \to (AX)(n)\]
as in \cref{phistarAphi}.
\end{itemize}
Composition is defined as
\[(\psi,g) \circ (\phi,f) = \big(\psi \phi, g \circ (\psi_* f)\big).\]
The \emph{identity morphism} of an object $(m,x)$ is the pair $(1_m,1_x)$ of identity morphisms.  
\item[Permutative Structure]
Interpreting \cref{def:permutativegroconst} in this context, $\cP X$ is a permutative category with strict \emph{monoidal unit} given by the pair 
\begin{equation}\label{PXmonoidalunit}
\big(\ang{},*\big) \in \cP X
\end{equation}
consisting of
\begin{itemize}
\item the empty tuple $\ang{}$ in $\cA$ and
\item the unique object $* \in (AX)\ang{} = X\ord{0} = \boldone$.
\end{itemize} 
The \emph{monoidal product} on objects is given by
\begin{equation}\label{PXmonoidalprod}
(m,x) \Box (n,y) = \begin{cases}
\big((m,n), (x,y)\big) & \text{if $p,q > 0$},\\
(m,x) & \text{if $q=0$, and}\\
(n,y) & \text{if $p=0$}
\end{cases}
\end{equation}
with
\[\left\{\begin{split}
(m,n) &= m \oplus n = \brb{\ang{m_i}_{i=1}^p, \ang{n_j}_{j=1}^q} \in \cA \andspace\\
(x,y) &= \brb{\ang{x_i}_{i=1}^p, \ang{y_j}_{j=1}^q} \in (AX)(m,n).
\end{split}\right.\]
The monoidal product on morphisms is defined similarly, using the fact that $\Aplus$ is a permutative category.  

The \emph{braiding} in $\cP X$ is given componentwise by the isomorphism\label{not:PXbetabox}
\[\begin{tikzpicture}[xscale=1,yscale=.9,vcenter]
\draw[0cell=.9]
(0,0) node (x11) {(m,x) \Box (n,y)}
(x11)++(4.5,0) node (x12) {(n,y) \Box (m,x)}
(x11)++(0,-1) node (x21) {\big((m,n), (x,y)\big)}
(x12)++(0,-1) node (x22) {\big((n,m), (y,x)\big)}
;
\draw[1cell=.9]  
(x11) edge node {\betabox_{(m,x), (n,y)}} (x12)
(x11) edge[equal, transform canvas={xshift=.2ex}] (x21)
(x12) edge[equal, transform canvas={xshift=-.2ex}] (x22)
(x21) edge node {\srb{\betaplus_{m,n},1}} node[swap] {\iso} (x22);
\end{tikzpicture}\]
with $\betaplus_{m,n}$ the additive braiding \cref{Abetaplusmn} in $\cA$.  This finishes the description of the permutative category $\cP X = \groa AX$ for a $\Ga$-category $X$.
\item[Morphisms]
Interpreting \cref{def:groconarityn} in the current context, for a morphism $F \cn X \to Y$ of $\Ga$-categories, the strict symmetric monoidal functor\label{not:PF}
\[\begin{tikzpicture}[xscale=1,yscale=1,baseline={(x11.base)}]
\draw[0cell=.9]
(0,0) node (x11) {\cP X = \groa AX}
(x11)++(4.2,0) node (x12) {\cP Y = \groa AY};
\draw[1cell=.9]  
(x11) edge node {\cP F = \groa AF} (x12);
\end{tikzpicture}\]
is given on objects by
\begin{equation}\label{PFmx}
\begin{split}
(\cP F)(m,x) &= \big(m, (AF)_m x\big)\\ 
&= \Big(m, \bang{F_{\ord{m_i}}(x_i)}_{i=1}^p \Big)
\end{split}
\end{equation}
and on morphisms by
\[\begin{split}
(\cP F)(\phi, f) &= \big(\phi, (AF)_n f\big)\\
&= \Big(\phi, \bang{F_{\ord{n_j}}(f_j)}_{j=1}^q \Big).
\end{split}\]
The \emph{unit constraint} is the identity morphism
\[(\cP F)^0 = 1_{(\ang{},*)} \inspace \cP Y.\]  
The \emph{monoidal constraint}
\[\begin{tikzcd}[column sep=huge]
(\cP F)(m,x) \Box (\cP F)(n,y) \ar{r}{(\cP F)^2} & (\cP F)\big( (m,x) \Box (n,y) \big)
\end{tikzcd}\]
is the identity morphism in $\cP Y$ of the following object, using \cref{PXmonoidalprod} to make the appropriate adjustment if either $p=0$ or $q=0$.
\[\begin{aligned}
&(\cP F)\big( (m,x) \Box (n,y) \big) &&\\
&= (\cP F) \big((m,n), (x,y)\big) & \phantom{=} & \text{by \cref{PXmonoidalprod}}\\
&= \Big((m,n) \scs \big( \bang{F_{\ord{m_i}}(x_i)}_{i=1}^p \scs \bang{F_{\ord{n_j}}(y_j)}_{j=1}^q  \big) \Big) & \phantom{=} & \text{by \cref{PFmx}}\\
&= \Big(m, \bang{F_{\ord{m_i}}(x_i)}_{i=1}^p \Big) \Box \Big(n, \bang{F_{\ord{n_j}}(y_j)}_{j=1}^q \Big) & \phantom{=} & \text{by \cref{PXmonoidalprod}}\\
&= (\cP F)(m,x) \Box (\cP F)(n,y) & \phantom{=} & \text{by \cref{PFmx}}
\end{aligned}\]
\end{description}
This finishes the description of the inverse $K$-theory functor $\cP = \groa A(-)$.
\end{explanation}

\begin{example}[Terminal $\Ga$-Category]\label{ex:Pterminal}\index{Ga-category@$\Ga$-category!terminal}
Suppose $X$ is the $\Ga$-category defined by
\[X\ord{n} = \boldone \foreachspace \ord{n} \in \Fskel.\]
By definition \cref{AXmdef} there is an isomorphism of categories 
\[(AX)(m) \iso \boldone \foreachspace m \in \cA.\]
By \cref{expl:PX} there is a strict symmetric monoidal isomorphism
\[\begin{tikzcd}[/tikz/column 1/.append style={anchor=base east},
/tikz/column 2/.append style={anchor=base west}, column sep=large, row sep=-.6ex]
\cP X \ar{r}{\iso} & \Aplus\\
(m,*) \ar[mapsto]{r} & m\\
(\phi, 1_*) \ar[mapsto]{r} & \phi
\end{tikzcd}\]
with $\Aplus = (\cA, \oplus, \ang{}, \betaplus)$ the additive structure of $\cA$ in \cref{ex:mandellcategory}.
\end{example}

\section{Pseudo Symmetric \texorpdfstring{$\Cat$}{Cat}-Multifunctorial Inverse \texorpdfstring{$K$}{K}-Theory}
\label{sec:invKpseudosym}

The purpose of this section is to prove that the inverse $K$-theory functor $\cP$ in \cref{mandellinvK} extends to a pseudo symmetric $\Cat$-multifunctor; see \cref{thm:invKpseudosym}.  More specifically we obtain two factorizations of the pseudo symmetric $\Cat$-multifunctor $\cP$.
\begin{enumerate}
\item In \cref{invKAgroa} $\cP$ factors as
\begin{itemize}
\item a $\Cat$-multifunctor $A$ followed by
\item a pseudo symmetric $\Cat$-multifunctor $\groa$ given by the Grothendieck construction on the small tight bipermutative category $\cA$.
\end{itemize} 
\item In the refinement \cref{invKAgroaU}, $\cP$ factors as the composite of
\begin{itemize}
\item a $\Cat$-multifunctor $A$,
\item a non-symmetric $\Cat$-\emph{multiequivalence} $\groa$, and 
\item a non-symmetric $\Cat$-multifunctor $\U$.
\end{itemize}  
\end{enumerate}
We emphasize that $\cP$ is \emph{not} a $\Cat$-multifunctor because it does not strictly preserve the symmetric group action.  See \cref{expl:invKrelated,expl:invKrelatedelm} for related work in the literature about the multifunctoriality of inverse $K$-theory.

Here is an outline of this section.
\begin{itemize}
\item \cref{thm:invKpseudosym} proves the pseudo symmetric $\Cat$-multifunctoriality of inverse $K$-theory $\cP$.
\item In \cref{cor:Ppreservesnonsym} we observe that $\cP$ preserves any algebraic structure parametrized by a non-symmetric $\Cat$-multifunctor, including the associative operad (\cref{ex:Preservesnonsym}).
\item In \cref{cor:Ppreservespseudosym} we observe that $\cP$ preserves any algebraic structure parametrized by a pseudo symmetric $\Cat$-multifunctor.  Examples related to the Barratt-Eccles $E_\infty$-operad are discussed in \cref{sec:barratteccles,sec:apppseudo}.
\end{itemize}

\subsection*{Main Theorem}

Recall from \cref{def:pseudosmultifunctor} that a pseudo symmetric $\Cat$-multifunctor
\begin{itemize}
\item strictly preserves the colored units \cref{enr-multifunctor-unit} and composition \cref{v-multifunctor-composition} and
\item is equipped with pseudo symmetry isomorphisms that satisfy the four coherence axioms \cref{pseudosmf-unity,pseudosmf-product,pseudosmf-topeq,pseudosmf-boteq}.
\end{itemize} 
Each $\Cat$-multifunctor is a pseudo symmetric $\Cat$-multifunctor with pseudo symmetry isomorphisms given by  identities (\cref{ex:idpsm}).  Composition of pseudo symmetric $\Cat$-multifunctors is well defined, associative, and unital (\cref{pseudosmfcomposite}).

Recall the $\Cat$-multicategories
\begin{itemize}
\item $\Gacat$ in \cref{thm:Gacatsmc},
\item $\ACat$ in \cref{thm:dcatcatmulticat} with $\cA$ in \cref{ex:mandellcategory}, and
\item $\permcatsg$ in \cref{thm:permcatenrmulticat}.
\end{itemize}

\begin{theorem}\label{thm:invKpseudosym}\index{inverse K-theory@inverse $K$-theory!pseudo symmetric Cat-multifunctor@pseudo symmetric $\Cat$-multifunctor}\index{pseudo symmetric!Cat-multifunctor@$\Cat$-multifunctor!inverse K-theory@inverse $K$-theory}
The following statements hold for the inverse $K$-theory functor $\cP$ in \cref{mandellinvK}.
\begin{enumerate}[label=(\roman*)]
\item\label{thm:invKpseudosym-i} $\cP$ extends to a pseudo symmetric $\Cat$-multifunctor
\begin{equation}\label{invKAgroa}
\begin{tikzpicture}[xscale=2.75,yscale=1,vcenter]
\draw[0cell=.9]
(0,0) node (x11) {\Gacat}
(x11)++(1,0) node (x12) {\permcatsg}
(x11)++(.42,-1) node (x21) {\ACat}
;
\draw[1cell=.9] 
(x11) edge node {\cP} (x12)
(x11) edge node[swap,pos=.3] {A} (x21)
(x21) edge node[swap,pos=.7] {\groa} (x12)
;
\end{tikzpicture}
\end{equation}
with
\begin{itemize}
\item $A$ the $\Cat$-multifunctor in \cref{Acatmultifunctor} and
\item $\groa$ the pseudo symmetric $\Cat$-multifunctor in \cref{grodpscatmultifunctor}.
\end{itemize}
\item\label{thm:invKpseudosym-ii} The pseudo symmetric $\Cat$-multifunctor $\cP$ in \cref{invKAgroa} is not a $\Cat$-multifunctor.
\item\label{thm:invKpseudosym-iii} As a \index{non-symmetric!Cat-multifunctor@$\Cat$-multifunctor!inverse K-theory@inverse $K$-theory}non-symmetric $\Cat$-multifunctor, $\cP$ in \cref{invKAgroa} factors as the boundary in
\begin{equation}\label{invKAgroaU}
\begin{tikzpicture}[xscale=3.3,yscale=1,vcenter]
\def\h{.1} \def\v{-1.2}
\draw[0cell=.9]
(0,0) node (x11) {\Gacat}
(x11)++(1,0) node (x12) {\permcatsg}
(x11)++(\h,\v) node (x21) {\ACat}
(x12)++(-\h,\v) node (x22) {\pfiba}
;
\draw[1cell=.85] 
(x11) edge[bend left=15] node {\cP} (x12)
(x11) edge node[swap,pos=.3] {A} (x21)
(x21) edge[transform canvas={xshift=-1ex}] node[pos=.4] {\groa} (x12)
(x21) edge[bend right=20] node[pos=.6] {\groa} node[swap,pos=.6] {\sim} (x22)
(x22) edge node[swap,pos=.7] {\U} (x12)
;
\end{tikzpicture}
\end{equation}
with the following data.
\begin{itemize}
\item $\pfiba$ is the non-symmetric $\Cat$-multicategory in \cref{ex:pfibdcatmulticat}.
\item The bottom $\groa$ is the non-symmetric $\Cat$-multiequivalence\index{non-symmetric!Cat-multiequivalence@$\Cat$-multiequivalence!inverse K-theory@inverse $K$-theory} in \cref{thm:dcatpfibdeq}.
\item $\U$ is the non-symmetric $\Cat$-multifunctor in \cref{pfibdtopermcatsg}.
\end{itemize}
\end{enumerate} 
\end{theorem}

\begin{proof}
In \cref{invKAgroa} $A$ is a $\Cat$-multifunctor, and $\groa$ is a pseudo symmetric $\Cat$-multifunctor.  Thus the composite $\groa \circ A$ is a pseudo symmetric $\Cat$-multifunctor by \cref{ex:idpsm,pseudosmfcomposite}.  Moreover, the following two statements hold.
\begin{itemize}
\item The $\Cat$-multifunctor $A$ extends the functor $A$ in the inverse $K$-theory functor by \cref{thm:Acatmultifunctor}.
\item The pseudo symmetric $\Cat$-multifunctor $\groa$ extends the functor $\groa$ in the inverse $K$-theory functor by definition \cref{grodone}.
\end{itemize}
So the composite pseudo symmetric $\Cat$-multifunctor $\groa \circ A$ extends the inverse $K$-theory functor $\cP$ in \cref{mandellinvK}.  This proves assertion \cref{thm:invKpseudosym-i}.

To prove assertion \cref{thm:invKpseudosym-ii}, consider an object $F \in \Gacat\scmap{\ang{X};Z}$ as in \cref{Fxonexntoz}.  The two strong $n$-linear functors
\[\begin{tikzpicture}[xscale=1,yscale=1,baseline={(x11.base)}]
\def\a{25}
\draw[0cell=.85]
(0,0) node (x11) {\txprod_{j=1}^n (\groa AX_{\sigma(j)})}
(x11)++(1,0) node (x12) {\phantom{Z}}
(x12)++(3.5,0) node (x13) {\phantom{Z}}
(x13)++(.3,.03) node (x14) {\groa AZ} 
;
\draw[1cell=.8]  
(x12) edge[transform canvas={yshift=.7ex}] node {\cP (F^\sigma) = \groa (AF)^\sigma} (x13)
(x12) edge[transform canvas={yshift=-.3ex}] node[swap] {(\cP F)^\sigma = (\groa AF)^\sigma} (x13)
;
\end{tikzpicture}\]
are different in general.  Indeed, by \cref{grodsigmadomain} applied to $\cD = \cA$, the object assignments of $\cP(F^\sigma)$ and $(\cP F)^\sigma$ are different in general because the multiplicative braiding in $\cA$, which is defined in \cref{Amtensorn}, is not the identity.  Thus $\cP$ does not strictly preserve the symmetric group action, so it is not a $\Cat$-multifunctor.

Assertion \cref{thm:invKpseudosym-iii} follows from the factorizations in \cref{invKAgroa} and \cref{ex:dcatgrodpfibdeq}.
\end{proof}

\subsection*{Preservation of Non-Symmetric Algebras}

Each pseudo symmetric $\Cat$-multifunctor (\cref{def:pseudosmultifunctor}) is, in particular, a non-symmetric $\Cat$-multifunctor (\cref{def:enr-multicategory-functor}) after forgetting the pseudo symmetry isomorphisms (\cref{ex:nscatmultifunctors}).  Thus, \cref{thm:invKpseudosym} \cref{thm:invKpseudosym-i} yields the following result, which is first recorded in \cite[Section 10]{johnson-yau-invK}.

\begin{corollary}\label{cor:Ppreservesnonsym}\index{inverse K-theory@inverse $K$-theory!preserves non-symmetric algebras}
Suppose
\begin{itemize}
\item $\Q$ is a non-symmetric $\Cat$-multicategory and
\item $F \cn \Q \to \Gacat$ is a non-symmetric $\Cat$-multifunctor.
\end{itemize}
Then the composite
\begin{equation}\label{QPnonsymmetric}
\begin{tikzcd}[column sep=large]
\Q \ar{r}{F} & \Gacat \ar{r}{\cP} & \permcatsg
\end{tikzcd}
\end{equation}
is a non-symmetric $\Cat$-multifunctor.
\end{corollary}

\cref{cor:Ppreservesnonsym} means that any algebraic structure in $\Gacat$ parametrized by a non-symmetric $\Cat$-multifunctor is preserved by inverse $K$-theory $\cP$.  

\begin{example}[Monoids to Ring Categories]\label{ex:Preservesnonsym}
Consider the \emph{non-symmetric associative operad}\index{non-symmetric!associative operad} $\Asns$ (\cref{def:enr-multicategory}) with each 
\[\Asns_n = \{*\},\]
a one-element set.  It is regarded as a non-symmetric $\Cat$-multicategory with one object and each $\Asns_n = \boldone$, the terminal category.  See \cite[Section 14.2]{yau-operad} and \cite[Section 11.1]{cerberusIII} for a detailed discussion of the free operad\index{associative operad}
\begin{equation}\label{Asoperad}
\As = \big\{ \Sigma_n \big\}_{n \geq 0}
\end{equation}
generated by $\Asns$.  In $\Gacat$ and $\permcatsg$, $\Asns$ parametrizes the following algebraic structures.
\begin{itemize}
\item By \cite[14.2.18]{yau-operad} or \cite[11.1.15]{cerberusIII} a monoid in $\Gacat$ is uniquely determined by a non-symmetric $\Cat$-multifunctor 
\[\Asns \to \Gacat.\]
\item By \cite[3.4]{elmendorf-mandell} or \cite[11.2.16]{cerberusIII} a small tight ring category (\cref{def:ringcat}) is uniquely determined by a non-symmetric $\Cat$-multifunctor 
\[\Asns \to \permcatsg.\]
\end{itemize}
These facts and \cref{cor:Ppreservesnonsym} imply that inverse $K$-theory $\cP$ sends each monoid\index{monoid!Ga-category@$\Ga$-category}\index{Ga-category@$\Ga$-category!monoid} in $\Gacat$ to a small tight ring category.\index{ring category!inverse K-theory@inverse $K$-theory}\index{inverse K-theory@inverse $K$-theory!ring category}  This is first proved in \cite[10.3]{johnson-yau-invK}.
\end{example}

\subsection*{Preservation of Pseudo Symmetric Algebras}

By \cref{pseudosmfcomposite} pseudo symmetric $\Cat$-multifunctors are closed under composition.  Thus, \cref{thm:invKpseudosym} \cref{thm:invKpseudosym-i} yields the following result.

\begin{corollary}\label{cor:Ppreservespseudosym}\index{inverse K-theory@inverse $K$-theory!preserves pseudo symmetric algebras}\index{pseudo symmetric!Cat-multifunctor@$\Cat$-multifunctor!preservation by inverse K-theory@preservation by inverse $K$-theory}
Suppose
\begin{itemize}
\item $\Q$ is a $\Cat$-multicategory and
\item $F \cn \Q \to \Gacat$ is a pseudo symmetric $\Cat$-multifunctor.
\end{itemize}
Then the composite
\begin{equation}\label{QPpseudosymmetric}
\begin{tikzcd}[column sep=large]
\Q \ar{r}{F} & \Gacat \ar{r}{\cP} & \permcatsg
\end{tikzcd}
\end{equation}
is a pseudo symmetric $\Cat$-multifunctor.
\end{corollary}

\cref{cor:Ppreservespseudosym} means that any algebraic structure in $\Gacat$ parametrized by a pseudo symmetric $\Cat$-multifunctor is preserved by inverse $K$-theory $\cP$.  The case where $\Q$ is the Barratt-Eccles $E_\infty$-operad is discussed in \cref{sec:barratteccles,sec:apppseudo}. 

\begin{example}[$\Cat$-Multifunctors]\label{ex1:Ppreservespseudosym}
By \cref{ex:idpsm} each $\Cat$-multifunctor becomes a pseudo symmetric $\Cat$-multifunctor when it is equipped with identity pseudo symmetry isomorphisms.  In this case, \cref{cor:Ppreservespseudosym} says that, if $F \cn \Q \to \Gacat$ is a $\Cat$-multifunctor, then the composite 
\[\begin{tikzcd}[column sep=huge]
\Q \ar{r}{\cP \circ F} & \permcatsg
\end{tikzcd}\]
in \cref{QPpseudosymmetric} is a pseudo symmetric $\Cat$-multifunctor.
\end{example}

\section{Unraveling the Pseudo Symmetric \texorpdfstring{$\Cat$}{Cat}-Multifunctorial Inverse $K$-Theory}
\label{sec:unravelingP}

In this section we explain in detail the pseudo symmetric $\Cat$-multifunctor $\cP$ in \cref{thm:invKpseudosym}.
\begin{itemize}
\item The arity 0 multimorphism functors of $\cP$ are discussed in \cref{expl:invKarityzero}.
\item The positive arity multimorphism functors of $\cP$ are discussed in \cref{expl:invKarityn}.
\item The pseudo symmetry isomorphisms $\cP_{\sigma}$ of $\cP$ are discussed in \cref{expl:Ppsisos}.  The displays \cref{groasigmaafmx,monensigmam} show that $\cP$ does \emph{not} strictly preserve the symmetric group action because the multiplicative braiding in $\cA$ is not the identity. 
\item \cref{expl:invKrelated,expl:invKrelatedelm} discuss related work in the literature about the multifunctoriality of $\cP$.
\end{itemize}

For each $\Ga$-category $X$, the inverse $K$-theory $\cP X = \groa AX$ as a small permutative category is discussed in detail in \cref{expl:PX}.  For $\Ga$-categories 
\[Z, X_1, \ldots, X_n \cn (\Fskel,\ord{0}) \to (\Cat,\boldone)\]
with $\angX = \ang{X_j}_{j=1}^n$ and $n \geq 0$, next we describe
\begin{itemize}
\item the multimorphism functor
\begin{equation}\label{invKmultifunctor}
\begin{tikzcd}[column sep=large]
\Gacat\scmap{\angX;Z} \ar{r}{\cP} & \permcatsg\scmap{\ang{\cP X};\cP Z}
\end{tikzcd}
\end{equation}
in \cref{expl:invKarityzero,expl:invKarityn} and
\item the pseudo symmetry isomorphisms in \cref{expl:Ppsisos}
\end{itemize}
for the pseudo symmetric $\Cat$-multifunctor $\cP = \groa \circ A$ in \cref{invKAgroa}.

\begin{explanation}[Inverse $K$-Theory in Arity 0]\label{expl:invKarityzero}
If $n=0$ then the multimorphism functor $\cP$ in \cref{invKmultifunctor} is given by the assignments 
\begin{equation}\label{Parityzero}
\begin{tikzpicture}[xscale=1,yscale=1,baseline={(x11.base)}]
\def\v{-.8} \def\u{-.7} \def\w{-.6}
\draw[0cell=.9]
(0,0) node (x11) {\Gacat\scmap{\ang{};Z}}
(x11)++(3,0) node (x12) {\ACat\scmap{\ang{};AZ}}
(x12)++(3.5,0) node (x13) {\permcatsg\scmap{\ang{};\cP Z}}
(x11)++(0,\v) node (x21) {\GMap(\ftu,Z)}
(x12)++(0,\v) node (x22) {(AZ)(1) = Z\ord{1}}
(x13)++(0,\v) node (x23) {\cP Z = \groa AZ}
(x21)++(0,\u) node (x31) {F}
(x22)++(0,\u) node (x32) {AF = F_{\ord{1}}(1)}
(x23)++(0,\u) node (x33) {\cP F = \brb{(1), F_{\ord{1}}(1)}}
(x31)++(0,\w) node (x41) {\theta}
(x32)++(0,\w) node (x42) {A\theta = \theta_{\ord{1},1}}
(x33)++(0,\w) node (x43) {\cP\theta = \brb{1_{(1)}, \theta_{\ord{1},1}}}
(x11)++(0,.6) node[inner sep=0pt] (s) {}
(x13)++(0,.6) node[inner sep=0pt] (t) {}
;
\draw[1cell=.85] 
(x11) edge node {A} node[swap] {\iso} (x12)
(x12) edge node {\groa} (x13)
(x11) edge[equal] (x21)
(x12) edge[equal] (x22)
(x13) edge[equal] (x23)
(x31) edge[|->] (x32)
(x32) edge[|->] (x33)
(x41) edge[|->] (x42)
(x42) edge[|->] (x43)
(x11) edge[-,shorten >=-1pt] (s)
(s) edge[-, shorten <=-1pt, shorten >=-1pt] node {\cP} (t)
(t) edge[shorten <=-1pt] (x13)
;
\end{tikzpicture}
\end{equation}
for objects $F$ and morphisms $\theta$ in $\Gacat\scmap{\ang{};Z}$.  The detail of \cref{Parityzero} is as follows.
\begin{itemize}
\item The equality under $\Gacat\scmap{\ang{};Z}$ is \cref{Gacatemptyz}.
\item The two equalities under $\ACat\scmap{\ang{};AZ}$ are from \cref{dcatbipermzeroary,AXmdef}.
\item The two equalities under $\permcatsg\scmap{\ang{};\cP Z}$ are from \cref{permcatemptyd,mandellinvK}.
\item The isomorphism $A$ is from \cref{Aarityzeroiso}.
\item The functor $\groa$ is from \cref{def:groconarityzero}.
\end{itemize}
The functor $\cP = \groa\circ A$ in \cref{Parityzero} is injective on objects and morphism sets.
\end{explanation}

\cref{not:cAtensorobjects} is used in \cref{expl:invKarityn,expl:Ppsisos} below.

\begin{explanation}[Inverse $K$-Theory in Positive Arity]\label{expl:invKarityn}
If $n > 0$ then the multimorphism functor $\cP = \groa \circ A$ in \cref{invKmultifunctor} is the composite
\begin{equation}\label{Parityn}
\begin{tikzpicture}[xscale=1,yscale=1,baseline={(x11.base)}]
\def\v{-.8} \def\u{-.7} \def\w{-.6}
\draw[0cell=.8]
(0,0) node (x11) {\Gacat\scmap{\ang{X};Z}}
(x11)++(3.2,0) node (x12) {\ACat\scmap{\ang{AX};AZ}}
(x12)++(3.8,0) node (x13) {\permcatsg\scmap{\ang{\cP X};\cP Z}}
(x11)++(0,\v) node (x21) {\GMap\brb{\txsma_{j=1}^n X_j,Z}}
(x11)++(0,.6) node[inner sep=0pt] (s) {}
(x13)++(0,.6) node[inner sep=0pt] (t) {}
;
\draw[1cell=.85] 
(x11) edge node {A} (x12)
(x12) edge node {\groa} (x13)
(x11) edge[equal] (x21)
(x11) edge[-,shorten >=-1pt] (s)
(s) edge[-, shorten <=-1pt, shorten >=-1pt] node {\cP} (t)
(t) edge[shorten <=-1pt] (x13)
;
\end{tikzpicture}
\end{equation}
with
\begin{itemize}
\item $\Gacat\scmap{\ang{X};Z} = \GMap\brb{\txsma_{j=1}^n X_j,Z}$ in \cref{Gacatnmorphism}, 
\item $\ACat\scmap{\ang{AX};AZ}$ in \cref{dcatnarycategory} for $\cD = \cA$ in \cref{ex:mandellcategory},
\item $\permcatsg\scmap{\ang{\cP X};\cP Z}$ in \cref{permcatsgangcd},
\item the functor $A$ in \cref{def:Aaritynfunctorobj,def:Aaritynfunctormor}, and
\item the functor $\groa$ in \cref{def:groconarityn,def:groconmodification}.
\end{itemize}
We first describe the object assignment of $\cP$ and then its morphism assignment.

\emph{Object Assignment}.  Suppose $F \in \Gacat\scmap{\ang{X};Z}$ is an object, which means a pointed natural transformation $F \cn \txsma_{j=1}^n X_j \to Z$ as in \cref{Fxonexntoz}.  By \cref{def:groconarityn,def:Aaritynfunctorobj,expl:AFangmformulas}, the strong $n$-linear functor (\cref{def:nlinearfunctor})
\[\begin{tikzcd}[column sep=huge]
\txprod_{j=1}^n \cP X_j = \txprod_{j=1}^n (\groa AX_j) \ar{r}{\cP F} & \cP Z = \groa AZ
\end{tikzcd}\]
is given as follows.  
\begin{description}
\item[Objects] $\cP F = \groa AF$ sends an object
\begin{equation}\label{angmjcaxj}
\lrang{\brb{m^j \in \cA , x^j \in (AX_j)m^j}} \in \txprod_{j=1}^n (\groa AX_j)
\end{equation}
to the object
\[(\cP F)\bang{(m^j , x^j)} = \left( m^{1,\ldots,n} \scs (AF)_{\angm} \ang{x^j}_{j=1}^n \right) \in \groa AZ\]
consisting of the object in \cref{monen} 
\[m^{1,\ldots,n} = \txotimes_{j=1}^n m^j \in \cA\]
and the object in \cref{AFangmangxj}
\[(AF)_{\angm} \ang{x^j}_{j=1}^n \in (AZ)m^{1,\ldots,n}.\] 
\item[Morphisms] $\cP F = \groa AF$ sends a morphism
\[\bang{(f_j,p_j)} \cn \bang{(m^j,x^j)} \to \bang{(k^j,w^j)} \inspace \txprod_{j=1}^n (\groa AX_j)\]
to the morphism
\[\begin{tikzpicture}[xscale=1,yscale=1,vcenter]
\draw[0cell=.85]
(0,0) node (x11) {(\cP F)\bang{(m^j , x^j)} = \left( m^{1,\ldots,n} \scs (AF)_{\angm} \ang{x^j}_{j=1}^n \right)}
(x11)++(0,-1.4) node (x21) {(\cP F)\bang{(k^j , w^j)} = \left( k^{1,\ldots,n} \scs (AF)_{\angk} \ang{w^j}_{j=1}^n \right)}
;
\draw[1cell=.85] 
(x11) edge[shorten <=-.2ex, shorten >=-.2ex, transform canvas={xshift=-2.5em}] node {= \brb{\txotimes_{j=1}^n f_j, (AF)_{\angk} \ang{p_j}_{j=1}^n}} node [swap] {(\cP F)\bang{(f_j \scs p_j)}} (x21)
;
\end{tikzpicture}\]
with the morphism in \cref{AFangmangxj}
\[(AF)_{\angk} \ang{p_j}_{j=1}^n \in (AZ)k^{1,\ldots,n}.\]
\item[Linearity Constraints] By \cref{grodphilinearity}, for each $j \in \{1,\ldots,n\}$, the $j$-th linearity constraint $(\cP F)^2_j$ of $\cP F$ is componentwise an isomorphism $(\lapinv,1)$ with $\lap$ the unique Laplaza coherence isomorphism in \cref{laplazaapp}.
\end{description}

\emph{Morphism Assignment}.  Suppose $\theta \cn F \to G$ is a morphism in $\Gacat\scmap{\ang{X};Z}$, which means a pointed modification with component pointed natural transformations as in \cref{thetaordmcomponent}.  By \cref{def:groconmodification} the $n$-linear natural transformation (\cref{def:nlineartransformation})\label{not:Ptheta}
\[\begin{tikzpicture}[xscale=1,yscale=1,baseline={(x1.base)}]
\def\a{25}
\draw[0cell=.9]
(0,0) node (x1) {\txprod_{j=1}^n (\groa AX_j)}
(x1)++(.8,.01) node (x2) {\phantom{D}}
(x2)++(3.3,0) node (x3) {\phantom{D}}
(x3)++(.3,0) node (x4) {\groa AZ}
;
\draw[1cell=.8]  
(x2) edge[bend left=\a] node[pos=.5] {\cP F = \groa AF} (x3)
(x2) edge[bend right=\a] node[swap,pos=.5] {\cP G = \groa AG} (x3)
;
\draw[2cell] 
node[between=x2 and x3 at .25, rotate=-90, 2label={above,\cP\theta = \groa A\theta}] {\Rightarrow}
;
\end{tikzpicture}\]
has, for each object $\bang{(m^j,x^j)}$ as in \cref{angmjcaxj}, a component morphism
\[\begin{tikzpicture}[xscale=1,yscale=1,vcenter]
\draw[0cell=.85]
(0,0) node (x11) {(\cP F)\bang{(m^j , x^j)} = \left( m^{1,\ldots,n} \scs (AF)_{\angm} \ang{x^j}_{j=1}^n \right)}
(x11)++(0,-1.4) node (x21) {(\cP G)\bang{(m^j , x^j)} = \left( m^{1,\ldots,n} \scs (AG)_{\angm} \ang{x^j}_{j=1}^n \right)}
;
\draw[1cell=.85] 
(x11) edge[shorten <=-.2ex, shorten >=-.2ex, transform canvas={xshift=-2.5em}] node {= \big(1_{m^{1,\ldots,n}} \scs \big((A\theta)_{\angm}\big)_{\ang{x^j}_{j=1}^n} \big)} node [swap] {(\cP \theta)_{\ang{(m^j,x^j)}}} (x21)
;
\end{tikzpicture}\]
in $\groa AZ$.  Here
\[\big((A\theta)_{\angm}\big)_{\ang{x^j}_{j=1}^n} \in (AZ)m^{1,\ldots,n}\] 
is the morphism in \cref{Athetaangmangxj}.  This finishes our description of the multimorphism functor $\cP = \groa \circ A$ in \cref{Parityn}.
\end{explanation}

\begin{explanation}[Pseudo Symmetry Isomorphisms of Inverse $K$-Theory]\label{expl:Ppsisos}
Here we describe the pseudo symmetry isomorphisms \cref{Fsigmaangccp} of the pseudo symmetric $\Cat$-multifunctor $\cP = \groa \circ A$ in \cref{invKAgroa}.  By definition \cref{pseudosymisopasting}, for each permutation $\sigma \in \Sigma_n$, the pseudo symmetry isomorphism $\cP_\sigma$ is given by the following pasting of natural isomorphisms.
\begin{equation}\label{Psigmapasting}
\begin{tikzpicture}[xscale=3,yscale=1.3,vcenter]
\def\d{.5}
\draw[0cell=.8]
(0,0) node (x11) {\Gacat\scmap{\angX;Z}}
(x11)++(1,0) node (x12) {\ACat\scmap{\ang{AX};AZ}}
(x12)++(1.3,0) node (x13) {\permcatsg\scmap{\ang{\cP X};\cP Z}}
(x11)++(0,-1) node (x21) {\Gacat\scmap{\angX\sigma;Z}}
(x12)++(0,-1) node (x22) {\ACat\scmap{\ang{AX}\sigma;AZ}}
(x13)++(0,-1) node (x23) {\permcatsg\scmap{\ang{\cP X}\sigma;\cP Z}}
(x11)++(0,\d) node[inner sep=0pt] (s) {}
(x13)++(0,\d) node[inner sep=0pt] (t) {}
(x21)++(0,-\d) node[inner sep=0pt] (s2) {}
(x23)++(0,-\d) node[inner sep=0pt] (t2) {}
;
\draw[1cell=.8]
(x11) edge node (f1) {A} (x12)
(x12) edge node (g1) {\groa} (x13)
(x21) edge node[swap] (f2) {A} (x22)
(x22) edge node[swap] (g2) {\groa} (x23)
(x11) edge node[swap] {\sigma} (x21)
(x12) edge node {\sigma} (x22)
(x13) edge node {\sigma} (x23)
(x11) edge[-,shorten >=-1pt] (s)
(s) edge[-, shorten <=-1pt, shorten >=-1pt] node {\cP} (t)
(t) edge[shorten <=-1pt] (x13)
(x21) edge[-,shorten >=-1pt] (s2)
(s2) edge[-, shorten <=-1pt, shorten >=-1pt] node {\cP} (t2)
(t2) edge[shorten <=-1pt] (x23)
;
\draw[2cell=.9]
node[between=f1 and f2 at .5] {=}
node[between=g1 and g2 at .58, 2label={above,(\groa)_{\sigma}}, 2label={below,\iso}] {\Longrightarrow}
;
\end{tikzpicture}
\end{equation} 
The detail of \cref{Psigmapasting} is as follows.
\begin{itemize}
\item By \cref{thm:Acatmultifunctor} $A$ is a $\Cat$-multifunctor, with the preservation of the symmetric group action proved in \cref{Apreservessymmetry}.  As in \cref{ex:idpsm}, $A$ is regarded as a pseudo symmetric $\Cat$-multifunctor with identity pseudo symmetry isomorphisms.  So the left square in \cref{Psigmapasting} is strictly commutative.
\item The pseudo symmetry isomorphism $(\groa)_{\sigma}$ is from \cref{def:gropseudosymmetry} for $\cD = \cA$ in \cref{ex:mandellcategory}.  By \cref{grodsigmacomponent} the component isomorphisms of $(\groa)_{\sigma}$ are always of the form $(\sigma,1)$. 
\end{itemize}

We unpack \cref{Psigmapasting} further.  Suppose $F \in \Gacat\scmap{\angX;Z}$ is an object, which means a pointed natural transformation $F \cn \txsma_{j=1}^n X_j \to Z$ as in \cref{Fxonexntoz}.  The component of the pseudo symmetry isomorphism $(\groa)_{\sigma}$ at the object 
\[AF \in \ACat\scmap{\ang{AX};AZ}\] 
as in \cref{angAXAFAZ} is an $n$-linear natural isomorphism (\cref{def:nlineartransformation}) as follows.
\[\begin{tikzpicture}[xscale=1,yscale=1,baseline={(x11.base)}]
\def\a{25}
\draw[0cell=.85]
(0,0) node (x11) {\txprod_{j=1}^n (\groa AX_{\sigma(j)})}
(x11)++(1,0) node (x12) {\phantom{Z}}
(x12)++(3,.03) node (x13) {\phantom{Z}}
(x13)++(.3,0) node (x14) {\groa AZ} 
;
\draw[1cell=.8]  
(x12) edge[bend left=\a] node {\cP (F^\sigma) = \groa (AF)^\sigma} (x13)
(x12) edge[bend right=\a] node[swap] {(\cP F)^\sigma = (\groa AF)^\sigma} (x13)
;
\draw[2cell=.9]
node[between=x12 and x13 at .35, rotate=-90, 2label={above,(\groa)_{\sigma,AF}}] {\Rightarrow}
;
\end{tikzpicture}\]
To describe this natural isomorphism explicitly, consider an object 
\[\bang{(m^j,x^j)} \in \txprod_{j=1}^n (\groa AX_{\sigma(j)})\]
with, for each $j \in \{1,\ldots,n\}$, objects $m^j \in \cA$ as in \cref{not:cAtensorobjects} and 
\[x^j = \bang{x^j_{i_j}}_{i_j=1}^{r_j} \in \txprod_{i_j=1}^{r_j} X_{\sigma(j)} \ord{m^j_{i_j}} = (AX_{\sigma(j)})m^j.\]
By \cref{grodsigmadomain,grodsigmacomponent} the component of $(\groa)_{\sigma, AF}$ at the object $\ang{(m^j,x^j)}$ is the following isomorphism in $\groa AZ$.
\begin{equation}\label{groasigmaafmx}
\begin{tikzpicture}[xscale=1,yscale=1,vcenter]
\draw[0cell=.85]
(0,0) node (x11) {\phantom{\cP (F^\sigma)\bang{(m^j , x^j)} = \left(m^{1,\ldots,n} \scs \sigmainv_* (AF)_{\sigma\angm} \bang{x^{\sigmainv(j)}}_{j=1}^n \right)}}
(x11)++(-.3,0) node (x10) {\cP (F^\sigma)\bang{(m^j , x^j)} = \left(m^{1,\ldots,n} \scs \sigmainv_* (AF)_{\sigma\angm} \bang{x^{\sigmainv(j)}}_{j=1}^n \right)}
(x11)++(0,-1.4) node (x21) {(\cP F)^\sigma\bang{(m^j , x^j)} = \left(m^{\sigmainv(1),\ldots,\sigmainv(n)} \scs (AF)_{\sigma\angm} \bang{x^{\sigmainv(j)}}_{j=1}^n \right)}
;
\draw[1cell=.85] 
(x11) edge[shorten <=-.5ex, shorten >=-.5ex, transform canvas={xshift=-5em}] node {= (\sigma,1)} node [swap] {(\groa)_{\sigma, AF, \ang{(m^j,x^j)}}} (x21)
;
\end{tikzpicture}
\end{equation}
In \cref{groasigmaafmx} the object
\[(AF)_{\sigma\angm} \bang{x^{\sigmainv(j)}}_{j=1}^n 
= \lrang{ \,\cdots\, \Bang{ \scalebox{.9}{$F_{\bang{m^{\sigmainv(j)}_{i_{\sigmainv(j)}}}_{j=1}^n} \bang{x^{\sigmainv(j)}_{i_{\sigmainv(j)}}}_{j=1}^n$} }_{i_{\sigmainv(1)}=1}^{r_{\sigmainv(1)}} \,\cdots\, }_{i_{\sigmainv(n)}=1}^{r_{\sigmainv(n)}}\]
belonging to
\[(AZ) m^{\sigmainv(1),\ldots,\sigmainv(n)} 
= \txprod_{i_{\sigmainv(n)}=1}^{r_{\sigmainv(n)}} \,\cdots\, \txprod_{i_{\sigmainv(1)}=1}^{r_{\sigmainv(1)}}\, \Big(Z \ord{m^{\sigmainv(1), \ldots, \sigmainv(n)}_{i_{\sigmainv(1)}, \ldots, i_{\sigmainv(n)}}} \Big)\]
is as in \cref{AFangmangxj} for the objects
\[\left\{\begin{split}
\sigma\angm &= \bang{m^{\sigmainv(j)}}_{j=1}^n \in \cA^n \andspace\\
\bang{x^{\sigmainv(j)}}_{j=1}^n & \in \txprod_{j=1}^n (AX_j) m^{\sigmainv(j)}.
\end{split}\right.\]
The isomorphism
\begin{equation}\label{monensigmam}
\begin{tikzcd}[column sep=large]
m^{1,\ldots,n} = \txotimes_{j=1}^n m^j \ar{r}{\sigma}[swap]{\iso} 
& \txotimes_{j=1}^n m^{\sigmainv(j)} = m^{\sigmainv(1),\ldots,\sigmainv(n)}
\end{tikzcd}
\end{equation}
in \cref{groasigmaafmx} permutes the $n$ factors $\ang{m^j}_{j=1}^n$ according to $\sigma \in \Sigma_n$ using the multiplicative braiding $\betate$ in $\cA$ \cref{Abetate}.  

The isomorphism $(\sigma,1)$ in \cref{groasigmaafmx} is the identity morphism if and only if $\sigma \in \cA$ in \cref{monensigmam} is the identity morphism.  Thus the pseudo symmetry isomorphism 
\begin{equation}\label{Psigmagroasigma}
\cP_\sigma = (\groa)_{\sigma, A(-)}
\end{equation}
for inverse $K$-theory $\cP$ is the identity natural transformation if and only if $\sigma \in \Sigma_n$ is the identity permutation.
\end{explanation}

\subsection*{Related Literature}

Next we discuss the relationship between \cref{thm:invKpseudosym} and \cite{elmendorf-multi-inverse,johnson-yau-invK} on the multifunctoriality of the inverse $K$-theory functor $\cP$.

\begin{explanation}[Related Work of Johnson-Yau]\label{expl:invKrelated}
The main result in \cite{johnson-yau-invK} says that inverse $K$-theory $\cP = \groa \circ A$ is a non-symmetric $\Cat$-multifunctor.  Both \cite{johnson-yau-invK} and \cref{thm:invKpseudosym} consider $\Cat$-enrichment.   \cref{thm:invKpseudosym} is an extension and a refinement of the main result in \cite{johnson-yau-invK} in the following ways:
\begin{enumerate}[label=(\roman*)]
\item While $\cP$ does not strictly preserve the symmetric group action, \cref{thm:invKpseudosym} shows that $\cP$ is a pseudo symmetric $\Cat$-multifunctor.  In particular, $\cP$ preserves the symmetric group action up to the pseudo symmetry isomorphisms $\cP_\sigma$ in \cref{Psigmapasting}, which satisfy the four coherence axioms in \cref{def:pseudosmultifunctor}. 
\item The proof in \cite{johnson-yau-invK} is a direct analysis of the functor $\cP$ and does not consider the two functors $A$ and $\groa$ separately.  On the other hand, \cref{thm:invKpseudosym} shows that each of the two constituent functors of $\cP = \groa \circ A$ is $\Cat$-multifunctorial: in the symmetric sense for $A$ and in the pseudo symmetric sense for $\groa$.
\item \cref{thm:invKpseudosym} contains the further factorization \cref{invKAgroaU} of $\cP$ involving the non-symmetric $\Cat$-\emph{multiequivalence} $\groa$ that takes into account the permutative opfibration structure of the Grothendieck construction on $\cA$.\defmark
\end{enumerate}
\end{explanation}

\begin{explanation}[Related Work of Elmendorf]\label{expl:invKrelatedelm}
The main claim in \cite{elmendorf-multi-inverse} is that inverse $K$-theory $\cP$ is a composite of three multifunctors in the symmetric sense, and, therefore, $\cP$ is also a multifunctor in the symmetric sense.  That work differs from \cref{thm:invKpseudosym} and \cite{johnson-yau-invK} in the following ways:
\begin{enumerate}[label=(\roman*)]
\item $\Cat$-enrichment is not considered in \cite{elmendorf-multi-inverse}.
\item The constructions in \cite{elmendorf-multi-inverse} involve an indexing category $\bbN$ that is different from $\cA$. 
\item The three multifunctors in \cite[Theorems 1--3]{elmendorf-multi-inverse} are not the ones in either \cref{invKAgroa} or \cref{invKAgroaU}. 
\item The description of the monoidal unit in $\Gacat$ in \cite[Section 4]{elmendorf-multi-inverse}, denoted $B$ there, is not correct.  The correct monoidal unit is $\ftu$ in \cref{Gacatunit,expl:Ftu}.
\end{enumerate}
Furthermore, the two main assertions in \cite{elmendorf-multi-inverse} are not correct, as we discuss next.  

Contrary to the claim in \cite[Corollary 5]{elmendorf-multi-inverse}, $\cP$ is \emph{not} a multifunctor in the symmetric sense because it does not strictly preserve the symmetric group action.  This is noted in \cite[Remark 8.1]{johnson-yau-invK} and described in detail in \cref{expl:Ppsisos}.  In particular, \cref{groasigmaafmx} shows how $\cP(-)^\sigma$ and $(\cP (-))^\sigma$ are different on objects.

Contrary to the claim in \cite[Observation 4]{elmendorf-multi-inverse}, the composite of the three constructions in \cite[Theorems 1--3]{elmendorf-multi-inverse} does \emph{not} yield $\cP$.  For a general $\Ga$-category $X$, that composite contains many objects that are not in $\cP X$ for the following reasons.
\begin{itemize}
\item The second construction in \cite[Theorem 2]{elmendorf-multi-inverse} is a wreath product that creates a basepoint corresponding to the object $0 \in \bbN$ and the unique object $* \in X\ord{0} = \boldone$. 
\item The third construction in \cite[Theorem 3]{elmendorf-multi-inverse} creates its own basepoint corresponding to the empty sequence, which should correspond to the monoidal unit in $\cP X$ as in \cref{PXmonoidalunit}.
\end{itemize}   
However, that composite construction also yields finite sequences $\ang{*}_{j=1}^n$ containing $n$ copies of $* \in X\ord{0}$ created by the wreath product.  The sequences $\ang{*}_{j=1}^n$ for $n \geq 1$ are not objects in $\cP X$ as in \cref{expl:PX}.  Similarly, finite sequences containing $*$ and other objects are also not supposed to be in $\cP X$.  For example, the composite construction in \cite[Theorems 1--3]{elmendorf-multi-inverse} is already different from $\cP X$ on objects for the terminal $\Ga$-category $X$ in \cref{ex:Pterminal}.
\end{explanation}

\section{The Barratt-Eccles Operad and $E_\infty$-Algebras}
\label{sec:barratteccles}

By \cref{cor:Ppreservespseudosym} inverse $K$-theory $\cP$ preserves any algebraic structure parametrized by a pseudo symmetric $\Cat$-multifunctor from a $\Cat$-multicategory.  In particular, this holds for pseudo symmetric $\Cat$-multifunctors from $E_\infty$-operads.  To explain in detail these pseudo symmetric $E_\infty$-algebras in $\Gacat$ and $\permcatsg$, in this section we first recall the categorical Barratt-Eccles $E_\infty$-operad.
\begin{itemize}
\item \cref{mapfrombarratteccles} is a coherence theorem for $E_\infty$-algebras.
\item Small permutative categories and small (tight) bipermutative categories are $E_\infty$-algebras in, respectively, $\Cat$, $\permcat$, and $\permcatsg$ (\cref{ex:Einftyalgebras}).
\item In \cref{ex:Einftygacat} we explain in detail $E_\infty$-algebras in $\Gacat$.
\end{itemize}
In \cref{sec:apppseudo} we define the pseudo symmetric variant of $E_\infty$-algebras using pseudo symmetric $\Cat$-multifunctors from the Barratt-Eccles operad.

The following definition of the Barratt-Eccles operad is from \cite[page 178]{elmendorf-mandell}.  For a detailed discussion, see \cite[Sections 11.4--11.6]{cerberusIII}, which uses the notation $E\As$ instead of $\BE$.

\begin{definition}\label{def:barratteccles}
The \emph{Barratt-Eccles operad}\index{Barratt-Eccles operad}\index{operad!Barratt-Eccles}
\[\BE = \big\{E\Sigma_n\big\}_{n \geq 0}\] 
is the $\Cat$-operad (\cref{def:enr-multicategory}) defined by the following data.
\begin{itemize}
\item It has one object $*$.
\item The category $E\Sigma_n$ has object set $\Sigma_n$.  For each pair of permutations $\sigma,\tau \in \Sigma_n$, the morphism set $E\Sigma_n(\sigma,\tau)$ is the one-element set $\{*\}$.
\item The \emph{$*$-colored unit} is the object $\id_1 \in \Sigma_1$.
\item The \emph{symmetric group action} on objects is the group multiplication 
\[\Sigma_n \times \Sigma_n \to \Sigma_n.\]
\item The \emph{composition} on objects is defined by group multiplication
\begin{equation}\label{Asgamma}
\ga\scmap{\sigma; \ang{\sigma_j}_{j=1}^n} 
= \sigmatil \cdot \sigmatimes \in \Sigma_{k_1 + \cdots + k_n}
\end{equation}
for $\sigma \in \Sigma_n$ and $\sigma_j \in \Sigma_{k_j}$ using the following notation.  The block permutation
\begin{equation}\label{sigmatilde}
\sigmatil = \sigma\ang{k_1,\ldots,k_n} \in \Sigma_{k_1 + \cdots + k_n}
\end{equation}
permutes $n$ consecutive blocks of lengths $k_1, \ldots, k_n$ according to $\sigma$ \cref{blockpermutation}.  The block sum
\begin{equation}\label{sigmatimes}
\sigmatimes = \sigma_1 \times \cdots \times \sigma_n \in \Sigma_{k_1 + \cdots + k_n}
\end{equation}
is as in \cref{blocksum}.
\item On morphisms the symmetric group action and composition are uniquely defined by the property that each morphism set in each $E\Sigma_n$ is a one-element set.
\end{itemize}
This finishes the definition of $\BE$.  

Moreover, for a $\Cat$-multicategory $\M$, a $\Cat$-multifunctor 
\[\BE \to \M\]
is called an \emph{$E_\infty$-algebra}\index{Einfty-algebra@$E_\infty$-algebra} in $\M$.
\end{definition}

The following coherence theorem for $E_\infty$-algebras is \cite[6.3.3]{fresse1}.  See also \cite[11.4.14]{cerberusIII} for a detailed proof and \cite[19.2.1]{yau-inf-operad} for an extension involving other types of equivariant structure.  The transposition that permutes $i$ and $j$ is denoted $(i,j)$, so $(1,2)$\label{not:onetwopermutation} is the non-identity permutation in $\Sigma_2$.

\begin{theorem}[$E_\infty$-Coherence]\label{mapfrombarratteccles}\index{Einfty-algebra@$E_\infty$-algebra!coherence}\index{operad!Barratt-Eccles!coherence}\index{Barratt-Eccles operad!coherence}
Suppose $(\M,\gamma,\operadunit)$ is a $\Cat$-multicategory.  Then an $E_\infty$-algebra in $\M$
\[\begin{tikzcd}[column sep=large]
\BE \ar{r}{F} & \M
\end{tikzcd}\]
is uniquely determined by
\begin{itemize}
\item the object $x = F(*)$ in $\M$;
\item the 0-ary 1-cell 
\[\tu = F(\id_0) \in \M\smscmap{\ang{};x};\]
\item the 2-ary 1-cell 
\[\mu = F(\id_2) \in \M\smscmap{x,x;x}; \andspace\]
\item the 2-ary 2-cell isomorphism 
\[\begin{tikzcd}[column sep=huge]
\mu \ar{r}{\xi \,=\, F(\tau)}[swap]{\iso} & \muop = \mu(1,2)
\end{tikzcd}  \inspace \M\scmap{x,x;x}\]
with $\tau \in E\Sigma_2\big(\id_2, (1,2)\big)$ the unique isomorphism. 
\end{itemize}
The above data are subject to the conditions \cref{eastopunityassociativity,eastopsymmetry,eastopunity,eastophexagon} below.
\begin{description}
\item[1-Cell Unity and Associativity] 
The following object equalities hold.
\begin{equation}\label{eastopunityassociativity}
\left\{\begin{split}
\ga\scmap{\mu;\tu, \operadunit_x} &= \operadunit_x = \ga\scmap{\mu;\operadunit_x,\tu} \in \M\smscmap{x;x}\\
\ga\scmap{\mu;\mu, \operadunit_x} &= \ga\scmap{\mu;\operadunit_x, \mu} \in \M\scmap{x,x,x;x}
\end{split}\right.
\end{equation}
\item[Symmetry] The following diagram in $\M\smscmap{x,x;x}$ is commutative.
\begin{equation}\label{eastopsymmetry}
\begin{tikzcd}[column sep=small, row sep=small]
\mu \ar{rr}{1_{\mu}} \ar[shorten >=-.5ex]{dr}[swap,pos=.3]{\xi} && \mu\\
& \muop \ar[shorten <=-.5ex]{ur}[swap,pos=.7]{\xi (1,2)} &
\end{tikzcd}
\end{equation}
\item[Unity] The following diagram in $\M\smscmap{x;x}$ is commutative.
\begin{equation}\label{eastopunity}
\begin{tikzcd}[column sep=huge, cells={nodes={scale=.85}},
every label/.append style={scale=.9}]
\operadunit_x = \ga\scmap{\mu;\operadunit_x,\tu} \ar[shift left]{r}{1_{\operadunit_x}} 
\ar[shift right]{r}[swap]{\ga\smscmap{\xi;1_{\operadunit_x},1_{\tu}}} &
\ga\scmap{\mu;\tu,\operadunit_x} = \ga\scmap{\muop;\operadunit_x,\tu} = \operadunit_x
\end{tikzcd}
\end{equation}
\item[Hexagon] The following diagram\index{hexagon axiom!Barratt-Eccles operad} in $\M\smscmap{x,x,x;x}$ is commutative.
\begin{equation}\label{eastophexagon}
\begin{tikzpicture}[xscale=4.5,yscale=1.2,baseline={(x21.base)}]
\def\h{1} \def\g{.05} \def\v{-1} \def\u{-.7}
\draw[0cell=.8]
(0,0) node (x11) {\ga\scmap{\mu;\mu,\operadunit_x}}
(x11)++(\h-2*\g,0) node (x12) {\ga\scmap{\muop;\mu,\operadunit_x}}
(x11)++(-\g,\u) node (x21) {\ga\scmap{\mu;\operadunit_x,\mu}}
(x21)++(\h,0) node (x22) {\ga\scmap{\mu;\muop,\operadunit_x} (2,3)}
(x21)++(4*\g,\v) node (x31) {\ga\scmap{\mu;\operadunit_x,\muop}}
(x31)++(\h-8*\g,0) node (x32) {\ga\scmap{\mu;\mu,\operadunit_x} (2,3)}
;
\draw[1cell=.8] 
(x21) edge[equal] (x11)
(x11) edge node {\ga\smscmap{\xi; 1_\mu, 1_{\operadunit_x}}} (x12)
(x12) edge[equal] (x22)
(x21) edge node[swap,pos=.2] {\ga\scmap{1_\mu; 1_{\operadunit_x},\xi}} (x31)
(x31) edge[equal] (x32)
(x32) edge node[swap,pos=.8] {\ga\scmap{1_\mu; \xi, 1_{\operadunit_x}} (2,3)} (x22)
; 
\end{tikzpicture}
\end{equation}
\end{description}
\end{theorem}

\begin{explanation}[$E_\infty$-Coherence]\label{expl:Einftyalgebras}
A key aspect of the proof of \cref{mapfrombarratteccles} is that, in addition to strictly preserving the unit and composition, a $\Cat$-multifunctor $F \cn \BE \to \M$ also strictly preserves the symmetric group action.  The last property is not available if $F$ is merely a pseudo symmetric $\Cat$-multifunctor (\cref{def:pseudosmultifunctor}).  Therefore, the pseudo symmetric variant of \cref{mapfrombarratteccles} in \cref{expl:psEinftyalgebra} below has more data and more axioms.
\end{explanation}

\begin{explanation}[Underlying Monoid of an $E_\infty$-Algebra]\label{expl:Einftymonoid}\index{Einfty-algebra@$E_\infty$-algebra!underlying monoid}\index{monoid!Einfty-algebra@$E_\infty$-algebra}
In \cref{mapfrombarratteccles} the triple $(x,\mu,\tu)$ is a monoid in $\M$ in the sense that the unity and associativity axioms \cref{eastopunityassociativity} are satisfied.  This is a consequence of the $\Cat$-multifunctor
\[\As \to \BE\]
from the associative operad \cref{Asoperad} to the Barratt-Eccles operad given by the identity function on $\Sigma_n$ for $n$-ary 1-cells.  The unit, symmetric group action, and composition in $\As$ are as in \cref{def:barratteccles}.
\end{explanation}

\begin{example}[(Bi)permutative Categories as $E_\infty$-Algebras]\label{ex:Einftyalgebras}
The following examples are applications of \cref{mapfrombarratteccles} to the $\Cat$-multicategories $\Cat$ (\cref{ex:cat}), $\permcat$, and $\permcatsg$ (\cref{thm:permcatenrmulticat}).
\begin{enumerate}
\item\label{ex:Einftyalgebras-i}\index{Einfty-algebra@$E_\infty$-algebra!in Cat@in $\Cat$}\index{permutative category!Einfty-algebra@$E_\infty$-algebra} By \cite[11.4.26]{cerberusIII} an $E_\infty$-algebra in $\Cat$ is precisely a small permutative category (\cref{def:symmoncat}).
\item\label{ex:Einftyalgebras-ii}\index{Einfty-algebra@$E_\infty$-algebra!in permcat@in $\permcat$}\index{bipermutative category!Einfty-algebra@$E_\infty$-algebra} By \cite[3.8]{elmendorf-mandell} or \cite[11.5.5]{cerberusIII} an $E_\infty$-algebra in $\permcat$ is precisely a small bipermutative category (\cref{def:embipermutativecat}).
\item\label{ex:Einftyalgebras-iii}\index{Einfty-algebra@$E_\infty$-algebra!in permcatsg@in $\permcatsg$} Restricting the proof of \cite[11.5.5]{cerberusIII} to $\permcatsg$ shows that an $E_\infty$-algebra in $\permcatsg$ is precisely a small \emph{tight} bipermutative category.\defmark
\end{enumerate}
\end{example}

\begin{example}[$E_\infty$-Algebras in $\Gacat$]\label{ex:Einftygacat}\index{Einfty-algebra@$E_\infty$-algebra!in Gacat@in $\Gacat$}\index{Ga-category@$\Ga$-category!Einfty-algebra@$E_\infty$-algebra}
Recall the $\Cat$-multicategory structure on $\Gacat$ (\cref{sec:Gacatmulticat}).  Applying \cref{mapfrombarratteccles}, an $E_\infty$-algebra in $\Gacat$ is a quadruple\label{not:EinfGacat}
\[(X,\tu,\mu,\xi)\]
consisting of the data \cref{Xgacat,oneoneone,mumncomponent,ximnbraiding} below.
\begin{description}
\item[Underlying Object] It consists of a $\Ga$-category (\cref{def:gammacategory}) 
\begin{equation}\label{Xgacat}
\begin{tikzcd}[column sep=large]
(\Fskel,\ord{0}) \ar{r}{X} & (\Cat,\boldone).
\end{tikzcd}
\end{equation}
\item[Unit] It is equipped with a unit object 
\[\tu \in \GMap(\ftu,X).\] 
This means a pointed natural transformation $\tu \cn \ftu \to X$ with $\ftu$ the monoidal unit diagram \cref{Gacatunit}.  As discussed in \cref{Gacatemptyz} through \cref{Fordnjziotaj}, such an object $\tu$ is uniquely determined by the object
\begin{equation}\label{oneoneone}
\tu_{\ord{1}}(1) \in X\ord{1}.
\end{equation}
\item[Multiplication] It is equipped with a 2-ary 1-cell
\[\mu \in \Gacat\scmap{X,X;X} = \GMap\brb{X \sma X, X}.\]
This means a pointed natural transformation (\cref{expl:pointedday})
\[\begin{tikzcd}[column sep=large]
X \sma X \ar{r}{\mu} & X
\end{tikzcd}\]
with $\sma$ the pointed Day convolution \cref{pointedday}.  It has component pointed functors
\begin{equation}\label{mumncomponent}
\begin{tikzcd}[column sep=large]
X\ord{m} \sma X\ord{n} \ar{r}{\mu_{\ord{m},\ord{n}}} & X\ord{mn}
\end{tikzcd}
\end{equation}
that are natural in $\ord{m}, \ord{n} \in \Fskel$ in the sense of the naturality diagram \cref{ptdaynaturality}.
\item[Braiding] It is equipped with a 2-ary 2-cell isomorphism 
\[\begin{tikzcd}[column sep=large]
\mu \ar{r}{\xi}[swap]{\iso} & \mu \circ \xiday
\end{tikzcd} \inspace \Gacat\scmap{X,X;X}\]
as follows, with $\xiday$ the braiding for the pointed Day convolution \cref{eq:Fxi}.
\[\begin{tikzpicture}[xscale=1,yscale=1,vcenter]
\def\h{2.5} \def\a{25}
\draw[0cell=.9]
(0,0) node (x11) {X \sma X}
(x11)++(\h,0) node (x12) {X \sma X}
(x11)++(\h/2,-1) node (x21) {X}
;
\draw[1cell=.85] 
(x11) edge node (xi) {\xiday} (x12)
(x12) edge[bend left=\a] node[pos=.3] {\mu} (x21)
(x11) edge[bend right=\a] node[swap,pos=.3] {\mu} (x21)
; 
\draw[2cell=1]
node[between=xi and x21 at .65, 2label={above,\xi}] {\Longrightarrow}
;
\end{tikzpicture}\]
As discussed in \cref{thetaordmcomponent}, this means that $\xi$ is a pointed modification with component pointed natural isomorphisms for $\ord{m}, \ord{n} \in \Fskel$ as follows.
\begin{equation}\label{ximnbraiding}
\begin{tikzpicture}[xscale=1,yscale=1,vcenter]
\def\h{3} \def\v{-1.3}
\draw[0cell=.9]
(0,0) node (x11) {X\ord{m} \sma X\ord{n}}
(x11)++(\h,0) node (x12) {X\ord{n} \sma X\ord{m}}
(x11)++(0,\v) node (x21) {X\ord{mn}}
(x12)++(0,\v) node (x22) {X\ord{nm}}
;
\draw[1cell=.85] 
(x11) edge node (s) {\xisma} (x12)
(x12) edge node {\mu_{\ord{n},\ord{m}}} (x22)
(x22) edge node (t) {X(\xisma_{\ord{m},\ord{n}})^\inv} (x21)
(x11) edge node[swap] {\mu_{\ord{m},\ord{n}}} (x21)
; 
\draw[2cell=.9]
node[between=s and t at .55, 2label={above,\xi_{\ord{m},\ord{n}}}, 2label={below,\iso}] {\Longrightarrow}
;
\end{tikzpicture}
\end{equation}
\begin{itemize}
\item The top $\xisma$ in \cref{ximnbraiding} is the braiding for the smash product \cref{eq:smash-pushout}, which is induced by the braiding for the Cartesian product.
\item The bottom 
\[\begin{tikzcd}[column sep=large]
\ord{mn} \ar{r}{\xisma_{\ord{m},\ord{n}}}[swap]{\iso} & \ord{nm}
\end{tikzcd}\]
in \cref{ximnbraiding} is the multiplicative braiding in $\Fskel$ \cref{xisma}.
\end{itemize}
These component pointed natural isomorphisms satisfy the modification axiom \cref{thetaordmodification}, which means naturality in $\ord{m}, \ord{n} \in \Fskel$.
\end{description}

In terms of the above data \cref{Xgacat,oneoneone,mumncomponent,ximnbraiding}, the conditions \cref{eastopunityassociativity,eastopsymmetry,eastopunity,eastophexagon} in \cref{mapfrombarratteccles} are as follows for $m,n,p \geq 0$.
\begin{description}
\item[1-Cell Unity and Associativity]
The 1-cell unity condition in \cref{eastopunityassociativity} means the equalities of functors
\[\mu_{\ord{1},\ord{n}}\brb{\tu_{\ord{1}}(1),-} 
= 1_{X\ord{n}} 
= \mu_{\ord{n},\ord{1}}\brb{-,\tu_{\ord{1}}(1)}.\] 
The 1-cell associativity condition in \cref{eastopunityassociativity} means that the following diagram of functors is commutative.
\begin{equation}\label{onecellassociativity}
\begin{tikzpicture}[xscale=1,yscale=1,vcenter]
\def\h{4} \def\v{-1.3}
\draw[0cell=.9]
(0,0) node (x11) {X\ord{m} \sma X\ord{n} \sma X\ord{p}}
(x11)++(\h,0) node (x12) {X\ord{m} \sma X\ord{np}}
(x11)++(0,\v) node (x21) {X\ord{mn} \sma X\ord{p}}
(x12)++(0,\v) node (x22) {X\ord{mnp}}
;
\draw[1cell=.85] 
(x11) edge node {1 \sma \mu_{\ord{n},\ord{p}}} (x12)
(x12) edge node {\mu_{\ord{m},\ord{np}}} (x22)
(x11) edge node[swap] {\mu_{\ord{m},\ord{n}} \sma 1} (x21)
(x21) edge node {\mu_{\ord{mn},\ord{p}}} (x22)
;
\end{tikzpicture}
\end{equation}
\item[Symmetry] 
The condition \cref{eastopsymmetry} means that the following pasting is equal to the identity natural transformation on $\mu_{\ord{m},\ord{n}}$.
\[\begin{tikzpicture}[xscale=1,yscale=1,vcenter]
\def\h{3} \def\v{-1.3} \def\u{.6} \def\a{25}
\draw[0cell=.8]
(0,0) node (x11) {X\ord{m} \sma X\ord{n}}
(x11)++(\h,0) node (x12) {X\ord{n} \sma X\ord{m}}
(x12)++(\h,0) node (x13) {X\ord{m} \sma X\ord{n}}
(x11)++(0,\v) node (x21) {X\ord{mn}}
(x12)++(0,\v) node (x22) {X\ord{nm}}
(x13)++(0,\v) node (x23) {X\ord{mn}}
(x11)++(0,\u) node[inner sep=0pt] (t1) {}
(x13)++(0,\u) node[inner sep=0pt] (t3) {}
(x21)++(0,-\u) node[inner sep=0pt] (b1) {}
(x23)++(0,-\u) node[inner sep=0pt] (b3) {}
;
\draw[1cell=.8] 
(x11) edge node (s) {\xisma} (x12)
(x12) edge node (s2) {\xisma} (x13)
(x11) edge node[swap] {\mu_{\ord{m},\ord{n}}} (x21)
(x12) edge node {\mu_{\ord{n},\ord{m}}} (x22)
(x13) edge node {\mu_{\ord{m},\ord{n}}} (x23)
(x11) edge[-,shorten >=-1pt] (t1)
(t1) edge[-,shorten <=-1pt,shorten >=-1pt] node {1} (t3)
(t3) edge[shorten <=-1pt] (x13)
(x23) edge[-,shorten >=-1pt] (b3)
(b3) edge[-,shorten <=-1pt,shorten >=-1pt] node[swap] {1} (b1)
(b1) edge[shorten <=-1pt] (x21)
; 
\draw[1cell=.7]
(x23) edge node (t2) {X(\xisma_{\ord{n},\ord{m}})^\inv} (x22)
(x22) edge node (t) {X(\xisma_{\ord{m},\ord{n}})^\inv} (x21)
; 
\draw[2cell=.9]
node[between=s and t at .55, 2label={above,\xi_{\ord{m},\ord{n}}}, 2label={below,\iso}] {\Longrightarrow}
node[between=s2 and t2 at .55, 2label={above,\xi_{\ord{n},\ord{m}}}, 2label={below,\iso}] {\Longrightarrow}
;
\end{tikzpicture}\]
The top and bottom regions are commutative by the symmetry axiom \cref{symmoncatsymhexagon} in, respectively, the symmetric monoidal categories
\begin{itemize}
\item $(\pCat,\sma,\xisma)$ in \cref{theorem:pC-monoidal} and
\item $(\Fskel,\sma,\xisma)$ in \cref{ex:Fskel}.
\end{itemize}
\item[Unity]
The condition \cref{eastopunity} means that the natural transformation
\[\begin{tikzpicture}[xscale=1,yscale=1,vcenter]
\def\a{25}
\draw[0cell=.9]
(0,0) node (x11) {X\ord{n}}
(x11)++(3.5,0) node (x12) {X\ord{n}}
;
\draw[1cell=.7] 
(x11) edge[bend left=\a] node {\mu_{\ord{n},\ord{1}}(-,\tu_{\ord{1}}(1)) = 1} (x12)
(x11) edge[bend right=\a] node[swap] {\mu_{\ord{1},\ord{n}}(\tu_{\ord{1}}(1),-) = 1} (x12)
;
\draw[2cell=.9]
node[between=x11 and x12 at .25, rotate=-90, 2label={above,(\xi_{\ord{n},\ord{1}})_{(-,\tu_{\ord{1}}(1))}}] {\Rightarrow}
;
\end{tikzpicture}\]
is equal to the identity on the identity functor $1_{X\ord{n}}$.
\item[Hexagon]
In the following diagram we use the notation
\[\xismaminus = (\xisma)^\inv \andspace \quad (-)' = X(-)\]
and abbreviate $\sma$ to concatenation.  For example, 
\[\ord{mn}' = X\ord{mn},\quad \ord{m}' \ord{n}' \ord{p}' = X\ord{m} \sma X\ord{n} \sma X\ord{p}, \andspace
(\xismaminus_{\ord{mn},\ord{p}})' = X(\xisma_{\ord{mn},\ord{p}})^\inv.\]
The condition \cref{eastophexagon} means the following pasting diagram equality.
\[\begin{tikzpicture}[xscale=1,yscale=1,vcenter]
\def\s{.7} \def\a{.6} \def\v{-1.4} \def\u{-1.4} \def\t{-1.4} \def\h{2.5}
\draw[0cell=\s]
(0,0) node (x11) {\ord{m}' \ord{n}' \ord{p}'}
(x11)++(1.7,0) node (x12) {\ord{p}' \ord{m}' \ord{n}'}
(x11)++(0,\v) node (x21) {\ord{mn}' \ord{p}'}
(x12)++(0,\v) node (x22) {\ord{p}' \ord{mn}'}
(x21)++(0,\v) node (x31) {\ord{mnp}'}
(x22)++(0,\v) node (x32) {\ord{pmn}'}
;
\draw[0cell=1.5]
(x22)++(1,0) node (eq) {=}
;
\draw[1cell=\s] 
(x11) edge node {\xisma} (x12)
(x21) edge node {\xisma} (x22)
(x32) edge node {(\xismaminus_{\ord{mn},\ord{p}})'} (x31)
(x11) edge node[swap] {\mu_{\ord{m},\ord{n}} 1} (x21)
(x21) edge node[swap] {\mu_{\ord{mn},\ord{p}}} (x31)
(x12) edge node {1 \mu_{\ord{m},\ord{n}}} (x22)
(x22) edge node {\mu_{\ord{p},\ord{mn}}} (x32)
;
\draw[2cell=.8]
node[between=x21 and x32 at .5, shift={(0,-.25)}, 2label={above,\xi_{\ord{mn},\ord{p}}}] {\Longrightarrow}
;
\draw[0cell=\s]
(eq)++(1,-\u) node (y11) {\ord{m}' \ord{n}' \ord{p}'}
(y11)++(\h,0) node (y12) {\ord{m}' \ord{p}' \ord{n}'}
(y12)++(\h,0) node (y13) {\ord{p}' \ord{m}' \ord{n}'}
(y11)++(0,\u) node (y21) {\ord{m}' \ord{np}'}
(y13)++(0,\u) node (y23) {\ord{pm}' \ord{n}'}
(y21)++(0,\u) node (y31) {\ord{mnp}'}
(y31)++(\h,0) node (y32) {\ord{mpn}'}
(y32)++(\h,0) node (y33) {\ord{pmn}'}
(y12)++(-\a,\t) node (z1) {\ord{m}' \ord{pn}'}
(y12)++(\a,\t) node (z2) {\ord{mp}' \ord{n}'}
;
\draw[1cell=\s] 
(y11) edge node {1 \xisma} (y12)
(y12) edge node {\xisma 1} (y13)
(y23) edge node {(\xismaminus_{\ord{m},\ord{p}})' 1} (z2)
(z1) edge node {1 (\xismaminus_{\ord{n},\ord{p}})'} (y21)
(y33) edge node {(\xismaminus_{\ord{m},\ord{p}} \sma 1)'} (y32)
(y32) edge node {(1 \sma \xismaminus_{\ord{n},\ord{p}})'} (y31)
(y11) edge node[swap] {1 \mu_{\ord{n},\ord{p}}} (y21)
(y21) edge node[swap] {\mu_{\ord{m},\ord{np}}} (y31)
(y13) edge node {\mu_{\ord{p},\ord{m}} 1} (y23)
(y23) edge node {\mu_{\ord{pm},\ord{n}}} (y33)
(y12) edge node[swap] {1 \mu_{\ord{p},\ord{n}}} (z1)
(z1) edge node[swap,pos=.4] {\mu_{\ord{m},\ord{pn}}} (y32)
(y12) edge node {\mu_{\ord{m},\ord{p}} 1} (z2)
(z2) edge node[pos=.4] {\mu_{\ord{mp},\ord{n}}} (y32)
;
\draw[2cell=.8]
node[between=y11 and z1 at .5, shift={(0,-.25)}, 2label={above, 1 \xi_{\ord{n},\ord{p}}}] {\Longrightarrow}
node[between=y13 and z2 at .5, shift={(0,-.25)}, 2label={above,\xi_{\ord{m},\ord{p}} 1}] {\Longrightarrow}
;
\end{tikzpicture}\]
The detail of this pasting diagram equality is as follows.
\begin{itemize}
\item The top square on the left-hand side is commutative by the naturality of the braiding $\xisma$ for the smash product \cref{eq:smash-pushout}, which is induced by the braiding for the Cartesian product.
\item Consider the right-hand side.
\begin{itemize}
\item The center diamond is commutative by the 1-cell associativity condition \cref{onecellassociativity}.
\item The two unlabeled trapezoids are commutative by the naturality of $\mu \cn X \sma X \to X$ \cref{ptdaynaturality}.
\end{itemize}
\item The left, respectively, right, vertical boundaries of the two sides are equal by the 1-cell associativity condition \cref{onecellassociativity}.
\item The top, respectively, bottom, horizontal boundaries of the two sides are equal by the hexagon axiom \cref{symmoncatsymhexagon} in, respectively, the symmetric monoidal categories
\begin{itemize}
\item $(\pCat,\sma,\xisma)$ in \cref{theorem:pC-monoidal} and
\item $(\Fskel,\sma,\xisma)$ in \cref{ex:Fskel}.
\end{itemize}
\end{itemize}
Equivalently, in terms of objects 
\[a \in X\ord{m}, \quad b \in X\ord{n}, \andspace c \in X\ord{p},\] 
the condition \cref{eastophexagon} means the morphism equality
\[\begin{split}
& \big( \xi_{\ord{mn},\ord{p}} \big)_{( \mu_{\ord{m},\ord{n}}(a,b) \scs c )} \\
&= X\left(1_{\ord{m}} \sma \xisma_{\ord{n},\ord{p}}\right)^\inv 
\left( \mu_{\ord{mp},\ord{n}} \big( (\xi_{\ord{m},\ord{p}})_{(a,c)} \scs 1_b \big) \right) 
\circ \mu_{\ord{m}, \ord{np}} \left( 1_a \scs (\xi_{\ord{n},\ord{p}})_{(b,c)} \right)
\end{split}\]
in $X\ord{mnp}$.
\end{description}
This finishes our description of an $E_\infty$-algebra in $\Gacat$.
\end{example}

\section{Inverse \texorpdfstring{$K$}{K}-Theory Preserves Pseudo Symmetric $E_\infty$-Algebras} 
\label{sec:apppseudo}

Inverse $K$-theory $\cP$ does \emph{not} preserve $E_\infty$-algebras in general, since $E_\infty$-algebras are parametrized by $\Cat$-multifunctors from the Barratt-Eccles operad $\BE$ (\cref{def:barratteccles}).  Composing such $\Cat$-multifunctors with the pseudo symmetric $\Cat$-multifunctor $\cP$ \cref{invKAgroa} yields pseudo symmetric $\Cat$-multifunctors (\cref{ex1:Ppreservespseudosym}), not $\Cat$-multifunctors in general.  In this section we
\begin{itemize}
\item define pseudo symmetric $E_\infty$-algebras (\cref{def:psEinftyalgebra}) and
\item observe that inverse $K$-theory $\cP$ preserves them (\cref{cor:Pefinitypsalg}).
\end{itemize}
Then we completely unpack the structure of a pseudo symmetric $E_\infty$-algebra in a general $\Cat$-multicategory, $\Gacat$, $\permcat$, and $\permcatsg$ (\cref{expl:psEinftyalgebra,expl:psEinftyGacat,expl:psEinftypermcat}).

\begin{definition}\label{def:psEinftyalgebra}\index{pseudo symmetric!Einfty-algebra@$E_\infty$-algebra}\index{Einfty-algebra@$E_\infty$-algebra!pseudo symmetric}
Suppose $\M$ is a $\Cat$-multicategory.  A pseudo symmetric $\Cat$-multifunctor (\cref{def:pseudosmultifunctor}) from the Barratt-Eccles operad (\cref{def:barratteccles})
\[\begin{tikzcd}[column sep=huge]
\BE \ar{r}{(F,\{F_\sigma\})} & \M
\end{tikzcd}\]
is called a \emph{pseudo symmetric $E_\infty$-algebra} in $\M$.
\end{definition}

Applied to the Barratt-Eccles operad $\BE$, \cref{cor:Ppreservespseudosym} has the following special case.

\begin{corollary}\label{cor:Pefinitypsalg}\index{inverse K-theory@inverse $K$-theory!preserves pseudo symmetric $E_\infty$-algebras}\index{pseudo symmetric!Einfty-algebra@$E_\infty$-algebra!preservation by inverse $K$-theory}\index{Einfty-algebra@$E_\infty$-algebra!pseudo symmetric!preservation by inverse $K$-theory}
The inverse $K$-theory pseudo symmetric $\Cat$-multifunctor in \cref{invKAgroa}
\[\begin{tikzcd}[column sep=large]
\Gacat \ar{r}{\cP} & \permcatsg
\end{tikzcd}\]
preserves pseudo symmetric $E_\infty$-algebras.
\end{corollary}

\begin{example}[$E_\infty$-Algebras]\label{ex:Peinftypseudo}
By \cref{ex:idpsm} each $E_\infty$-algebra (\cref{def:barratteccles}) is a pseudo symmetric $E_\infty$-algebra when it is equipped with identity pseudo symmetry isomorphisms.  By \cref{cor:Pefinitypsalg} inverse $K$-theory $\cP$ sends each $E_\infty$-algebra in $\Gacat$ (\cref{ex:Einftygacat}) to a pseudo symmetric $E_\infty$-algebra in $\permcatsg$.  
\end{example}

Next we explain in detail the structure of a pseudo symmetric $E_\infty$-algebra in
\begin{itemize}
\item a general $\Cat$-multicategory (\cref{expl:psEinftyalgebra}),
\item $\Gacat$ (\cref{expl:psEinftyGacat}), and
\item $\permcat$ and $\permcatsg$ (\cref{expl:psEinftypermcat}).
\end{itemize}
\cref{expl:psEinftyalgebra} below should be compared with the Coherence \cref{mapfrombarratteccles} for $E_\infty$-algebras.

\begin{explanation}[Pseudo Symmetric $E_\infty$-Algebras]\label{expl:psEinftyalgebra}\index{pseudo symmetric!Einfty-algebra@$E_\infty$-algebra}\index{Einfty-algebra@$E_\infty$-algebra!pseudo symmetric}
Interpreting \cref{def:pseudosmultifunctor} with domain the Barratt-Eccles operad $\BE$ (\cref{def:barratteccles}), a pseudo symmetric $E_\infty$-algebra in a $\Cat$-multicategory $(\M,\ga,1)$
\[\begin{tikzcd}[column sep=huge]
\BE \ar{r}{(F,\{F_\sigma\})} & \M
\end{tikzcd}\]
consists of the data \cref{xFstar,musigmaonecell,mutausigma,Fsigmatau} below for $n \geq 0$ and $\sigma,\tau \in \Sigma_n$.
\begin{description}
\item[Object] The object assignment of $F$ yields an object 
\begin{equation}\label{xFstar}
x = F(*) \in \M.
\end{equation}
\item[1-Cells in Arity $n$] The $n$-ary 1-cell assignment of $F$ yields an object
\begin{equation}\label{musigmaonecell}
\mu_\sigma = F(\sigma) \in \M\scmap{x^n;x}
\end{equation}
with $x^n = \angx_{j=1}^n$ the $n$-tuple consisting of $n$ copies of $x$.
\item[2-Cells in Arity $n$] The $n$-ary 2-cell assignment of $F$ yields an isomorphism
\begin{equation}\label{mutausigma}
\begin{tikzcd}[column sep=large]
\mu_\sigma \ar{r}{\mu_\sigma^\tau}[swap]{\iso} & \mu_\tau
\end{tikzcd} \inspace \M\scmap{x^n;x}
\end{equation}
corresponding to the unique isomorphism in $E\Sigma_n(\sigma,\tau) = \{*\}$.
\item[Pseudo Symmetry Isomorphisms] The $\tau$-component of the pseudo symmetry isomorphism $F_\sigma$ \cref{Fsigmaangccp} is an $n$-ary 2-cell isomorphism
\begin{equation}\label{Fsigmatau}
\begin{tikzcd}[column sep=large]
\mu_{\tau\sigma} \ar{r}{F_{\sigma; \tau}}[swap]{\iso} & \mu_\tau \cdot \sigma
\end{tikzcd} \inspace \M\scmap{x^n;x}
\end{equation}
with $\mu_\tau \cdot \sigma$ the right $\sigma$-action on $\mu_\tau$.
\end{description}
The above data are subject to the conditions \cref{Fcomponentfunctor} through \cref{Fbotequivariance} below for $n \geq 0$ and $\sigma,\tau,\pi \in \Sigma_n$.
\begin{description}
\item[Component Functoriality] The functoriality of each functor 
\[\begin{tikzcd}[column sep=large]
E\Sigma_n \ar{r}{F} & \M\scmap{x^n;x}
\end{tikzcd}\]
means the following morphism equalities in $\M\scmap{x^n;x}$.
\begin{equation}\label{Fcomponentfunctor}
\left\{\begin{split}
\mu_\sigma^\sigma &= 1_{\mu_\sigma}\\
\mu_\tau^\pi \circ \mu_\sigma^\tau &= \mu_\sigma^\pi
\end{split}\right.
\end{equation}
\item[Pseudo Symmetry Naturality] The naturality of $F_\sigma$ in \cref{Fsigmatau} means that the following diagram in $\M\smscmap{x^n;x}$ is commutative.
\begin{equation}\label{Fsigmanaturality}
\begin{tikzpicture}[xscale=1,yscale=1,vcenter]
\def\v{-1.2}
\draw[0cell=.9]
(0,0) node (x11) {\mu_{\tau\sigma}}
(x11)++(2.5,0) node (x12) {\mu_\tau \cdot \sigma}
(x11)++(0,\v) node (x21) {\mu_{\pi\sigma}}
(x12)++(0,\v) node (x22) {\mu_\pi \cdot \sigma}
;
\draw[1cell=.9] 
(x11) edge node {F_{\sigma; \tau}} (x12)
(x21) edge node {F_{\sigma; \pi}} (x22)
(x11) edge node[swap] {\mu^{\pi\sigma}_{\tau\sigma}} (x21)
(x12) edge node {\mu^{\pi}_{\tau} \cdot \sigma} (x22)
; 
\end{tikzpicture}
\end{equation}
\item[Unit] The unit axiom \cref{enr-multifunctor-unit} is the object equality
\begin{equation}\label{muidone}
\mu_{\id_1} = 1_x \in \M\smscmap{x;x}.
\end{equation}
\item[Composition] Using the notation in \cref{Asgamma,sigmatilde,sigmatimes}, the composition axiom \cref{v-multifunctor-composition} on objects is the following object equality in $\M\smscmap{x^{k_1+\cdots+k_n};x}$.
\begin{equation}\label{compaxiomobject}
\mu_{\sigmatil \cdot \sigmatimes} = \ga\scmap{\mu_\sigma; \ang{\mu_{\sigma_j}}_{j=1}^n}
\end{equation}
The composition axiom \cref{v-multifunctor-composition} on morphisms is the following morphism equality for $\sigma, \tau \in \Sigma_n$ and $\sigma_j, \tau_j \in \Sigma_{k_j}$.
\begin{equation}\label{compaxiommorphism}
\mu_{\sigmatil \cdot \sigmatimes}^{\tautil \cdot \tautimes} 
= \ga\scmap{\mu_\sigma^\tau; \ang{\mu_{\sigma_j}^{\tau_j}}_{j=1}^n}
\end{equation}
\end{description}
The axioms \cref{pseudosmf-unity,pseudosmf-product,pseudosmf-topeq,pseudosmf-boteq} in \cref{def:pseudosmultifunctor} are as follows.
\begin{description}
\item[Unit Permutation] The axiom \cref{pseudosmf-unity} is the following morphism equality in $\M\smscmap{x^n;x}$.
\begin{equation}\label{Fidntau}
F_{\id_n; \tau} = 1_{\mu_\tau} \cn \mu_\tau \to \mu_\tau
\end{equation}
\item[Product Permutation] The axiom \cref{pseudosmf-product} means that the following diagram in $\M\smscmap{x^n;x}$ is commutative, with the equality given by the symmetric group action axiom \cref{enr-multicategory-symmetry} in $\M$.
\begin{equation}\label{Fproductpermutation}
\begin{tikzpicture}[xscale=1,yscale=1,vcenter]
\def\v{-1.2}
\draw[0cell=.9]
(0,0) node (x11) {\mu_{\pi\sigma\tau}}
(x11)++(3,0) node (x12) {\mu_\pi \cdot \sigma\tau}
(x11)++(0,\v) node (x21) {\mu_{\pi\sigma} \cdot \tau}
(x12)++(0,\v) node (x22) {(\mu_\pi \cdot \sigma) \cdot \tau}
;
\draw[1cell=.9] 
(x11) edge node {F_{\sigma\tau; \pi}} (x12)
(x21) edge node {F_{\sigma; \pi} \cdot \tau} (x22)
(x11) edge node[swap] {F_{\tau; \pi\sigma}} (x21)
(x12) edge[equal] (x22)
; 
\end{tikzpicture}
\end{equation}
\item[Top Equivariance] The axiom \cref{pseudosmf-topeq} means that, for $\sigma,\tau \in \Sigma_n$ and $\tau_j \in \Sigma_{k_j}$ with $j \in \{1,\ldots,n\}$, the following diagram in $\M\smscmap{x^{k_1+\cdots+k_n};x}$ is commutative.
\begin{equation}\label{Ftopequivariance}
\begin{tikzpicture}[xscale=1,yscale=1,vcenter]
\def\v{-1}
\draw[0cell=.9]
(0,0) node (x11) {\mu_{\tautil \cdot \tautimes \cdot \sigmabar}}
(x11)++(6,0) node (x12) {\mu_{\tautil \cdot \tautimes} \cdot \sigmabar}
(x11)++(0,\v) node (x21) {\mu_{\til{\tau\sigma} \cdot (\smallprod_j \tau_{\sigma(j)})}}
(x12)++(0,\v) node (x22) {\ga\scmap{\mu_\tau; \ang{\mu_{\tau_j}}_{j=1}^n} \cdot \sigmabar}
(x21)++(0,\v) node (x31) {\ga\scmap{\mu_{\tau\sigma}; \ang{\mu_{\tau_{\sigma(j)}}}_{j=1}^n}}
(x22)++(0,\v) node (x32) {\ga\scmap{\mu_\tau \cdot \sigma; \ang{\mu_{\tau_{\sigma(j)}}}_{j=1}^n}}
;
\draw[1cell=.8] 
(x11) edge node {F_{\sigmabar;\, \tautil \cdot \tautimes}} (x12)
(x31) edge node {\ga\scmap{F_{\sigma; \tau}; \ang{1_{\mu_{\tau_{\sigma(j)}}}}_{j=1}^n}} (x32)
(x11) edge[equal] node[swap] {\spadesuit} (x21)
(x22) edge[equal] node {\spadesuit} (x32)
(x21) edge[equal] node[swap] {\clubsuit} (x31)
(x12) edge[equal] node {\clubsuit} (x22)
; 
\end{tikzpicture}
\end{equation}
\begin{itemize}
\item The block permutation
\[\sigmabar = \sigma\ang{k_{\sigma(1)}, \ldots, k_{\sigma(n)}} \in \Sigma_{k_1 + \cdots + k_n}\]
is induced by $\sigma$ as in \cref{blockpermutation}.
\item The two equalities labeled $\spadesuit$ use the top equivariance axiom \cref{enr-operadic-eq-1} in, respectively, $\BE$ and $\M$. 
\item The two equalities labeled $\clubsuit$ use the composition axiom \cref{compaxiomobject}.
\end{itemize}
\item[Bottom Equivariance] The axiom \cref{pseudosmf-boteq} means that the following diagram in $\M\smscmap{x^{k_1+\cdots+k_n};x}$ is commutative for $\sigma \in \Sigma_n$ and $\sigma_j, \tau_j \in \Sigma_{k_j}$ with $j \in \{1,\ldots,n\}$.
\begin{equation}\label{Fbotequivariance}
\begin{tikzpicture}[xscale=1,yscale=1,vcenter]
\def\v{-1}
\draw[0cell=.9]
(0,0) node (x11) {\mu_{\sigmatil \cdot \sigmatimes \cdot \tautimes}}
(x11)++(5.5,0) node (x12) {\mu_{\sigmatil \cdot \sigmatimes} \cdot \tautimes}
(x11)++(0,\v) node (x21) {\mu_{\sigmatil \cdot (\smallprod_j \sigma_j \tau_j)}}
(x12)++(0,\v) node (x22) {\ga\scmap{\mu_\sigma; \ang{\mu_{\sigma_j}}_{j=1}^n} \cdot \tautimes}
(x21)++(0,\v) node (x31) {\ga\scmap{\mu_{\sigma}; \ang{\mu_{\sigma_j \tau_j}}_{j=1}^n}}
(x22)++(0,\v) node (x32) {\ga\scmap{\mu_\sigma; \ang{\mu_{\sigma_j} \cdot \tau_j}_{j=1}^n}}
;
\draw[1cell=.8] 
(x11) edge node {F_{\tautimes;\, \sigmatil \cdot \sigmatimes}} (x12)
(x31) edge node {\ga\scmap{1_{\mu_\sigma}; \ang{F_{\tau_j;\, \sigma_j}}_{j=1}^n}} (x32)
(x11) edge[equal] node[swap] {\heartsuit} (x21)
(x22) edge[equal] node {\heartsuit} (x32)
(x21) edge[equal] node[swap] {\lozenge} (x31)
(x12) edge[equal] node {\lozenge} (x22)
; 
\end{tikzpicture}
\end{equation}
\begin{itemize}
\item The two equalities labeled $\heartsuit$ use the bottom equivariance axiom \cref{enr-operadic-eq-2} in, respectively, $\BE$ and $\M$. 
\item The two equalities labeled $\lozenge$ use the composition axiom \cref{compaxiomobject}.
\end{itemize}
\end{description}
This finishes the explicit description of a pseudo symmetric $E_\infty$-algebra in a $\Cat$-multicategory $(\M,\ga,1)$.
\end{explanation}

\cref{expl:psEinftyGacat} below should be compared with the description of an $E_\infty$-algebra in $\Gacat$ (\cref{ex:Einftygacat}).

\begin{explanation}[Pseudo Symmetric $E_\infty$-Algebras in $\Gacat$]\label{expl:psEinftyGacat}\index{pseudo symmetric!Einfty-algebra@$E_\infty$-algebra!in Gacat@in $\Gacat$}\index{Einfty-algebra@$E_\infty$-algebra!pseudo symmetric!in Gacat@in $\Gacat$}
Applying \cref{expl:psEinftyalgebra} to the $\Cat$-multicategory $\M = \Gacat$ (\cref{sec:Gacatmulticat}), a pseudo symmetric $E_\infty$-algebra in $\Gacat$\label{not:psEinfGacat}
\[\big(X, \{\mu_\sigma\}, \{\mu_\sigma^\tau\}, \{F_\sigma\} \big) \cn \BE \to \Gacat\]
consists of the data \cref{underlyinggacat,musigmaxnx,muidzeroone,musigmatauxnx,FsigmatauGacat} below for $n \geq 0$ and $\sigma,\tau \in \Sigma_n$.
\begin{description}
\item[Underlying $\Ga$-Category] It consists of a $\Ga$-category (\cref{def:gammacategory})
\begin{equation}\label{underlyinggacat}
\begin{tikzcd}[column sep=large]
(\Fskel,\ord{0}) \ar{r}{X} & (\Cat,\boldone).
\end{tikzcd}
\end{equation}
\item[1-Cells in Arity $n$] It has a pointed natural transformation \cref{Fxonexntoz}
\begin{equation}\label{musigmaxnx}
\begin{tikzcd}[column sep=large]
\txsma_{j=1}^n X \ar{r}{\mu_\sigma} & X
\end{tikzcd} 
\inspace \GMap\brb{\txsma_{j=1}^n X,X}.
\end{equation}
For each $n$-tuple of objects $\angm = \bang{\ord{m_j}}_{j=1}^n$ in $\Fskel$, $\mu_\sigma$ consists of a component pointed functor \cref{Fordmcomponents}
\begin{equation}\label{musigmaangmmxn}
\begin{tikzcd}[column sep=huge]
\txsma_{j=1}^n X\ord{m_j} \ar{r}{\mu_{\sigma; \angm}} & X\big(\ord{m_1 \cdots\, m_n}\big)
\end{tikzcd}
\end{equation}
that satisfies the naturality diagram \cref{Fordmnaturality} with respect to $\bang{\ord{m_j}}_{j=1}^n$.  

If $n=0$, then $\txsma_{j=1}^n X$ is interpreted as the monoidal unit diagram $\ftu$ in \cref{Gacatunit}.  As discussed in the paragraph containing \cref{Fordnjziotaj}, the pointed natural transformation
\[\begin{tikzcd}[column sep=large]
\ftu \ar{r}{\mu_{\id_0}} & X
\end{tikzcd}
\inspace \GMap\brb{\ftu,X}\]
is determined by the object 
\begin{equation}\label{muidzeroone}
(\mu_{\id_0})_{\ord{1}}(1) \in X\ord{1}.
\end{equation}
\item[2-Cells in Arity $n$] It has an invertible pointed modification
\begin{equation}\label{musigmatauxnx}
\begin{tikzpicture}[xscale=1,yscale=1,baseline={(x11.base)}]
\def\d{25} 
\draw[0cell=.9]
(0,0) node (x11) {\txsma_{j=1}^n X}
(x11)++(.35,.07) node (a) {\phantom{Z}}
(a)++(2.2,0) node (b) {X}
;
\draw[1cell=.8] 
(a) edge[bend left=\d] node {\mu_\sigma} (b)
(a) edge[bend right=\d] node[swap] {\mu_\tau} (b)
;
\draw[2cell=.9]
node[between=a and b at .4, rotate=-90, 2label={above,\mu_\sigma^\tau}] {\Rightarrow}
;
\end{tikzpicture}
\inspace \GMap\brb{\txsma_{j=1}^n X,X}.
\end{equation}
For each $n$-tuple of objects $\angm = \bang{\ord{m_j}}_{j=1}^n$ in $\Fskel$, $\mu_\sigma^\tau$ consists of a component pointed natural isomorphism \cref{thetaordmcomponent}
\begin{equation}\label{musigmatauangm}
\begin{tikzpicture}[xscale=1,yscale=1,baseline={(x11.base)}]
\def\d{27} 
\draw[0cell=.9]
(0,0) node (x11) {\txsma_{j=1}^n X\ord{m_j}}
(x11)++(.55,.06) node (a) {\phantom{Z}}
(a)++(2.5,0) node (b) {\phantom{Z}}
(b)++(.7,0) node (x12) {X\big(\ord{m_1 \cdots\, m_n}\big)}
;
\draw[1cell=.8] 
(a) edge[bend left=\d] node {\mu_{\sigma;\angm}} (b)
(a) edge[bend right=\d] node[swap] {\mu_{\tau;\angm}} (b)
;
\draw[2cell=.9]
node[between=a and b at .3, rotate=-90, 2label={above,(\mu_\sigma^\tau)_{\angm}}] {\Rightarrow}
;
\end{tikzpicture}
\end{equation}
that satisfies the modification axiom \cref{thetaordmodification}.

If $n=0$, then by \cref{Fcomponentfunctor} $\mu_{\id_0}^{\id_0}$ consists of the identity morphism on the object determined by $\mu_{\id_0}$ in \cref{muidzeroone}.
\item[Pseudo Symmetry Isomorphisms] It has an invertible pointed modification
\begin{equation}\label{FsigmatauGacat}
\begin{tikzpicture}[xscale=1,yscale=1,baseline={(x11.base)}]
\def\d{25} 
\draw[0cell=.9]
(0,0) node (x11) {\txsma_{j=1}^n X}
(x11)++(.35,.07) node (a) {\phantom{Z}}
(a)++(2.2,0) node (b) {X}
;
\draw[1cell=.8] 
(a) edge[bend left=\d] node {\mu_{\tau\sigma}} (b)
(a) edge[bend right=\d] node[swap] {\mu_\tau \cdot \sigma} (b)
;
\draw[2cell=.9]
node[between=a and b at .37, rotate=-90, 2label={above,F_{\sigma; \tau}}] {\Rightarrow}
;
\end{tikzpicture}
\inspace \GMap\brb{\txsma_{j=1}^n X,X}.
\end{equation}
For each $n$-tuple of objects $\angm = \bang{\ord{m_j}}_{j=1}^n$ in $\Fskel$, $F_{\sigma; \tau}$ consists of a component pointed natural isomorphism \cref{thetaordmcomponent}
\begin{equation}\label{FsigmatauGacatcomp}
\begin{tikzpicture}[xscale=1,yscale=1,baseline={(x11.base)}]
\def\d{27} 
\draw[0cell=.9]
(0,0) node (x11) {\txsma_{j=1}^n X\ord{m_j}}
(x11)++(.55,.06) node (a) {\phantom{Z}}
(a)++(2.5,0) node (b) {\phantom{Z}}
(b)++(.7,0) node (x12) {X\big(\ord{m_1 \cdots\, m_n}\big)}
;
\draw[1cell=.8] 
(a) edge[bend left=\d] node {\mu_{\tau\sigma;\angm}} (b)
(a) edge[bend right=\d] node[swap] {(\mu_\tau \cdot \sigma)_{\angm}} (b)
;
\draw[2cell=.9]
node[between=a and b at .3, rotate=-90, 2label={above,F_{\sigma; \tau,\angm}}] {\Rightarrow}
;
\end{tikzpicture}
\end{equation}
that satisfies the modification axiom \cref{thetaordmodification}.  By \cref{Fsigmacomponent} the pointed functor $(\mu_\tau \cdot \sigma)_{\angm}$ in \cref{FsigmatauGacatcomp} is the following composite.
\begin{equation}\label{mutausigmaangmcomposite}
\begin{tikzpicture}[xscale=1,yscale=1.4,vcenter]
\draw[0cell=.9]
(0,0) node (x11) {\txsma_{j=1}^n X \ord{m_j}}
(x11)++(4.5,0) node (x12) {X\big(\ord{m_1 \cdots\, m_n} \big)}
(x11)++(0,-1) node (x21) {\txsma_{j=1}^n X \ord{m_{\sigmainv(j)}}}
(x12)++(0,-1) node (x22) {X\big(\ord{m_{\sigmainv(1)} \cdots\, m_{\sigmainv(n)}} \big)}
;
\draw[1cell=.9]  
(x11) edge node {(\mu_\tau \cdot \sigma)_{\angm}} (x12)
(x11) edge node {\iso} node[swap] {\sigma} (x21)
(x21) edge node {\mu_{\tau; \sigma\angm}} (x22)
(x22) edge[shorten <=-.4ex] node {\iso} node [swap] {X(\sigmainv)} (x12)
;
\end{tikzpicture}
\end{equation}

If $n=0$, then by \cref{Fidntau} $F_{\id_0;\id_0}$ is the identity morphism on the object determined by $\mu_{\id_0}$ in \cref{muidzeroone}.
\end{description}
This finishes the description of the list of data of a pseudo symmetric $E_\infty$-algebra in $\Gacat$.

The conditions \cref{Fcomponentfunctor} through \cref{Fbotequivariance} for $n \geq 0$ and $\sigma,\tau,\pi \in \Sigma_n$ are interpreted componentwise using the data \cref{musigmaangmmxn,musigmatauangm,FsigmatauGacatcomp}, with appropriate adjustment using \cref{muidzeroone} for 0-ary inputs.  Below we state explicitly the three most nontrivial conditions from that list.
\begin{description}
\item[Composition] Consider objects
\[\ord{m_{ji}} \in \Fskel \qquad \text{for $1\leq j \leq n$ and $1 \leq i \leq k_j$}\]
along with the following notation.
\begin{equation}\label{angmji}
\left\{\begin{aligned}
\bang{\ord{m_j}} &= \bang{\ord{m_{ji}}}_{i=1}^{k_j} \in \Fskel^{k_j} 
& m_{j\bdot} &= \txprod_{i=1}^{k_j} m_{ji}\\
\angm &= \bang{\bang{\ord{m_j}}}_{j=1}^n \in \Fskel^{k_1+\cdots+k_n} &&
\end{aligned}\right.
\end{equation}
Then the composition axiom on objects \cref{compaxiomobject} is the following commutative diagram of functors.
\begin{equation}\label{musigmasigmaangm}
\begin{tikzpicture}[xscale=1,yscale=1,vcenter]
\draw[0cell=.9]
(0,0) node (x11) {\txsma_{j=1}^n \txsma_{i=1}^{k_j} X\ord{m_{ji}}}
(x11)++(4.5,.0) node (x12) {X\big(\ord{\txprod_{j=1}^n m_{j\bdot}}\big)}
(x11)++(2.3,-1.2) node (x2) {\txsma_{j=1}^n X \ord{m_{j\bdot}}}
;
\draw[1cell=.9] 
(x11) edge node {\mu_{\sigmatil \cdot \sigmatimes; \angm}} (x12)
(x11) edge[transform canvas={xshift=-.5em}] node[swap] {\wedge_j\, \mu_{\sigma_j; \ang{\ord{m_j}}}} (x2)
(x2) edge[transform canvas={xshift=.5em}] node[swap] {\mu_{\sigma; \ang{\ord{m_{j\bdot}}}_{j=1}^n}} (x12)
;
\end{tikzpicture}
\end{equation}
The composition axiom on morphisms \cref{compaxiommorphism} is the analog of \cref{musigmasigmaangm} with each arrow replaced by its 2-cell variant in \cref{musigmatauangm}.
\item[Top Equivariance] We use the notation in \cref{angmji} together with
\[\angm\sigmabar = \bang{\bang{\ord{m_{\sigma(j)}}}}_{j=1}^n \in \Fskel^{k_{\sigma(1)} + \cdots + k_{\sigma(n)}}.\]
The top equivariance axiom \cref{Ftopequivariance} means that the component natural isomorphism
\[\begin{tikzpicture}[xscale=1,yscale=1,baseline={(x11.base)}]
\def\d{20} 
\draw[0cell=.9]
(0,0) node (x11) {\txsma_{j=1}^n \txsma_{i=1}^{k_{\sigma(j)}} X\ord{m_{\sigma(j),i}}}
(x11)++(1.2,0) node (a) {\phantom{Z}}
(a)++(3.5,0) node (b) {\phantom{Z}}
(b)++(.9,-.04) node (x12) {X\big(\ord{\txprod_{j=1}^n m_{\sigma(j),\bdot}}\big)}
;
\draw[1cell=.8] 
(a) edge[bend left=\d] node {\mu_{\tautil \cdot \tautimes \cdot \sigmabar; \angm\sigmabar}} (b)
(a) edge[bend right=\d] node[swap] {(\mu_{\tautil \cdot \tautimes} \cdot \sigmabar)_{\angm\sigmabar}} (b)
;
\draw[2cell=.9]
node[between=a and b at .3, rotate=-90, 2label={above,F_{\sigmabar;\, \tautil \cdot \tautimes,\angm\sigmabar}}] {\Rightarrow}
;
\end{tikzpicture}\]
is equal to the following whiskered natural isomorphism.
\[\begin{tikzpicture}[xscale=1,yscale=1,vcenter]
\def\d{20} \def\v{-1.4}
\draw[0cell=.8]
(0,0) node (x11) {\txsma_{j=1}^n \txsma_{i=1}^{k_{\sigma(j)}} X\ord{m_{\sigma(j),i}}}
(x11)++(2.5,-\v/2) node (x12) {\txsma_{j=1}^n X\ord{m_{\sigma(j),\bdot}}}
(x12)++(3.8,0) node (x13) {X\big(\ord{\txprod_{j=1}^n m_{\sigma(j),\bdot}}\big)}
(x12)++(0,\v) node (x22) {\txsma_{j=1}^n X\ord{m_{j \bdot}}}
(x13)++(0,\v) node (x23) {X\big(\ord{\txprod_{j=1}^n m_{j \bdot}}\big)}
;
\draw[1cell=.8] 
(x11) edge[shorten <=-5ex, transform canvas={yshift=.5em}] node[pos=.3] {\wedge_j \, \mu_{\tau_{\sigma(j)}; \ang{\ord{m_{\sigma(j)}}}}} (x12)
(x12) edge node {\mu_{\tau\sigma; \ang{\ord{m_{\sigma(j),\bdot}}}_{j=1}^n}} (x13)
(x12) edge node[swap] {\sigma} (x22)
(x22) edge node[swap] {\mu_{\tau; \ang{\ord{m_{j \bdot}}}_{j=1}^n}} (x23)
(x23) edge node[swap] {X(\sigmainv)} (x13)
;
\draw[2cell=.9]
node[between=x12 and x23 at .5, shift={(-.8,0)}, rotate=-90, 2label={above,F_{\sigma; \tau, \ang{\ord{m_{\sigma(j),\bdot}}}_{j=1}^n}}] {\Longrightarrow}
;
\end{tikzpicture}\]
In the previous two diagrams, the top, respectively, bottom, boundary composites are equal by the left, respectively, right, vertical equalities in \cref{Ftopequivariance}, along with \cref{mutausigmaangmcomposite}.
\item[Bottom Equivariance] We use \cref{angmji} along with the following notation.
\[\left\{\begin{aligned}
\bang{\ord{m_j}} \tau_j &= \bang{\ord{m_{j, \tau_j(i)}}}_{i=1}^{k_j} \in \Fskel^{k_j} 
& m_{j, \tau_j \bdot} &= \txprod_{i=1}^{k_j} m_{j, \tau_j(i)}\\
\angm\tautimes &= \bang{\bang{\ord{m_j}} \tau_j}_{j=1}^n \in \Fskel^{k_1+\cdots+k_n} &&
\end{aligned}\right.\]
The bottom equivariance axiom \cref{Fbotequivariance} means that the component natural isomorphism
\[\begin{tikzpicture}[xscale=1,yscale=1,baseline={(x11.base)}]
\def\d{20} 
\draw[0cell=.9]
(0,0) node (x11) {\txsma_{j=1}^n \txsma_{i=1}^{k_j} X\ord{m_{j, \tau_j(i)}}}
(x11)++(1.15,0) node (a) {\phantom{Z}}
(a)++(3.5,0) node (b) {\phantom{Z}}
(b)++(.85,-.04) node (x12) {X\big(\ord{\txprod_{j=1}^n m_{j, \tau_j \bdot}}\big)}
;
\draw[1cell=.8] 
(a) edge[bend left=\d] node {\mu_{\sigmatil \cdot \sigmatimes \cdot \tautimes; \angm\tautimes}} (b)
(a) edge[bend right=\d] node[swap] {(\mu_{\sigmatil \cdot \sigmatimes} \cdot \tautimes)_{\angm\tautimes}} (b)
;
\draw[2cell=.9]
node[between=a and b at .3, rotate=-90, 2label={above,F_{\tautimes;\, \sigmatil \cdot \sigmatimes, \angm\tautimes}}] {\Rightarrow}
;
\end{tikzpicture}\]
is equal to the following whiskered natural isomorphism.
\[\begin{tikzpicture}[xscale=1,yscale=1,vcenter]
\def\v{-1.4}
\draw[0cell=.8]
(0,0) node (x11) {\txsma_{j=1}^n \txsma_{i=1}^{k_j} X\ord{m_{j, \tau_j(i)}}}
(x11)++(3.8,0) node (x12) {\txsma_{j=1}^n X\ord{m_{j, \tau_j \bdot}}}
(x12)++(2.5,\v/2) node (x13) {X\big(\ord{\txprod_{j=1}^n m_{j, \tau_j \bdot}}\big)}
(x11)++(0,\v) node (x21) {\txsma_{j=1}^n \txsma_{i=1}^{k_j} X\ord{m_{ji}}}
(x12)++(0,\v) node (x22) {\txsma_{j=1}^n X\ord{m_{j \bdot}}}
;
\draw[1cell=.8] 
(x11) edge node {\wedge_j\, \mu_{\sigma_j \tau_j; \ang{\ord{m_j}} \tau_j}} (x12)
(x12) edge[shorten >=-2ex, transform canvas={yshift=.5em}] node[pos=.3] {\mu_{\sigma; \ang{\ord{m_{j, \tau_j \bdot}}}_{j=1}^n}} (x13)
(x11) edge node[swap] {\wedge_j\, \tau_j} (x21)
(x21) edge node[swap] {\wedge_j\, \mu_{\sigma_j; \ang{\ord{m_j}}}} (x22)
(x22) edge[transform canvas={xshift=-.5em}] node[swap] {\wedge_j\, X(\tau_j^\inv)} (x12)
;
\draw[2cell=.9]
node[between=x11 and x22 at .5, shift={(-.8,-.1)}, rotate=-90, 2label={above,\txsma_j F_{\tau_j;\, \sigma_j, \ang{\ord{m_j}} \tau_j}}] {\Longrightarrow}
;
\end{tikzpicture}\]
In the previous two diagrams, the top, respectively, bottom, boundary composites are equal by the left, respectively, right, vertical equalities in \cref{Fbotequivariance}, along with \cref{mutausigmaangmcomposite}.
\end{description}
This finishes the description of a pseudo symmetric $E_\infty$-algebra in $\Gacat$.
\end{explanation}

Recall from \cref{ex:Einftyalgebras} that $E_\infty$-algebras in $\permcat$ and $\permcatsg$ are, respectively, small bipermutative categories and their tight variants (\cref{def:embipermutativecat}).  \cref{expl:psEinftypermcat} below also applies to $\permcatsg$ \cref{permcatsgangcd} by restricting to \emph{strong} $n$-linear functors, which have invertible linearity constraints \cref{musigmatwoj}.

\begin{explanation}[Pseudo Symmetric $E_\infty$-Algebras in $\permcat$]\label{expl:psEinftypermcat}\index{pseudo symmetric!Einfty-algebra@$E_\infty$-algebra!in permcat@in $\permcat$}\index{Einfty-algebra@$E_\infty$-algebra!pseudo symmetric!in permcat@in $\permcat$}
Applying \cref{expl:psEinftyalgebra} to the $\Cat$-multicategory $\M = \permcat$ (\cref{thm:permcatenrmulticat}), a pseudo symmetric $E_\infty$-algebra in $\permcat$\label{not:psEinfpermcat}
\[\big(\C, \{\mu_\sigma\}, \{\mu_\sigma^\tau\}, \{F_\sigma\} \big) \cn \BE \to \permcat\]
consists of the data \cref{Coplusebeta,musigmanlinear,musigmatwoj,musigmataunlinear,Fsigmanlinear} below for $n \geq 0$ and $\sigma,\tau \in \Sigma_n$.
\begin{description}
\item[Underlying Permutative Category] It consists of a small permutative category (\cref{def:symmoncat})
\begin{equation}\label{Coplusebeta}
(\C,\oplus,e,\beta).
\end{equation}
\item[1-Cells in Arity $n$] It has an $n$-linear functor (\cref{def:nlinearfunctor})
\begin{equation}\label{musigmanlinear}
\begin{tikzcd}[column sep=2.5cm]
\txprod_{j=1}^n \C = \C^n \ar{r}{\left(\mu_\sigma \scs \ang{\mu^2_{\sigma; j}}_{j=1}^n\right)} & \C.
\end{tikzcd}
\end{equation}
The natural transformation
\begin{equation}\label{musigmatwoj}
\begin{tikzcd}[column sep=large]
\mu_\sigma \ang{X \compj X_j} \oplus \mu_\sigma \ang{X \compj X_j'} \ar{r}{\mu^2_{\sigma;j}} & 
\mu_\sigma \ang{X \compj (X_j \oplus X_j')}
\end{tikzcd}
\end{equation}
for objects $\angX \in \C^n$ and $X_j' \in \C$ is the $j$-th linearity constraint \cref{laxlinearityconstraints} of $\mu_\sigma$.

If $\M = \permcatsg$, then $\mu_\sigma$ is a \emph{strong} $n$-linear functor.  This means that each linearity constraint $\mu^2_{\sigma;j}$ is a natural isomorphism.
\item[2-Cells in Arity $n$] It has an $n$-linear natural isomorphism (\cref{def:nlineartransformation})
\begin{equation}\label{musigmataunlinear}
\begin{tikzpicture}[xscale=1,yscale=1,baseline={(x11.base)}]
\def\d{25} 
\draw[0cell=.9]
(0,0) node (x11) {\C^n}
(x11)++(2.5,0) node (b) {\C}
;
\draw[1cell=.8] 
(x11) edge[bend left=\d] node {\left(\mu_\sigma, \ang{\mu^2_{\sigma;j}}_{j=1}^n\right)} (b)
(x11) edge[bend right=\d] node[swap] {\left(\mu_\tau, \ang{\mu^2_{\tau;j}}_{j=1}^n\right)} (b)
;
\draw[2cell=.9]
node[between=x11 and b at .4, rotate=-90, 2label={above,\mu_\sigma^\tau}] {\Rightarrow}
;
\end{tikzpicture}
\inspace \permcat\scmap{\angC_{j=1}^n;\C}.
\end{equation}
\item[Pseudo Symmetry Isomorphisms] It has an $n$-linear natural isomorphism 
\begin{equation}\label{Fsigmanlinear}
\begin{tikzpicture}[xscale=1,yscale=1,baseline={(x11.base)}]
\def\d{25} 
\draw[0cell=.9]
(0,0) node (x11) {\C^n}
(x11)++(2.5,0) node (b) {\C}
;
\draw[1cell=.8] 
(x11) edge[bend left=\d] node  {\left(\mu_{\tau\sigma}, \ang{\mu^2_{\tau\sigma; j}}_{j=1}^n\right)} (b)
(x11) edge[bend right=\d] node[swap] {\left(\mu_\tau, \ang{\mu^2_{\tau; j}}_{j=1}^n\right) \cdot \sigma} (b)
;
\draw[2cell=.9]
node[between=x11 and b at .4, rotate=-90, 2label={above,F_{\sigma; \tau}}] {\Rightarrow}
;
\end{tikzpicture}
\inspace \permcat\scmap{\angC_{j=1}^n;\C}.
\end{equation}
Here $\big(\mu_\tau, \ang{\mu^2_{\tau; j}}_{j=1}^n\big) \cdot \sigma$ is the right $\sigma$-action on $\big(\mu_\tau, \ang{\mu^2_{\tau; j}}_{j=1}^n\big)$ as defined in \cref{permcatsigmaaction,fsigmatwoj}.
\end{description}
This finishes the description of the list of data of a pseudo symmetric $E_\infty$-algebra in $\permcat$.

The conditions \cref{Fcomponentfunctor} through \cref{Fbotequivariance} for $n \geq 0$ and $\sigma,\tau,\pi \in \Sigma_n$ are interpreted using the data \cref{Coplusebeta,musigmanlinear,musigmatwoj,musigmataunlinear,Fsigmanlinear}.  Below we state explicitly the three most nontrivial conditions from that list.
\begin{description}
\item[Composition] The composition axiom on objects \cref{compaxiomobject} means that the following diagram of functors is commutative.
\begin{equation}\label{prodmusigmajmusigma}
\begin{tikzpicture}[xscale=1,yscale=1,vcenter]
\def\d{25} 
\draw[0cell=.9]
(0,0) node (x11) {\C^{k_1 + \cdots + k_n}}
(x11)++(3,0) node (x12) {\C}
(x11)++(1.5,-1) node (x2) {\C^n}
;
\draw[1cell=.8] 
(x11) edge node {\mu_{\sigmatil \cdot \sigmatimes}} (x12)
(x11) edge[shorten <=-.5ex] node[swap,pos=.3] {\txprod_{j=1}^n \, \mu_{\sigma_j}} (x2)
(x2) edge node[swap] {\mu_\sigma} (x12)
;
\end{tikzpicture}
\end{equation}
Moreover, their corresponding linearity constraints are equal.  The linearity constraints for the composite $\mu_\sigma \circ (\txprod_{j=1}^n\, \mu_{\sigma_j})$ are defined in \cref{ffjlinearity}.

The composition axiom on morphisms \cref{compaxiommorphism} is the analog of \cref{prodmusigmajmusigma} with each functor replaced by its 2-cell variant in \cref{musigmataunlinear}.
\item[Top Equivariance]
With the notation
\[k = k_1 + \cdots + k_n \andspace k_\sigma = k_{\sigma(1)} + \cdots + k_{\sigma(n)},\]
the top equivariance axiom \cref{Ftopequivariance} is the following equality of natural isomorphisms.
\[\begin{tikzpicture}[xscale=1,yscale=1,vcenter]
\def\d{10} \def\h{1.9} \def\g{1.4} \def\v{.7}
\draw[0cell=.9]
(0,0) node (x11) {\C^{k_{\sigma}}}
(x11)++(\h,\v) node (x12) {\C}
(x11)++(\h,-\v) node (x2) {\C^k}
;
\draw[1cell=.8]
(x11) edge[bend left=\d] node[pos=.3] {\mu_{\tautil \cdot \tautimes \cdot \sigmabar}} (x12)
(x11) edge[bend right=\d] node[swap,pos=.3] {\sigmabar} (x2)
(x2) edge node[swap] (m) {\mu_{\tautil \cdot \tautimes}} (x12)
;
\draw[2cell=.9]
node[between=x11 and m at .3, rotate=-90, 2label={above, F_{\sigmabar;\, \tautil \cdot \tautimes}}] {\Rightarrow}
;
\draw[0cell=1.5]
(x11)++(\h+1,0) node (eq) {=}
;
\draw[0cell=.9]
(eq)++(1,0) node (y11) {\C^{k_{\sigma}}}
(y11)++(2.4,0) node (y12) {\C^n}
(y12)++(\g,\v) node (y13) {\C}
(y12)++(\g,-\v) node (y2) {\C^n}
;
\draw[1cell=.8]
(y11) edge node {\txprod_{j=1}^n \, \mu_{\tau_{\sigma(j)}}} (y12)
(y12) edge[bend left=\d] node[pos=.3] {\mu_{\tau\sigma}} (y13)
(y12) edge[bend right=\d] node[swap,pos=.3] {\sigma} (y2)
(y2) edge node[swap] (mt) {\mu_{\tau}} (y13)
;
\draw[2cell=.9]
node[between=y12 and mt at .33, rotate=-90, 2label={above, F_{\sigma;\, \tau}}] {\Rightarrow}
;
\end{tikzpicture}\]
In the above equality, the top, respectively, bottom, boundary composites are equal by the left, respectively, right, vertical equalities in \cref{Ftopequivariance}, along with \cref{permcatsigmaaction,permcatcomposite}.
\item[Bottom Equivariance]
The bottom equivariance axiom \cref{Fbotequivariance} is the following equality of natural isomorphisms.
\[\begin{tikzpicture}[xscale=1,yscale=1,vcenter]
\def\d{10} \def\h{1.9} \def\g{1.4} \def\v{.7}
\draw[0cell=.9]
(0,0) node (x11) {\C^k}
(x11)++(\h,-\v) node (x12) {\C}
(x11)++(0,-2*\v) node (x2) {\C^k}
;
\draw[1cell=.8]
(x11) edge[bend left=\d] node[pos=.5] {\mu_{\sigmatil \cdot \sigmatimes \cdot \tautimes}} (x12)
(x11) edge node[swap] (t) {\tautimes} (x2)
(x2) edge[bend right=\d] node[swap] {\mu_{\sigmatil \cdot \sigmatimes}} (x12)
;
\draw[2cell=.9]
node[between=t and x12 at .27, rotate=-90, 2label={above, F_{\tautimes;\, \sigmatil \cdot \sigmatimes}}] {\Rightarrow}
;
\draw[0cell=1.5]
(x12)++(1,0) node (eq) {=}
;
\draw[0cell=.9]
(eq)++(1,\v) node (y11) {\C^k}
(y11)++(1.4*\h,-\v) node (y12) {\C^n}
(y11)++(0,-2*\v) node (y2) {\C^k}
(y12)++(\g,0) node (y13) {\C}
;
\draw[1cell=.8]
(y11) edge[bend left=\d] node {\smallprod_{j=1}^n \, \mu_{\sigma_j \tau_j}} (y12)
(y11) edge node[swap] (tt) {\tautimes} (y2)
(y2) edge[bend right=\d] node[swap] {\smallprod_{j=1}^n \, \mu_{\sigma_j}} (y12)
(y12) edge node {\mu_\sigma} (y13)
;
\draw[2cell=.9]
node[between=tt and y12 at .25, rotate=-90, 2label={above, \txprod_{j=1}^n \, F_{\tau_j;\, \sigma_j}}] {\Rightarrow}
;
\end{tikzpicture}\]
In the above equality, the top, respectively, bottom, boundary composites are equal by the left, respectively, right, vertical equalities in \cref{Fbotequivariance}, along with \cref{permcatsigmaaction,permcatcomposite}.
\end{description}
This finishes the description of a pseudo symmetric $E_\infty$-algebra in $\permcat$.
\end{explanation}

\part*{Bibliography and Indices}
\appendix
\chapter{Open Questions}
\label{ch:questions}

\chapquote{``I enjoy questions that seem honest, even when they admit or reveal confusion, in preference to questions that appear designed to project sophistication."}{Bill Thurston}{\href{https://mathoverflow.net/users/9062/bill-thurston}{MathOverflow user profile}}

In this chapter, we discuss open questions related to the topics of this work: enriched multicategories, bipermutative categories, the Grothendieck construction, and inverse $K$-theory.

\begin{question}[Multi-Bicategories]\label{qu:multibicat}
Each $\Cat$-multicategory has an underlying 1-ary 2-category (\cref{ex:unarycategory}).  A \emph{bicategory}\index{bicategory} \cite[2.1.3]{johnson-yau} is a pseudo version of a 2-category.
\begin{itemize}
\item What is a pseudo variant\index{multicategory!pseudo variant} of a $\Cat$-multicategory?  Do theorems for bicategories have multicategorical variants?
\end{itemize}
More concretely, define
\begin{itemize}
\item \emph{multi-bicategories}\index{multi-bicategory}\index{bicategory!multi-} as the bicategorical variants of $\Cat$-multicategories and the multiple-input variants of bicategories, 
\item \emph{lax $\Cat$-multifunctors} between multi-bicategories, generalizing \cref{def:pseudosmultifunctor,expl:laxsmf}, as the lax variants of $\Cat$-multifunctors and the multiple-input variants of lax functors \cite[4.1.2]{johnson-yau},
\item \emph{lax $\Cat$-multinatural transformations} between lax $\Cat$-multifunctors as the lax variants of pseudo symmetric $\Cat$-multinatural transformations (\cref{def:psmftransformation}) and the multiple-input variants of lax transformations \cite[4.2.1]{johnson-yau}, and so on.  
\end{itemize} 
Then use these definitions to prove multi-bicategorical variants of
\begin{itemize}
\item the Bicategorical Pasting Theorem\index{Bicategorical!Pasting Theorem} \cite[3.6.6]{johnson-yau},
\item the Bicategorical Quillen Theorem A\index{Bicategorical!Quillen Theorem A} \cite[7.3.1]{johnson-yau},
\item the Whitehead Theorem for Bicategories\index{Bicategorical!Whitehead Theorem} \cite[7.4.1]{johnson-yau}, generalizing \cref{thm:multiwhitehead,thm:psmultiwhitehead},
\item the Bicategorical Yoneda Lemma\index{Bicategorical!Yoneda Lemma} \cite[8.3.16]{johnson-yau}, and
\item the Bicategorical Coherence Theorem\index{Bicategorical!Coherence Theorem} \cite[8.4.1]{johnson-yau}.
\end{itemize} 
There may be a conceptual approach to some of these questions via 2-monads \cite[6.5]{johnson-yau}.
\begin{itemize}
\item Is there a 2-monad\index{2-monad} with
\begin{itemize}
\item $\Cat$-multicategories with a given object set $X$ as strict algebras and
\item multi-bicategories with object set $X$ as pseudo algebras? 
\end{itemize} 
\end{itemize}
If the answer is yes, then the coherence theorems in \cite{bkp,power-coherence} may apply to this 2-monad.
\end{question}

\begin{question}[Pseudo Symmetry]\label{qu:pseudo}\index{pseudo symmetric!Cat-multifunctor@$\Cat$-multifunctor}\index{Cat-multifunctor@$\Cat$-multifunctor!pseudo symmetric}
Recall the notion of a pseudo symmetric $\Cat$-multifunctor $(F,\{F_\sigma\})$ in \cref{def:pseudosmultifunctor}.
\begin{itemize}
\item Is it possible for $F$---consisting of its object assignment and component functors---to admit several different pseudo symmetric $\Cat$-multifunctor structures $\{F_\sigma\}$?
\end{itemize}
In \cref{ex:idpsm} we equip a $\Cat$-multifunctor $F$ with identity natural transformations $F_\sigma = \Id$ to make $(F,\{\Id\})$ into a pseudo symmetric $\Cat$-multifunctor.  Here is a special case of the previous question.
\begin{itemize}
\item Is it possible for a $\Cat$-multifunctor $F$ to admit generally non-identity pseudo symmetry isomorphisms $\{F_\sigma\}$ such that $(F,\{F_\sigma\})$ is a pseudo symmetric $\Cat$-multifunctor different from the one in \cref{ex:idpsm}?\defmark
\end{itemize}
\end{question}

\begin{question}[Permutative Opfibrations]\label{qu:permopfibration}\index{permutative opfibration}
By the Grothendieck Fibration Theorem \cite[9.5.6]{johnson-yau}, there is a 2-monad on $\Cat/\B$ whose pseudo, respectively, strict, algebras are in canonical bijections with cloven, respectively, split, opfibrations over $\B$.
\begin{itemize}
\item Prove the analogous theorem for \emph{permutative} opfibrations and the cloven variant in \cref{def:permutativefibration}. 
\end{itemize}
In a permutative opfibration $P \cn \E \to \cD$, the functor $P$ is a strict symmetric monoidal functor between permutative categories.  Therefore, instead of $\Cat/\B$, one should consider the over-category of small permutative categories and strict symmetric monoidal functors over $\cD$.
\end{question}

\begin{question}[Pseudo Bipermutative-Indexed Categories]\label{qu:dcatpseudo}\index{diagram category!generalization}
For a small tight bipermutative category $\cD$, the $\Cat$-multicategory $\DCat$ (\cref{thm:dcatcatmulticat}) has objects, $n$-ary 1-cells, and $n$-ary 2-cells given by, respectively, additive symmetric monoidal functors, additive natural transformations, and additive modifications (\cref{def:additivesmf,def:additivenattr,def:additivemodification}).
\begin{itemize}
\item If the objects, $n$-ary 1-cells, and $n$-ary 2-cells are replaced by appropriate versions of, respectively, pseudo functors, strong transformations, and modifications \cite[4.1.2, 4.2.1, and 4.4.1]{johnson-yau}, does $\DCat$ still form a $\Cat$-multicategory or a multi-bicategory (\cref{qu:multibicat})?
\end{itemize}
One can simply replace an additive symmetric monoidal functor in \cref{def:additivesmf} by a pseudo functor.  However, each of \cref{def:additivenattr,def:additivemodification} has a unity axiom and an additivity axiom.  Their pseudo variants are, therefore, more subtle to formulate.
\end{question}

\begin{question}[Pseudo Permutative Opfibrations]\label{qu:pfibdpseudo}\index{permutative opfibration!generalization}
For a small tight bipermutative category $\cD$, the non-symmetric $\Cat$-multicategory $\pfibd$ (\cref{thm:pfibdmulticat}) has objects, $n$-ary 1-cells, and $n$-ary 2-cells given by, respectively, small permutative opfibrations over $\cD$, opcartesian strong $n$-linear functors, and opcartesian $n$-linear transformations (\cref{def:permutativefibration,def:opcartnlinearfunctor,def:opcartnlineartr}).  Similar to \cref{qu:dcatpseudo}, we ask the following.
\begin{itemize}
\item If the objects, $n$-ary 1-cells, and $n$-ary 2-cells are replaced by, respectively, cloven opfibrations and some pseudo versions of $n$-ary 1-cells and $n$-ary 2-cells, does $\pfibd$ still form a non-symmetric $\Cat$-multicategory or a non-symmetric multi-bicategory (\cref{qu:multibicat})?\defmark
\end{itemize}
\end{question}

\begin{question}[Pseudo Grothendieck $\Cat$-Multiequivalences]\label{qu:grodpseudo}\index{Grothendieck construction!generalization}
For each small tight bipermutative category $\cD$, there is a non-symmetric $\Cat$-multiequivalence given by the Grothendieck construction (\cref{thm:dcatpfibdeq}) 
\[\begin{tikzpicture}[xscale=2.75,yscale=1.2,baseline={(dc.base)}]
\draw[0cell=.9]
(0,0) node (dc) {\DCat}
(dc)++(1,0) node (pf) {\pfibd.};
\draw[1cell=.9]  
(dc) edge node {\grod} node[swap] {\sim} (pf);
\end{tikzpicture}\]
If \cref{qu:dcatpseudo,qu:pfibdpseudo} can both be answered in the positive, then we ask the following.
\begin{itemize}
\item Can the non-symmetric $\Cat$-multiequivalence $\grod$ be extended to the pseudo setting?
\end{itemize}
In fact, one can regard this question as a litmus test for the correctness of the definitions of the pseudo variants of $\DCat$ and $\pfibd$ in \cref{qu:dcatpseudo,qu:pfibdpseudo}.
\end{question}

\begin{question}[Global Grothendieck $\Cat$-Multiequivalences]\label{qu:globular}\index{non-symmetric!Cat-multiequivalence@$\Cat$-multiequivalence}\index{Cat-multiequivalence@$\Cat$-multiequivalence!global}
In the construction $\grod$ the small tight bipermutative category $\cD$ is fixed.  So a natural question is whether $\cD$ can be used as a variable.
\begin{itemize}
\item Is there a global version of the non-symmetric $\Cat$-multiequivalence $\grod$ where $\cD$ varies in tight bipermutative categories?
\end{itemize}
To answer this question, one would need
\begin{itemize}
\item a notion of functors between bipermutative categories, which should be analogous to symmetric bimonoidal functors in \cite[Section 5.1]{cerberusI};
\item a global version of the $\Cat$-multicategory $\DCat$ (\cref{thm:dcatcatmulticat}) where $\cD$ varies in tight bipermutative categories; and
\item a global version of the non-symmetric $\Cat$-multicategory $\pfibd$ (\cref{thm:pfibdmulticat}) where $\cD$ varies in tight bipermutative categories.
\end{itemize} 
Moreover, if \cref{qu:grodpseudo} has a positive answer, then one can also ask whether there is a global variant where $\cD$ varies.
\end{question}

\begin{question}[Operadic Grothendieck $\Cat$-Multiequivalences]\label{qu:grooperad}
Small tight bipermutative categories are $E_\infty$-algebras in $\permcatsg$ (\cref{ex:Einftyalgebras} \eqref{ex:Einftyalgebras-iii}) with respect to the Barratt-Eccles operad $\BE$ (\cref{def:barratteccles}).  
\begin{itemize}
\item Is there a generalization of the non-symmetric $\Cat$-multiequivalence $\grod$ where $\cD$ is a fixed $\Q$-algebra in $\permcatsg$ for a categorical operad $\Q$?
\end{itemize}
In addition to the Barratt-Eccles $E_\infty$-operad $\BE$, other categorical operads of interest include the $E_n$-operads in \cite[13.2.1]{cerberusIII} and the $\G$-monoidal category operads for an action operad $\G$ in \cite[19.2.1]{yau-inf-operad}.  If the answer to this question is yes, then we can extend \cref{qu:globular} to this setting.
\begin{itemize}
\item Is there a global version of the non-symmetric $\Cat$-multiequivalence $\grod$ where $\cD$ varies among $\Q$-algebras in $\permcatsg$?
\end{itemize}
One can extend this question further by
\begin{itemize}
\item varying the categorical operad $\Q$ and
\item considering it in the context of \cref{qu:grodpseudo}.\defmark
\end{itemize} 
\end{question}

\begin{question}[$\Cat$-Multifunctorial Inverse $K$-Theory]\label{qu:symmetricP}\index{inverse K-theory@inverse $K$-theory!multifunctorial variant}
Inverse $K$-theory $\cP$ is a pseudo symmetric $\Cat$-multifunctor but \emph{not} a $\Cat$-multifunctor in the symmetric sense (\cref{thm:invKpseudosym}).  
\begin{itemize}
\item Is there an equivalent variant of $\cP$ that is a $\Cat$-multifunctor in the symmetric sense? 
\end{itemize}
More precisely, this question asks for the existence of
\begin{itemize}
\item a $\Cat$-multifunctor 
\[\begin{tikzcd}[column sep=large]
\Gacat \ar{r}{\cP'} & \permcatsg
\end{tikzcd}\]
together with
\item a natural transformation $\theta \cn \cP \to \cP'$ or $\cP' \to \cP$ that is componentwise a stable equivalence in $\permcatsg$.  
\end{itemize}  
A morphism in $\permcatsg$ is a \emph{stable equivalence}\index{stable equivalence} if its image under Segal $K$-theory\index{Segal K-theory@Segal $K$-theory} \cite{segal} is a stable equivalence; see \cite{gjo1,johnson-yau-invK,johnson-yau-multiK,bousfield_friedlander,mandell_inverseK,thomason} for more discussion of stable equivalences.  Since $\cP$ is a homotopy inverse of Segal $K$-theory \cite{mandell_inverseK}, $\theta$ ensures that $\cP'$ is also a homotopy inverse of Segal $K$-theory.  Moreover, being a $\Cat$-multifunctor, $\cP'$ preserves any algebraic structure---such as $E_\infty$-algebras and $E_n$-algebras---parametrized by a $\Cat$-multifunctor, strengthening \cref{cor:Ppreservesnonsym,cor:Ppreservespseudosym}.
\end{question}

\backmatter
\bibliographystyle{sty/amsalpha3}
\bibliography{references}

\newcommand{\fact}[2]{\noindent({#2}) {#1}}

\newcommand{\chapNumName}[1]{\medskip\begin{center}\textbf{\Cref{#1}}.  \nameref{#1}\end{center}}
\newcommand{\thm}[1]{\textbf{#1}.}



\chapter*{List of Main Facts}

\chapNumName{ch:prelim}

\fact{Each monoidal category satisfies $\lambda_{\tu} = \rho_{\tu}$.}{\ref{lambda=rho}}

\fact{Each monoidal category satisfies the left and right unity axioms.}{\ref{moncat-other-unit-axioms}}

\fact{Each braided monoidal category satisfies $\rho = \lambda \xi_{-,\tu}$ and $\lambda = \rho \xi_{\tu,-}$}{\ref{braidedunity}}

\fact{A symmetric monoidal category is precisely a braided monoidal category that satisfies the symmetry axiom.}{\ref{rk:smcat}}

\fact{$\Cat$ is a symmetric monoidal closed category.}{\ref{ex:cat}}

\fact{$\Cat$ is a 2-category.}{\ref{ex:catastwocategory}}

\fact{$\VCat$ is a 2-category.}{\ref{ex:vcatastwocategory}}

\fact{A locally small 2-category is precisely a $\Cat$-category.}{\ref{locallysmalltwocat}}

\fact{The diagram category $\DV$ is symmetric monoidal closed.}{\ref{thm:Day}}

\fact{$\DV$ is enriched, tensored, and cotensored over $\V$.}{\ref{thm:DV}}

\chapNumName{ch:bipermutative}

\fact{In a symmetric bimonoidal category, 12 of the 24 axioms follow from the other 12 axioms.}{\ref{thm:sbc-redundant-axioms}}

\fact{There is a canonical bijection between tight bipermutative categories and tight symmetric bimonoidal categories with permutative additive structure, permutative multiplicative structure, $\lambdadot = 1$, and $\rhodot = 1$}{\ref{thm:smcbipermuative}}

\fact{There is a small tight bipermutative category $\Finsk$ with unpointed finite sets $\ufs{n}$ for $n \geq 0$ as objects and functions as morphisms.}{\ref{ex:finsk}}

\fact{There is a small tight bipermutative category $\Fset$ with the same objects as $\Finsk$ and permutations as morphisms.}{\ref{ex:Fset}}

\fact{There is a small tight bipermutative category $\Fskel$ with pointed finite sets $\ord{n}$ for $n \geq 0$ as objects and pointed functions as morphisms.}{\ref{ex:Fskel}}

\fact{Mandell's category $\cA$ for inverse $K$-theory is a small tight bipermutative category.}{\ref{ex:mandellcategory}}

\fact{In a symmetric bimonoidal category $\C$ in which the value of each $\delta$-prime edge is a monomorphism, any two parallel paths in $\grx$ with a regular domain have the same value in $\C$.}{\ref{thm:laplaza-coherence-1}}

\fact{Laplaza's Coherence \cref{thm:laplaza-coherence-1} applies to tight bipermutative categories.}{\ref{ex:laplazacoherence}}

\chapNumName{ch:multicat}

\fact{The terminal multicategory $\Mterm$ consists of a single object and a single $n$-ary operation for each $n$.}{\ref{definition:terminal-operad-comm}}

\fact{Each object in a $\V$-multicategory has an endomorphism $\V$-operad.}{\ref{example:enr-End}}

\fact{Each permutative category has an endomorphism multicategory.}{\ref{ex:endc}}

\fact{Each $\V$-multicategory has an underlying $\V$-category.}{\ref{ex:unarycategory}}

\fact{There is a 2-category with small $\V$-multicategories as objects.}{\ref{v-multicat-2cat}}

\fact{$\VMulticat$ and $\VMulticatns$ are complete and cocomplete.}{\ref{vmulticatbicomplete}}

\fact{The initial $\V$-multicategory has an empty set of objects.  The terminal $\V$-multicategory has one object and each multimorphism object given by the terminal object in $\V$.}{\ref{ex:vmulticatinitialterminal}}

\fact{A $\Cat$-multinatural transformation consists of component 1-ary 1-cells that satisfy two $\Cat$-naturality conditions, one for objects and one for morphisms.}{\ref{expl:catmultitransformation}}

\fact{Each $\Cat$-multiequivalence is essentially surjective on objects and an isomorphism on each multimorphism category.}{\ref{catmultieqesff}}

\fact{A $\Cat$-multifunctor is a $\Cat$-multiequivalence if and only if it is essentially surjective on objects and an isomorphism on each multimorphism category.}{\ref{thm:multiwhitehead}}

\fact{A 2-functor between 2-categories is a 2-equivalence if and only if it is essentially surjective on objects and an isomorphism on each hom category.}{\ref{thm:iiequivalence}}

\chapNumName{ch:pseudosymmetry}

\fact{A pseudo symmetric $\Cat$-multifunctor strictly preserves the colored units and composition and has pseudo symmetry isomorphisms that satisfy four coherence axioms.}{\ref{def:pseudosmultifunctor}}

\fact{There is a canonical bijection between $\Cat$-multifunctors and pseudo symmetric $\Cat$-multifunctors with identity pseudo symmetry isomorphisms.}{\ref{ex:idpsm}}

\fact{A pseudo symmetric $\Cat$-multifunctor yields a non-symmetric $\Cat$-multifunctor.}{\ref{ex:nscatmultifunctors}}

\fact{Pseudo symmetric $\Cat$-multifunctors are closed under composition, which is associative and unital.}{\ref{pseudosmfcomposite}}

\fact{A pseudo symmetric $\Cat$-multinatural transformation satisfies two $\Cat$-naturality conditions and the pseudo symmetry preservation axiom.}{\ref{def:psmftransformation}}

\fact{Identities and vertical composites are well defined for pseudo symmetric $\Cat$-multinatural transformations.}{\ref{pstvcomp}}

\fact{Pseudo symmetric $\Cat$-multinatural transformations are closed under horizontal composition.}{\ref{psthcomp}}

\fact{There is a 2-category $\catmulticatps$ with small $\Cat$-multicategories as objects, pseudo symmetric $\Cat$-multifunctors as 1-cells, and pseudo symmetric $\Cat$-multinatural transformations as 2-cells.}{\ref{thm:catmulticatps}}

\fact{Each pseudo symmetric $\Cat$-multiequivalence is essentially surjective on objects and an isomorphism on each multimorphism category.}{\ref{pscatmultieqesff}}

\fact{A pseudo symmetric $\Cat$-multifunctor is a pseudo symmetric $\Cat$-multiequivalence if and only if it is essentially surjective on objects and an isomorphism on each multimorphism category.}{\ref{thm:psmultiwhitehead}}

\chapNumName{ch:diagram}

\fact{$\DV$ is a $\V$-multicategory.}{\ref{thm:permindexedcat}}

\fact{An additive natural transformation is a natural transformation that satisfies the unity axiom and the additivity axiom.}{\ref{def:additivenattr}}

\fact{An additive modification is a modification that satisfies the unity axiom and the additivity axiom.}{\ref{def:additivemodification}}

\fact{$\DCat\smscmap{\angX;Z}$ is a subcategory of $\Dtecat\smscmap{\angX;Z}$.}{\ref{dcatnarycategory}}

\fact{$\DCat$ is not a symmetric monoidal category.}{\ref{expl:dcatbipermzero}}

\fact{$\DCat$ is a $\Cat$-multicategory.}{\ref{thm:dcatcatmulticat}}

\fact{There is a forgetful $\Cat$-multifunctor $\DCat \to \Dtecat$.}{\ref{dcatforgetdtecat}}

\chapNumName{ch:multigro}

\fact{A 1-linear functor is precisely a strictly unital symmetric monoidal functor.}{\ref{ex:onelinearfunctor}}

\fact{The multiplicative structure in a tight bipermutative category is a strong $n$-linear functor relative to the additive structure.}{\ref{ex:bipermnlinearfunctor}}

\fact{$\permcat$ and $\permcatsg$ are $\Cat$-multicategories.}{\ref{thm:permcatenrmulticat}}

\fact{The Grothendieck construction of an indexed category is a category.}{\ref{groconstcategory}}

\fact{The Grothendieck construction of a symmetric monoidal functor to $\Cat$ is a permutative category.}{\ref{grodxpermutative}}

\fact{The Grothendieck construction of an additive natural transformation is a strong $n$-linear functor.}{\ref{grodphifunctor}}

\fact{The Grothendieck construction of an additive modification is an $n$-linear natural transformation.}{\ref{grodPhinlinear}}

\fact{The Grothendieck construction defines a functor from $\DCat\smscmap{\angX;Z}$ to $\permcatsg\smscmap{\ang{\grod X}; \grod Z}$.}{\ref{grodmultimorphismfunctor}}

\fact{The Grothendieck construction preserves colored units.}{\ref{grodpreservesunits}}

\fact{The Grothendieck construction preserves composition.}{\ref{grodpreservescomposition}}

\fact{The Grothendieck construction has pseudo symmetry isomorphisms.}{\ref{grodsigmanaturaliso}}

\fact{The Grothendieck construction $\grod$ is a pseudo symmetric $\Cat$-multifunctor from $\DCat$ to $\permcatsg$.  It is a $\Cat$-multifunctor if the multiplicative braiding in $\cD$ is the identity.  Conversely, if $\grod$ is a $\Cat$-multifunctor, then $x \otimes y = y \otimes x$ for all objects $x,y \in \cD$.}{\ref{thm:grocatmultifunctor}}

\chapNumName{ch:permfib}

\fact{A permutative opfibration is a split opfibration with additional properties but no extra structure.}{\ref{expl:permutativefibration}}

\fact{An opcartesian $n$-linear functor is an $n$-linear functor with additional properties but no extra structure.}{\ref{expl:opcartnlinearfunctor} \eqref{expl:opcartnlinearfunctor-one}}

\fact{An opcartesian 0-linear functor is an object in the $P_0$-preimage of $\tu \in \cD$.}{\ref{expl:opcartnlinearfunctor} \eqref{expl:opcartnlinearfunctor-0}}

\fact{In the definition of an opcartesian $n$-linear functor, axiom \cref{def:opcartnlinearfunctor-i} follows from axiom \cref{def:opcartnlinearfunctor-iii}.}{\ref{preserveopcmorphisms}}

\fact{An opcartesian 1-linear functor $F$ is precisely a strict symmetric monoidal functor that satisfies $P_1 = P_0 F$ and preserves chosen opcartesian lifts.}{\ref{ex:idoponelinear}}

\fact{The multiplicative structure of a small tight bipermutative category defines an opcartesian strong $n$-linear functor.}{\ref{ex:bipermopnlinear}}

\fact{An opcartesian $n$-linear transformation is an $n$-linear natural transformation with an additional property but no extra structure.}{\ref{expl:opcartnlineartr} \eqref{expl:opcartnlineartr-i}}

\fact{An opcartesian 0-linear transformation is a morphism in the $P_0$-preimage of $\tu \in \cD$.}{\ref{expl:opcartnlineartr} \eqref{expl:opcartnlineartr-ii}}

\fact{$\pfibd\smscmap{\ang{P};P_0}$ is a subcategory of $\permcatsg\smscmap{\angC;\C_0}$.}{\ref{pfibdcomponentwelldef}}

\fact{Composition in $\pfibd$ is well defined on opcartesian strong multilinear functors.}{\ref{pfibdcompobjwelldef}}

\fact{Composition in $\pfibd$ is well defined on opcartesian multilinear transformations.}{\ref{pfibdcompmorwelldef}}

\fact{If the multiplicative braiding in $\cD$ is the identity, then the symmetric group action in $\pfibd$ is well defined on opcartesian strong multilinear functors.  Conversely, if the symmetric group action in $\pfibd$ is well defined on opcartesian strong multilinear functors, then $\otimes = \tensorop$ in $\cD$.}{\ref{pfibdsymwelldefobj}}

\fact{If the multiplicative braiding in $\cD$ is the identity, then the symmetric group action in $\pfibd$ is well defined on opcartesian multilinear transformations.  Conversely, if the symmetric group action in $\pfibd$ is well defined on opcartesian multilinear transformations, then $\otimes = \tensorop$ in $\cD$.}{\ref{pfibdsymwelldefmor}}

\fact{$\pfibd$ is a non-symmetric $\Cat$-multicategory.  If the multiplicative braiding in $\cD$ is the identity, then $\pfibd$ is a $\Cat$-multicategory.  Conversely, if $\pfibd$ is a $\Cat$-multicategory, then $\otimes = \tensorop$ in $\cD$.}{\ref{thm:pfibdmulticat}}

\fact{There exists a non-symmetric $\Cat$-multifunctor $\U \cn \pfibd \to \permcatsg$.  If the multiplicative braiding in $\cD$ is the identity, then $\U$ is a $\Cat$-multifunctor.}{\ref{cor:pfibdtopermcatsg}}

\fact{The first-factor projection is a small permutative opfibration over $\cD$.}{\ref{grodxpermopfib}}

\fact{$\grod \cn \DCat\smscmap{\ang{};Z} \to \pfibd\smscmap{\ang{};U_Z}$ is an isomorphism.}{\ref{grodfactorarityzero}}

\fact{$\grod \cn \DCat\smscmap{\angX;Z} \to \pfibd\smscmap{\ang{U_X};U_Z}$ is a functor.}{\ref{grodfactorarityn}}

\fact{$\grod \cn \DCat \to \pfibd$ is a non-symmetric $\Cat$-multifunctor, which is, furthermore, a $\Cat$-multifunctor if the multiplicative braiding in $\cD$ is the identity.}{\ref{thm:dcatpfibd}}

\chapNumName{ch:gromultequiv}

\fact{$\grod \cn \DCat \to \pfibd$ is essentially surjective on objects.}{\ref{grodesssurjective}}

\fact{$\grod \cn \DCat\smscmap{\angX;Z} \to \pfibd\smscmap{\ang{U_X};U_Z}$ is injective on objects.}{\ref{grodarityninjobj}}

\fact{$\grod \cn \DCat\smscmap{\angX;Z} \to \pfibd\smscmap{\ang{U_X};U_Z}$ is injective on morphism sets.}{\ref{grodarityninjmor}}

\fact{$\grod \cn \DCat\smscmap{\angX;Z} \to \pfibd\smscmap{\ang{U_X};U_Z}$ is surjective on objects.}{\ref{grodobjectfullfphi}}

\fact{$\grod \cn \DCat\smscmap{\angX;Z} \to \pfibd\smscmap{\ang{U_X};U_Z}$ is surjective on morphism sets.}{\ref{grodPhiomega}}

\fact{$\grod \cn \DCat \to \pfibd$ is a non-symmetric $\Cat$-multiequivalence.  It is a $\Cat$-multiequivalence if the multiplicative braiding in $\cD$ is the identity.}{\ref{thm:dcatpfibdeq}}

\chapNumName{ch:multifunctorA}

\fact{$\pCat$ is a complete and cocomplete symmetric monoidal closed category.}{\ref{theorem:pC-monoidal}}

\fact{$\Gacat$ is the category of pointed functors $(\Fskel,\ord{0}) \to (\Cat,\boldone)$ and pointed natural transformations.}{\ref{def:gammacategory}}

\fact{Each $\Ga$-category is canonically a pointed functor $(\Fskel,\ord{0}) \to (\pCat,\boldone)$.}{\ref{pointedgammacat}}

\fact{$\Gacat$ is complete, cocomplete, and symmetric monoidal closed.  It is enriched, tensored, and cotensored over $\pCat$.  It is a $\pCat$-multicategory.}{\ref{thm:Gacatsmc}}

\fact{The monoidal unit diagram $\ftu$ satisfies $\ftu\ord{n} \iso \ord{n}$.}{\ref{expl:Ftu}}

\fact{The multimorphism categories of $\Gacat$ are given by $\GMap\brb{\txsma_j X_j,Z}$.}{\ref{Gacatnmorphism}}

\fact{Objects and morphisms in $\Gacat\smscmap{\ang;Z}$ are determined by their values at $1 \in \ord{1}$.}{\ref{Fordnjziotaj} and \ref{thetaordoneone}}

\fact{The symmetric group action in $\Gacat$ corresponds to permuting the domain factors.}{\ref{Gacatsigmaaction}}

\fact{Composition in $\Gacat$ is given by the pointed Day convolution product and composition in $\GMap$.}{\ref{Gacatgammacomp}}

\fact{Each $\Ga$-category $X$ yields a strictly unital strong symmetric monoidal functor $AX \cn \cA \to \Cat$.}{\ref{AXsmf}}

\fact{Each $\Ga$-category morphism $F$ yields a monoidal natural transformation $AF$.}{\ref{AFmonoidal}}

\fact{The underlying category of a $\Cat$-multicategory has the same objects and 1-ary 1-cells as morphisms.}{\ref{def:underlyingcategory}}

\fact{$A \cn \Gacat \to \ACat$ is a functor.}{\ref{Afunctor}}

\fact{$A \cn \Gacat\smscmap{\ang{};Z} \fto{\iso} \ACat\smscmap{\ang{};AZ}$ is an isomorphism.}{\ref{Aarityzeroiso}}

\fact{$A \cn \Gacat\smscmap{\angX;Z} \to \ACat\scmap{\ang{AX};AZ}$ is a functor that is injective on objects and morphism sets.}{\ref{Amultimorfunctor}}

\fact{$A$ preserves the symmetric group action.}{\ref{Apreservessymmetry}}

\fact{$A$ preserves the multicategorical composition.}{\ref{Apreservescomp}}

\fact{$A \cn \Gacat \to \ACat$ is a $\Cat$-multifunctor.}{\ref{thm:Acatmultifunctor}}

\chapNumName{ch:invK}

\fact{The inverse $K$-theory functor $\cP$ is the composite $\groa \circ A$.}{\ref{mandellinvK}}

\fact{$\cP$ extends to a pseudo symmetric $\Cat$-multifunctor $\groa \circ A$ that is not a $\Cat$-multifunctor.  As a non-symmetric $\Cat$-multifunctor, $\cP$ factors as $\U \circ \groa \circ A$.}{\ref{thm:invKpseudosym}}

\fact{$\cP$ preserves algebraic structures parametrized by non-symmetric $\Cat$-multifunctors.}{\ref{cor:Ppreservesnonsym}}

\fact{$\cP$ sends monoids in $\Gacat$ to small tight ring categories.}{\ref{ex:Preservesnonsym}}

\fact{$\cP$ preserves algebraic structures parametrized by pseudo symmetric $\Cat$-multifunctors.}{\ref{cor:Ppreservespseudosym}}

\fact{$\cP(-)^\sigma$ and $(\cP(-))^\sigma$ are different on objects in general.}{\ref{groasigmaafmx}}

\fact{The pseudo symmetry isomorphism $\cP_\sigma$ is the identity if and only if $\sigma$ is the identity permutation.}{\ref{Psigmagroasigma}}

\fact{An $E_\infty$-algebra is a $\Cat$-multifunctor from the Barratt-Eccles operad $\BE$.}{\ref{def:barratteccles}}

\fact{An $E_\infty$-algebra in a $\Cat$-multicategory is uniquely determined by an object, a $0$-ary 1-cell, a 2-ary 1-cell $\mu$, and a 2-ary 2-cell isomorphism $\mu \fto{\iso} \muop$, subject to several axioms.}{\ref{mapfrombarratteccles}}

\fact{Each $E_\infty$-algebra has an underlying monoid.}{\ref{expl:Einftymonoid}}

\fact{$E_\infty$-algebras in $\Cat$, $\permcat$, and $\permcatsg$ are, respectively, small permutative categories and small (tight) bipermutative categories.}{\ref{ex:Einftyalgebras}}

\fact{A pseudo symmetric $E_\infty$-algebra is a pseudo symmetric $\Cat$-multifunctor from $\BE$.}{\ref{def:psEinftyalgebra}}

\fact{$\cP$ preserves pseudo symmetric $E_\infty$-algebras.}{\ref{cor:Pefinitypsalg}}

\fact{A pseudo symmetric $E_\infty$-algebra consists of an object, an $n$-ary 1-cell $\mu_\sigma$ for each permutation $\sigma$, and invertible $n$-ary 2-cells $\mu_\sigma^\tau \cn \mu_\sigma \fto{\iso} \mu_\tau$ and $F_{\sigma;\tau} \cn \mu_{\tau\sigma} \fto{\iso} \mu_\tau \cdot \sigma$ for each pair of permutations $(\sigma,\tau)$, subject to several axioms.}{\ref{expl:psEinftyalgebra}}

\chapter*{List of Notations}

\newcommand{\entry}[3]{\scalebox{.75}{#1} \> \> \scalebox{.75}{\pageref{#2}} \> 
  \> \scalebox{.75}{#3}\\[-.2ex]}
\newcommand{\entryNoPage}[3]{\scalebox{.75}{#1} \> \> \> \> \scalebox{.75}{#3}\\[-.2ex]} 
\newcommand{\blob}{\> \> \> \> }
\newcommand{\stuff}[1]{\blob \scalebox{.75}{#1}\\[-.2ex]}
\newcommand{\Header}{\>\> \textbf{Page} \>\> \textbf{Description}\\}
\newcommand{\newchHeader}[1]{\blob\\ \textbf{\Cref{#1}}\Header}


\begin{tabbing}
\phantom{\textbf{Notation}} \= \hspace{1.5cm}\= \phantom{\textbf{Page}}\= \hspace{.6cm}\= \phantom{\textbf{Description}} \\

\blob\\
\textbf{Standard Notations}\blob\textbf{Description}\\
\entryNoPage{$\Ob(\C)$, $\Ob\C$}{not:objects}{objects in a category $\C$}
\entryNoPage{$\C(X,Y)$, $\C(X;Y)$}{not:morphisms}{set of morphisms $X \to Y$}
\entryNoPage{$1_X$}{not:idmorphism}{identity morphism}
\entryNoPage{$g \circ f$, $gf$}{notation:morphism-composition}{composition of morphisms}
\entryNoPage{$\iso$, $\fto{\iso}$}{not:iso}{an isomorphism}
\entryNoPage{$F \cn \C \to \D$}{def:functors}{a functor}
\entryNoPage{$\Id_{\C}$, $1_{\C}$}{not:idc}{identity functor}
\entryNoPage{$\boldone$}{ex:terminal-category}{terminal category}
\entryNoPage{$(\Set,\times,*)$}{notation:set}{category of sets and functions}
\entryNoPage{$\theta_X$}{thetax}{a component of a natural transformation $\theta$}
\entryNoPage{$1_F$}{not:idf}{identity natural transformation}
\entryNoPage{$\phi\theta$}{not:vcomp}{vertical composition of natural transformations}
\entryNoPage{$\theta' * \theta$}{not:hcomp}{horizontal composition of natural transformations}
\entryNoPage{$(L,R,\phi)$, $L \dashv R$}{notation:adjunction}{an adjunction}
\entryNoPage{$\eta$, $\epz$}{adjunction-unit}{unit and counit of an adjunction}
\entryNoPage{$\varnothing$, $\varnothing^{\C}$}{not:initialobj}{an initial object}
\entryNoPage{$\coprod$, $\amalg$}{not:coprod}{a coproduct}
\entryNoPage{$\prod$, $\smallprod$}{not:coprod}{a product}
\entryNoPage{$\Sigma_n$}{not:Sigman}{symmetric group on $n$ letters}

\newchHeader{ch:prelim}

\entry{$\cU$}{not:universe}{a Grothendieck universe}
\entry{$\otimes$}{notation:monoidal-product}{monoidal product}
\entry{$\tu$}{not:monoidalunit}{monoidal unit}
\entry{$\alpha$}{not:associativityiso}{associativity isomorphism}
\entry{$\lambda$, $\rho$}{not:unitisos}{unit isomorphisms}
\entry{$(X,\mu,\eta)$}{notation:monoid}{a monoid}
\entry{$\xi$}{notation:symmetry-iso}{symmetry isomorphism or braiding}
\entry{$\oplus$, $e$}{not:opluse}{monoidal product and monoidal unit in a permutative category}
\entry{$[-,-]$, $\Hom(-,-)$}{notation:internal-hom}{internal hom}
\entry{$\DC$}{not:DC}{diagram category}
\entry{$(\Cat, \times, \boldone,[,])$}{not:cat}{category of small categories}
\entry{$(F,F^2,F^0)$}{def:monoidalfunctor}{a monoidal functor}
\entry{$\Rightarrow$}{twocellnotation}{a natural transformation}
\entry{$m$}{not:enrcomposition}{composition in a $\V$-category}
\entry{$i_X$}{not:enridentity}{identity in a $\V$-category}
\entry{$\C \otimes \D$}{definition:vtensor-0}{tensor product of $\V$-categories}
\entry{$\ximid$}{not:ximid}{interchanging the middle two factors}
\entry{$\B_0$, $\B_1$, $\B_2$}{def:twocategory}{objects, 1-cells, and 2-cells in a 2-category $\B$}
\entry{$\alpha'\alpha$}{not:vcompiicell}{vertical composition of 2-cells}
\entry{$gf$}{not:hcompicell}{horizontal composition of 1-cells}
\entry{$\beta * \alpha$}{not:hcompiicell}{horizontal composition of 2-cells}
\entry{$\B(X,Y)$}{not:homcat}{a hom category}
\entry{$\Phi_X$}{not:componentiicell}{a component 2-cell of a modification $\Phi$}
\entry{$\Psi\Phi$}{modificationvcomp}{vertical composition of modifications}
\entry{$\Phi' * \Phi$}{modificationhcomp}{horizontal composition of modifications}
\entry{$(f,g,\eta,\epz)$}{def:equivalences-ii}{an adjunction in a 2-category}
\entry{$X \otimes Y$}{eq:dayconvolution}{Day convolution}
\entry{$\Dhom$}{eq:dayhom}{hom diagram}
\entry{$J$}{eq:dayunit}{unit diagram}
\entry{$X \otimes A$}{not:tensored}{an object of a tensored structure}
\entry{$X^A$}{not:cotensored}{an object of a cotensored structure}
\entry{$\ev_e$}{not:evaluation}{evaluation at $e$}

\newchHeader{ch:bipermutative}

\entry{$\Cplus$}{not:sbcadditive}{additive structure of a symmetric bimonoidal category}
\entry{$\Cte$}{not:sbcmultiplicative}{multiplicative structure of a symmetric bimonoidal category}
\entry{$\lambdadot$, $\rhodot$}{sbc-multiplicative-zero}{multiplicative zeros}
\entry{$\deltal$, $\deltar$}{sbc-distributivity}{distributivity morphisms}
\entry{$\Cplus$, $(\C,\oplus,\zero,\xiplus)$}{not:Cplusring}{additive structure of a ring category}
\entry{$\Cte$, $(\C,\otimes,\tu)$}{not:Ctering}{multiplicative structure of a ring category}
\entry{$\fal$, $\far$}{ringcatfactorization}{factorization morphisms in a ring category}
\entry{$\Cte$, $(\C,\otimes,\tu,\xitimes)$}{not:Ctebiperm}{multiplicative structure of a bipermutative category}
\entry{$\Finsk$}{ex:finsk}{bipermutative category of finite sets $\ufs{n}$ and functions}
\entry{$\ufs{n}$}{unpointedfs}{unpointed finite set $\{1,\ldots,n\}$ with $\ufs{0} = \emptyset$}
\entry{$\Fset$}{ex:Fset}{bipermutative category of permutations}
\entry{$\Fskel$}{ex:Fskel}{bipermutative category of pointed finite sets $\ord{n}$ and pointed functions}
\entry{$\ord{n}$}{ordn}{pointed finite set $\{0,\ldots,n\}$}
\entry{$\cA$}{ex:mandellcategory}{bipermutative category for inverse $K$-theory}
\entry{$\ang{}$, $\emptyset$}{not:emptyseq}{empty sequence}
\entry{$\mathtt{m}$}{summn}{$m_1 + \cdots + m_p$}
\entry{$\ang{x_i}_{i=1}^p$}{Amtensorn}{$(x_1,\ldots,x_p)$}
\entry{$\grelx$}{not:grelx}{elementary graph}
\entry{$\eelx$}{not:eelx}{set of elementary edges}
\entry{$\eelfrx$}{not:eelfrx}{free $\plustimes$-algebra of $\eelx$}
\entry{$\pedgex$}{not:pedgex}{set of prime edges}
\entry{$\grx$}{not:grx}{graph of $X$}
\entry{$\stx$}{not:stx}{strict $\plustimes$-algebra of $X$}
\entry{$\supp$}{support}{support}

\newchHeader{ch:multicat}

\entry{$\Prof(D)$}{notation:profs}{class of $D$-profiles}
\entry{$\angd = (d_1,\ldots,d_n)$}{notation:us}{a $D$-profile of length $n$}
\entry{$\oplus$}{not:concat}{concatenation of profiles}
\entry{$\smscmap{\angd;d'}$}{notation:duc}{an element in $\Prof(D) \times D$}
\entry{$\M\smscmap{\angx;x'}$}{notation:enr-cduc}{$n$-ary operation object}
\entry{$\angx\sigma$}{enr-notation:c-sigma}{right symmetric group action by $\sigma$}
\entry{$\operadunit_x$}{notation:enr-unit-c}{$x$-colored unit}
\entry{$\gamma$}{notation:enr-multicategory-composition}{composition in a $\V$-multicategory}
\entry{$\sigma\ang{k_{\sigma(1)},\ldots,k_{\sigma(n)}}$}{blockpermutation}{block permutation induced by $\sigma$}
\entry{$\tau_1 \times \cdots \times \tau_n$}{blocksum}{block sum}
\entry{$\M_n$}{not:nthobject}{$n$-ary operation object of a $\V$-operad}
\entry{$\Mterm$}{definition:terminal-operad-comm}{terminal multicategory}
\entry{$\End(x)$}{example:enr-End}{endomorphism $\V$-operad of an object $x$}
\entry{$\EndC$}{ex:endc}{endomorphism multicategory of a permutative category}
\entry{$\VMulticat$}{v-multicat-2cat}{2-category of small $\V$-multicategories}
\entry{$\VMulticatns$}{not:vmulticatns}{2-category of small non-symmetric $\V$-multicategories}
\entry{$\bt$}{not:bt}{a terminal object in $\V$}
\entry{$\angx \to y$}{not:naryonecell}{an $n$-ary 1-cell}
\entry{$\Rightarrow$}{twocellmulticat}{an $n$-ary 2-cell}
\entry{$F\angc$}{Fangcthetaangc}{$(Fc_1,\ldots,Fc_k)$}
\entry{$\theta_{\angc}$}{Fangcthetaangc}{$(\theta_{c_1},\ldots,\theta_{c_k})$}
\entry{$1_{\theta_{\angc}}$}{catmultinaturalityiicell}{$(1_{\theta_{c_1}}, \ldots, 1_{\theta_{c_k}})$}
\entry{$\catmulticat$}{not:catmulticat}{2-category of small $\Cat$-multicategories}
\entry{$\catmulticatns$}{not:catmulticatns}{non-symmetric variant of $\catmulticat$}

\newchHeader{ch:pseudosymmetry}

\entry{$\left(F, \{F_{\sigma,\angc,c'}\}\right)$, $(F,\{F_\sigma\})$}{def:pseudosmultifunctor}{a pseudo symmetric $\Cat$-multifunctor}
\entry{$F_{\sigma,\angc,c'}$, $F_\sigma$}{Fsigmaangccp}{a pseudo symmetry isomorphism}
\entry{$\sigmabar$}{not:topeq}{block permutation $\sigma\langle k_{\sigma(1)}, \ldots , k_{\sigma(n)} \rangle$}
\entry{$\tautimes$}{not:boteq}{block sum $\tau_1 \times \cdots \times \tau_n$}
\entry{$F_{\sigma;p}$, $F_{\sigma,p}$}{not:Fsubsigmap}{$p$-component of $F_\sigma$}
\entry{$p^\sigma$, $p \cdot \sigma$}{not:ptosigma}{right $\sigma$-action}
\entry{$GF_\sigma$, $(GF)_\sigma$}{not:GFsigma}{a pseudo symmetry isomorphism of $GF$}
\entry{$\catmulticatps$}{thm:catmulticatps}{pseudo symmetric variant of $\catmulticat$}
\entry{$\catmulticatlax$}{not:catmulticatlax}{lax symmetric variant of $\catmulticat$}

\newchHeader{ch:diagram}

\entry{$\Dcat$, $\Dtecat$}{not:Dcat}{$\Cat$-multicategory of $\cD$-indexed categories with $\cD$ permutative}
\entry{$\sigma\anga$}{sigmaangaasigmainv}{$\ang{a_{\sigmainv(j)}}_{j=1}^n$}
\entry{$a_{\sigmainv}$}{sigmaangaasigmainv}{$\txotimes_{j=1}^n a_{\sigmainv(j)}$}
\entry{$(X,X^2,X^0)$}{def:additivesmf}{an additive symmetric monoidal functor}
\entry{$X^0*$}{additivesmfunitobj}{unit object of an additive symmetric monoidal functor}
\entry{$X^2_{a,b}$}{additivesmfconstraint}{$(a,b)$-component of the monoidal constraint $X^2$}
\entry{$\lap$}{laplazaapp}{Laplaza coherence isomorphism that distributes over a sum}
\entry{$\varrho$}{varrhofactor}{diagonal in all but one factors followed by a transpose permutation}
\entry{$\DCat$}{def:dcatcatmulticat}{$\Cat$-multicategory of additive symmetric monoidal functors}
\entry{$\DCat\smscmap{\angX;Z}$}{dcatangxz}{a subcategory of $\Dtecat\smscmap{\angX;Z}$}
\entry{$\DCat\smscmap{\ang{};Z}$}{dcatbipermzeroary}{$Z\tu$}
\entry{$U$}{dcatforgetdtecat}{forgetful $\Cat$-multifunctor $\DCat \to \Dtecat$}

\newchHeader{ch:multigro}

\entry{$\ang{x \compi x_i'}$, $\ang{x} \compi x_i'$}{compnotation}{$\angx$ with its $i$-th entry replaced by $x_i'$}
\entry{$\ang{x \compi x_i' \comp_\ell x'_\ell}$}{compcompnotation}{$\ang{x \compi x_i'}$ with its $\ell$-th entry replaced by $x'_\ell$}
\entry{$(F,\ang{F^2_j}_{j=1}^n)$, $(F,\ang{F^2_j})$}{nlinearfunctorangcd}{an $n$-linear functor}
\entry{$F^2_j$}{laxlinearityconstraints}{$j$-th linearity constraint of $F$}
\entry{$\permcat\smscmap{\angC;\D}$}{permcatangcd}{category of $n$-linear functors and natural transformations}
\entry{$\permcatsg\smscmap{\angC;\D}$}{permcatsgangcd}{strong variant of $\permcat\smscmap{\angC;\D}$}
\entry{$\permcat\smscmap{\ang{};\D}$}{permcatemptyd}{$\D$}
\entry{$\permcatsg\smscmap{\ang{};\D}$}{permcatemptyd}{$\D$}
\entry{$\permcat$}{thm:permcatenrmulticat}{$\Cat$-multicategory of small permutative categories}
\entry{$\permcatsg$}{not:permcatsgcatm}{strong variant of $\permcat$}
\entry{$\grod X$}{def:grothendieckconst}{Grothendieck construction of a functor $X \cn \cD \to \Cat$}
\entry{$(\zero, X^0*)$}{groconunit}{monoidal unit of $\grod X$}
\entry{$\gbox$}{gboxobjects}{monoidal product of $\grod X$}
\entry{$\betabox$}{groconbeta}{braiding of $\grod X$}
\entry{$\grod \phi$}{groadditivenattr}{Grothendieck construction of an additive natural transformation}
\entry{$(\grod \phi)^2_i$}{grodphilinearity}{$i$-th linearity constraint of $\grod \phi$}
\entry{$\grod \Phi$}{grodPhi}{Grothendieck construction of an additive modification}
\entry{$(\grod)_{\sigma}$}{gropseudosymmetry}{a pseudo symmetry isomorphism of $\grod$}
\entry{$(\grod)_{\sigma,\phi}$, $(\grod)_{\sigma,\phi,\ang{(a_j,x_j)}}$}{not:grodsigmaphi}{components of $(\grod)_\sigma$}
\entry{$\grod$}{grodpscatmultifunctor}{pseudo symmetric $\Cat$-multifunctor $\DCat \to \permcatsg$}

\newchHeader{ch:permfib}

\entry{$\fl{y}{f}$}{def:opfibration-i}{a fore-lift with object $y$ and morphism $f$}
\entry{$\fbar \cn y \to y_f$}{def:opfibration-ii}{a lift of a fore-lift $\fl{y}{f}$}
\entry{$\fr{g}{h}{f}$}{def:opfibration-iii}{a fore-raise}
\entry{$\ftil$}{def:opfibration-iv}{a raise of a fore-raise $\fr{g}{h}{f}$}
\entry{$P^\op \cn \E^\op \to \B^\op$}{expl:opfibterminology-i}{opposite functor of $P \cn \E \to \B$}
\entry{$\ang{P} \to P_0$}{Fangcczero}{an opcartesian $n$-linear functor}
\entry{$\pfibd\smscmap{\ang{P};P_0}$}{not:pfibdangppzero}{a subcategory of $\permcatsg\smscmap{\angC;\C_0}$}
\entry{$\pfibd$}{def:pfibdmulticat}{ns $\Cat$-multicategory of small permutative opfibrations over $\cD$}
\entry{$\tensorop$}{oppositetensor}{opposite monoidal product}
\entry{$\U$}{pfibdtopermcatsg}{non-symmetric $\Cat$-multifunctor $\pfibd \to \permcatsg$}
\entry{$U_X$}{Usubx}{first-factor projection $\grod X \to \cD$}
\entry{$\ang{\grod X}$}{not:anggrodx}{$\ang{\grod X_j}_{j=1}^n$}
\entry{$\ang{U_X}$}{not:anggrodx}{$\ang{U_{X_j}}_{j=1}^n$}
\entry{$\grod$}{groddcatpfibd}{non-symmetric $\Cat$-multifunctor $\DCat \to \pfibd$}

\newchHeader{ch:gromultequiv}

\entry{$P^\inv(a)$}{Xapinva}{preimage subcategory}
\entry{$F_2\angx$}{Fpreimageobject}{second entry in $F\ang{(a_j,x_j)}$}
\entry{$F_2\angp$}{Fangoneajpj}{second entry in $F\ang{(1_{a_j},p_j)}$}

\newchHeader{ch:multifunctorA}

\entry{$(\C,*)$, $(\C,\iota)$}{def:pointedcategory}{a pointed category $\C$ with basepoint $\iota \cn \boldone \to \C$}
\entry{$\pCat$}{not:pCat}{category of small pointed categories and pointed functors}
\entry{$E$}{smashunit}{smash unit $\boldone \bincoprod \boldone$}
\entry{$\C \wed \D$}{not:wedgeproduct}{wedge product}
\entry{$\C \sma \D$}{eq:smash-pushout}{smash product}
\entry{$\pr$}{eq:smash-pushout}{universal arrow $\C \times \D \to \C \sma \D$}
\entry{$[\C,\D]_*$}{eq:pHom}{pointed diagram category}
\entry{$\Fpunc\scmap{\ord{m};\ord{n}}$}{nonzeromorphism}{set of nonzero morphisms $\ord{m} \to \ord{n}$ in $\Fskel$}
\entry{$\Gacat$}{def:gammacategory}{pointed diagram category $[\Fskel,\Cat]_*$}
\entry{$\iota_n$}{not:iotan}{unique pointed function $\ord{0} \to \ord{n}$}
\entry{$\Fhat$}{def:Fhat}{$\pCat$-category with hom objects $\Fhat \smscmap{\ord{m};\ord{n}}$}
\entry{$\Fhat \smscmap{\ord{m};\ord{n}}$}{Fhathom}{pointed category $\bigvee_{\Fpunc\smscmap{\ord{m};\ord{n}}} E$}
\entry{$X \Gotimes Y$}{pointedday}{pointed Day convolution}
\entry{$\ftu$}{Gacatunit}{monoidal unit diagram $\Fhat \smscmap{\ord{1};-}$}
\entry{$\GHom(X,Y)$}{Fpointedhom}{pointed hom diagram}
\entry{$\xiday$}{eq:Fxi}{braiding for the pointed Day convolution}
\entry{$(L_{\ord{1}} \scs \ev_{\ord{1}})$}{not:Lordone}{evaluation at $\ord{1}$ and its left adjoint}
\entry{$\GMap(X,Y)$}{Fpointedmap}{pointed mapping object}
\entry{$\Gacat\smscmap{\angX;Z}$}{Gacatnmorphism}{a multimorphism category of $\Gacat$}
\entry{$\nu_j$}{iotajordoneordn}{morphism $\ord{1} \to \ord{n}$ in $\Fskel$ with $\nu_j(1) = j$}
\entry{$AX$}{functorAX}{functor $\cA \to \Cat$ associated to a $\Ga$-category $X$}
\entry{$(AX)(m)$}{AXmdef}{$\prod_{i=1}^p X\ord{m_i}$}
\entry{$\phi_*$}{phistarAphi}{$(AX)(\phi)$}
\entry{$\phi_{i,j}$}{phiijpointed}{pointed function induced by $\phi$}
\entry{$(AX,(AX)^2,(AX)^0)$}{AXsmf}{strictly unital strong symmetric monoidal structure of $AX$}
\entry{$AF$}{natAF}{monoidal natural transformation associated to $F$}
\entry{$(AF)_m$}{AFm}{$\prod_{i=1}^p F_{\ord{m_i}}$}
\entry{$A$}{AGacatACati}{functor $\Gacat \to \ACat$ of underlying categories}
\entry{$m^j$}{mjangmji}{$\bang{m^j_{i_j}}_{i_j=1}^{r_j}$}
\entry{$m^{1,\ldots,n}_{i_1,\ldots,i_n}$}{monenionein}{$\txprod_{j=1}^n m^j_{i_j}$}
\entry{$m^{1,\ldots,n}$}{monen}{$\txotimes_{j=1}^n m^j$}
\entry{$\ang{m}$}{angmangmj}{$\ang{m^j}_{j=1}^n$}
\entry{$(AF)_{\angm}$}{AFangm}{$\angm$-component functor of $AF$}
\entry{$\pr_?$}{AFangm}{coordinate projection}
\entry{$x^j$}{angxjij}{$\bang{x^j_{i_j}}_{i_j=1}^{r_j}$}
\entry{$A\theta$}{AFAangXAZ}{additive modification associated to $\theta$}
\entry{$(A\theta)_{\angm}$}{Athetaangm}{$\angm$-component natural transformation of $A\theta$}
\entry{$A$}{Acatmultifunctor}{$\Cat$-multifunctor $\Gacat \to \ACat$}

\newchHeader{ch:invK}

\entry{$\cP$}{mandellinvK}{inverse $K$-theory functor $\groa \circ A$}
\entry{$\cP X$}{mandellPX}{inverse $K$-theory of a $\Ga$-category $X$}
\entry{$(\ang{},*)$}{PXmonoidalunit}{monoidal unit in $\cP X$}
\entry{$\Box$}{PXmonoidalprod}{monoidal product in $\cP X$}
\entry{$\betabox$}{not:PXbetabox}{braiding in $\cP X$}
\entry{$\cP F$}{not:PF}{inverse $K$-theory of a $\Ga$-category morphism $F$}
\entry{$\cP$}{invKAgroa}{pseudo symmetric $\Cat$-multifunctorial inverse $K$-theory}
\entry{$\Asns$}{ex:Preservesnonsym}{non-symmetric associative operad}
\entry{$\As$}{Asoperad}{associative operad}
\entry{$\cP \theta$}{not:Ptheta}{inverse $K$-theory of a pointed modification $\theta$}
\entry{$\BE$}{def:barratteccles}{Barratt-Eccles operad}
\entry{$\sigmatil$}{sigmatilde}{block permutation induced by $\sigma$}
\entry{$(1,2)$}{not:onetwopermutation}{non-identity permutation in $\Sigma_2$}
\entry{$(X,\tu,\mu,\xi)$}{not:EinfGacat}{an $E_\infty$-algebra in $\Gacat$}
\entry{$\big(X, \{\mu_\sigma\}, \{\mu_\sigma^\tau\}, \{F_\sigma\} \big)$}{not:psEinfGacat}{a pseudo symmetric $E_\infty$-algebra in $\Gacat$}
\entry{$\big(\C, \{\mu_\sigma\}, \{\mu_\sigma^\tau\}, \{F_\sigma\} \big)$}{not:psEinfpermcat}{a pseudo symmetric $E_\infty$-algebra in $\permcat$}

\end{tabbing}

\printindex
\end{document}